\documentclass[11pt]{amsbook}
\usepackage{amscd}
\usepackage{amssymb}
\usepackage{hyperref}
\usepackage{mathrsfs}
\usepackage{verbatim}

\newtheorem{thm}{Theorem}[section]
\newtheorem{lem}[thm]{Lemma}

\newtheorem{cor}[thm]{Corollary}

\newtheorem{prop}[thm]{Proposition}

\newtheorem{mainthm}{Theorem}

\theoremstyle{definition}

\renewcommand{\thecase}{}

\newtheorem{conj}[thm]{Conjecture}
\newtheorem{mainconj}[mainthm]{Conjecture}

\newtheorem{defn}[thm]{Definition}
\newtheorem{exmp}[thm]{Example}
\newtheorem{hyp}[thm]{Hypothesis}

\newtheorem{notn}[thm]{Notation} 
\newtheorem{rmk}[thm]{Remark}
\renewcommand{\thestep}{}

\theoremstyle{remark}

\makeatletter
\def\alphenumi{
  \def\theenumi{\alph{enumi}}
  \def\p@enumi{\theenumi}
  \def\labelenumi{(\@alph\c@enumi)}}
\makeatother

\makeatletter
\def\thecase{\@arabic\c@case}
\makeatother

\numberwithin{equation}{section}
\numberwithin{section}{chapter}

\makeatletter
\def\thestep{\@arabic\c@step}
\makeatother

\renewcommand\emptyset{\varnothing}

\newcommand\embed{\hookrightarrow}

\newcommand\barM{{\bar{M}}}

\newcommand\barmu{{\bar\mu}}

\newcommand\ubarRR{{\underline{\mathbb{R}}}}

\newcommand\AAA{\mathbb{A}}

\newcommand\CC{\mathbb{C}}

\newcommand\FF{\mathbb{F}}

\newcommand\LL{\mathbb{L}}

\newcommand\NN{\mathbb{N}}

\newcommand\PP{\mathbb{P}}
\newcommand\QQ{\mathbb{Q}}
\newcommand\RR{\mathbb{R}}

\newcommand\ZZ{\mathbb{Z}}

\newcommand\beps{{\boldsymbol{\varepsilon}}}

\newcommand\bga{{\boldsymbol{\gamma}}}
\newcommand\bgamma{{\boldsymbol{\gamma}}}
\newcommand\bchi{{\boldsymbol{\chi}}}

\newcommand\bla{{\boldsymbol{\lambda}}}

\newcommand\bA{{\mathbf{A}}}

\newcommand\bB{{\mathbf{B}}}

\newcommand\bD{{\mathbf{D}}}

\newcommand\bh{{\mathbf{h}}}

\newcommand\bL{{\mathbf{L}}}

\newcommand\bQ{{\mathbf{Q}}}
\newcommand\br{{\mathbf{r}}}

\newcommand\bTheta{{\mathbf{\Theta}}}

\newcommand\bv{{\mathbf{v}}}

\newcommand\bx{{\mathbf{x}}}

\newcommand\by{{\mathbf{y}}}

\newcommand\bz{{\mathbf{z}}}

\newcommand{\cov}{\nabla}

\newcommand{\rd}{\partial}

\newcommand\thalf{{\textstyle{\frac{1}{2}}}}
\newcommand\tquarter{{\textstyle{\frac{1}{4}}}}

\newcommand\tthreequarter{{\textstyle{\frac{3}{4}}}}

\newcommand\teighth{{\textstyle{\frac{1}{8}}}}

\newcommand\half{{{\frac{1}{2}}}}

\newcommand\quarter{{{\frac{1}{4}}}}
\newcommand\threehalf{{{\frac{3}{2}}}}

\newcommand\fg{{\mathfrak{g}}}

\newcommand\fo{{\mathfrak{o}}}

\newcommand\fs{{\mathfrak{s}}}
\newcommand\fS{{\mathfrak{S}}}
\newcommand\ft{{\mathfrak{t}}}

\newcommand\fU{{\mathfrak{U}}}
\newcommand\fV{{\mathfrak{V}}}

\newcommand\eps{\varepsilon}
\newcommand\ga{\gamma}
\newcommand\Ga{\Gamma}
\newcommand\la{\lambda}
\newcommand\La{\Lambda}
\newcommand\ka{\kappa}
\newcommand\om{\omega}
\newcommand\Om{\Omega}
\newcommand\si{\sigma}
\newcommand\Si{\Sigma}

\newcommand\so{{\mathfrak{s}\mathfrak{o}}}
\newcommand\su{{\mathfrak{s}\mathfrak{u}}}
\newcommand\fu{{\mathfrak{u}}}

\newcommand\BG{\operatorname{BG}}

\newcommand\BS{\operatorname{BS}}
\newcommand\BSO{\operatorname{BSO}}

\newcommand\BSpinu{{\operatorname{BSpin}}^u}
\newcommand\BtH{\operatorname{B{\widetilde H}}}

\newcommand\BU{\operatorname{BU}}

\newcommand\EG{\operatorname{EG}}

\newcommand\ES{\operatorname{ES}}
\newcommand\ESO{\operatorname{ESO}}

\newcommand\ESpinu{{\operatorname{ESpin}}^u}

\newcommand\EtH{\operatorname{E{\widetilde H}}}

\newcommand\GL{\operatorname{GL}}

\newcommand\SO{\operatorname{SO}}

\newcommand\Spin{\operatorname{Spin}}

\newcommand\SU{\operatorname{SU}}
\newcommand\U{\operatorname{U}}

\newcommand\less{\setminus}

\newcommand{\8}{\infty}

\newcommand\ad{{\operatorname{ad}}}
\newcommand\Ad{{\operatorname{Ad}}}

\newcommand\Aut{\operatorname{Aut}}

\newcommand\CCl{\operatorname{{\mathbb{C}\ell}}}

\newcommand\Coker{\operatorname{Coker}}

\newcommand\dist{\operatorname{dist}}

\newcommand\End{\operatorname{End}}

\newcommand\Fr{\operatorname{Fr}}

\newcommand\Gl{\operatorname{Gl}}

\newcommand\Hom{\operatorname{Hom}}

\newcommand\ind{\operatorname{Index}}
\newcommand\Ind{\operatorname{Index}}
\newcommand\Int{\operatorname{int}}
\newcommand\Imag{\operatorname{Im}}
\newcommand\Isom{\operatorname{Isom}}

\newcommand\Ker{\operatorname{Ker}}

\newcommand\loc{\operatorname{loc}}
\newcommand\Map{\operatorname{Map}}

\newcommand\Ob{\operatorname{Ob}}

\newcommand\PD{\operatorname{PD}}

\newcommand\Ran{\operatorname{Ran}}
\newcommand\rank{\operatorname{rank}}

\newcommand\red{\operatorname{red}}

\newcommand\Stab{\operatorname{Stab}}

\newcommand\SW{SW}
\newcommand\Sym{\operatorname{Sym}}

\newcommand\Tor{\operatorname{Tor}}
\newcommand\Th{\operatorname{Th}}

\newcommand\Tr{\operatorname{Tr}}

\newcommand\vol{\operatorname{vol}}

\newcommand\asd{\textsc{asd}}

\newcommand\cl{{\mathrm{cl}}}

\newcommand\even{{\mathrm{even}}}

\newcommand\id{{\mathrm{id}}}

\newcommand\odd{{\mathrm{odd}}}

\newcommand\round{{\mathrm{round}}}
\newcommand\sing{\mathrm{sing}}

\newcommand\spinc{\text{$\mathrm{spin}^c$ }}
\newcommand\spinu{\text{$\mathrm{spin}^u$ }}
\newcommand\Spinc{\text{$\mathrm{Spin}^c$}}
\newcommand\Spinu{\text{$\mathrm{Spin}^u$}}
\newcommand\spl{\textsc{spl}}
\newcommand\std{{\mathrm{std}}}

\newcommand\vir{\mathrm{vir}}

\newcommand\sA{{\mathscr{A}}}
\newcommand\sB{{\mathscr{B}}}
\newcommand\sC{{\mathscr{C}}}

\newcommand\sF{{\mathscr{F}}}
\newcommand\sG{{\mathscr{G}}}

\newcommand\sI{{\mathscr{I}}}

\newcommand\sK{{\mathscr{K}}}
\newcommand\sL{{\mathscr{L}}}
\newcommand\sM{{\mathscr{M}}}

\newcommand\sO{{\mathscr{O}}}
\newcommand\sP{{\mathscr{P}}}

\newcommand\sS{{\mathscr{S}}}

\newcommand\sU{{\mathscr{U}}}
\newcommand\sV{{\mathscr{V}}}
\newcommand\sW{{\mathscr{W}}}
\newcommand\sX{{\mathscr{X}}}

\newcommand\tsC{{\tilde\sC}}

\newcommand\tG{{\widetilde G}}

\newcommand\tH{{\widetilde H}}

\newcommand\tM{{\tilde M}}
\newcommand\tN{{\tilde N}}

\newcommand\tXi{{\widetilde\Xi}}

\begin{document}

\frontmatter
\title[Cobordism formula]
{An SO(3)-monopole cobordism formula relating Donaldson and Seiberg--Witten
invariants}
\author[Paul M. N. Feehan]{Paul M. N. Feehan}
\author[Thomas G. Leness]{Thomas G. Leness}

\begin{abstract}
We prove an analogue of the Kotschick--Morgan Conjecture in the context of $\SO(3)$ monopoles, obtaining a formula relating the Donaldson and Seiberg--Witten invariants of smooth four-manifolds using the $\SO(3)$-monopole cobordism. The main technical difficulty in the $\SO(3)$-monopole program relating the Seiberg--Witten and Donaldson invariants has been to compute intersection pairings on links of strata of reducible $\SO(3)$ monopoles, namely the moduli spaces of Seiberg--Witten monopoles lying in lower-level strata of the Uhlenbeck compactification of the moduli space of $\SO(3)$ monopoles \cite{FL1}. In this monograph, we prove --- modulo a gluing theorem which is an extension of our earlier work in \cite{FL3} --- that these intersection pairings can be expressed in terms of topological data and Seiberg--Witten invariants of the four-manifold. Our proofs that the $\SO(3)$-monopole cobordism yields both the Superconformal Simple Type Conjecture of Moore, Mari\~no, and Peradze \cite{MMPhep,MMPdg} and Witten's Conjecture \cite{Witten} in full generality for all closed, oriented, smooth four-manifolds with $b_1=0$ and odd $b^+\ge 3$ appear in \cite{FL7, FL8}.
\end{abstract}

\dedicatory{Paul Feehan dedicates this monograph to the memory of his parents, Martin and Odilla Feehan
and to his wife, Julia Stanat, and daughter, Odilla Yseult Tyler Feehan.\\ Thomas Leness dedicates this monograph to
his wife, Kirsten Wood.}



\maketitle

\makeatletter
\newcommand\@dotsep{3.0}
\def\@tocline#1#2#3#4#5#6#7{\relax
  \ifnum #1>\c@tocdepth % then omit
  \else
    \par \addpenalty\@secpenalty\addvspace{#2}%
    \begingroup \hyphenpenalty\@M
    \@ifempty{#4}{%
      \@tempdima\csname r@tocindent\number#1\endcsname\relax
    }{%
      \@tempdima#4\relax
    }%
    \parindent\z@ \leftskip#3\relax \advance\leftskip\@tempdima\relax
    \rightskip\@pnumwidth plus1em \parfillskip-\@pnumwidth
    #5\leavevmode\hskip-\@tempdima #6\relax
    \leaders\hbox{$\m@th
      \mkern \@dotsep mu\hbox{.}\mkern \@dotsep mu$}\hfill
    \hbox to\@pnumwidth{\@tocpagenum{#7}}\par
    \nobreak
    \endgroup
  \fi}
\makeatother

\tableofcontents

\chapter*{Preface}
The $\SO(3)$-monopole cobordism formula is an equality between the Donaldson invariant of a closed, smooth four-manifold, $X$, and a universal expression involving the Seiberg--Witten invariants of $X$, the intersection form of $X$, and other homotopy invariants of $X$. In this monograph, we prove that this formula follows from certain properties of the gluing map constructed in \cite{FL3}.

Our proof that the $\SO(3)$-monopole cobordism formula implies Witten's Conjecture \cite{Witten} relating the Donaldson and Seiberg--Witten invariants for four-manifolds of simple type which  satisfy a `geography inequality' or which are `abundant' (as defined in \cite{FKLM}) appeared in \cite{FL6}.
Our proofs that the $\SO(3)$-monopole cobordism yields both the Superconformal Simple Type Conjecture of Moore, Mari\~no, and Peradze \cite{MMPhep,MMPdg} and Witten's Conjecture \cite{Witten} in full generality  for all closed, oriented, smooth four-manifolds with $b_1(X)=0$, odd $b^+(X)\ge 3$, and Seiberg--Witten simple type have appeared in \cite{FL7} and \cite{FL8}, respectively.

The monopole cobordism formula has an important feature in common with that conjectured by Kotschick and Morgan \cite{KotschickMorgan} concerning the wall-crossing property for Donaldson invariants: both formulae give an equality relating the invariants that involves unknown coefficients depending only on the homotopy type of the manifold.
G\"ottsche \cite{Goettsche} and G\"ottsche and Zagier \cite{GoettscheZagier} assumed the validity of the Kotschick--Morgan Conjecture to derive a wall-crossing formula, with all coefficients explicitly determined, and the resulting structure of Donaldson invariants for four-manifolds with
$b^+=1$.
We applied a broadly similar strategy in our articles \cite{FL7, FL8}, but refer interested readers to those articles for further details.

The main technical difficulty arising in the $\SO(3)$-monopole program to prove Witten's Conjecture has been to compute intersection pairings on the links of strata of reducible $\SO(3)$ monopoles, namely the moduli spaces of Seiberg--Witten monopoles lying in lower-level strata of the Uhlenbeck compactification of the moduli space of $\SO(3)$ monopoles \cite{FL1}. Our local gluing map \cite{FL3} parameterizes a neighborhood of one of these strata.  In this monograph, we solve the `overlap problem' described in \cite{FLMcMaster} and show how these local gluing maps fit together to describe a neighborhood of the union of these strata. This description allows us to prove that the desired intersection numbers can be expressed as a universal expression in the relevant Seiberg--Witten invariant, the intersection form of the manifold, and other homotopy invariants of the manifold appearing in the $\SO(3)$-monopole cobordism formula.

Because the $\SO(3)$-monopole cobordism formula has certain key features in common with that conjectured by Kotschick and Morgan \cite{KotschickMorgan}, not only in the form of the result but also in the nature of the intersection pairings to be computed, we can use the techniques developed herein to prove that the Kotschick--Morgan Conjecture also follows by similar arguments and the form of our gluing theorem for anti-self-dual connections \cite{FLKM1}.

Since the first version of this monograph was circulated, applications of our results have appeared in the proof of Property P for knots by Kronheimer and Mrowka \cite{KMPropertyP} and work of Sivek on Donaldson invariants for symplectic four-manifolds \cite{Sivek_2015imrn}. An alternative approach to Witten's Conjecture, inspired by results in physics, has been pursued by G\"ottsche, Nakajima, and Yoshioka \cite{Goettsche_Nakajima_Yoshioka_2009, Goettsche_Nakajima_Yoshioka_2008, Goettsche_Nakajima_Yoshioka_2011}, based in part on work of Mochizuki \cite{Mochizuki_2009}.

\section*{Acknowledgments}
We began work on the problem of giving a mathematically rigorous proof of Witten's Conjecture in 1995. We are most grateful for the longstanding support of Ronald Fintushel, Peter Kronheimer, and Cliff Taubes for our research in gauge theory and applications. We are extremely grateful to Tom Mrowka for his constant encouragement over the years and in particular to some very fruitful conversations where he urged us to develop the concept of `spliced-ends moduli spaces' of instantons over the four-sphere, $S^4$. He explained this concept to us in March 1999, and proposed it as a key idea in this part of our work on Witten's Conjecture and the Kotschick--Morgan Conjecture. Some of our research was conducted while visiting Columbia University, the Institute for Advanced Study, Princeton, the Max Planck Institute for Mathematics, Bonn, the Simons Center for Geometry and Physics, Stony Brook, Trinity College Dublin, and the University of Pennsylvania: we warmly thank those institutions for their kind hospitality. 
We are very grateful to the editorial and publishing staff of the American Mathematical Society for their assistance and patience during the preparation of our manuscript.
We thank the National Science Foundation for its support of our research, most recently via the grants DMS-1266145 (Feehan) and DMS-1510064 (Feehan), DMS-0905786 (Leness) and DMS-1510063 (Leness).

Paul Feehan pays tribute to the memory of his parents and their dedication, faith, love, and encouragement throughout their lives and thanks his wife, Julia Stanat, for her love and support.

Thomas Leness thanks his wife, Kirsten Wood, and his children, George and Anne, for their patience and support.
\bigskip

\rightline{December 10, 2018\footnote{This version: December 26, 2018, incorporating final galley proof corrections and corrections of minor typographical errors noticed since publication online. \emph{Memoirs of the American Mathematical Society} \textbf{256} (2018), no. 1226, \url{https://doi.org/10.1090/memo/1226}}}
\bigskip

\leftline{Paul M. N. Feehan}
\leftline{Department of Mathematics}
\leftline{Rutgers, The State University of New Jersey}
\leftline{Piscataway, NJ 08854-8019, United States}
\bigskip

\leftline{\texttt{feehan@math.rutgers.edu}}
\leftline{\url{math.rutgers.edu/~feehan}}

\bigskip

\leftline{Thomas G. Leness}
\leftline{Department of Mathematics}
\leftline{Florida International University}
\leftline{Miami, FL 33199, United States}
\bigskip

\leftline{\texttt{lenesst@fiu.edu}}
\leftline{\url{fiu.edu/~lenesst}}

\mainmatter

\chapter{Introduction}
\label{chap:Introduction}

\section{Summary of main results}
\label{sec:Summary_main_results}
In \cite{Witten}, Witten
defined the Seiberg--Witten invariants of smooth four-manifolds with $b_1(X)=0$ and odd $b^+(X)>0$ (based in part on earlier joint work with Seiberg \cite{SeibergWittenElecMag})
and stated a conjecture, based on arguments from quantum field
theory, relating the Seiberg--Witten and
Donaldson invariants.
Moore and Witten extended
this conjecture to a formula for four-manifolds with
$b_1(X)>0$ and $b^+(X)=0$ in \cite{MooreWitten}.
All computed examples
of Donaldson and Seiberg--Witten invariants satisfy these
conjectures.
However, these conjectures were
based on mathematically non-rigorous arguments from quantum field theory
and so the need for a mathematical explanation remained.

In \cite{PTLocal}, Pidstrigach and
Tyurin outlined an $\SO(3)$-monopole program
\label{Pidstrigach-Tyurin_SO(3)-monopole_cobordism}
with the goal of
giving a mathematically rigorous proof of Witten's Conjecture by
using a moduli space of $\SO(3)$ monopoles.
After the foundations of the $\SO(3)$-monopole program were
developed in \cite{FL1,FL2a,FL2b,FeehanGenericMetric,TelemanGenericMetric},
it soon became evident that the $\SO(3)$-monopole program exhibited a fundamental difficulty --- that of computing intersection numbers for links of certain
singularities in the lower strata of the Uhlenbeck compactification of
the moduli space of $\SO(3)$ monopoles.
This difficulty echoes that arising in attempts \cite{FLKM1, Goettsche} to prove the conjecture of Kotschick and Morgan \cite{KotschickMorgan}
\label{Introduction_Kotschick-Morgan_conjecture}
concerning wall-crossing
formulas for Donaldson invariants on four-manifolds with $b^+(X)=1$,
although there are many additional challenges in the case of the $\SO(3)$-monopole program.

In this monograph, we complete the topological
computation of these intersection numbers in the sense that we reduce
the computation to a \emph{local gluing theorem} which extends
that of \cite{FL3}, stated here as Hypothesis \ref{hyp:Gluing}.

\begin{mainthm}[$\SO(3)$-monopole cobordism formula]
\label{thm:MainThm}
Assume Hypothesis \ref{hyp:Gluing} holds.
Let $X$ be a closed, connected, oriented smooth four-manifold with $b_1(X)=0$, odd $b^+(X)>1$, Euler characteristic $\chi$, and signature $\si$. Let $\La,w\in H^2(X;\ZZ)$ obey $w-\La\equiv w_2(X)\pmod 2$. Let $\delta,m$ be non-negative integers for which $m\leq [\delta/2]$, where $[\,\cdot\,]$ denotes the greatest integer function, and $\delta\equiv -w^2-\frac{3}{4}(\chi+\sigma)\pmod{4}$, with
$\Lambda$ and $\delta$ obeying $\delta<i(\La)$, where $i(\La)=\La^2-\frac{1}{4}(3\chi+7\si)$. Then, for any $h\in H_2(X;\RR)$ and
generator $x\in H_0(X;\ZZ)$, one has the following expression for the
Donaldson invariant,
\begin{equation}
\label{eq:MainEquation}
\begin{aligned}
D^w_X(h^{\delta-2m}x^m)
&=
\sum_{\fs\in{\textrm{\em \Spinc}}(X)}
(-1)^{\frac{1}{4}(w-\La+c_1(\fs))^2}SW_X(\fs)
\\
&\qquad\times\sum_{i=0}^{ \min(\ell,[\delta/2]-m)}
\left(
p_{\delta,\ell,m,i}(c_1(\fs),\La)Q_X^i
\right)(h),
\end{aligned}
\end{equation}
where $Q_X$ is the intersection form on $H_2(X;\RR)$, and
$\ell=\frac{1}{4}(\delta+(c_1(\fs)-\La)^2+\frac{3}{4}(\chi+\si))$, and $p_{\delta,\ell,m,i}(\cdot,\cdot)$ is a
homogeneous polynomial of degree $\delta-2m-2i$ with coefficients which are universal functions of
$\chi$, $\sigma$, $c_1(\fs)^2$, $\La^2$, $c_1(\fs)\cdot\La$, $\delta$, $m$,
and $\ell$.
\end{mainthm}

Many significant results in low-dimensional topology can be proved using Theorem \ref{thm:MainThm} and related results in this monograph.  Foremost among these are our proof in \cite{FL8} of the Superconformal Simple Type Conjecture\label{Superconformal_simple_type_conjecture} of Moore, Mari\~no, and Peradze \cite{MMPhep,MMPdg}, our proof of Witten's Conjecture \cite{Witten} for a large class of four-manifolds in \cite[Corollary 7]{KMPropertyP} and \cite{FL6} and, finally, our proof in \cite{FL7} of Witten's Conjecture in full generality for all closed, oriented, smooth four-manifolds with $b_1(X)=0$, odd $b^+(X) > 1$, and Seiberg--Witten simple type.

In addition, Kronheimer and Mrowka apply our Theorem \ref{thm:MainThm} to give a proof of Property P for knots\label{Property_P_knots} in \cite{KMPropertyP}.
Recall that Property P is the statement that $+1$ surgery on a non-trivial
knot $K$ in $S^3$ yields a manifold which is not a homotopy sphere.
In \cite[Theorem 6]{KMPropertyP}, Kronheimer and Mrowka employ
Theorem \ref{thm:MainThm} to prove that Witten's Conjecture holds
for a large family of manifolds\footnote{There is a slight difference between Theorem \ref{thm:MainThm} as presented here and the version used in \cite{KMPropertyP}. Specifically, the bound on the degree of the Donaldson invariant used in \cite[Theorem 6]{KMPropertyP} uses $i(\La)=-\frac{1}{4}(\chi+\si)$
instead of the correct value of $i(\La)=-\frac{1}{4}(3\chi+7\si)$.
We believe that this difference is due to a typographical error in an earlier draft of this monograph.  The version given here is correct, but the error in no
way effects any of the results in \cite{KMPropertyP}.}.
They then argue that a counterexample to Property P would allow them to construct a four-manifold with non-trivial Seiberg--Witten invariants but trivial Donaldson invariants. As such a four-manifold would contradict Theorem \ref{thm:MainThm},
there can be no counterexample to Property P.

Another result following from the methods in this monograph is the `Multiplicity Conjecture', \cite[Conjecture 3.1]{FKLM}, stated here as Theorem \ref{thm:Multiplicity}. Recall from \cite{FKLM} that a four-manifold is `abundant' if the orthogonal complement of the Seiberg--Witten basic classes, with respect to the intersection form, contains a hyperbolic summand. In collaboration with Kronheimer and Mrowka, we showed in \cite{FKLM} that, for abundant four-manifolds, this conjecture verified the ideas of Mari{\~n}o, Moore, and Peradze from \cite{MMPdg,MMPhep} in which they proposed a lower bound on the multiplicity of the vanishing of the Seiberg--Witten series at zero in terms of $c(X)=-\frac{1}{4}(7\chi+11\si)$ and from this bound derived a lower bound on the number of basic classes of a four-manifold also in terms of $c(X)$.

Sivek has applied Theorem \ref{thm:MainThm} to show that symplectic four-manifolds with $b_1=0$ and odd $b^+>1$ have non-vanishing Donaldson invariants, and that the canonical class is always a Kronheimer--Mrowka basic class \cite{Sivek_2015imrn}.

Finally, the Kotschick--Morgan Conjecture \cite{KotschickMorgan} for the wall-crossing formulas for Donaldson invariants of manifolds with $b^+=1$ also follows from the results in this monograph as described in Section \ref{subsec:IntroKMConj}.

Witten's Conjecture for the relation between the Donaldson and Seiberg--Witten invariants is stated below.  The definitions of all terms in its statement appear in Sections \ref{sec:Donaldsonseries} and \ref{subsec:SWInvariants}.

\begin{mainconj}[Witten's Conjecture]
\label{conj:WC}
\cite{Witten}
Let $X$ be a closed, connected, oriented, smooth four-manifold with $b_1(X)=0$ and odd $b^+(X)>1$.  Assume that $X$ has Seiberg--Witten simple type. Then $X$ has Kronheimer--Mrowka simple type, the Kronheimer--Mrowka basic classes coincide with the Seiberg--Witten basic classes, and the Donaldson series of $X$ is given by
\begin{equation}
\bD^w_X(h)= 2^{2-c(X)}e^{Q_X(h)/2}\sum_{\fs} (-1)^{\frac{1}{2}(w^2+w\cdot c_1(\fs))}
\SW_X(\fs)e^{\langle c_1(\fs),h\rangle},
\end{equation}
where $c(X)=-\frac{1}{4}(7\chi+11\si)$.
\end{mainconj}

Detailed introductions of the $\SO(3)$-monopole program have appeared in
\cite{FLGeorgia,FLGeorgiaII,FLMcMaster}.  We now give a brief summary
of these ideas here.
Recall that a \spinc structure
\label{Introduction_spinc_structure}
on $X$ is determined by a complex rank-four vector bundle, $W\to X$, and a Clifford multiplication map, $\rho:T^*X\to \Hom(W)$. A \spinu structure\label{Introduction_spinu_structure} $\ft$ on $X$, as defined in \cite{FL2a} or Section \ref{subsubsec:SpincuStr} here, is given by $(\rho_V,V)$, where $V=W\otimes E$, and $\rho_V=\rho\otimes\id_E$, and $(\rho,W)$ is a \spinc structure, and $E\to X$ is a complex rank-two vector bundle,

An $\SO(3)$ monopole for a \spinu structure, $\ft=(\rho_V,W\otimes E)$, is a pair,
$(A,\Phi)$, where $A$ is a unitary connection on $E$ and $\Phi$ is a section of $V^+:=W^+\otimes E$, satisfying equations which can be thought of as a higher-rank version of the Seiberg--Witten equations, and $W = W^+\oplus W^-$ is the usual splitting of $W$ into the complex rank-two vector bundles of positive and negative spinors. The moduli space, $\sM_{\ft}$, is the space of $\SO(3)$ monopoles
\label{Introduction_moduli_space_SO(3)_monopoles}
modulo gauge equivalence. For generic perturbations, $\sM_{\ft}$ is a smooth manifold away from two types of singular subspaces which are also identified as fixed points of an $S^1$ action on $\sM_{\ft}$. The first type of singular subspace is identified with the moduli space of anti-self-dual connections \cite[Equation (3.5)]{FL2a}. The second type of singular subspace, that of reducible $\SO(3)$ monopoles, is identified in \cite[Lemma 3.13]{FL2a} with the
Seiberg--Witten moduli space\label{Introduction_moduli_space_Seiberg-Witten_monopoles}, $M_{\fs}$, where the \spinu structure $\ft$ admits a splitting $\ft=\fs\oplus \fs\otimes L$.  The $\SO(3)$-monopole program aims to use $\sM_{\ft}$ as a cobordism between the links of these singularities. Because $\sM_{\ft}$ is not compact in general, this cobordism does not provide any useful homological information.

The moduli space, $\sM_{\ft}$, admits an Uhlenbeck compactification,
\label{Introduction_Uhlenbeck_compactification} $\bar\sM_{\ft}$, which is contained in the space
of ideal monopoles,
$$
\bigcup_{\ell=0}^N\ \sM_{\ft(\ell)}\times \Sym^\ell(X),
$$
where $\sM_{\ft(\ell)}$ is a moduli space of $\SO(3)$ monopoles for a `lower charge' \spinu structure, $\ft(\ell)$, and $\Sym^\ell(X)$ is the $\ell$-th symmetric product of $X$. The space, $\bar\sM_{\ft}/S^1$, provides a compact and oriented cobordism between a link of the moduli space of anti-self-dual connections and links, $\bar\bL_{\ft,\fs}$, of subspaces of reducible $\SO(3)$ monopoles. However, additional reducible $\SO(3)$ monopoles, in the form of subspaces of ideal Seiberg--Witten monopoles,
\begin{equation}
\label{eq:IntroLowerLevelReducibles}
M_{\fs}\times\Sym^\ell(X) \subset \sM_{\ft(\ell)}\times \Sym^\ell(X),
\end{equation}
appear in the lower levels of the Uhlenbeck compactification.  The intersection number of certain geometric representatives of cohomology classes with the link of the anti-self-dual connections yields a multiple of the Donaldson invariant.  The main technical result of this monograph, Theorem \ref{thm:LinkPairing}, is a formula for the intersection number of these geometric representatives with the link, $\bar\bL_{\ft,\fs}$, of the subspace \eqref{eq:IntroLowerLevelReducibles} of $\bar\sM_{\ft}/S^1$. This formula expresses the intersection number in terms of
universal functions which, while not explicit, do not depend on the manifold, $X$.
The cobordism provided by $\bar\sM_{\ft}/S^1$ then yields the equality \eqref{eq:MainEquation} between the aforementioned multiple of the Donaldson invariant and the sum, over \spinc structures $\fs$, of these intersection numbers.

We note that a generalization of Theorem \ref{thm:LinkPairing} can also
be proved by the methods of this monograph, where the assumptions that $b_1(X)=0$ or $z=h^{\delta-2m}x^m$ are relaxed, but the resulting expression becomes considerably more complicated.

\section{Outline of the argument}
As described above, the problem we address in this monograph
is the computation of the intersection numbers,
\begin{equation}
\label{eq:IntroIntersectionExpression}
\#
\left(
\bar\sV(z)\cap\bar\sW^{\eta}\cap \bar\bL_{\ft,\fs}
\right),
\end{equation}
where $\bar\sV(z)$ and $\bar\sW$ are the geometric representatives
mentioned above and $\bar\bL_{\ft,\fs}$ is the link in
$\bar\sM_{\ft}/S^1$ of the Seiberg--Witten
singularities,
$M_{\fs}\times\Sym^\ell(X)$, appearing in \eqref{eq:IntroLowerLevelReducibles}.
A summary of our approach to this computation has appeared
in \cite{FLMcMaster}.  We give a short
sketch of our method here.

For the simplest case, when $M_{\fs}\subset\sM_{\ft}$, computations of the intersection number \eqref{eq:IntroIntersectionExpression} appeared in \cite{FL2b} (where $\ell=0$)
while computations of \eqref{eq:IntroIntersectionExpression} for singularities of the form $M_{\fs}\times X$ appeared in \cite{FLLevelOne} (where $\ell=1$).   While it should be possible to adapt the techniques of \cite{LenessWC} to directly compute \eqref{eq:IntroIntersectionExpression} for singularities of the form $M_{\fs}\times \Sym^2(X)$,
the difficulty increases rapidly as $\ell$ becomes larger and direct calculations appear intractable when $\ell \geq 3$.

Our gluing result \cite[Theorem 1.1]{FL3} can be used to define \emph{local gluing maps}
which parameterize neighborhoods of the strata,
\[
\left( M_{\fs}\times\Si \right)\cap \bar\sM_{\ft}
\subset
\left( M_{\fs}\times\Sym^\ell(X) \right)\cap \bar\sM_{\ft},
\]
in $\bar\sM_{\ft}$, where $\Si$ is a stratum of the symmetric product, $\Sym^\ell(X)$.
A local gluing map is then a composition of a \emph{local splicing map} and a certain perturbation of the image of the splicing map,
called the \emph{solution map}.
The local splicing map is defined, roughly, by patching together solutions of the $\SO(3)$-monopole equations over
$X$ and over the four-sphere, $S^4$, with its standard round metric of radius one using cut-off functions.
Let $[A_0,\Phi_0]\in M_{\fs}$ be a gauge-equivalence class of $\SO(3)$ monopoles over $X$ and let $\bx\in\Si \subset \Sym^\ell(X)$ be represented by a set of distinct points with multiplicity, $\{x_1,\ldots,x_r\}$, for $1 \leq r \leq \ell$. We refer to $(A_0,\Phi_0)$ as the \emph{background pair} and $x_i$ as the \emph{splicing points}.
Let $(A_i,0)$ for $i=1,\ldots,r$ be $\SO(3)$ monopoles over $S^4$; such solutions are given by an anti-self-dual connection, $A_i$, on an $\SU(2)$-bundle over $S^4$,
with second Chern (or `instanton') number equal to the multiplicity of $x_i$,
and the zero section.
Then, using local trivializations of the bundle supporting the pair $(A_0,\Phi_0)$
near the points $x_i$ and cut-off functions, one can patch together the $\SO(3)$-monopoles, $(A_0,\Phi_0)$ and $(A_i,0)$, to form a pair, $(A',\Phi')$, which is equal to $(A_0,\Phi_0)$ away from the points $x_i$ and equal to $(A_i,0)$ near $x_i$. We denote this pair schematically by
\begin{equation}
\label{eq:IntroSplicingMap}
(A',\Phi')
=
(A_0,\Phi_0)\#_{x_1}(A_1,0)\#_{x_2}(A_2,0)\#\cdots\#_{x_r}(A_r,0).
\end{equation}
This construction of the pair $(A',\Phi')$ is used to define the splicing map. Roughly speaking, the domain of this splicing map is then given by a fiber bundle,
\begin{equation}
\label{eq:IntroGluingFiberBundle}
\Gl(\Si)=\Fr(\Si)\times_{G(\Si)}M(\Si)
\to
M_{\fs}\times \Si,
\end{equation}
with a well-understood principal bundle, $\Fr(\Si)$, and structure group, $G(\Si)$.  The fiber, $M(\Si)$, is a product of moduli spaces of anti-self-dual connections over $S^4$.
In the preceding description of the spliced pair, $(A',\Phi')$, the points $[A_0,\Phi_0]$ and $\bx$ lie in the base of the bundle \eqref{eq:IntroGluingFiberBundle} and the points $[A_i,0]$ lie in the fiber, $M(\Si)$. The bundle $\Fr(\Si)$ arises from the trivializations of the bundles used in the construction of $(A',\Phi')$.

The pair $(A',\Phi')$ is only approximately a solution to the $\SO(3)$-monopole equations. There is a section, $\bchi_\Si$, of an \emph{obstruction pseudo-bundle}\footnote{See Definition \ref{defn:PseudoBundle} for an explanation of the term `pseudo-bundle'.},
$\Upsilon_\Si\to\Gl(\Si)$, whose fibers are finite-dimensional vector spaces, such that the
\emph{solution map} defined by
$$
(A',\Phi')\mapsto (A',\Phi')+(a,\phi),
$$
where $(a,\phi):=\wp(A',\Phi')$ is the solution to a system of partial differential equations given by
$$
\fS((A',\Phi')+\wp(A',\Phi')) \in \Upsilon|_{(A',\Phi')},
$$
where $(A,\Phi) \mapsto \fS(A,\Phi)$ is the map defined by the expression on the left-hand side of $\SO(3)$-monopole equations \eqref{eq:PerturbedSO3MonopoleEquations}.
The \emph{local obstruction section} is defined by
\[
\bchi_\Si(A',\Phi')=\fS((A',\Phi')+\wp(A',\Phi')),
\]
and hence the
solution map
identifies the zero locus of the obstruction section, $\bchi_\Si^{-1}(0)$, with a subset of $\sM_{\ft}$.

There are two main sources of difficulty in using the parametrization of a neighborhood of $M_{\fs}\times\Si$ provided by the gluing map in computing the intersection pairings \eqref{eq:IntroIntersectionExpression} of singularities of the form \eqref{eq:IntroLowerLevelReducibles} when $\ell>1$. The first difficulty is the presence of `higher-charge' moduli spaces of anti-self-dual connections on $S^4$ in the fiber
$M(\Si)$.  It would be an interesting question to see if the
work of \cite{NakajimaInstLect,NakajimaInstCountI,LiQinWang}
could be adapted to give the necessary information for the $G(\Si)$-equivariant cohomology ring of $M(\Si)$, but one which is beyond the scope of this monograph.  Such computations are carried out for the case $\ell=2$ in \cite{LenessWC}.
We do not address that problem in this monograph.
Instead, we use a pushforward-pullback argument (see Section \ref{sec:Computations}) to isolate the topology of these fibers as universal polynomials of the type appearing in \eqref{eq:MainEquation}.
The second difficulty, which we do address in this monograph, is the `overlap problem' arising from the presence of more than one stratum $\Si\subset\Sym^\ell(X)$ and the resulting need for more than one gluing map to parameterize a neighborhood of $M_{\fs}\times\Sym^\ell(X)$.

\subsection{Problem of overlaps}
\label{subsec:Overlaps}
When $\ell\ge 2$, the singularity, $M_{\fs}\times\Sym^\ell(X)$, has more than one stratum.  Hence, to compute the intersection pairing \eqref{eq:IntroIntersectionExpression} when $\ell\ge 2$, we must understand the overlap of the images of the gluing maps with domains given by the fiber bundles $\Gl(\Si)$ and $\Gl(\Si')$ of
\eqref{eq:IntroGluingFiberBundle} for different strata $\Si$ and $\Si'$ of $\Sym^\ell(X)$.
Merely adding up the intersection numbers in the open sets parameterized by each gluing map could yield the wrong answer because these open sets overlap and we might be counting each intersection point more than once.  Examples of this type of problem have appeared in \cite{OzsvathBlowUp,LenessBlowUp,LenessWC,WW}.

Our approach is to describe the overlaps of the images of the gluing maps
(compositions of splicing and solution maps)
by defining so-called \emph{crude splicing maps}\label{intro_crude_splicing_maps}. Because we can write the splicing map explicitly while the
solution map
is defined by the Implicit Function Theorem,
it is easier to compare two splicing maps than to compare two gluing maps.
Because the $\SO(3)$-monopole gluing maps defined in \cite{FL3} only identify the zero locus of the obstruction section in $\Gl(\Si)$ with an open subspace of $\sM_{\ft}$, the intersection of the images of $\Gl(\Si)$ and $\Gl(\Si')$ under the gluing map could be quite complicated and need not be given by an open subspace of $\Gl(\Si)$ or $\Gl(\Si')$ and need not share any of the fiber bundle properties which we wish to use in our computation. Instead, using a deformation of the splicing map and of the fiber
of the bundle $\Gl(\Si)$ to construct the crude splicing map, we are able to ensure that the overlap of the images of two splicing maps is a subbundle of each image.

Given two strata $\Si$ and $\Si'$ of $\Sym^{\ell}(X)$, with $\Si\subset \cl\, \Si'$
(where $\cl\,\Sigma$ denotes the closure in $\Sym^{\ell}(X)$ of a stratum $\Sigma \subset \Sym^{\ell}(X)$), and crude splicing maps
$$
\bga''_{\Si'}:\Gl(\Si')\to \bar\sC_{\ft},
\quad
\bga''_{\Si}:\Gl(\Si)\to \bar\sC_{\ft},
$$
we control the overlap of their images by defining a space of \emph{overlap data},
$\Gl(\Si,\Si')$, and maps
$$
\rho^u_{\Si,\Si'}: \Gl(\Si,\Si')\to \Gl(\Si'),
\quad
\rho^d_{\Si,\Si'}: \Gl(\Si,\Si')\to \Gl(\Si),
$$
such that the diagram
\begin{equation}
\label{eq:IntroOverlapCD}
\begin{CD}
\Gl(\Si,\Si') @> \rho^u_{\Si,\Si'} >> \Gl(\Si')
\\
@V \rho^d_{\Si,\Si'} VV @V \bga''_{\Si'} VV
\\
\Gl(\Si) @> \bga''_{\Si} >> \bar\sC_{\ft}
\end{CD}
\end{equation}
commutes and such that
$$
\Imag(\bga''_{\Si'})\cap\Imag(\bga''_{\Si})
=
\Imag(\bga''_{\Si'}\circ\rho^u_{\Si,\Si'})
=
\Imag(\bga''_{\Si}\circ\rho^d_{\Si,\Si'}).
$$
Moreover, the upwards and downwards overlap maps, $\rho^u_{\Si,\Si'}$ and $\rho^d_{\Si,\Si'}$, are both fiber bundle maps. The diagram \eqref{eq:IntroOverlapCD} is then used to define the space of {\em global splicing data\/} as a pushout of the spaces of local splicing data, $\Gl(\Si)$. Once this is accomplished, the rest of the computation of the intersection number \eqref{eq:IntroIntersectionExpression} is largely a formal accounting of choices of cohomology classes with compact support and a use of the familiar pullback-pushforward technique.

\subsection{Overlap space and overlap maps}
\label{subsec:IntroOverlapSpaceAndMaps}
We now sketch the construction of the overlap space, $\Gl(\Si,\Si')$, and overlap maps, $\rho^u_{\Si,\Si'}$ and $\rho^d_{\Si,\Si'}$. First, we note that the overlap of the images will be empty unless $\Si\subset\cl\, \Si'$ or
$\Si'\subset\cl\, \Si$.  We will assume that the first case holds,
thus $\Si\subset\cl\, \Si'$, and will refer to $\Si$ as the \emph{lower} stratum.

There is an open neighborhood,
$\nu(\Si,\Si') \subset \Sigma'$, of $\Si$
which can be thought of as
a tubular neighborhood
of $\Si$ in $\Si'$.  We then define the space of \emph{overlap data}, $\Gl(\Si,\Si')$, to be the restriction of the bundle $\Gl(\Si')$ in \eqref{eq:IntroGluingFiberBundle} to
$$
M_{\fs}\times\nu(\Si,\Si')\subset M_{\fs}\times \Si'.
$$
The map $\rho^u_{\Si,\Si'}$ is the inclusion of the bundles.  The definition of $\rho^d_{\Si,\Si'}$ will require a redefinition of the fibers $M(\Si)$.

The \emph{overlap map}, $\rho^d_{\Si,\Si'}$, will map the fiber of the composition,
\begin{equation}
\label{eq:IntroOverlapLowerFiber}
\Gl(\Si,\Si')\to M_{\fs}\times \nu(\Si,\Si')\to M_{\fs}\times\Si,
\end{equation}
to the fiber of the projection $\Gl(\Si)\to M_{\fs}\times\Si$ from \eqref{eq:IntroGluingFiberBundle}. The fiber of the map $\nu(\Si,\Si')\to\Si$ is, up to a local trivialization, a collection of points in $\RR^4$.  As an example, the  fiber of the normal bundle of the diagonal $\Delta\subset X^2$ is, after specifying a trivialization of the tangent bundle of $X$ at the point given by the point in $\Delta$, a pair of points
in $\RR^4$ with a mass-centering condition.  The fiber of the composition \eqref{eq:IntroOverlapLowerFiber} is then given by $M(\Si')$ and a collection of points in $\RR^4$. Hence, the map $\rho^d_{\Si,\Si'}$ must take the points in $\RR^4$ and the connections in $M(\Si')$ and map them to elements of $M(\Si)$.

The strata of $\Sym^\ell(X)$ are in bijective correspondence with partitions of $\ell$.  A stratum $\Si$ is \emph{lower} than a stratum $\Si'$ if and only if the partition determining $\Si'$ is a refinement of the partition determining $\Si$. By \emph{refinement}, we mean that if $\Si$ is given by a partition $\ell=\ka_1+\dots+\ka_s$, then $\Si'$ is given by a partition $\sum_{i,j}\ka_{i,j}$, where $\ka_i=\sum_j \ka_{i,j}$. If
$\Sigma$ and $\Sigma'$
are given by the partitions above, then a point in the fiber of $\nu(\Si,\Si')\to \Si$ over $\bx\in\Si$ consists of clusters of points $y_{i,j}\in X$ near each point $x_i\in X$ defining $\bx$. The points $y_{i,j}$ correspond to the weights $\ka_{i,j}$ making up the partition $\sum_j \ka_{i,j}$ of $\ka_i$. With a trivialization of the tangent bundle at each $x_i$, the points $y_{i,j}$ can be considered as points  in $\RR^4$.  Consequently, the data in the fiber of \eqref{eq:IntroOverlapLowerFiber}, points in $\RR^4$, and connections in $M(\Si')$, comprise the domain of a gluing map parameterizing a neighborhood
(in the Uhlenbeck compactification of a moduli space of framed, mass-centered, anti-self-dual connections on $S^4$) of the points,
\begin{equation}
\label{eq:IntroTrivialStrata}
[\Theta]\times\Si(\ka_{i,j}), 
\end{equation}
where $\Theta$ is the
product connection and 
$$
\Si(\ka_{i,j})\subset\Sym^{\ka_i}(\RR^4)-\{c\},\quad\text{where } c=[0,\ldots,0],
$$
is the stratum given by the partition $\kappa_i = \sum_j \kappa_{i,j}$.
We refer to the strata \eqref{eq:IntroTrivialStrata} as \emph{product connection strata}.
Hence, it would seem natural to choose $\rho^d_{\Si,\Si'}$ to be the gluing map and identify the data in the fiber of \eqref{eq:IntroOverlapLowerFiber} with an element of $M(\Si)$.
However, the above definition will \emph{not} make the diagram \eqref{eq:IntroOverlapCD} commute.

\subsection{Associativity of splicing maps}
\label{subsec:IntroAssocSplic}
The commutativity of the diagram \eqref{eq:IntroOverlapCD} can be understood as an associativity property
of the splicing construction as follows.
Analogous properties of gluing maps for gradient flow trajectories are discussed in \cite[Theorem 2.1]{CohenJonesSegal_1995}, \cite{LizhengQin_2011}, and \cite{Wehrheim_2012}.

Let $\Si$ and $\Si'$ be two strata of $\Sym^\ka(\RR^4)$. Strata of such a symmetric product are specified by
partitions of $\ka$.  If $\Si$ is the lower stratum, that is, $\Si\subset\cl\, \Si'$, then the partition defining $\Si'$ is a refinement of the partition $\ka=\ka_1+\dots +\ka_s$ defining $\Si$ in the sense that  $\ka_i=\ka_{i,1}+\cdots+\ka_{i,r_i}$, where $\Si'$ is given by the partition $\ka=\sum_{i,j}\ka_{i,j}$.

Let $(A_0,\Phi_0)$ be the background pair and $A_{i,j}$ be the connections over $S^4$ given by the point in the fiber $M(\Si')$. Let $\bx\in \Si$ be given by the set of points $\{x_1,\dots,x_r\}$. Let  $\bz\in\nu(\Si,\Si')\subset\Si'$ lying over $\bx$ be given by points $z_{i,j}\in X$ corresponding to points $y_{i,j}\in\RR^4$. Let $\Theta$ be the product connection on $S^4$ and identify the points $y_{i,j}$ with points in $S^4$ by stereographic projection. We define connections with charge $\ka_i$ on $S^4$ by applying the splicing construction of \eqref{eq:IntroSplicingMap} just to the connection,
\begin{equation}
\label{eq:IntroConnSplicing}
A_i'=
\Theta\#_{y_{i,1}}A_{i,1}\#_{y_{i,2}}A_{i,2}\#\cdots \#_{y_{i,r_i}}A_{i,r_i}.
\end{equation}
If we consider the splicing construction as
defining a binary operation and the pair $(\Theta,0)$ as a unit for this binary operation, then the following can be seen as an associativity equality.
\begin{equation}
\label{eq:IntroSplicingAssoc}
\begin{aligned}
{}&
(A_0,\Phi_0)\#_{x_1}(A_1',0)\#_{x_2}(A_2',0)\#\cdots\#_{x_r}(A_r',0)
\\
{}&\quad
=
(A_0,\Phi_0)\#_{z_{1,1}}(A_{1,1},0)\#\cdots\#_{z_{1,r_1}}(A_{1,r_1},0)
\#\cdots
\#_{z_{s,1}}A_{s,1}\#\cdots\#_{z_{s,r_s}}(A_{s,r_s},0).
\end{aligned}
\end{equation}
We wish to define the map $\rho^d_{\Si,\Si}$ in the diagram \eqref{eq:IntroOverlapCD} by the construction
of the connections $A_i'$.
We can see as follows that the diagram \eqref{eq:IntroOverlapCD} will commute if the equality \eqref{eq:IntroSplicingAssoc} holds. The composition $\bga_{\Si'}''\circ\rho^u_{\Si,\Si'}$
corresponds to the right-hand-side of the equality \eqref{eq:IntroSplicingAssoc} and is equal to
the map obtained by
splicing the pairs at the
points in $X$ defined by an element of
$\nu(\Si,\Si')$.  The composition $\bga_{\Si}''\circ\rho^d_{\Si,\Si'}$ corresponds to the left-hand-side of the equality \eqref{eq:IntroSplicingAssoc} and is equal to the composition of two splicing maps.

There are two problems with the preceding strategy.
First, the connections $A_i'$ are not anti-self-dual and thus
do not define an element of $M(\Si)$.  Second, the equality \eqref{eq:IntroSplicingAssoc} does not hold in general.
In the sequel, we describe how we correct these problems.

\subsection{Instanton moduli space with spliced ends}
\label{subsec:Introduction_spliced-ends_moduli_space}
To overcome the first difficulty mentioned above, we shall change the fiber $M(\Si)$ as follows.
Recall that $M^s_\ka(S^4)$ denotes the moduli space of charge-$\ka$, mass-centered, anti-self-dual, framed $\SU(2)$-connections over $S^4$ and that $\bar M^s_\ka(S^4)$ denotes its Uhlenbeck compactification. We shall redefine the fiber,
$$
M(\Si)=\prod_{i=1}^r \bar M^s_{\ka_i}(S^4),
$$
by replacing $\bar M^s_{\ka_i}(S^4)$ in the preceding product with an \emph{instanton moduli space with spliced ends}, $\bar M^s_{\spl,\ka_i}$, whose construction we now indicate.

Recall that for each stratum, $\Si(\ka_{i,j})\subset\Sym^{\ka_i}(\RR^4)-\{c\}$, there is a gluing map,
\begin{equation}
\label{eq:IntroGluingToTrivial}
\prod_{j} \bar M^s_{\ka_{i,j}}(S^4)\times \Si(\ka_{i,j})
\to
\bar M^s_{\ka_i}(S^4),
\end{equation}
parameterizing a neighborhood of the stratum \eqref{eq:IntroTrivialStrata} in $\bar M^s_{\ka_i}(S^4)$. We define the instanton moduli space with spliced ends, $\bar M^s_{\spl, \ka_i}$, using induction on $\ka_i$ by replacing the image of the above gluing map by the image of the corresponding splicing map with domain,
\[
\bar M^s_{\spl,\ka_{i,j}}\times \Si(\ka_{i,j}).
\]
We now describe how the images of these splicing maps fit together with each other and with the rest of $\bar M^s_{\ka_i}(S^4)$ to form a smoothly stratified space.
By \emph{smoothly stratified space}, we mean a stratified space in the sense of \cite[Definition 3.4]{TauL2} where each stratum is a smooth manifold.
Such a space is called a Whitney pre-stratified space satisfying the condition of the frontier in \cite[p. 480]{Mather_2012}.

For these splicing maps, the background connection is
the product connection, $\Theta$, and the metric is flat.  Under these conditions, the analogue of the associative equality \eqref{eq:IntroSplicingAssoc} for splicing connections holds and the analogue of the diagram \eqref{eq:IntroOverlapCD} commutes. This implies that the union, over strata in $\Sym^{\ka_i}(\RR^4)-\{c\}$, of the images of these splicing maps forms a smoothly stratified space $W(\ka_i)$ whose image, $W'(\ka_i)$, under the gluing deformation is a neighborhood of $\{\Theta\}\times\Sym^{\ka_i}(\RR^4)-\{c\}$ in $\bar M^s_{\ka_i}(S^4)$. The gluing deformation defines a smoothly stratified isotopy $R_t$, $0\le t\le 1$, in the sense of 
\cite[p. 482]{Mather_2012}, between $W(\ka_i)$ and $W'(\ka_i)$.
For neighborhoods\footnote{If $U, U'$ are subsets of a topological space, $S$, we write $U' \sqsubset U$ if there is a continuous function $f:S\to [0,1]$ such that $f(U')=1$ and $f(S\setminus U)=0$.}
$W_1(\ka_i)\sqsubset W_2(\ka_i)\sqsubset W(\ka_i)$ of the
product connection strata \eqref{eq:IntroTrivialStrata}
and $f:W(\ka_i)\to [0,1]$ satisfying $W_1(\ka_i)\subset f^{-1}(0)$ and $W(\ka_i) \setminus W_2(\ka_i)\subset f^{-1}(1)$,  let $W''(\ka_i)$ be the image of $W(\ka_i)$ under the map $R_{f(\cdot)}(\cdot)$. Defining $\bar M^s_{\spl,\ka_i}$ by
\[
\left(\bar M^s_{\ka_i}(S^4) \setminus W'(\ka_i)\right)
\cup
W''(\ka_i),
\]
then yields the desired smoothly stratified space with an Uhlenbeck neighborhood of the
product connection strata given by the image of a splicing map.
By redefining $M(\Si)$ with the aid of the above construction of the instanton moduli space with spliced ends, we can define $\rho^d_{\Si,\Si'}$ using a splicing map.

\subsection{Space of global splicing data}
\label{subsec:SpaceOfGlobalSplicingData}
With the redefinition of $M(\Si)$ described in Section \ref{subsec:Introduction_spliced-ends_moduli_space},
the only remaining problem in making the diagram \eqref{eq:IntroOverlapCD} commute is the possible failure of the equality \eqref{eq:IntroSplicingAssoc} when the Riemannian metric on $X$ and the background connection, $A_0$, are not flat.

These problems can be overcome by performing the splicing operation with respect to perturbations
of the metric and background connection which are flat near the splicing points.  Because the metric and connection $A_0$ need not be globally flat, this deformation will vary with the splicing point. This process defines the crude splicing maps, $\bga_\Si''$, appearing in \eqref{eq:IntroOverlapCD}.

The space of global splicing data, $\bar\sM^{\vir}_{\ft,\fs}$, is the pushout of the spaces $\Gl(\Si)$ by the diagram \eqref{eq:IntroOverlapCD}.  That is, $\bar\sM^{\vir}_{\ft,\fs}$ is the union of the spaces $\Gl(\Si)$ subject to the
relation that for $\bA\in\Gl(\Si,\Si')$, we have $\rho^d_{\Si,\Si'}(\bA)=\rho^u_{\Si,\Si'}(\bA)$.

\subsection{Definition of link of a subspace of a moduli space of ideal Seiberg--Witten monopoles}
The commutativity of the diagram \eqref{eq:IntroOverlapCD} allows us to define a
smoothly stratified space, $\bar\sM^{\vir}_{\ft,\fs}$, as the union of the images of the crude splicing maps discussed above.  This space is then the union of the cone bundle neighborhoods defined by the domains \eqref{eq:IntroGluingFiberBundle} of the crude splicing maps.  These local cone bundle structures allow us to define a subspace, $\bL^{\vir}_{\ft,\fs}\subset \bar\sM^{\vir}_{\ft,\fs}/S^1$, as the union of
smoothly stratified, codimension-one subspaces,
$\bL(\Si)\subset\Gl(\Si)/S^1$, each of which is a subbundle of the restriction of $\Gl(\Si)/S^1$
to a compact subspace, $K_\Si\Subset\Si$, namely
\begin{equation}
\label{eq:IntroSubBundle}
\bL(\Si)=\Fr(\Si)|_{M_{\fs}\times K_\Si}\times_{G(\Si)\times S^1}\rd M(\Si),
\end{equation}
where $\rd M(\Si)\subset M(\Si)$ is a
smoothly stratified, codimension-one subspace.
The link, $\bar\bL_{\ft,\fs}$, will be defined, roughly, as the intersection of $\bL^{\vir}_{\ft,\fs}$ with the zero-locus of a
section of a vector bundle, referred to as the \emph{obstruction bundle}, over $\bar\sM^{\vir}_{\ft,\fs}$. The intersection,
$\bL(\Si)\cap\bL(\Si')$, will be a smoothly stratified,
codimension-one subspace of $\bL^{\vir}_{\ft,\fs}$ and can be described as
\begin{equation}
\label{eq:IntroPieceIntersection}
\begin{aligned}
\bL(\Si)\cap\bL(\Si')
{}&=
\Fr(\Si')|_{M_{\fs}\times \rd_\Si K_{\Si'}}
\times_{G(\Si')\times S^1}
\rd M(\Si')
\\
{}&=
\Fr(\Si)|_{K_\Si}\times_{G(\Si)\times S^1}\rd_{\Si'}\rd M(\Si),
\end{aligned}
\end{equation}
where $\rd_\Si K_{\Si'}\subset K_{\Si'}$ is a codimension-one
boundary admitting a fibration, $\rd_\Si K_{\Si'}\to K_\Si$,
and $\rd_{\Si'}\rd M(\Si)\subset  \rd M(\Si)$ is a
smoothly stratified, codimension-one subspace.

\subsection{Computation of intersection numbers with the link of the
  moduli space of ideal Seiberg--Witten monopoles}
\label{subsubsec:IntroComp}
Duality arguments allow us to convert the intersection number \eqref{eq:IntroIntersectionExpression} into a pairing of cohomology classes, which we write schematically as
$$
\langle \barmu\smile e,[\bL^{\vir}_{\ft,\fs}]\rangle,
$$
where $\barmu$ is a cohomology class to which the
intersection of the geometric representatives are dual and $e$ denotes
the Euler class of the obstruction bundle. We observe that the cohomology classes $\barmu$ and $e$ are generated by cohomology classes pulled back from $M_{\fs}\times K_\Si$ and by the first Chern class of the $S^1$ action. We wish to write this pairing as a sum,
$$
\sum_i \langle \barmu\smile e,[\bL(\Si_i)]\rangle,
$$
over the
subspaces $\bL(\Si_i)\subset \bL^{\vir}_{\ft,\fs}$ and apply a pushforward-pullback argument to the diagram,
\begin{equation}
\label{eq:IntroPushPullDiag}
\begin{CD}
\bL(\Si)=\Fr(\Si)\times_{G(\Si)\times S^1} \rd M(\Si)
@>>>
\EG(\Si)\times_{G(\Si)\times S^1}\rd M(\Si)
\\
@VVV @VVV
\\
M_{\fs}\times K_\Si @>>> \BG(\Si)
\end{CD}
\end{equation}
To do so, we must select a representative of the cohomology class $\barmu\smile e$ with compact support away from the boundaries of $\bL(\Si)$ that are given by $\bL(\Si)\cap\bL(\Si')$.  We specify such a choice in a manner similar to that introduced in \cite{OzsvathBlowUp}. We define quotients $q_\Si:\bL(\Si)\to\widehat\bL(\Si)$ with the following properties:
\begin{enumerate}
\item
The map $q_\Si:\bL(\Si)\to\widehat\bL(\Si)$ is injective on the interior of $\bL(\Si)$.
\item
The restrictions of $q_\Si$ and $q_{\Si'}$ to the intersection $\bL(\Si)\cap\bL(\Si')$ are equal.
\item
There is a cohomology class $\widehat\mu\smile\widehat e$ on $\widehat\bL(\Si)$ such that the restriction of
$\barmu\smile e$ to $\bL(\Si)$ equals $q_{\Si}^*(\widehat\mu\smile\widehat e)$.
\item
Each quotient admits an (orbifold) fiber bundle structure,
$$
\widehat\bL(\Si)\to M_{\fs}\times\cl(\Si),
$$
with the same structure group $G(\Si)$.
\end{enumerate}
We construct the quotient $\widehat\bL(\Si)$ by replacing the restricted bundle $\Fr(\Si)|_{M_{\fs}\times K_\Si}$ with the extension of
that bundle over $M_{\fs}\times\cl(\Si)$. The equality \eqref{eq:IntroPieceIntersection} allows us
to do this and satisfy the second requirement above simultaneously by replacing the fibers $\rd M(\Si)$ with
a quotient of this fiber. Then, we can write
$$
\langle \barmu\smile e,[\bL^{\vir}_{\ft,\fs}]\rangle
=
\sum_i
\langle \widehat\mu\smile\widehat e,[\widehat\bL(\Si_i)]\rangle
$$
and apply the pushforward-pullback argument to the diagram analogous to  \eqref{eq:IntroPushPullDiag}
for $\widehat\bL(\Si)$.  Because the topology of the fiber bundle $\Fr(\Si)$ depends only on quantities
described following \eqref{eq:MainEquation}, such a pushforward-pullback argument will yield the desired result.

\section{Kotschick--Morgan Conjecture}
\label{subsec:IntroKMConj}
We now describe how the Kotschick--Morgan Conjecture \cite{Goettsche, KotschickMorgan} also follows from the methods in this monograph.

If $b^+(X)=1$, the Donaldson invariant depends on the metric.  If $\om^+(g)$ indicates the unique (once
an orientation for $H^+(X)$ is specified) unit-length, self-dual, harmonic two-form, the Donaldson invariant will change
when this harmonic form `crosses a wall' for reasons we now describe.

For manifolds with $b^+(X)>1$, one proves that the Donaldson invariant is independent of the metric
(see \cite[Theorem 9.2.12]{DK})
by considering the cobordism defined by
$$
\bar M^w_\ka(g_I)
=
\{[A,\bx,t]: [A,\bx]\in \bar M^w_\ka(g_t)\ \text{and $t\in [-1,1]$}\},
$$
where $\bar M^w_\ka(g_t)$ is the Uhlenbeck compactification of the moduli space of  connections which are anti-self-dual with respect to the Riemannian metric $g_t$ and $g_t$ is a smooth path of Riemannian metrics on $X$ connecting two different Riemannian metrics, $g_{-1}$ and $g_1$. A Donaldson invariant defined by the metric $g_t$ is given by a pairing
$$
\langle \barmu, [\bar M^w_\ka(g_t)]\rangle.
$$
The cohomology classes, namely the $\mu$-classes \cite{DK, KMStructure}, used in defining the Donaldson invariants extend over the cobordism $\bar M^w_\ka(g_I)$. Because the boundaries of this cobordism are $\bar M^w_\ka(g_{\pm 1})$, we can write
\begin{equation}
\label{eq:DInvariantWallCross}
\langle \barmu, [\bar M^w_\ka(g_1)]\rangle
-
\langle \barmu, [\bar M^w_\ka(g_{-1})]\rangle
=
\langle \delta^*\barmu, [\bar M^w_\ka(g_I)]\rangle
=
0,
\end{equation}
proving that the Donaldson invariant for the metric $g_1$ is equal to that for the metric $g_{-1}$.

If $b^+(X)=1$, the above argument fails because the $\mu$-classes do not extend over the cobordism. The cobordism may contain  singularities of the form
\begin{equation}
\label{eq:ReducibleASDConnInParamModuli}
\{[A_L]\}\times\Sym^\ell(X) \subset \bar M^w_\ka(g_t),
\end{equation}
where $A_L$ is a connection which is reducible with respect to a reduction of an $\SO(3)$-bundle, $\fg=\ubarRR\oplus L$, and $L$ is a complex line bundle with $c_1(L)\smile\om^+(g_t)=0$.  For each of these families of reducible singularities appearing in the cobordism, the change in the Donaldson invariant, as expressed by the left-hand-side of \eqref{eq:DInvariantWallCross}, will be given by
\begin{equation}
\label{eq:DWCPairing}
\langle \barmu, [\partial\bar U^w_\ka(L)]\rangle,
\end{equation}
where $\partial\bar U^w_\ka(L)$ is the link of the set of singularities \eqref{eq:ReducibleASDConnInParamModuli} in the
cobordism $\bar M^w_\ka(g_I)$.

In \cite[Conjectures 6.2.1 and 6.2.2]{KotschickMorgan}, Kotschick and Morgan proposed that the pairing \eqref{eq:DWCPairing} was given by a polynomial in $c_1(L)$ and the intersection form $Q_X$ and that the coefficients of this polynomial depended only on the homotopy type of $X$. Assuming this conjecture, G{\"o}ttsche was able to compute an explicit formula for the pairing \eqref{eq:DWCPairing} in \cite{Goettsche}.

This conjecture follows almost immediately from the methods of this monograph
and an analogue of Hypothesis \ref{hyp:Gluing} for the anti-self-dual equations.
Our gluing theorem for anti-self-dual connections \cite{FLKM1} (based in turn on \cite{DK, DonConn, MorganMrowkaTube, TauSelfDual, TauIndef,TauFrame,TauStable}), parameterize neighborhoods of the strata $[A_L]\times\Si$, by bundles that are almost identical to the bundles, $\Gl(\Si)$, described in
\eqref{eq:IntroGluingFiberBundle}.
The argument described in the preceding sections, where we employ the overlap spaces and instanton moduli space with spliced ends to create a space of global splicing data, can be adapted --- essentially without change --- to define a link, $\partial\bar U^w_\ka(L)$, to which we can apply the quotient arguments in Section \ref{subsubsec:IntroComp}, thus proving the Kotschick--Morgan Conjecture.

\section{Outline of the monograph}
We give a brief outline of our monograph. A review of notation and definitions from our previous
articles in this sequence appears in
Chapter \ref{chap:prelim}. Chapter \ref{chap:Diagonals} contains a description of a stratification of $\Sym^\ell(X)$, and definitions of normal bundles of the strata of $\Sym^\ell(X)$. In Chapter \ref{chap:diagTM}, we show that the projection maps of the normal bundles introduced in Chapter \ref{chap:Diagonals} satisfy Thom--Mather type equalities on the overlap of their images. In Chapter \ref{chap:SplicedEnd}, we define the spliced-ends moduli space of anti-self-dual connections over $S^4$ and this is used in Chapter \ref{chap:GlobalSplicingData} to define the space of global splicing data.  An analogous construction of an obstruction bundle over the space of global splicing data is carried out in Chapter \ref{chap:obstr}. In Chapter \ref{chap:Link}, we define the link
of a Seiberg--Witten moduli space contained in a lower level of the Uhlenbeck compactification of the moduli space of $\SO(3)$ monopoles
and describe the fiber bundle structure of the
subspaces of the link described in \eqref{eq:IntroSubBundle}.
We use cohomological computations in Chapter \ref{chap:Cohom} to compute the pullbacks of the relevant cohomology classes to the space of global splicing data.
In Chapter \ref{chap:Comp},  we perform the remaining computations required to prove the main results of the monograph. Finally, in Chapter \ref{chap:KMconj}, we show how the arguments of this monograph yield a proof of the Kotschick--Morgan Conjecture.

\chapter{Preliminaries}
\label{chap:prelim}
In this chapter, we introduce the notation and review the definitions from our preceding work on $\SO(3)$ monopoles, specifically \cite{FL2a,FL2b}. We begin in Section \ref{sec:ModuliSpace} by recalling the definition of the moduli space of $\SO(3)$ monopoles and its basic properties \cite{FL1, FeehanGenericMetric}. In Section \ref{sec:ASDsingularities}, we describe the stratum of zero-section monopoles, that is, anti-self-dual connections. In Section \ref{sec:Reducibles}, we discuss the strata of reducible, or
Seiberg--Witten monopoles, together with their `virtual' neighborhoods and normal bundles.  In Section \ref{sec:Cohomology}, we define the cohomology classes which will be
paired with the links of the anti-self-dual and Seiberg--Witten moduli spaces.  In Section \ref{sec:Donaldsonseries}, we review the definition of the Donaldson series.  Lastly, in Section \ref{sec:LinkCobordism}, we describe the basic relation between the pairings with links of the anti-self-dual and Seiberg--Witten moduli spaces provided by the $\SO(3)$-monopole cobordism.

\section{The moduli space of $\SO(3)$ monopoles}
\label{sec:ModuliSpace}
Throughout this monograph, $(X,g)$ will denote a closed, connected, oriented, smooth, Riemannian four-manifold.

\subsection{Clifford modules}
\label{subsubsec:SpincuStr}
Let $V$ be a Hermitian vector bundle over $(X,g)$ and let $\rho:T^*X\to
\End_\CC(V)$ be a real-linear map satisfying
\begin{equation}
\label{eq:CliffordMapDefn}
\rho(\alpha)^2 = -g(\alpha,\alpha)\id_{V}
\quad\text{and}\quad
\rho(\alpha)^\dagger = -\rho(\alpha),
\quad \alpha \in C^\8(T^*X).
\end{equation}
The map $\rho$ uniquely extends to a linear isomorphism, $\rho:\Lambda^{\bullet}(T^*X)\otimes_\RR\CC\to\End_\CC(V)$, and gives $V$ the structure of a Hermitian Clifford module for the complex Clifford algebra $\CCl(T^*X)$.  There is a splitting $V=V^+\oplus V^-$, where $V^\mp$ are the
$\pm 1$ eigenspaces of $\rho(\vol)$. A unitary connection $A$ on $V$ is \emph{spin} if
\begin{equation}
\label{eq:SpinConnection}
[\nabla_A,\rho(\alpha)] =\rho(\nabla\alpha)
\quad\text{on }C^\8(V),
\end{equation}
for any $\alpha\in C^\8(T^*X)$, where $\nabla$ is the Levi-Civita connection.

A Hermitian Clifford module $\fs=(\rho,W)$ is a \emph{\spinc structure}\label{spinc_structure} when $W$ has complex rank four; it defines a class
\begin{equation}
\label{eq:DefineChernClassOfSpinc}
c_1(\fs)=c_1(W^+),
\end{equation}
which is an integral lift of the second Stiefel-Whitney class.  If $L\to X$ is a complex line bundle, we write $\fs\otimes L$ for the \spinc structure $(\rho\otimes\id_L,W\otimes L)$.

We call a Hermitian Clifford module $\ft=(\rho,V)$ a \emph{\spinu structure}\label{spinu_structure} when $V$ has complex rank eight. Recall that $\fg_{\ft}\subset\su(V)$ is the
$\SO(3)$ subbundle given by the span of the sections of the bundle $\su(V)$ which commute with the action of $\CCl(T^*X)$ on $V$. We obtain a splitting,
\begin{equation}
\label{eq:EndSplitting}
\su(V^+)
\cong
\rho(\Lambda^+)\oplus i\rho(\Lambda^+)\otimes_\RR\fg_{\ft}
\oplus \fg_{\ft},
\end{equation}
and similarly for $\su(V^-)$. The fibers $V_x^+$ define complex lines whose tensor-product square is $\det(V^+_x)$ and thus a complex line bundle over
$X$,
\begin{equation}
\label{eq:CliffordDeterminantBundle}
{\det}^{\frac{1}{2}}(V^+).
\end{equation}
A \spinu structure $\ft$ thus defines characteristic classes,
\begin{equation}
\label{eq:SpinUCharacteristics}
c_1(\ft)=\textstyle{\frac{1}{2}} c_1(V^+),
\quad
p_1(\ft) = p_1(\fg_{\ft}),
\quad
\text{and}\quad
w_2(\ft)=w_2(\fg_{\ft}).
\end{equation}
Given a \spinc bundle $W$, one has an isomorphism $V\cong W\otimes_\CC E$ of Hermitian Clifford modules, where $E$ is a complex rank-two Hermitian vector bundle
\cite[Lemma 2.3]{FL2a}; then
\begin{equation}
\label{eq:SpinAssociatedBundles}
\fg_{\ft}
=
\su(E)
\quad\text{and}\quad
{\det}^{\frac{1}{2}}(V^+)
=
\det(W^+)\otimes_\CC\det(E).
\end{equation}
If for $\ft=(\rho,V)$, there is a \spinc structure $\fs=(\rho_W,W)$ and the bundle $V$ admits a splitting $V\cong W\oplus W\otimes L$ and $\rho=\rho_W\oplus \rho_W\otimes\id_L$, then we write $\ft=\fs\oplus \fs\otimes L$.

\subsection{$\SO(3)$ monopoles}
\label{subsubsec:PU2Monopoles}
We fix a smooth unitary connection $A_\La$ on the line bundle ${\det}^{\frac{1}{2}}(V^+)$, let $k\geq 2$ be an integer, and let $\sA_{\ft}$ be the affine space of $L^2_k$ spin connections\footnote{We adopt the notation of Freed and Uhlenbeck \cite{FU} for Sobolev spaces.} on $V$ which induce the connection $2A_\La$ on $\det(V^+)$. If $A$ is a spin connection on $V$, then it defines an $\SO(3)$ connection, $\hat A$, on the subbundle $\fg_{\ft}\subset\su(V)$
\cite[Lemma 2.5]{FL2a}; conversely, every $\SO(3)$ connection on $\fg_{\ft}$ lifts to a unique spin connection on $V$ inducing the connection $2A_\Lambda$ on $\det(V^+)$ \cite[Lemma 2.11]{FL2a}.

Let $\sG_{\ft}$\label{Gauge_transformations_on_SO(3)_pairs} denote the group of $L^2_{k+1}$ unitary automorphisms of $V$ which commute with $\CCl(T^*X)$ and which have Clifford-determinant one (see \cite[Definition 2.6]{FL2a}). Define\label{Configuration_space_SO(3)_pairs}
\begin{equation}
\label{eq:SpinUPreConfiguration}
\tsC_\ft = \sA_\ft\times L^2_k(X,V^+)
\quad\text{and}\quad
\sC_\ft = \tsC_\ft/\sG_\ft.
\end{equation}
The action of $\sG_{\ft}$ on $V$ induces an adjoint action on $\End_\CC(V)$, acting as the identity on $\rho(\Lambda^\bullet_\CC)\subset\End_\CC(V)$ and inducing an adjoint
action on $\fg_{\ft}\subset\End_\CC(V)$ (see \cite[Lemma 2.7]{FL2a}). The space $\tsC_\ft$ carries a circle action\label{Circle_action_on_SO(3)_pairs} induced by scalar multiplication on $V$:
\begin{equation}
\label{eq:S1ZAction}
S^1\times V \to V,
\quad (e^{i\theta},\Phi)\mapsto e^{i\theta}\Phi.
\end{equation}
Because this action commutes with that of $\sG_{\ft}$, the action \eqref{eq:S1ZAction} also defines an action on $\sC_{\ft}$. Note that $-1\in S^1$ acts trivially on $\sC_{\ft}$.  Let
$\sC^0_\ft\subset\sC_\ft$ be the subspace represented by pairs whose spinor components are not identically zero, let $\sC^*_\ft\subset\sC_\ft$ be the subspace represented by pairs where the induced $\SO(3)$ connections on $\fg_{\ft}$ are irreducible, and let $\sC^{*,0}_\ft=\sC^*_\ft\cap \sC^0_\ft$.

We call a pair $(A,\Phi)$ in $\tsC_\ft$ an \emph{$\SO(3)$ monopole}\label{SO(3)_monopole} if
\begin{equation}
\label{eq:PerturbedSO3MonopoleEquations}
\fS(A,\Phi)
:=
\begin{pmatrix}
\ad^{-1}(F^+_{\hat A}) - \tau\rho^{-1}(\Phi\otimes\Phi^*)_{00}
\\
D_A\Phi + \rho(\vartheta)
\end{pmatrix}
= 0.
\end{equation}
Here, $D_A=\rho\circ \cov_A:C^\8(X,V^+)\to C^\8(X,V^-)$ is the Dirac operator\label{Dirac_operator}; the isomorphism $\ad:\fg_{\ft}\to \so(\fg_{\ft})$ identifies the self-dual component of the curvature $F_{\hat A}^+$, a section of $\La^+\otimes\so(\fg_{\ft})$, with $\ad^{-1}(F^+_{\hat A})$, a section of $\La^+\otimes\fg_{\ft}$; the section $\tau$ of $\GL(\Lambda^+)$ is a
perturbation close to the identity; the perturbation $\vartheta$ is a complex one-form close to zero; $\Phi^* \in \Hom(V^+,\CC)$ is the pointwise Hermitian dual $\langle\cdot,\Phi\rangle$
of $\Phi$; and $(\Phi\otimes\Phi^*)_{00}$ is the component of the section $\Phi\otimes\Phi^*$ of $i\fu(V^+)$ lying in $\rho(\Lambda^+)\otimes\fg_{\ft}$ with respect to the splitting $\fu(V^+)=i\underline{\RR}\oplus\su(V^+)$ and decomposition \eqref{eq:EndSplitting} of $\su(V^+)$.

Equation \eqref{eq:PerturbedSO3MonopoleEquations} is invariant under the action of $\sG_\ft$. We let $\sM_{\ft}\subset \sC_\ft$\label{moduli_space_SO(3)_monopoles} be the subspace represented by pairs satisfying
equation
\eqref{eq:PerturbedSO3MonopoleEquations} and write
\begin{equation}
\label{eq:PU2MonopoleSubspaces}
\sM^*_{\ft} = \sM_{\ft}\cap\sC^*_\ft, \quad
\sM^0_{\ft}= \sM_{\ft}\cap\sC^0_\ft, \quad\text{and}\quad
\sM^{*,0}_\ft = \sM_{\ft}\cap\sC^{*,0}_\ft.
\end{equation}
Since equation \eqref{eq:PerturbedSO3MonopoleEquations} is invariant under the circle action induced by scalar multiplication on $V$, the subspaces \eqref{eq:PU2MonopoleSubspaces} of $\sC_{\ft}$ are also invariant under this action.

\begin{thm}
\label{thm:Transv}
(See Feehan \cite[Theorem 1.1]{FeehanGenericMetric} and Teleman \cite{TelemanGenericMetric}.)
Let $X$ be a closed, oriented, smooth four-manifold and let $V$ be a complex rank-eight, Hermitian vector bundle over $X$. Then for parameters $(\rho,g,\tau,\vartheta)$, which are generic in the sense of \cite{FeehanGenericMetric}, and $\ft=(\rho,V)$, the space $\sM^{*,0}_{\ft}$ is a smooth manifold of the expected dimension,
\begin{equation}
\label{eq:Transv}
\begin{aligned}
\dim \sM^{*,0}_{\ft}
=
d(\ft)
= d_a(\ft)+2n_a(\ft),
\quad\text{where }
d_a(\ft)
&=-2p_1(\ft)- \frac{3}{2}(\chi(X)+\sigma(X)),
\\
n_a(\ft)
&= \frac{1}{4}(p_1(\ft)+c_1(\ft)^2-\sigma(X)),
\end{aligned}
\end{equation}
and $\chi(X)$ is the Euler characteristic and $\si(X)$ is the signature of $X$.
\end{thm}

For the remainder of the article, we assume that the \emph{perturbation parameters}\label{perturbation_parameters} in \eqref{eq:PerturbedSO3MonopoleEquations} are chosen as indicated in
Theorem \ref{thm:Transv}. For each non-negative integer $\ell$, let $\ft(\ell)=(\rho,V_{\ell})$ be the \spinu structure characterized by
\begin{equation}
\label{eq:DefineLowerChargeSpinuStr}
c_1(V_\ell)=c_1(V),
\quad
p_1(\fg_{\ft(\ell)})=p_1(\fg_{\ft})+4\ell,
\quad\text{and}\quad
w_2(\fg_{\ft(\ell)})=w_2(\fg_{\ft}).
\end{equation}
We let $\bar\sM_{\ft}$ denote the closure of $\sM_{\ft}$ in the space of \emph{ideal $\SO(3)$ monopoles}\label{ideal_SO(3)_monopoles},
\begin{equation}
\label{eq:idealmonopoles}
I\sM_{\ft} = \bigsqcup_{\ell=0}^\8(\sM_{\ft(\ell)}\times\Sym^\ell(X)),
\end{equation}
with respect to an \emph{Uhlenbeck topology}\label{Uhlenbeck_topology} \cite[Definition 4.19]{FL1} defined by the following notion of convergence.

\begin{defn}
\label{defn:UhlenbeckConvergence}
A sequence $\{[A_\alpha,\Phi_\alpha]\}_{\alpha\in\NN}
\subset\sC_{\ft}$ converges to an ideal pair $[A_0,\Phi_0,\bx]\in\sC_{\ft(\ell)}\times\Sym^\ell(X)$ if
\begin{itemize}
\item
There is a sequence of $L^2_{k+1,\loc}$ \spinu bundle isomorphisms, $u_\alpha: V_\ft|_{X-\bx} \to V_{\ft(\ell)}|_{X-\bx}$, such that the sequence of pairs $u_\alpha(A_\alpha,\Phi_\alpha)$ converges as $\alpha\to\infty$ to $(A_0,\Phi_0)$ in $L^2_{k,\loc}$ over $X-\bx$ and
\item
The sequence of measures, $|F_{\hat A_\alpha}|^2\,d\vol$ converges as $\alpha\to\infty$ in the weak-star topology to the measure $|F_{\hat A}|^2\,d\vol + 8\pi^2\sum_{x\in \bx} \delta_x$, where $\delta_x$ is the Dirac delta
measure centered at $x \in X$.
\end{itemize}
We refer to neighborhoods in the topology defined by Uhlenbeck convergence as \emph{Uhlenbeck neighborhoods}.
\end{defn}
We call the intersection of $\bar\sM_{\ft}$ with $\sM_{\ft(\ell)}\times \Sym^\ell(X)$ its {\em $\ell$-th level\/}.
We define the space $\bar\sC_\ft$ to be the set
\begin{equation}
\label{eq:IdealPairs}
\bar\sC_{\ft}: = \bigsqcup_{\ell=0}^\8(\sC_{\ft(\ell)}\times\Sym^\ell(X)),
\end{equation}
with the topology given by Definition \ref{defn:UhlenbeckConvergence}.

\begin{thm}
\label{thm:Compactness}
(See Feehan and Leness \cite[Theorem 1.1]{FL1}.)
Let $X$ be a closed Riemannian four-manifold with \spinu structure $\ft$.
Then there is a positive integer $N$, depending at most on
the scalar curvature of $X$,
the curvature of the chosen unitary connection on $\det(V^+)$,
and $p_1(\ft)$, such that the Uhlenbeck closure $\bar\sM_{\ft}$ of $\sM_{\ft}$ in $\sqcup_{\ell=0}^N(\sM_{\ft(\ell)}\times\Sym^\ell(X))$ is a second-countable, compact, Hausdorff space. The space $\bar\sM_{\ft}$ carries a continuous circle action, which restricts to the circle action defined on $\sM_{\ft_\ell}$ for each level.
\end{thm}

\section{Stratum of anti-self-dual or zero-section solutions}
\label{sec:ASDsingularities} {}From equation \eqref{eq:PerturbedSO3MonopoleEquations}, we see that the stratum of $\sM_{\ft}$ represented by pairs with zero spinor is identified with
\begin{equation}
\label{eq:ASDModuliSpace}
\{A\in\sA_{\ft}: F_{\hat A}^+ = 0\}/\sG_{\ft} \cong M_\kappa^w(X,g),
\end{equation}
the moduli space of $g$-anti-self-dual connections on the $\SO(3)$ bundle, $\fg_{\ft}$, where $\ka=-\quarter p_1(\ft)$ and $w\equiv w_2(\ft)\pmod 2$. For a generic Riemannian metric $g$, the space $M_\kappa^w(X,g)$ is a smooth manifold of the expected dimension,  $-2p_1(\ft) - \threehalf(\chi+\sigma)=d_a(\ft)$.

As explained in \cite[Section 3.4.1]{FL2a}, it is desirable to choose $w\pmod{2}$ so as to exclude points in $\bar\sM_{\ft}$ with associated flat
$\SO(3)$ connections, so we have a \emph{disjoint} union,
\begin{equation}
\label{eq:StratificationCptPU(2)Space}
\bar\sM_{\ft}
\cong
\bar\sM_{\ft}^{*,0} \sqcup \bar M_\kappa^w \sqcup \bar\sM_{\ft}^{\red},
\end{equation}
where $\bar\sM_{\ft}^*\subset\bar\sM_{\ft}$ is the subspace represented by triples whose associated $\SO(3)$ connections are irreducible, $\bar\sM_{\ft}^0\subset\bar\sM_{\ft}$ is the subspace represented by triples whose spinors are not identically zero, $\bar\sM_{\ft}^{*,0} = \bar\sM_{\ft}^{*}\cap\bar\sM_{\ft}^{0}$, while $\bar\sM_{\ft}^{\red}\subset\bar\sM_{\ft}$
is the subspace $\bar\sM_{\ft}-\bar\sM_{\ft}^*$ represented by triples whose associated $\SO(3)$ connections are reducible. We recall the

\begin{defn}
\label{defn:Good}
(See Feehan and Leness \cite[Definition 3.20]{FL2b}.)
A class $v\in H^2(X;\ZZ/2)$ is {\em good\/} if no integral lift of $v$ is
torsion.
\end{defn}

If $w\pmod{2}$ is good, then the union \eqref{eq:StratificationCptPU(2)Space} is disjoint, as desired. In practice, rather than constraining $w\pmod{2}$ itself, we use the
blow-up trick of \cite{MorganMrowkaPoly}, replacing $X$ with the blow-up, $X\#\overline{\CC\PP}^2$, and replacing $w$ by $w+\PD[e]$ (where $e\in H_2(X;\ZZ)$ is the exceptional class and $\mathrm{PD}[e]$ denotes its Poincar{\'e} dual), noting that $w+\PD[e]\pmod{2}$ is always good, and define gauge-theoretic invariants of $X$ in terms of moduli spaces on $X\#\overline{\CC\PP}^2$. When $w\pmod{2}$ is good, we
define \cite[Definition 3.7]{FL2a} the link of $\barM^w_{\ka}$ in $\bar\sM_{\ft}/S^1$ by
\begin{equation}
\label{eq:DefineASDLink}
\bar\bL^{w}_{\ft,\ka}
=
\{[A,\Phi,\bx]\in\bar\sM_{\ft}/S^1: \|\Phi\|_{L^2}^2=\eps\},
\end{equation}
where $\eps$ is a small positive constant; for generic $\eps$, the link, $\bar\bL^w_{\ft,\ka}$, is a smoothly-stratified, codimension-one subspace of $\bar\sM_{\ft}/S^1$.

\section{Strata of Seiberg--Witten or reducible solutions}
\label{sec:Reducibles}
We call a pair $(A,\Phi)\in\tsC_\ft$ \emph{reducible} if the connection $A$ on $V$ respects a splitting,
\begin{equation}
\label{eq:BasicSplitting}
V = W\oplus W\otimes L = W\otimes (\underline{\CC}\oplus L),
\end{equation}
for some \spinc structure $\fs=(\rho,W)$ and complex line bundle $L$, in which case $c_1(L)=c_1(\ft)-c_1(\fs)$. A spin connection $A$ on $V$ is reducible with respect to the splitting
\eqref{eq:BasicSplitting} if and only if $\hat A$ is reducible with respect to the splitting $\fg_{\ft}\cong \underline{\RR}\oplus L$, \cite[Lemma 2.9]{FL2a}. If $A$ is reducible, we can write $A=B\oplus B\otimes A_L$, where $B$ is a spin connection on $W$ and $A_L$ is a unitary connection on $L$; then $\hat A= d_{\RR}\oplus A_L$ and $A_L = A_\La\otimes (B^{\det})^*$, where $B^{\det}$ is the connection on $\det(W^+)$ induced by $B$ on $W$
and $d_\RR$ is the product connection on the line bundle $\underline{\RR} := X\times\RR$.

\subsection{Seiberg--Witten monopoles}
\label{subsec:SWMonopoles}
Given a \spinc structure $\fs=(\rho,W)$ on $X$, let $\sA_{\fs}$ denote the affine space of $L^2_k$ spin connections on $W$. Let $\sG_{\fs}$ denote the group of $L^2_{k+1}$ unitary automorphisms of $W$, commuting with $\CCl(T^*X)$, which we identify with $L^2_{k+1}(X,S^1)$. We then define
\begin{equation}
\label{eq:SpincPreConfig}
\tsC_{\fs}= \sA_{\fs}\times L^2_k(W^+)
\quad\text{and}\quad
\sC_{\fs}=\tsC_{\fs}/\sG_{\fs},
\end{equation}
where $\sG_{\fs}$ acts on $\tsC_{\fs}$ by
\begin{equation}
\label{eq:SWGaugeGroupAction}
(s,(B,\Psi))
\mapsto
s(B,\Psi)
=
(B- (s^{-1}ds)\id_W, s\Psi).
\end{equation}
We call a pair $(B,\Psi)\in\tsC_{\fs}$ a Seiberg--Witten monopole if
\begin{equation}
\label{eq:SeibergWitten}
\begin{aligned}
\Tr(F^+_B) - \tau\rho^{-1}(\Psi\otimes\Psi^*)_{0} - F^+(A_{\Lambda})
&=0,
\\
D_B\Psi + \rho(\vartheta)\Psi
&=0,
\end{aligned}
\end{equation}
where $\Tr:\fu(W^+)\to i\ubarRR$ is defined by the trace on $2\times 2$ complex matrices, $(\Psi\otimes\Psi^*)_0$ is the component of the section $\Psi\otimes\Psi^*$ of $i\fu(W^+)$ contained in $i\su(W^+)$, $D_B:C^\8(W^+)\to C^\8(W^-)$ is the Dirac operator, and $A_\Lambda$ is a unitary connection on a line bundle with first Chern class $\Lambda\in
H^2(X;\ZZ)$.  The perturbations are chosen so that solutions to equation \eqref{eq:SeibergWitten} are identified with reducible solutions to \eqref{eq:PerturbedSO3MonopoleEquations} when $c_1(\ft)=\Lambda$. Let $\tM_{\fs}\subset\tsC_{\fs}$ be the subspace cut out by equation \eqref{eq:SeibergWitten} and denote the moduli space of Seiberg--Witten monopoles by $M_{\fs}=\tM_{\fs}/\sG_{\fs}$.

\subsection{Seiberg--Witten invariants}
\label{subsec:SWInvariants}
We let $\sC_{\fs}^0\subset\sC_{\fs}$ be the open subspace represented by pairs whose spinor components are not identically zero and define a complex line bundle over $\sC_\fs^0\times X$ by
\begin{equation}
\label{eq:DefineSWUniversal}
\LL_{\fs}= \tsC_\fs^0\times_{\sG_\fs}\underline{\CC},
\end{equation}
where $\underline{\CC}=X\times\CC$ and $s\in \sG_{\fs}$ acts on $(B,\Psi)\in \tsC_{\fs}$ and $(x,\zeta)\in \underline{\CC}$ by
\begin{equation}
\label{eq:DefineSWUniversalS1Action}
((B,\Psi),(x,\zeta))
\mapsto
(s(B,\Psi),(x,s(x)^{-1}\zeta)).
\end{equation}
Define
\begin{equation}
\label{eq:DefineAAA2}
\AAA_2(X) = \Sym\left( H_{0}(X;\RR)\right)\otimes\La^\bullet(H_{1}(X;\RR))
\end{equation}
to be the graded algebra, with $z=\beta_1\beta_2\cdots\beta_r$ having total degree $\deg(z) = \sum_p(2-i_p)$, when $\beta_p\in H_{i_p}(X;\RR)$. Then the map,
\begin{equation}
\label{eq:SWMuMap}
\mu_{\fs}:H_\bullet(X;\RR)\to H^{2-\bullet}(\sC^0_{\fs};\RR),
\quad
\beta\mapsto c_1(\LL_{\fs})/\beta,
\end{equation}
extends in the usual way to a homomorphism of graded real algebras,
$$
\mu_{\fs}:\AAA_2(X) \to H^{\bullet}(\sC^0_\fs;\RR).
$$
If $x\in H_0(X;\ZZ)$ denotes the positive generator, we set
\begin{equation}
\label{eq:SWClass}
\mu_{\fs} = c_1(\LL_{\fs})/x \in H^2(\sC_\fs^0;\ZZ).
\end{equation}
Equivalently, $\mu_{\fs}$ is the first Chern class of the $S^1$ base-point fibration over $\sC_\fs^0$. If $b_1(X)=0$, then $c_1(\LL_{\fs})=\mu_{\fs}\times 1$ by \cite[Lemma 2.14]{FL2a}.

For $b^+(X)>0$ and generic Riemannian metrics on $X$, the space $M_{\fs}$ contains zero-section pairs if and only if $c_1(\fs)-\Lambda$ is a torsion class by \cite[Proposition 6.3.1]{MorganSWNotes}. If $M_{\fs}$ contains no zero-section pairs then, for generic perturbations, it is a compact, orientable, smooth manifold of dimension
\begin{equation}
\label{eq:DimSW}
d_s(\fs)
=
\dim M_{\fs}
=
\frac{1}{4}(c_1(\fs)^2 -2\chi -3\sigma).
\end{equation}
Let $\tilde X=X\#\overline{\CC\PP}^2$ denote the blow-up of $X$ with exceptional class $e\in H_2(\tilde X;\ZZ)$ and denote its Poincar\'e dual by $\PD[e]\in H^2(\tilde X;\ZZ)$. Let $\fs^\pm=(\tilde\rho,\tilde W)$ denote the \spinc structure on $\tilde X$ with $c_1(\fs^\pm)=c_1(\fs)\pm \PD[e]$ obtained by splicing the \spinc structure $\fs=(\rho,W)$ on $X$
with the \spinc structure on $\overline{\CC\PP}^2$ with first Chern class $\pm \PD[e]$. (See \cite[Section 4.5]{FL2b} for an explanation of the relation between \spinc structures on $X$ and $\tilde X$.) Now
$$
c_1(\fs)\pm \PD[e]-\Lambda \in H^2(\tilde X;\ZZ)
$$
is not a torsion class and so --- for $b^+(X)>0$, generic Riemannian metrics on $X$ and related metrics on the connected sum $\tilde X$ --- the moduli spaces $M_{\fs^\pm}(\tilde X)$ contain no zero-section pairs. Thus, for our choice of generic perturbations, the moduli spaces $M_{\fs^\pm}(\tilde X)$ are compact, oriented, smooth manifolds, both of dimension $\dim M_{\fs}(X)$.

For $b_1(X)=0$ and odd $b^+(X)>1$, we define the \emph{Seiberg--Witten invariant} by \cite[Section 4.1]{FL2b}
\begin{equation}
\label{eq:DefSW}
  SW_{X}(\fs)
=
\langle\mu_{\fs^+}^{d},[M_{\fs^+}(\tilde X)]\rangle
=
\langle\mu_{\fs^-}^{d},[M_{\fs^-}(\tilde X)]\rangle,
\end{equation}
where $2d=d_s(\fs)=d_s(\fs^\pm)$. When $b^+_2(X)=1$ the pairing on the right-hand side of definition \eqref{eq:DefSW} depends on the chamber in the positive cone of $H^2(\tilde X;\RR)$ determined by the period point of the Riemannian metric on $\tilde X$. The definition of the Seiberg--Witten invariant for this case is also given in \cite[Section 4.1]{FL2b}: we assume that the class $w_2(X)-\Lambda\pmod{2}$ is good to avoid technical difficulties involved in chamber specification. Since $w\equiv w_2(X)-\Lambda \pmod{2}$, this coincides with the constraint we use to define the Donaldson invariants in Section \ref{sec:Donaldsonseries} when $b^+(X)=1$. We refer to \cite[Lemma 4.1 and Remark 4.2]{FL2b} for a comparison of the chamber structures required for the definition of Donaldson and Seiberg--Witten invariants when $b^+(X)=1$.

We say that $c_1(\fs)$ is a \emph{Seiberg--Witten basic class} if the map $\mu_{\fs}$ is non-trivial. The manifold $X$ is said to have \emph{Seiberg--Witten simple type} if all basic classes satisfy $d_s(\fs)=0$.

\subsection{Reducible $\SO(3)$ monopoles}
\label{subsec:RedPU2Monopole}
If $\ft=(\rho,V)$ and $\fs=(\rho,W)$ with $V=W\oplus W\otimes L$, then there is a smooth embedding
\begin{equation}
\label{eq:DefnOfIota}
\tsC_{\fs} \embed \tsC_{\ft},
\quad
(B,\Psi)\mapsto (B\oplus B\otimes A_\La\otimes B^{\det,*},\Psi\oplus 0),
\end{equation}
which is gauge-equivariant with respect to the homomorphism
\begin{equation}
\label{eq:GaugeGroupInclusion}
\varrho:\sG_{\fs}\embed \sG_{\ft},
\quad
s\mapsto \id_W\otimes\begin{pmatrix}s & 0 \\ 0 & s^{-1}\end{pmatrix}.
\end{equation}
According to \cite[Lemma 3.13]{FL2a}, the map \eqref{eq:DefnOfIota} defines a topological embedding $M_{\fs}^0\embed\sM_{\ft}$, where $M_{\fs}^0=M_{\fs}\cap\sC_{\fs}^0$ and an embedding of $M_{\fs}$ if $w_2(\ft)\neq 0$ or $b_1(X)=0$; its image in $\sM_{\ft}$ is represented by pairs which are reducible with respect to the splitting $V=W\oplus W\otimes L$. Henceforth, we shall not distinguish between $M_\fs$ and its image in $\sM_{\ft}$ under this embedding.

If $\fg_{\ft(\ell)}\cong i\ubarRR\oplus L$ where $c_1(L)=c_1(\ft)-c_1(\fs)$, then $p_1(\fg_{\ft(\ell)})=(c_1(\ft)-c_1(\fs))^2$.  Hence, for
\begin{equation}
\label{eq:ReducibleLevel}
\ell(\ft,\fs)=\frac{1}{4}\left((c_1(\ft)-c_1(\fs))^2- p_1(\ft)\right),
\end{equation}
the embedding \eqref{eq:DefnOfIota} gives an inclusion of $M_{\fs}\times\Sym^\ell(X)$ into $I\sM_{\ft}$, where $\ell=\ell(\ft,\fs)$.

\subsection{Circle actions}
\label{subsec:S1Actions}
When $V=W\oplus W\otimes L$ and $\ft=(\rho,V)$, the space $\tsC_{\ft}$ inherits a circle action defined by
\begin{equation}
\label{eq:DefineS1LAction}
S^1\times V \to V,
\quad (e^{i\theta},\Psi\oplus\Psi')\mapsto \Psi\oplus e^{i\theta}\Psi',
\end{equation}
where $\Psi\in C^\8(W)$ and $\Psi'\in C^\8(W\otimes L)$. With respect to the splitting $V=W\oplus W\otimes L$, the actions \eqref{eq:DefineS1LAction} and \eqref{eq:S1ZAction} are related by
\begin{equation}
\label{eq:RelateS1Actions}
\begin{pmatrix}1 & 0 \\ 0 & e^{i2\theta}\end{pmatrix}
=
e^{i\theta}u,
\text{ where }
u
= \varrho_\fs(e^{-i\theta})
:= \begin{pmatrix}e^{-i\theta} & 0 \\ 0 & e^{i\theta}\end{pmatrix}
\in\sG_{\ft},
\end{equation}
and so, when we pass to the induced circle actions on the quotient $\sC_{\ft}=\tsC_{\ft}/\sG_{\ft}$, the actions \eqref{eq:DefineS1LAction} and \eqref{eq:S1ZAction} on
$\sC_{\ft}$ differ only in their multiplicity. Recall \cite[Lemma 3.11]{FL2a} that the image  in $\tsC_{\ft}$ of the map \eqref{eq:DefnOfIota} contains all pairs which are fixed by the circle action \eqref{eq:DefineS1LAction}.

When $V=W\oplus W\otimes L$, the bundle $\fg_{\ft}$ defined prior to \eqref{eq:EndSplitting} admits a splitting
\begin{equation}
\label{eq:EndReduction}
\fg_{\ft}\cong i\underline{\RR}\oplus L,
\end{equation}
where $\underline{\RR}=X\times \RR$.
The action \eqref{eq:DefineS1LAction} induces an $S^1$ action on $\fg_{\ft}$ given, for $z\in L$ and $\zeta\in i\underline{\RR}$, by
\begin{equation}
\label{eq:S1LActionOnAdj}
\left( e^{i\theta},(\zeta,z)\right)
\mapsto
(\zeta, e^{-2i\theta}z),
\end{equation}
as described in \cite[Equation 3.56]{FLLevelOne}.

\subsection{The virtual normal bundle of the Seiberg--Witten moduli space}
\label{subsubsec:ThickenedNeighborhood}
Suppose $\ft=(\rho,V)$ and $\fs=(\rho,W)$, with $V=W\oplus W\otimes L$, so we have a topological embedding $M_{\fs}\embed \sM_{\ft}$; we assume $M_{\fs}$ contains no zero-section monopoles.  Recall from \cite[Section 3.5]{FL2a} that there exist finite-rank, complex vector bundles,
\begin{equation}
\label{eq:VirtualNormalandObstBundles}
\pi_{\Xi}:\Xi_{\ft,\fs}\to M_{\fs}
\quad\text{and}\quad
\pi_{N}:N_{\ft,\fs}\to M_{\fs},
\end{equation}
with $\Xi_{\ft,\fs}\cong M_{\fs}\times\CC^{r_\Xi}$, called the \emph{obstruction bundles}\label{Obstruction_bundle_Seiberg-Witten_moduli_space} and \emph{virtual normal bundles} of $M_{\fs}\embed\sM_{\ft}$, respectively.  For a small enough positive radius $\eps$, there are a topological embedding \cite[Theorem 3.21]{FL2a} of an open tubular neighborhood,
\begin{equation}
\label{eq:BackgroundConfigEmbedding}
\bga_{\fs}:N_{\ft,\fs}(\eps)
\embed
\sC_{\ft},
\end{equation}
and a smooth section $\bchi_{\fs}$ of the pulled-back complex vector bundle,
\begin{equation}
\label{eq:PulledBackObstructionBundle}
\pi_{N}^*\Xi_{\ft,\fs} \to N_{\ft,\fs}(\eps),
\end{equation}
such that the restriction of $\bga_{\fs}$ yields a homeomorphism
\begin{equation}
\label{eq:HomeoSWEmbedding}
\bga_{\fs}: \bchi_{\fs}^{-1}(0)\cap N_{\ft,\fs}(\eps)
\cong
\sM_{\ft}\cap \bga_{\fs}(N_{\ft,\fs}(\eps)),
\end{equation}
restricting to a diffeomorphism on the complement of $M_{\fs}$ and identifying $M_{\fs}$ with its image in $\sM_{\ft}$ under the embedding \eqref{eq:DefnOfIota}. We often refer to
the image $\bga_{\fs}(N_{\ft,\fs}(\eps))$ as a \emph{virtual moduli space}.\label{Virtual_Seiberg-Witten_moduli_space}

Our terminology is loosely motivated by that of \cite{GraberPand} and \cite{RuanSW}, where the goal (translated to our setting) would be to construct a \emph{virtual fundamental class}\label{Virtual_fundamental_class} for $\sM_{\ft}$, given by the cap product of the fundamental class of a virtual space containing $\sM_{\ft}$ with the Euler class of a vector bundle over this virtual space, where $\sM_{\ft}$ is the zero locus of a (possibly non-transversally vanishing) section. Here, $\bga_{\fs}(N_{\ft,\fs}(\eps))$ plays the role of the virtual space and (the pullback of) $\Xi_{\ft,\fs}$ the vector bundle with zero section yielding (an open neighborhood in) $\sM_{\ft}$.  Then, $N_{\ft,\fs}$ is the normal bundle of $M_{\fs}\embed \bga_{\fs}(N_{\ft,\fs}(\eps))$, while $[N_{\ft,\fs}]-[\Xi_{\ft,\fs}]$ would more properly be called the `virtual normal bundle' of $M_{\fs}\embed \sM_{\ft}$, in the language of $K$-theory.

By \cite[Equations (2.47) and (3.35)]{FL2a}, the negative of the index of the $\SO(3)$-monopole elliptic deformation complex  at a reducible solution can be written as
\begin{equation}
\label{eq:SWDimRelations}
\dim\sM_{\ft}
=
2n_s(\ft,\fs) + d_{s}(\fs),
\end{equation}
where $d_s(\fs)$ is the expected dimension of the Seiberg--Witten moduli space $M_{\fs}$ (see equation \eqref{eq:DimSW}), while $n_s(\ft,\fs)=n_s'(\ft,\fs)+n_s''(\ft,\fs)$ is minus the complex index of the normal deformation operator \cite[Equations (3.71) and (3.72)]{FL2a}), with
\begin{equation}
\label{eq:NormalComponentDims}
\begin{aligned}
n_s'(\ft,\fs)
&= -(c_1(\ft)-c_1(\fs))^2-\frac{1}{2}(\chi+\sigma),
\\
n_s''(\ft,\fs)
&=
\frac{1}{8}(c_1(\fs)-2c_1(\ft))^2-\sigma).
\end{aligned}
\end{equation}
If $r_\Xi$ is the complex rank of $\Xi_{\ft,\fs}\to M_{\fs}$, and $r_N(\ft,\fs)$ is the complex rank of $N_{\ft,\fs}\to M_{\fs}$, then
\begin{equation}
\label{eq:SWBundleRankRelations}
r_N(\ft,\fs)
=
n_s(\ft,\fs) + r_\Xi,
\end{equation}
as we can see from the dimension relation \eqref{eq:SWDimRelations} and the topological model \eqref{eq:HomeoSWEmbedding}.

The map \eqref{eq:BackgroundConfigEmbedding} is $S^1$-equivariant when $S^1$ acts trivially on $M_{\fs}$, by scalar multiplication on the fibers of $N_{\ft,\fs}(\eps)$, and by the action
\eqref{eq:DefineS1LAction} on $\sC_{\ft}$. The bundle \eqref{eq:PulledBackObstructionBundle} and section $\bchi_{\fs}$ are $S^1$-equivariant if $S^1$ acts on $N_{\ft,\fs}$ and the fibers
of $\bga_{\fs}^*\Xi_{\ft,\fs}$ by scalar multiplication.

Let $\tN_{\ft,\fs}\to \tM_{\fs}$ be the pullback of $N_{\ft,\fs}$ by the projection $\tM_{\fs}\to M_{\fs}=\tM_{\fs}/\sG_{\fs}$, so $\tN_{\ft,\fs}$ is a $\sG_{\fs}$-equivariant bundle,
where $\sG_{\fs}$ acts on the base $\tM_{\fs}$ by the usual gauge group action \eqref{eq:SWGaugeGroupAction} and the induced action on the total space,
\begin{equation}
\label{eq:VirtualNormalBundleTotalSpace}
\tN_{\ft,\fs}
\subset
\tM_{\fs}\times
L^2_k(\Lambda^1\otimes_\RR L) \oplus L^2(W^+\otimes L)
\subset
\tM_{\fs}\times
L^2_k(\Lambda^1\otimes_\RR \fg_{\ft}) \oplus L^2(V^+),
\end{equation}
via the embedding \eqref{eq:GaugeGroupInclusion} of $\sG_{\fs}$ into $\sG_{\ft}$ and the splittings $\fg_{\ft}\cong\underline{\RR}\oplus L$ \cite[Lemma 3.10]{FL2a} and
$V=W\oplus W\otimes L$.  Thus, $s\in\sG_{\fs}$ acts by scalar multiplication by $s^{-2}$ on sections of $\Lambda^1\otimes_\RR L$ and by $s^{-1}$ on sections of $W^+\otimes L$ \cite[Section 3.5.4]{FL2a}.

For a small enough positive $\eps$, there is a smooth embedding \cite[Section 3.5.4]{FL2a} of the open tubular neighborhood $\tN_{\ft,\fs}(\eps)$,
\begin{equation}
\label{eq:BackgroundPreConfigEmbedding}
\tilde\bga_{\fs}:\tN_{\ft,\fs}(\eps)
\to
\tsC_{\ft},
\end{equation}
which is gauge-equivariant with respect to the preceding action of $\sG_{\fs}$, and covers the topological embedding \eqref{eq:BackgroundConfigEmbedding}.  The map
\eqref{eq:BackgroundPreConfigEmbedding} is $S^1$-equivariant, where $S^1$ acts trivially on $\tM_{\fs}$, by scalar multiplication on the fibers of $\tN_{\ft,\fs}(\eps)$, and by the action
\eqref{eq:DefineS1LAction} on $\sC_{\ft}$. We note that the map \eqref{eq:BackgroundPreConfigEmbedding} is also $S^1$-equivariant with respect to the action \eqref{eq:S1ZAction} on $\sC_{\ft}$, if $S^1$ acts on $\tN_{\ft,\fs}$ by
\begin{equation}
\label{eq:S1ZActionOnN}
(e^{i\theta},(B,\Psi,\beta,\psi))
\mapsto
\varrho(e^{i\theta})( B,\Psi,e^{2i\theta}\beta,e^{2i\theta}\psi)
=
(B,e^{i\theta}\Psi,\beta,e^{i\theta}\psi),
\end{equation}
where $(B,\Psi)\in \tM_{\fs}$ and $(\beta,\psi)\in L^2_k(\Lambda^1\otimes_\RR L) \oplus L^2(W^+\otimes L)$, so $(B,\Psi,\beta,\psi)\in \tN_{\ft,\fs}$, and
$\varrho:\sG_{\fs}\to\sG_{\ft}$ is the homomorphism \eqref{eq:GaugeGroupInclusion}. This equivariance follows from the relation \eqref{eq:RelateS1Actions} between the actions \eqref{eq:DefineS1LAction} and \eqref{eq:S1ZAction}.

We note that the Chern character of the bundle $N_{\ft,\fs}$ is computed in \cite[Theorem 3.29]{FL2a} while the Segre classes of this bundle are computed, under some additional assumptions, in \cite[Lemma 4.11]{FL2b}.

\section{Cohomology classes on the moduli space of $\SO(3)$ monopoles}
\label{sec:Cohomology}
The identity \eqref{eq:MainEquation} arises as an equality between pairings of suitable cohomology classes with a link in $\bar\sM_{\ft}^{*,0}/S^1$ of the anti-self-dual moduli subspace and with the links of the Seiberg--Witten moduli subspaces in $\bar\sM_{\ft}/S^1$. We now review the definitions of these cohomology classes and their dual geometric representatives given in \cite[Section 3]{FL2b}.

The first kind of cohomology class is defined on $\sM^*_{\ft}/S^1$, via the associated $\SO(3)$ bundle,
\begin{equation}
\label{eq:UniversalU2Bundle}
\FF_{\ft}
=
\tsC^{*}_\ft/S^1\times_{\sG_\ft}\fg_{\ft}
\to
\sC^{*}_\ft/S^1\times X.
\end{equation}
The group $\sG_{\ft}$ acts diagonally in \eqref{eq:UniversalU2Bundle}, with $\sG_{\ft}$ acting on the left on $\fg_{\ft}$. We define \cite[Section 3.1]{FL2b}
\begin{equation}
\label{eq:DefnMuClasses}
\mu_p: H_\bullet(X;\RR)\to H^{4-\bullet}(\sC^*_\ft/S^1;\RR),
\quad
\beta\mapsto -\frac{1}{4} p_1(\FF_{\ft})/\beta.
\end{equation}
On restriction to $M^w_{\ka}\hookrightarrow \sM_{\ft}$, the cohomology classes $\mu_p(\beta)$ coincide with those used in the definition of Donaldson invariants \cite[Lemma 3.1]{FL2b}. Define
\begin{equation}
\label{eq:DefineAAA}
\AAA(X) = \Sym\left( H_{\even}(X;\RR)\right)\otimes\La^\bullet(H_{\odd}(X;\RR))
\end{equation}
to be the graded algebra, with $z=\beta_1\beta_2\cdots\beta_r$ having total degree $\deg(z) = \sum_p(4-i_p)$, when $\beta_p\in H_{i_p}(X;\RR)$. Then $\mu_p$ extends in the usual way to
a homomorphism of graded real algebras,
\begin{equation}
\label{eq:MuPMap}
\mu_p:\AAA(X) \to H^{\bullet}(\sC^*_\ft/S^1;\RR),
\end{equation}
which preserves degrees. Next, we define a complex line bundle over $\sC^{*,0}_{\ft}/S^1$,
\begin{equation}
\label{eq:DefineDetLineBundle}
\LL_{\ft}= \sC^{*,0}_{\ft}\times_{(S^1,\times -2)}\CC,
\end{equation}
where the $S^1$ action is given, for $[A,\Phi]\in\sC^{*,0}_{\ft}$ and $\zeta\in\CC$, by
\begin{equation}
\label{eq:DeterminantS1Action}
([A,\Phi],\zeta) \mapsto ([A,e^{i\theta}\Phi], e^{2i\theta}\zeta).
\end{equation}
Then we define the second kind of cohomology class on $\sM_{\ft}^{*,0}/S^1$ by
\begin{equation}
\label{eq:DefineMuC1}
\mu_c = c_1(\LL_{\ft}) \in H^2(\sC_{\ft}^{*,0}/S^1;\RR).
\end{equation}
For monomials $z\in\AAA(X)$, we constructed \cite[Section 3.2]{FL2b} geometric representatives\label{geometric_representatives} $\sV(z)$ dual to $\mu_p(z)$ (following the discussion in \cite{KMStructure})
and $\sW$ dual to $\mu_c$, defined on $\sM_{\ft}^*/S^1$ and $\sM_{\ft}^{*,0}/S^1$, respectively; their closures in $\bar\sM_{\ft}/S^1$ are denoted by $\bar\sV(z)$ and $\bar\sW$ \cite[Definition 3.14]{FL2b}. When
\begin{equation}
\label{eq:DimensionConditionPrelim}
\deg(z) +2\eta
=
\dim(\sM^{*,0}_{\ft}/S^1) -1,
\end{equation}
and $\deg(z)\ge \dim M^w_\ka$ it follows from \cite[Section 3.3]{FL2b} that the intersection
\begin{equation}
\label{eq:GeomReprIntersection}
\bar\sV(z) \cap \bar\sW^{\eta} \cap \bar\sM^{*,0}_{\ft}/S^1,
\end{equation}
is an oriented one-manifold (not necessarily connected) whose closure in $\bar\sM_{\ft}/S^1$ can only intersect $(\bar\sM_{\ft}-\sM_{\ft})/S^1$ at points in $\sM_{\ft}^{\red}\cong\cup(M_{\fs}\times\Sym^\ell(X))$, where the union is over $\ell\geq 0$ and $\fs\in\Spinc(X)$ \cite[Corollary 3.18]{FL2b}.

\section{Donaldson invariants}
\label{sec:Donaldsonseries}
We first recall the definition \cite[Section 2]{KMStructure} of the Donaldson series when $b_1(X)=0$ and $b^+(X)>1$ is odd, so that, in this case, $\chi(X)+\sigma(X)\equiv 0\pmod{4}$. See also Section 3.4.2 in \cite{FL2b}, especially for a definition of the Donaldson invariants when $b^+(X)=1$. For any choice of $w\in H^{2}(X;\ZZ)$, the Donaldson invariant is a linear function,
$$
D^{w}_{X}:\AAA(X) \to \RR.
$$
Let $\tilde X=X\#\overline{\CC\PP}^2$ be the blow-up of $X$ and let $e\in H_2(\tilde X;\ZZ)$ be the exceptional class, with Poincar\'e dual $\PD[e]\in H^2(\tilde X;\ZZ)$.  If $z\in\AAA(X)$ is a monomial, we define $D_X^w(z)=0$ unless
\begin{equation}
    \mathrm{deg}(z)\equiv -2w^{2}-\frac{3}{2}(\chi+\sigma)\pmod{8}.
    \label{mod8}
\end{equation}
If $\deg(z)$ obeys equation \eqref{mod8}, we let $\kappa\in\frac{1}{4}\ZZ$ be defined by
$$
\deg(z)=8\ka - \frac{3}{2}\left(\chi+\si\right).
$$
There exists an $\SO(3)$ bundle over $\tilde X$ with first Pontrjagin number $-4\ka -1$ and second Stiefel-Whitney class $w+\PD[e]\pmod{2}$. One then defines the Donaldson invariant on monomials by
\begin{equation}
\label{eq:DefineDonaldson}
D^w_X(z)
=
\#\left( \bar\sV(ze)\cap\barM^{w+\PD[e]}_{\ka+1/4}(\tilde X)\right),
\end{equation}
and extends $D_X^w$ to a real linear function on $\AAA(X)$. Note that $w+\PD[e]\pmod{2}$ is good in the sense of Definition \ref{defn:Good}. If $w'\equiv w \pmod{2}$, then \cite{DonOrient}
\begin{equation}
\label{eq:DonaldsonsSignChange}
D^{w'}_{X}=(-1)^{\frac{1}{4}(w'-w)^{2}}D^{w}_{X}.
\end{equation}
The Donaldson series is a formal power series,
\begin{equation}
\label{eq:DefineDonaldsonSeries}
\bD^{w}_{X}(h) = D^{w}_{X}((1+\textstyle{\frac{1}{2}} x)e^{h}),
\quad h \in H_{2}(X;\RR).
\end{equation}
By equation~\eqref{mod8}, the series $\bD^{w}_{X}$ is even if
$$
-w^{2}-\frac{3}{4}(\chi+\sigma)\equiv  0 \pmod 2,
$$
and odd otherwise. A four-manifold has Kronheimer--Mrowka simple type\label{Kronheimer-Mrowka_simple_type} if for some $w$ and all $z\in \AAA(X)$,
$$
D^{w}_{X}(x^{2}z)=4D^{w}_{X}(z).
$$
According to \cite[Theorem 1.7]{KMStructure}, when $X$ has Kronheimer--Mrowka simple type the series $\bD^{w}_{X}(h)$ is an analytic function of $h$ and there are finitely many characteristic cohomology classes, $K_{1},\ldots,K_{m}$ in $H^2(X;\ZZ)$ (the Kronheimer--Mrowka basic classes)\label{Kronheimer-Mrowka_basic_classes}, and non-zero rational numbers, $a_{1},\ldots,a_{m}$ (independent of $w$), so that\label{Kronheimer-Mrowka_structure_theorem}
$$
\bD^{w}_X(h)
=
e^{\half h \cdot h}
\sum_{i=1}^{r}(-1)^{\half(w^{2}+K_{i}\cdot w)}a_{i}e^{\langle
K_{i},h\rangle}.
$$
Witten's Conjecture \cite{Witten} then relates the Donaldson and Seiberg--Witten series for four-manifolds of simple type.

When $b^+(X)=1$, the pairing on the right-hand side of definition \eqref{eq:DefineDonaldson} depends on the chamber in the positive cone of $H^2(\tilde X;\RR)$ determined by the period point of the Riemannian metric on $\tilde X$, just as in the case of Seiberg--Witten invariants described in Section \ref{subsec:SWInvariants}. We refer to Section 3.4.2 in \cite{FL2b} for a detailed discussion of this case and, as in \cite{FL2b}, we assume that the class $w\pmod{2}$ is good in order to avoid technical difficulties involved in chamber specification.

\section{Links and the cobordism}
\label{sec:LinkCobordism}
Since the ends of the components of the one-manifold \eqref{eq:GeomReprIntersection} either lie near $M_\kappa^w$ or $M_{\fs}\times\Sym^\ell(X)$, for some $\fs$ and $\ell(\ft,\fs)\geq 0$ for which $\ft(\ell) = \fs\oplus\fs\otimes L$, we have when $w_2(\ft)$ is good in the sense of Definition \ref{defn:Good} that (see Theorem \ref{thm:CobordismThm})
\begin{equation}
\label{eq:RawCobordismSum}
\#\left(\bar\sV(z)\cap \bar\sW^{n_a-1}\cap \bar\bL^w_{\ft,\ka}\right)
=
-\sum_{\fs\in\Spinc(X)}
\# \left( \bar\sV(z)\cap \bar\sW^{n_a-1}\cap \bar\bL_{\ft,\fs}\right),
\end{equation}
where $\bar\bL^w_{\ft,\ka}$ is the link of $\bar M^w_\ka$ in $\bar\sM_{\ft}/S^1$ (see \cite[Definition 3.7]{FL2a}) and where $\bar\bL_{\ft,\fs}$ is empty if $\ell(\ft,\fs)<0$
and is the boundary of an open neighborhood of the Seiberg--Witten stratum $M_{\fs}\times\Sym^\ell(X)$ in $\bar\sM_{\ft}/S^1$ if $\ell=\ell(\ft,\fs)\ge 0$.  By construction, the intersection of $\bar\bL_{\ft,\fs}$ with the top stratum of $\bar\sM_{\ft}/S^1$ is a smoothly-stratified space and the intersection of $\bar\bL_{\ft,\fs}$ with the geometric representatives is in 
the top stratum of $\bar{\mathbf{L}}_{\mathfrak{t},\mathfrak{s}}$. The precise definition of $\bar\bL_{\ft,\fs}$ is given in \cite[Definition 3.22]{FL2a} for $\ell=0$ and in Definition \ref{defn:DefineLink} for $\ell\ge 1$.

When $\deg(z) = \dim M^w_{\ka}$ and $n_a(\ft)>0$, the intersection of the one-manifold \eqref{eq:GeomReprIntersection} with the link $\bar\bL^w_{\ft,\ka}$ is given by
\cite[Lemma 3.30]{FL2b}
\begin{equation}
\label{eq:ASDPairing}
2^{1-n_a}\#\left( \bar\sV(z)\cap \bar\sW^{n_a-1}\cap \bar\bL^w_{\ft,\ka}\right)
=
\#\left(\bar\sV(z)\cap\barM^w_{\ka}\right).
\end{equation}
Applying this identity to the blow-up, $X\#\overline{\CC\PP}^2$, when $n_a(\ft)>0$, we recover the Donaldson invariant $D_X^w(z)$ on the right-hand side of \eqref{eq:ASDPairing} via definition \eqref{eq:DefineDonaldson}.

\begin{rmk}
If in \eqref{eq:ASDPairing} we have $\deg(z)>\dim M^w_\ka$ and we replace $\bar\sW^{n_a-1}$ with $\bar\sW^\eta$ where $\eta$ satisfies \eqref{eq:DimensionConditionPrelim}, then the intersection number vanishes. If $n_a(\ft)\le 0$ and $\deg(z)$ still satisfies \eqref{eq:DimensionConditionPrelim} with $\eta=0$, then the intersection number in \eqref{eq:ASDPairing} is a constant times the spin polynomial invariant \cite{PTDirac}.
\end{rmk}

\chapter{Diagonals of symmetric products of manifolds}
\label{chap:Diagonals}
To describe the link of the family of singularities, $M_{\fs}\times\Sym^\ell(X)$, in $\bar\sM_{\ft}$, we will need to describe the strata of the symmetric product $\Sym^\ell(X)$, the normal bundles of the strata, and their incidence relations.

\section{Definitions}
\label{subsec:DiagDefs}
In this section, we introduce some vocabulary for describing the action of the symmetric group on $X^\ell$, the fixed point sets of this action, and the resulting strata of the symmetric product $\Sym^{\ell}(X)$. For basic definitions of the vocabulary of group actions, we refer to \cite[pp. 2--8]{tomDieck}.

\subsection{Subgroups of the symmetric group}
\label{subsubsec:SymGrps}
For any set $P$, let $\fS_P$ be the group of bijections,  $\si:P\to P$. For $P=N_\ell=\{1,\dots,\ell\}$, we will write $\fS_{\ell}$ for the symmetric group on $\ell$ elements.
We will write \emph{partitions} of $N_\ell$ as $\sP=\{P_1,\dots,P_{r(\sP)}\}$ where $P_i\subset N_\ell$.  We refer to the number of sets in $\sP$ as the \emph{length} of $\sP$. Each such partition of $N_\ell$ gives a partition of $\ell$: $\ell=|P_1|+\dots +|P_{r(\sP)}|$. The symmetric group $\fS_{\ell}$ acts on the set of partitions of $N_\ell$ by
\begin{equation}
\label{eq:ActionOnPartitions}
\left( \si,\sP\right)
\to \si(\sP):=
\{\si(P_1),\dots,\si(P_{r(\sP)})\}.
\end{equation}
The orbits of this action are distinguished by the partitions of $\ell$. Let $\Ga(\sP)<\fS_\ell$ be the stabilizer of a partition $\sP$, considered as an ordered set, with respect to the action \eqref{eq:ActionOnPartitions}.

\begin{lem}
\label{lem:CharacterizeGaP}
Let $\sP$ be a partition of $N_\ell$. Then the subgroup $\Ga(\sP)$ of $\fS_\ell$ is given by the image of the inclusion,
\begin{equation}
\label{eq:StabilizerInclusion}
\prod_{P\in\sP} \fS_{P} \to \fS_{\ell}.
\end{equation}
\end{lem}

\begin{proof}
Let $\Ga'(\sP)$ be the image of the homomorphism \eqref{eq:StabilizerInclusion}. There is then an inclusion $\Ga'(\sP)\to\Ga(\sP)$. This map is surjective because any $\si\in\fS_\ell$
which preserves the subsets $P_1,\dots,P_{r(\sP)}$ is a permutation of each of these subsets and thus in $\Ga'(\sP)$.
\end{proof}

\begin{lem}
\label{lem:CharNormalizerOfGaP}
Let $\fS(\sP)\le \fS_\ell$ be the normalizer of $\Ga(\sP)$.  Then for all $\si\in\fS(\sP)$ and all $P\in\sP$, one has that $\si(P)\in\sP$.
\end{lem}

\begin{proof}
Assume that there are
$\si\in\fS(\sP)$ and $P\in\sP$ with $\si(P)\notin\sP$. If $|P|=1$, then there is a $P'\in\sP$ with $\si(P)\subsetneqq P'$.  The strictness of the preceding inclusion implies that $|P'|>1$.
Thus, there are $a,b\in P'$ with $a\in\si(P)$ and $b\notin \si(P)$. By Lemma \ref{lem:CharacterizeGaP}, the transposition $(a\ b)$ is in $\Ga(\sP)$ and by the definition of $\fS(\sP)$ as the normalizer of $\Ga(\sP)$, we have $\si^{-1}(a\ b)\si=(c\ d)\in\Ga(\sP)$.  However, $c=\si^{-1}(a)\in P$ while $d=\si^{-1}(b)\notin P$, so $(c\ d)\notin\Ga(\sP)$.  This contradiction proves
the lemma if $|P|=1$.

If $|P|>1$, then there are $a,b\in \si(P)$ with $a\neq b$ and $P'\in\sP$ with $a\in P'$ and $b\notin P'$.  By Lemma \ref{lem:CharacterizeGaP}, the transposition $(a\ b)$ is not in $\Ga(\sP)$.  However, because $\si^{-1}(a),\si^{-1}(b)\in P$, the transposition $\si^{-1}(a\ b)\si=(\si^{-1}(a)\ \si^{-1}(b) )$ is in $\Ga(\sP)$, contradicting the definition of $\fS(\sP)$ as the normalizer of $\Ga(\sP)$.
\end{proof}

Lemma \ref{lem:CharNormalizerOfGaP} shows that $\fS(\sP)$ acts via \eqref{eq:ActionOnPartitions} as a permutation of $\sP$.  If we define
\begin{equation}
\label{eq:WhOfGaP}
W(\sP):=\fS(\sP)/\Ga(\sP),
\end{equation}
then $W(\sP)$ will act freely on $\sP$.

\subsection{Definition of the diagonals}
\label{subsec:DefiningTheDiagonals}
We define the strata of $\Sym^\ell(X)$ by the quotients of diagonals in $X^\ell$ by the permutation action of $\fS_\ell$ on $X^\ell$.
For a partition $\sP$ of $N_\ell$, define
\begin{equation}
\label{eq:DefineDelta}
\begin{aligned}
\Delta^\circ(X^\ell,\sP) & := \{ (x_1,\dots,x_{\ell}) \in X^\ell:
 x_i=x_j \
\text{if and only if $\exists\, P\in\sP$ with $i,j\in P$ } \},
\\
\Delta(X^\ell,\sP) & := \{ (x_1,\dots,x_{\ell}) \in X^\ell:
 x_i=x_j \
\text{if $\exists\, P\in\sP$ with $i,j\in P$ } \}.
\end{aligned}
\end{equation}
The preceding subspaces are related by
\begin{equation}
\label{eq:Closure}
\cl_{X^\ell}\Delta^\circ(X^\ell,\sP)
=
\Delta(X^\ell,\sP).
\end{equation}
The \emph{big (respectively, small) diagonal}  in $X^\ell$ comprises the set of points $(x_1,\dots,x_\ell)\in X^\ell$ with $x_i=x_j$ for at least one (respectively, all) pair $i,j$. We can identify $\Delta^\circ(X^\ell,\sP)$ with the complement of the big diagonal in $X^r$  (where $r$ is the length of $\sP$) as follows. Write $\sP$ as an ordered collection $P_1,\dots,P_{r(\sP)}$. Then define an embedding,
\begin{equation}
\label{eq:DefineDiagonalComplementEmbedding}
\iota_{\sP}:
X^r\less \ \bigcup_{i<j}\ \{x_i=x_j\} \to \Delta^\circ(X^\ell,\sP)\subset X^\ell,
\end{equation}
by $\iota(x_1,\dots,x_r)=(y_1,\dots,y_\ell)$, where $y_k=x_i$ for $k\in P_i$.  The map $\iota_{\sP}$ is a diffeomorphism.
\begin{lem}
\label{lem:DiagonalStabilizers}
Let $\sP$ be a partition of $N_\ell$. Let $\fS_\ell$ act on $X^\ell$ by the permutation action,
\begin{equation}
\label{eq:SymmetricGroupActionOnXell}
\left( \si,(x_1,\dots,x_n)\right)\mapsto
(x_{\si(1)},\dots,x_{\si(n)}).
\end{equation}
Then the following hold:
\begin{enumerate}
\item
\label{item:DiagonalStabilizers1}
The diagonal $\Delta(X^\ell,\sP)$ is the fixed-point set of the group $\Ga(\sP)$;
\item
\label{item:DiagonalStabilizers2}
For $\bx\in X^\ell$, one has $\bx\in\Delta^\circ(X^\ell,\sP)$ if and only if $\Stab_{\bx}=\Ga(\sP)$;
\item
\label{item:DiagonalStabilizers3}
If $\si\in\fS_\ell$ then $\si(\Delta^\circ(X^\ell,\sP))\subset \Delta^\circ(X^\ell,\sP)$ if and only if $\si\in\fS(\sP)$.
\end{enumerate}
\end{lem}

\begin{proof}
That $\Delta(X^\ell,\sP)$ is contained in the fixed-point set of $\Ga(\sP)$ follows immediately from the definitions. The fixed-point set of $\fS_i$ acting on $X^i$ is the small diagonal.
Thus, by the characterization of  $\Ga(\sP)$ in Lemma \ref{lem:CharacterizeGaP}, if $\bx\in X^\ell$ has the property that $\si(\bx)=\bx$ for all $\si\in\Ga(\sP)$, then $x_i=x_j$ for all $i,j\in P$ for every $P\in\sP$. Hence, the fixed-point set of $\Ga(\sP)$ is contained in $\Delta(X^\ell,\sP)$, proving Item \eqref{item:DiagonalStabilizers1}.

Item \eqref{item:DiagonalStabilizers1} implies that $\Ga(\sP)\subset\Stab_\bx$ for $\bx\in\Delta^\circ(X^\ell,\sP)$. If $\si\in\Stab_\bx$ and $\si\notin\Ga(\sP)$, then there are $P\in\sP$ and $i\in P$ such that $\si(i)\notin P$. But $\si(\bx)=\bx$ then implies that $x_{\si(i)}=x_i$, which contradicts the fact that $\bx\in\Delta^\circ(X^\ell,\sP)$. Therefore, $\bx\in\Delta^\circ(X^\ell,\sP)$ implies that $\Stab_\bx=\Ga(\sP)$. If $\Stab_x=\Ga(\sP)$, then Item \eqref{item:DiagonalStabilizers1} implies that $\bx\in\Delta(X^\ell,\sP)$.  If $\bx\notin\Delta^\circ(X^\ell,\sP)$ then there are $i\in P_a$ and $j\in P_b$, $P_a\neq P_b\in\sP$, with $x_i=x_j$.  Consequently, for $\si=(i\ j)\in\fS_\ell$, one has $\si(\bx)=\bx$ but $\si\notin \Ga(\sP)$, a contradiction, and this proves Item \eqref{item:DiagonalStabilizers2}.

Item \eqref{item:DiagonalStabilizers3} follows immediately from the definition of $\fS(\sP)$ as the normalizer of $\Ga(\sP)$, Item \eqref{item:DiagonalStabilizers2}, and the relation $\Stab_{\si(\bx)}=\si\Stab_\bx\si^{-1}$
(see \cite[p. 3]{tomDieck}).
\end{proof}

\subsection{Strata of the symmetric product}
Let $\tilde\pi_{\ell}:X^\ell\to\Sym^\ell(X)=X^\ell/\fS_\ell$ be the projection. Define the \emph{stratum of the symmetric product $\Sym^\ell(X)$ corresponding to a partition $\sP$} by
\begin{equation}
\label{eq:DefineStratumOfPartition}
\Si(X^\ell,\sP) := \tilde\pi_{\ell}\left(\Delta^\circ(X^\ell,\sP)\right).
\end{equation}
Lemma \ref{lem:DiagonalStabilizers} implies that
\begin{equation}
\label{eq:StratumQuotientOfDiagonal}
\Si(X^\ell,\sP)
= \Delta^\circ(X^\ell,\sP)/\fS(\sP)
= \Delta^\circ(X^\ell,\sP)/W(\sP).
\end{equation}
The result below follows immediately from the definitions.

\begin{lem}
\label{lem:PreImageOfSi}
Let $\sP$ be a partition of $N_\ell$ and let $[\sP]$ denote the orbit of $\sP$ under $\fS_\ell$. Then,
$$
\tilde\pi_{\ell}^{-1}\left(\Si(X^\ell,\sP)\right)
=
\bigcup_{\sP'\in[\sP]}\ \Delta^\circ(X,\sP').
$$
\end{lem}

\section{Incidence relations among diagonals and strata}
\label{subsec:IncidenceRelations}
We now describe incidence relations among the diagonals of $X^\ell$ and the resulting relations among the strata $\Si$ of $\Sym^\ell(X)$.

\begin{defn}
\label{defn:DefineRefinement}
Let $\sP$ and $\sP'$ be partitions of $N_\ell$. We say $\sP'$ is a \emph{refinement} of $\sP$, or $\sP<\sP'$, if for every $P'\in\sP'$ there is a $P\in\sP$ with $P'\subseteqq P$.
\end{defn}

\begin{lem}
\label{lem:IncidenceRelation}
If $\sP$ and $\sP'$ are partitions of $N_{\ell}$ then $\Delta^\circ(X^\ell,\sP)\subset \Delta(X,\sP')$ if and only if $\sP<\sP'$.
\end{lem}

\begin{proof}
If $\sP'$ is a refinement of $\sP$ and $\bx\in\Delta^\circ(X^\ell,\sP)$, then for every $P'\in\sP'$ and $i,j\in P'$, there is a $P\in\sP$ with $i,j\in P'\subset P$.  Since $i,j\in P\in\sP$, we have $x_i=x_j$ so $\bx\in\Delta(X,\sP')$, proving one implication.

Conversely, if $\Delta^\circ(X^\ell,\sP)\subset \Delta(X,\sP')$ and $P'\in\sP'$ then, because every $\bx=(x_1,\dots,x_\ell) \in\Delta^\circ(X^\ell,\sP)$ is also in $\Delta(X,\sP')$, we see that for every $i,j\in P'$, we have $x_i=x_j$.  This implies (by the only if in the definition of $\Delta^\circ(X^\ell,\sP)$) that there must be a $P\in\sP$ with $i,j\in P$. Thus, $P'\subset P$ and $\sP<\sP$.
\end{proof}

Now, we compare the incidence relations with the strata of $\Sym^\ell(X)$.  Recall that for a partition $\sP$ of $N_\ell$, we use $[\sP]$ to denote the orbit of $\sP$ under the action \eqref{eq:ActionOnPartitions} of $\fS_\ell$ on the set of partitions of $N_\ell$. For partitions $\sP$ and $\sP'$ of $N_\ell$, define
\begin{equation}
\label{eq:ConjugateRefinements}
[\sP<\sP']
:=
\{\sP''\in [\sP']: \sP< \sP''\}.
\end{equation}
We now give an example to show that while $\fS(\sP)$ acts  on $[\sP<\sP']$ by the action \eqref{eq:ActionOnPartitions}, this action need  not be transitive.

\begin{exmp}
\label{exmp:ActionOnRelatedPartitionOrbit}
Consider the partitions $\sP=\{P_1,P_2\}$, where $P_1=\{1,2\}$ and $P_2=\{3,4,5\}$, and $\sP'=\{Q_1,Q_2,Q_3,Q_4\}$, where $Q_1=\{1\}$, $Q_2=\{2\}$, $Q_3=\{3\}$, and $Q_4=\{4,5\}$.
Observe that  $\sP<\sP'$ because $Q_1,Q_2\subset P_1$ and $Q_3,Q_4\subset P_2$. Define $\si:=(1\ 4)(2\ 5)\in \fS_5$.  Then $\si(\sP')=\sP''=\{R_1,R_2,R_3,R_4\}$, where $R_1=\{1, 2\}$, $R_2=\{3\}$, $R_3=\{4\}$ and $R_4=\{5\}$.  Because $R_1=P_1$ and $R_2,R_3,R_4\subset P_2$, we have $\sP<\sP''$.  Thus, $\sP',\si(\sP')\in [\sP<\sP']$ but $\si\notin \fS(\sP)$ because $\si$ does not satisfy the conclusion of Lemma \ref{lem:CharNormalizerOfGaP}.
\end{exmp}

The next lemma follows from the definition of $\Si(X^\ell,\sP)$ (see \eqref{eq:DefineStratumOfPartition}) as well as Lemma \ref{lem:IncidenceRelation} and the third assertion of Lemma \ref{lem:PreImageOfSi}.

\begin{lem}
\label{lem:IncidenceRelationInSymmetricProduct}
If $\sP$ and $\sP'$ are partitions of $N_{\ell}$, then
$$
\Si(X^\ell,\sP)\subset \cl_{\Sym^\ell(X)} \Si(X^\ell,\sP')
$$
if and only if there is a $\sP''\in [\sP']$ with $\sP<\sP''$.
\end{lem}

\section{Normal bundles of diagonals and strata}
\label{sec:DiagonalNormals}
We now introduce the normal bundles of the diagonals $\Delta^\circ(X^\ell,\sP)$ of $X^{\ell}$.

\begin{lem}
Let $\sP$ be any partition of $N_{\ell}$. The tangent bundle of the submanifold $\Delta^\circ(X^\ell,\sP)$ is given by
\begin{equation}
\label{eq:TangentToDiagonal}
T\Delta^\circ(X^\ell,\sP)
=
\left\{(v_1,\dots,v_{\ell})\in TX^\ell: v_i=v_j \iff \exists\, P\in\sP \text{ with } i,j\in P \right\}.
\end{equation}
If $g$ is any metric on $X$, then the orthogonal complement of $T\Delta^\circ(X^\ell,\sP)$ with respect to the metric on $X^\ell$ given by the $\ell$-fold product of $g$
is the following normal bundle of $\Delta^\circ(X^\ell,\sP)$,
\begin{equation}
\label{eq:DiagonalNormal}
\tilde\nu(X^\ell,\sP)
:=
\left\{(v_1,\dots,v_{\ell})\in TX^\ell|_{\Delta^\circ(X^\ell,\sP)}: \sum_{i\in P}v_i=0 \text{ for all } P\in\sP \right\}.
\end{equation}
\end{lem}

\begin{proof}
Differentiate any path in $\Delta^\circ(X^\ell,\sP)$ to see the form of $T\Delta^\circ(X^\ell,\sP)$ appearing in equation \eqref{eq:TangentToDiagonal}. For $\vec v=(v_1,\dots,v_\ell)\in T\Delta^\circ(X^\ell,\sP)$, define $v_P=v_i$ for any $i$ with $i\in P$. If $\vec v=(v_1,\dots,v_\ell)\in T\Delta^\circ(X^\ell,\sP)$ and $\vec w=(w_1,\dots,w_\ell)\in \tilde\nu(X^\ell,\sP)$, then the equality,
$$
\vec v\cdot \vec w=\sum_i v_i\cdot w_i
=
\sum_{P\in\sP} v_P\cdot\left( \sum_{i\in P} w_i\right)
=0,
$$
which holds for any product metric on $X^\ell$, implies that the bundle in equation \eqref{eq:DiagonalNormal} is contained in the orthogonal complement of $T\Delta^\circ(X^\ell,\sP)$.
Since the above inclusion must hold for any value of $v_P$, we see that any element of the orthogonal complement of $T\Delta^\circ(X^\ell,\sP)$ must satisfy the equations defining
the bundle in \eqref{eq:DiagonalNormal}.
\end{proof}

Lemma \ref{eq:TangentToDiagonal} leads to the following description of a neighborhood of $\Delta^\circ(X^\ell,\sP)$.

\begin{lem}
\label{lem:Normal}
There is an open neighborhood, $\tilde\sO(X^\ell,\sP)\subset \tilde\nu(X^\ell,\sP)$, of the zero-section and an $\fS(\sP)$-equivariant exponential map identifying
$\tilde\sO(X^\ell,\sP)$ with an open neighborhood, $\tilde\sU(X^\ell,\sP)$, of $\Delta^\circ(X^\ell,\sP)$ in $X^\ell$.
\end{lem}

\begin{rmk}
\label{rmk:XlTopStratum}
If $\sP$ is the partition of $N_\ell$ given by $\ell=1+1+\dots +1$, so $\Delta^\circ(X^\ell,\sP)$ is the complement of the
big diagonal in $X^\ell$, then the normal bundle \eqref{eq:DiagonalNormal} is trivial in the sense that it is equal to $\Delta^\circ(X^\ell,\sP)$. The diagonal $\Delta^\circ(X^\ell,\sP)$ is an open, dense subspace of $X^\ell$ and equal to the open neighborhood $\tilde\sU(X^\ell,\sP)$ of Lemma \ref{lem:Normal}.
\end{rmk}

There is an obvious  $\fS(\sP)$ action on the normal bundle $\tilde\nu(X^\ell,\sP)$.  We will write $\nu(X^\ell,\sP)$ for the quotient $\tilde\nu(X^\ell,\sP)/\fS(\sP)$. Because the exponential map is $\fS(\sP)$-equivariant, equation \eqref{eq:StratumQuotientOfDiagonal} yields

\begin{lem}
\label{lem:NormalInSymmProd}
Continue the notation of Lemma \ref{lem:Normal}. If $\sO(X^\ell,\sP)=\tilde\sO(X^\ell,\sP)/\fS(\sP)$, then a neighborhood,
$$
\sU(X^\ell,\sP) := \tilde\sU(X^\ell,\sP)/\fS(\sP),
$$
of $\Si(X^\ell,\sP)$ in $\Sym^\ell(X)$ is homeomorphic to $\sO(X^\ell,\sP)$.
\end{lem}

We now identify the intersection of $\Delta^\circ(X,\sP')$ with $\tilde\sU(X^\ell,\sP)$.

\begin{lem}
\label{lem:DiagonalInDiagonalNormal}
If $\sP<\sP'$ are partitions of $N_{\ell}$ and
\begin{multline}
\label{eq:DoubleDiagonalNormal}
\tilde\nu(X^\ell,\sP\to\sP')
\\
:=
\left\{(v_1,\dots,v_\ell)\in\tilde\nu(X^\ell,\sP): v_i=v_j \iff \exists\, P'\in\sP'\text{ with } i,j\in P' \right\},
\end{multline}
then the intersection of $\Delta^\circ(X^\ell,\sP')$ with $\tilde\sU(X^\ell,\sP)$ is homeomorphic to
\[
\tilde\sO(X^\ell,\sP\to\sP',g_{\sP}) := \tilde\nu(X^\ell,\sP\to\sP')\cap\tilde\sO(X^\ell,g_{\sP}).
\]
\end{lem}

\begin{proof}
This follows immediately from the equivariance of the exponential map.
\end{proof}

Thus, a neighborhood of the lower diagonal, $\Delta^\circ(X^\ell,\sP)$, in $\Delta^\circ(X^\ell,\sP')$  can be described by the bundle $\tilde\nu(X^\ell,\sP\to\sP')$.  However,
the end of the stratum $\Si(X^\ell,\sP')$ near the lower stratum, $\Si(X^\ell,\sP)$, can be more complicated than a quotient of $\tilde\nu(X^\ell,\sP\to\sP')$ by a subgroup of $\fS(\sP)$.
Indeed, $\tilde\nu(X^\ell,\sP\to\sP')$ need not be an invariant subspace of $\tilde\nu(X^\ell,\sP)$ under the action of $\fS(\sP)$. Recall from Lemma \ref{lem:PreImageOfSi}
that $\tilde\pi_{\ell}^{-1}(\Si(X^\ell,\sP'))$ comprises
not just $\Delta^\circ(X^\ell,\sP')$ but rather all the diagonals given by a partition in $[\sP']$, the orbit of $\sP'$ under the action
\eqref{eq:ActionOnPartitions}. Therefore, while we can cover $\Si(X^\ell,\sP)$ by $\Delta^\circ(X^\ell,\sP)$, to account for all the ends of $\Delta^\circ(X^\ell,\sP')$ which
get mapped by $\tilde\pi_\ell$ to a neighborhood of the lower stratum, $\Si(X^\ell,\sP)$, we must consider all diagonals $\Delta^\circ(X^\ell,\sP'')$ where
$\sP''\in[\sP']$ and $\sP< \sP''$.  That is, we must consider all partitions $\sP''\in [\sP<\sP']$, where $[\sP<\sP']$ is defined in \eqref{eq:ConjugateRefinements}. The next lemma follows from Lemmas \ref{lem:Normal} and \ref{lem:DiagonalInDiagonalNormal}.

\begin{lem}
\label{lem:EndOfUpperStratum}
Let $\sP<\sP'$ be partitions of $N_{\ell}$. Let $[\sP<\sP']$ be the set of partitions defined in \eqref{eq:ConjugateRefinements}. Then a neighborhood of $\Si(X^\ell,\sP)$ in $\Si(X^\ell,\sP')$ is homeomorphic to a neighborhood of the zero-section in
\begin{equation}
\label{eq:EndOfUpperStratum}
\nu(X^\ell,\sP\to [\sP'])
:=
\tilde\nu(X^\ell,\sP\to [\sP'])/\fS(\sP)
\end{equation}
where
$$
\tilde\nu(X^\ell,\sP\to [\sP'])
:=
\bigsqcup_{\sP''\in [\sP<\sP']}\tilde\nu(X^\ell,\sP\to\sP'').
$$
\end{lem}

\begin{rmk}
As noted in Example \ref{exmp:ActionOnRelatedPartitionOrbit}, $\fS(\sP)$ need not act transitively on $[\sP<\sP']$. The orbits of $\fS(\sP)$ in $[\sP<\sP']$ will enumerate
the components of the end of $\Si(X^\ell,\sP')$ in a neighborhood of $\Si(X^\ell,\sP)$.
\end{rmk}

\section{Enumeration of the strata}
\label{subsec:EnumStrata}
We will often need to give arguments by induction on strata of $\Sym^\ell(X)$. To that end, we now give a method of enumerating these strata. Choose a representative, $\sP$, from each orbit of partitions, $[\sP]$, under the action \eqref{eq:ActionOnPartitions}.  We refer to the partition $\sP_0=\{N_\ell\}$ as the \emph{crudest} partition.
Enumerate the remaining representatives in such a way that
$$
\Si(X^\ell,\sP_i)\subset\cl\left(\Si(X^\ell,\sP_j)\right) \quad\text{only if}\quad i<j.
$$
Write these partitions as $\sP_0,\sP_2,\dots,\sP_r$ so that $\Si(X^\ell,\sP_0)$ is the lowest stratum and $\Si(X^\ell,\sP_r)$ is the highest.

\chapter{A partial Thom--Mather structure on symmetric products}
\label{chap:diagTM}
In this chapter, we define neighborhoods of the strata of $\Sym^\ell(X)$ 
defined in Section \ref{subsec:DefiningTheDiagonals}, projections of these neighborhoods onto the strata,
and vector-valued tubular distance functions on these neighborhoods.

\section{Introduction}
\label{sec:Introduction_diagTM}
A topological space $Z$ is \emph{smoothly stratified}\label{smoothly_stratified_space} if it admits a decomposition into a locally finite collection of locally closed subspaces, $Z=\cup_i \Si_i$, each of which is a smooth manifold and which satisfy the \emph{condition of the frontier}\label{condition_frontier}: $\Si\cap\cl(\Si_j)\neq \emptyset$ if and only if $\Si\subset \cl(\Si_j)$.  Such spaces are called Whitney pre-stratified in \cite[p. 480]{Mather_2012} (see also, \cite[p. 36]{GorMacPh}).
A \emph{Thom--Mather} stratification\label{Thom-Mather_stratification} on a smoothly-stratified space $Z$ with strata $\Si_i$ consists of
\emph{tubular neighborhood structures}\label{tubular_neighborhood_structure}, by which we mean neighborhoods, $U_i$ of $\Si_i$, with \emph{tubular neighborhood projections}\label{tubular_neighborhood_projection}, $\pi_i:U_i\to\Si_i$, and \emph{tubular distance functions}\label{tubular_distance_function}, $t_i:U_i\to [0,1)$, satisfying the following properties:
\begin{enumerate}
\item
For all $i$, $t_i^{-1}(0)=\Si_i$.
\item
For $\Si_i\subset\cl_Z(\Si_j)$ the \emph{compatibility conditions}\label{compatibility_condition}
\begin{equation}
\label{eq:TMProperties}
\begin{aligned}
\pi_i\circ\pi_j&= \pi_i,
\\
t_i\circ \pi_j&=t_i,
\end{aligned}
\end{equation}
wherever these compositions are defined.
\item
The map $\pi_i\times t_i:U_i\to\Si_i\times[0,1)$ is a submersion.
\end{enumerate}
These conditions are called the \emph{control conditions}\label{control_condition} of a Thom--Mather stratification (see \cite[p. 42]{GorMacPh}). Because $\Sym^\ell(X)$ is the quotient of the smooth manifold $X^\ell$ by the action of
a finite group, it admits the structure of a \emph{Whitney stratification}\label{Whitney_stratification} (for example, see \cite[Theorem 4.3.7]{Pflaum})  and hence a Thom--Mather stratification (see \cite[p. 42]{GorMacPh}).

In this chapter, we give an explicit construction of a \emph{partial tubular neighborhood structure}\label{partial_tubular_neighborhood_structure}
on the stratification  of $\Sym^\ell(X)$ defined in Section \ref{subsec:DefiningTheDiagonals}.  By partial tubular neighborhood structure on a smoothly-stratified space, $Z$, with strata, $\Si_i$, we mean neighborhoods, $U_i$ of $\Si_i$, with projection maps, $\pi_i:U_i\to\Si_i$, and \emph{vector-valued tubular distance functions}\label{vector-valued_tubular_distance_function}, $\vec t_i:U_i\to [0,1)^{n_i}$, with $\Si_i=\vec t_i^{-1}(\vec 0)$. We show these maps satisfy the first compatibility condition in \eqref{eq:TMProperties} and how to control the failure of the second compatibility condition.

We shall need the explicit construction of these projection and distance functions to compare them with the functions of a partial tubular neighborhood structure on neighborhoods of $M_{\fs}\times\Sym^\ell(X)$. In addition, this presentation serves as an introduction to a technique for constructing partial tubular neighborhood structures which we shall employ again in later chapters in less familiar settings.

To construct the partial tubular neighborhood structure on $\Sym^\ell(X)$, we begin by constructing one on $X^\ell$ for the stratification given by the diagonals of $X^\ell$. In Section \ref{sec:R4Diags}
we define subspaces of $(\RR^4)^\ell$ which will appear as the fibers of the normal bundles defined in Section \ref{sec:DiagonalNormals}.  We also define stratifications of these subspaces of $(\RR^4)^\ell$ and tubular neighborhoods of these stratifications. To use the normal bundles of Section \ref{sec:DiagonalNormals} to construct partial tubular neighborhood structures on $\Sym^\ell(X)$, we define an exponential map with a varying metric in Section \ref{sec:FamiliesOfMetric}. We describe the overlap of the exponential maps from these normal bundles in Section \ref{sec:OverlapMaps}.  In Section \ref{sec:TMMapsOnSymX} we construct families of metrics so that the projection maps defined by the associated exponential maps satisfy the first equation in \eqref{eq:TMProperties}. Then, we show how this partial tubular neighborhood structure on $X^\ell$ descends to one on $\Sym^\ell(X)$ in Section \ref{sec:NormalBundlesOfSymlXStrata}. In Section \ref{sec:SymlXTubDistance}, we define the vector-valued tubular distance functions and describe their failure to satisfy the second equation in \eqref{eq:TMProperties}.  This description and the analogous description for the vector-valued tubular distance
functions for $M_\fs\times\Sym^\ell(X)$ will allow us to define the link of $M_\fs\times\Sym^\ell(X)$
in \eqref{eq:DefineLinkStratum} as a union of fiber bundles as described in Proposition \ref{prop:LinkPieceFiberBundleStructure}
which will enable the pushforward-pullback arguments of Sections \ref{sec:BundlesPushforwards},
\ref{sec:Computations} and \ref {sec:Proofs}.
Finally, in Section \ref{sec:SymlXStrataDecomp}, we use the partial tubular neighborhood structure to produce a partition of $\Sym^\ell(X)$ into compact subsets of these tubular neighborhoods.  These compact subsets will
be used in our construction of the link of $M_\fs\times\Sym^\ell(X)$ and to give a criterion for patching together maps to $\Sym^\ell(X)$ which are homotopic but not equal.

\section{Diagonals in products of $\RR^4$}
\label{sec:R4Diags}
The fiber of the normal bundle of the diagonal $\Delta^\circ(X^\ell,\sP)$ is given
by the product, over $P\in\sP$, of the subspace of mass-centered points in $\oplus_{i\in P}\RR^4$ (see \eqref{eq:DiagonalNormal}).
Hence,
for any subset $P$ of $N_\ell$, define
\begin{equation}
\label{eq:DefineZ}
z_P: \bigoplus_{i\in P} \RR^4 \ni (v_i)_{i\in P} \mapsto \sum_i v_i \in \RR^4
\quad\text{and}\quad
Z_P := z_P^{-1}(0)\subset \bigoplus_{i\in P} \RR^4.
\end{equation}
We will write $z_\ell$ for $z_{N_\ell}$.  We define the \emph{cone point},
\begin{equation}
\label{eq:DefineR4DiagonalConePoint}
c_P\in Z_P,
\end{equation}
to be the zero vector in $\oplus_{i\in P}\ \RR^4$. Note that if $|P|=1$, then $Z_P=\{c_P\}$. The fiber of $\tilde\nu(X^\ell,\sP)$ is then given by (compare \eqref{eq:DiagonalNormal})
\begin{equation}
\label{eq:DiagonalNormalFiber}
Z(\sP)= \prod_{P\in\sP} Z_P.
\end{equation}
For partitions $\sP<\sP'$ of $N_\ell$, define a partition of each $P\in\sP$ by
\begin{equation}
\label{eq:DefineRefinementPartition}
\sP'_P := \{P'\in\sP': P'\subseteqq P\}.
\end{equation}
If $P\in\sP\cap\sP'$, then $\sP'_P=\{P\}$. We define
\begin{equation}
\Delta^\circ(Z_P,\sP'_P)
:=
\{(v_i)_{i\in P}\in Z_P: \text{$v_i=v_j$ if and only if $\exists\, P'\in\sP'_P$ with $i,j\in P'$} \},
\label{eq:DefineCenteredR4Diagonal}
\end{equation}
and observe that the fiber of the bundle $\tilde\nu(X^\ell,\sP\to\sP')$ in \eqref{eq:DoubleDiagonalNormal} is given by
\begin{equation}
\label{eq:DoubleDiagNormalFiber}
Z(\sP,\sP')
=
\prod_{P\in\sP} \Delta^\circ(Z_P,\sP'_P).
\end{equation}
Note that if $\sP'_P$ is given by one set, $P$, as is the case when $P\in\sP\cap\sP'$, then $\Delta^\circ(Z_P,\sP'_P)$ is a single point, given by $|P|$ copies of the zero vector in $\RR^4$.

We define a \emph{cone parameter function} on $Z_P$ by the map
\begin{equation}
\label{eq:DefineInitialConeParameter}
t_{P,0}:Z_P \ni (v_i)_{i\in P} \mapsto \sum_{i\in P} |v_i|^2 \in [0,\infty).
\end{equation}
Then $t_{P,0}^{-1}(0)=c_P$ and there is a deformation retraction of $Z_P$ to $c_P$ given by
\begin{equation}
\label{eq:ConeRetractionMap}
r_P: Z_P\times [0,1] \ni ((v_i)_{i\in P},s) \mapsto (sv_i)_{i\in P} \in Z_P.
\end{equation}
Observe that $r_P(\cdot,1)$ is the identity map, while $r_P(\cdot,0)=c_P$.  Moreover,
\begin{equation}
t_{P,0}\circ r_P(\cdot,s)=s^2 t_{P,0}(\cdot).
\end{equation}
Hence, there is a homeomorphism between $Z_P$ and the cone on $t_{P,0}^{-1}(1)$.

The following lemma describes the normal bundle of $\Delta^\circ(Z_P,\sP'_P)\subset z_P^{-1}(0)$.

\begin{lem}
\label{lem:RelatingTrivialStrata}
Let $\sP'_P$ be a partition of $P$. Then the restriction of the map
\begin{equation}
\label{eq:NormalOfR4Diagonal}
e(Z_P,\sP'_P):
\Delta^\circ(Z_P,\sP'_P) \times
\prod_{P'\in \sP'_P} Z_{P'}
\to
Z_P
\end{equation}
defined by
\begin{equation}
\label{eq:DefineR4DiagonalExponential}
\left(
(v_i)_{i\in P},(z_j)_{\begin{subarray}{c}j\in P', \\ P'\in\sP_P' \end{subarray}}
\right)
\mapsto (v_i+z_i)_{i\in P}
\end{equation}
to a sufficiently small neighborhood $\sO(Z_P,\sP'_P)$ of the subspace
$$
\Delta^\circ(Z_P,\sP'_P) \times \prod_{P'\in \sP'_P}  \{c_{P'}\}
$$
is a homeomorphism onto its image, $\sU(Z_P,\sP'_P)$.
\end{lem}

\begin{proof}
We observe that the map \eqref{eq:DefineR4DiagonalExponential}
can be inverted by the \emph{center of mass map},
\begin{equation}
\label{eq:R4CenterOfMassProjection}
\begin{aligned}
(w_i)_{i\in P} \to & \left(
(v_i)_{i\in P},(z_i)_{\begin{subarray}{c}j\in P', \\ P'\in\sP_P' \end{subarray}}
\right),
\\
v_{i}&=\frac{1}{|P'|}\sum_{j\in P'}w_j,
\quad
\text{for $i\in P'$},
\\
z_j&=w_j-v_j.
\end{aligned}
\end{equation}
(Note that the above map is an inverse because $\sum_{j\in P'} z_j=0$
by the definition of $Z_{P'}$.)  However, the map
\eqref{eq:R4CenterOfMassProjection} only takes values
in
$$
\Delta^\circ(Z_P,\sP'_P) \times \prod_{P'\in \sP'_P} Z_{P'}
$$
if we restrict to a
sufficiently small neighborhood $\sO(Z_P,\sP'_P)$.
\end{proof}

\begin{rmk}
\label{rmk:ZPTopStratumCase}
As in Remark \ref{rmk:XlTopStratum}, if $\sP'_P$ is the finest partition of $P$, made of subsets all of size one, then for all $P'\in\sP'_P$, the point $Z_{P'}$ is the cone point and $e(Z_P,\sP'_P)$ is the identity map.  In this case, $\Delta^\circ(Z_P,\sP'_P)=\sO(Z_P,\sP'_P)$ is an open and dense subset of $Z_P$.
\end{rmk}

\begin{rmk}
\label{rmk:EquivarianceOfConeExponential}
Observe that the map $e(Z_P,\sP'_P)$ is equivariant with respect to the $\RR^+$ actions given by scalar multiplication on all factors in the domain,
$$
\left(
(v_i)_{i\in P},(z_j)_{\begin{subarray}{c}j\in P', \\ P'\in\sP_P' \end{subarray}}
\right)
\mapsto
\left(
(sv_i)_{i\in P},(sz_j)_{\begin{subarray}{c}j\in P', \\ P'\in\sP_P' \end{subarray}}
\right),
$$
and the map $r_P(\cdot,s)$ on $Z_P$.
\end{rmk}

Because $Z_P$ is a completely normal topological space (see \cite[Section 32, Example 4]{Munkres2ndEdition}), we can find open neighborhoods $U(Z_P,\sP'_P)$ of $\Delta^\circ(Z_P,\sP'_P)$ for every partition $\sP'_P$ of $P$ such that
$U(Z_P,\sP'_P)\cap U(Z_P,\sP''_P)$ is empty unless $\sP'_P<\sP''_P$ or $\sP''_P<\sP'_P$.  By Remark \ref{rmk:EquivarianceOfConeExponential}, we can assume that the neighborhoods $U(Z_P,\sP'_P)$ are closed under the maps $r_P(\cdot,s)$ for $s\in [0,1]$. The normal bundle structure given by the map \eqref{eq:DefineR4DiagonalExponential} then defines a tubular neighborhood projection,
\begin{equation}
\label{eq:R4DiagonalProjection}
\pi(Z_P,\sP'_P): U(Z_P,\sP'_P)\to \Delta^\circ(Z_P,\sP'_P),
\end{equation}
given by the $v_i$ component of the map \eqref{eq:R4CenterOfMassProjection}:
\begin{equation}
\label{eq:DefineR4DiagonalProjection}
(w_i)_{i\in P} \mapsto (v_i)_{i\in P},
\quad\text{where } v_i := \frac{1}{|P'|}\sum_{j\in P'}w_j \text{ for } i\in P'\in\sP'_P.
\end{equation}
In the case of the top stratum described in Remark \ref{rmk:ZPTopStratumCase}, the map $\pi(Z_P,\sP'_P)$ is just the identity map. By construction, the maps $\pi(Z_P,\sP'_P)$ commute with the map $r_P$ of \eqref{eq:ConeRetractionMap} in the sense that
\begin{equation}
\label{eq:ConeRetractionCommuteWithProjection}
r_P(\pi(Z_P,\sP'_P)(\cdot),s)
=
\pi(Z_P,\sP'_P)(r_P(\cdot,s)),
\end{equation}
for all $s\in [0,1]$.

The following lemma shows that the maps $\pi(Z_P,\sP'_P)$ satisfy the first equality in \eqref{eq:TMProperties}.

\begin{lem}
If $\sP'_P<\sP''_P$ are partitions of $P$, then on the intersection
$U(Z_P,\sP'_P)\cap U(Z_P,\sP''_P)$, the following equality holds:
\begin{equation}
\label{eq:R4DiagonalTM1}
\pi(Z_P,\sP'_P)\circ \pi(Z_P,\sP''_P) = \pi(Z_P,\sP'_P).
\end{equation}
\end{lem}

\begin{proof}
The lemma follows immediately from the definition
\eqref{eq:DefineR4DiagonalProjection}.
\end{proof}

We use the above set of tubular neighborhoods to define a deformation of the scale function,
$t_{P,0}$, which is constant on the fibers of the maps $\pi(Z_P,\sP'_P)$.

\begin{lem}
\label{lem:R4TMScale}
There are
\begin{enumerate}
\item
\label{item:R4TMScale_function}
A smooth function, $t(Z_P):t_{P,0}^{-1}([0,1])\subset Z_P\to [0,\infty)$,
\item
\label{item:R4TMScale_neighborhoods}
Open neighborhoods, $U'(Z_P,\sP'_P)\sqsubset U(Z_p,\sP'_P)\cap t_{P,0}^{-1}([0,1])$, of the diagonals, $\Delta^\circ(Z_P,\sP'_P)\cap t_{P,0}^{-1}([0,1])$,
\end{enumerate}
such that the following hold:
\begin{enumerate}
\item
\label{item:R4TMScale_one}
If $c_P\in Z_P$ is the cone point in \eqref{eq:DefineR4DiagonalConePoint}, then $t(Z_P)^{-1}(0)=c_P$.
\item
\label{item:R4TMScale_two}
If $r_P$ is the map in \eqref{eq:ConeRetractionMap} then, for all $s\in [0,1]$,
\begin{equation}
\label{eq:RetractionAndScale}
t(Z_P)(r_P(\cdot,s))=s^2t(Z_P)(\cdot).
\end{equation}
\item
\label{item:R4TMScale_three}
For every partition $\sP'_P$ of $P$ of length greater than one,
\begin{equation}
\label{eq:R4ZPTM2}
t(Z_P)\circ \pi(Z_P,\sP'_P)=t(Z_P) \quad\text{on } U'(Z_P,\sP'_P).
\end{equation}
\item
\label{item:R4TMScale_four}
The function,
$t(Z_P)$,
is invariant with respect to the permutation action of $\fS_P$ and the diagonal rotation action of
$\SO(4)$ on $Z_P$.
\end{enumerate}
\end{lem}

\begin{proof}
Enumerate the conjugacy classes of partitions of $P$, $[\sP_0],\dots,[\sP_n]$ in the manner described in
Section \ref{subsec:EnumStrata}.  That is, for $\sP\in [\sP_i]$ and $\sP'\in [\sP_j]$, we have
$\Delta(Z_P,\sP)\subset \cl\, \Delta(Z_P,\sP')$ only if $i<j$. We will construct a sequence of functions, $t_0,\dots,t_n$, such that $t_i$ satisfies \eqref{eq:R4ZPTM2} with respect to $\pi(Z_P,\sP)$ for any partition $\sP\in [\sP_j]$ with $j<i$.

Begin with $t_0=t_{P,0}$.  This function is invariant with respect to the actions of $\SO(4)$ and $\fS_P$
and satisfies \eqref{eq:RetractionAndScale}.

Assume we have defined functions, $t_i:Z_P\to [0,\infty)$, satisfying Items \eqref{item:R4TMScale_one}, \eqref{item:R4TMScale_two}, and \eqref{item:R4TMScale_four} of the lemma and satisfying \eqref{eq:R4ZPTM2} for all partitions $\sP \in [\sP_j]$  with $j<i$.

There are neighborhoods
$$
\sU_1\sqsubset \sU_2 \sqsubset \sU(Z_P,\sP_i)
$$
of $\Delta^\circ(Z_P,\sP_P)$ which are closed under the $\SO(4)$ action and under the action of $\RR^{+}$ given by $r_P(\cdot,s)$.  Let $\sV_i$ be the orbit of $\sU_i$ under the action of $\fS_P$.  By shrinking $\sU_i$ if necessary, we can assume that either $\sU_i= \si(\sU_i)$ or $\sU_i\cap \si(\sU_i)=\emptyset$ for all $\si\in\fS_P$. Let $\beta_i:Z_P\to [0,1]$ be a function, invariant under the actions of $\SO(4)$, $\fS_P$, and $\RR^{+}$, supported on $\sV_2$, and satisfying $\sV_1\subset \beta^{-1}(1)$. Denote
$$
\pi_i := \coprod_{\sP\in [\sP_i]} \pi(Z_P,\sP):\sV_i\to \bigsqcup_{\sP\in [\sP_i]}\Delta^\circ(Z_P,\sP).
$$
Observe that the map, $\pi_i$, is equivariant with respect to the actions of $\SO(4)$, $\fS_P$, and $\RR^{+}$.
On $\sV_2$, we define
$$
t_{i+1}(\cdot)=\beta_i(\cdot)(t_i\circ \pi_i)(\cdot)+(1-\beta_i(\cdot))t_i(\cdot).
$$
Extend $t_{i+1}$ to $Z_P$ by setting $t_{i+1} := t_i$ on $Z_P-\sV_2$.

By the invariance of $\beta_i$ and $t_i$ with respect to the actions of $\SO(4)$ and $\fS_P$, and the equivariance of $\pi_i$ with respect to these actions, the function $t_{i+1}$ is invariant with respect to these actions.
The same argument shows that $t_{i+1}$ satisfies \eqref{eq:RetractionAndScale}, that is, $t_{i+1}(r_P(\cdot,s))=s^2 t_{i+1}(\cdot)$.

For $j<i$ and $\sP\in [\sP_j]$, the function $t_i$ is constant on the fibers of $\pi(Z_P,\sP)$ by our induction hypothesis. The relation \eqref{eq:R4DiagonalTM1} implies that on $\sV_2\cap \sU'(X_P,\sP)$, the fibers of $\pi_i$ are contained in those of $\pi(Z_P,\sP)$ and thus $t_i\circ \pi_i$ equals $t_i$ on $\sV_2\cap \sU'(X_P,\sP)$. This implies that on $\sV_2\cap \sU'(X_P,\sP)$ any convex linear combination of $t_i\circ\pi_i$ and $t_i$ equals $t_i$.
Hence for any $j<i$ and any $\sP\in [\sP_j]$, the restriction of $t_{i+1}$ to $\sU(X_P,\sP)$ equals the restriction of $t_i$. By definition, this equality also holds on the complement of $\sV_2$. Thus, $t_{i+1}$ satisfies \eqref{eq:R4ZPTM2} for any $\sP\in [\sP_j]$ and $j<i$. We further observe that on $\sV_1\subset \beta_i^{-1}(1)$, we have $t_{i+1}=t_i\circ\pi_i$ while the equality
$\pi_i\circ\pi_i=\pi_i$ implies that
$$
(t_{i+1}\circ\pi_i)(\cdot)
=
(\beta_i\circ\pi_i)(\cdot)(t_i\circ \pi_i)(\cdot)
+
(1-(\beta_i\circ\pi_i)(\cdot))(t_i\circ \pi_i)(\cdot)
=
(t_i\circ \pi_i)(\cdot).
$$
Hence, if we define $\sU'(Z_P,\sP)=\sV_1\cap \sU(Z_P,\sP)$ for $\sP\in [\sP_i]$, then \eqref{eq:R4ZPTM2} holds, completing the induction.
\end{proof}

\section{Families of metrics}
\label{sec:FamiliesOfMetric}
A Riemannian metric, $g$, on $X$ determines a product Riemannian metric on $X^\ell$ and thus exponential maps of the normal bundles of the diagonal. We will use a non-standard diffeomorphism of a neighborhood of the zero section in
the normal bundle $\tilde\nu(X^\ell,\sP)$ with a neighborhood of $\Delta^\circ(X^\ell,\sP)$ in $X^\ell$.  This diffeomorphism will be defined by varying the Riemannian metric on $X$ and hence the resulting product metric on $X^\ell$. We begin by noting that it is possible to define such an exponential map with a varying metric.

\begin{lem}
\label{lem:VaryingMetricExponentialMap}
If $g_{\sP}$ is a smooth family of metrics on $X$ parameterized by $\bx\in\Delta^\circ(X^\ell,\sP)$, then there is an open neighborhood $\tilde\sO(X^\ell,g_{\sP})\subset\tilde\nu(X^\ell,\sP)$ of the zero section and a smooth embedding,
\begin{equation}
\label{eq:VaryingMetricExponentialMap}
e(X^\ell,g_{\sP}):\tilde\sO(X^\ell,g_{\sP})\subset\tilde\nu(X^\ell,\sP)\to X^\ell,
\end{equation}
parameterizing an open neighborhood, $\tilde\sU(X^\ell,g_{\sP})$, of $\Delta^\circ(X^\ell,\sP)$ in $X^\ell$ such that the restriction of the map $e(X^\ell,g_{\sP})$ to a fiber $\tilde\sO(X^\ell,g_{\sP})|_\bx$ is given by the exponential map defined by the $\ell$-fold product of the metric $g_{\sP,\bx}$.
\end{lem}

\begin{proof}
We claim that the derivative of the map $e(X^\ell,g_{\sP})$ at a point $p$ on the zero-section of $\tilde\nu(X^\ell,\sP)$ is the identity.  We identify this zero section with $\Delta^\circ(X^\ell,\sP)$ throughout this discussion. Over a point $p$ in this zero section, the tangent bundle of $\tilde\nu(X^\ell,\sP)$ admits a decomposition into the direct sum of $\tilde\nu(X^\ell,\sP)|_p$ and $T_p\Delta^\circ(X^\ell,\sP)$. By definition of the exponential map, $D e(X^\ell,g_{\sP})_p(v)=v$ for any $v\in\tilde\nu(X^\ell,\sP)|_p$. To evaluate $D e(X^\ell,g_{\sP})_p(w)$ for $w\in T_p\Delta^\circ(X^\ell,\sP)$, let $\ga:[0,1]\to \Delta^\circ(X^\ell,\sP)$ be a smooth path with $\ga'(0)=w$.  Then for all $t\in [0,1]$, we have $e(X^\ell,g_{\sP})(\ga(t))=\ga(t)$.  Differentiating the preceding equality yields $D e(X^\ell,g_{\sP})_p(w)=w$, so $D e(X^\ell,g_{\sP})_p$ is
the identity.  Hence, $e(X^\ell,g_{\sP})$ defines an embedding of a neighborhood of the zero section.
\end{proof}

We can assume that the neighborhood $\tilde\sO(X^\ell,g_{\sP})$
is  small enough to have the following property.

\begin{lem}
\label{lem:DisjointBallsInVaryingMetricExpMap}
Continue the hypotheses of Lemma \ref{lem:VaryingMetricExponentialMap}. The open neighborhood $\tilde\sO(X^\ell,g_{\sP})$ can be chosen so that the images of $\pi_i\circ e(X^\ell,g_{\sP})$ and
$\pi_j\circ e(X^\ell,g_{\sP})$ are disjoint for $i\neq j$, where $\pi_i:X^\ell\to X$ is projection onto the $i$-th factor.
\end{lem}

\begin{proof}
For $(x_P)_{P\in\sP}\in \Delta^\circ(X^\ell,\sP)$, the points $x_P$ are distinct.  Hence, we can chose balls around the origin in $T_{x_P}X$ whose image under the exponential map are disjoint, thus satisfying the conclusion of the lemma.
\end{proof}

The image of the exponential map in \eqref{eq:VaryingMetricExponentialMap} is a tubular neighborhood of $\Delta^\circ(X^\ell,\sP)$. Such a tubular neighborhood, $\tilde\sU(X^\ell,g_{\sP})\subset X^\ell$, defines a
projection,
\begin{equation}
\label{eq:DefineXlTubularProjection}
\pi(X^\ell,g_{\sP}):\tilde\sU(X^\ell,g_{\sP})\to \Delta^\circ(X^\ell,\sP),
\end{equation}
by the composition of  $e(X^\ell,g_{\sP})^{-1}$ and the projection $\tilde\nu(X^\ell,\sP)\to\Delta^\circ(X,\sP)$.

\begin{rmk}
\label{rmk:XlTopStratumExpMap}
If $\sP$ is the finest partition as described in Remark \ref{rmk:XlTopStratum}, then the exponential map $e(X^\ell,g_{\sP})$ and the projection map $\pi(X^\ell,g_{\sP})$ both equal the identity map.
\end{rmk}

For any smooth manifold, $M$, of dimension $n$, let $\Fr_{\GL}(TM)$ denote the $\GL(n)$-bundle of frames of tangent vectors. Define
$$
\Fr_{\GL}(TX^\ell,\sP)
:=
\{(F_1,\dots,F_\ell)\in \Fr(TX^\ell)|_{\Delta^\circ(X^\ell,\sP)}:
F_i=F_j\ \text{for all $i,j\in P\in\sP$}\}.
$$
If $g_{\sP}$ is a family of metrics parameterized by $\Delta^\circ(X^\ell,\sP)$, then we define
$\Fr(TX^\ell,\sP,g_{\sP})$ to be the subbundle of $\Fr_{\GL}(TX^\ell,\sP)$ of frames such that if
$(F_1,\dots,F_\ell)\in \Fr(TX^\ell,\sP,g_{\sP})|_{\bx}$, then $F_i$ is orthonormal with respect to $g_{\sP,\bx}$. The structure group of $\Fr(TX^\ell,\sP,g_{\sP})$ is
\begin{equation}
\label{eq:TgBundleDiagonalStructureGrp}
\tG(T,\sP):=\SO(4)^{r(\sP)}=\prod_{P\in\sP}\SO(4),
\quad\text{where $\sP=\{P_1,\dots,P_{r(\sP)}\}$}.
\end{equation}
Recalling the definition of $Z(\sP)$ from \eqref{eq:DiagonalNormalFiber},
we have the identity,
\begin{equation}
\label{eq:DiagonalNormalAssocBundle}
\tilde\nu(X^\ell,\sP)
=
\Fr(TX^\ell,\sP,g_{\sP})\times_{\tG(T,\sP)}Z(\sP).
\end{equation}
Lemma \ref{lem:DiagonalInDiagonalNormal} implies the following.

\begin{lem}
\label{lem:RestricedExponentialMapDomain}
Let $\sP<\sP'$ be partitions of $N_\ell$. Let $\tilde\nu(X^\ell,\sP\to\sP')$ be the bundle defined in \eqref{eq:DoubleDiagonalNormal}. Then the restriction of the exponential map $e(X^\ell,g_{\sP})$ to
\begin{equation}
\label{eq:RestrictedExponentialMapDomain}
\tilde\sO(X^\ell,g_{\sP})
\cap
\tilde\nu(X^\ell,\sP\to\sP')
\end{equation}
takes values in $\Delta(X,\sP')$.
\end{lem}

We have the following relation between the bundles $\Fr(TX^\ell,\sP,g_{\sP})$ for different partitions, $\sP$.

\begin{lem}
\label{lem:RelateTXFiberBundles}
Let $\sP<\sP'$ be partitions of $N_\ell$. Let $\tilde\sU(X,\sP\to\sP',g_{\sP})\subset\Delta^\circ(X^\ell,\sP')$ be the image of the restriction of the exponential map $e(X^\ell,g_{\sP})$ to the domain in
\eqref{eq:RestrictedExponentialMapDomain}. Let $g_{\sP}$ and $g_{\sP'}$ be smooth families of metrics parameterized
by $\Delta^\circ(X^\ell,\sP)$ and $\Delta^\circ(X,\sP')$ respectively. Assume that $g_{\sP',\bx'}=g_{\sP,\bx}$ if $\pi(X^\ell,g_{\sP})(\bx')=\bx$. Then,
\begin{multline}
\label{eq:TangentBundleOnEnd}
\Fr(TX^\ell,\sP',g_{\sP'})|_{\tilde\sU(X,\sP\to\sP',g_{\sP})}
\\
\cong
\Fr(TX^\ell,\sP,g_{\sP})
\times_{\tilde G(T,\sP)}
\prod_{P\in\sP}
\left(
\Delta^\circ(Z_P,\sP'_P)
\times
\tilde G(T,\sP'_P)
\right).
\end{multline}
\end{lem}

\begin{proof}
By Lemma \ref{lem:RestricedExponentialMapDomain}, a neighborhood of $\Delta^\circ(X^\ell,\sP)$ in
$\Delta^\circ(X^\ell,\sP')$ is diffeomorphic to a neighborhood of the zero section in
\begin{equation}
\label{eq:NormalData}
\tilde\nu(X^\ell,\sP\to\sP')
\cong
\Fr(TX^\ell,\sP,g_{\sP})
\times_{\tilde G(T,\sP)}
\prod_{P\in\sP}
\Delta^\circ(Z_P,\sP'_P)
\subset
\nu(X^\ell,\sP,g_{\sP})
\end{equation}
by the exponential map $e(X^\ell,g_{\sP})$. The isomorphism \eqref{eq:TangentBundleOnEnd} is then defined by the obvious parallel translation of the frames in $\Fr(TX^\ell,\sP,g_{\sP})$ to the point in $\tilde\sU(X^\ell,\sP\to\sP')$ given by the data in \eqref{eq:NormalData}.
\end{proof}

\section{Overlap maps}
\label{sec:OverlapMaps}
Let $\sP<\sP'$ be partitions of $N_\ell$. To define tubular neighborhoods and projection maps which satisfy the conditions \eqref{eq:TMProperties}, we define a space to control the overlap of the tubular neighborhoods, $\sU(X,\sP,g_{\sP})$ and $\sU(X,\sP',g_{\sP'})$, as follows. We begin by defining an overlap fiber bundle,
\begin{equation}
\label{eq:OverlapNormalBundle}
\tilde\nu(X^\ell,\sP,\sP')
:=
\Fr(TX^\ell,\sP,g_{\sP})
\times_{\tilde G(T,\sP)}
\prod_{P\in\sP}
\left(
\Delta^\circ(Z_P,\sP'_P)
\times
\prod_{P'\in\sP'_P} Z_{P'}
\right).
\end{equation}
The cone points, $c_{P'}\in Z_{P'}$, define a subspace of $\tilde\nu(X^\ell,\sP,\sP')$ which is identified with
$\tilde \nu(X^\ell,\sP\to\sP')$ by \eqref{eq:NormalData}. Neighborhoods in $\Delta^\circ(Z_P,\sP_P')$ of the cone points $c_P\in \Delta(Z_P,\sP'_P)$ define `neighborhoods of the zero section' in $\tilde\nu(X^\ell,\sP,\sP')$.

We will define an open neighborhood, $\tilde\sO(X^\ell,\sP,\sP',g_{\sP})$, of the zero section in $\tilde\nu(X^\ell,\sP,\sP')$
and maps, $\rho^{X,u}_{\sP,\sP'}$ and $\rho^{X,d}_{\sP,\sP'}$, so that the diagram,
\begin{equation}
\label{eq:OverlapEquation}
\begin{CD}
\tilde\sO(X^\ell,\sP,\sP',g_{\sP})
@> \rho^{X,u}_{\sP,\sP'} >>
\tilde\sO(X^\ell,g_{\sP'})
\\
@V \rho^{X,d}_{\sP,\sP'} VV
@V e(X^\ell,g_{\sP'}) VV
\\
\tilde\sO(X^\ell,g_{\sP})
@> e(X^\ell,g_{\sP}) >>
X^\ell
\end{CD}
\end{equation}
commutes and controls the overlaps of the images of the exponential maps, $e(X^\ell,g_{\sP})$ and $e(X^\ell,g_{\sP'})$.

\subsection{The downwards overlap map}
The downwards overlap map,
\begin{equation}
\label{eq:DefineXlDownwardsOverlap}
\rho^{X,d}_{\sP,\sP'}:\tilde\sO(X^\ell,\sP,\sP',g_{\sP})
\to
\tilde\sO(X^\ell,g_{\sP})
\end{equation}
is given by the fiberwise inclusions defined by the map in equation \eqref{eq:NormalOfR4Diagonal},
$$
e(Z_P,\sP'_P):
\Delta^\circ(Z_P,\sP'_P)\times\prod_{P'\in\sP'_P} Z_{P'}
\to
Z_P.
$$
Note that if $P\in\sP\cap\sP'$, then $\Delta^\circ(Z_P,\sP'_P)$ is a point and $e(Z_P,\sP'_P)$ is the identity map on $Z_P$. Because $e(Z_P,\sP'_P)$  is $\SO(4)$-equivariant, the product,
$$
\prod_{P\in\sP}
e(Z_P,\sP'_P):
\prod_{P\in\sP}
\left(\Delta^\circ(Z_P,\sP'_P)\times\prod_{P'\in\sP'_P} Z_{P'} \right)
\to
\prod_{P\in\sP} Z_P
$$
is $G(T,\sP)$-equivariant and thus extends to the domain and range given in \eqref{eq:DefineXlDownwardsOverlap}. The open subspace, $\tilde\sO(X^\ell,\sP,\sP',g_{\sP})$, must satisfy
\begin{equation}
\label{eq:XlOverlapCondition1}
\tilde\sO(X^\ell,\sP,\sP',g_{\sP})
\subseteqq
(\rho^{X,d}_{\sP,\sP'})^{-1}
\left(
\tilde\sO(X^\ell,g_{\sP})
\right)
\end{equation}
for the composition in the diagram \eqref{eq:OverlapEquation} to be defined.

\subsection{The upwards overlap map}
\label{subsubsec:XUpwardsOverlap}
We define the upwards overlap map,
\begin{equation}
\label{eq:DefineXlUpwardsOverlap}
\rho^{X,u}_{\sP,\sP'}:\tilde\nu(X^\ell,\sP,\sP')
\to
\tilde\sO(X^\ell,g_{\sP'}),
\end{equation}
as follows. Let
\begin{equation}
\label{eq:ProjectionToOverlapNormal}
\pi_1: \tilde\nu(X^\ell,\sP,\sP')
\to
\tilde\nu(X^\ell,\sP\to\sP')
\end{equation}
be the projection which is obvious from the definitions \eqref{eq:NormalData} and \eqref{eq:OverlapNormalBundle} (delete the factors of $Z_{P'}$).  We assume that the open set $\tilde\sO(X^\ell,\sP,\sP',g_{\sP})$ satisfies
\begin{equation}
\label{eq:XlOverlapCondition2}
\tilde\sO(X^\ell,\sP,\sP',g_{\sP})
\subseteqq
\pi_1^{-1}\left(
\tilde\sO(X^\ell,g_{\sP})
\right).
\end{equation}
Then the composition $e(X^\ell,g_{\sP})\circ \pi_1$ is defined on $\tilde\sO(X^\ell,\sP,\sP',g_{\sP})$ and takes values in $\Delta^\circ(X^\ell,\sP')$.

We now define the map $\rho^{X,u}_{\sP,\sP'}$. We make the assumption on the families of metrics $g_{\sP}$ and $g_{\sP'}$ appearing in
Lemma \ref{lem:RelateTXFiberBundles}. Given a point in $\tilde\nu(X^\ell,\sP,\sP')$ lying over $\by=(y_P)\in \Delta^\circ(X^\ell,\sP)$, define a point $\by'=(y_{P'})_{P'\in\sP'}\in\Delta^\circ(X^\ell,\sP')$ by the composition,
$$
e(X^\ell,g_{\sP})\circ \pi_1:
\tilde\nu(X^\ell,\sP,\sP')
\to
\Delta^\circ(X^\ell,\sP').
$$
Then, for $P'\in\sP'_P$ (where $P\in\sP$), define $F_{P'}$ to be the parallel translation of the frame $F_P\in \Fr(TX)|_{y_P}$ from $y_P$ to $y'_{P'}$ using the Levi--Civita connection defined by the metric $g_{\sP,\by}$. Leave the data in $Z_{P'}$ alone.  The result is an element of
$$
\Fr(TX^\ell,\sP',g_{\sP'})|_{\by'}\times_{G(T,\sP')} \prod_{P'\in\sP'} Z_{P'}
=
\tilde\nu(X^\ell,\sP')|_{\by'},
$$
thus defining the map $\rho^{X,u}_{\sP,\sP'}$. Observe that if $P\in\sP\cap\sP'$, the parallel translation for the $P$-th component will be trivial because $\Delta^\circ(Z_P,\sP'_P)$ is $|P|$-copies of the zero vector.

\subsection{Commuting overlap maps}
\label{subsec:Commuting_overlap_map}
We now define conditions on the metrics $g_{\sP}$ and $g_{\sP'}$ which ensure that the diagram \eqref{eq:OverlapEquation} commutes.
The diagram will commute if we chose the families of metrics $g_{\sP}$ and $g_{\sP'}$ in such a way as to eliminate holonomy. To do this, we introduce the notion of a \emph{locally flattened metric}.

\begin{defn}
\label{defn:Locally_flattened_metric}
If $A\subset X^\ell$, then the \emph{support} of $A$ in $X$ is the union of the images of $A$ under the projection maps,
$X^\ell\to X$. Let $g_{\sP}$ be a smooth family of metrics on $X$ parameterized by $\Delta^\circ(X^\ell,\sP)$. Let $\tilde\sU'\subseteqq\tilde\sU(X^\ell,g_{\sP})$ be a neighborhood of $\Delta^\circ(X^\ell,\sP)$. The family of metrics $g_{\sP}$ is \emph{locally flat with respect to $\tilde\sU'$ at} $\bx\in\Delta^\circ(X^\ell,\sP)$ if the metric $g_{\sP,\bx}$ is flat on the support of $\tilde\sU'\cap\pi(X^\ell,g_{\sP})^{-1}(\bx)$ in $X$. The family of metrics is \emph{locally flat with respect to} $\tilde\sU'$ if this holds for all $\bx\in \Delta^\circ(X^\ell,\sP)$.
\end{defn}

The following lemma shows that the diagram \eqref{eq:OverlapEquation} commutes when the metrics are locally flat in the sense of Definition \ref{defn:Locally_flattened_metric}.

\begin{lem}
\label{lem:CommutingExponentialMaps}
Let $\sP<\sP'$ be partitions of $N_\ell$. Assume that the smooth families of metrics $g_{\sP}$ and $g_{\sP'}$ satisfy
\begin{enumerate}
\label{item:CommutingExponentialMaps1}
\item
The metrics $g_{\sP}$ are flat with respect to a neighborhood $\tilde\sU'\subseteqq\tilde\sU(X^\ell,g_{\sP})$,
\item
\label{item:CommutingExponentialMaps2}
For $\by'\in \tilde\sU(X^\ell,g_{\sP})\cap\Delta^\circ(X^\ell,\sP')$ with $\pi(X^\ell,g_{\sP})(\by')=\by$, the metrics satisfy
$g_{\sP,\by}=g_{\sP',\by'}$.
\end{enumerate}
Then there is a neighborhood $\tilde\sO'(X,\sP,\sP',g_{\sP})$ of the zero section in the bundle $\tilde\nu(X^\ell,\sP,\sP')$ defined in
\eqref{eq:OverlapNormalBundle} such that the diagram \eqref{eq:OverlapEquation} commutes when restricted to $\tilde\sO'(X,\sP,\sP',g_{\sP})$.
\end{lem}

\begin{proof}
We observe that for the compositions
$$
e(X^\ell,g_{\sP})\circ\rho^{X,d}_{\sP,\sP'}
\quad\text{and}\quad
e(X^\ell,g_{\sP'})\circ\rho^{X,u}_{\sP,\sP'}
$$
to be defined, the neighborhood $\tilde\sO'(X,\sP,\sP',g_{\sP})$ must satisfy the inclusion relations \eqref{eq:XlOverlapCondition1} and \eqref{eq:XlOverlapCondition2}. However, we replace the condition \eqref{eq:XlOverlapCondition1} with the stronger requirement that
\begin{equation}
\label{eq:XOverlapAssumption}
\tilde\sO'(X,\sP,\sP',g_{\sP})
\subseteqq
(e(X^\ell,g_{\sP})\circ\rho^{X,d}_{\sP,\sP'})^{-1}( \tilde\sU').
\end{equation}
On any open subset satisfying the above condition, the composition
$$
e(X^\ell,g_{\sP})\circ\rho^{X,d}_{\sP,\sP'}
$$
is defined by the exponential map for the metric $g_{\sP,\by}$ of the vectors obtained by appropriately adding  vectors in the fiber of $\tilde\nu(X^\ell,\sP,\sP',g_{\sP'})$. The composition,
$$
e(X^\ell,g_{\sP'})\circ\rho^{X,u}_{\sP,\sP'},
$$
is defined by first parallel translating the vectors defined by the elements of $Z_{P'}$ from $\by$ to the nearby points in the support of $\by'$, and then exponentiating.
By \eqref{eq:XOverlapAssumption} and the first assumption on the metrics $g_{\sP}$ and $g_{\sP'}$, the metric $g_{\sP,\by}$ is flat where this parallel translation takes place and hence the two compositions are equal.
\end{proof}

\subsection{The projection maps}
We now show that the commutativity of the diagram \eqref{eq:OverlapEquation} will imply that the first Thom--Mather property holds.

\begin{lem}
\label{lem:DiagonalsTM1}
Continue the hypotheses and notation of Lemma \ref{lem:CommutingExponentialMaps}. Then, upon restriction to $\tilde\sU'\cap \pi(X^\ell,g_{\sP'})^{-1}\left(\Delta^\circ(X^\ell,\sP')\right)$, the following equality holds:
\begin{equation}
\label{eq:DiagonalsTM1}
\pi(X^\ell,g_{\sP})\circ\pi(X^\ell,g_{\sP'})=\pi(X^\ell,g_{\sP}).
\end{equation}
\end{lem}

\begin{proof}
The subspace $\tilde\sU'\cap \pi(X^\ell,g_{\sP'})^{-1}\left(\Delta^\circ(X^\ell,\sP')\right)$ is in the image of the composition
$e(X^\ell,g_{\sP})\circ\rho^{X,d}_{\sP,\sP'}$. The equality follows from the observation that the pullback of $\pi(X^\ell,g_{\sP'})$ by $\rho^{X,d}_{\sP,\sP'}$ is given by the projection $\pi_1$ of \eqref{eq:ProjectionToOverlapNormal}.
\end{proof}

\section{Construction of the families of locally flattened metrics}
\label{sec:TMMapsOnSymX}
We now construct the families of metrics satisfying the conditions of Lemma \ref{lem:CommutingExponentialMaps} and thus yielding projection maps $\pi(X^\ell,g_{\sP})$ which satisfy \eqref{eq:DiagonalsTM1}.  The first step is to construct a locally flat family of metrics.

\begin{lem}
\label{lem:FlatteningMetrics}
If $g_{\sP}$ is a family of metrics smoothly parameterized by $\Delta^\circ(X^\ell,\sP)$, then there is a family of smooth metrics
$\bar g_{\sP}$ smoothly parameterized by $\Delta^\circ(X^\ell,\sP)$ such that the following hold:
\begin{enumerate}
\item
\label{item:FlatteningMetrics_1}
$\bar g_{\sP}$ is locally flat with respect to a tubular neighborhood
$\tilde\sU'\subseteqq\tilde\sU(X^\ell,g_{\sP})$.
\item
\label{item:FlatteningMetrics_2}
For every $\bx\in \Delta^\circ(X^\ell,\sP)$, the metric $\bar g_{\sP, \bx}$ is equal to
$g_{\sP,\bx}$ on the complement of the support of
$\tilde\sU(X^\ell,g_{\sP})\cap\pi(X^\ell,g_{\sP})^{-1}(\bx)$.
\item
\label{item:FlatteningMetrics_3}
The families $g_{\sP}$ and $\bar g_{\sP}$ are $C^1$ close within the support of $\tilde\sU'$.
\item
\label{item:FlatteningMetrics_4}
The exponential maps $e(X^\ell,g_{\sP})$ and $e(X^\ell,\bar g_{\sP})$ are equal.
\item
\label{item:FlatteningMetrics_5}
If $g_{\sP}$ is already flat with respect to $\tilde\sU(X^\ell,g_{\sP})$ at $\bx\in \Delta^\circ(X^\ell,\sP)$,
then $\bar g_{\sP, \bx}=g_{\sP,\bx}$.
\end{enumerate}
\end{lem}

\begin{proof}
Let $\tilde\sU(X^\ell,g_{\sP})$ be the tubular neighborhood of $\Delta^\circ(X^\ell,\sP)$ defined by the family of metrics, $g_{\sP}$, so
\[
e(X^\ell,g_{\sP}):\tilde\sO(X^\ell,g_{\sP})\subset\tilde\nu(X^\ell,\sP)\cong\tilde\sU(X^\ell,g_{\sP}).
\]
For $P\in\sP$, let $s_P:\Delta^\circ(X^\ell,\sP)\to[0,\infty)$ be a smooth function such that
\begin{equation}
\label{eq:DefineSepFunctions}
\{(v_1,\dots,v_\ell)\in\tilde\nu(X^\ell,\sP)|_{\bx}: |v_i|\le 4s_P(\bx)
\quad\text{for all $i\in P$ and $P\in\sP$}\}
\subset
\sO(\sP).
\end{equation}
For $\bx\in\Delta^\circ(\sP)$ and $(F_P)_{P\in\sP}\in \Fr(TX^\ell,\sP,g_{\sP})|_{\bx}$, the frame $F_P$ and
the exponential map of the metric $g_{\sP,\bx}$ identify $B(0,4 s_P(\bx))\subset \RR^4$ with
$B_{g_{\sP,\bx}}(x_P,4s_P(\bx)) \subset X$. Let $\delta_{x_P}$ denote the flat metric on $B_{g_{\sP,\bx}}(x_P,4s_P(\bx))$ given by the pushforward of the Euclidean metric on $B(0,4 s_P(\bx))$ by the exponential map.
For $\la\in (0,1]$, let $\chi_\la:\RR\to \RR$ be a smooth function such that $\chi_\la(t)=0$ for $t\le \half \la$ and $\chi_\la(t)=1$ for $t\ge \la$. We define $\bar g_{\sP, \bx}$ by
\begin{equation}
\label{eq:FlatFamily1}
\bar g_{\sP, \bx} :=
  \begin{cases}
    g_{\sP,\bx} & \text{on $X\less \cup_i B_{g_{\sP,\bx}}(x_P,4s_P(\bx))$},
    \\
    \chi_{4s_P(\bx)}(|x|)g_{\sP,\bx} + (1-\chi_{4s_P(\bx)}(|x|))\delta_{x_P}
      &\text{on $\Omega_{g_{\sP,\bx}}(x_P,2s_P(\bx),4s_P(\bx)))$},
    \\
    \delta_{x_P} & \text{on $B_{g_{\sP,\bx}}(x_P,2s_P(\bx))$},
  \end{cases}
\end{equation}
where $x(\cdot)$ denotes the local coordinate chart defined by the exponential map for $g_{\sP,\bx}$.

The resulting family of metrics, $\bar g_\sP$, is locally flat with respect to the tubular neighborhood $\tilde\sU'$ defined by replacing $4s_P(\bx)$ with $2s_P(\bx)$ in \eqref{eq:DefineSepFunctions}.
Item \eqref{item:FlatteningMetrics_2} follows immediately from the construction
\eqref{eq:FlatFamily1} while Item \eqref{item:FlatteningMetrics_3} follows from the expression
for a metric in geodesic normal coordinates
\cite[Definition 1.24, Proposition 1.25, and Corollary 1.32]{Aubin}.

For $x$ in the support of $\bx$, the radial geodesics from $x$ are the same for $g_{\sP,\bx}$ and for $\bar g_{\sP,\bx}$.
These radial geodesics from $x$ define the exponential maps $e(X^\ell,g_{\sP})$ and $e(X^\ell,\bar g_{\sP})$ and
hence these maps are equal.

Because the metric $g_{\sP,\bx}$ has not been changed at $x_P$, the frame bundles do not change.  If the metric $g_{\sP,\bx}$ is flat on $B(x_P,4s_P(\bx))$, then it is given by the Euclidean metric in geodesic normal coordinates by \cite[Theorem 7.11]{Spivak2}.  Thus, $\delta_{x_P}=g_{\sP,\bx}$ and the interpolation in \eqref{eq:FlatFamily1} does not change the metric.
\end{proof}

We now construct the desired families of metrics.

\begin{lem}
\label{lem:FlatteningMetrics2}
For each partition $\sP$ of $N_\ell$, there is a smooth family of metrics $g_{\sP}$ on $X$ defining tubular neighborhoods $\tilde\sU(X^\ell,g_{\sP})$ by their exponential maps, such that the following hold:
\begin{enumerate}
\item
\label{item:FlatteningMetrics2_1}
The family $g_{\sP}$ is locally flat with respect to a tubular neighborhood
$\tilde\sU'(X^\ell,g_{\sP})\subset \tilde\sU(X^\ell,g_{\sP})$.
\item
\label{item:FlatteningMetrics2_2}
If $\sP<\sP'$, then for all $\bx'\in\Delta^\circ(X^\ell,\sP')\cap\tilde\sU'(X^\ell,g_{\sP})$, we have
$g_{\sP',\bx'}=g_{\sP,\bx}$, where $\bx=\pi(X^\ell,g_{\sP})(\bx')$.
\item
\label{item:FlatteningMetrics2_3}
For $\si\in \fS_\ell$, we have $g_{\sP,\bx}=g_{\si(\sP),\si(\bx)}$.
\end{enumerate}
\end{lem}

\begin{proof}
The proof is by induction on the partitions, $\sP$. We will construct the family of metrics, $g_{\sP}$, for
one partition in the orbit, $[\sP]$, and then use Item \eqref{item:FlatteningMetrics2_3} in the lemma to define $g_{\sP}$ for all other $\sP'\in [\sP]$.

Choose a Riemannian metric, $g$, on $X$.  To every partition $\sP$ of $N_\ell$, define an initial family of metrics $g_{\sP}$ by setting $g_{\sP,\bx}:=g$.

Let $\sP_0$ be the crudest partition. Redefine $g_{\sP_0}$ by applying the flattening procedure of Lemma \ref{lem:FlatteningMetrics}. The redefined family of metrics $g_{\sP_0}$ is locally flat with respect to the tubular neighborhood $\tilde\sU'(X^\ell,\sP_0,g_{\sP_0})$ given in Lemma \ref{lem:FlatteningMetrics}.

Assume that families of metrics, $g_{\sP}$, and tubular neighborhoods, $\tilde\sU'(X^\ell,g_{\sP})$, satisfying the
conclusions of the lemma have been defined for all $\sP<\sP'$.  There are two families of metrics
on $X$ parameterized by the intersection,
\begin{equation}
\Delta^\circ(X^\ell,\sP')\cap \sU(X^\ell,g_{\sP}),
\end{equation}
namely, the initial family, $g_{\sP'}$, and the family, $\pi(X^\ell,g_{\sP})^*g_{\sP}$, which is defined by
$$
\pi(X^\ell,g_{\sP})^*g_{\sP,\bx'}:=g_{\sP,\bx},
\quad\text{where}\quad
\pi(X^\ell,g_{\sP})(\bx')=\bx.
$$
For $\sP_1<\sP'$ and $\sP_2<\sP'$ and by shrinking the tubular neighborhoods when necessary, observe that if
the intersection,
\begin{equation}
\label{eq:TripleOverlap}
\Delta^\circ(X^\ell,\sP')\cap \tilde\sU(X^\ell,g_{\sP_1})
\cap \tilde\sU(X^\ell,g_{\sP_2}),
\end{equation}
is non-empty, then there is a relation $\sP_1<\sP_2<\sP'$. Consequently, the relations,
$$
\pi(X^\ell,g_{\sP_1})\circ \pi(X^\ell,g_{\sP_2})
=
\pi(X^\ell,g_{\sP_1}),
$$
(which follow by induction and Lemma \ref{lem:DiagonalsTM1}) and the identity,
$$
g_{\sP_2}=\pi(X^\ell,g_{\sP_1})^*g_{\sP_1},
$$
imply that the families of metrics,
\begin{equation}
\label{eq:PullbackMetrics}
\pi(X^\ell,g_{\sP_1})^*g_{\sP_1}
\quad
\text{and}
\quad
\pi(X^\ell,g_{\sP_2})^*g_{\sP_2},
\end{equation}
are equal on \eqref{eq:TripleOverlap}.  Thus we can form a new family of metrics by interpolating between $g_{\sP'}$ on
$$
\Delta^\circ(X,\sP') \setminus \bigcup_{\sP<\sP'}\tilde\sU'(X^\ell,g_{\sP}),
$$
and the pullback metrics \eqref{eq:PullbackMetrics}. By further shrinking the tubular neighborhoods,
$\tilde\sU'(X^\ell,g_{\sP})$, for $\sP<\sP'$ and taking a sufficiently small neighborhood,
$\tilde\sU'(X^\ell,g_{\sP'})$, we can assume this interpolated family is locally flat with respect
to $\tilde\sU'(X^\ell,g_{\sP'})$ on
\begin{equation}
\label{eq:NghOfLowerDiagonals}
\Delta^\circ(X^\ell,\sP')\cap\left(\bigcup_{\sP<\sP'}\tilde\sU(X^\ell,g_{\sP})\right).
\end{equation}
Now apply the flattening procedure of Lemma \ref{lem:FlatteningMetrics} to this interpolated family.  Because this family is already flat with respect to $\tilde\sU'(X^\ell,g_{\sP'})$ for $\bx$ in the open space \eqref{eq:NghOfLowerDiagonals}, the flattening procedure does not change the interpolated family there.  Hence, the new flattened family satisfies Items \eqref{item:FlatteningMetrics2_1} and \eqref{item:FlatteningMetrics2_2}, completing the induction.
\end{proof}

\section{Normal bundles of strata of $\Sym^\ell(X)$}
\label{sec:NormalBundlesOfSymlXStrata}
We now use the $\fS_\ell$-equivariance of the construction of the previous section to define normal bundles of the strata $\Si(X^\ell,\sP)$ in $\Sym^\ell(X)$.  Recall from Lemma \ref{lem:EndOfUpperStratum} that to describe the end of $\Si(X^\ell,\sP')$ near $\Si(X^\ell,\sP)$, we have to consider the ends of the diagonals $\Delta^\circ(X^\ell,\sP'')$ near $\Delta^\circ(X^\ell,\sP)$ for all $\sP''\in [\sP<\sP']$.

The analogue of the bundle \eqref{eq:OverlapNormalBundle} appropriate for this discussion is
\begin{equation}
\label{eq:OverlapNormalBundleSymm}
\begin{aligned}
{}&
\tilde\nu(X^\ell,\sP,[\sP'],g_{\sP})
\\
{}&\quad:=
\Fr(TX^\ell,\sP,g_{\sP})
\times_{\tilde G(T,\sP)}
\bigsqcup_{\sP''\in [\sP<\sP']}
\prod_{P\in\sP}
\left(
\Delta^\circ(Z_P,\sP''_P)
\times
\prod_{Q\in\sP''_P} Z_{Q}
\right),
\end{aligned}
\end{equation}
where we use $\sqcup$ to denote disjoint union.
The symmetric group $\fS(\sP)$ acts on the spaces,
$$
\Fr(TX^\ell,\sP,g_{\sP}),
\quad
\bigsqcup_{\sP''\in [\sP<\sP']}
\prod_{P\in\sP}\Delta^\circ(Z_P,\sP''_P),
\quad\text{and}\quad
\bigsqcup_{\sP''\in [\sP<\sP']} \prod_{P\in\sP} \prod_{Q\in\sP''_P} Z_{Q},
$$
by permuting the factors.  The diagonal action of $\fS(\sP)$ on the preceding three spaces defines an action on $\tilde\nu(X^\ell,\sP,[\sP'],g_{\sP})$. We define
\begin{equation}
\label{eq:OverlapNormalBundleSymm1}
\begin{aligned}
{}&
\nu(X^\ell,\sP,[\sP'],g_{\sP})
\\
{}&\quad :=
\Fr(TX^\ell,\sP,g_{\sP})
\times_{G(T,\sP)}
\bigsqcup_{\sP''\in [\sP<\sP']}
\prod_{P\in\sP}
\left(
\Delta^\circ(Z_P,\sP''_P)
\times
\prod_{Q\in\sP''_P} Z_{Q}
\right)
\end{aligned}
\end{equation}
where, if $\tilde G(T,\sP)$ as defined as in \eqref{eq:TgBundleDiagonalStructureGrp},
$$
G(T,\sP) := \tilde G(T,\sP)\rtimes \fS(\sP),
$$
where $\rtimes$ denotes the semi-direct product, \cite[p. 59]{LangAlgebra} (compare the structure group
appearing on \cite[p. 252]{FrM}).
Thus,
$$
\nu(X^\ell,\sP,[\sP'],g_{\sP})
=
\tilde \nu(X^\ell,\sP,[\sP'],g_{\sP})/\fS(\sP).
$$
We define
\begin{equation}
\label{eq:TubNghBundleOverlaps}
\begin{aligned}
{}&
\tilde\sO(X^\ell,\sP,[\sP'],g_{\sP})
\subseteqq
\tilde \nu(X^\ell,\sP,[\sP'],g_{\sP}),
\\
{}&
\sO(X^\ell,\sP,[\sP'],g_{\sP})
=
\tilde\sO(X^\ell,\sP,[\sP'],g_{\sP})/\fS
\subseteqq
\nu(X^\ell,\sP,[\sP'],g_{\sP})
\end{aligned}
\end{equation}
by analogy with $\tilde\sO(X^\ell,\sP,\sP',g_{\sP})$ in \eqref{eq:XlOverlapCondition1}.

We define a downwards overlap map,
\begin{equation}
\label{eq:DefineXellDownwardsOverlapMap}
\rho^d_{\sP,[\sP']}:
\nu(X^\ell,\sP,[\sP'],g_{\sP})
\to
\sO(X^\ell,g_{\sP})=
\tilde
\sO(X^\ell,g_{\sP})/\fS(\sP),
\end{equation}
by extending the map of fibers given by
$$
\coprod_{\sP''\in [\sP<\sP']}
\prod_{P\in\sP}\
e(Z_P,\sP''_P),
$$
over the domain.
We have used the convention that if $f_\alpha:X_\alpha\to Y_\alpha$ is a family of continuous maps, then
$\coprod_\alpha f_\alpha:\sqcup_\alpha X_\alpha\to\sqcup_\alpha Y_\alpha$ is the map defined on
the disjoint unions by this family.
Note that the quotient in \eqref{eq:DefineXellDownwardsOverlapMap} could be taken to be by $W(\sP)$ instead of $\fS(\sP)$, as $\Ga(\sP)$ acts trivially.  This phenomenon will occur frequently, but as it will not affect the discussion we will not mention it again.

To define an upwards overlap map, we define
\begin{equation}
\label{eq:XlTubNghOverlapUpwardsSym}
\sO(X^\ell,[\sP<\sP'],g_{\sP''})
:=
\bigsqcup_{\sP''\in [\sP<\sP']}
\tilde\sO(X^\ell,\sP'',g_{\sP''})/\fS(\sP).
\end{equation}
The upwards transition map is then given by
\begin{equation}
\label{eq:XlUpwardsOverlapMapSym}
\rho^u_{\sP,[\sP']}:
\sO(X^\ell,\sP,[\sP'],g_{\sP})
\to
\sO(X^\ell,[\sP<\sP',g_{\sP''}),
\quad
\rho^u_{\sP,[\sP']}=\coprod_{\sP''\in [\sP<\sP']} \rho^u_{\sP,[\sP']}.
\end{equation}
The equivariance of the exponential map $e(X^\ell,g_{\sP})$ with respect to the symmetric group action implies that $e(X^\ell,g_{\sP})$ defines an embedding,
\begin{equation}
\label{eq:SymExponentialMap}
e(X^\ell,g_{\sP}):
\sO(X^\ell,g_{\sP})
\to
\Sym^\ell(X).
\end{equation}
This equivariance and
Item \eqref{item:FlatteningMetrics2_3} of Lemma \ref{lem:FlatteningMetrics2} then imply that the exponential maps, $e(X^\ell,g_{\sP''})$, define an embedding,
$$
e(X^\ell,g_{[\sP<\sP']})=
\coprod_{\sP''\in [\sP<\sP']} e(X^\ell,g_{\sP'}):
\sO(X^\ell,[\sP<\sP',g_{\sP''})
\to
\Sym^\ell(X).
$$
We then have the

\begin{prop}
\label{prop:SymmetricProductOverlapDiagram}
For each partition $\sP$ of $N_\ell$, let $g_{\sP}$ be the family of metrics constructed in Lemma \ref{lem:FlatteningMetrics2}. Then the image of the map \eqref{eq:SymExponentialMap} defines a homeomorphism onto a neighborhood $\sU(X^\ell,g_{\sP})$ of $\Si(X^\ell,\sP)$. For partitions $\sP<\sP'$ of $N_\ell$, there is a commutative diagram,
\begin{equation}
\label{eq:SymOverlapEquation}
\begin{CD}
\sO(X^\ell,\sP,[\sP'],g_{\sP})
@> \rho^{X,u}_{\sP,[\sP']} >>
\sO(X^\ell,[\sP<\sP'],g_{\sP''})
\\
@V \rho^{X,d}_{\sP,[\sP']} VV
@V e(X^\ell,g_{[\sP<\sP']}) VV
\\
\sO(X^\ell,g_{\sP})
@> e(X^\ell,g_{\sP'}) >>
\Sym^\ell(X)
\end{CD}
\end{equation}
such that any point in $\sU(X^\ell,g_{\sP})\cap\sU(X^\ell,g_{\sP'})$ is in the image of the compositions
$e(X^\ell,g_{\sP})\circ \rho^d_{\sP,[\sP']}$ and $e(X^\ell,g_{[\sP<\sP']})\circ \rho^u_{\sP,[\sP']}$.
\end{prop}

\begin{proof}
The proposition follows immediately from the $\fS_\ell$-equivariance
of the constructions of Sections \ref{sec:OverlapMaps} and \ref{sec:TMMapsOnSymX}.
\end{proof}

We note that the projection maps, $\pi(X^\ell,g_{\sP}):\tilde \sU(X^\ell,g_{\sP})
\to\Delta^\circ(X^\ell,\sP)$, are $\fS(\sP)$-equivariant and thus define projection
maps,
$$
\pi(X^\ell,g_{\sP}):
\sU(X^\ell,g_{\sP})
\to
\Si(X^\ell,\sP),
$$
on the quotient. These maps on the quotient still satisfy \eqref{eq:DiagonalsTM1}.

\section{The tubular distance function}
\label{sec:SymlXTubDistance}
Rather than a single tubular distance function, $t_{\sP}:\sU(X^\ell,g_{\sP})\to [0,\infty)$,
we find it convenient to define several functions, one for each $P\in\sP$, whose common zero locus is $\Si(X^\ell,\sP)$.
The functions, $t(Z_P)$, defined in Lemma \ref{lem:R4TMScale} are $\SO(4)$-equivariant and thus define functions,
$$
t(Z_P):
\Fr(TX^\ell,\sP,g_{\sP})\times_{\tilde G(T,\sP)}\prod_{P\in\sP}Z_P
\to [0,\infty).
$$
For $I^{\sP}=\prod_{P\in\sP} [0,1]$, we define
\begin{equation}
\label{eq:DefineDiagonalTubularDistFunction}
\vec t(X^\ell,g_{\sP}) :=\left(t(Z_P)\circ e(X^\ell,g_{\sP})^{-1}\right)_{P\in\sP}
:\tilde\sU(X^\ell,g_{\sP}) \to I^{\sP}.
\end{equation}
If $\sP$ is the finest partition described in Remark \ref{rmk:XlTopStratum}, then the function, $\vec t(X^\ell,g_{\sP})$, is the zero map.

To define hypersurfaces and disk bundles in $\sU(X^\ell,g_{\sP})$, we will use the following `square' sets,
\begin{equation}
\label{eq:DefineSquares}
\begin{aligned}
D(\sP,\eps)
{}&=
\{(t_P)_{P\in\sP}: 0\le t_P< \eps\ \text{for all $P\in\sP$}\},
\\
\bar D(\sP,\eps)
{}&=
\{(t_P)_{P\in\sP}: 0\le t_P\le \eps\ \text{for all $P\in\sP$}\},
\\
\rd \bar D(\sP,\eps)
{}&=
\{(t_P)_{P\in\sP}\in \bar D(\sP,\eps): t_P=\eps\ \text{for some $P\in\sP$}\},
\end{aligned}
\end{equation}
to define tubular neighborhoods of the diagonals. The symmetric group, $\fS(\sP)$, acts on $I^{\sP}$ and thus on the sets \eqref{eq:DefineSquares} by permuting the factors. Observe that $\vec t(X^\ell,g_{\sP})$ is equivariant
with respect to the action of $\fS(\sP)$ and thus defines a map,
$$
\vec t(X^\ell,g_{\sP}): \sU(X^\ell,g_{\sP})
\to
I^{\sP}/\fS(\sP).
$$
For $s\in [0,1]$, we will write $s\vec t(X^\ell,g_{\sP})$ for the function $(s t_P)_{P\in\sP}$ and $\vec t^{\,P}(X^\ell,g_{\sP})$ for the $P$-th component of $\vec t(X^\ell,g_{\sP})$. The identification of $\sU(X^\ell,g_{\sP})$ with an open subspace of the bundle $\nu(X^\ell,\sP)$ then gives the following result.

\begin{lem}
\label{lem:XDefRetraction}
For every partition $\sP$ of $N_\ell$, there is a smooth map,
$$
r(X^\ell,g_{\sP}):
\tilde\sU''(X^\ell,g_{\sP})
\times
[0,1]
\to
\tilde\sU''(X^\ell,g_{\sP}),
$$
such that $r(X^\ell,g_{\sP})(\cdot,s)$ is a homeomorphism for all $s\in (0,1]$ and $r(X^\ell,g_{\sP})$ satisfies
the following properties:
\begin{enumerate}
\item
\label{item:XDefRetraction1}
$r(X^\ell,g_{\sP})(\cdot,1)$ is the identity map.
\item
\label{item:XDefRetraction2}
$r(X^\ell,g_{\sP})(\cdot,0)=\pi(X^\ell,g_{\sP})(\cdot)$.
\item
\label{item:XDefRetraction3}
For all $s\in [0,1]$,
\begin{equation}
\label{eq:XDefScaleCondition}
\vec t(X^\ell,g_{\sP})\circ r(X^\ell,g_{\sP})(\cdot,s)
=
s^2 \vec t(X^\ell,g_{\sP})(\cdot).
\end{equation}
\end{enumerate}
If $\sP<\sP'$, then on $\tilde\sU''(X^\ell,g_{\sP})\cap \tilde\sU''(X^\ell,g_{\sP'})$ we have
\begin{equation}
\label{eq:XOverlapRetractionCondition}
\pi(X^\ell,g_{\sP})\circ r(X^\ell,g_{\sP'})(\cdot,s)
=
\pi(X^\ell,g_{\sP})(\cdot),
\end{equation}
for all $s\in [0,1]$.
\end{lem}

\begin{proof}
If $\sP$ is the finest partition, as described in Remark \ref{rmk:XlTopStratum}, set $r(X^\ell,g_{\sP})(\cdot,s)$
equal to the identity map for all $s$. Otherwise, we define the map $r(X^\ell,g_{\sP})$ by pushing forward
a deformation retraction, $r_\nu(X^\ell,\sP)$, on $\nu(X^\ell,\sP)$ to $\sU(X^\ell,\sP)$ by the exponential map, $e(X^\ell,g_{\sP})$. Let $r_P: Z_P\times [0,1] \to Z_P$ be the deformation retraction defined in \eqref{eq:ConeRetractionMap}.  If $\Delta_{\sP}:[0,1]\to \prod_{P\in\sP}[0,1]$ is the diagonal inclusion
and $\id_{Z,P}:Z_P\to Z_P$ is the identity map, then the composition
$$
\begin{CD}
(\prod_{P\in\sP} Z_P)\times [0,1]
\\
@V \prod_{P\in\sP}\id_{Z,P}\times \Delta_{\sP} VV
\\
\prod_{P\in\sP} (Z_P\times [0,1])
\\
@V \prod_{P\in\sP} r_P VV
\\
\prod_{P\in\sP} Z_P
\end{CD}
$$
defines a $G(T,\sP)$-equivariant deformation retraction and hence a deformation retraction,
$$
r_\nu(X^\ell,\sP):\nu(X^\ell,\sP) \times [0,1]\to \nu(X^\ell,\sP).
$$
If $r(X^\ell,g_{\sP})$ is defined by pushing $r_\nu(X^\ell,\sP)$ forward by $e(X^\ell,g_{\sP})$, then $r(X^\ell,g_{\sP})$ immediately satisfies Items \eqref{item:XDefRetraction1} and \eqref{item:XDefRetraction2}
of the lemma.  The property \eqref{eq:XDefScaleCondition} follows from the definitions of $r(X^\ell,g_{\sP})$  and $\vec t(X^\ell,g_{\sP})$ and
\eqref{eq:RetractionAndScale}. Finally, \eqref{eq:XOverlapRetractionCondition} follows immediately from
the definition of $r(X^\ell,g_{\sP})$ as a map of fibers of $\nu(X^\ell,\sP)$.
\end{proof}

The following result describes the failure of the second Thom--Mather conditions in \eqref{eq:TMProperties} for these maps.

\begin{lem}
\label{lem:SymProjTubDistRelations1}
Let $\sP<\sP'$ be partitions of $N_\ell$ and
let $r(X^\ell,g_{\sP'})$ be the map defined in Lemma \ref{lem:XDefRetraction}.
If $\bx\in\sU(X^\ell,g_{\sP})\cap \sU(X^\ell,g_{\sP'})$ and
$s\in [0,1]$, then the following hold:
\begin{enumerate}
\item
\label{item:SymProjTubDistRelations1_1}
If $P\notin\sP\cap \sP'$, then $\vec t^{\,P}(X^\ell,g_{\sP})(\bx) =
\vec t^{\,P}(X^\ell,g_{\sP})(r(X^\ell,g_{\sP'})(\bx,s))$.
\item
\label{item:SymProjTubDistRelations1_2}
If $P\in\sP\cap\sP'$, then $s^2\vec t^{\,P}(X^\ell,g_{\sP})(\bx) =
\vec t^{\,P}(X^\ell,g_{\sP})(r(X^\ell,g_{\sP'})(\bx,s))$.
\item
\label{item:SymProjTubDistRelations1_3}
If $P\in\sP\cap\sP'$, then $\vec t^{\,P}(X^\ell,g_{\sP})(\bx)=\vec t^{\,P}(X^\ell,g_{\sP'})(\bx)$.
\end{enumerate}
\end{lem}

\begin{proof}
We prove Item \eqref{item:SymProjTubDistRelations1_1} by showing that if $P\notin\sP\cap\sP'$, then the functions $t(Z_P)$ are constant on the fibers of $\pi(X^\ell,g_{\sP'})$ which contain the paths $s\mapsto r(X^\ell,g_{\sP'})(\cdot,s)$. The definition of the upwards overlap map $\rho^u_{\sP,[\sP']}$ implies that pulling the map $\pi(X^\ell,g_{\sP'})$ back by $\rho^u_{\sP,[\sP']}$ gives the projection map,
$$
\bigsqcup_{\sP''\in [\sP<\sP']}
\prod_{P\in\sP}
\Delta^\circ(Z_P,\sP''_P)
\times
\prod_{Q\in\sP''_P} Z_{Q}
\to
\bigsqcup_{\sP''\in [\sP<\sP']}
\prod_{P\in\sP} \Delta^\circ(Z_P,\sP''_P).
$$
From the definition of the downwards overlap map, $\rho^d_{\sP,[\sP']}$, in \eqref{eq:DefineXellDownwardsOverlapMap},
we see that the image of the restriction of $\rho^d_{\sP,[\sP']}$ to the above fiber is given by
$$
\bigsqcup_{\sP''\in [\sP<\sP']}
\prod_{P\in\sP} \Imag(e(Z_P,\sP''_P))
=
\bigsqcup_{\sP''\in [\sP<\sP']}
\prod_{P\in\sP} \sU(Z_P,\sP''_P)
\subset
\prod_{P\in\sP} Z_P,
$$
where the set $\sU(Z_P,\sP'_P)$ is defined prior to \eqref{eq:R4DiagonalProjection}. We focus on the restriction of $\pi(X^\ell,g_{\sP'})\circ \rho^u_{\sP,[\sP']}$ to one of the components in the domain above and see that the result holds for each such restriction. Pushing this restriction of  $\pi(X^\ell,g_{\sP'})\circ \rho^u_{\sP,[\sP']}$
forward by $\rho^d_{\sP,[\sP']}$ to $\prod_{P\in\sP} \sU(Z_P,\sP'_P)$ gives the map,
$$
\prod_{P\in\sP} \pi(Z_P,\sP'_P):
\prod_{P\in\sP}   \sU(Z_P,\sP'_P) \to \prod_{P\in\sP}  Z_P,
$$
where $\pi(Z_P,\sP'_P)$ is defined in \eqref{eq:DefineR4DiagonalProjection}. The property \eqref{eq:R4ZPTM2} then implies that $t(Z_P)$ is constant on the fibers of $\pi(Z_P,\sP'_P)$ for all $P\notin \sP\cap\sP'$.

If $P\in \sP\cap \sP'$, the equality \eqref{eq:RetractionAndScale}, namely, $t(Z_P)(r_P(\cdot,s))=s^2t(Z_P)(\cdot)$, then yields Item
\eqref{item:SymProjTubDistRelations1_2}.

Item \eqref{item:SymProjTubDistRelations1_3}
follows by observing that if $P\in\sP\cap\sP'$, then both $\rho^u_{\sP,[\sP']}$  and $\rho^d_{\sP,[\sP']}$ are defined by the identity map on $Z_P$ and the $P$-th components of
$\vec t(X^\ell,g_{\sP})\circ \rho^d_{\sP,[\sP']}$ and $\vec t(X^\ell,g_{\sP'})\circ \rho^u_{\sP,[\sP']}$ are both defined by $t(Z_P)$ and hence are equal.
\end{proof}

Although the preceding lemma showed that the second Thom--Mather identity in \eqref{eq:TMProperties} does not hold, the following lemma shows that for $\sP<\sP'$, the restriction of the projection $\pi(X^\ell,g_{\sP'})$ to a sufficiently small tubular neighborhood of $\Delta^\circ(X^\ell,\sP')$ preserves the tubular neighborhood defined by $\vec t(X^\ell,g_{\sP})$.

\begin{lem}
\label{lem:SymProjTubDistRelations}
Let $\sP<\sP'$ be partitions of $N_\ell$. Let $\hat D(\sP,\eps)$ denote either $D(\sP,\eps)$ or $\bar D(\sP,\eps)$.
If $\bx\in\sU(X^\ell,g_{\sP})\cap \sU(X^\ell,g_{\sP'})$, then the following hold:
\begin{enumerate}
\item
\label{item:SymProjTubDistRelations_1}
If $\vec t(X^\ell,g_{\sP})(\bx)\in \hat D(\sP,\eps)$, then
$(\vec t(X^\ell,g_{\sP}) \circ \pi(X^\ell,g_{\sP'}))(\bx)\in \hat D(\sP,\eps)$.
\item
\label{item:SymProjTubDistRelations_2}
If $\eps'<\eps$, and $\vec t(X^\ell,g_{\sP'})(\bx)\in \bar D(\sP',\eps')$,
and $(\vec t(X^\ell,g_{\sP}) \circ \pi(X^\ell,g_{\sP'}))(\bx)\in \hat D(\sP,\eps)$,
then  $\vec t(X^\ell,g_{\sP})(\bx)\in \hat D(\sP,\eps)$.
\end{enumerate}
\end{lem}

\begin{proof}
Item \eqref{item:SymProjTubDistRelations_1} follows immediately from
Items \eqref{item:SymProjTubDistRelations1_1} and \eqref{item:SymProjTubDistRelations1_2} of Lemma \ref{lem:SymProjTubDistRelations1}.
To prove Item \eqref{item:SymProjTubDistRelations_2}, we note that if $(\vec t(X^\ell,g_{\sP}) \circ \pi(X^\ell,g_{\sP'}))(\bx)\in \hat D(\sP,\eps)$ and $\vec t(X^\ell,g_{\sP})(\bx)\notin \hat D(\sP,\eps)$, then Item \eqref{item:SymProjTubDistRelations1_1} of Lemma \ref{lem:SymProjTubDistRelations1} implies that there must be a $P\in\sP\cap\sP'$ such that the $P$-th component of $\vec t(X^\ell,g_{\sP})(\bx)$ is greater than or equal to $\eps$.  However, Item \eqref{item:SymProjTubDistRelations1_3} of Lemma \ref{lem:SymProjTubDistRelations1} and the assumption $\vec t(X^\ell,g_{\sP'})(\bx)\in \bar D(\sP',\eps')$ imply that the $P$-th component of $\vec t(X^\ell,g_{\sP})(\bx)$ is less than or equal to $\eps'$, which is strictly less than $\eps$. This contradiction proves Item \eqref{item:SymProjTubDistRelations_2}.
\end{proof}

\section{Decomposition of the strata}
\label{sec:SymlXStrataDecomp}
We now construct decompositions of $\Sym^\ell(X)$ and $\Si(X^\ell,\sP)$ to be used in the construction of the link
of $M_{\fs}\times\Sym^\ell(X)$ appearing in equation \eqref{eq:RawCobordismSum}.

\begin{lem}
\label{lem:DecomposeSymProduct}
Enumerate the strata of $\Sym^\ell(X)$ as described in
Section \ref{subsec:EnumStrata} using partitions $\sP_0,\dots,\sP_n$.
To each partition $\sP_i$, there is a small, positive
parameter $\eps_i$ satisfying $\eps_i>\eps_j$ for $i<j$
such that if we define
\begin{align*}
T_i
{}&:=
\sU(X^\ell,\sP_i)
\cap
\vec t(X^\ell,g_{\sP_i})^{-1}(D(\sP_i,\eps_i))
\setminus
\bigcup_{j<i}\ \vec t(X^\ell,g_{\sP_j})^{-1}(D(\sP_j,\eps_j)),
\quad\text{for $i\neq n$},
\\
T_n
{}&:=
\Si(X^\ell,\sP_n)
\setminus
\bigcup_{j<n}\ \vec t(X^\ell,g_{\sP_j})^{-1}(D(\sP_j,\eps_j)),
\end{align*}
then for $i\neq n$, we have $T_i\Subset \sU(X^\ell,g_{\sP_i})$, and $\Sym^\ell(X)=\cup_i T_i$, and this union is disjoint.
\end{lem}

\begin{proof}
The lemma follows immediately from the definitions.
\end{proof}

The decomposition of $\Sym^\ell(X)$ in the previous lemma leads to decompositions of the strata $\Si(X^\ell,\sP_i)$.

\begin{lem}
\label{lem:CompactSubsetsOfSi}
Continue the hypotheses of Lemma \ref{lem:DecomposeSymProduct} and let $\eps_i$ be the parameter constructed in
Lemma \ref{lem:DecomposeSymProduct}. For one of the partitions, $\sP_k$, in the given enumeration, define
$$
T_{k,j} := \cl\left(\Si(X^\ell,\sP_k)\right) \cap T_j,
\quad\text{and}\quad
K_k
:=
\cl\left(\Si(X^\ell,\sP_k)\right)
\setminus
\bigcup_{i<k}\ T_{k,i}.
$$
Then the following hold.
\begin{enumerate}
\item
\label{item:CompactSubsetsOfSi_1}
There is an equality,
$$
\cl\left(\Si(X^\ell,\sP_k)\right) = K_k \cup \left(\bigcup_{j<k} T_{k,j}\right),
$$
where the union on the right-hand side is disjoint.
\item
\label{item:CompactSubsetsOfSi_2}
$K_k$ is compact.
\item
\label{item:CompactSubsetsOfSi_3}
The restriction of $\pi(X^\ell,\sP_j)$ to $T_{k,j}$ takes values in $K_j$.
\end{enumerate}
\end{lem}

\begin{proof}
The first assertion follows immediately from Lemma \ref{lem:DecomposeSymProduct}. We note that $K_k$ is compact by observing that it is a closed subset of the compact set $\cl\, \Si(X^\ell,\sP_k)$.

We now prove Item \eqref{item:CompactSubsetsOfSi_3}. If $\bx\in T_{k,j}$, then
Item \eqref{item:SymProjTubDistRelations_1} in Lemma \ref{lem:SymProjTubDistRelations}
implies that $\bx'=\pi(X^\ell,\sP_k)(\bx)\in \vec t(X^\ell,\sP_j)^{-1}(D(\sP_k,\eps_j))$.
If $\pi(X^\ell,\sP_k)(\bx)\notin K_j$, then there is an index $i<j$ such that
$\vec t(X^\ell,\sP_i)(\bx')\in D(\sP_i,\eps_i)$.  However, $\bx\in T_{k,j}$ implies that $\vec t(X^\ell,\sP_i)(\bx')\notin D(\sP_i,\eps_i)$ contradicting
Item \eqref{item:SymProjTubDistRelations_2} in Lemma \ref{lem:SymProjTubDistRelations}.
\end{proof}

The spaces $T_i\subset \Sym^\ell(X)$ admit deformation retractions to $K_i\subset\Si(X^\ell,\sP_i)$ as described in the following
lemma.

\begin{lem}
\label{lem:SymmetricProductRetractions}
Continue the hypotheses and notation of Lemma \ref{lem:DecomposeSymProduct}.
If the parameters $\eps_0,\dots,\eps_n$ defining $T_i$
in Lemma \ref{lem:DecomposeSymProduct}, are sufficiently small then, for each partition $\sP_i\neq \sP_n$, there is a map,
$$
r_i:\Sym^\ell(X)\to\Sym^\ell(X),
$$
satisfying the following properties, for $\pi_i=\pi(X^\ell,g_{{\mathscr{P}}_i})$:
\begin{enumerate}
\item
\label{item:SymmetricProductRetractions_1}
$r_i$ is homotopic to the identity.
\item
\label{item:SymmetricProductRetractions_2}
The restriction of $r_i$ to $T_i$ is equal to $\pi(X^\ell,g_{\sP_i})$.
\item
\label{item:SymmetricProductRetractions_3}
For each partition $\sP_j$, we have $r_i(\Si(X^\ell,\sP_j)) \subset \cl\, \Si(X^\ell,\sP_j)$.
\item
\label{item:SymmetricProductRetractions_4}
For $j<i$, we have $\pi(X^\ell,g_{\sP_j})\circ r_i=\pi(X^\ell,g_{\sP_j})$.
\item
\label{item:SymmetricProductRetractions_5}
For $j<i$, the restriction of $r_i$ to the complement of
$\vec t(X^\ell,g_{\sP_i})^{-1}(D(\sP_i,\eps_j))$ is equal to the identity map.
\item
\label{item:SymmetricProductRetractions_6}
For $j<i$ and $A_j$ equal to $\vec t(X^\ell,g_{\sP_j})^{-1}(\bar D(\sP_j,\eps_j))$
or $\vec t(X^\ell,g_{\sP_j})^{-1}(D(\sP_j,\eps_j))$, we have
$r_i(\bx)\in A_j$ if and only if $x\in A_j$.
\item
\label{item:SymmetricProductRetractions_7}
For $j<i$, we have $r_i(T_j)\subset T_j$.
\item
\label{item:SymmetricProductRetractions_8}
$r_i (K_k\cup \cup_{j=k-1}^{i+1} T_{k,j}) = K_k\cup \cup_{j=k-1}^{i} T_{k,j}$.
\item
\label{item:SymmetricProductRetractions_9}
The restriction of $r_i$ to the complement of $\cup_{j\le i}\vec t(X^\ell,g_{\sP_j})^{-1}(\bar D(\sP_j,\eps_j))$
is injective.
\end{enumerate}
\end{lem}

\begin{proof}
If the parameters $\eps_i$ are sufficiently small, then we can find neighborhoods,
$$
\sO_1(X^\ell,g_{\sP_i})\sqsubset\sO_2(X^\ell,g_{\sP_i})
\sqsubset
\sO(X^\ell,g_{\sP_i}),
$$
of $\Si(X^\ell,\sP_i)$ such that
$$
T_i=
e(X^\ell,g_{\sP_i})(\sO_1(X^\ell,g_{\sP_i}))
 - \bigcup_{j<i} \vec t(X^\ell,g_{\sP_j})^{-1} D(\sP_j,\eps_j),
$$
and for all $j<i$,
\begin{equation}
\label{eq:OuterOpenSet}
e(X^\ell,g_{\sP_i})(\sO_2(X^\ell,g_{\sP_i}))
\sqsubset
\vec t(X^\ell,g_{\sP_i})^{-1} D(\sP_i,\eps_j).
\end{equation}
There is a continuous  function, $p:\sO(X^\ell,g_{\sP_i})\to [0,1]$, with $p^{-1}(0)=\sO_1(X^\ell,g_{\sP_i})$ and $p^{-1}(1)=\sO(X^\ell,g_{\sP_i})-\sO_2(X^\ell,g_{\sP_i})$.
Setting $p_e:=p\circ e(X^\ell,g_{\sP_i})^{-1}$,
define $r_i$ on $\sU(X^\ell,g_{\sP_i})$ by
$$
r_i(\cdot) := r(X^\ell,g_{\sP_i})(\cdot, p_e(\cdot)),
$$
and extend it as the identity map on $\Sym^\ell(X)-\sU(X^\ell,g_{\sP_i})$. Items
\eqref{item:SymmetricProductRetractions_1} and \eqref{item:SymmetricProductRetractions_2}
follow immediately from the definition.  Item  \eqref{item:SymmetricProductRetractions_3}
follows from noting that the maps $r_P$ and thus $r(X^\ell,g_{\sP})$ preserve the diagonals.  Item \eqref{item:SymmetricProductRetractions_4} follows from \eqref{eq:XOverlapRetractionCondition}. Item \eqref{item:SymmetricProductRetractions_5} follows from the construction of $r_i$ and the inclusion \eqref{eq:OuterOpenSet}.

We now prove Item \eqref{item:SymmetricProductRetractions_6}. If $\bx\in A_j$, then
Items \eqref{item:SymProjTubDistRelations1_1} and \eqref{item:SymProjTubDistRelations1_2}  of Lemma \ref{lem:SymProjTubDistRelations1} imply that
$r(X^\ell,g_{\sP_i})(\bx,s)\in A_j$ for all $s\in [0,1]$ and thus $r_i(\bx)\in A_j$. If $r_i(\bx)\in A_j$ and $\bx\notin A_j$, then $\vec t(X^\ell,g_{\sP_j})(\bx)\neq  \vec t(X^\ell,g_{\sP_j})(r_i(\bx))$.
Item \eqref{item:SymProjTubDistRelations1_1} of Lemma \ref{lem:SymProjTubDistRelations1} implies that
$$
\vec t^{\,P}(X^\ell,g_{\sP_j})(\bx)=\vec t^{\,P}(X^\ell,g_{\sP_j})(r_i(\bx))\le\eps_j,
$$
for all $P\in\sP_i\cap\sP_j$.  Thus $\bx\notin A_j$ implies that there is a $P\in\sP_i\cap\sP_j$ with
$\vec t^{\,P}(X^\ell,g_{\sP_j})(\bx)\ge \eps_j$.
Item \eqref{item:SymProjTubDistRelations1_3} in Lemma \ref{lem:SymProjTubDistRelations1} then implies that
$\vec t^{\,P}(X^\ell,g_{\sP_i})(\bx)\ge \eps_j$, so $\bx\notin \vec t(X^\ell,g_{\sP_i})^{-1}( D(\sP_i,\eps_j))$
and thus, by Item \eqref{item:SymmetricProductRetractions_5}, $r_i(\bx)=\bx$, a contradiction.

Item \eqref{item:SymmetricProductRetractions_7} follows from the definition of $T_i$ and Item \eqref{item:SymmetricProductRetractions_6}. Item \eqref{item:SymmetricProductRetractions_8} follows from the definition of $T_{k,j}$ and Item \eqref{item:SymmetricProductRetractions_6}. Item
\eqref{item:SymmetricProductRetractions_9} follows from the definition of $r_i$ and Item \eqref{item:SymmetricProductRetractions_6}.
\end{proof}

In the following lemma, we show how to use the maps $r_i$ to construct a global map to $\Sym^\ell(X)$ from a collection of maps to $\Sym^\ell(X)$ which are approximately equal.

\begin{lem}
\label{lem:GlobalProjToSymCrit}
Continue the hypotheses of Lemma \ref{lem:SymmetricProductRetractions}. Let $\{U_i\}_{i=0}^n$ be an open cover of a topological space $Y$ with $U_i\cap U_j=\emptyset$ unless $i\le j$.  For $i=0,\dots,n$,
let $f_i:U_i\to \Sym^\ell(X)$ be a continuous map satisfying the following properties:
\begin{enumerate}
\item
\label{item:GlobalProjToSymCritAssumption_1}
If $j<k$, then $f_k(U_k\cap U_j) \subset \cup_{i\le j} T_i$,
\item
\label{item:GlobalProjToSymCri2Assumption_2}
If $i<j$, then $\pi_i\circ f_j|_{U_i\cap U_j}=f_i|_{U_i\cap U_j}$.
\end{enumerate}
Then there is a map $F:Y\to\Sym^\ell(X)$ such that $F|_{U_i}$ is homotopic to $f_i$.
\end{lem}

\begin{proof}
Define $m_i:=r_0\circ r_1\circ\dots \circ r_i$ and let $m_{-1}$ be the identity map. Each map $m_i$ is homotopic to the identity by
Item \eqref{item:SymmetricProductRetractions_1} in Lemma \ref{lem:SymmetricProductRetractions}. If $g_i=m_{i-1}\circ f_i$, so $g_0=f_0$, then $g_i$ is homotopic to $f_i$.  We claim that the maps $g_i$ and $g_j$ are equal on $U_i\cap U_j$ and thus define the global map $F$.

Arguing by upwards induction, assume that $g_i|_{U_i\cap U_j}=g_j|_{U_i\cap U_j}$ for all $i,j<k$. We now prove that if $j<k$, then $g_k|_{U_j\cap U_k}=g_j|_{U_j\cap U_k}$. If $x\in U_j\cap U_k$, then by property \eqref{item:GlobalProjToSymCritAssumption_1} there is an index $i\le j$ such that $f_k(x)\in T_i$.

Item \eqref{item:SymmetricProductRetractions_7} in Lemma \ref{lem:SymmetricProductRetractions} and the fact that $f_k(x)\in T_i$ imply that
\begin{equation}
\label{eq:RetractionInclusion}
(r_{i+1}\circ \dots \circ r_{k-1}\circ f_k)(x)\in T_i.
\end{equation}
Hence,
\begin{align*}
g_k(x)
{}&=
(m_{i-1}\circ r_i\circ r_{i+1}\circ \dots \circ r_{k-1}\circ f_k)(x))
\\
{}&=
(m_{i-1}\circ \pi_i\circ r_{i+1}\circ \dots \circ r_{k-1}\circ f_k)(x)
\quad\text{(by Item \eqref{item:SymmetricProductRetractions_2} in Lemma \ref{lem:SymmetricProductRetractions} and
\eqref{eq:RetractionInclusion})}
\\
{}&=
(m_{i-1}\circ \pi_i\circ f_k)(x)
\quad\text{(by Item \eqref{item:SymmetricProductRetractions_4} in Lemma \ref{lem:SymmetricProductRetractions})}
\\
{}&
=(m_{i-1}\circ f_i)(x)
\quad\text{(by property \eqref{item:GlobalProjToSymCri2Assumption_2})}
\\
{}&=g_i(x).
\end{align*}
By induction, we have $g_i(x)=g_j(x)$ which, together with the previous equality $g_k(x)=g_i(x)$, gives the
equality $g_j(x)=g_k(x)$ required to complete the proof.
\end{proof}

\chapter{The instanton moduli space with spliced ends}
\label{chap:SplicedEnd}
The gluing maps defined in \cite{FL3} provide a tubular
neighborhood structure for each stratum $M_{\fs}\times\Si$
of $M_{\fs}\times\Sym^\ell(X)$ just as the exponential map
$e(X^\ell,g_{\sP})$ did for the stratum $\Si(X^\ell,\sP)$.
In this analogy, $\bar{M}^{s,\natural}_{\ka}(S^4,\delta)$ plays the same role as the space $Z_P$,
making up the fiber of the tubular neighborhood.
Here, $\bar{M}^{s,\natural}_{\ka}(S^4,\delta)$  denotes
the Uhlenbeck compactification of the moduli
space of framed anti-self-dual connections on $S^4$ which are
mass-centered and have scale less than $\delta$, where the mass and scale are defined in
\eqref{eq:CenterOfMass} and \eqref{eq:ConnectionScale} respectively.
However, it will require much more work to establish the commutativity
of the diagram analogous to \eqref{eq:OverlapEquation} for the
tubular neighborhood structure of $M_{\fs}\times\Sym^\ell(X)$.
As we described in our Introduction, we will work with splicing maps
instead of gluing maps to establish such an equality because
the explicit definition of the splicing map allows us to verify equalities
analogous to \eqref{eq:OverlapEquation}.
In order to use splicing maps instead of gluing maps,
we must introduce a deformation of $\bar{M}^{s,\natural}_{\ka}(S^4,\delta)$,
which we call the \emph{instanton moduli space with spliced ends}.
The new space can still be used to define the domain of the gluing map of Hypothesis \ref{hyp:Gluing}
because the connections in this new space satisfy the same crucial
estimates (for example, Item \eqref{item:ExistenceOfSplicedEndsModuli4} in
Theorem \ref{thm:ExistenceOfSplicedEndsModuli}) as the connections in $\bar{M}^{s,\natural}_{\ka}(S^4,\delta)$.

We will be working with topological spaces which are \emph{smoothly stratified} as in the definition given in
Section \ref{sec:Introduction_diagTM}.
A continuous map between such spaces is smoothly stratified if it maps each smooth stratum into a smooth
stratum and if the restriction of the map to each smooth stratum of the domain is a smooth map.

\section{Introduction}
\label{sec:SplicedEnd_introduction}
The instanton moduli space with spliced ends  will be identical
to $\barM^{s,\natural}_{\ka}(S^4,\delta)$
except on an Uhlenbeck neighborhood of the \emph{punctured, centered symmetric product},
\begin{equation}
\label{eq:PuncturedTrivialStrata}
[\Theta]\times \left(\Sym^{\ka,\natural}_\delta(\RR^4)\less c_{\ka}\right),
\end{equation}
where $\Theta$ is the product connection on $S^4\times\SU(2)$, and
$\Sym^{\ka,\natural}_\delta(\RR^4)$ is the subspace of mass-centered points in
the $\ka$-th symmetric product of $\RR^4$ with a scale constraint (see \eqref{eq:SingularPointScale}),
and $c_\ka=[0,\dots,0]\in \Sym^{\ka}(\RR^4)$ is given by $\ka$ copies of the origin.

A neighborhood of the strata \eqref{eq:PuncturedTrivialStrata} in
$\barM^{s,\natural}_\ka(S^4,\delta)$ is parameterized by the union of the images
of the gluing maps $\bga_{\Theta,\Si}$, as $\Si$ varies over the strata
in \eqref{eq:PuncturedTrivialStrata}.
As described in Section \ref{subsec:R4Diagonals}, the space $\Sym^{\ka,\natural}_\delta(\RR^4)$
is the quotient of a subspace $Z_\ka(\delta)\subset (\RR^4)^\ka$ by the symmetric
group action.
If $\Si$ is a stratum of $\Sym^{\ka,\natural}_\delta(\RR^4)$
and $\sP$ is a partition of $N_\ka$ corresponding to $\Si$, so
$\Si=\Delta^\circ((\RR^4)^\ka,\sP)\cap Z_\ka(\delta)/\fS(\sP)$ as in \eqref{eq:DiagonalsOfZP},
then the gluing map is a smoothly stratified embedding \cite{Feehan_Leness_monopolegluingbook},
\begin{equation}
\label{eq:ProductStratumGluingOnS4}
\bga_{\Theta,\sP}:
\left(\Delta^\circ(\RR^{4\ka},\sP)\cap Z_\ka(\delta)\right)\times_{\fS(\sP)}
\prod_{P\in\sP} \bar M^{s,\natural}_{|P|}(S^4,\delta_P)
\to \barM^{s,\natural}_{\ka}(S^4,\delta).
\end{equation}
The instanton moduli space with spliced ends, $\barM^{s,\natural}_{\spl,\ka}(S^4,\delta)$, will be defined by replacing the image of $\bga_{\Theta,\sP}$ with the image of the splicing map $\bga'_{\Theta,\sP}$ on the interior of this tubular neighborhood.
The connections defining this deformation of the instanton moduli space will be almost anti-self-dual as
measured by the norm $\|\cdot \|_{L^{\sharp,2}(X)}$, which on a four-manifold $X$ was defined
in \cite[Equation (6.1a)]{TauFrame} and \cite[Equation (4.3)]{FeehanSlice} by
\begin{equation}
\label{eq:LSharpNorm}
\| a\|_{L^{\sharp,2}(X)}
:=
\|a\|_{L^2(X)}
+
\sup_{x\in X}\|\dist^{-2}(x,\cdot)|a|\|_{L^1(X)}.
\end{equation}
We construct $\barM^{s,\natural}_{\spl,\ka}(S^4,\delta)$ and summarize its properties in the following

\begin{thm}
\label{thm:ExistenceOfSplicedEndsModuli}
Let $\bar\sB^s_\ka(S^4,\delta)$ be the Uhlenbeck extension of the
framed quotient space defined in
\eqref{eq:ExtendedFramedQuotientSpace}.
There exists a
smoothly-stratified subspace,
$\barM^{s,\natural}_{\spl,\ka}(S^4,\delta)\subset\bar\sB^s_\ka(S^4,2\delta)$,
closed under the
$\SO(3)\times\SO(4)$ action defined in \eqref{eq:ExtendedFrameAction} and \eqref{eq:ExtendedRotationAction}, satisfying the following properties:
\begin{enumerate}
\item
\label{item:ExistenceOfSplicedEndsModuli1}
There is a smoothly-stratified, $\SO(3)\times\SO(4)$-equivariant
homeomorphism,
\begin{equation}
\label{eq:HomeomorphismFromSplEndsToInstModuli}
\barM^{s,\natural}_{\spl,\ka}(S^4,\delta)\cong \barM^{s,\natural}_{\ka}(S^4,\delta),
\end{equation}
which is the identity on the levels,
$$
[\Theta]\times \Sym^{\ka,\natural}_\delta(\RR^4),
$$
where $\Sym^{\ka,\natural}_\delta(\RR^4)\subset\Sym^\ka(\RR^4)$ is
defined in \eqref{eq:SingularPointScale}.
\item
\label{item:ExistenceOfSplicedEndsModuli2}
There is an Uhlenbeck neighborhood in the sense of Definition
\ref{defn:UhlenbeckConvergence}, namely $W_{\ka}$, in $\bar{\sB}^s_{\ka}(S^4,2\delta)$ of the punctured,
centered symmetric product \eqref{eq:PuncturedTrivialStrata}
such that
$$
\barM^{s,\natural}_{\spl,\ka}(S^4,\delta)\less W_{\ka}
 =
\barM^{s,\natural}_{\ka}(S^4,\delta) \less W_{\ka}.
$$
\item
\label{item:ExistenceOfSplicedEndsModuli3}
For each stratum $[\Theta]\times\Si$ in \eqref{eq:PuncturedTrivialStrata}, there is
an Uhlenbeck neighborhood, $W(\Si)\subset\bar{\sB}^s_{\ka}(S^4,\delta)$, of the
stratum
\begin{equation}
\label{eq:TrivialStratumP}
[\Theta]\times \Si
\end{equation}
and an open neighborhood, $\sO^{\asd}_1(\Theta,\sP,\delta)\subset
\nu(\Theta,\sP,\delta)$,
of the stratum \eqref{eq:DefineASDSplicingDomain}
in the domain of the splicing map $\bga'_{\Theta,\sP}$
defined in \eqref{eq:R4TrivialSplicing},
such that
\begin{equation}
\label{eq:EqualsSplicingImage}
\barM^{s,\natural}_{\spl,\ka}(S^4,\delta)\cap
W(\Si) = \bga'_{\Theta,\sP}\left( \sO^{\asd}_1(\Theta,\sP,\delta)\right).
\end{equation}
\item
\label{item:ExistenceOfSplicedEndsModuli4}
For each integer $\kappa \geq 0$, there is a constant, $C = C(\ka)$, such that for all $[A,\bx]\in \barM^{s,\natural}_{\spl,\ka}(S^4,\delta)$, one has
\[
\|F^+_{A}\|_{L^{\sharp,2}(S^4)}\le C\delta,
\]
where the norm $\|\cdot\|_{L^{\sharp,2}}$ is defined in \eqref{eq:LSharpNorm}.
\end{enumerate}
\end{thm}

We construct $\barM^{s,\natural}_{\spl,\ka}(S^4,\delta)$ by
induction on $\ka \in \NN$.
We first construct what we call the `spliced end' of
$\barM^{s,\natural}_{\spl,\ka}(S^4,\delta)$.
The \emph{spliced end} will be an Uhlenbeck neighborhood of the punctured, centered symmetric product
\eqref{eq:PuncturedTrivialStrata} in $\barM^{s,\natural}_{\spl,\ka}(S^4,\delta)$.
We shall define the spliced end of $\barM^{s,\natural}_{\spl,\ka}(S^4,\delta)$
to be the union of the images of splicing maps \eqref{eq:R4TrivialSplicing},
$\bga'_{\Theta,\sP}$, where
$\sP$ is a partition of $N_\ka$ corresponding to a stratum $\Si$
in $\Sym^{\ka,\natural}_\delta(\RR^4)-\{c_{\ka}\}$
and where the domain of the splicing map is not defined by
the moduli spaces, $\barM^{s,\natural}_{\ka_i}(S^4,\delta_i)$, with smaller instanton number
but rather by the instanton moduli spaces with spliced ends,
$\barM^{s,\natural}_{\spl,\ka_i}(S^4,\delta_i)$, with smaller instanton number.
A technical result on the equality of two compositions of splicing maps,
Proposition \ref{prop:CommutingSplicingSymmQuotient},
shows that the images of the different splicing maps intersect in open sets,
proving that the union of these images, $W_\ka$, is a smoothly stratified space.

After reviewing the properties of connections on $S^4$ in
Section \ref{subsec:SpaceOfConn}, we adapt the notation of
Chapter \ref{chap:Diagonals} to describe the
product-connection  strata, \eqref{eq:PuncturedTrivialStrata}, in Section \ref{subsec:R4Diagonals}.
In Section \ref{subsec:SplicingMapsR4}, we define the relevant splicing maps.
In Section \ref{subsec:OverlapSpliceMapsR4}, we prove the crucial
technical result, Proposition \ref{prop:CommutingSplicingSymmQuotient},
that provides control over the overlap of two splicing maps.
In Section \ref{subsec:SplicedEnd}, we finish the first stage of
the proof of Theorem \ref{thm:ExistenceOfSplicedEndsModuli}
by constructing $W_\ka$ in Proposition \ref{prop:ExistenceOfSplicedEnd}.

The second stage of the proof of Theorem \ref{thm:ExistenceOfSplicedEndsModuli} appears in
Section \ref{subsec:Collar}, where we
construct an isotopy of
the complement of a neighborhood of the punctured, centered symmetric
product in $W_\ka$ to the actual moduli space,
$\barM^{s,\natural}_{\ka}(S^4,\delta)$.
This isotopy is defined by the composition of the isotopy
defined by the gluing
map and the isotopy defined by the centering map.  That such an isotopy exists
follows immediately from the estimates on
$\|F^+_{A}\|_{L^{\sharp,2}(S^4)}$ for a connection, $A$, in the spliced end; most of the
work in Section \ref{subsec:Collar} lies in constructing the parameter
on the spliced end for the isotopy.
Finally, in Section \ref{subsec:PropsOfSplicedEnd} we construct a
cone parameter on the
instanton moduli space with spliced end
and prove that
this space has the structure
of a
\emph{Whitney stratified space}\label{Instanton_moduli_space_spliced_end_Whitney_stratified_space} in the sense of \cite[Section 1.2]{GorMacPh}.

\section{Connections over the four-dimensional sphere}
\label{subsec:SpaceOfConn}
We begin by reviewing some standard definitions concerning
connections over $S^4$; similar definitions appear in
\cite{FLKM1,FL3}.

For an integer $\kappa \geq 1$, let $\sB_\ka(S^4)$ denote the quotient, $\sA_\kappa/\sG_\kappa$, of the affine space, $\sA_\kappa$, of $L_2^2$ $\mathrm{SU}(2)$ connections on a Hermitian smooth vector bundle $E_\ka\to S^4$ with $c_2(E_\ka)=\ka$, modulo the action of the group, $\sG_\kappa = \Aut(E_\kappa)$, of $L^2_3$ $\mathrm{SU}(2)$  gauge transformations of $E_\kappa$.
Let $\sB^s_\ka(S^4)\to\sB_\ka(S^4)$ be the
principal
$\SO(3)$-bundle defined by
\[
\sB^s_\ka(S^4) := \left(\sA_\ka(S^4) \times P_\kappa|_s\right)/\sG_\kappa,
\]
where $P_\kappa = \Fr(E_\kappa)$, the principal $\SU(2)$-bundle of $\mathrm{SU}(2)$  frames for $E_\kappa$ and $s \in S^4$ is the South Pole, identified with the point at infinity in $S^4 = \RR^4\cup\{\infty\}$.
Let
\begin{equation}
\label{eq:IdealConnS4}
\bar\sB_{\ka}(S^4)
=
\bigsqcup_{\ell=0}^\ka \ \left( \sB_{\ka-\ell}(S^4)\times \Sym^\ell(\RR^4)\right),
\end{equation}
denote the space of ideal connections on
$E_\kappa$
given the topology
induced by Uhlenbeck convergence as defined in Definition \ref{defn:UhlenbeckConvergence}.
We use $\Sym^\ell(\RR^4)$ in \eqref{eq:IdealConnS4} rather than $\Sym^\ell(S^4)$ to simplify the definitions
\eqref{eq:CenterOfMass} and \eqref{eq:ConnectionScale} below.  We will be working with ideal connections with
finite scale so using $\Sym^\ell(\RR^4)$ does not omit relevant ideal connections.
We write elements of
$\bar\sB_{\ka}(S^4)$
as $[A,\bx]$, where $[A]\in \sB_{\ka-\ell}(S^4)$ and $\bx\in \Sym^\ell(\RR^4)$.

If $n\in S^4$ denotes the North Pole, identified with the origin in $S^4 = \RR^4\cup\{\infty\}$,
let $y(\cdot) = \varphi_n^{-1}: S^4\less\{s\}\to\RR^4$ be the coordinate chart given by a stereographic
projection from the South Pole with $\varphi_n(0) = n \in S^4$.
For a point
$[A,\bx]\in \bar\sB_{\ka}(S^4)$ with $\bx=[x_1,\dots,x_\ell]$,
we define its \emph{center of mass}
by
\begin{equation}
\label{eq:CenterOfMass}
z[A,\bx]
:=
\left(\int_{\RR^4} |\varphi_n^*F_A|^2 d^4y \right)^{-1} \int_{\RR^4} y |\varphi_n^*F_A|^2 d^4y
+ \sum_{i=1}^\ell x_i \in \RR^4,
\end{equation}
and the {\em scale\/} by
\begin{equation}
\label{eq:ConnectionScale}
\la[A,\bx]^2
:=
\left(\int_{\RR^4} |\varphi_n^*F_A|^2 d^4y \right)^{-1}\int_{\RR^4} | y - z[A, \bx]|^2 |\varphi_n^*F_A|^2 d^4y
+
\sum_{i=1}^\ell |x_i|^2 \in [0,\infty).
\end{equation}
For a constant $\eps > 0$, we define
$$
\bar\sB_{\ka}(S^4,\eps) := \la^{-1}([0,\eps]).
$$
We continue to denote by $z[\cdot]$ and $\la[\cdot]$ the pull-back of those functions from $\sB_{\ka}(S^4)$ to $\sB^s_{\ka}(S^4)$.
If $[A,\bx]\in \bar\sB_{\ka}(S^4,\eps)$ and $\eps\ll 1$, then the support of $\bx$ is disjoint from the south pole.
Hence, the quotient space, $\sB^s_{\ka}(\eps)$, is continuously embedded in the quotient space of pairs of \emph{ideal} connections and frames,
\begin{equation}
\label{eq:ExtendedFramedQuotientSpace}
\begin{aligned}
\bar\sB^{s}_{\ka}(S^4,\eps)
&:=
\la^{-1}[0,\eps]
\cap
\left(
\bigsqcup_{\ell=0}^\ka \ \left( \sB^s_{\ka-\ell}(S^4)\times \Sym^\ell(\RR^4)\right)\right),
\\
\bar\sB^{s,\natural}_{\ka}(S^4,\eps)
&=
z^{-1}(0)\cap \bar\sB^{s}_{\ka}(S^4,\eps),
\end{aligned}
\end{equation}
and $\sB^{s,\natural}_{\ka}(S^4,\eps)$ is continuously embedded in the quotient space of pairs of ideal \emph{mass-centered} connections and frames. If $\Theta$ is the
product connection, we call
$[\Theta,F^s,c_\ka]\in \sB^s_0(S^4)\times\Sym^{\ka}(\RR^4)$
the \emph{cone point} of $\bar\sB^{s,\natural}_{\ka}(S^4,\eps)$.

The group $\SU(2)$ acts on $\bar\sB^{s,\natural}_{\ka}(S^4,\eps)$ by
the action of $\SU(2)$ on the frame in $\Fr(E_\ka)|_s$.
(Because connections that are mass-centered at the North Pole and have uniformly positive scale cannot bubble over the South Pole
--- see \cite{FeehanGeometry} for an explanation based on the Chebychev Inequality --- the bundle fibers $E_{\ka-\ell}$ for $0\leq \ell \leq \ka$
can all be identified with $E_\ka|_s$.)
The stabilizer of any point in $\bar\sB^{s,\natural}_{\ka}(S^4,\eps)$
contains
$\{\pm\id\}\subset\SU(2)$,
so this action of $\SU(2)$ on
$\bar\sB^{s,\natural}_{\ka}(S^4,\eps)$ descends to an action
of $\SO(3)=\SU(2)/\{\pm\id\}$ on $\bar\sB^{s,\natural}_{\ka}(S^4,\eps)$,
\begin{equation}
\label{eq:ExtendedFrameAction}
\SO(3)\times \bar\sB^{s,\natural}_{\ka}(S^4,\eps)
\to
\bar\sB^{s,\natural}_{\ka}(S^4,\eps)
\end{equation}
which we refer to as the
action on the frame.

The group $\SO(4)$ acts on $S^4$ by 
pull-back of its action on $\RR^4$ via stereographic projection, $S^4\less\{s\}\to\RR^4$.
The action of $\SO(4)$ on $S^4$ induces an action on
$\bar\sB^{s,\natural}_{\ka}(S^4,\eps)$,
\begin{equation}
\label{eq:ExtendedRotationAction}
\SO(4)\times \bar\sB^{s,\natural}_{\ka}(S^4,\eps)
\to
\bar\sB^{s,\natural}_{\ka}(S^4,\eps),
\end{equation}
 by pull-back of
the connection and frame for $E_\ka|_s$ and the rotation action on the points in
$\Sym^\ka(\RR^4)$ (see \cite[Section 4, p. 343]{TauFrame}
or \cite[Section 3.2]{FeehanGeometry}).

\section{Strata containing the product connection}
\label{subsec:R4Diagonals}
We now adapt
the framework of
Chapter \ref{chap:Diagonals} to
describe the product-connection strata  in
\eqref{eq:PuncturedTrivialStrata}.
Define $\Sym^{\ka,\natural}(\RR^4)$ to be the quotient
of the zero locus,
\[
Z_\kappa := z_{\ka}^{-1}(0),
\]
of the center-of-mass map, $z_\ka:\oplus_{i\in N_{\ka}}\RR^4\to\RR^4$
(defined following \eqref{eq:DefineZ}), by the symmetric group, $\fS_\ka$,
$$
\Sym^{\ka,\natural}(\RR^4) := Z_\kappa/\fS_\ka = z_\ka^{-1}(0)/\fS_\ka.
$$
We define a scale for elements of $\Sym^\ka(\RR^4)$ by
\begin{equation}
\label{eq:SingularPointScale}
\la[v_1,\dots,v_\ka]^2 := \sum_{i=1}^\kappa |v_i|^2,
\end{equation}
and set
\begin{align*}
Z_\ka(\delta) &:= z_{\ka}^{-1}(0)\cap \la^{-1}([0,\delta)),
\\
\Sym^{\ka,\natural}_\delta(\RR^4) &:= Z_\ka(\delta)/\fS_\ka.
\end{align*}
We will define strata of $\Sym^{\ka,\natural}_\delta(\RR^4)$
to be quotients of the intersections of the diagonals in $(\RR^4)^{\ka}$ with $Z_\ka(\delta)$.
For any partition $\sP$ of $N_\ka$, we define
\begin{equation}
\label{eq:DiagonalsOfZP}
\begin{aligned}
\Delta^\circ(Z_\ka(\delta),\sP)
&:=
\Delta^\circ(Z_\ka,\sP)\cap Z_{\ka}(\delta),
\\
\Si(Z_\ka(\delta),\sP)
&:=
\Delta^\circ(Z_\ka(\delta),\sP)/\fS(\sP)
\subset\Sym^{\ka,\natural}_\delta(\RR^4),
\end{aligned}
\end{equation}
where $\Delta^\circ(Z_\ka,\sP)$ is defined in
\eqref{eq:DefineCenteredR4Diagonal}.
Let $\Si(Z_\ka(\delta),\sP)$ be the image of $\Delta^\circ(Z_\ka(\delta),\sP)$
under the projection $Z_{\ka}(\delta)\to\Sym^{\ka,\natural}_\delta(\RR^4)$.
We then have the

\begin{lem}
\label{lem:DiagonalsInZ}
If $\kappa\geq 1$ is an integer and
$\sP$ and $\sP'$ are partitions of $N_{\ka}$ and
$\tilde\pi_{\ka}: Z_{\ka}\subset
\oplus_{i\in N_{\ka}}\RR^4\to\Sym^{\ka,\natural}(\RR^4)$ is
the projection, then the following hold:
\begin{enumerate}
\item
\label{item:DiagonalsInZ1}
$\Si(Z_\ka(\delta),\sP)
= \Delta^\circ(Z_\ka(\delta),\sP)/\fS(\sP)
= \Delta^\circ(Z_\ka(\delta),\sP)/W(\sP)$.
\item
\label{item:DiagonalsInZ2}
$\tilde\pi_{\ka}^{-1}\left( \Si(Z_\ka(\delta),\sP)\right)
=
\sqcup_{\sP'\in[\sP]} \Delta^\circ(Z_\ka(\delta),\sP')$.
\item
\label{item:DiagonalsInZ3}
$\Si(Z_\ka(\delta),\sP) \subset \cl\left( \Si(Z_{\ka}(\delta),\sP')\right)$
if and only if there is $\sP''\in [\sP]$ such that $\sP''<\sP$.
\end{enumerate}
\end{lem}

If $(v_1,\dots,v_{\ka})\in
\Delta^\circ(Z_{\ka},\sP)$ and $P\in\sP$, then $v_i=v_j$ for $i,j\in
P$.  Thus, we will write $v_P=v_i$ for any $i\in P$.  With this
notation, we will write $(y_1,\dots,y_\ka)=(y_P)_{P\in\sP}$ for an
element of $\Delta^\circ(Z_{\ka},\sP)$.

\subsection{Tubular neighborhoods}
\label{subsubsec:TubNgh}

Because $Z_\ka(\delta)$ is an open subspace of $Z_\ka$,
the tubular neighborhood of the diagonals in
\eqref{eq:DiagonalsOfZP} can be described by
the restriction of the tubular neighborhoods
described in Lemma \ref{lem:RelatingTrivialStrata},
as we formally state in the following lemma.

\begin{lem}
\label{lem:R4DiagonalsNormal}
Let $\sP$ be a partition of $N_\ka$.
For $P\subseteq N_\ka$,
let $Z_P\subset \oplus_{i\in P}\RR^4$ be the
subspace of mass-centered sets of points
defined in \eqref{eq:DefineZ}.
The normal bundle of the diagonal $\Delta^\circ(Z_\ka(\delta),\sP)$
in $Z_{\ka}$ is
\begin{equation}
\label{eq:R4DiagonalNormal}
\tilde\nu(Z_{\ka}(\delta),\sP)
=
\Delta^\circ(Z_\ka(\delta),\sP)
\times\prod_{P\in\sP} Z_P.
\end{equation}
There is an open neighborhood, $\tilde\sO(Z_{\ka}(\delta),\sP)$,
of the zero section in $\tilde\nu(Z_{\ka}(\delta),\sP)$
such that the restriction of
the exponential map $e(Z_{\ka},\sP)$ defined in \eqref{eq:NormalOfR4Diagonal}
to $\tilde\sO(Z_{\ka}(\delta),\sP)$ is injective and $\fS(\sP)$-equivariant.
\end{lem}

Let $\tilde\sU(Z_{\ka}(\delta),\sP)$ be the image of
$\tilde\sO(Z_{\ka}(\delta),\sP)$ under the exponential map
\eqref{eq:NormalOfR4Diagonal}.
By the $\fS(\sP)$-equivariance of the exponential map
$e(Z_{\ka}(\delta),\sP)$, there is a homeomorphism,
\begin{equation}
\label{eq:R4DiagonalTubNghSymmQuotient}
\sO(Z_{\ka}(\delta),\sP)\cong \sU(Z_{\ka}(\delta),\sP),
\end{equation}
where $\sO(Z_{\ka}(\delta),\sP)=\tilde\sO(Z_{\ka}(\delta),\sP)/\fS(\sP)$
and $\sU(Z_{\ka}(\delta),\sP)=\tilde\sU(Z_{\ka}(\delta),\sP)/\fS(\sP)$.

To describe the pre-image of $\Delta^\circ(Z_\ka(\delta),\sP')$
under $e(Z_{\ka},\sP)$, where $\sP<\sP'$, we introduce the
following notation.

\begin{lem}
\label{lem:R4DiagonalExpPreImage}
Let $\sP<\sP'$ be partitions of $N_\ka$.
For $P\in\sP$, let $\sP'_P$ be the partition of $P$
defined in \eqref{eq:DefineRefinementPartition}.
Then there is an inclusion,
\begin{equation}
\label{eq:DefineR4NormalPP'}
\tilde\nu(Z_{\ka}(\delta),\sP\to\sP')
=
\Delta^\circ(Z_\ka(\delta),\sP)
\times
\prod_{P\in\sP} \Delta^\circ(Z_P,\sP'_P)
\hookrightarrow
\tilde\nu(Z_{\ka}(\delta),\sP),
\end{equation}
with the following significance. If we define
\begin{equation}
\label{eq:DefineTubNghR4Normal}
\tilde\sO(Z_{\ka}(\delta),\sP\to\sP')
:=
\tilde\sO(Z_{\ka}(\delta),\sP)
\cap
\tilde\nu(Z_{\ka}(\delta),\sP\to\sP'),
\end{equation}
then
\[
\tilde\sO(Z_{\ka}(\delta),\sP\to\sP')
=
e(Z_{\ka},\sP)^{-1}(\Delta^\circ(Z_\ka(\delta),\sP')).
\]
\end{lem}

\begin{proof}
We observe that
$$
\Delta^\circ(Z_P,\sP'_P)
=\left\{(v_i)_{i\in P}\in Z_P: v_i=v_j
\iff \exists P'\in\sP'_P \text{ with } i,j\in P'\right\}.
$$
Because $\sP$ is a partition of $N_{\ka}$, there
is an isomorphism,
$$
\bigoplus_{P\in\sP}\bigoplus_{i\in P}\RR^4
\cong
\bigoplus_{i\in N_{\ka}}\RR^4,
$$
which induces an inclusion,
\begin{equation}
\label{eq:NormalBundleInclusion}
\prod_{P\in\sP}\Delta^\circ(Z_P,\sP'_P)
\ni ((v_i)_{i\in P,P\in\sP})
\mapsto (v_i)_{i\in N_{\ka}}
\in Z_{\ka}.
\end{equation}
By the $\fS(\sP)$-equivariance of the exponential map,
$e(Z_{\ka},\sP)$, an element of the pre-image of
$\Delta^\circ(Z_\ka(\delta),\sP')$ under $e(Z_{\ka},\sP)$
can be written as
$$
((x_1,\dots,x_{\ka}),(v_1,\dots,v_{\ka}))\in
\Delta^\circ(Z_\ka,\sP)\times
\oplus_{i\in N_{\ka}}\RR^4
$$
where $(v_1,\dots,v_{\ka})\in \oplus_{i\in N_{\ka}}\RR^4$
satisfies
\begin{enumerate}
\item
$\sum_{i\in P} v_i=0$ for all $P\in\sP$.
\item
For all $i,j$, we have $v_i=v_j$ if and only if there exists
$P'\in\sP'$ with $i,j\in P'$.
\end{enumerate}
Then one observes that the set of elements
$(v_1,\dots,v_{\ka})\in \oplus_{i\in N_{\ka}}\RR^4$
satisfying these two conditions is precisely the
image of the inclusion \eqref{eq:NormalBundleInclusion}.
\end{proof}

By Items \eqref{item:DiagonalsInZ2} and \eqref{item:DiagonalsInZ3} in
Lemma \ref{lem:DiagonalsInZ},
to describe the end of  $\Si(Z_{\ka}(\delta),\sP')$ near
$\Si(Z_\ka(\delta),\sP)$, we must describe not only the
end of $\Delta^\circ(Z_\ka(\delta),\sP')$ near
$\Delta^\circ(Z_\ka(\delta),\sP)$, as is done in Lemma
\ref{lem:R4DiagonalExpPreImage}, but also the ends of
$\Delta^\circ(Z_\ka(\delta),\sP'')$ near
$\Delta^\circ(Z_\ka(\delta),\sP)$ for all $\sP''\in [\sP']$
with $\sP<\sP''$.

Thus, define
\begin{equation}
\label{eq:R4DiagonalOverlap}
\tilde\nu(Z_{\ka}(\delta),[\sP<\sP'])
:=
\Delta^\circ(Z_\ka(\delta),\sP)
\times
\bigsqcup_{\sP''\in [\sP<\sP']}
\left(
\prod_{P\in\sP}
\Delta^\circ(Z_P,\sP''_P)
\right).
\end{equation}
There is an action of $\fS(\sP)$ on $\tilde\nu(Z_{\ka},[\sP<\sP'])$
defined by the standard action on $\Delta^\circ(Z_\ka,\sP)$
and the action of $\fS(\sP)$ on $[\sP<\sP']$.
Although $\Delta^\circ(Z_\ka(\delta),\sP)$
is not a subset of $\tilde\nu(Z_{\ka}(\delta),[\sP<\sP'])$,
we use the phrase {\em a neighborhood of $\Delta^\circ(Z_\ka(\delta),\sP)$
in $\nu(Z_{\ka}(\delta),[\sP<\sP'])$\/} to refer to the obvious parallel with
subspaces of $\nu(Z_{\ka}(\delta),\sP)$.

\begin{lem}
\label{lem:SymmProdNeigh}
Let $\sP$ and $\sP'$ be partitions of $N_{\ka}$ with
$\sP<\sP'$.  Then a neighborhood of $\Si(Z_\ka(\delta),\sP)$ in
$\Si(Z_{\ka}(\delta),\sP')$ is homeomorphic to a
neighborhood, $\sO(Z_\ka(\delta),[\sP<\sP'])$,
of the zero section in
$$
\tilde\nu(Z_{\ka}(\delta),[\sP<\sP'])
/\fS(\sP).
$$
\end{lem}

\begin{proof}
The bundle $\tilde\nu(Z_{\ka}(\delta),[\sP<\sP'])$ is the union
of the bundles $\tilde\nu(Z_{\ka}(\delta),\sP\to\sP'')$ defined
in \eqref{eq:DefineR4NormalPP'} as $\sP''$ varies in
$[\sP<\sP']$ along the subspace $\Delta^\circ(Z_\ka(\delta),\sP)$
of each of these bundles.
Therefore, a neighborhood of
$\Delta^\circ(Z_\ka(\delta),\sP)$ in
\begin{equation}
\label{eq:UnionOfDiags}
\bigcup_{\sP''\in [\sP<\sP']} \Delta^\circ(Z_\ka(\delta),\sP'')
\end{equation}
is homeomorphic to a neighborhood of the zero-section in
$\nu(Z_{\ka}(\delta),[\sP<\sP'])$.  By the characterization of
the pre-image of $\Si(Z_\ka(\delta),\sP)$ under the projection
$Z_{\ka}(\delta)\to\Sym^{\natural,\ka}_\delta(\RR^4)$ in
Item \eqref{item:DiagonalsInZ2} of Lemma
\ref{lem:DiagonalsInZ}, a neighborhood of
$\Si(Z_\ka(\delta),\sP)$ in $\Si(Z_\ka(\delta),\sP')$
is homeomorphic to the quotient of \eqref{eq:UnionOfDiags}
by $\fS(\sP)$, thus proving the lemma.
\end{proof}

\section{The splicing map with the product connection over $\RR^4$}
\label{subsec:SplicingMapsR4}
In this section, we define the
basic operation of splicing
onto the
product connection
as used in \cite{FL3}.  Our main
result here is in Lemma \ref{lem:CommutingSplicing}, where we show that
the composition of two splicing maps is equal to a single
splicing map on suitably small open sets.

Let $\sP$ be a partition of $N_{\ka}$. We define the splicing map
on an open subspace of
$$
\Delta^\circ(Z_\ka(\delta),\sP) \times_{\fS(\sP)}
\prod_{P\in\sP}\bar{\sB}^{s,\natural}_{|P|}(S^4,\delta_P),
$$
where the constants $\delta_P$ are invariant under the action of
$\fS(\sP)$ and
where $\fS(\sP)$ acts on the generalized connections by permuting
them, sending an element of $\bar{\sB}^{s,\natural}_{|P|}(S^4,\delta_P)$
to an element of $\bar{\sB}^{s,\natural}_{|\si(P)|}(S^4,\delta_{|\si(P)|})$.
Note that elements of $\Delta^\circ(Z_\ka(\delta),\sP)$
can be written as $(y_P)_{P\in\sP}$, where
$y_P\in\RR^4$, while elements of
$\prod_{P\in\sP}\bar{\sB}^{s,\natural}_{|P|}(S^4,\delta_P)$ can be written
as $([A_P,F^s_P,\bx_P])_{P\in\sP}$, where
$[A_P,F^s_P,\bx_P]\in \bar{\sB}^{s,\natural}_{|P|}(S^4,\delta_P)$.
The open subspace on which the splicing map will be defined is
\begin{multline}
\label{eq:DefineSplicingDataR4}
\tilde\sO(\Theta,\sP,\delta)
:=
\left\{ \left( y_P,[A_P,F^s_P,\bx_P]\right)_{P\in\sP}\in
\Delta^\circ(Z_\ka(\delta),\sP) \times \prod_{P\in\sP}\bar{\sB}^s_{|P|}(S^4,\delta_P):\right.
\\
\left. 8\sqrt{\la[A_P,\bx_P]} + 8\sqrt{\la[A_{P'},\bx_{P'}]} < \dist(y_P,y_{P'}),
\ \forall\, P\neq P' \right\}.
\end{multline}
We define the splicing map,
\begin{equation}
\label{eq:R4TrivialSplicing}
\bga'_{\Theta,\sP}:\sO(\Theta,\sP,\delta)
:=
\tilde\sO(\Theta,\sP,\delta)/\fS(\sP)\to \bar{\sB}^s_{\ka}(S^4,2\delta),
\end{equation}
as follows. For $P\in\sP$, let $E_P\to S^4$ be the
Hermitian vector
bundle supporting the connection
$A_P$ and let $\Fr(E_P)$ be the bundle of unitary frames of $E_P$.
The framed connection $[A_P,F^s_P]$ defines a section $\phi(A_P,F^s_P)$ of
$\Fr(E_P)$ over the complement of the North Pole
by parallel translation of $F^s_P$ with
respect to the connection $A_P$,
along great circles from the South Pole. Let
$\Theta(A_P,F^s_P)$ be the product
connection defined by this section.  Note that if $A_P$ is flat on the domain of
$\phi(A_P,F^s_P)$, then $\Theta(A_P,F^s_P)=A_P$.

If $\Theta$ denotes the
product
connection on
$S^4\times \CC^2$,
then the product connection $\Theta(A_P,F^s_P)$ is
identified with the restriction of $\Theta$ to the
domain of $\phi(A_P,F^s_P)$ as follows. The section
$\phi(A_P,F^s_P)$ gives a trivialization of the bundle $E$ and
thus identifies it with the restriction of the
product
bundle
$S^4\times \CC^2$
to the domain of the section. Under this
identification the connection $\Theta$ is identified with
$\Theta(A_P,F^s_P)$.

We also note that given any two connections, $A_1$ and $A_2$, a
convex linear combination of them, $tA_1 + (1-t)A_2$, also defines
a connection.  For example, if $A_1$ is a connection which equals the product connection
in a particular trivialization and $A_2=A_1+a_2$, we
would write the connection one form for $tA_1 + (1-t)A_2$ in this
trivialization as $(1-t)a_2$.

Given $(y,\la) \in \RR^4\times (0,\infty)$,
let $c_{y,\la}:\RR^4\to \RR^4$ be defined by
\[
c_{y,\la}(z) := (z-y)/\la,
\]
and so $c_{y,\la}$ maps the ball $B(y,\la)$ onto $B(0,1)$.

Let $\beta:\RR\to [0,1]$ be a smooth function satisfying $\beta(x)=0$ for $x\le 1/2$ and $\beta(x)=1$ for $x\ge 1$. Define
$\chi:\RR^4\to [0,1]$ by $\chi(x)=\beta(|x|)$ and
$\chi_{y,\la}:\RR^4\to [0,1]$ by
$\chi_{y,\la}(x)=(c_{y,\la}^*\chi)(x)=\beta(|x-y|/\la)$.
The function $1-\chi_{y,\la}$ is thus supported on the ball
$B(y,\la)$ and equal to one on the ball $B(y,\frac{1}{2}\la)$.

For $\by=(y_P)_{P\in\sP}\in \Delta^\circ(Z_\ka(\delta),\sP)$,
we define
the splicing map \eqref{eq:R4TrivialSplicing} by setting
\begin{equation}
\label{eq:SplicedGenConn}
\bga'_{\Theta,\sP}
\left(
\by,([A_P,F^s_P,\bx_P])_{P\in\sP}
\right)
:=
[A',F^s,\bx'],
\end{equation}
where, for
$\la_P=\la([A_P,F^s_P,\bx_P])$,
we define
the framed connection $[A',F^s]$  by
\begin{align}
\label{eq:SplicedConnection}
A'
&=
\begin{cases}
\Theta & \text{on } \RR^4\less \cup_{P\in\sP}B(y_P,4 \sqrt{\la_P}),
 \\
(1-\chi_{y_P,4 \sqrt{\la_P}})c_{y_P,1}^*A_P +
\chi_{y_P,4 \sqrt{\la_P}}c_{y_P,1}^*\Theta(A_P,F^s_P)
{}&\text{on }\Omega(y_P;2\sqrt{\la_P},4\sqrt{\la_P}),
\\
c_{y_P,1}^*A_P & \text{on } B(y_P,2\sqrt{\la_P}),
\end{cases}
\end{align}
with the frame $F^s$ in $[A',F^s]$ being given by the
canonical frame for the trivialization associated with the product connection, $\Theta$.
If the points $\bx_P$ in \eqref{eq:SplicedGenConn} are given by points $x_{P,i}\in\RR^4$
with multiplicities $\ka_{P,i}$, then the point  $\bx'$ in the expression \eqref{eq:SplicedGenConn} is defined by
the points $c_{y_P,1}^{-1}(x_{P,i})$ with multiplicities $\ka_{P,i}$.

Note that we have not computed the behavior of the function
$\la\circ\bga'_{\Theta,\sP}$ precisely, so we cannot assert
that the image of $\bga'_{\Theta,\sP}$ is contained in
$\sB^s_{\ka}(2\delta)$.  However, that containment will follow
from the continuity of $\la$ and of $\bga'_{\Theta,\sP}$
with respect to Uhlenbeck limits if we shrink the domain
$\sO(\Theta,\sP,\delta)$ by requiring the scales $\la_P$ to
be sufficiently small.

We have the following result regarding the behavior of this
splicing map at the cone point of $\bar\sB^s_{|P|}(\delta_P)$.

\begin{lem}
\label{lem:UhlebeckLimitCont}
For any partition $\sP$ of $N_{\ka}$, the map $\bga'_{\Theta,\sP}$ is smoothly stratified.
If $c_{|P|}\in\Sym^{|P|,\natural}_{\delta}(\RR^4)$
denotes the cone point, then for any
$\by\in\Delta^\circ(Z_\ka(\delta),\sP)$,
\begin{equation}
\label{eq:SplicingConePoint}
\bga'_{\Theta,\sP}
\left(
\by,([\Theta,F^s,c_{|P|})_{P\in\sP})
\right)
=
[\Theta,F^s,\by].
\end{equation}
\end{lem}
\begin{proof}
The fact that the map  $\bga'_{\Theta,\sP}$ is smooth when restricted to the subspaces defined by
the strata of $\bar\sB^{s,\natural}_{|P|}(\delta_P)$ is clear from the definition.
If a sequence
$$
\{[A_{P,\alpha},F^s_{P,\alpha},\bx_{P,\alpha}]\}_{\alpha=1}^\infty
\subset
\bar\sB^{s,\natural}_{|P|}(\delta_P)
$$
converges in the Uhlenbeck topology to
$[A_{P,0},F^s_{P,0},\bx_{P,0}]\in \bar\sB^{s,\natural}_{|P|}(\delta_P)$, then by the definition of
Uhlenbeck convergence,
$$
\lim_{\alpha\to\infty} \la [A_{P,\alpha},F^s_{P,\alpha},\bx_{P,\alpha}] =
\la [A_{P,0},F^s_{P,0},\bx_{P,0}]:=\la_{P,0}.
$$
By \eqref{eq:ConnectionScale},
none of the points in $\bx_{P,0}$ can lie outside of the ball $B(0,\la_{P,0})\subset B(0,2\sqrt{\la_{P,0}})$
and so the sequence of connections $\{A_{P,\alpha}\}$ cannot ``bubble'' outside of the ball $B(0,2\sqrt{\la_{P,0}})$.
These observations and the construction of $\bga'_{\Theta,\sP}$ in \eqref{eq:SplicedConnection} imply
that $\bga'_{\Theta,\sP}$ is continuous with respect to these Uhlenbeck limits.

The equality \eqref{eq:SplicingConePoint} follows immediately from the construction \eqref{eq:SplicedConnection}.
\end{proof}

It is important that our constructions in Section \ref{subsec:SplicedEnd} be equivariant with respect to the
$\SO(3)\times\SO(4)$ actions on the space of connections
defined in \eqref{eq:ExtendedFrameAction} and \eqref{eq:ExtendedRotationAction}, so that
the instanton moduli space with spliced ends can be used in the space of gluing data.
The action of $\SO(3)\times\SO(4)$ on the domain \eqref{eq:DefineSplicingDataR4} is defined
as follows.  For $R\in\SO(4)$, let $\tilde R_P$ be any automorphism of
$\Fr(E_P)$ covering $R:S^4\to S^4$. For $A\in\SO(3)$ and $R\in \SO(4)$, define
\begin{equation}
\label{eq:SO3SO4Actions}
\left( \left((y_P)_{P\in\sP},
[A_P,F^s_P,\bx_P]_{P\in\sP}\right),A,R \right)
\mapsto
\left( (Ry_P)_{P\in\sP}, [(\tilde R_P^{-1})^*A_P,\tilde R_P F^s_PA^{-1},R\bx_P]_{P\in\sP} \right).
\end{equation}
We have the following lemma.

\begin{lem}
\label{lem:GroupEquiv}
Let $\sP$ be a partition of $N_{\ka}$.  The splicing map $\bga_{\Theta,\sP}'$ defined in
\eqref{eq:SplicedGenConn} is equivariant with respect to the $\SO(3)\times\SO(4)$ actions
defined in \eqref{eq:SO3SO4Actions} on its domain and in
\eqref{eq:ExtendedFrameAction} and \eqref{eq:ExtendedRotationAction} on its image.
\end{lem}

\begin{proof}
The equivariance of $\bga_{\Theta,\sP}'$ with respect to the $\SO(3)$ actions follows immediately from
the construction \eqref{eq:SplicedConnection}.

The construction of the splicing map in \eqref{eq:SplicedConnection} defines,
for each point
$$
\bA=\left( (y_P)_{P\in\sP},[A_P,F^s_P,\bx_P]_{P\in\sP}\right)
\in\tilde\sO(\Theta,\sP,\delta),
$$
an embedding of principal bundles,
$$
\iota_\bA: \bigsqcup_{P\in\sP} \Fr(E_P)|_{B(0,4 \sqrt{\la_P})}\to \Fr(E_\ka),
$$
covering the embedding $B(0,4 \sqrt{\la_P})\to S^4$ given by $x_P\mapsto y_P+x_P$ on the $P$-th component.
If $\tilde R_\ka:\Fr(E_\ka)\to\Fr(E_\ka)$ is a bundle map
covering $R:S^4\to S^4$,
then $\tilde R\circ\iota_A$ equals, up to a
gauge transformation, the embedding $\iota_{R\bA}$, where
$$
R\bA=\left( (Ry_P)_{P\in\sP}, [(\tilde R_P^{-1})^*A_P,\tilde R_P F^s_P,R\bx_P]_{P\in\sP} \right).
$$
The desired equivariance, $(\tilde R_\ka^{-1})^*\bga_{\Theta,\sP}'(\bA)=\bga_{\Theta,\sP}'(R\bA)$, then
follows by considering the images of the horizontal distributions  defining the connections under these
embeddings.
\end{proof}

We will also require the following observation about centering
maps.

\begin{lem}
\label{lem:ComposeCenterings}
Let $c_{y,\la}:\RR^4\to \RR^4$ be defined by
$c_{y,\la}(z)=(z-y)/\la$.
For any $\la_1,\la_2>0$ and
$x_1,x_2\in\RR^4$, we have:
\begin{align*}
c_{x_2,\la_2}\circ c_{x_1,\la_1} &= c_{x_1+\la_1x_2,\la_1\la_2},
\\
c_{x_2,\la_2}^*\chi_{x_1,\la_1} &= \chi_{x_2+\la_2x_1,\la_1\la_2}.
\end{align*}
\end{lem}

If one splices
a mass-centered connection, of charge $\ka$
and scale $\la$ at a point $y\in\RR^4$,
with the product connection
it is not immediately clear what the scale of the resulting connection will
be. However, the limit of the scale of the spliced connection
(as $\la$ tends to zero)  will be $\sqrt{\ka |y|^2}$. The following
exploitation of this fact will be used in Lemma
\ref{lem:CommutingSplicing} to prove  that when composing two splicing maps--- for a suitable open set
of splicing data --- the annulus on which the second splicing map
interpolates between connections (see the second line of \eqref{eq:SplicedConnection})
does not intersect the annuli of the first splicing map.

\begin{lem}
\label{lem:ScaleLimit}
Let $\Theta$ be the product connection and $c_{|P|}\in\Sym^{|P|,\natural}_{\delta}(\RR^4)$ the cone point.
If
\begin{equation}
\label{eq:DefineTrivialStrataInR4GluingSpace}
T(\Theta,\sP,\delta)
:=
\left\{\left(y_P,[\Theta,F^s_P,c_{|P|}]\right)_{P\in\sP}\in \sO(\Theta,\sP,\delta)\right\},
\end{equation}
then, for $\Si(Z_\ka(\delta),\sP)$ as defined in \eqref{eq:DiagonalsOfZP}, we have
\begin{equation}
\label{eq:SplicingOnTrivialStratum}
\bga'_{\Theta,\sP}(T(\Theta,\sP,\delta))
=
[\Theta]|\times \Si(Z_\ka(\delta),\sP),
\end{equation}
and there is an open neighborhood, $D(\Theta,\sP,\delta)$, of
$T(\Theta,\sP,\delta)$ in $\sO(\Theta,\sP,\delta)$, such that for every
$$
\bA:= \left( y_P,[A_P,F^s_P,\bx_P] \right)_{P\in\sP} \in
D(\Theta,\sP,\delta)
$$
and for every $P\in\sP$, we have
\begin{equation}
\label{eq:BallContainment}
B\left(y_P,4\sqrt{\la_P}\right)\Subset B(0,2\sqrt{\la}),
\end{equation}
and
\begin{equation}
\label{eq:BallContainment2}
B\left(y_P,\tquarter\la_P^{1/3}\right)\Subset B(0,\teighth\la^{1/3}),
\end{equation}
where  $\la_P=\la([A_P,F^s_P,\bx_P])$ and $\la=\la(\bga'_{\Theta,\sP}(\bA))$.
\end{lem}

\begin{proof}
The equality \eqref{eq:SplicingOnTrivialStratum} follows immediately
from Lemma \ref{lem:UhlebeckLimitCont}.

The functions $\la_P$ vanish on $T(\Theta,\sP,\delta)$
while by Lemma \ref{lem:UhlebeckLimitCont}, the composition
$\la\circ\bga'_{\Theta,\sP}$ is equal to the square root of
$\sum_{P\in\sP}|P||x_P|^2$ and thus is non-zero on $T(\Theta,\sP,\delta)$.
The
continuity of the functions involved then yields the
existence of the desired open subspace.
\end{proof}

\section{Composition of splicing maps}
\label{subsec:OverlapSpliceMapsR4}
To construct the instanton moduli space with spliced ends, we need to show
that the union of the images of the splicing maps $\bga'_{\Theta,\sP}$
form a smoothly-stratified
space. Let $\sP,\sP'$ be partitions of $N_{\ka}$ with $\sP<\sP'$.
For each $P\in\sP$, let $\sP'_P$ be the partition of $P$
defined in \eqref{eq:DefineRefinementPartition}. We will need to
understand the overlaps of the images of the maps $\bga'_{\Theta,\sP}$ and
$\bga'_{\Theta,\sP'}$. To that end,  we define a space
$\tilde\sO(\Theta,\sP,\sP',\delta)$ with maps to $\tilde\sO(\Theta,\sP,\delta)$
and $\tilde\sO(\Theta,\sP',\delta)$ so that the following diagram commutes:
\begin{equation}
\label{eq:CommutingSplice0}
\begin{CD}
\tilde\sO(\Theta,\sP,\sP',\delta)
@>\rho^{\Theta,d}_{\sP,\sP'} >>
\tilde\sO(\Theta,\sP,\delta)
\\
@V \rho^{\Theta,u}_{\sP,\sP'}   VV @V \bga'_{\Theta,\sP} VV
\\
\tilde\sO(\Theta,\sP',\delta)
@> \bga'_{\Theta,\sP'} >>
\bar{\sB}^s_{\ka}(S^4,2\delta)
\end{CD}
\end{equation}
The superscripts $d$ and $u$ appearing in the maps
$\rho^{\Theta,d}_{\sP,\sP'}$ and $\rho^{\Theta,d}_{\sP,\sP'}$
suggest ``down'' and ``up'', respectively, as
$\sP$ and $\sP'$ correspond to the lower and upper strata.

\subsection{Definition of the overlap data}
\label{subsec:Definition_overlap data}
We now define the objects appearing in the diagram
\eqref{eq:CommutingSplice0}.
The open set
$\tilde\sO(\Theta,\sP,\sP',\delta)$ will be defined to be a subspace of
\begin{equation}
\label{eq:R4OverlapSpace}
\tilde\nu(\Theta,\sP,\sP',\delta)
:=
\Delta^\circ(Z_\ka(\delta),\sP) \times
\prod_{P\in\sP}
\left(
    \Delta^\circ(Z_{|P|}(\delta_P),\sP'_P)
    \times
    \prod_{P'\in\sP_P'} \bar{\sB}^s_{|P'|}(S^4,\delta_{P'})
\right),
\end{equation}
where the constants $\delta_P$ and $\delta_{P'}$ will not need to be
defined explicitly.
The open set $\tilde\sO(\Theta,\sP,\sP',\delta)$  will be defined by the conditions
\eqref{eq:OverlapCondition1}, \eqref{eq:OverlapCondition2}, \eqref{eq:OverlapCondition3}, and
\eqref{eq:OverlapCondition4}.

An element $\si\in\fS(\sP)$ with $\si(\sP')=\sP''$
defines  a bijection,
\begin{equation}
\label{eq:SymmGrpActionOnR4OverlapSpace}
\si:
\tilde\nu(\Theta,\sP,\sP',\delta)
\to
\tilde\nu(\Theta,\sP,\sP'',\delta),
\end{equation}
by the diagonal action on each of the factors.  That is,
$\si$ defines bijections,
\begin{align*}
{}&
\Delta^\circ(Z_\ka(\delta),\sP)
\to
\Delta^\circ(Z_\ka(\delta),\sP),
\\
{}&
\Delta^\circ(Z_{|P|}(\delta_P),\sP'_P)
\to
\Delta^\circ(Z_{|\si(P)|}(\delta_{\si(P)},\sP''_{\si(P)}),
\\
{}&
\bar{\sB}^s_{|P'|}(S^4,\delta_{P'})
\to
\bar{\sB}^s_{|\si(P')|}(S^4,\delta_{\si(P')}),
\end{align*}
which define the desired bijection \eqref{eq:SymmGrpActionOnR4OverlapSpace}.

The following maps
are given by the obvious projections onto the factors of
$\tilde\nu(\Theta,\sP,\sP',\delta)$:
\begin{equation}
\label{eq:DefineProjections}
\begin{aligned}
 \pi_{\sP,\bx} &: \tilde\nu(\Theta,\sP,\sP',\delta) \to
 \Delta^\circ(Z_\ka(\delta),\sP),
\\
\pi_{\sP'_P,\bx} &: \tilde\nu(\Theta,\sP,\sP',\delta) \to
\Delta^\circ(Z_{|P|}(\delta_P),\sP'_P),
\\
\pi_{\bx} &: \tilde\nu(\Theta,\sP,\sP',\delta) \to
\Delta^\circ(Z_\ka(\delta),\sP)\times
    \prod_{P\in\sP}\Delta^\circ(Z_{|P|}(\delta_P),\sP'_P),
\\
\pi_{P} &: \tilde\nu(\Theta,\sP,\sP',\delta) \to
\Delta^\circ(Z_{|P|}(\delta_P),\sP'_P)\times \prod_{P'\in\sP'_P}\bar{\sB}^s_{|P'|},
\\
\pi_{P'}&: \tilde\nu(\Theta,\sP,\sP',\delta) \to \bar{\sB}^s_{|P'|}(\delta_P).
\end{aligned}
\end{equation}
We use the inclusion given in \eqref{eq:DefineR4NormalPP'},
$$
\tilde\nu(Z_\ka(\delta),\sP\to\sP')
=\Delta^\circ(Z_\ka(\delta),\sP)
\times
\prod_{P\in\sP} \Delta^\circ(Z_P(\delta_P),\sP'_P)
\to
\tilde\nu(Z_{\ka}(\delta),\sP),
$$
to view
$\tilde\nu(Z_\ka(\delta),\sP\to\sP')$ as a subspace of
$\tilde\nu(Z_\ka(\delta),\sP)$.  By Lemma \ref{lem:R4DiagonalExpPreImage},
the image of the restriction of the exponential map
$e(Z_\ka,\sP)$ in \eqref{eq:NormalOfR4Diagonal}
to elements of $\tilde\nu(Z_\ka(\delta),\sP\to\sP')$
is contained in $\Delta^\circ(Z_{\ka}(\delta),\sP')$.
We would like to define the map $\rho^{\Theta,u}_{\sP,\sP'}$
appearing in the diagram
\eqref{eq:CommutingSplice0} by
\begin{equation}
\label{eq:R4TrivialSplicingUpwardsOverlapMap}
\rho^{\Theta,u}_{\sP,\sP'}
:=
(e(Z_\ka,\sP)\circ\pi_{\bx}) \times
\prod_{P'\in\sP'} \pi_{P'}.
\end{equation}
That is, the map $\rho^{\Theta,u}_{\sP,\sP'}$ leaves the connection
data in $\nu(\Theta,\sP,\sP',\delta)$
unchanged and maps the points in $\tilde\nu(Z_{\ka}(\delta),\sP\to\sP')$ to
the diagonal $\Delta^\circ(Z_\ka(\delta),\sP')$.

For the composition $e(Z_\ka,\sP)\circ\pi_{\bx}$
in \eqref{eq:R4TrivialSplicingUpwardsOverlapMap}
to be defined, we must restrict the domain of
$\rho^{\Theta,u}_{\sP,\sP'}$ to an open subspace
$\tilde\sO(\Theta,\sP,\sP',\delta)$ of $\tilde\nu(\Theta,\sP,\sP',\delta)$ satisfying
\begin{equation}
\label{eq:OverlapCondition1}
\tilde\sO(\Theta,\sP,\sP',\delta) \subset
\pi_{\bx}^{-1}\left( \tilde\sO(Z_{\ka}(\delta),\sP\to\sP') \right),
\end{equation}
where $\tilde\sO(Z_{\ka}(\delta),\sP\to\sP')$ is defined in \eqref{eq:DefineTubNghR4Normal}.

To define the composition
$\bga'_{\Theta,\sP'}\circ\rho^{\Theta,u}_{\sP,\sP'}$
appearing in the diagram \eqref{eq:CommutingSplice0},
the open subspace $\tilde\sO(\Theta,\sP,\sP',\delta)$ must also satisfy
\begin{equation}
\label{eq:OverlapCondition2}
\tilde\sO(\Theta,\sP,\sP',\delta) \subset
(\rho^{\Theta,u}_{\sP,\sP'})^{-1}\left( \tilde\sO(\Theta,\sP',\delta)\right),
\end{equation}
where $\tilde\sO(\Theta,\sP',\delta)$, defined in \eqref{eq:DefineSplicingDataR4}, is the domain of $\bga'_{\Theta,\sP'}$.

Before we proceed to define $\rho^{\Theta,d}_{\sP,\sP'}$, we note that
if $P\in\sP\cap\sP'$, then the partition $\sP'_P$ is just the set $P$
and $\Delta^\circ(Z_{|P|},\sP'_P)$ is
a single point (the zero vector).
In this case, we define the splicing map,
$$
\bga'_{\Theta,\sP'_P}:
\Delta^\circ(Z_{|P|}(\delta),\sP'_P)\times \bar\sB^s_{|P|}(S^4,\delta_P)
\to \bar\sB^s_{|P|}(S^4,2\delta),
$$
to be the projection onto
$\bar\sB^s_{|P|}(S^4,\delta_P)\subset\bar\sB^s_{|P|}(S^4,2\delta)$.
The map $\rho^{\Theta,d}_{\sP,\sP'}$ of \eqref{eq:CommutingSplice0} will
then be defined by
\begin{equation}
\label{eq:DefineUpwardTransition}
\rho^{\Theta,d}_{\sP,\sP'}
:=
\pi_{\sP,\bx}\times
\left(\prod_{P\in\sP}\bga'_{\Theta,\sP'_P}\circ\pi_{P}\right).
\end{equation}
That is, the map $\rho^{\Theta,d}_{\sP,\sP'}$ leaves the point
in $\Delta^\circ(Z_\ka,\sP)$ unchanged and for each $P\in\sP$
applies the splicing map $\bga'_{\Theta,\sP'_P}$ to the space
$$
\Delta^\circ(Z_{|P|}(\delta_P),\sP'_P)\times
\prod_{P'\in\sP'_P} \bar\sB^s_{|P'|}(S^4,\delta_{P'}).
$$
For the map
$\rho^{\Theta,d}_{\sP,\sP'}$
to be defined on $\tilde\sO(\Theta,\sP,\sP',\delta)$, we need
the image of $\pi_P$ to be contained in the
domain, $\sO(\Theta,\sP'_P,\delta_P)$, of the splicing map $\bga'_{\Theta,\sP'_P}$.
However, we will make a stronger requirement (to be used in Lemma
\ref{lem:CommutingSplicing}),
\begin{equation}
\label{eq:OverlapCondition3}
\tilde\sO(\Theta,\sP,\sP',\delta) \subset
\pi_P^{-1}\left( D(\Theta,\sP'_P,\delta_P)\right) \quad \text{for all
$P\in\sP\less\sP'$},
\end{equation}
where $D(\Theta,\sP'_P,\delta_P)$ is the open subspace of
$\tilde\sO(\Theta,\sP'_P,\delta_P)$ defined in Lemma \ref{lem:ScaleLimit}.  Since
\[
D(\Theta,\sP'_P,\delta_P)\subset \sO(\Theta,\sP'_P,\delta_P),
\]
the condition
\eqref{eq:OverlapCondition3} implies that $\rho^{\Theta,d}_{\sP,\sP'}$ is
defined on $\tilde\sO(\Theta,\sP,\sP',\delta)$.

Finally, for the composition $\bga'_{\Theta,\sP}\circ \rho^{\Theta,d}_{\sP,\sP'}$
to be well-defined,
$\tilde\sO(\Theta,\sP,\sP',\delta)$ must satisfy
\begin{equation}
\label{eq:OverlapCondition4}
\tilde\sO(\Theta,\sP,\sP',\delta) \subset
\left(\rho^{\Theta,d}_{\sP,\sP'}\right)^{-1} \left( \tilde\sO(\Theta,\sP,\delta) \right).
\end{equation}
In the proof of Lemma \ref{lem:OverlappingTrivial} below, we will show that the image of the restriction
of $\bga'_{\Theta,\sP'}\circ\rho^{\Theta,u}_{\sP,\sP'}$ to $\tilde\sO(\Theta,\sP,\sP',\delta)$  contains
\begin{equation}
\label{eq:TrivialStrataInImageOfDiagonalOverlapMap}
\{[\Theta,\bx]: \bx\in e(Z_\ka,\sP)(\sO(Z_\ka(\delta),\sP,\sP'))\}.
\end{equation}
Hence, the image of $\bga'_{\Theta,\sP'}\circ\rho^{\Theta,u}_{\sP,\sP'}$
will contain the intersection of the stratum $[\Theta]\times\Si(Z_\ka(\delta),\sP')$
with the image of $\bga'_{\Theta,\sP}$.

\begin{lem}
\label{lem:OverlappingTrivial}
Let $\sP,\sP'$ be partitions of
$N_{\ka}$ satisfying $\sP<\sP'$. Define a
subspace of $\tilde\nu(\Theta,\sP,\sP',\delta)$ by
\begin{equation}
\label{eq:DefineTrivialOverlap}
T(\Theta,\sP,\sP',\delta)
:=
\pi_{\bx}^{-1}\left(\tilde\sO(Z_{\ka}(\delta),\sP\to\sP')\right) \cap \left(
\bigcap_{P'\in\sP'} \pi_{P'}^{-1}([\Theta,c_{|P'|}])\right).
\end{equation}
Then the following hold:
\begin{enumerate}
\item
\label{item:OverlappingTrivial_1}
$(\bga'_{\Theta,\sP'}\circ\rho^{\Theta,u}_{\sP,\sP'} )
\left( T(\Theta,\sP,\sP',\delta) \right)
=
\left\{[\Theta,\bx]: \bx\in
e(Z_\ka,\sP)\left(\tilde\sO(Z_{\ka}(\delta),\sP\to\sP')\right)\right\}$,
\item
\label{item:OverlappingTrivial_2}
$(\bga'_{\Theta,\sP}\circ \rho^{\Theta,d}_{\sP,\sP'})
\left( T(\Theta,\sP,\sP',\delta)\right)
=
 \left([\Theta]\times
e(Z_\ka,\sP)(\tilde\sO(Z_{\ka}(\delta),\sP\to\sP'))\right)$.
\end{enumerate}
Moreover,
there is an open neighborhood, $\tilde\sO(\Theta,\sP,\sP',\delta)$, of
$T(\Theta,\sP,\sP',\delta)$ in
$\tilde\nu(\Theta,\sP,\sP',\delta)$ satisfying the conditions
\eqref{eq:OverlapCondition1}, \eqref{eq:OverlapCondition2},
\eqref{eq:OverlapCondition3}, and \eqref{eq:OverlapCondition4}.
\end{lem}

\begin{proof}
Item \eqref{item:OverlappingTrivial_1}
follows from the definitions of the spaces and maps
and the value of the splicing map given
in \eqref{eq:SplicingConePoint}.
Equation \eqref{eq:SplicingConePoint} implies that
$$
\left.\rho^{\Theta,d}_{\sP,\sP'}\right|_{T(\Theta,\sP,\sP',\delta)}
=
\pi_{\bx}|_{T(\Theta,\sP,\sP',\delta)}.
$$
The preceding identity and the equality,
\begin{equation}
\label{eq:piXofT}
\pi_{\bx}(T(\Theta,\sP,\sP',\delta)) =\tilde \sO(Z_\ka(\delta),\sP\to\sP'),
\end{equation}
which follows immediately from the definition of $T(\Theta,\sP,\sP',\delta)$,
yields Item \eqref{item:OverlappingTrivial_2}.

Equation \eqref{eq:piXofT}
implies that there is an open neighborhood of
$T(\Theta,\sP,\sP',\delta)$ in the space
$\tilde\nu(\Theta,\sP,\sP',\delta)$ satisfying
\eqref{eq:OverlapCondition1}. For every $P'\in\sP'$, we have
$$
\pi_{P'}(T(\Theta,\sP,\sP',\delta))=[\Theta,c_{|P'|}],
$$
so
$\pi_P(T(\Theta,\sP,\sP',\delta))\subset T(\Theta,\sP'_P,\delta_P)
\subset D(\Theta,\sP,\delta)$ for
all $P\in\sP$ and
there is an open neighborhood of
$T(\Theta,\sP,\sP',\delta)$ in $\tilde\nu(\Theta,\sP,\sP',\delta)$ satisfying
\eqref{eq:OverlapCondition3}.

Item \eqref{item:OverlappingTrivial_1} implies that
$\rho^{\Theta,u}_{\sP,\sP'}(T(\Theta,\sP,\sP',\delta)) \subset
\tilde\sO(\Theta,\sP',\delta)$, so there is an open neighborhood of
$T(\Theta,\sP,\sP',\delta)$ in $\tilde\nu(\Theta,\sP,\sP',\delta)$ satisfying
\eqref{eq:OverlapCondition2}.

Item \eqref{item:OverlappingTrivial_2}  implies
that $\rho^{\Theta,d}_{\sP,\sP'}(T(\Theta,\sP,\sP',\delta))$ is contained in
the domain of $\bga'_{\Theta,\sP'}$,
namely, $\sO(\Theta,\sP,\delta)$, so there is an open neighborhood of
$T(\Theta,\sP,\sP',\delta)$ in $\tilde\nu(\Theta,\sP,\sP',\delta)$ satisfying
\eqref{eq:OverlapCondition4}.

The intersection of these open neighborhoods of
$T(\Theta,\sP,\sP',\delta)$ in $\tilde\nu(\Theta,\sP,\sP',\delta)$, yields
the desired open neighborhood $\tilde\sO(\Theta,\sP,\sP',\delta)$.
\end{proof}

\subsection{Equality of splicing maps}
\label{subsec:Equality of splicing maps}
We now prove that the diagram \eqref{eq:CommutingSplice0} commutes. The key point in this proof is the
restriction \eqref{eq:OverlapCondition3} on the domain $\tilde\sO(\Theta,\sP,\sP',\delta)$
which, as shown in Lemma \ref{lem:ScaleLimit}, ensures that
the connections given by $\rho^{\Theta,d}_{\sP,\sP'}$ are \emph{already} equal to
the product connection on the annuli on which the splicing map
$\bga'_{\Theta,\sP}$ interpolates between them and the product
connection.  Hence, in the iterated splicing construction defining the composition
$\bga'_{\Theta,\sP}\circ \rho^{\Theta,d}_{\sP,\sP'}$, the interpolation in the
definition of $\bga'_{\Theta,\sP}$ does not change the connection.
This iterated splicing procedure is thus equivalent to a single splicing, that given by the composition
$\bga'_{\Theta,\sP'}\circ \rho^{\Theta,u}_{\sP,\sP'}$.

\begin{lem}
\label{lem:CommutingSplicing}
Let $\sP$ and $\sP'$ be partitions
of $N_{\ka}$ with $\sP<\sP'$. If $\tilde\sO(\Theta,\sP,\sP',\delta)$ is the
open subset defined in Lemma \ref{lem:OverlappingTrivial}
of the space $\tilde\nu(\Theta,\sP,\sP',\delta)$
defined in \eqref{eq:R4OverlapSpace}, then the following diagram
commutes:
\begin{equation}
\label{eq:CommutingSplicing}
\begin{CD}
\tilde\sO(\Theta,\sP,\sP',\delta)
@> \rho^{\Theta,u}_{\sP,\sP'} >>
\tilde\sO(\Theta,\sP',\delta)
\\
@V \rho^{\Theta,d}_{\sP,\sP'} VV @V \bga'_{\Theta,\sP'} VV
\\
\tilde\sO(\Theta,\sP,\delta)
@> \bga'_{\Theta,\sP} >>
\bar{\sB}^s_{\ka}(S^4,2\delta)
\end{CD}
\end{equation}
\end{lem}

\begin{proof}
To clarify the notation, we write $Q$ for
an element of the partition $\sP'$. Denote
$$
\bA
:=
\left( (y_P)_{P\in\sP},
\left((x_Q)_{Q\in\sP'_P},[A_Q,F^s_Q]_{Q\in\sP'_P}\right)_{P\in\sP}\right)
\in \tilde\sO(\Theta,\sP,\sP',\delta),
$$
where  $(y_P)_{P\in\sP}\in \Delta^\circ(Z_\ka,\sP)$ and, by \eqref{eq:OverlapCondition3},
$((x_Q)_{Q\in\sP'_P},[A_Q,F^s_Q]_{Q\in\sP'_P})\in
D(\Theta,\sP'_P,\delta_P)$.
(For simplicity of exposition, we assume that
$[A_Q,F^s_Q]\in\bar\sB^s_{|Q|}(S^4,\delta_Q)$ has no ideal points;
the proof is no more difficult without that assumption, but
the notation becomes more opaque.)

By \eqref{eq:R4TrivialSplicingUpwardsOverlapMap},
\begin{align*}
\rho^{\Theta,u}_{\sP,\sP'}\left( \bA\right)
{}&=
\rho^{\Theta,u}_{\sP,\sP'}
\left( (y_P)_{P\in\sP},
\left((x_Q)_{Q\in\sP'_P},[A_Q,F^s_Q]_{Q\in\sP'_P}\right)_{P\in\sP}\right)
\\\
{}&=
\left( (x_Q')_{Q\in\sP'}, \left([A_Q,F^s_Q]_{Q\in\sP'_P}\right)\right),
\end{align*}
where
\begin{equation}
\label{eq:DefinePsi}
x'_Q=y_P+x_Q \quad \text{if $Q\in\sP'_P$.}
\end{equation}
(Note that
$x'_Q=y_{P}$ if $Q\in\sP\cap\sP'$ as $x_Q$ is then the zero vector.)
If we denote $\la_Q:=\la([A_Q])$, then
\begin{multline}
\label{eq:SplicedConn0}
\bga'_{\Theta,\sP'}\circ\rho^{\Theta,u}_{\sP,\sP'}(\bA)
\\
=
\begin{cases}
\Theta & \text{on $\RR^4 \less \cup_{P\in\sP,Q\in\sP'_P} B(x_Q',4\sqrt{\la_Q})$},
\\
(1-\chi_{x'_Q,4\sqrt{\la_Q}})c_{x'_Q,1}^*A_Q
\\
+
\chi_{x'_Q,4\sqrt{\la_Q}}c_{x'_Q,1}^*\Theta(A_Q,F^s_Q) & \text{on
$\Omega(x'_Q;2\sqrt{\la_Q},4\sqrt{\la_Q})$},
\\
c_{x'_Q,1}^*A_Q & \text{on $B(x_Q',2\sqrt{\la_Q})$}.
\end{cases}
\end{multline}
We now compare \eqref{eq:SplicedConn0} with
$\bga'_{\Theta,\sP}\circ\rho^{\Theta,d}_{\sP,\sP'}(\bA)$.

For $P\in\sP\less\sP'$, we denote
$$
[A_P,F^s_P]
:=
\bga'_{\Theta,\sP'_P}(\pi_P(\bA))
=
\bga'_{\Theta,\sP'_P} \left(
(x_Q)_{Q\in\sP'_P},[A_Q,F^s_Q]_{Q\in\sP'_P} \right),
$$
and
$[A_P,F^s_P]:=[A_Q,F^s_Q]$ for $P=Q\in\sP\cap\sP'$. Then by \eqref{eq:DefineUpwardTransition},
$$
\rho^{\Theta,d}_{\sP,\sP'}(\bA)
=
\left(
(y_P)_{P\in\sP},([A_P,F^s_P])_{P\in\sP}
\right)
\in \sO(\Theta,\sP,\delta).
$$
If we denote $\la_P':=\la([A_P])$, then
\begin{equation}
\label{eq:SplicedConn1}
\bga'_{\Theta,\sP}\circ\rho^{\Theta,d}_{\sP,\sP'}\left(\bA\right)
=
\begin{cases}
\Theta & \text{on $\RR^4\less \cup_{P\in\sP} B(y_P,4 \sqrt{\la_P'})$},
\\
(1-\chi_{y_P,4 \sqrt{\la_P'}}) c_{y_P,1}^*A_P
\\
+ \chi_{y_P,4 \sqrt{\la_P'}}c_{y_P,1}^*\Theta(A_P,F^s_P)) & \text{on $\Omega(y_P:2 \sqrt{\la_P'},4 \sqrt{\la_P'})$},
\\
c_{y_P,1}^*A_{P} & \text{on $B(y_P,4 \sqrt{\la_P'})$}.
\end{cases}
\end{equation}
Over the balls $B(x_P,4 \sqrt{\la_P'})$, for
$P\in\sP\cap\sP'$, and over $\RR^4\less\cup_{P\in\sP}B(x_P,4 \sqrt{\la_P'})$,
the connections \eqref{eq:SplicedConn0} and \eqref{eq:SplicedConn1} are
identical.  We thus focus our attention on the balls
$B(x_P,4 \sqrt{\la_P'})$ for $P\in\sP\less \sP'$.  For such a $P$, denoting
$\la_Q:=\la([A_Q])$ as done prior to \eqref{eq:SplicedConn0},  we have
\begin{align}
\label{eq:SplicedConn2}
A_P
{}&=\bga'_{\Theta,\sP'_P}(\pi_P(\bA))
\\
\notag &=
\begin{cases}
\Theta & \text{on $\RR^4\less \cup_{Q\in\sP'_P} B(x_Q,4\sqrt{\la_Q})$},
\\
(1-\chi_{x_Q,4 \sqrt{\la_Q}})c_{x_Q,1}^*A_Q
\\
+ \chi_{x_Q,4 \sqrt{\la_Q}}c_{x_Q,1}^*\Theta(A_Q,F^s_Q)
{}& \text{on $\Omega(x_Q,2\sqrt{\la_Q},4\sqrt{\la_Q})$},
\\
c_{x_Q,1}^*A_Q & \text{on $B(x_Q,2\sqrt{\la_Q})$}.
\end{cases}
\end{align}
Because the points $(x_Q)_{Q\in\sP'_P}\in
\Delta^{\circ}(Z_P(\delta),\sP'_P)$ and the connections
$[A_Q,F^s_Q]\in\bar{\sB}^s_{|Q|}(S^4)$ form a point in
$D(\Theta,\sP'_P,\delta_P)$ by the condition \eqref{eq:OverlapCondition3}
in the construction
of $\tilde\sO(\Theta,\sP,\sP',\delta)$,
Lemma \ref{lem:ScaleLimit} and equation \eqref{eq:BallContainment}
imply that for all $Q\in\sP'_P$ we have
$$
B\left(x_Q,4\sqrt{\la_Q}\right) \subset B\left(0,2\sqrt{\la_P'}\right),
$$
where $\la_P'$ is defined prior to \eqref{eq:SplicedConn1}.
Therefore \eqref{eq:SplicedConn2} implies that $A_P$
is already equal to the
product connection $\Theta$ over
$$
\RR^4\less B\left(0,2 \sqrt{\la_P'}\right),
$$
so $c_{y_P,1}^*A_P$ is equal to $\Theta$ over
$$
\RR^4\less B\left(y_P,2 \sqrt{\la_P'}\right).
$$
Hence, the convex combination in the second line of
\eqref{eq:SplicedConn1} is equal to $\Theta$.
Denoting $x'_Q=y_P+x_Q$ as in \eqref{eq:DefinePsi} and
$\la_Q:=\la([A_Q])$ as before \eqref{eq:SplicedConn2},
we examine the
final line of \eqref{eq:SplicedConn1} by applying Lemma
\ref{lem:ComposeCenterings} to rewrite
\eqref{eq:SplicedConn2} as
\begin{equation}
\label{eq:SplicedConn3}
c_{y_P,1}^*A_P
=
\begin{cases}
\Theta & \text{on $\RR^4 \less \cup_{Q\in\sP'_P}=B(x_Q',4\sqrt{\la_Q})$},
\\
(1-\chi_{x_Q',4 \sqrt{\la_Q}})c_{x_Q',1}^*A_Q
\\
+ \chi_{x_Q',4 \sqrt{\la_Q}}c_{x_Q',1}^*\Theta(A_Q,F^s_Q) &
\text{on $\Omega(x'_Q;2\sqrt{\la_Q},4\sqrt{\la_Q})$},
\\
c_{x'_Q,1}^*A_Q & \text{on $B(x_Q',2\sqrt{\la_Q})$}.
\end{cases}
\end{equation}
Consequently, we see that $\bga'_{\Theta,\sP}\circ\rho^{\Theta,d}_{\sP,\sP'}(\bA)$
is given by
\begin{equation}
\label{eq:SplicedConn4}
\begin{cases}
\Theta & \text{on $\RR^4 \less \cup_{P\in\sP,Q\in\sP'_P} B(x_Q',4\sqrt{\la_Q})$},
\\
(1-\chi_{x_Q',4 \sqrt{\la_Q}})c_{x_Q',1}^*A_Q
\\
+ \chi_{x_Q',4 \sqrt{\la_Q}}c_{x_Q',1}^*\Theta(A_Q,F^s_Q)  &
\text{on $\Omega(x'_Q;2\sqrt{\la_Q},4\sqrt{\la_Q})$},
\\
c_{x'_Q,1}^*A_Q & \text{on $B(x_Q',2\sqrt{\la_Q})$}.
\end{cases}
\end{equation}
The conclusion of Lemma \ref{lem:CommutingSplicing} now follows by comparing \eqref{eq:SplicedConn4} with
\eqref{eq:SplicedConn0}.
\end{proof}

\subsection{Symmetric group actions and quotients}
\label{subsubsec:SymGrpQuots}
We now make some observations on the action of the symmetric
group necessary to define
a quotient of the diagram \eqref{eq:CommutingSplice0} by the
symmetric group.
For $\sP<\sP'$, we recall from Lemma \ref{lem:SymmProdNeigh}
that to describe a neighborhood of
$\Si(Z_\ka(\delta),\sP)$ in $\Si(Z_\ka,\sP)$,
we must describe all the diagonals
$\Delta^\circ(Z_\ka(\delta),\sP'')$ where $\sP''\in [\sP<\sP']$.
We define
\begin{equation}
\label{eq:SplicingOverlapSpaceR4SymmQuotient}
\begin{aligned}
\tilde\nu(\Theta,\sP,[\sP'],\delta)
{}&:=
\bigsqcup_{\sP''\in [\sP<\sP']}\tilde\nu(\Theta,\sP,\sP'',\delta),
\\
\nu(\Theta,\sP,[\sP'],\delta)
{}&:=\tilde\nu(\Theta,\sP,[\sP'],\delta)/\fS(\sP),
\end{aligned}
\end{equation}
where the group $\fS(\sP)$ acts on $\tilde\nu(\Theta,\sP,[\sP'],\delta)$ by
the action defined in \eqref{eq:SymmGrpActionOnR4OverlapSpace}.
If the open subspaces $\tilde\sO(\Theta,\sP,\sP',\delta)$
defined in Lemma \ref{lem:OverlappingTrivial} are sufficiently
small, then their union
$\tilde\sO(\Theta,\sP,[\sP'],\delta)\subset\tilde\nu(\Theta,\sP,[\sP'],\delta)$
will be closed under the action of $\fS(\sP)$ and we define
\begin{equation}
\label{eq:SplicingOverlapSpaceR4SymmQuotientDiskBundles}
\sO(\Theta,\sP,[\sP'],\delta):=\tilde \sO(\Theta,\sP,[\sP'],\delta)/\fS(\sP).
\end{equation}
By analogy with the definition of
$T(\Theta,\sP,\sP',\delta)$ in \eqref{eq:DefineTrivialOverlap}, we define
\begin{equation}
\begin{aligned}
\label{eq:TrivialStrataOverlapSymmQuotient}
\tilde T(\Theta,\sP,[\sP'],\delta)
&{}:=
\left\{
\left(
    (y_P)_{P\in\sP},
    (x_Q)_{Q\in\sP''},
    ([\Theta,c_{|Q|}])_{Q\in\sP''}
\right) \in \tilde\nu(\Theta,\sP,[\sP'],\delta): \right.
\\
{}&\qquad\qquad
\left. \left(
    (y_P)_{P\in\sP},
    (x_Q)_{Q\in\sP''}
\right)\in \sO(Z_\ka(\delta),[\sP<\sP'])
\right\},
\\
T(\Theta,\sP,[\sP'],\delta)
{}&:=
\tilde T(\Theta,\sP,[\sP'],\delta)/\fS(\sP).
\end{aligned}
\end{equation}
where $\sO(Z_\ka(\delta),[\sP<\sP'])$ is defined in
Lemma \ref{lem:SymmProdNeigh}.

Recall that we denote the disjoint union or coproduct of sets $A_\alpha$ by $\sqcup_\alpha A_\alpha$
and for maps $f_\alpha:A_\alpha\to B_\alpha$, the coproduct map
$\coprod f_\alpha:\sqcup_\alpha A_\alpha \to \sqcup B_\alpha$ is defined by $f_\alpha$ on the subset
$A_\alpha\subset\sqcup_\alpha A_\alpha$.
The coproduct of the maps $\rho^{\Theta,u}_{\sP,\sP''}$,
$$
\coprod_{\sP''\in [\sP<\sP']} \rho^{\Theta,u}_{\sP,\sP''}:
\bigsqcup_{\sP''\in [\sP<\sP']}\tilde\sO(\Theta,\sP,\sP',\delta)
\to
\bigsqcup_{\sP''\in [\sP<\sP']} \nu(\Theta,\sP'',\delta)
$$
is equivariant with respect to the $\fS(\sP)$ actions and thus defines a map on $\fS(\sP)$ quotients,
\begin{equation}
\label{eq:SymmetricProductNormalBundlesOfHigherStratum}
\begin{aligned}
\rho^{\Theta,u}_{\sP,[\sP']}:&
\sO(\Theta,\sP,[\sP'],\delta)
\to
\nu(\Theta,[\sP<\sP'],\delta)
\\
{}&\nu(\Theta,[\sP<\sP'],\delta)
:=
\left.\left(\bigsqcup_{\sP''\in [\sP<\sP']} \nu(\Theta,\sP'',\delta)\right)\right/\fS(\sP).
\end{aligned}
\end{equation}
Similarly,
the coproduct of the maps $\rho^{\Theta,d}_{\sP,\sP''}$ is equivariant with respect to the action of $\fS(\sP)$ and thus
defines a map of the $\fS(\sP)$ quotients,
\begin{equation}
\label{eq:SymmetricProductNormalBundlesOfLowerStratum}
\rho^{\Theta,d}_{\sP,[\sP']}:
\sO(\Theta,\sP,[\sP'],\delta)
\to
\sO(\Theta,\sP,\delta).
\end{equation}
The coproduct of the splicing maps $\bga'_{\Theta,\sP''}$,
$$
\coprod_{\sP''\in [\sP<\sP']}\bga'_{\Theta,\sP''}:
\bigsqcup_{\sP''\in [\sP<\sP']}\sO(\RR^4,\sP'',\delta)
\to
\bar{\sB}^s_{\ka}(S^4,2\delta)
$$
is invariant under the action of $\fS(\sP)$
and thus defines a map
on the $\fS(\sP)$ quotient
\begin{equation}
\label{eq:SplicingOverSeveralPartitions}
\bga'_{\Theta,\sP,[\sP']} :
\left.\left(\bigsqcup_{\sP''\in [\sP<\sP']}\sO(\RR^4,\sP'',\delta)\right)\right/\fS(\sP)
\to
\bar{\sB}^s_{\ka}(S^4,2\delta),
\end{equation}
which has the same image as $\bga'_{\Theta,\sP'}$.
Lemma \ref{lem:CommutingSplicing} and the preceding discussion of symmetric group actions
then yields the

\begin{prop}
\label{prop:CommutingSplicingSymmQuotient}
Let $\sP$ and $\sP'$ be partitions
of $N_{\ka}$ with $\sP<\sP'$.
If $\sO(\Theta,\sP,[\sP'],\delta)$ is the open subspace defined
in \eqref{eq:SplicingOverlapSpaceR4SymmQuotientDiskBundles},
then one obtains the commutative diagram,
\begin{equation}
\label{eq:CommutingSplicingQuotient}
\begin{CD}
\sO(\Theta,\sP,[\sP'],\delta) @> \rho^{\Theta,u}_{\sP,[\sP']} >>
\left(\bigsqcup_{\sP''\in [\sP<\sP']}\sO(\RR^4,\sP'')\right)/\fS(\sP)
\\
@V \rho^{\Theta,d}_{\sP,[\sP']} VV @V \bga'_{\Theta,\sP,[\sP']} VV
\\
\sO(\Theta,\sP,\delta) @> \bga'_{\Theta,\sP} >> \bar{\sB}^s_{\ka}(S^4,2\delta)
\end{CD}
\end{equation}
where
$\rho^{\Theta,d}_{\sP,[\sP']}$ is defined in \eqref{eq:SymmetricProductNormalBundlesOfLowerStratum},
$\bga'_{\Theta,\sP,[\sP']}$ in \eqref{eq:SplicingOverSeveralPartitions},
and $\rho^{\Theta,u}_{\sP,[\sP']}$ in \eqref{eq:SymmetricProductNormalBundlesOfHigherStratum}.
\end{prop}

\section{The spliced end of the instanton moduli space}
\label{subsec:SplicedEnd}
We now construct the deformation of an Uhlenbeck neighborhood of the strata
\eqref{eq:PuncturedTrivialStrata} in $\bar M^{s,\natural}_\ka(S^4,\delta)$ by taking the union of the
images of the splicing maps $\bga'_{\Theta,\sP}$ restricted to appropriate subspaces of $\sO(\Theta,\sP,\delta)$.

This construction is inductive.  We therefore assume that
the instanton moduli space with spliced ends, $\barM^{s,\natural}_{\spl,|P|}(\delta)$,
with the properties stated in Theorem \ref{thm:ExistenceOfSplicedEndsModuli}
has already been constructed for $|P|<\ka$.
As described in Section \ref{sec:SplicedEnd_introduction}, the space $\barM^{s,\natural}_{\spl,\ka}(\delta)$
differs from $\barM^{s,\natural}_{\ka}(\delta)$ only on Uhlenbeck neighborhoods of the
strata containing the product connection described in \eqref{eq:PuncturedTrivialStrata}.
Such strata are determined by partitions $\sP=\{P_1,\dots,P_r\}$ of $N_\ka$ with $r>1$
so $|P_i|<\ka$ for all $i$.
Hence, the Uhlenbeck neighborhoods of these strata
are parameterized by the gluing maps \eqref{eq:ProductStratumGluingOnS4}
whose domain includes moduli spaces $\barM^{s,\natural}_{|P|}(\delta)$ where $|P|<\ka$.
The inductive hypothesis will  thus allow us to replace these moduli spaces with
$\barM^{s,\natural}_{\spl,|P|}(\delta)$.
There are no partitions of $\ka=1$ with length greater than one,
so defining
$\barM^{s,\natural}_{\spl,1}(\delta)=\barM^{s,\natural}_{1}(\delta)$
will complete the initial step of the induction.
Therefore, the
inductive hypothesis suffices to complete this construction.

The splicing maps will be restricted to the following spaces
(where the precise choice of the constants $\delta_P>0$ will not be relevant
to the rest of the discussion),
\begin{equation}
\label{eq:DefineASDSplicingDomain}
\begin{aligned}
\tilde\sO^{\asd}(\Theta,\sP,\delta)
&:=
\tilde\sO(\Theta,\sP,\delta)\cap \left(
    \Delta^\circ(Z_\ka(2\delta),\sP)
    \times
    \prod_{P\in\sP} M^{s,\natural}_{\spl,|P|}(\delta_P)
\right),
\\
\sO^{\asd}(\Theta,\sP,\delta)&:=\tilde\sO^{\asd}(\Theta,\sP,\delta)/\fS(\sP).
\end{aligned}
\end{equation}
Observe that $\tilde T(\Theta,\sP,\delta)\subset \tilde\sO^{\asd}(\Theta,\sP,\delta)$,
where $\tilde T(\Theta,\sP,\delta)$
is defined in \eqref{eq:DefineTrivialStrataInR4GluingSpace} and
that the dimension of $\sO^{\asd}(\Theta,\sP,\delta)$ is equal to that of
$\barM^{s,\natural}_{\ka}(S^4,\delta)$.

The spliced end of the instanton moduli space
is constructed in the following proposition.

\begin{prop}
\label{prop:ExistenceOfSplicedEnd}
Assume that the instanton moduli space with spliced ends, $\barM^{s,\natural}_{\spl,\ka'}(\delta)$,
with the properties enumerated in Theorem \ref{thm:ExistenceOfSplicedEndsModuli} has already been constructed
for all $\ka'<\ka$.
Then for every partition $\sP$ of $N_\ka$ with length greater than one,
there is an $\SO(3)\times\SO(4)$-invariant
open neighborhood, $\sO^{\asd}_1(\Theta,\sP,\delta)$, of
$T(\Theta,\sP,\delta)$ in $\sO^{\asd}(\Theta,\sP,\delta)$ such that
\begin{equation}
\label{eq:DefineSplicedEnd}
W_\ka := \bigcup_{\sP} \ \bga'_{\Theta,\sP}\left(\sO^{\asd}_1(\Theta,\sP,\delta)\right)
\end{equation}
is a smoothly-stratified subspace of $\bar\sB^s_{\ka}(S^4,\delta)$,
with every point $[A,F^s,\bx]\in W_\ka$ satisfying the estimate
$\|F^+_{A}\|_{L^{\sharp,2}(S^4)}\le \eps$,
where $\|\cdot\|_{L^{\sharp,2}}$ is the norm defined in \eqref{eq:LSharpNorm}
and $\eps$ is the constant appearing in \cite[Proposition 7.6]{FLKM1}.
\end{prop}

We will prove Proposition \ref{prop:ExistenceOfSplicedEnd} by using
Proposition \ref{prop:CommutingSplicingSymmQuotient} to
describe the intersections of the images of the splicing maps.
The transition maps $\rho^{\Theta,u}_{\sP,[\sP']}$ and
$\rho^{\Theta,d}_{\sP,[\sP']}$ will be restricted to the
following subspace of $\sO(\Theta,\sP,[\sP'],\delta)$,
\begin{subequations}
\label{eq:SplicedModuliOverlapR4}
\begin{multline}
\label{eq:SplicedModuliOverlapR4_prequotient}
\tilde\sO^{\asd} (\Theta,\sP,[\sP'],\delta)
:=
\tilde\sO(\Theta,\sP,[\sP'],\delta)
\cap \left(
    \Delta^\circ(Z_\ka(\delta),\sP)
    \right.
\\
    \left.
    \times
    \bigsqcup_{\sP''\in [\sP<\sP']}
    \prod_{P\in\sP}
            \left(
            \Delta^\circ(Z_P(\delta_P),\sP_P'')
            \times
            \prod_{Q\in\sP_P''} \barM^{s,\natural}_{\spl,|Q|}(\delta_Q)
            \right)
\right),
\end{multline}
\begin{equation}
\label{eq:SplicedModuliOverlapR4_quotient}
\sO^{\asd}(\Theta,\sP,[\sP'],\delta) := \tilde\sO^{\asd} (\Theta,\sP,[\sP'],\delta)/\fS(\sP).
\end{equation}
\end{subequations}
(Compare the definition
\eqref{eq:SplicingOverlapSpaceR4SymmQuotient}.)
Observe that $T(\Theta,\sP,[\sP'],\delta)\subset
\sO^{\asd} (\Theta,\sP,[\sP'],\delta)$,
where $T(\Theta,\sP,[\sP'],\delta)$ is well-defined in
\eqref{eq:TrivialStrataOverlapSymmQuotient} because
$[\Theta,c_{|Q|}]\in \barM^{s,\natural}_{\spl,|Q|}(\delta_Q)$ for
all $|Q|<\ka$.

We must first verify that
the transition maps $\rho^{\Theta,u}_{\sP,[\sP']}$ and
$\rho^{\Theta,d}_{\sP,[\sP']}$ map
$\sO^{\asd} (\Theta,\sP,[\sP'],\delta)$ to
$\sO^{\asd}(\Theta,\sP',\delta)$ and $\sO^{\asd}(\Theta,\sP,\delta)$ respectively.

\begin{lem}
\label{lem:SplicedEndsEmbeddingR4d}
Assume that $\barM^{s,\natural}_{\spl,\ka'}(\delta)$
has been constructed so that it satisfies the properties enumerated in
Theorem \ref{thm:ExistenceOfSplicedEndsModuli} for all
$\ka'<\ka$.  Let $T(\Theta,\sP,[\sP'],\delta)\subset \sO^{\asd} (\Theta,\sP,[\sP'],\delta)$
be as defined in
\eqref{eq:TrivialStrataOverlapSymmQuotient}.
Then there is an open neighborhood
$\sO^{\asd}_1(\Theta,\sP,[\sP'],\delta)$
of $T(\Theta,\sP,[\sP'],\delta)$ in $\sO^{\asd} (\Theta,\sP,[\sP'],\delta)$
such that  the restriction of $\rho^{\Theta,d}_{\sP,[\sP']}$
to $\sO^{\asd}_1(\Theta,\sP,[\sP'],\delta)$ is an open embedding
of $\sO^{\asd}_1(\Theta,\sP,[\sP'],\delta)$ into
$\sO^{\asd}(\Theta,\sP,\delta)$.
\end{lem}

\begin{proof}
Induction and the property \eqref{eq:EqualsSplicingImage} in
Theorem \ref{thm:ExistenceOfSplicedEndsModuli} imply
that for each $P\in\sP\less\sP'$, there is an open
neighborhood $\sO^{\asd}_1(\Theta,\sP'_P,\delta_P)$ of
$T(\Theta,\sP'_P,\delta)$ in
$\sO^{\asd}(\Theta,\sP'_P,\delta)$ such that
$$
\bga'_{\Theta,\sP'_P}\left( \sO^{\asd}_1(\Theta,\sP'_P,\delta_P)\right) \subset
\barM^{s,\natural}_{\spl,|P|}(\delta_P).
$$
Therefore, in the definition \eqref{eq:SplicedModuliOverlapR4}of $\sO^{\asd} (\Theta,\sP,[\sP'],\delta)$, if we replace the factor
$$
\Delta^\circ(Z_P(\delta_P),\sP_P'')
            \times
            \prod_{Q\in\sP_P''} \barM^{s,\natural}_{\spl,|Q|}(\delta_Q)
$$
with the open subspace,
$$
\tilde\sO^{\asd}_1(\Theta,\sP'_P,\delta_P)
 \subset
\Delta^\circ(Z_P(\delta_P),\sP_P'')
            \times
            \prod_{Q\in\sP_P''} \barM^{s,\natural}_{\spl,|Q|}(\delta_Q),
$$
then the resulting open subspace $\sO^{\asd}_1(\Theta,\sP,[\sP'],\delta)$
of $\sO^{\asd} (\Theta,\sP,[\sP'],\delta)$ will satisfy
\begin{equation}
\label{eq:RestrictedUpwardsTransitionInclusion}
\rho^{\Theta,d}_{\sP,[\sP']}
\left(\sO^{\asd}_1(\Theta,\sP,[\sP'],\delta)\right)
\subset
\sO^{\asd}(\Theta,\sP,\delta).
\end{equation}
(Note that we are not concerned with the
subsets $P\in\sP\cap\sP'$
because $\rho^{\Theta,d}_{\sP,\sP'}$ is the identity on these factors
and so they do not affect the inclusion
\eqref{eq:RestrictedUpwardsTransitionInclusion}.)

The fact that the restriction of $\rho^{\Theta,d}_{\sP,[\sP']}$
to $\sO^{\asd}_1(\Theta,\sP,[\sP'],\delta)$ gives
an open embedding into $\sO^{\asd}(\Theta,\sP,\delta)$
follows from the assumption that
$\bga'_{\Theta,\sP'_P}$ gives an open embedding of
$\sO^{\asd}_1(\Theta,\sP'_P,\delta_P)$ into
$\barM^{s,\natural}_{\spl,|P|}(\delta_P)$.
\end{proof}

We now prove an analogue of Lemma \ref{lem:SplicedEndsEmbeddingR4d}
for the map $\rho^{\Theta,u}_{\sP,[\sP']}$.
First, we must define the appropriate range of $\rho^{\Theta,u}_{\sP,[\sP']}$,
to take into account the appearance of the conjugate partitions
$\sP''\in [\sP<\sP']$.
Let
\begin{multline}
\label{eq:UpperPartitionR4ASDNormal}
\sO^{\asd}(\Theta,[\sP<\sP'],\delta)
:=
\left.\left(\bigsqcup_{\sP''\in [\sP<\sP']}\sO^{\asd}(\Theta,\sP'',\delta)
\right)
\right/\fS(\sP)
\\
\subset
\nu(\Theta,[\sP<\sP'],\delta),
\end{multline}
where $\nu(\Theta,[\sP<\sP'],\delta)$
is defined in \eqref{eq:SymmetricProductNormalBundlesOfHigherStratum}.
We then have the

\begin{lem}
\label{lem:SplicedEndsEmbeddingR4u}
Continue the assumptions and hypotheses of Lemma
\ref{lem:SplicedEndsEmbeddingR4d}.
Then the restriction of the map
$\rho^{\Theta,u}_{\sP,[\sP']}$
 to $\sO^{\asd} (\Theta,\sP,[\sP'],\delta)$
is an open embedding into $\sO^{\asd}(\Theta,[\sP<\sP'],\delta)$.
\end{lem}

\begin{proof}
The restriction of the map $\rho^{\Theta,u}_{\sP,\sP'}$ defines a map
$$
\begin{CD}
\tilde\sO(\Theta,\sP,\sP',\delta)
\cap
\left(
\prod_{P\in\sP}
\left(
\Delta^\circ(Z_P(\delta_P),\sP'_P)
\times
\prod_{Q\in\sP'_P}  \barM^{s,\natural}_{\spl,|Q|}(\delta_Q)
\right)
\right)
\\
@VVV
\\
\Delta^\circ(Z_P(\delta_P),\sP')
\times
\prod_{P\in\sP}\prod_{Q\in\sP'_P} \barM^{s,\natural}_{\spl,|Q|}(\delta_Q)
\end{CD}
$$
By the definition of $\rho^{\Theta,u}_{\sP,\sP'}$ in \eqref{eq:R4TrivialSplicingUpwardsOverlapMap},
this restriction is given by the product of the exponential map $e(Z_\ka,\sP)$ with the projection maps
$\pi_Q$ onto $\barM^{s,\natural}_{\spl,|Q|}(\delta_Q)$ defined in \eqref{eq:DefineProjections} and is thus an open embedding.
The conclusion of the lemma then follows immediately from the definition
of  $\rho^{\Theta,u}_{\sP,[\sP']}$ in terms of the maps $\rho^{\Theta,u}_{\sP,\sP'}$.
\end{proof}

The following lemma yields the existence of suitably small,
$\SO(3)\times\SO(4)$-invariant
neighborhoods of $T(\Theta,\sP,\delta)$.
For any subspace $V$ of $\bar{\sB}^s_{\ka}(S^4)$, let
$\cl(V)$ denote the closure of $V$ in the Uhlenbeck topology.

\begin{lem}
\label{lem:SeparatingNgh}
If $V$ is any  set in
$\bar{\sB}^s_{\ka}(S^4)$ satisfying
$$
\cl(V) \cap \bga'_{\Theta,\sP}\left(
T(\Theta,\sP,\delta)\right)=\emptyset,
$$
then there are  $\SO(3)\times\SO(4)\times \fS(\sP)$-invariant
neighborhoods,
$$
\sO^{\asd}_{f,3}(\Theta,\sP,\delta)\sqsubset
\sO^{\asd}_{f,2}(\Theta,\sP,\delta))\sqsubset
\sO^{\asd}_{f,1}(\Theta,\sP,\delta),
$$
of $T(\Theta,\sP,\delta)$ in
$\sO^{\asd}(\Theta,\sP,\delta)$ such that the intersection,
$$
\cl(V) \cap \bga'_{\Theta,\sP}\left( \sO^{\asd}_{f,1}(\Theta,\sP,\delta)\right),
$$
is empty and there is an inclusion,
$$
\cl(V) \cap \cl
  \left(
    \bga'_{\Theta,\sP}\left( \sO^{\asd}_{f,1}(\Theta,\sP,\delta)\right)
  \right)
\subset \cl(V) \cap \cl
  \left(
     \bga'_{\Theta,\sP}\left( T(\Theta,\sP,\delta)\right)
  \right).
$$
\end{lem}

\begin{proof}
For $\bA\in \sO^{\asd}(\Theta,\sP,\delta)$ given by
$\bA=\left((x_P)_{P\in\sP},([A_P,F^s_P,\by_P])_{P\in\sP}\right)$,
denote $\la_P(\bA):=\la(A_P,\by_P)$. For
$\bx:=(x_P)_{P\in\sP}\in \Delta^\circ(Z_\ka,\sP)$, define
$$
V(\bx,\sP)
:=
\{\left(\bx,([A_P,F^s_P,\by_P])_{P\in\sP}\right)
  \in (\bga'_{\Theta,\sP})^{-1}\left(\cl( V) \right)\}
$$
and define $\bla_{V}:\Delta^\circ(Z_\ka,\sP)\to (0,1]$ by
$$
\bla_{V}(\bx)^2
:=
\min\left\{\min_{\bA\in
V(\bx,\sP)}\sum_{P\in\sP}\la_P(\bA)^2, 1\right\}.
$$
By the hypothesis that $\cl(V)$ is disjoint from
$\bga'_{\Theta,\sP}\left( T(\Theta,\sP,\delta)\right)$
and the observation that $T(\Theta,\sP,\delta)$ is the zero locus of $\sum_{P\in\sP}\la_P(\cdot)$, we see that
$\bla_{V}$ is
always positive on $\Delta^\circ(Z_\ka(\delta),\sP)$.
Let $f:\Delta^\circ(Z_\ka(\delta),\sP)\to (0,\eps)$ be a continuous
function with
\begin{equation}
\label{eq:fBoundedByScaleDistance}
f((x_P)_{P\in\sP}) < \half \bla_{V}((x_P)_{P\in\sP})^2,
\end{equation}
for all $(x_P)_{P\in\sP}\in \Delta^\circ(Z_\ka(\delta),\sP)$.
Then for $j=1,2,3$, we define the neighborhood $\sO^{\asd}_{f,j}(\Theta,\sP,\delta)$ by
$$
\sO^{\asd}_{f,j}(\Theta,\sP,\delta)
:=
\left\{\bA=\left
(\bx,([A_P,F^s_P,\by_P])_{P\in\sP}\right)\in\sO^{\asd}(\Theta,\sP,\delta):
  \sum_{P\in\sP} \la_P(\bA)^2< f(\bx)/j\right\}.
$$
These sets have the desired invariance and their images do not
intersect $\cl(V)$.
The existence of a smoothly-stratified function $g_j$ supported in
${\mathscr{O}}^{\textsc{asd}}_{f,j}(\Theta,\mathscr{P},\delta)$
with
$g_j({\mathscr{O}}^{\textsc{asd}}_{f,j+1}(\Theta,\mathscr{P},\delta))=1$
(needed to prove ${\mathscr{O}}^{\textsc{asd}}_{f,j+1}(\Theta,\mathscr{P},\delta)\sqsubset {\mathscr{O}}^{\textsc{asd}}_{f,j}(\Theta,\mathscr{P},\delta)$) follows from the
definition of these sets.

To prove the final inclusion, let
$\{\bA(\alpha)\}_{\alpha=1}^{\infty}\subset\sO^{\asd}_{f,1}(\Theta,\sP,\delta)$ be a sequence
with
\begin{equation}
\label{eq:LimitInV}
\lim_{\alpha\to\infty} \bga'_{\Theta,\sP}\left( \bA(\alpha)\right)
\in \cl(V).
\end{equation}
If
$$
\bA(\alpha)=[ \bx(\alpha), ([A_P(\alpha),F^s_P(\alpha),\by_P(\alpha)])_{P\in\sP}],
$$
where $\bx(\alpha)\in \Delta^\circ(Z_\ka(\delta),\sP)$, then \eqref{eq:LimitInV} implies
that
\begin{equation}
\label{eq:LimitOfNgh1}
\lim_{\alpha\to\infty} \sum_{P\in\sP}\la_P(\bA(\alpha))^2
=
\lim_{\alpha\to\infty}\bla_V(\bx(\alpha))^2.
\end{equation}
Because $\bA(\alpha)\in \sO^{\asd}_{f,1}(\Theta,\sP,\delta)$, we have
\begin{equation}
\label{eq:LimitOfNgh2}
\sum_{P\in\sP}\la_P(\bA(\alpha))^2<
\frac{1}{2}\bla_V(\bx(\alpha))^2.
\end{equation}
Combining \eqref{eq:LimitOfNgh1} and \eqref{eq:LimitOfNgh2} yields
$$
\lim_{\alpha\to\infty} \sum_{P\in\sP}\la_P^2(\bA(\alpha))=0,
$$
and so for all $P\in\sP$,
$$
\lim_{\alpha\to\infty} [A_P(\alpha),F^s_P(\alpha),\by_P(\alpha)]
=
[\Theta,F^s_P(\infty),c_P],
$$
and thus,
$$
\lim_{\alpha\to\infty}\bga'_{\Theta,\sP}\left( \bA(\alpha)\right)
=
\lim_{\alpha\to\infty}
\bga'_{\Theta,\sP}
\left( \bx(\alpha),\left( [\Theta,F^s_P(\infty),c_P]\right)_{P\in\sP}\right)
\in
\cl \left(\bga'_{\Theta,\sP}(T(\Theta,\sP,\delta))\right),
$$
as asserted.
\end{proof}

We now construct the subspaces $\sO^{\asd}_1(\Theta,\sP,\delta)$
referred to in Proposition \ref{prop:ExistenceOfSplicedEnd}.
We will construct these subspaces to be $\fS_\ka$-invariant
in the sense that if $\sP''\in [\sP]$, then the open
subspaces $\sO^{\asd}_1(\Theta,\sP,\delta)$ and $\sO^{\asd}_1(\Theta,\sP'',\delta)$
will be identified by the natural
action of $\fS_\ka$ on $\bigsqcup_{\sP_i\in[\sP_1]}\nu(\Theta,\sP_i)$.
Thus, defining $\sO^{\asd}_1(\Theta,\sP',\delta)$ for one $\sP'\in [\sP']$
suffices to define the space
\begin{equation}
\label{eq:DefineSymmetricQuotientOfR4UpperASDNgh}
\sO^{\asd}_1(\Theta,[\sP<\sP'])
=
\left.\left(
\bigsqcup_{\sP''\in [\sP<\sP']} \sO^{\asd}_1(\Theta,\sP'',\delta)
\right)\right/\fS(\sP)
\subset
\sO^{\asd}(\Theta,[\sP<\sP'],\delta),
\end{equation}
where $\sO^{\asd}(\Theta,[\sP<\sP'],\delta)$ is defined in
\eqref{eq:UpperPartitionR4ASDNormal}.

\begin{lem}
\label{lem:ConstructCollar}
Assume that $\barM^{s,\natural}_{\spl,\ka'}(\delta)$
has been constructed so that it satisfies the properties enumerated in
Theorem \ref{thm:ExistenceOfSplicedEndsModuli} for all
$\ka'<\ka$. Then for every partition $\sP$ of $N_\ka$
of length greater than one, there  are
open neighborhoods
$$
\sO^{\asd}_3(\Theta,\sP,\delta)
\sqsubset
\sO^{\asd}_2(\Theta,\sP,\delta)
\sqsubset
\sO^{\asd}_1(\Theta,\sP,\delta)
$$
of $T(\Theta,\sP,\delta)$ in $\sO^{\asd}(\Theta,\sP,\delta)$ such that the following hold for $j=1,2,3$.
\begin{enumerate}
\item
\label{item:ConstructCollar1}
The neighborhoods
$\sO^{\asd}_j(\Theta,\sP,\delta)$
are closed  under the
$\SO(3)\times\SO(4)$ action given in \eqref{eq:SO3SO4Actions}.
\item
\label{item:ConstructCollar2}
If the partitions $\sP_1$ and $\sP_2$ are conjugate under the action of
$\fS_\ka$, then the neighborhoods $\sO^{\asd}_j(\Theta,\sP_1,\delta)$
and $\sO^{\asd}_j(\Theta,\sP_2,\delta)$  are identified by the natural
action of $\fS_\ka$ on $\bigsqcup_{\sP_i\in[\sP_1]}\nu(\Theta,\sP_i)$.
\item
\label{item:ConstructCollar3}
If there are no partitions $\sP_1\in [\sP]$ and $\sP_2\in [\sP']$
with $\sP_1<\sP_2$ or $\sP_2<\sP_1$,  then
$$
\bga'_{\Theta,\sP}(\sO^{\asd}_1(\Theta,\sP,\delta))
  \cap
\bga'_{\Theta,\sP'}(\sO^{\asd}_1(\Theta,\sP',\delta))=\emptyset.
$$
\item
\label{item:ConstructCollar4}
If $\sP< \sP'$ then there exists an open neighborhood
$\sO^{\asd}_1(\Theta,\sP,[\sP'],\delta)$ of the subspace
$T(\Theta,\sP,[\sP'],\delta)$ of
$\sO^{\asd} (\Theta,\sP,[\sP'],\delta)$ such that
\begin{equation}
\label{eq:CompleteOverlapsR4}
\begin{aligned}
\bga'_{\Theta,\sP} & \left( \sO^{\asd}_1(\Theta,\sP,\delta)\right) \cap
\bga'_{\Theta,\sP'}\left( \sO^{\asd}_1(\Theta,\sP',\delta)\right)
\\
&\subseteq \bga'_{\Theta,\sP'}\left(
\rho^{\Theta,u}_{\sP,[\sP']}\left(\sO^{\asd}_1(\Theta,\sP,[\sP'],\delta)\right)\right)
\\
&=
\bga'_{\Theta,\sP}\left(
\rho^{\Theta,d}_{\sP,[\sP']}\left(\sO^{\asd}_1(\Theta,\sP,[\sP'],\delta)\right)\right).
\end{aligned}
\end{equation}
\end{enumerate}
\end{lem}

\begin{proof}
Enumerate the strata $\Si(Z_\ka(\delta),\sP)$ in the manner described in
Section \ref{subsec:EnumStrata} with conjugacy classes $[\sP_0],\dots,[\sP_n]$ so that
\begin{equation}
\label{eq:ClosureOfStrataInPartitionEnumeration}
\Si(Z_\ka(\delta),\sP_i)\subset\cl\left(\Si(Z_\ka(\delta),\sP_j)\right) \quad\text{only if}\quad i<j.
\end{equation}
Because the length of the partition $\sP_0$ is one, we begin our induction with $\sP_1$ and
using induction on $k$, construct  neighborhoods $\sO^{\asd}_j(\Theta,\sP_k,\delta)$ satisfying items
\eqref{item:ConstructCollar1} and \eqref{item:ConstructCollar2}.  In addition, these inductively constructed
neighborhoods will, for all $i<k$, satisfy Items \eqref{item:ConstructCollar3} and \eqref{item:ConstructCollar4}
with $[\sP_i]=[\sP]$ and $[\sP_k]=[\sP']$.  We note that the symmetry between $\sP$ and $\sP'$
in Items \eqref{item:ConstructCollar3} and \eqref{item:ConstructCollar4} implies that this
will suffice to complete the induction and prove the lemma.

In each step of the induction, to ensure this construction satisfies Item \eqref{item:ConstructCollar2}, it suffices to construct $\sO^{\asd}_j(\Theta,\sP_k,\delta)$ for
a single partition $\sP'\in [\sP_k]$ and for any other
partition,
$\sP''\in [\sP_k]$, define $\sO^{\asd}_j(\Theta,\sP'',\delta)$
to be the image of $\sO^{\asd}_j(\Theta,\sP',\delta)$ under the action of
$\fS_\ka$.

For $k=1$, take $\sO^{\asd}_1(\Theta,\sP_1,\delta)=\sO^{\asd}(\Theta,\sP_1,\delta)$.
The neighborhoods $\sO^{\asd}_j(\Theta,\sP_1,\delta)$ for $j=2,3$ can be constructed
by replacing the coefficient $8$ in the definition \eqref{eq:DefineSplicingDataR4} with
slightly larger coefficients.
Item \eqref{item:ConstructCollar2} holds from the definition of $\sO^{\asd}(\Theta,\sP_1,\delta)$
while Items \eqref{item:ConstructCollar3} and \eqref{item:ConstructCollar4} hold trivially.

By induction, assume that the open sets $\sO^{\asd}_1(\Theta,\sP_i,\delta)$ satisfying Items \eqref{item:ConstructCollar1}-
\eqref{item:ConstructCollar4} have been
defined for all $i<k$; we will complete the induction by constructing $\sO^{\asd}_1(\Theta,\sP_k,\delta)$.
Note that the conclusion still holds if we
shrink the sets $\sO^{\asd}_1(\Theta,\sP_i,\delta)$ to smaller neighborhoods of
$T(\Theta,\sP_i,\delta) $.
We now construct $\sO^{\asd}_1(\Theta,\sP_k,\delta)$.

For each $i<k$, if there are no
$\sP'\in [\sP_k]$ and $\sP\in [\sP_i]$ with $\sP<\sP'$, then
by \eqref{eq:SplicingOnTrivialStratum},
\eqref{eq:ClosureOfStrataInPartitionEnumeration}, and Lemma \ref{lem:IncidenceRelationInSymmetricProduct},
\begin{align}
\label{eq:UnrelatedTrivIntersection1}
\cl\left(\bga'_{\Theta,\sP_i}\left( T(\Theta,\sP_i,\delta)\right)\right)
\cap \bga'_{\Theta,\sP_k}\left( T(\Theta,\sP_k,\delta)\right) & =
\emptyset, \\
\label{eq:UnrelatedTrivIntersection2}
\bga'_{\Theta,\sP_i}\left(T(\Theta,\sP_i,\delta)\right)
\cap
\cl\left(\bga'_{\Theta,\sP_k}\left( T(\Theta,\sP_k,\delta)\right)\right)
{}&=\emptyset .
\end{align}
Lemma \ref{lem:SeparatingNgh} and
\eqref{eq:UnrelatedTrivIntersection1} imply that there  are $\SO(3)\times\SO(4)$-invariant
neighborhoods,
$$
\sO^{\asd}_3(\Theta,\sP_k,\delta)\sqsubset
\sO^{\asd}_2(\Theta,\sP_k,\delta)\sqsubset
\sO^{\asd}_1(\Theta,\sP_k,\delta),
$$
of $T(\Theta,\sP_k,\delta)$ such that
$$
\cl
 \left(  \bga'_{\Theta,\sP_i}\left( T(\Theta,\sP_i,\delta)\right)\right)
\cap \bga'_{\Theta,\sP_k}\left( \sO^{\asd}_1(\Theta,\sP_k,\delta)\right) =\emptyset.
$$
Lemma \ref{lem:SeparatingNgh} also ensures that
\begin{align}
\notag \cl {}& \left(
      \bga'_{\Theta,\sP_i}\left( T(\Theta,\sP_i,\delta)\right)
    \right)
\cap \cl \left(
      \bga'_{\Theta,\sP_k}\left( \sO^{\asd}_1(\Theta,\sP_k,\delta)\right)
    \right)
\\
\label{eq:UnrelatedClosureContainment}
{}&\subseteq \cl \left(
      \bga'_{\Theta,\sP_i}\left( T(\Theta,\sP_i,\delta)\right)
    \right)
\cap \cl \left(
      \bga'_{\Theta,\sP_k}\left( T(\Theta,\sP_k,\delta)\right)
    \right).
\end{align}
Equations \eqref{eq:UnrelatedClosureContainment} and
\eqref{eq:UnrelatedTrivIntersection2} yield
\begin{align*}
 \bga'_{\Theta,\sP_i}& \left( T(\Theta,\sP_i,\delta)\right)
\cap \cl \left(
      \bga'_{\Theta,\sP_k}\left( \sO^{\asd}_1(\Theta,\sP_k,\delta)\right)
    \right) \\
& \subset
 \bga'_{\Theta,\sP_i}  \left( T(\Theta,\sP_i,\delta)\right)
\cap \cl \left(
      \bga'_{\Theta,\sP_k}\left( T(\Theta,\sP_k,\delta)\right)
    \right)=\emptyset.
\end{align*}
Thus, by Lemma \ref{lem:SeparatingNgh} there are
(smaller) $\SO(3)\times\SO(4)$-invariant neighborhoods
$\sO^{\asd}_{j,2}(\Theta,\sP_i,\delta)\subset \sO^{\asd}_j(\Theta,\sP_i,\delta)$ of $T(\Theta,\sP_i,\delta)$
with $j=1,2,3$ satisfying
$$
\sO^{\asd}_{3,2}(\Theta,\sP_i,\delta)\sqsubset
\sO^{\asd}_{2,2}(\Theta,\sP_i,\delta)\sqsubset
\sO^{\asd}_{1,2}(\Theta,\sP_i,\delta)
$$
and
$$
 \bga'_{\Theta,\sP_i}\left( \sO^{\asd}_{1,2}(\Theta,\sP_i,\delta)\right)
\cap \cl \left(
      \bga'_{\Theta,\sP_k}\left( \sO^{\asd}_1(\Theta,\sP_k,\delta)\right)
    \right)=\emptyset,
$$
as desired.
Replace the open neighborhoods $\sO^{\asd}_j(\Theta,\sP_i,\delta)$
with $\sO^{\asd}_{j,2}(\Theta,\sP_i,\delta)$.
For all other $\sP''\in [\sP_i]$, replace $\sO^{\asd}_j(\Theta,\sP'',\delta)$
with the image of $\sO^{\asd}_{j,2}(\Theta,\sP_i,\delta)$  under the action of
$\fS_\ka$. Then $\sO^{\asd}_1(\Theta,\sP_k,\delta)$ satisfies
Items \eqref{item:ConstructCollar3} and \eqref{item:ConstructCollar4}
 for $\sP_i$ with $i<k$ and no
$\sP'\in [\sP_k]$ and $\sP\in [\sP_i]$ with $\sP<\sP'$.

If $[\sP_i]<\sP_k$, by conjugating $\sP_i$ if necessary and using the
invariance of the neighborhoods $\sO^{\asd}_j(\Theta,\sP_i,\delta)$ under the
action of $\fS_{\ka}$, we can assume that $\sP_i<\sP_k$.
By the equality \eqref{eq:SplicingOnTrivialStratum}, the closure of
$$
[\Theta] \times \left( \Si(Z_\ka(\delta),\sP_k) \setminus \Imag
(e(Z_\ka,\sP_i)) \right),
$$
(where $e(Z_\ka,\sP_i)$ is the exponential map defined in
equation \eqref{eq:NormalOfR4Diagonal})
does not intersect the image $\bga'_{\Theta,\sP_i}(T(\Theta,\sP_i,\delta))$.
Therefore, by applying Lemma \ref{lem:SeparatingNgh} and shrinking
$\sO^{\asd}_j(\Theta,\sP_i,\delta)$, for $j=1,2,3$, as above, we can assume that
$$
\cl \left(
\bga'_{\Theta,\sP_i}\left( \sO^{\asd}_1(\Theta,\sP_i,\delta)\right)
    \right)
\cap
\left(
[\Theta] \times \left( \Si(Z_\ka,\sP'') \setminus \Imag
(e(Z_\ka,\sP)) \right)
\right)
=\emptyset
$$
for all $\sP''\in [\sP_i<\sP_k]$
and hence
\begin{multline}
\label{eq:ClosureLowerOpenLessTrivOverlap}
\cl \left(
\bga'_{\Theta,\sP}\left( \sO^{\asd}_1(\Theta,\sP_i,\delta)\right)
    \right)
\cap
\left(
\left.\left(
\bigcup_{\sP''\in [\sP_i<\sP_k]}
    \bga'_{\Theta,\sP''}( T(\Theta,\sP'',\delta))
\right)\right/\fS(\sP)\right)
\\
\sqsubset [\Theta] \times
e(Z_\ka(\delta),\sP_i)\left(\sO(Z_\ka(\delta),[\sP_i<\sP_k])\right),
\end{multline}
where $\sO(Z_\ka(\delta),[\sP_i<\sP_k])$ is defined in
Lemma \ref{lem:SymmProdNeigh}.

Let the spaces $T(\Theta,\sP_i,[\sP_k],\delta)$
and $\sO^{\asd}_1(\Theta,\sP_i,[\sP_k],\delta)$
be as defined in \eqref{eq:TrivialStrataOverlapSymmQuotient}.
The equality and inclusion,
\begin{align*}
{}& [\Theta]\times e(Z_\ka,\sP_i)\left(\sO(Z_\ka(\delta),[\sP_i<\sP_k])\right)
\\
&\quad =
\bga'_{\Theta,[\sP_i<\sP_k]}
  \left(
    \rho^{\Theta,u}_{\sP_i,[\sP_k]}\left( T(\Theta,\sP_i,[\sP_k],\delta)\right)
  \right)
\\
&\quad \subset \bga'_{\Theta,[\sP_i<\sP_k]}
  \left(
    \rho^{\Theta,u}_{\sP_i,[\sP_k]}\left( \sO^{\asd}_1(\Theta,\sP_i,[\sP_k],\delta)\right)
  \right),
\end{align*}
together with \eqref{eq:ClosureLowerOpenLessTrivOverlap}
then imply that
\begin{equation}
\label{eq:LowerR4SplicingMinusOverlap}
\cl \left( \bga'_{\Theta,\sP_i}\left( \sO^{\asd}_1(\Theta,\sP_i,\delta)\right)
    \right)
\setminus
\bga'_{\Theta,[\sP_i<\sP_k]}
  \left(
    \rho^{\Theta,u}_{\sP_i,[\sP_k]}\left( \sO^{\asd}_1(\Theta,\sP_i,[\sP_k],\delta)\right)
  \right)
\end{equation}
is disjoint from $\bga'_{\Theta,\sP''}(T(\Theta,\sP'',\delta))$ for all
$\sP''\in [\sP_i<\sP_k]$.
Lemma \ref{lem:SeparatingNgh} gives
$\SO(3)\times\SO(4)\times\fS(\sP)$-invariant neighborhoods
$$
\sO^{\asd}_3(\Theta,[\sP_i<\sP_k],\delta)\sqsubset
\sO^{\asd}_2(\Theta,[\sP_i<\sP_k],\delta)\sqsubset
\sO^{\asd}_1(\Theta,[\sP_i<\sP_k],\delta),
$$
of
$$
\left.\left(\bigsqcup_{\sP''\in [\sP_i<\sP_k]}T(\Theta,\sP'',\delta)\right)\right/\fS(\sP_i)
\quad\text{in}\quad
\left.\left(\bigsqcup_{\sP''\in [\sP_i<\sP_k]}\sO^{\asd}(\Theta,\sP'',\delta)\right)\right/\fS(\sP_i)
$$
whose image under $\bga'_{\Theta,[\sP_i<\sP_k]}$ is
disjoint from the
set
\eqref{eq:LowerR4SplicingMinusOverlap}.
This disjointness yields the inclusion,
\begin{align*}
{}&
\bga'_{\Theta,\sP_i}\left( \sO^{\asd}_1(\Theta,\sP_i,\delta)\right) \cap
\bga'_{\Theta,\sP''}\left( \sO^{\asd}_1(\Theta,\sP'',\delta)\right)
\\
{}&\quad\subset
\bga'_{\Theta,[\sP_i<\sP_k]}
  \left(
    \rho^{\Theta,u}_{\sP_i,[\sP_k]}\left( \sO^{\asd}_1(\Theta,\sP_i,[\sP_k],\delta)\right)
  \right),
\end{align*}
for all
$\sP''\in [\sP_i<\sP_k]$, as required  to show that Items
\eqref{item:ConstructCollar3} and \eqref{item:ConstructCollar4} hold.  This completes the induction and hence the proof of the lemma.
\end{proof}

\begin{proof}[Proof of Proposition \ref{prop:ExistenceOfSplicedEnd}]
For each partition $\sP$ of $N_\ka$ with length greater than one, let $\sO^{\asd}_1(\Theta,\sP,\delta)$
be the $\SO(3)\times\SO(4)$ invariant subspace of $\sO^{\asd}(\Theta,\sP,\delta)$
constructed in Lemma \ref{lem:ConstructCollar}.
For any two partitions $\sP$ and $\sP'$, the images
$\bga'_{\Theta,\sP}(\sO^{\asd}_1(\Theta,\sP,\delta))$
and
$\bga'_{\Theta,\sP'}(\sO^{\asd}_1(\Theta,\sP',\delta))$
are either disjoint or are given  by the open subspaces described in
Item \eqref{item:ConstructCollar4} in Lemma \ref{lem:ConstructCollar}.
Hence, the overlaps are open subspaces and the union of the images,
$W_\ka$, has the desired properties.

The arguments in \cite[Proposition 5.10]{FLKM1} imply that if
$A'$ is defined by splicing connections $A_P$ with scales
$\delta_P$ for $P\in\sP$
onto the
product
connection
and if $\delta_P\le \delta$,
then
$$
\|F^+_{A'}\|_{L^{\sharp,2}(S^4)}
\le
\delta + \sum_{P\in\sP} \|F^+_{A_P}\|_{L^{\sharp,2}(S^4)}.
$$
The bound on $\|F^+_{A'}\|_{L^{2,\sharp}(S^4)}$ for all $A'\in W_\ka$ then
follows by induction.
\end{proof}

\section{Tubular neighborhoods of the instanton moduli space with spliced ends}
\label{subsec:TMOfSplicedEnd}
We now introduce the partial tubular neighborhood structure of the spliced
end $W_\ka$ defined in \eqref{eq:DefineSplicedEnd}.  This structure
will satisfy the
first of the conditions in \eqref{eq:TMProperties}.
In the following lemma, we introduce the
projection map, analogous to the map $\pi_i$ in \eqref{eq:TMProperties}.

\begin{lem}
\label{lem:TMProjectionOnSplicedEndR4}
For each partition $\sP$ of $N_\ka$ with
$|\sP|>1$, let $\sU(\Theta,\sP)\subset W_\ka$
be the image of the open subspace $\sO^{\asd}_1(\Theta,\sP,\delta)$
defined in Lemma \ref{lem:ConstructCollar}.
Then there is an $\SO(3)\times\SO(4)$-equivariant map,
\begin{equation}
\label{eq:SplicedEndR4Projection}
\pi(\Theta,\sP):
\sU(\Theta,\sP) \to [\Theta]\times\Si(Z_\ka(\delta),\sP),
\end{equation}
such that the following hold.
\begin{enumerate}
\item
\label{item:TMProjectionOnSplicedEndR41}
If $\sP'\in [\sP]$, then
$\sU(\Theta,\sP)=\sU(\Theta,\sP')$, and
$\Si(Z_\ka(\delta),\sP)=\Si(Z_\ka(\delta),\sP')$, and
$\pi(\Theta,\sP)=\pi(\Theta,\sP')$.
\item
\label{item:TMProjectionOnSplicedEndR42}
On the overlap $\sU(\Theta,\sP)\cap\sU(\Theta,\sP')$, we have
\begin{equation}
\label{eq:R4SplicedEndProjectionTM}
\pi(\Theta,\sP)\circ \pi(\Theta,\sP')=\pi(\Theta,\sP).
\end{equation}
\end{enumerate}
\end{lem}

\begin{proof}
The map $\pi(\Theta,\sP)$ is defined to be the
composition of the projection,
$$
\pi_{Z,\sP}:
\sO^{\asd}_1(\Theta,\sP,\delta)\to\Si(Z_\ka(\delta),\sP),
$$
with the inverse of the splicing map, $\bga'_{\Theta,\sP}$.
Item \eqref{item:TMProjectionOnSplicedEndR41} follows from this definition.

The equality \eqref{eq:R4SplicedEndProjectionTM} follows from the assertion that
the intersection of $\sU(\Theta,\sP)\cap\sU(\Theta,\sP')$ with a fiber of $\pi(\Theta,\sP')$ is contained in a fiber
of $\pi(\Theta,\sP)$.
To prove this, first note that by \eqref{eq:CompleteOverlapsR4},
$$
\sU(\Theta,\sP)\cap\sU(\Theta,\sP') =
\Imag\left(\bga'_{\Theta,\sP}\circ \rho^{\Theta,d}_{\sP,[\sP']} \right)
=
\Imag\left(\bga'_{\Theta,\sP'}\circ\rho^{\Theta,u}_{\sP,[\sP']} \right).
$$
Consider the following diagram,
$$
\begin{CD}
\tilde\sO(\Theta,\sP,\sP',\delta)
@> \rho^{\Theta,u}_{\sP,[\sP']} >>
\sO^{\asd}_1(\Theta,\sP',\delta)
\\
@V \pi_{\bx} VV @V \pi_{Z,\sP'} VV
\\
\Delta^\circ(Z_\ka(\delta),\sP)
\times
\prod_{P\in\sP}\Delta^\circ(Z_{|P|}(\delta_P),\sP'_P)
@> e(Z_\ka(\delta),\sP) >>
\Delta(Z_\ka(\delta),\sP')
\\
@V p VV
\\
\Delta^\circ(Z_\ka(\delta),\sP)
\end{CD}
$$
where the map $\pi_\bx$ is defined in  \eqref{eq:DefineProjections} and the map $p$ is the obvious projection.
Note that $p\circ \pi_\bx=\pi_{\sP,\bx}$, where $\pi_{\sP,\bx}$ is defined in \eqref{eq:DefineProjections}.
The diagram commutes by the definition of $\rho^{\Theta,u}_{\sP,[\sP']}$ in \eqref{eq:R4TrivialSplicingUpwardsOverlapMap}.
By the commutativity of this diagram and the injectivity of $e(Z_\ka(\delta),\sP)$, for $\by\in \Delta(Z_\ka(\delta),\sP')$,
and $\bx\in \Delta^\circ(Z_\ka(\delta),\sP)$, and $\tilde\bx\in p^{-1}(\bx)$ with $e(Z_\ka(\delta),\sP) (\tilde \bx)=\by$,
\begin{equation}
\label{eq:FiberInclusion1}
\Imag\left(  \rho^{\Theta,u}_{\sP,[\sP']}\right)
\cap \pi_{Z,\sP'}^{-1}(\by)
=
\rho^{\Theta,u}_{\sP,[\sP']}\left(\pi_{\bx}^{-1}(\tilde \bx) \right)
 \subset
\rho^{\Theta,u}_{\sP,[\sP']}\left( (p\circ \pi_\bx)^{-1}(\bx)\right).
\end{equation}
Next, we observe that by the definition of $\rho^{\Theta,d}_{\sP,[\sP']}$
in \eqref{eq:DefineUpwardTransition}, for $\bx\in \Delta^\circ(Z_\ka(\delta),\sP)$,
\begin{equation}
\label{eq:FiberInclusion2}
\Imag\left(  \rho^{\Theta,d}_{\sP,[\sP']}\right)
\cap \pi_{Z,\sP}^{-1}(\bx)
=
\rho^{\Theta,d}_{\sP,[\sP']}\left( (p\circ \pi_\bx)^{-1}(\bx)\right).
\end{equation}
Combining the inclusion \eqref{eq:FiberInclusion1} with the equalities \eqref{eq:FiberInclusion2}
and \eqref{eq:CommutingSplicingQuotient} yields
$$
\bga'_{\Theta,\sP'}\circ\rho^{\Theta,u}_{\sP,[\sP']}\left( \pi_{\bx}^{-1}(\tilde \bx) \right)
\subset
\bga'_{\Theta,\sP'}\circ\rho^{\Theta,u}_{\sP,[\sP']}\left((p\circ \pi_\bx)^{-1}(\bx) \right)
=
\bga'_{\Theta,\sP}\circ\rho^{\Theta,d}_{\sP,[\sP']}\left( (\pi_{Z,\sP})^{-1}(\bx)\right),
$$
which gives the inclusion of fibers required to prove \eqref{eq:R4SplicedEndProjectionTM}.
\end{proof}

The obvious identification,
\begin{equation}
\label{eq:TrivStratIdent}
\{[\Theta]\}\times Z_\ka/\fS_\ka
\cong
Z_\ka/\fS_\ka,
\end{equation}
allows us to view the projection maps $\pi(Z_\ka,\sP)$
defined in Section \ref{sec:R4Diags} as maps  on
the
strata of $W_\ka$ containing the product connection.
We identify the restriction of the projection
$\pi(\Theta,\sP)$ with such a map in the following lemma.

\begin{lem}
\label{lem:R4TrivialProjectionIdentification}
Under the identification \eqref{eq:TrivStratIdent},
the restriction of the projection map $\pi(\Theta,\sP)$
in Lemma \ref{lem:TMProjectionOnSplicedEndR4} to the
intersection,
$$
\left(\{[\Theta]\}\times Z_\ka/\fS_\ka\right)
\cap
\sU(\Theta,\sP),
$$
is equal to the map
$$
\pi(Z_\ka,\sP):U(Z_\ka,\sP)/\fS(\sP)\to
\Delta^\circ(Z_\ka,\sP)=\Si(Z_\ka,\sP)/\fS(\sP)
$$
defined in
\eqref{eq:R4DiagonalProjection}.
\end{lem}

\begin{proof}
The conclusion follows by noting that the restriction of the
splicing map $\bga'_{\Theta,\sP}$ to points where the connections being spliced
in are
the product connection is exactly the map $e(Z_\ka,\sP)$ in
\eqref{eq:NormalOfR4Diagonal}.
\end{proof}

\section{Isotopy of the spliced end of the instanton moduli space}
\label{subsec:Collar}
The final step in the proof of Theorem \ref{thm:ExistenceOfSplicedEndsModuli}
is to use the isotopy provided by the gluing and centering maps to
attach the spliced end, $W_\ka$, to the original moduli space
$\barM^{s,\natural}_{\ka}(S^4,\delta)$.

We first construct the isotopy of $W_\ka$ in $\bar\sB^s_{\ka}(S^4,\delta)$.

\begin{lem}
\label{lem:SplicedEndIsotopy}
For $\delta$ sufficiently small,
there is a continuous, smoothly-stratified map,
\begin{equation}
\label{eq:Isotopy}
R: (-\infty,1] \times W_\ka \to \bar\sB^s_\ka(S^4,2\delta),
\end{equation}
such that the following hold.
\begin{enumerate}
\item
\label{item:SplicedEndIsotopy1}
For all $t\in (-\infty,\frac{1}{2}]$
and $[A,F^s,\bx]\in W_\ka$, we have
$R(t,[A,F^s,\bx])=[A,F^s,\bx]$.
\item
\label{item:SplicedEndIsotopy2}
For all $t\in [\tthreequarter,1]$
and $[A,F^s,\bx]\in W_\ka$, we have
$R(t,[A,F^s,\bx])\in\barM^s_\ka(\delta)$.
\item
\label{item:SplicedEndIsotopy2a}
If we define
$$
\| [A,F^s,\bx]\|_{L^{\sharp,2}(S^4)}:=
\| F^+(A)\|_{L^{\sharp,2}(S^4)}
$$
then for all $t\in (-\infty,1]$ and $[A,F^s,\bx]\in W_\ka$,
$$
\|R(t,[A,F^s,\bx])\|_{L^{\sharp,2}(S^4)} \le \|[A,F^s,\bx]\|_{L^{\sharp,2}(S^4)}.
$$
\item
\label{item:SplicedEndIsotopy3}
For all $[A,F^s,\bx]\in W_\ka$, we have
$R(1,[A,F^s,\bx])\in \barM^{s,\natural}_\ka(\delta)$.
\item
\label{item:SplicedEndIsotopy4}
For all $t\in (-\infty,1]$, the map $R(t,\cdot):W_\ka\to\bar\sB^s_\ka(S^4,\delta)$
is $\SO(3)\times\SO(4)$-equivariant.
\end{enumerate}
\end{lem}

\begin{proof}
The bound on the $L^{\sharp,2}(S^4)$-norm
of the $F^+_A$ for all
$[A,\bx]\in W_\ka$ appearing in Proposition \ref{prop:ExistenceOfSplicedEnd}
and \cite[Proposition 7.6]{FLKM1} imply that for $\delta$ small enough to satisfy
\cite[Proposition 7.6]{FLKM1},
there is a continuous, smoothly-stratified,
$\SO(3)\times\SO(4)$-equivariant embedding,
\begin{equation}
\label{eq:S4GluingMap}
G:[0,1]\times W_\ka \ni (t,[A,F^s,\bx]) \mapsto [A+a_t(A),F^s,\bx] \in \bar\sB^{s,\natural}_\ka(S^4,2\delta),
\end{equation}
where $a_t(A)=d^{+,*}_Av_t(A)$ satisfies
\begin{equation}
\label{eq:S4GluingMapProperty}
F^+(A+a_t(A))=(1-t) F^+_A.
\end{equation}
Observe that $a_0(A)=0$ while $A\mapsto A+a_1(A)$ is the gluing map defined by \cite[Proposition 7.6]{FLKM1}.
For $t\in [\frac{1}{2},\tthreequarter]$,
we define
\begin{equation}
\label{eq:GluingIsotopy}
R(t,[A,F^s,\bx]) := [A+a_{4(t-1/2)}(A),F^s,\bx].
\end{equation}
This smoothly-stratified isotopy satisfied Items \eqref{item:SplicedEndIsotopy1} and
\eqref{item:SplicedEndIsotopy2}.
Although $R(3/4,[A,F^s,\bx])$ will be anti-self-dual, it need not be mass-centered.  Thus, define
$$
w(A,F^s,\bx):=z[A+a_1(A),F^s,\bx].
$$
Because, by \cite[Lemma 4.4]{Feehan_yangmillsenergy_lojasiewicz4d} and the definition \eqref{eq:CenterOfMass},
$$
z[ c_{y,1}^*A,F^s,c_{y,1}^{-1}(\bx)]
=
z[A,F^s,\bx]+y,
$$
if we define
$$
R(t,[A,F^s,\bx]):=
[ c_{(3-4t)w(A,F^s,\bx),1}^*A,F^s,c_{(3-4t)w(A,F^s,\bx),1}^{-1}(\bx)],
$$
for $t\in [3/4,1]$,
then $R(3/4,[A,F^s,\bx])=[A,F^s,\bx]$ and
$$
z\left( R(1,[A,F^s,\bx])\right)=0.
$$
Consequently, the isotopy $R$ satisfies Item \eqref{item:SplicedEndIsotopy3}.

Item \eqref{item:SplicedEndIsotopy2a} follows for $t\in [1/2,3/4]$ by observing that
equation \eqref{eq:S4GluingMapProperty} implies that
$$
\|R(t,[A,F^s,\bx])\|_{L^{\sharp,2}(S^4)}
=
\| F^+(A+a_t(A))\|_{L^{\sharp,2}(S^4)}
=
(1-t)\| F^+(A)\|_{L^{\sharp,2}(S^4)}, \quad\text{for $t\in [1/2,3/4]$.}
$$
For $t\in [3/4,1]$, Item \eqref{item:SplicedEndIsotopy2a} follows by the invariance of
the $L^{\sharp,2}(S^4)$ norm under pullback by translations.

We now prove that Item \eqref{item:SplicedEndIsotopy4} holds.
The $\SO(3)$ equivariance of $R$ is immediate.  For $S\in \SO(4)$, the $\SO(4)$ equivariance of
\eqref{eq:S4GluingMapProperty} yields
\begin{align*}
F^+(S^*A+S^*a_t(A))
{}&=
S^*F^+(A+a_t(A))
\\
{}&=
(1-t)S^*F^+(A)
\\
{}&=
(1-t)F^+(S^*A)
\\
{}&=
F^+(S^*A+a_t(S^*A)).
\end{align*}
The uniqueness of the solution $a_t(A)$ to \eqref{eq:S4GluingMapProperty} then implies that
$S^*a_t(A)=a_t(S^*A)$.  This gives the $\SO(4)$ equivariance of $R$ for $t\in [1/2,3/4]$.
To prove the equivariance for $t\in [3/4,1]$, we first note that
$$
z[ S^*A,\tilde S^{-1}F^s,S^{-1}\bx]
=
S^{-1} z[A,F^s,\bx],
$$
and so, by the $\SO(4)$ equivariance of $A\mapsto a_1(A)$,
\begin{equation}
\label{eq:EquivarianceOfCoM}
w(S^*A,\tilde S^{-1}F^s,S^{-1}\bx)=S^{-1}w(A,F^s,\bx).
\end{equation}
The equality $S\circ c_{S^{-1}w,1}=c_{w,1}\circ S$ implies that $c_{S^{-1}w,1}^*S^*=S^*c_{w,1}^*$.
Combining this with \eqref{eq:EquivarianceOfCoM} in the definition of $R(t,\cdot)$ for $3/4\le t\le 1$
completes the proof of the $\SO(4)$ equivariance of $R$ and hence the proof of Item \eqref{item:SplicedEndIsotopy4}.
\end{proof}

\begin{proof}[Proof of Theorem \ref{thm:ExistenceOfSplicedEndsModuli}]
Let
$$
\sU_3\sqsubset \sU_2 \sqsubset W_\ka
$$
be the union, over partitions $\sP$ of $N_\ka$ of length greater than one,
of the neighborhoods,
\begin{equation}
\label{eq:SmallestNgh}
\bga'_{\Theta,\sP}\left( \sO^{\asd}_j(\Theta,\sP,\delta)\right),
\end{equation}
constructed in Lemma \ref{lem:ConstructCollar}, where $W_\ka$ is defined in \eqref{eq:DefineSplicedEnd}.
These sets are closed under the actions of $\SO(3)$ and $\SO(4)$ and form a neighborhood of
$[\Theta]\times (Z_P-\{c_P\})/\fS_P$ in $W_\ka$.
Then, from the proof of \cite[Theorem 1.3.13]{Pflaum},
there is a smoothly-stratified map\footnote{A map that is continuous and whose restriction to each stratum is smooth.}
$\beta:W_\ka\to [0,1]$ with $\sU_3\subset \beta^{-1}(0)$
and $W_\ka\setminus\sU_2\subset \beta^{-1}(1)$.

Define $R(W_\ka)$ to be the image of $W_\ka$ under the map $R(\beta(\cdot),\cdot)$:
$$
R(W_\ka)
:=
\{R(\beta([A,F^s,\bx]),[A,F^s,\bx]): [A,F^s,\bx]\in W_\ka\}.
$$
Define the instanton moduli space with spliced ends by
\begin{equation}
\label{eq:DefineSplicedEndsModuliSpace}
\barM^{s,\natural}_{\spl,\ka}(S^4,\delta)
:=
\left(
\barM^{s,\natural}_{\ka}(S^4,\delta)
\setminus
R(1,W_\ka)
\right)
\cup
R(W_\ka).
\end{equation}
The $\SO(3)\times\SO(4)$ equivariance of $R$,
provided by Item \eqref{item:SplicedEndIsotopy2a}, implies that the space \eqref{eq:DefineSplicedEndsModuliSpace}
is closed under the $\SO(3)\times\SO(4)$ action.  The isotopy $R$ defines the $\SO(3)\times\SO(4)$-equivariant homeomorphism
in Item \eqref{item:ExistenceOfSplicedEndsModuli1} of the theorem.
The neighborhood $W_\ka$ appearing in Item \eqref{item:ExistenceOfSplicedEndsModuli2} is given by $R(W_\ka)$.
The Uhlenbeck neighborhoods $W(\Si)$ appearing in Item \eqref{item:ExistenceOfSplicedEndsModuli3} are given by
the sets \eqref{eq:SmallestNgh} with $j=3$.
By the construction of the
solution
$a_t(A)$ in
\eqref{eq:S4GluingMapProperty}, the $\|\cdot\|_{L^{\sharp,2}(S^4)}$-norm of
the curvature decreases along the isotopy $R(t,\cdot)$ as $t$ increases.  The bounds
appearing in \eqref{prop:ExistenceOfSplicedEnd} on this norm then imply that the space
\eqref{eq:DefineSplicedEndsModuliSpace} satisfies Item \eqref{item:ExistenceOfSplicedEndsModuli4}.
\end{proof}

\section{Properties of the instanton moduli space with spliced ends}
\label{subsec:PropsOfSplicedEnd}
We now prove two properties of the
instanton moduli space with spliced ends, $\barM^{s,\natural}_{\spl,\ka}(S^4,\delta)$, defined in \eqref{eq:DefineSplicedEndsModuliSpace}
which will be used in the Definition \ref{defn:DefineLink} of the link $\bar\bL_{\ft,\fs}$ of the lower-level reducible
$\SO(3)$ monopoles \eqref{eq:IntroLowerLevelReducibles} and of a fundamental class for $\bar\bL_{\ft,\fs}$.
To construct the link, we will require the following perturbation of the scale function that is constant along the fibers of the
projections $\pi(\Theta,\sP)$.

\begin{lem}
\label{lem:TMTubularDistanceR4}
There are a continuous, $\SO(3)\times\SO(4)$ invariant map,
\begin{equation}
\label{eq:TMScaleR4}
\tilde\la_\ka: \barM^{s,\natural}_{\spl,\ka}(S^4,\delta)\to [0,2\delta),
\end{equation}
which is smooth on each stratum, and
 neighborhoods
$\sU'(\Theta,\sP)\sqsubset\sU(\Theta,\sP)$ of
$[\Theta]\times\Si(Z_\ka(\delta),\sP)$ with the following properties.
\begin{enumerate}
\item
\label{item:TMTubularDistanceR41}
The function $\tilde\la_\ka$ is equal to the
scale function $\la$ defined in \eqref{eq:ConnectionScale} on the
complement of $W_\ka$.
\item
\label{item:TMTubularDistanceR42}
Under the identification \eqref{eq:TrivStratIdent},
the restriction of $\tilde\la_\ka$ to the
strata containing the product connection
in \eqref{eq:TrivStratIdent} is equal to the function
$t(Z_\ka)$ defined in Lemma \ref{lem:R4TMScale}.
\item
\label{item:TMTubularDistanceR43}
For all partitions $\sP$ of $N_\ka$ with $|\sP|>1$, we have
$\tilde\la_\ka=\tilde\la_\ka\circ\pi(\Theta,\sP)$ on
$\sU'(\Theta,\sP)$.
\end{enumerate}
\end{lem}

\begin{proof}
The function $\tilde\la_\ka$ is constructed in
a manner identical to that of $t(Z_P)$ in
Lemma \ref{lem:R4TMScale} with the scale function
$\la$ defined in \eqref{eq:ConnectionScale}
playing the roles of $t_{P,0}$ and the neighborhoods
$\sU(\Theta,\sP)$ and maps $\pi(\Theta,\sP)$ playing
the role of $\sU(Z_P,\sP)$ and $\pi(Z_P,\sP)$.
\end{proof}

We will use the following lemma when constructing a fundamental class
for the link of the
strata \eqref{eq:IntroLowerLevelReducibles} of gauge-equivalence
classes of reducible ideal $\SO(3)$ monopoles.

\begin{lem}
\label{lem:SplicedEndNDR}
The stratified-space $\barM^{s,\natural}_{\spl,\ka}(S^4,\delta)$ is
a Whitney-stratified space in the sense of \cite[Section 1.2]{GorMacPh}.
If $\Si\subset \barM^{s,\natural}_{\spl,\ka}(S^4,\delta)$ is the complement of the top stratum,
then there is a neighborhood $U$ of $\Si$ in $\barM^{s,\natural}_{\spl,\ka}(S^4,\delta)$
of which $\Si$ is a retraction.
\end{lem}

\begin{proof}
The smoothly-stratified homeomorphism \eqref{eq:HomeomorphismFromSplEndsToInstModuli}
between $\barM^{s,\natural}_{\spl,\ka}(S^4,\delta)$
and $\barM^{s,\natural}_{\ka}(S^4,\delta)$ preserves strata.
Hence, the lemma will follow
by proving that the pair $(\barM^{s,\natural}_{\ka}(S^4,\delta),\Si')$ has the
required properties, where $\Si'$ is the complement of the top stratum of $\barM^{s,\natural}_{\ka}(S^4,\delta)$.
The space $\barM^{s,\natural}_{\ka}(S^4,\delta)$ is
identified with a subspace of a completion of the ADHM data in \cite[pp. 470--471]{Maciocia}
or \cite[Corollary 3.4.10]{DK}.
In \cite[Theorem 1.2]{LenessADHM}, it is shown that this completion
of the ADHM data is a
semi-algebraic space (that is, defined by a finite set of polynomial inequalities on a Euclidean space) and hence Whitney-stratified
(see \cite[Section 1.2]{GorMacPh}).  The same is true for the unframed space,
$\barM^{\natural}_{\ka}(S^4,\delta)$.
In \cite[Proposition 5]{Goresky_Triang_Strat_Objects},
Goresky proves that Whitney-stratified spaces can be given the structure of a simplicial
complex in such a way that the complement of the top stratum is a subcomplex.
The existence of the neighborhood $U$ then follows from \cite[Corollary 3.3.11]{Spanier}.
\end{proof}

\chapter{The space of global splicing data}
\label{chap:GlobalSplicingData}

\section{Introduction}
\label{sec:Introduction_space_global_splicing_data}
Let $\ft$ be a \spinu structure and $\fs$ a \spinc structure such that the intersection
\begin{equation}
\label{eq:IntroLowerLevelReducibles1}
\left(M_\fs\times\Sym^\ell(X)\right)\cap \bar\sM_\ft/S^1,
\end{equation}
is non-empty, where $M_{\mathfrak{s}}$ is 
the subspace of reducible $\SO(3)$ monopoles appearing in \eqref{eq:IntroLowerLevelReducibles}.
In this
chapter
we construct the space of global splicing data
discussed in Section \ref{subsec:SpaceOfGlobalSplicingData} and which will be used to construct the link
of the subspace \eqref{eq:IntroLowerLevelReducibles1}.  This space will be the union of
spaces of splicing data, $\bar\sU(\ft,\fs,\sP)$, associated with diagonals
$\Si(X^\ell,\sP)$ defined in \eqref{eq:DefineStratumOfPartition} of $\Sym^\ell(X)$
as $\sP$ varies over partitions of $N_{\ell}=\{1,\dots,\ell\}$.
(Of course, we will do this equivariantly with respect to the
symmetric group.) To each partition $\sP$, we will define a {\em
crude splicing map\/}\label{global_splicing_data_crude_splicing_map} which  will be identical to the splicing map
defined in \cite[Equation (3.27)]{FL3} except in the following respects:
\begin{enumerate}
\item
\label{item:Introduction_space_global_splicing_data1}
The connections over $S^4$ being spliced are elements of the
instanton moduli space with spliced ends defined in Theorem \ref{thm:ExistenceOfSplicedEndsModuli}
rather than the moduli space of
anti-self-dual connections on $S^4$; and
\item
\label{item:Introduction_space_global_splicing_data2}
Each background pair,
$(A_0,\Phi_0)$, is `flattened' in the sense that
the connection $A_0$ is
replaced with one which is flat while the section $\Phi_0$ is
multiplied by a cut-off function that is zero on balls of fixed
radius around the splicing point rather than on balls whose radius
tends to zero as the connection on $S^4$ becomes more concentrated.
\end{enumerate}
The space of global splicing data is the union, over partitions of $N_\ell$,
of the images of the crude splicing maps.  To ensure that this union is a smoothly-stratified space, we require the
images of the crude splicing maps for different partitions to overlap in open subsets, as
in the construction of the spliced end of the instanton moduli space, $W_\ka$, in Section \ref{subsec:SplicedEnd}.
In Proposition \ref{prop:XOverlapControl}, we prove an analogue of
Proposition \ref{prop:CommutingSplicingSymmQuotient}, the result
describing
the overlaps of images of different splicing maps when splicing
connections
on $S^4$
with the product connection on $S^4$.
However, the proof of Proposition
\ref{prop:CommutingSplicingSymmQuotient} relied on the flatness of the
Riemannian metric
on the manifold
and of the product connection.
We will use the method from Section \ref{sec:TMMapsOnSymX}
to produce the needed flatness of the metric on manifold.
The
condition \eqref{item:Introduction_space_global_splicing_data2}
above
on the crude splicing map gives us the flatness of the connection needed to
prove Proposition \ref{prop:XOverlapControl}.

After defining the domain of the splicing map in Section
\ref{sec:SplicingData}, we show how to flatten pairs as
described above in Section \ref{sec:FlatteningPairs} and thus define
the crude splicing maps in Section \ref{sec:CrudeSplicingMap}. Given
the analysis of the overlaps of the images of the crude
splicing maps developed in Section \ref{sec:OverlapSpacesModuli}, we
construct the space of global splicing data in Section
\ref{sec:GlobalSplicingData}. In Section \ref{sec:TMStr}, we
construct the partial Thom--Mather structure (as described in Section \ref{sec:Introduction_diagTM})
which will be used in
the construction of the link $\bar\bL_{\ft,\fs}$ of the subspace \eqref{eq:IntroLowerLevelReducibles1}.
In Section \ref{sec:GlobalSplicingMap}, we describe a global splicing map
whose image will contain the subspace \eqref{eq:IntroLowerLevelReducibles1}. Finally, in Section
\ref{sec:ProjectionToSym}, we construct a map from
the space of global splicing data to $\Sym^\ell(X)$.

\section{Splicing data}
\label{sec:SplicingData}
We begin by recalling the domains of the splicing maps from \cite{FL3}.
Assume that $M_{\fs}\times\Sym^\ell(X)$ appears in the $\ell$-th level
of the space of ideal $\SO(3)$ monopoles on $\ft$, namely
$I\sM_{\ft}$.  Let $\ft(\ell)$ be the \spinu structure defined prior to
equation \eqref{eq:idealmonopoles}.
By generalizing the splicing map defined in \eqref{eq:SplicedGenConn},
we will obtain a splicing map that will attach connections on $S^4$ to the `background pair'
$(A_0,\Phi_0)\in\tsC_{\ft(\ell)}$ to create a pair
$(A,\Phi)\in\tsC_{\ft}$.

\subsection{Background pairs}
\label{subsec:Background pairs}
Let $\tilde M_{\fs}\subset\tsC_{\ft(\ell)}$ denote the image of the
embedding \eqref{eq:DefnOfIota} of the space of solutions to the
perturbed Seiberg--Witten equations from \cite[Lemma 3.12]{FL2a}.
The quotient, $M_{\fs}=\tilde M_{\fs}/\sG_{\fs}$, of $\tilde M_{\fs}$ by the group $\sG_{\fs}$ of gauge transformations embeds in
$\sC_{\ft(\ell)}$.  Let $N_{\ft(\ell),\fs}\to M_{\fs}$
be the virtual normal bundle discussed in Section \ref{subsubsec:ThickenedNeighborhood}
with an embedding of a $\delta$-disk subbundle
$N_{\ft(\ell),\fs}(\delta)\subset N_{\ft(\ell),\fs}\to \sC_{\ft(\ell)}$, as
defined in \eqref{eq:BackgroundConfigEmbedding}.
Let $\tilde N_{\ft(\ell),\fs}\to\tilde M_{\fs}$ be the pullback of
$N_{\ft(\ell),\fs}$ by the projection $\tilde M_{\fs}\to M_{\fs}$.
As in \eqref{eq:BackgroundPreConfigEmbedding}, there is an embedding
$\tilde N_{\ft(\ell),\fs}(\delta)\to\tsC_{\ft(\ell)}$ covering the embedding
$N_{\ft(\ell),\fs}(\delta)\to \sC_{\ft(\ell)}$.

\subsection{Riemannian metrics}
\label{subsubsec:FiberBundleMetrics}
Recall that in Lemma \ref{lem:FlatteningMetrics2}, we constructed for each partition
$\sP$ of $N_\ell$ a smooth family of metrics $g_{\sP}$ on $X$,
parameterized by
$\Delta^\circ(X^\ell,\sP)\subset X^\ell$ and a tubular neighborhood
$\tilde\sU(X^\ell,g_{\sP})\subset X^\ell$ of $\Delta^\circ(X^\ell,\sP)$.  The map
$\pi(X^\ell,\sP):\tilde\sU(X^\ell,g_{\sP})\to\Delta^\circ(X^\ell,\sP)$
denotes the projection map of this tubular neighborhood. We list here
the relevant properties
of this family of metrics:
\begin{enumerate}
\item
For each $\by\in\Delta^\circ(X^\ell,\sP)$, the metric
$g_{\sP,\by}$ is flat on the
support of
\[
\tilde\sU(X^\ell,g_{\sP})\cap
\pi(X^\ell,\sP)^{-1}(\by),
\]
where support is defined in Definition \ref{defn:Locally_flattened_metric}.
\item
For each pair of partitions $\sP<\sP'$ of
$N_\ell$ and point $\by'\in
\tilde\sU(X^\ell,g_{\sP})\cap\Delta^\circ(X^\ell,\sP')$ with
$\pi(X^\ell,\sP)(\by')=\by$, the equality
$g_{\sP,\by}=g_{\sP',\by'}$ holds.
\item
The metrics are $\fS_\ell$-invariant in the sense that $g_{\sP,\by}=g_{\si(\sP),\si(\by)}$
for all $\si\in\fS_\ell$.
\end{enumerate}

\subsection{Frame bundles}
For a partition $\sP$ of $N_{\ell}$, let
$$
\Fr(TX^\ell,\sP,g_\sP)\to \Delta^\circ(X^\ell,\sP),
$$
be the fiber bundle defined just before
\eqref{eq:DiagonalNormalAssocBundle}. Similarly, we define
\begin{multline}
\label{eq:DefineDiagonalOfBackgroundFrameBundle}
\Fr(\fg_{\ft(\ell)},\sP) := \left\{(F^{\fg}_1,\dots,F^{\fg}_\ell) \in
\left.\left(\prod_{i=1}^\ell \Fr(\fg_{\ft(\ell)})\right)\right|_{\Delta^\circ(X^\ell,\sP)}: \right.
\\
\left. F^{\fg}_i=F^{\fg}_j \iff \exists\, P\in\sP \text{ with } i,j\in P \right\}.
\end{multline}
We define the \emph{gluing data bundle} by
\begin{equation}
\label{eq:DefineGluingDataBundle}
\Fr(\ft,\fs,\sP)
:=
\Fr(TX,\sP,g_{\sP})\times_{\Delta^\circ(X^\ell,\sP)}\Fr(\fg_{\ft(\ell)},\sP).
\end{equation}
Recall that we denote elements of $\Delta^\circ(X^\ell,\sP)$ by
$(y_P)_{P\in\sP}$, where $y_P\in X$. Similarly, we denote
elements of the
total space of the
bundle $\Fr(\ft,\fs,\sP)$ lying over
$(y_P)_{P\in\sP}$ by $(F^T_P,F^{\fg}_P)_{P\in \sP}$, where
$F^T_P\in \Fr(TX)|_{y_P}$ and
$F^{\fg}_P\in\Fr(\fg_{\ft(\ell)})|_{y_P}$. We also write
$(F^{T,\fg}_P)_{P\in\sP}$ for such an element of
$\Fr(\ft,\fs,\sP,g_{\sP})$.

The structure groups of
the bundle
$\Fr(X,\sP)$ over
$\Delta^\circ(X^\ell,\sP)$ and over $\Si(X^\ell,\sP)$ are
$\tG(\sP)$ and $G(\sP)/\Ga(\sP)$, respectively, where
\begin{subequations}
\label{eq:DefineGluingDataBundleStructureGroup}
\begin{multline}
\label{eq:DefineGluingDataBundleStructureGroup_first_part}
\tG(\sP)
:=
\{\left((R_1,M_1),\ldots,(R_\ell,M_\ell)\right)\in
\left( \SO(4)\times\SO(3)\right)^\ell:
\\
R_i=R_j \text{ and } M_i=M_j,\
\forall\, i,j\in P\in\sP\},
\end{multline}
\begin{equation}
\label{eq:DefineGluingDataBundleStructureGroup_second_part}
G(\sP) := \tG(\sP) \rtimes \fS(\sP),
\end{equation}
\end{subequations}
where $\rtimes$ denotes the semi-direct product \cite[p. 59]{LangAlgebra} (compare the structure group
appearing on \cite[p. 252]{FrM}).
Observe that for $\sP<\sP'$, there is an inclusion,
\begin{equation}
\label{eq:StructureGroupInclusion}
\widetilde G(\sP)\to \widetilde G(\sP'),
\end{equation}
given by the inclusion of the diagonals. We can write elements of
$\tG(\sP)$ as $(R_P,M_P)_{P\in\sP}$. With this notation, the action
of $G(\sP)$ on
\begin{equation}
\label{eq:DefineGluingDataFiber}
\barM(\sP)
:=
\prod_{P\in\sP} \barM^{s,\natural}_{\spl,|P|}(\delta)
\end{equation}
is given by the factor $(R_P,M_P)$ acting by the standard
$\SO(3)\times\SO(4)$ action on
$\barM^{s,\natural}_{\spl,|P|}(\delta)$ and by the natural permutation
action of $\fS(\sP)$.

\subsection{Group actions on the frame bundles}
\label{subsubsec:GrpActionOnFrameBundle}
In addition to the action
of the structure groups $G(\sP)$ and
$\widetilde G(\sP)$ on the frame
bundle $\Fr(\ft,\fs,\sP)$, there are also actions of the
gauge group $\sG_{\fs}$ and of $S^1$ on $\Fr(\ft,\fs,\sP)$.
Both of these actions are diagonal actions defined on the factors
of $\Fr(\fg_{\ft(\ell)})$ in $\Fr(\ft,\fs,\sP)$ by the $S^1$ action
on $\fg_{\ft(\ell)}$ defined by the reduction of $\fg_{\ft(\ell)}$ to an $S^1$ bundle
as in \eqref{eq:EndReduction}.

\subsection{Space of splicing data}
\label{subsubsection:GrpActionOnFrames}
The splicing map will be defined on a subspace
$$
\tilde N_{\ft(\ell),\fs}(\delta)\times_{\sG_{\fs}} \sO(\ft,\fs,\sP)
\subset
\tilde N_{\ft(\ell),\fs}(\delta)\times_{\sG_{\fs}}\bar\Gl(\ft,\fs,\sP),
$$
where
\begin{equation}
\label{eq:GluingData}
\bar\Gl(\ft,\fs,\sP) := \Fr(\ft,\fs,\sP) \times_{G(\sP)} \barM(\sP).
\end{equation}
The subspace $\sO(\ft,\fs,\sP)$ of $\bar\Gl(\ft,\fs,\sP)$ is an open neighborhood
of
\begin{equation}
\label{eq:DefineTrivialStratumXConePoints}
\Si(\ft,\fs,\sP)
:=
\left\{ \left( (F^{T,\fg}_P)_{P\in\sP},([\Theta,F^s_P,c_{|P|}])_{P\in\sP} \right)
\subset
\bar\Gl(\ft,\fs,\sP)\right\},
\end{equation}
where $\Theta$ is the
product connection over $S^4$
and $c_{|P|}$ is the cone point of $\Sym^{|P|,\natural}(\RR^4)$.
The notation in \eqref{eq:DefineTrivialStratumXConePoints} is
motivated by the observation that because $\SO(3)\times\SO(4)$
acts trivially on $[\Theta,F^s_P,c_{|P|}]$, there is an
identification,
\begin{equation}
\label{eq:IdentifyConePointStratum} \Si(\ft,\fs,\sP) \cong
\Si(X^\ell,\sP),
\end{equation}
where $\Si(X^\ell,\sP)\subset\Sym^\ell(X)$ is defined in \eqref{eq:StratumQuotientOfDiagonal}.
We define a
subspace of $\bar\Gl(\ft,\fs,\sP)$
containing $\Si(\ft,\fs,\sP)$ by
\begin{equation}
\label{eq:DefineTrivialStratumX}
T(\ft,\fs,\sP)
:=
\left\{\left( (F^{T,\fg}_P)_{P\in\sP},([\Theta,F^s_P,\bv_P])_{P\in\sP}\right)
\in
\bar\Gl(\ft,\fs,\sP)\right\}.
\end{equation}
One requirement in the definition of $\sO(\ft,\fs,\sP)$ will be
that for any point,
$$
\left(
(F^T_P,F^{\fg}_P)_{P\in\sP},([A_P,F^s_P,\bv_P])_{P\in\sP} \right)
\in\sO(\ft,\fs,\sP),
$$
where $F^T_P$ and $F^{\fg}_P$ lie over $y_P\in X$ with $P'\neq P$, we have
\begin{equation}
\label{eq:SeparatingCondition1}
8\la([A_P,F^s_P,\bv_P]) < \min_{P'\in\sP,P'\neq P} d(y_P,y_P'),
\end{equation}
where $\la$ is the scale parameter
\eqref{eq:ConnectionScale}.  Finally, we define
\begin{equation}
\label{eq:DefineSplicingDataToDiagonalProjection} \pi_{\Si}:
\tilde N_{\ft(\ell),\fs}(\delta)\times_{\sG_{\fs}}\bar\Gl(\ft,\fs,\sP)
 \to \Si(X^\ell,\sP)
\end{equation}
to be the obvious projection.  The $S^1$ action on $\tilde
N_{\ft(\ell),\fs}$ defined in \eqref{eq:S1ZActionOnN} defines an $S^1$
action,
\begin{equation}
\label{eq:S1ActionOnPieceHalfWeight}
S^1\times \tilde N_{\ft(\ell),\fs}(\delta)\times_{\sG_{\fs}}\bar\Gl(\ft,\fs,\sP)
\to
\tilde N_{\ft(\ell),\fs}(\delta)\times_{\sG_{\fs}}\bar\Gl(\ft,\fs,\sP),
\end{equation}
and we will see that all the constructions of this
chapter are
equivariant with respect to this $S^1$ action.
The following analogue of \cite[Lemma 3.6]{FLLevelOne} gives
two descriptions of the action \eqref{eq:S1ActionOnPieceHalfWeight}.

\begin{lem}
\label{lem:S^1ActionOnSplicingData}
The following two circle
actions on
$\tilde N_{\ft(\ell),\fs}(\delta)\times_{\sG_{\fs}}\bar\Gl(\ft,\fs,\sP)$ are
equivalent:
\begin{enumerate}
\item The action \eqref{eq:S1ZActionOnN} on $\tN_{\ft(\ell),\fs}(\delta)$
and the trivial action on $\bar\Gl(\ft,\fs,\sP)$.
\item The diagonal action with weight two on the fibers of $\tilde
N_{\ft(\ell),\fs}(\delta)\to \tilde M_{\fs}$ and with the action on
$\bar\Gl(\ft,\fs,\sP)$ described in Section
\ref{subsubsec:GrpActionOnFrameBundle} with weight two.
\end{enumerate}
\end{lem}

\begin{proof}
The proof of this lemma is identical to that of \cite[Lemma
3.6]{FLLevelOne} using the equality in \eqref{eq:S1ZActionOnN} of
$S^1$ actions.
\end{proof}

We will need to refer to the $S^1$ action,
\begin{equation}
\label{eq:S1ActionOnPiece}
S^1\times
\tilde N_{\ft(\ell),\fs}(\delta)\times_{\sG_{\fs}}\bar\Gl(\ft,\fs,\sP)
\to
\tilde N_{\ft(\ell),\fs}(\delta)\times_{\sG_{\fs}}\bar\Gl(\ft,\fs,\sP),
\end{equation}
given by the diagonal action with weight one on the fibers of
$\tilde N_{\ft(\ell),\fs}(\delta)\to \tilde M_{\fs}$ and with the action on
$\bar\Gl(\ft,\fs,\sP)$ defined by the $S^1$ action described in Section
\ref{subsubsec:GrpActionOnFrameBundle} with weight one.
The action described in
Lemma \ref{lem:S^1ActionOnSplicingData} is the action
\eqref{eq:S1ActionOnPiece} with weight two.

Each $\si\in\fS_\ell$ defines a smoothly-stratified diffeomorphism $\si:\bar\Gl(\ft,\fs,\sP)\to\bar\Gl(\ft,\fs,\si(\sP))$,
where $\si(\sP)$ is the partition defined in \eqref{eq:ActionOnPartitions} by
the following relabelling map,
\begin{multline}
\label{eq:SymmetricGroupActionOnSplicingData1}
\bA:=\left( (F^{\fg}_P,F^T_P),[A_P,F^s_P,\bv_P]\right)_{P\in\sP}
\\
\mapsto
\si(\bA):=
\left( (F^{\fg}_{\si(P)},F^T_{\si(P)}),[A_{\si(P)},F^s_{\si(P)},\bv_{\si(P)}]\right)_{P\in\sP}.
\end{multline}
Because the stabilizer $\fS(\sP)$ of $\sP$ under the action \eqref{eq:ActionOnPartitions}
is in the structure group $G(\sP)$ in \eqref{eq:DefineGluingDataBundleStructureGroup_second_part},
if $\si\in\fS(\sP)$, then the map defined by $\si$ by \eqref{eq:SymmetricGroupActionOnSplicingData1} is trivial.

\section{The flattening map on pairs}
\label{sec:FlatteningPairs}
The commutativity of the diagram
\eqref{eq:CommutingSplicingQuotient} for the instanton moduli space with spliced ends
depended on the flatness of the metric on $\RR^4$ and on the
flatness of the
product connection $\Theta$ around the splicing points.
The locally flat metrics $g_{\sP}$ described in Section
\ref{subsubsec:FiberBundleMetrics} have this property. We now
introduce a method for achieving the same kind of local flatness
for the background connection $A_0$. This method also yields a way
of cutting off the background section consistently from one
stratum of the subspace \eqref{eq:IntroLowerLevelReducibles1}
to another.

\begin{defn}
\label{defn:SeparatingFunctions} A collection of smooth functions
$s_P:\Delta^\circ(X^\ell,\sP)\to (0,1)$, indexed by $P\in\sP$, is
a \emph{separating family} if the following hold:
\begin{enumerate}
\item The functions are $\fS(\sP)$-invariant in the sense that
$s_{P}=s_{\si(P)}$  and $s_{P}=s_P\circ\si$, for all
$\si\in\fS(\sP)$.
\item For all
$\by=(y_P)_{P\in\sP}\in\Delta^\circ(X^\ell,\sP)$ and for all
$P\neq P'\in\sP$, the balls $B(y_P,4s_P(\by)^{1/2})$ and
$B(y_{P'},4s_{P'}(\by)^{1/2})$ are disjoint.
\end{enumerate}
\end{defn}

\begin{defn}
Let $\{s_P\}$ be a separating family of smooth functions on
$\Delta^\circ(X^\ell,\sP)$. A connection
$\hat A$ on a vector bundle $\fg\to X$
is \emph{flat
with respect to $\{s_P\}$} at
$\by=(y_P)_{P\in\sP}\in\Delta^\circ(X^\ell,\sP)$ if the connection
$\hat A$ is flat on
\begin{equation}
\label{eq:SeparatingFunctionBallS^1}
\bigcup_{P\in\sP} B(y_P,2s_P(\by)^{1/2}).
\end{equation}
A pair $(A,\Phi)\in\tilde\sC_{\ft(\ell)}$ is \emph{flat with respect to
$\{s_P\}$} if $\hat A$ is flat with respect to $\{s_P\}$ and
$\Phi$ vanishes on
\begin{equation}
\label{eq:SeparatingFunctionBalls2}
\bigcup_{P\in\sP} B(y_P,4s_P(\by)^{1/3}).
\end{equation}
\end{defn}

Lemma \ref{lem:FlatteningPairs} below provides the analogue of the flattening construction
for metrics given in Lemma \ref{lem:FlatteningMetrics}.
Recall from Section \ref{subsubsec:PU2Monopoles}
that for a \spinu structure $\ft$ and \spinu
connection $A$ on the \spinu bundle $V_\ft\to X$
then $\hat A$ denotes the unique connection induced by $A$ on the $\so(3)$ bundle $\fg_{\ft}$
appearing in \eqref{eq:EndSplitting}.

\begin{lem}
\label{lem:FlatteningPairs}
Let $\sP$ be a partition of $N_\ell$.
Let $\{s_P\}$ be a separating family of smooth functions on
$\Delta^\circ(X^\ell,\sP)$. Then there is an $S^1$-equivariant
map,
\begin{equation}
\label{eq:DefineFlatteningMap}
\bTheta_{\sP}': \tN_{\ft(\ell),\fs}(\delta)\times
\Delta^\circ(X^\ell,\sP) \to \tsC_{\ft(\ell)},
\end{equation}
which is $\sG_{\fs}$-equivariant and descends to an $S^1$-equivariant
map,
$$
N_{\ft(\ell),\fs}(\delta)\times \Delta^\circ(X^\ell,\sP) \to \sC_{\ft(\ell)}.
$$
In addition, for all
$\by=(y_P)_{P\in\sP}\in\Delta^\circ(X^\ell,\sP)$, if
$(A',\Phi')=\bTheta'_{\sP}((A_0,\Phi_0),\by)$, then
\begin{enumerate}
\item The connection $\hat A'$ is flat on
$$
\bigcup_{P\in\sP} B(y_P,2s_P(\by)^{1/2})
$$
and equal to $\hat A_0$ on the complement of
\begin{equation}
\label{eq:ConnToBeFlat}
\bigcup_{P\in\sP} B(y_P,4s_P(\by)^{1/2}).
\end{equation}
\item The section $\Phi'$ is identically zero on
$$
\bigcup_{P\in\sP} B(y_P,4s_P(\by)^{1/3})
$$
and is equal to $\Phi$ on the complement of
\begin{equation}
\label{eq:SectionToVanish}
\bigcup_{P\in\sP} B(y_P,8s_P(\by)^{1/3}).
\end{equation}
\end{enumerate}
If the connection $\hat A_0$ is already flat on the space
\eqref{eq:ConnToBeFlat}, then $A'=A_0$. If the section $\Phi$
is already identically zero on the space \eqref{eq:SectionToVanish}, then
$\Phi'=\Phi_0$.
\end{lem}

\begin{proof}
The map $\bTheta_{\sP'}$ is initially defined as a map,
$$
\bTheta_{\sP}':\tilde N_{\ft(\ell),\fs}(\delta)\times\Fr(\ft,\fs,\sP) \to
\tsC_{\ft(\ell)},
$$
as described in \cite[Section 3.2]{FL3}.  Specifically, the frames and
radial parallel translation define a
product connection on the
space \eqref{eq:ConnToBeFlat}.  If $\hat A_0$ is flat on the space
\eqref{eq:ConnToBeFlat}, this product connection is equal to $\hat
A_0$.  We now use a cut-off function to interpolate between the
product connection and $\hat A_0$ on the annuli around the points
$y_P$, defining a new connection $\hat A'$ which is flat on the
inner balls and equal to $\hat A_0$ outside the space
\eqref{eq:ConnToBeFlat}. Observe that the resulting connection
$\hat A'$ is independent of the choice of the frames used as it is
equal to the connection obtained by splicing in the
product connection on the space \eqref{eq:ConnToBeFlat} which has stabilizer
$\SO(3)$ and is invariant under the rotation action on the tangent frame
(see \cite[Proposition 7.2.9]{DK}). Hence, the
map $\bTheta_{\sP}'$ descends to the domain stated in the lemma.

The section $\Phi'$ is defined by multiplying $\Phi_0$ by a
cut-off function equal to one on the space \eqref{eq:SectionToVanish}
and vanishing on the balls of radius
$4s_P(\by)^{1/3}$.
\end{proof}

To define the crude splicing maps in such a way that an analogue
of Proposition \ref{prop:CommutingSplicingSymmQuotient}
describes the overlaps of their images,
we must redefine the flattening maps
in such a way that for $\by'\in\Delta^\circ(X^\ell,\sP')$ near
$\by\in\Delta^\circ(X^\ell,\sP)$, the following equality holds:
\[
\bTheta_{\sP'}(A_0,\Phi_0,\by')=\bTheta_{\sP}(A_0,\Phi_0,\by).
\]
(Compare the preceding equality with the second property of the metrics $g_{\sP}$ appearing in
Section \ref{subsubsec:FiberBundleMetrics}.) To this end, we introduce
a refinement of the flattening map constructed in Lemma
\ref{lem:FlatteningPairs}.

\begin{lem}
\label{lem:CompatibleFlattening}
Continue the assumptions of Lemma
\ref{lem:FlatteningPairs}. Then, for every partition $\sP$ of
$N_\ell$ there are separating functions
$\{\tilde s_P\}$ and an
$S^1$-equivariant map,
\begin{equation}
\label{eq:CompatibleFlattening}
\bTheta_{\sP}: \tilde N_{\ft(\ell),\fs}(\delta)\times \Delta^\circ(X^\ell,\sP) \to \tsC_{\ft(\ell)},
\end{equation}
which is $\fS(\sP)$-invariant, $\sG_{\fs}$-equivariant,
descends to the quotient,
$$
N_{\ft(\ell),\fs}(\delta)\times \Delta^\circ(X^\ell,\sP) \to \sC_{\ft(\ell)},
$$
and has the properties that if $(A',\Phi')=\bTheta_{\sP}(A_0,\Phi_0,\by)$, where
$\by=(y_P)_{P\in\sP}\in\Delta^\circ(X^\ell,\sP)$, then the following hold:
\begin{enumerate}
\item
\label{item:CompatFlatteningSymmetricEquiv}
If the symmetric group $\fS_\ell$ acts on the partitions of $N_\ell$ by the action \eqref{eq:ActionOnPartitions}
and on $X^\ell$ by the obvious permutation action, then for $\si\in\fS_\ell$,
$$
\bTheta_{\si(\sP)}(A_0,\Phi_0,\si(\by))=\bTheta(A_0,\Phi_0,\by)
$$
and the separating functions $\{\tilde s_{\si(P)}\}_{\si(P)\in\si(\sP)}$ associated to
$\si(\sP)$ satisfy $\tilde s_{\si(P)}=\tilde s_P\circ \si$.
\item
\label{item:BallOnWhichFlattenedConnIsFlat}
The connection $\hat A'$ is flat on
$$
\bigcup_{P\in\sP} B(y_P,2\tilde s_P(\by)^{1/2}).
$$
\item The section $\Phi'$ vanishes on
$$
\bigcup_{P\in\sP} B(y_P,4\tilde s_P(\by)^{1/3}).
$$
\end{enumerate}
Furthermore, if $\sP<\sP'$ and the tubular neighborhood
$\tilde\sU(X^\ell,g_{\sP})$ of $\Delta^\circ(X^\ell,\sP)$ is suitably small,
$\by'\in\tilde\sU(X^\ell,g_{\sP})\cap\Delta^\circ(X^\ell,\sP')$, and
$\pi(X^\ell,g_{\sP})(\by')=\by\in\Delta^\circ(X^\ell,\sP), $ then
\begin{equation}
\label{eq:ConsistentFlattening}
\bTheta_{\sP}(A_0,\Phi_0,\by)=\bTheta_{\sP'}(A_0,\Phi_0,\by').
\end{equation}
\end{lem}

\begin{rmk}
The separating functions $\{\tilde s_P\}$ appearing in the statement of Lemma \ref{lem:CompatibleFlattening}
differ from those appearing in Lemma \ref{lem:FlatteningPairs} in that the functions $\{\tilde s_P\}$ satisfy
the additional condition \eqref{eq:BallInclusions} relating the separating functions of one partition to another
and which is used to derive \eqref{eq:ConsistentFlattening}.
\end{rmk}

\begin{proof}
The proof is similar to the construction of the consistent flat
families of metrics in Lemma \ref{lem:FlatteningMetrics2}. One
constructs $\bTheta_{\sP}$ by induction on $\sP$.  As described in
Section \ref{subsec:EnumStrata}, the induction uses one partition from each orbit of the action of $\fS_\ell$
on the set of partitions of $N_\ell$ and extends the results to all partitions in the orbit
so that Item \eqref{item:CompatFlatteningSymmetricEquiv} holds.  Thus, Item \eqref{item:CompatFlatteningSymmetricEquiv}
follows trivially from the method of induction.

For the crudest partition, $\sP_0:=\{N_\ell\}$, of $N_\ell$, define
$\bTheta_{\sP_0}:=\bTheta_{\sP_0}'$ (as constructed in Lemma
\ref{lem:FlatteningPairs}).

Assume that $\bTheta_{\sP}$ satisfying the conclusions of the lemma has been constructed for all $\sP<\sP'$. To
construct $\bTheta_{\sP'}$, we first define
$$
\bTheta_{\sP',\sP}:\tilde N_{\ft(\ell),\fs}(\delta)\times
\left(\sU(X^\ell,g_{\sP})\cap \Delta^\circ(X^\ell,\sP')\right) \to \tsC_{\ft(\ell)},
$$
by
$$
\bTheta_{\sP',\sP}(A_0,\Phi_0,\by')
:=
\bTheta_{\sP}(A_0,\Phi_0,\pi(X^\ell,g_{\sP})(\by')).
$$
Observe that if $\sP_1<\sP$ and $\sP_2<\sP$ and
$\sU(X^\ell,g_{\sP_1})\cap\sU(X^\ell,g_{\sP_2})$ is non-empty, then we can
assume $\sP_1<\sP_2<\sP$.  The Thom--Mather
condition \eqref{eq:DiagonalsTM1}
on the
projections $\pi(X^\ell,g_{\sP})$,
$$
\pi(X^\ell,g_{\sP_1})\circ\pi(X^\ell,g_{\sP_2})=\pi(X^\ell,g_{\sP_1}),
$$
and the inductive assumption imply that for
$\by'\in\sU(X^\ell,g_{\sP_1})\cap\sU(X^\ell,g_{\sP_2})\cap\Delta^\circ(X^\ell,\sP')$ we have
\begin{align*}
{}&\bTheta_{\sP',\sP_2}(A_0,\Phi_0,\by')
\\
{}&\quad
=
\bTheta_{\sP_2}(A_0,\Phi_0,\pi(X^\ell,g_{\sP_2})(\by'))
\\
{}&\quad
=
\bTheta_{\sP_1}(A_0,\Phi_0,\pi(X^\ell,g_{\sP_1})\circ\pi(X^\ell,g_{\sP_2})(\by'))
\quad\text{(by inductive hypothesis and \eqref{eq:ConsistentFlattening})}
\\
{}&\quad
=
\bTheta_{\sP_1}(A_0,\Phi_0,\pi(X^\ell,g_{\sP_1})(\by'))
\quad\text{(by Thom--Mather property \eqref{eq:DiagonalsTM1})}
\\
&\quad
=
\bTheta_{\sP',\sP_1}(A_0,\Phi_0,\by')
\quad\text{(by the definition of $\bTheta_{\sP',\sP_1}$)}
.
\end{align*}
Thus, if we define the ends of the diagonal $\Delta^\circ(X^\ell,\sP')$ by
\begin{equation}
\label{eq:EndOfDiagonal}
\sU(\sP')
:=
\bigcup_{\sP<\sP'}\sU(X^\ell,g_{\sP})\cap\Delta^\circ(X^\ell,\sP'),
\end{equation}
then the map
\begin{align*}
{}&\bTheta_{\sP',\sU}: \tilde N_{\ft(\ell),\fs}(\delta)\times \sU(\sP')\to \tsC_{\ft(\ell)},
\\
{}&\quad\bTheta_{\sP',\sU}(A_0,\Phi_0,\by')
:=
\bTheta_{\sP',\sP}(A_0,\Phi_0,\by'), \quad\text{for } \by'\in\sU(X^\ell,g_{\sP}),
\end{align*}
is defined consistently for $\by'\in\sU(X^\ell,g_{\sP_1})\cap\sU(X^\ell,g_{\sP_2})\cap\Delta^\circ(X^\ell,\sP')$.
We next define
$$
\bTheta_{\sP',1}:\tilde N_{\ft(\ell),\fs}(\delta)\times \Delta^\circ(X^\ell,\sP')
\to \tsC_{\ft(\ell)}
$$
by
$$
\bTheta_{\sP',1}(A_0,\Phi_0,\by')
:=
(1-\bchi(\by'))(A_0,\Phi_0) + \bchi(\by')\bTheta_{\sP',\sU}(A_0,\Phi_0,\by'),
$$
where $\bchi:\Delta^\circ(X^\ell,\sP)\to [0,1]$ is an $\fS(\sP)$-invariant function supported in \eqref{eq:EndOfDiagonal} and equal
to one on a neighborhood $\sU_2$ of the end of
$\Delta^\circ(X^\ell,\sP)$ (hence $\sU_2$ is a proper subspace of
$\sU(\sP')$). By shrinking the tubular neighborhoods
$\sU(X^\ell,g_{\sP})$, we can find a family of separating functions
$\{\tilde s_{P'}\}_{P'\in\sP'}$, such that for
$\by'=(y'_{P'})_{P'\in\sP'}
\in\sU(X^\ell,g_{\sP})\cap\Delta^\circ(X^\ell,\sP')$ with
$\pi(X^\ell,g_{\sP})(\by')=\by=(y_P)_{P\in\sP}$, we have
\begin{equation}
\label{eq:BallInclusions}
\begin{aligned}
B(y'_{P'},8\tilde s_{P'}(\by')^{1/3}) &\Subset B(y_P,4\tilde
s_{P}(\by)^{1/3}),
\\
B(y'_{P'},4\tilde s_{P'}(\by')^{1/2}) &\Subset
B(y_P,2\tilde s_{P}(\by)^{1/2}), \quad \forall\, P'\in\sP'_P.
\end{aligned}
\end{equation}
From the inclusions \eqref{eq:BallInclusions}, the inductive
hypothesis, and the final statement of Lemma
\ref{lem:FlatteningPairs}, we have the equality,
$$
\bTheta_{\sP'}' \left(
\bTheta_{\sP,1}(A_0,\Phi_0,\by'),\by'\right)
=
\bTheta_{\sP,1}(A_0,\Phi_0,\by'), \quad\forall\, \by'\in \sU_2.
$$
Thus, for $\by'\in\sU_2$, the pair
$\bTheta_{\sP,1}(A_0,\Phi_0,\by')$ is already flat on the relevant
balls and so the flattening construction of Lemma
\ref{lem:FlatteningPairs} does not change
$\bTheta_{\sP,1}(A_0,\Phi_0,\by')$.
Therefore, if we define
$$
\bTheta_{\sP'}(A_0,\Phi_0,\by')
:=
\bTheta_{\sP'}'(\bTheta_{\sP,1}(A_0,\Phi_0,\by'),\by'),
$$
this flattening map will satisfy the condition \eqref{eq:ConsistentFlattening},
completing the construction of the map $\bTheta_{\sP'}$ and hence
completing the induction.
\end{proof}

\section{The crude splicing map}
\label{sec:CrudeSplicingMap}
The crude splicing map,
\begin{equation}
\label{eq:DefineCrudeSplicing}
\bga''_{\ft,\fs,\sP}:
\tilde N_{\ft(\ell),\fs}(\delta)\times_{\sG_{\fs}}
\sO(\ft,\fs,\sP)
\subset
\tilde N_{\ft(\ell),\fs}(\delta)\times_{\sG_{\fs}}
\bar\Gl(\ft,\fs,\sP)
 \to
\bar\sC_{\ft},
\end{equation}
where $\bar\sC_{\ft}$ is defined in \eqref{eq:IdealPairs}
and
$\bar\Gl(\ft,\fs,\sP)$ is defined in  \eqref{eq:GluingData},
has a domain defined by a \emph{separation condition},
\begin{equation}
\label{eq:DefineCrudeSplicingDomain}
\begin{aligned}
{}& \sO(\ft,\fs,\sP)
\\
{}&\quad= \{(F^{T,\fg}_P,[A_P,F^s_P,\bv_P])_{P\in\sP} \in
\bar\Gl(\ft,\fs,\sP): 8\la([A_P,F^s_P,\bv_P])^{1/3}\le
s_P(F^{T,\fg}_P)\},
\end{aligned}
\end{equation}
where $\{s_P\}$ is the pullback of the family of separating
functions $\{\tilde s_P\}$ constructed in Lemma
\ref{lem:CompatibleFlattening} by the projection map,
$\bar\Gl(\ft,\fs,\sP)\to\Si(X^\ell,\sP)$.
We will denote elements of the domain of the crude splicing map by $[(A_0,\Phi_0),\bA]$
and use the following notation.
\begin{notn}
\label{notn:SplicingDataFactors}
A point $[(A_0,\Phi_0),\bA]$ in the domain of the map \eqref{eq:DefineCrudeSplicingDomain} is given by
\begin{enumerate}
\item
\label{item:SplicingDataFactors1}
A pair $(A_0,\Phi_0)\in\tilde N_{\ft(\ell),\fs}(\delta)$,
\item
\label{item:SplicingDataFactors2}
A point $(y_P)_{P\in\sP}\in \Delta^\circ(X^\ell,\sP)$,
\item
\label{item:SplicingDataFactors3}
Frames $F^{\fg}_P\in \Fr(\fg_{\ft(\ell)})|_{y_P}$,
\item
\label{item:SplicingDataFactors4}
Frames $F^T_P\in \Fr(TX)|_{y_P}$,
\item
\label{item:SplicingDataFactors5}
Gauge-equivalence classes of
framed, mass-centered
connections, $[A_P,F^s_P]\in \barM^{s,\natural}_{\spl,|P|}$, over $S^4$.
\end{enumerate}
\end{notn}
The crude splicing
map will be $S^1$-equivariant with respect to the $S^1$ action on the
domain given in Lemma \ref{lem:S^1ActionOnSplicingData} and the
$S^1$ action on the range given by the action
\eqref{eq:S1ZAction}.  Before defining the crude splicing map,
we review the definition of the splicing map defined
in \cite[Section 3.2]{FL3}.

\subsection{The standard splicing map}
\label{subsubsec:StdSplicingMap}
For clarity, we shall refer to the splicing map defined in \cite[Section 3.2]{FL3}
as the \emph{standard} splicing map.
The \emph{crude} splicing map will be similar to the
standard splicing map, with some differences described below.
We now sketch the construction in \cite[Section 3.2]{FL3}
of a pair $(A',\Phi')=\bga'_{\ft,\fs,\sP}([(A_0,\Phi_0),\bA])$, where $[(A_0,\Phi_0),\bA]$ is a point in the domain of the crude splicing map, \eqref{eq:DefineCrudeSplicingDomain}, given by the data in
Notation \ref{notn:SplicingDataFactors}.

Recall from Section \ref{subsubsec:PU2Monopoles} that for a \spinu connection $A$ on a \spinu bundle $V_\ft$,
we write $\hat A$ for the connection induced by $A$ on the $\so(3)$ bundle $\fg_{\ft}$
appearing in \eqref{eq:EndSplitting}.  If we restrict to the \spinu connections inducing a fixed connection on
$\det(V_\ft)$, then the map $A\mapsto \hat A$ is a bijection.

For $[A_P,F^s_P]\in  \barM^{s,\natural}_{\spl,|P|}$ as in Item \eqref{item:SplicingDataFactors5} of Notation \ref{notn:SplicingDataFactors},
the connection $\hat A_P$ will be an orthogonal connection on the
real rank-three, Riemannian vector bundle, $\fg_P\to S^4$.
Similarly, for $(A_0,\Phi_0)\in\tilde N_{\ft(\ell),\fs}(\delta) $ as in Item \eqref{item:SplicingDataFactors1} of Notation \ref{notn:SplicingDataFactors},
$\hat A_0$ will be a connection on
$\fg_{\ft(\ell)} \to X$.
Define
\begin{equation}
\label{eq:ScaleOfSplicingData}
\la_P:=\la([A_p,F^s_P])
\end{equation}
to be the  scale of the connection $A_P$ as defined in
\eqref{eq:ConnectionScale}.
For each $P\in\sP$, let
$\tau(A_P,F^s_P)$ be the trivialization of
$\fg_P$ over $S^4\less\{n\}$
(where $n$ is the North Pole)
defined by parallel translation of the frame
$F^s_P$ along great circles.  This trivialization defines a
product (and hence flat) connection, $\Theta(A_P,F^s_P)$.
Recall from Section \ref{subsec:SpaceOfConn}, that the map $\varphi_n:\RR^n\to S^4\less\{s\}$ was defined
by stereographic projection.
Let $B(r)\subset \RR^4$ be the open ball
of radius $r$ centered at the
origin
and let
$\Om(r_1,r_2):=B(r_2)\setminus \bar B(r_1)$.
More generally, for $x\in X$, let $\Omega(x,r_1,r_2):=B(x,r_2)\setminus \bar{B}(x,r_1)$,
where $B(x,r)\subset X$ denotes the open ball with center $x$ and radius $r$ defined by
the Riemannian metric, $g$.
For a smooth function $\beta:\mathbb{R}\to [0,1]$ obeying $\beta(t)=1$ for $t\le 1$ and $\beta(t)=0$ for $t\ge 2$,
constant $\varepsilon>0$, and point $x\in\RR^4$, define $\chi_{x,\varepsilon}:\RR^4\to [0,1]$ by
$\chi_{x,\varepsilon}:=\beta(|\cdot -x|/\varepsilon)$ and define $\chi_{n,\varepsilon}:S^4\to [0,1]$ by
$\chi_{n,\varepsilon}:=\chi_{0,\varepsilon}\circ\varphi_n^{-1}$.
Thus, $\chi_{n,\frac{1}{4}\lambda^{1/2}}$ is a smooth cut-off function on $S^4$ that is equal to one on
$\varphi_n(B(\frac{1}{4}\lambda^{1/2}))$ and is supported in $\varphi_n(B(\frac{1}{2}\lambda^{1/2}))$.

Define a connection on
$\fg_P \to S^4$ by
\begin{equation}
\label{eq:FlattenedS4Conn}
\hat A'_P
:=
\begin{cases}
\Theta(A_P,F^s_P)
{}& \text{on } S^4\setminus \varphi_n(B(\frac{1}{2} \la_P^{1/2})),
\\
(1-\chi_{n,\frac{1}{4}\lambda_P^{1/2}})) \Theta(A_P,F^s_P) +  \chi_{n,\frac{1}{4}\lambda_P^{1/2}}\hat A_P
{}& \text{on }\varphi_n(\Om(\frac{1}{4} \la_P^{1/2},\frac{1}{2}\la_p^{1/2})),
\\
\hat A_P
{}& \text{on }\varphi_n(B(\frac{1}{4}\la_P^{1/2})).
\end{cases}
\end{equation}
For $x\in X$, define $\chi_{x,\varepsilon}:X\to [0,1]$ by $\chi_{x,\varepsilon}:=\beta(\dist_g(\cdot,x)/\varepsilon)$.
Thus, $1-\chi_{x,2\lambda^{1/2}}$ is a smooth cut-off function on $X$ that is equal to one on $X\setminus B(x,4\lambda^{1/2})$
and is equal to zero on $B(x,2\lambda^{1/2})$.
Parallel translation with respect to the connection $\hat A_0$
of the frame $F^{\fg}_P$
(from Item \eqref{item:SplicingDataFactors3} of Notation \ref{notn:SplicingDataFactors})
along radial geodesics from the point
$y_P$ defines a trivialization
$\tau(A_0,F^{\fg}_P)$ of
the bundle $\fg_{\ft(\ell)}$ over a ball
in $X$ centered at the point $y_P$.
Let $\Theta(A_0,F^{\fg}_P)$ be the
product connection
defined by
this trivialization.
We then define
an orthogonal connection on $\fg_{\ft(\ell)}$ by
\begin{equation}
\label{eq:CutOffBackgroundConnection}
\hat A_0'
:=
\begin{cases}
\hat A_0
{}&\text{on } X\setminus \cup_P B(x_P,4\la_P^{1/2}),
\\
(1-\chi_{y_P,2\lambda_P^{1/2}})\hat A_0
+
\chi_{y_P,2\lambda^{1/2}}\Theta(A_P,F^{\fg}_P)
{}&\text{on } \Omega(y_P,2\lambda_P^{1/2},4\lambda_P^{1/2}),
\\
\Theta(A_0,F^{\fg}_P)
{}&\text{on }B(y_P,2\la_P^{1/2}).
\end{cases}
\end{equation}
The frame $F^T_P$ from Item \eqref{item:SplicingDataFactors4} of Notation \ref{notn:SplicingDataFactors}, 
the exponential map around $x_P$,
and stereographic projection identify the annuli
$\Om(x_P,\frac{1}{2} \la_P^{1/2},2\la_P^{1/2})\subset X$ with
the annuli
$\varphi_n(\Om(\frac{1}{2} \la_P^{1/2},2\la_P^{1/2}))\subset S^4$.
Identifying the trivializations
$\tau(A_0,F^{\fg}_P)$
and $\tau(A_P,F^s_P)$ over these annuli allows
us to
glue the bundle $\fg_{\ft(\ell)}$
and bundles $\fg_P$
(for $P\in\sP$)
to create a bundle isomorphic to $\fg_{\ft}$,
\begin{equation}
\label{eq:ConstructionOfSplicedBundle}
\fg_{\ft}
\cong
\fg_{\ft(\ell)}|_{X\setminus\cup_P B(y_P,\thalf \la_P^{1/2})}
\cup
\bigcup_{P\in\sP}\
\fg_P|_{\varphi_n(B(2\la_P^{1/2}))}
.
\end{equation}
Identifying these trivializations gives an
identification of the
product
connections
$\Theta(A_0,F^{\fg}_P)$ and
$\Theta(A_P,F^s_P)$.
Let $\Theta(A_0,F^{\fg}_P,F^T_P,A_P,F^s_P)$ be the
product
connection over $\Om(x_P,\frac{1}{2} \la_P^{1/2},2\la_P^{1/2})\subset X$
given by
this identification.
We define
an orthogonal connection on $\fg_\ft$ by
\begin{equation}
\label{eq:StandardSplicingConn}
\hat A'
:=
\begin{cases}
\hat A_0 '
{}&\text{on } X\setminus\cup_P B(x_P,4\la_P^{1/2}),
\\
\Theta(A_0,F^{\fg}_P,F^T_P,A_P,F^s_P)
{}&\text{on } \Om(y_P,\frac{1}{2}\la_P^{1/2},2\la_P^{1/2}),
\\
\hat A_P'
{}&\text{on } B(y_P,2\la_P^{1/2}).
\end{cases}
\end{equation}
If $x\in X$, then $\chi_{x,4\lambda^{1/3}}:X\to [0,1]$ is a smooth cut-off function such that
$1-\chi_{x,4\lambda^{1/3}}$ is equal to one on $X\setminus B(x,8\lambda^{1/3})$ and is equal to zero
on $B(x,4\lambda^{1/3})$.
We define
a section of $V_\ft$ by
\begin{equation}
\label{eq:StandardCutOffSection}
\Phi'
:=
\begin{cases}
\Phi
{}&\quad\text{on } X\setminus\cup_P B(y_P,8\la_P^{1/3}),
\\
(1-\chi_{y_P,4\lambda_P^{1/3}})\Phi
{}&\quad\text{on } \Om(y_P,4\la_P^{1/3},8\la_P^{1/3}),
\\
0
{}&\quad\text{on } \cup_P B(y_P,4\la_P^{1/3}).
\end{cases}
\end{equation}
The \emph{standard splicing map}
has the same domain and range as the crude splicing map in \eqref{eq:DefineCrudeSplicing}.
For the data $\bA$  given in  Notation \ref{notn:SplicingDataFactors}, the standard splicing map is defined by
\[
\bga_{\ft,\fs,\sP}'([(A_0,\Phi_0),\bA]) :=[A',\Phi'],
\]
where $A'$
is the \spinu connection with $\hat A'$
defined by \eqref{eq:StandardSplicingConn}
and $\Phi'$ is defined by \eqref{eq:StandardCutOffSection}.

\subsection{Construction of the crude splicing map}
\label{subsubsec:ConstrCrudeSplice}
The definition of the crude splicing map, $\bga''_{\ft,\fs,\sP}$, will differ
from the preceding description of the standard splicing map, $\bga_{\ft,\fs,\sP}'$, in the following ways:
\begin{enumerate}
\item
\label{item:SplicingMapDifference1}
The metric on $X$ used to identify geodesic balls in $X$ with balls
centered at
the North Pole of $S^4$ is fixed for the standard splicing
map.
The metric on $X$ so used in the definition of the crude splicing map
varies in the smooth family of metrics $g_{\sP}$ parameterized by the
splicing point $\by\in\Delta^\circ(X^\ell,\sP)$
constructed in Lemma \ref{lem:FlatteningMetrics2}.
Note that for $\by=(y_P)$, the balls of radius $r$ centered at $y_P$ defined by the Riemannian metrics $g$
and $g_{\sP,\by}$ are equal by Item \eqref{item:FlatteningMetrics_4} of Lemma \ref{lem:FlatteningMetrics}.
\item
\label{item:SplicingMapDifference2}
For the definition of the crude splicing map,
the trivializations of the background bundle, $\fg_{\ft(\ell)}$,
are defined by parallel translation with
respect to the flattened connection, $\hat A_0''$ where
$(A_0'',\Phi'')=\Theta_{\sP}(A_0,\Phi_0,\bx)$, instead of with
respect to $\hat A_0$.
\item
\label{item:SplicingMapDifference3}
In the definition of the standard splicing map in \cite[\S
3.2]{FL3}, the background connection $A_0$ is flattened
over the balls
$$
B(y_P,2\la_P^{1/2}),
$$
while the background section is multiplied by a cut-off function
that is identically zero on
$$
B(y_P,4\la_P^{1/3}),
$$
where $\la_P$ is the scale of the connection on $S^4$ being spliced
in at $x_P\in X$. The crude splicing map is defined instead by
replacing the background pair $(A_0,\Phi_0)$ with the background
pair $\bTheta_{\sP}(A_0,\Phi_0,\by)$ constructed in Lemma
\ref{lem:CompatibleFlattening}.  Hence, the background pair is
flattened over balls whose radius does not depend on the scale of
the connections over $S^4$.
\end{enumerate}

Let $[(A_0,\Phi_0),\bA]$ be a point in the domain of $\bga_{\ft,\fs,\sP}''$
in \eqref{eq:DefineCrudeSplicing}.  We will continue to use the notation from
Notation \ref{notn:SplicingDataFactors} for the data given by $\bA$.
For the flattening map $\bTheta_{\sP}$ in \eqref{eq:CompatibleFlattening}, let
\begin{equation}
\label{eq:CrudeSplicingFlattenedBackgroundPair}
(A''_0,\Phi''):=
\bTheta_{\sP}(A_0,\Phi_0,(y_P)_{P\in\sP}),
\end{equation}
where $(A_0,\Phi_0)$ and $(y_P)_{P\in\sP}$ are as in Items  \eqref{item:SplicingDataFactors1}
and \eqref{item:SplicingDataFactors2}  of Notation  \ref{notn:SplicingDataFactors}, respectively.
Let $\hat A_P$ and $\la_P$ be as defined in \eqref{eq:FlattenedS4Conn}  and \eqref{eq:ScaleOfSplicingData}.
Using the same notation as in \eqref{eq:StandardSplicingConn}, we define
\begin{equation}
\label{eq:CrudeSplicingConn}
\hat A''
:=
\begin{cases}
\hat A_0''
{}&\text{on } X\setminus \cup_P B(x_P,4\la_P^{1/2}),
\\
\Theta(A_0'',F^{\fg}_P,F^T_P,A_P,F^s_P)
{}&\text{on } \Om(y_P,\frac{1}{2}\la_P^{1/2},2\la_P^{1/2}),
\\
\hat A_P'
{}&\text{on } B(y_P,2\la_P^{1/2}).
\end{cases}
\end{equation}
Note, however, that the identification of $\Om(y_P,\frac{1}{2}\la_P^{1/2},2\la_P^{1/2})$
and
$\varphi_n(\Om(\frac{1}{2}\la_P^{1/2},2\la_P^{1/2}))$ is obtained using 
$\varphi_n$ and the exponential map around $y_P$ given by 
the
tangent
frame $F^T_P$
and the flattened metric $g_{\sP,\by}$
constructed in Lemma \ref{lem:FlatteningMetrics2}.
We then define the \emph{crude splicing map} by
\begin{equation}
\label{eq:CrudeSplicingMapExplicit}
\bga''_{\ft,\fs,\sP}([(A_0,\Phi_0),\bA]) := [A'',\Phi''],
\end{equation}
where $A''$ is the \spinu connection inducing the connection $\hat A''$ appearing in \eqref{eq:CrudeSplicingConn}.

\subsection{Properties of the crude splicing map}
\label{subsec:Properties of the crude splicing map}
The following equivariance of the crude splicing maps with respect to the symmetric group
map \eqref{eq:SymmetricGroupActionOnSplicingData1} holds because the point $[A'',\Phi'']$
in \eqref{eq:CrudeSplicingMapExplicit}
is independent of the labelling of the elements of the subsets $P\subset N_\ell$ making up the partition $\sP$.

\begin{lem}
\label{lem:CrudeSplicingSymmetricEquiv}
If $\si\in\fS_\ell$ and $\bA\in\sO(\ft,\fs,\sP)$, then
$$
\bga_{\ft,\fs,\sP}''\left( [(A_0,\Phi_0),\bA]\right)
=
\bga_{\ft,\fs,\si(\sP)}''\left( [(A_0,\Phi_0),\si(\bA])\right),
$$
where $\si(\bA)$ is the image of $\bA$ under the map \eqref{eq:SymmetricGroupActionOnSplicingData1}.
\end{lem}

The following lemma gives a more explicit description
of the crude splicing map when restricted to a subspace
of its domain.

\begin{lem}
\label{lem:CrudeSplicingOnTrivials}
Let $T(\ft,\fs,\sP)\subset\bar\Gl(\ft,\fs,\sP)$
be the subspace defined in \eqref{eq:DefineTrivialStratumX}.
Then there is a
smoothly-stratified diffeomorphism
$$
k_{\sP}: T(\ft,\fs,\sP) \cong \tilde\nu(X^\ell,\sP)/\fS(\sP),
$$
where $\tilde\nu(X^\ell,\sP)$ is the normal bundle defined in
\eqref{eq:DiagonalNormalAssocBundle}. The restriction of the crude splicing map $\bga''_{\ft,\fs,\sP}$
to
$$
\tilde N_{\ft(\ell),\fs}(\delta)\times_{\sG_{\fs}}T(\ft,\fs,\sP)
$$
is given by
$$
\bga''_{\ft,\fs,\sP}
\left(
[(A_0,\Phi_0),\bA]
\right)
=
\left(
\bTheta_{\sP}\left( (A_0,\Phi_0),(y_P)_{P\in\sP}\right),
e(X^\ell,g_{\sP})(k_{\sP}(\bA)\right),
$$
where $(A_0,\Phi_0)\in \tilde N_{\ft(\ell),\fs}(\delta)$, and $\bA\in T(\ft,\fs,\sP)$ lies over $(y_P)_{P\in\sP}\in\Delta^\circ(X^\ell,\sP)$,
and $\bTheta_{\sP}$ is the flattening map
\eqref{eq:CompatibleFlattening}, and
$e(X^\ell,g_{\sP})$ is the exponential map defined in \eqref{eq:VaryingMetricExponentialMap},
and
$e(X^\ell,g_{\sP})(k_{\sP}(\bA))$
represents an element of $\Sym^\ell(X)$.
\end{lem}

\begin{proof}
By the definition \eqref{eq:DefineTrivialStratumX}, points in $T(\ft,\fs,\sP)$ can be written as
$$
\bA=\left( (F^T_P,F^{\fg}_P)_{P\in\sP},([\Theta,F^s_P,\bv_P])_{P\in\sP}
\right),
$$
where
$F^T_P$ and $F^{\fg}_P$ are the frames appearing in Notation \ref{notn:SplicingDataFactors} and
$\Theta$ is the
product connection on $S^4\times \so(3)$ with
$[\Theta,F^s_P,\bv_P]\in \barM^{s,\natural}_{\spl,|P|}(\delta)$.
The mass-centering condition on $[\Theta,F^s_P,\bv_P]$ implies that $\bv_P\in Z_P(\delta)$,
where $Z_P(\delta)$ is defined following \eqref{eq:SingularPointScale}.  Because the $\SO(3)$ action on
the frame is trivial on such elements of $\barM^{s,\natural}_{\spl,|P|}(\delta)$,
we obtain a smoothly-stratified diffeomorphism by setting
$$
k_{\sP}
\left( (F^T_P,F^{\fg}_P)_{P\in\sP},([\Theta,F^s_P,\bv_P])_{P\in\sP}
\right)
:=
\left( (F^T_P,v_P)_{P\in\sP}\right),
$$
and comparing the definition of $\tilde\nu(X^\ell,\sP)$ in \eqref{eq:DiagonalNormalAssocBundle}.
The assertion regarding the crude splicing map follows immediately from its definition
in \eqref{eq:CrudeSplicingMapExplicit}, \eqref{eq:CrudeSplicingFlattenedBackgroundPair} and
\eqref{eq:CrudeSplicingConn}.
\end{proof}

\begin{rmk}
\label{rmk:LocalSplicingMap}
In contrast to the conclusion of Lemma \ref{lem:CrudeSplicingOnTrivials}, the standard splicing map would
satisfy, for $\bA\in \Si(\ft,\fs,\sP)$, the
subspace defined in \eqref{eq:DefineTrivialStratumXConePoints},
$$
\bga'_{X,\sP}([(A_0,\Phi_0),\bA]) =\left([
A_0,\Phi_0],e(X^\ell,g_{\sP})(k_{\sP}(\bA))\right),
$$
that is, the background pair would not be flattened.
\end{rmk}

\section{Overlap spaces and maps}
\label{sec:OverlapSpacesModuli}
As described in Section \ref{subsec:Overlaps}, we
describe the intersection of the images of the crude splicing maps $\bga''_{\ft,\fs,\sP}$  and
$\bga''_{\ft,\fs,\sP'}$,
for partitions $\sP<\sP'$  of $N_\ell$, by
introducing spaces of overlap data, $\Gl(\ft,\fs,\sP,[\sP'])$,
and defining maps from the space of overlap data to the domains of these crude splicing maps in such a way
that the intersection of the images of $\bga''_{\ft,\fs,\sP}$  and
$\bga''_{\ft,\fs,\sP'}$ is described by a diagram similar to \eqref{eq:IntroOverlapCD}.
The construction is analogous to that of the spliced end $W_\ka$ defined in \eqref{eq:DefineSplicedEnd},
specifically the proof of the commutativity of the diagram \eqref{eq:CommutingSplicingQuotient}.

\subsection{The overlap space}
\label{subsubsec:OverlapSpaceModuli}
Recall from  \eqref{eq:DefineGluingDataBundle} and \eqref{eq:GluingData}
that the space of gluing data is a fiber bundle  $\bar\Gl(\ft,\fs,\sP')\to\Si(X^\ell,\sP')$.
For partitions $\sP<\sP'$  of $N_\ell$,
one can think of the space of overlap data associated to these partitions as the restriction of the fiber bundle $\bar\Gl(\ft,\fs,\sP')\to\Si(X^\ell,\sP')$
to the intersection of the tubular neighborhood  $\sU(X^\ell,\sP)$ of $\Si(X^\ell,\sP)$ with
the stratum $\Si(X^\ell,\sP')$. This construction will be similar to that of the space
$\sO^{\asd}(\Theta,\sP,[\sP'],\delta)$ defined in \eqref{eq:SplicedModuliOverlapR4_quotient}.
Because our constructions must respect the action of the symmetric group,
the space of overlap data must
describe the intersection of the image of $\bga''_{\ft,\fs,\sP}$  with
$\bga''_{\ft,\fs,\sP''}$ for all $\sP''\in [\sP<\sP']$, where $[\sP<\sP']$ is defined
in \eqref{eq:ConjugateRefinements} to be the partitions conjugate under the symmetric group to $\sP'$
that are refinements of $\sP$.
(See the discussion prior to Lemma \ref{lem:EndOfUpperStratum} for more on this requirement.)
Then, we define the space of overlap data by
\begin{equation}
\label{eq:OverlapGluingDataSpace}
\begin{aligned}
{}& \bar\Gl(\ft,\fs,\sP,[\sP'])
\\
&\quad:=
\Fr(\ft,\fs,\sP)\times_{G(\sP)}
\bigsqcup_{\sP''\in[\sP<\sP']} \prod_{P\in\sP} \left(
\Delta^\circ(Z_{|P|}(\delta_P),\sP''_P) \times
\barM(\sP''_P)\right),
\end{aligned}
\end{equation}
where
$\Fr(\ft,\fs,\sP$ is the gluing-data bundle defined in
\eqref{eq:DefineGluingDataBundle}, and
$\barM(\sP''_P)$ is the product of instanton moduli spaces with spliced ends defined
in \eqref{eq:DefineGluingDataFiber}, and
$\sP''_P$ is the partition of $P$ given by the subsets of $P$ in $\sP''$,
and the diagonal $\Delta^\circ(Z_{|P|}(\delta_P),\sP''_P)$ is defined in \eqref{eq:DiagonalsOfZP}.
Recall that we denote elements of
$\Delta^\circ(Z_{|P|}(\delta_P),\sP''_P)$ by
$(v_Q)_{Q\in\sP''_P}$,
where $v_Q\in\RR^4$ (as described after Lemma \ref{lem:DiagonalsInZ}). We then denote elements of
$\bar\Gl(\ft,\fs,\sP,[\sP'])$ by
\begin{equation}
\label{eq:NotationForXOverlapSplicingData}
\bA'
=
\left(
(F^T_P,F^{\fg}_P)_{P\in\sP}, \left( (v_Q)_{Q\in\sP''_P},
([A_Q,F^s_Q,\bv_Q])_{Q\in\sP''_P}\right)_{P\in\sP} \right)
\in
\bar\Gl(\ft,\fs,\sP,[\sP']),
\end{equation}
where $(F^T_P,T^{\fg}_P)_{P\in\sP}\in \Fr(\ft,\fs,\sP)$,
and $(v_Q)_{Q\in\sP''_P}\in \Delta^\circ(Z_{|P|},\sP''_P)$, and
$[A_Q,F^s_Q,\bv_Q]\in \barM^{s,\natural}_{\spl,|Q|}(\delta)$.

The projection $\Fr(\ft,\fs,\sP)\to\Si(X^\ell,\sP)$
induces a projection,
\begin{equation}
\label{eq:ProjectionGluingOverlapToDiagonalInX}
\pi_{\Si}:
\bar\Gl(\ft,\fs,\sP,[\sP']) \to \Si(X^\ell,\sP).
\end{equation}
Suppose $\Theta$ is the product connection on $S^4\times \so(3)$
and $c_{|Q|}$ is the cone point defined following \eqref{eq:PuncturedTrivialStrata}.
Because $[\Theta,F^s_Q,c_{|Q|}]\in\barM^{s,\natural}_{\spl,|Q|}$ is a fixed point of the action of $\SO(3)\times\SO(4)$
on $\barM^{s,\natural}_{\spl,|Q|}$
the subspace
\begin{equation}
\label{eq:DefineTrivialInOverlap}
\begin{aligned}
{}& T(\ft,\fs,\sP,[\sP'])
\\
&\quad := \{ \left( (F^T_P,F^{\fg}_P)_{P\in\sP}, \left(
(v_Q)_{Q\in\sP''_P},
([\Theta,F^s_Q,c_{|Q|}])_{Q\in\sP''_P}\right)_{P\in\sP} \right)
\in \bar\Gl(\ft,\fs,\sP,[\sP'])\},
\end{aligned}
\end{equation}
is identified with the
subspace $\nu(X^\ell,\sP\to [\sP'])$ of the normal bundle $\nu(X^\ell,\sP)$ of $\Delta^\circ(X^\ell,\sP)$,
\begin{equation}
\label{eq:IdentifyNormal}
k_{\Si,\sP,[\sP']}:
T(\ft,\fs,\sP,[\sP']) \cong \nu(X^\ell,\sP\to [\sP'])\subset\nu(X^\ell,\sP),
\end{equation}
where $\nu(X^\ell,\sP\to [\sP'])$ was defined in \eqref{eq:EndOfUpperStratum} and which by \eqref{eq:NormalData} can be presented as
\begin{equation}
\label{eq:IdentifyingNormalBundleOfDiagonalInSymmProduct1}
\nu(X^\ell,\sP\to [\sP'])
\cong
\Fr(TX^\ell,\sP,g_{\sP}) \times_{G(T,\sP)}
\bigsqcup_{\sP''\in [\sP<\sP']} \prod_{P\in\sP}\Delta^\circ(Z_{|P|}(\delta_P),\sP''_P).
\end{equation}
Comparing the preceding presentation with \eqref{eq:OverlapGluingDataSpace}, we see that
there is a projection,
\begin{equation}
\label{eq:ProjectionGluingOverlapToNormalInX}
\pi_{\Si,\sP,[\sP']}:
\bar\Gl(\ft,\fs,\sP,[\sP']) \to
T(\ft,\fs,\sP,[\sP']),
\end{equation}
given by (using the same notation as in \eqref{eq:DefineTrivialInOverlap})
\begin{align*}
{}& \pi_{\Si,\sP,[\sP']} \left( (F^T_P,F^{\fg}_P)_{P\in\sP}, \left(
(v_Q)_{Q\in\sP''_P},
([A_Q,F^s_Q,\bx_Q])_{Q\in\sP''_P}\right)_{P\in\sP} \right)
\\&\quad=
\left( (F^T_P,F^{\fg}_P)_{P\in\sP}, \left( (v_Q)_{Q\in\sP''_P},
([\Theta,F^s_Q,c_{|Q|}])_{Q\in\sP''_P}\right)_{P\in\sP} \right).
\end{align*}
The composition of the projection \eqref{eq:ProjectionGluingOverlapToNormalInX} and
the identification \eqref{eq:IdentifyNormal} defines a map,
\begin{equation}
\label{eq:ProjectionGluingOverlapToNormalInX1}
\pi_{\nu(\Si)} = k_{\Si,\sP,[\sP']}\circ\pi_{\Si,\sP,[\sP']}:
\bar\Gl(\ft,\fs,\sP,[\sP']) \to
\nu(X^\ell,\sP\to [\sP']),
\end{equation}
which we call the projection onto the normal bundle.

The $S^1$ action on $\bar\Gl(\ft,\fs,\sP)$
discussed in Section \ref{subsubsec:GrpActionOnFrameBundle}
is defined by the $S^1$ action on $\Fr(\ft,\fs,\sP)$
given in Section \ref{subsubsec:GrpActionOnFrameBundle}.
This circle action then gives an action of $S^1$ on
$\bar\Gl(\ft,\fs,\sP,[\sP'])$ through
its action on the factor
$\Fr(\ft,\fs,\sP)$ in the
definition \eqref{eq:OverlapGluingDataSpace} of $\bar\Gl(\ft,\fs,\sP,[\sP'])$.
Hence, there is an $S^1$ action,
\begin{equation}
\label{eq:S1ActionOnOverlap}
S^1\times
\tilde N_{\ft(\ell),\fs}(\delta)\times_{\sG_{\fs}} \bar\Gl(\ft,\fs,\sP,[\sP'])
\to
\tilde N_{\ft(\ell),\fs}(\delta)\times_{\sG_{\fs}} \bar\Gl(\ft,\fs,\sP,[\sP']),
\end{equation}
defined in the
same manner as the circle action \eqref{eq:S1ActionOnPiece}.

\subsection{The upwards overlap map}
\label{subsec:Global_splicing_data_upwards overlap map}
We now construct the upwards overlap map described in the diagram
\eqref{eq:IntroOverlapCD}.
This construction will be analogous to that of the map $\rho^{\Theta,u}_{\sP,[\sP']}$ defined in \eqref{eq:SymmetricProductNormalBundlesOfHigherStratum}.
As mentioned in the beginning of the previous section, the image of this map will be contained in
the union of the gluing data bundles
corresponding to all partitions  $\sP''\in[\sP<\sP']$,
\begin{equation}
\label{eq:SimultUpwardsOverlapSpaces}
\bar\Gl(\ft,\fs,[\sP<\sP'])
:=
\bigsqcup_{\sP''\in[\sP<\sP']}\bar\Gl(\ft,\fs,\sP'')/\fS(\sP),
\end{equation}
where $\bar\Gl(\ft,\fs,\sP'')$ is defined in \eqref{eq:GluingData} (compare the space \eqref{eq:XlTubNghOverlapUpwardsSym} used to define the upwards transition map for the subspaces of $\Sym^\ell(X)$)
and $\fS(\sP)$ acts by permuting elements of $[\sP<\sP']$.
For the space of overlap data $\bar\Gl(\ft,\fs,\sP,[\sP'])$  defined in \eqref{eq:OverlapGluingDataSpace}, set
\begin{equation}
\label{eq:XUpwardsTransitionDomain1}
\sO(\ft,\fs,\sP,[\sP']):= \pi_{\nu(\Si)}^{-1} \left(\sO(X^\ell,g_{\sP})\cap \nu(X^\ell,\sP\to [\sP']) \right) \subset\bar\Gl(\ft,\fs,\sP,[\sP']).
\end{equation}
We will
define a subspace,
$$
\tilde N(\delta)\times_{\sG_{\fs}}
\sO(\ft,\fs,\sP,[\sP']) \subset
\tilde N(\delta)\times_{\sG_{\fs}}\bar\Gl(\ft,\fs,\sP,[\sP']),
$$
and an \emph{upwards overlap map},
\begin{equation}
\label{eq:XUpwardsTransitionMap}
\rho^{\ft,\fs,u}_{\sP,[\sP']}:
\tilde N_{\ft(\ell),\fs}(\delta)\times_{\sG_{\fs}}
\sO(\ft,\fs,\sP,[\sP']) \to
\tilde N_{\ft(\ell),\fs}(\delta)
\times_{\sG_{\fs}}
\bar\Gl(\ft,\fs,[\sP<\sP']),
\end{equation}
will be defined, using the notation for points in the domain from \eqref{eq:NotationForXOverlapSplicingData},
to act as the identity on the background pair $(A_0,\Phi_0)\in\tilde N_{\ft(\ell),\fs}(\delta)$,
the identity
on the $S^4$ connections,
and parallel translation of the
frames $(F^T_P,F^{\fg}_P)$ from $y_P\in X$ to each $x_Q\in X$ for $Q\in\sP''_P$,
where $(x_Q)_{Q\in\sP''}\in\Delta^\circ(X^\ell,\sP'')$ is the image of $\bA'$ under the
composition $e(X^\ell,g_{\sP})\circ\pi_{\nu(\Si)}$.
This parallel translation is done with respect to the locally
flattened metric and background connection $A_0''$.
If $P\in\sP\cap\sP'$, then because $\Delta^\circ(Z_P(\delta),\sP'_P)$
is given by
$|P|$-copies of the zero vector (see the sentence following \eqref{eq:DoubleDiagNormalFiber}),
we would have $x_Q=y_P$ and
the parallel translation for this case would be given by the identity map.

We now give a
formal definition of the upwards overlap map.
Assume that $x$ is a point in a geodesic ball
centered at
$y\in X$. Then for
any connection $A$ on a bundle over $X$ and
frame $F$ of this
bundle lying over $y$, let $T^A_{x,y}(F)$ denote parallel
translation, respect to the connection $A$, of the
frame $F$ from $y$ to $x$ along the radial geodesic. We let
$T^{g}_{x,y}$ denote parallel translation
with respect to the
Levi-Civita connection of a metric $g$ on $X$.

Using the notation \eqref{eq:NotationForXOverlapSplicingData}
for
$\bA'\in\sO(\ft,\fs,\sP,[\sP'])$, write
\begin{equation}
\label{eq:UpwardTransitionBasePoints}
\pi_{\Si}(\bA')=\by=(y_P)_{P\in\sP} \quad\text{and}\quad
e(X^\ell,g_{\sP})\circ \pi_{\nu(\Si)}(\bA')= \by'=(x_Q)_{Q\in\sP''}.
\end{equation}
For $(A_0,\Phi_0)\in \tilde N_{\ft(\ell),\fs}(\delta)$, let
$\bTheta_{\sP}(A_0,\Phi_0,\by)=(A_0',\Phi_0')$ be as
appearing
in
\eqref{eq:CompatibleFlattening}.
Then,
\begin{equation}
\label{eq:UnparamUpwardsTransition}
\begin{aligned}
{}& \rho^{\ft,\fs,u}_{\sP,[\sP']}\left([(A_0,\Phi_0),\bA']\right)
\\
{}= & \rho^{\ft,\fs,u}_{\sP,[\sP']} \left(\left[ (A_0,\Phi_0),(F^T_P,F^{\fg}_P)_{P\in\sP},
\left( (v_Q)_{Q\in\sP''_P},
([A_Q,F^s_Q,\bx_Q])_{Q\in\sP''_P}\right)_{P\in\sP}\right] \right)
\\
{}&\qquad :=  \left[(A_0,\Phi_0),
      \left( (T^{g_{\sP,\by}}_{x_Q,y_P}(F^T_P),T^{\hat A_0'}_{x_Q,y_P}(F^{\fg}_P))_{Q\in\sP''},
([A_Q,F^s_Q,\bx_Q])_{Q\in\sP''_P}\right) \right],
\end{aligned}
\end{equation}
where the indices $P$ and $Q$ appearing in the parallel
translations $T^{g_{\sP,\by}}_{x_Q,y_P}$ and $T^{\hat A_0'}_{x_Q,y_P}$ satisfy
$Q\subseteqq P$. The metric $g_{\sP,\by}$ defining the parallel
translation
$T^{g_{\sP,\by}}_{x_Q,y_P}$ is the locally flattened metric on $X$
defined in Lemma \ref{lem:FlatteningMetrics2}.
The proof of the following lemma is then straightforward.

\begin{lem}
\label{lem:UpwardsTransition}
Let $\sP<\sP'$ be partitions of
$N_\ell$.
Then the map $\rho^{\ft,\fs,u}_{\sP,[\sP']}$ defined above
is an open embedding which is equivariant with respect
to the $S^1$ actions defined by \eqref{eq:S1ActionOnPiece}
and \eqref{eq:S1ActionOnOverlap}.
\end{lem}

\subsection{Downwards overlap map}
\label{subsec:Global_splicing_data_downwards overlap map}
The downwards overlap map is an $S^1$-equivariant map,
$$
\rho^{\ft,\fs,d}_{\sP,\sP'}:
\tilde N_{\ft(\ell),\fs}(\delta)
\times_{\sG_{\fs}}
\sO_d(\ft,\fs,\sP,[\sP'])
\to
\tilde N_{\ft(\ell),\fs}(\delta)
\times_{\sG_{\fs}}\bar\Gl(\ft,\fs,\sP),
$$
defined similarly to the map $\rho^{\Theta,d}_{\sP,[\sP']}$ defined in \eqref{eq:SymmetricProductNormalBundlesOfLowerStratum}
and where the domain,
$$
\sO_d(\ft,\fs,\sP,[\sP'])\subset
\sO(\ft,\fs,\sP,[\sP']) \subset
\bar\Gl(\ft,\fs,\sP,[\sP']),
$$
will be defined in \eqref{eq:XUpwardsTransitionDeomainRequirement}.
We wish to define, for $\bA'\in \bar\Gl(\ft,\fs,\sP,[\sP'])$,
using the notation of \eqref{eq:NotationForXOverlapSplicingData},
and $(A_0,\Phi_0)\in\tilde N_{\ft(\ell),\fs}(\delta)$,
\begin{equation}
\label{eq:XDownwardTransition}
\begin{aligned}
{}& \rho^{\ft,\fs,d}_{\sP,[\sP']}([(A_0,\Phi_0),\bA'])
\\
{}= &
\rho^{\ft,\fs,d}_{\sP,[\sP']} \left( \left[(A_0,\Phi_0),(F^T_P,F^{\fg}_P)_{P\in\sP},
\left( (v_Q)_{Q\in\sP''_P},
([A_Q,F^s_Q,\bx_Q])_{Q\in\sP''_P}\right)_{P\in\sP}\right] \right)
\\
&:=
\left[ (A_0,\Phi_0),(F^T_P,F^{\fg}_P)_{P\in\sP}, \left( \bga'_{\Theta,\sP''_P}
\left( (v_Q)_{Q\in\sP''_P}, ([A_Q,F^s_Q,\bx_Q])_{Q\in\sP''_P}
\right) \right) \right],
\end{aligned}
\end{equation}
where $\bga'_{\Theta,\sP''_P}$ is the splicing map defined in \eqref{eq:SplicedGenConn}
and we use the convention that $\bga'_{\Theta,\sP''_P}$ is the identity map when $P\in\sP\cap\sP''$.
To motivate the definition of the domain $\sO_d(\ft,\fs,\sP,[\sP'])$ and to show that $\rho^{\ft,\fs,d}_{\sP,[\sP']}$
is a bundle map, we give the following description.
If we define a map $\rho^{\ft,\fs,d}_{f,\sP,[\sP']}$
on a subspace of the fiber of the map $\pi_\Si$ defined in \eqref{eq:ProjectionGluingOverlapToDiagonalInX}
by
\begin{equation}
\label{eq:DownwardsInclusionFiberMap}
\rho^{\ft,\fs,d}_{f,\sP,[\sP']}
:=
\bigsqcup_{\sP''\in [\sP<\sP']}\prod_{P\in\sP}\bga'_{\Theta,\sP''_P},
\end{equation}
(the $f$ subscript standing for `fiber') with domain and range,
\begin{equation}
\label{eq:DowardsOverlapFiberMap}
\begin{CD}
\sO\subset
\bigsqcup\limits_{\sP''\in [\sP<\sP']}
\prod\limits_{P\in\sP}
\left(
\Delta^\circ(Z_P(\delta_P),\sP''_P)
    \times
    \prod\limits_{Q\in\sP''_P} \barM^{s,\natural}_{\spl,|Q|}(\delta_Q)
\right)
\\
@V \rho^{\ft,\fs,d}_{f,\sP,[\sP']} VV
\\
\prod\limits_{P\in\sP}\barM^{s,\natural}_{\spl,|P|}(\delta_P)
\end{CD}
\end{equation}
then by the $\SO(3)\times\SO(4)$ equivariance of $\bga'_{\Theta,\sP''_P}$
given by Lemma \ref{lem:GroupEquiv}
and the $\fS(\sP)$ equivariance of the construction,
we see
that the map $\rho^{\ft,\fs,d}_{f,\sP,[\sP']}$ is $G(\sP)$-equivariant.
The map $\rho^{\ft,\fs,d}_{\sP,[\sP']}$ from \eqref{eq:XDownwardTransition} is then
defined by the property that its restriction to the fiber of $\pi_\Si$ equals $\rho^{\ft,\fs,d}_{f,\sP,[\sP']}$.

The preceding description of $\rho^{\ft,\fs,d}_{\sP,[\sP']}$
as the extension from the fiber of $\rho^{\ft,\fs,d}_{f,\sP,[\sP']}$
implies that
the domain of $\rho^{\ft,\fs,d}_{\sP,[\sP']}$
must be contained in
the set of points $\bA'\in \bar\Gl(\ft,\fs,\sP,[\sP'])$ which lie in the subset
$\sO$ appearing in \eqref{eq:DowardsOverlapFiberMap}
of the fiber of the map $\pi_\Si$ defined in \eqref{eq:ProjectionGluingOverlapToDiagonalInX}.
More  precisely,
we must assume that domain of
$\rho^{\ft,\fs,d}_{\sP,[\sP']}$
satisfies
\begin{multline}
\label{eq:XUpwardsTransitionDeomainRequirement}
\sO_d(\ft,\fs,\sP,[\sP'])
\\
\subseteqq \left\{ \left( (F^T_P,F^{\fg}_P)_{P\in\sP}, \left(
(v_Q)_{Q\in\sP''_P},
([A_Q,F^s_Q,\bx_Q])_{Q\in\sP''_P}\right)_{P\in\sP} \right) \in
\bar\Gl(\ft,\fs,\sP,[\sP']): \right.
\\
\left.\left( (v_Q)_{Q\in\sP''_P},
([A_Q,F^s_Q,\bx_Q])_{Q\in\sP''_P}\right) \in
\sO^{\asd}_1(\Theta,\sP''_P,\delta_P) \quad\text{for all $P\in\sP$}\right\},
\end{multline}
where $\sO^{\asd}_1(\Theta,\sP''_P,\delta_P)$ is the domain of $\bga'_{\Theta,\sP''_P}$
defined in Lemma \ref{lem:ConstructCollar}.

\begin{lem}
\label{lem:XUpwardsTransition}
If
$\sO_d(\ft,\fs,\sP,[\sP'])\subseteqq
\bar\Gl(\ft,\fs,\sP,[\sP'])$ is any open neighborhood of
the subspace
$T(\ft,\fs,\sP,[\sP'])$ satisfying
\eqref{eq:XUpwardsTransitionDeomainRequirement}, then the
restriction  of  $\rho^{\ft,\fs,d}_{\sP,[\sP']}$ to
$\tilde N_{\ft(\ell),\fs}(\delta)\times_{\sG_{\fs}}\sO_d(\ft,\fs,\sP,[\sP'])$
is an $S^1$-equivariant
open embedding, whose image is an open neighborhood of $N_{\ft(\ell),\fs}(\delta)\times\Si(\ft,\fs,\sP)$ in
$\tilde N_{\ft(\ell),\fs}(\delta)\times_{\sG_{\fs}}\bar\Gl(\ft,\fs,\sP)$ and where the $S^1$ actions are given by \eqref{eq:S1ActionOnPiece} and \eqref{eq:S1ActionOnOverlap},
\end{lem}

\begin{proof}
The
conclusion
follows immediately from
Item \eqref{item:ExistenceOfSplicedEndsModuli3} in Theorem \ref{thm:ExistenceOfSplicedEndsModuli}
that the splicing
maps $\bga'_{\Theta,\sP''_P}$ define open embeddings onto neighborhoods
of the strata \eqref{eq:TrivialStratumP} in $\barM^{s,\natural}_{\spl,|P|}(\delta_P)$.
\end{proof}

\subsection{Equality of splicing maps}
\label{subsec:Global_splicing_data_quality_splicing_maps}
We now show how the constructions of the overlap maps,
\eqref{eq:XUpwardsTransitionMap} and \eqref{eq:XDownwardTransition},
and crude splicing maps \eqref{eq:CrudeSplicingMapExplicit}
lead to a commutative diagram similar to the diagram
\eqref{eq:CommutingSplicingQuotient} used to construct the
spliced end $W_\ka$ defined in \eqref{eq:DefineSplicedEnd}.

For $\sP''\in [\sP<\sP']$, let
$\sO(\ft,\fs,\sP'')\subset \bar\Gl(\ft,\fs,\sP'')$
be the open subspace defined in \eqref{eq:DefineCrudeSplicingDomain} by
the separation condition.  Because this separation condition is invariant
under the action of $\fS(\sP)$ by Item \eqref{item:CompatFlatteningSymmetricEquiv} in
Lemma \ref{lem:CompatibleFlattening}, the collection of open subsets
 $\{\sO(\ft,\fs,\sP'')\}_{\sP''\in [\sP<\sP']}\}$
is $\fS(\sP)$-invariant.
Define
\begin{equation}
\label{eq:DefineSymmQuotientOfCrudeSplicingDomains}
\sO(\ft,\fs,[\sP<\sP'])
:=
\left.\left(\bigsqcup_{\sP''\in [\sP<\sP']} \sO(\ft,\fs,\sP'')\right)\right/\fS(\sP)
\subset
\bar\Gl(\ft,\fs,[\sP<\sP']).
\end{equation}
The $\fS(\sP)$-invariance of the map,
$$
\bigsqcup_{\sP''\in [\sP<\sP']}\bga''_{\ft,\fs,\sP''}:
\bigsqcup_{\sP''\in [\sP<\sP']}\tilde N_{\ft(\ell),\fs}(\delta)\times_{\sG_{\fs}} \sO(\ft,\fs,\sP'')
\to\bar\sC_{\ft},
$$
given by Lemma \ref{lem:CrudeSplicingSymmetricEquiv} then
implies that it defines a map on the $\fS(\sP)$-quotient,
\begin{equation}
\label{eq:UpperStratumCrudeSplicing}
\bga_{\ft,\fs,[\sP<\sP']}'':
\tilde N_{\ft(\ell),\fs}(\delta)\times_{\sG_{\fs}}\sO(\ft,\fs,[\sP<\sP'])
\to \bar\sC_{\ft}.
\end{equation}
Given open subspaces, $\sO(\ft,\fs,\sP)\subset \bar\Gl(\ft,\fs,\sP)$ as in \eqref{eq:DefineCrudeSplicingDomain},
and an $\fS(\sP)$-invariant collection of subspaces,
$\sO(\ft,\fs,\sP'')$ for $\sP''\in [\sP<\sP']$ as above, we
define an open subspace,
\begin{equation}
\label{eq:DefineOverlapDomain}
\sO_1(\ft,\fs,\sP,[\sP'])
\subset \bar\Gl(\ft,\fs,\sP,[\sP']),
\end{equation}
to be the
subset of
points satisfying the conditions
\eqref{eq:XUpwardsTransitionDomain1}, \eqref{eq:XUpwardsTransitionDeomainRequirement}, and,
$$
\sO_1(\ft,\fs,\sP,[\sP'])
\subseteqq
\left(\rho^{\ft,\fs,u}_{\sP,[\sP']}\right)^{-1} \left( \sO(\ft,[\sP<\sP'])
\right) \cap \left(\rho^{\ft,\fs,d}_{\sP,[\sP']}\right)^{-1} \left(
\sO(\ft,\fs,\sP) \right).
$$
We then have the

\begin{prop}
\label{prop:XOverlapControl}
Let $\sP<\sP'$ be partitions of
$N_\ell$. Assume that the families of metrics $g_{\sP''}$
on $X$
satisfy the
conditions in Section \ref{subsubsec:FiberBundleMetrics}.
Let $\sO(\ft,\fs,\sP)\subset\bar\Gl(\ft,\fs,\sP)$ and
$\sO(\ft,\fs,\sP'')\subset\bar\Gl(\ft,\fs,\sP'')$ for $\sP''\in [\sP<\sP']$
be the open
subspaces defined in \eqref{eq:DefineCrudeSplicingDomain}.
Let $\sO(\ft,\fs,[\sP<\sP'])$ and $\sO_1(\ft,\fs,\sP,[\sP'])$
be the
spaces defined in \eqref{eq:DefineSymmQuotientOfCrudeSplicingDomains}
and \eqref{eq:DefineOverlapDomain}, respectively.
Then the following diagram commutes,
\begin{equation}
\label{eq:GlobalSplicingCD}
\begin{CD}
\tilde
N_{\ft(\ell),\fs}(\delta)\times_{\sG_{\fs}}\sO_1(\ft,\fs,\sP,[\sP']) @>
\rho^{\ft,\fs,u}_{\sP,[\sP']} >> \tilde
N_{\ft(\ell),\fs}(\delta)\times_{\sG_{\fs}}\sO(\ft,\fs,[\sP<\sP'])
\\
@V \rho^{\ft,\fs,d}_{\sP,[\sP']} VV @V
\bga''_{\ft,\fs,[\sP<\sP']} VV
\\
\tilde N_{\ft(\ell),\fs}(\delta)\times_{\sG_{\fs}}\sO(\ft,\fs,\sP) @>
\bga''_{\ft,\fs,\sP} >> \bar\sC_{\ft}
\end{CD}
\end{equation}
and the maps are $S^1$-equivariant with respect to the action
\eqref{eq:S1ActionOnOverlap} with weight two on the domain of the overlap maps,
the action \eqref{eq:S1ActionOnPiece} with weight two on the domains of the crude splicing maps,
and the action \eqref{eq:S1ZAction} on $\bar\sC_{\ft}$.
\end{prop}

\begin{proof}
The conditions in the definition of the subspace
$\sO_1(\ft,\fs,\sP,[\sP'])$ imply that the compositions in the diagram \eqref{eq:GlobalSplicingCD},
\begin{equation}
\label{eq:CompositionsInSplicingDiagram}
\bga''_{\ft,\fs,[\sP<\sP']}\circ \rho^{\ft,\fs,u}_{\sP,[\sP']}
\quad\text{and}\quad
\bga''_{\ft,\fs,\sP}\circ(\rho^{\ft,\fs,d}_{\sP,[\sP']}),
\end{equation}
are defined on $\tilde
N_{\ft(\ell),\fs}(\delta)\times_{\sG_{\fs}}\sO_1(\ft,\fs,\sP,[\sP'])$.

The $S^1$-equivariance of the overlap maps
appears in Lemmas \ref{lem:UpwardsTransition}
and \ref{lem:XUpwardsTransition}.
The action \eqref{eq:S1ZAction} is just
scalar multiplication on the section.
Hence, the crude splicing maps will be equivariant
if $S^1$ acts on their domains by scalar multiplication
on the section of $V^+_0$ (where $\ft(\ell)=(\rho,V^\pm_0)$).
By Lemma \ref{lem:S^1ActionOnSplicingData}, the action in \eqref{eq:S1ActionOnPiece} with weight two
is scalar multiplication on the section of $V^+_0$, giving the desired equivariance.

To see that the diagram \eqref{eq:GlobalSplicingCD} commutes,
begin by restricting the compositions in \eqref{eq:CompositionsInSplicingDiagram} to
$$
\sG_{\fs}(A_0,\Phi_0)\times_{\sG_{\fs}}
\sO_1(\ft,\fs,\sP,[\sP'])|_{\pi_{\Si}^{-1}(\by)}
\subset
\tN_{\ft(\ell),\fs}(\delta)\times_{\sG_{\fs}} \sO(\ft,\fs,\sP,[\sP']),
$$
for $(A_0,\Phi_0)\in \tN_{\ft(\ell),\fs}(\delta)$ and $\by\in\Delta^\circ(X^\ell,\sP)$.
For $\bA'\in \sO_1(\ft,\fs,\sP,[\sP'])|_{\pi_{\Si}^{-1}(\by)}$,
we write
\begin{equation}
\label{eq:OverlapCompositionUpAndDownImages}
\begin{aligned}
(A''_u,\Phi''_u) &:=
(\bga''_{\ft,\fs,[\sP<\sP']}\circ \rho^{\ft,\fs,u}_{\sP,[\sP']})\left( [(A_0,\Phi_0),\bA']\right),
\\
(A''_d,\Phi''_d) &:=
(\bga''_{\ft,\fs,\sP}\circ\rho^{\ft,\fs,d}_{\sP,[\sP']})\left( [(A_0,\Phi_0),\bA']\right).
\end{aligned}
\end{equation}
Thus, if $[A''_u,\Phi''_u]=[A''_d,\Phi''_d]$, then the diagram \eqref{eq:GlobalSplicingCD} commutes.
The definitions
of the overlap maps
in \eqref{eq:UnparamUpwardsTransition} and \eqref{eq:XDownwardTransition}
and of the crude splicing maps
in \eqref{eq:CrudeSplicingConn} and \eqref{eq:CrudeSplicingMapExplicit}
imply that $A''_d$ and $A''_u$
are defined respectively by splicing in
connections on $S^4$
at the points in $X$ given by $(y_P)_{P\in\sP}=\by\in \Si(X,\sP)$  and at the points in $X$ given by $(x_Q)_{Q\in\sP'}=\by'\in\Si(X,\sP')$  as defined in \eqref{eq:UpwardTransitionBasePoints}.
Because $\pi(X^\ell,g_{\sP})(\by')=\by$  by \eqref{eq:UpwardTransitionBasePoints},
Item \eqref{item:CommutingExponentialMaps2} in Lemma \ref{lem:CommutingExponentialMaps}
and the consistency
of the families of metrics $g_{\sP}$ and
$g_{\sP'}$ imply that $g_{\sP,\by}=g_{\sP',\by'}$.
Similarly, the consistency condition \eqref{eq:ConsistentFlattening} for
the flattening maps constructed in Lemma \ref{lem:CompatibleFlattening} implies that
\[
(A'_0,\Phi'_0) = \bTheta_{\sP}(A_0,\Phi_0,\by) = \bTheta_{\sP'}(A_0,\Phi_0,\by').
\]
The definition in \eqref{eq:CrudeSplicingMapExplicit} and
\eqref{eq:CrudeSplicingFlattenedBackgroundPair} of the
section given by the crude splicing maps then imply that $\Phi''_u=\Phi_0=\Phi''_d$.

Property \eqref{eq:XUpwardsTransitionDomain1} of
$\sO(\ft,\fs,\sP,[\sP'])$ implies that $\by'$ lies in the
suitably small tubular neighborhood $\sU(X^\ell,g_{\sP})$ defined in
Lemma \ref{lem:CompatibleFlattening}. Thus, from
\eqref{eq:BallInclusions} the inclusions
$$
B(x_Q,4\tilde s_{P'}(\by')^{1/2}) \Subset B(y_P,2\tilde s_{P}(\by)^{1/2}),\quad\text{for $Q\subseteqq P$},
$$
hold for the balls where the splicing takes place.
The definition of the crude splicing map in \eqref{eq:CrudeSplicingConn} implies that
$A''_u$ and $A''_d$ both equal $A_0'$, and thus each other, on the complement of the balls
$B(y_P,2\tilde s_{P}(\by)^{1/2})$.
Item \eqref{item:BallOnWhichFlattenedConnIsFlat} of Lemma \ref{lem:CompatibleFlattening}
implies that the connection $\hat A_0'$ is flat on each ball $B(y_P,2\tilde s_{P}(\by)^{1/2})$
while Lemma \ref{lem:FlatteningMetrics2} implies that the metric $g_{\sP,\by}$ is flat on this set.
Because the metric and the connection $\hat A_0'$ are flat on the union of the balls
$B(y_P,2\tilde s_{P}(\by)^{1/2})$, the equality between the restrictions of $A''_u$ and $A''_d$
to these balls  follows from Lemma \ref{lem:CommutingSplicing}.
\end{proof}

\section{Construction of the space of global splicing data}
\label{sec:GlobalSplicingData}
The construction of the space of
global splicing data follows from  Proposition
\ref{prop:XOverlapControl} and the arguments on shrinking the
neighborhoods $\sO(\ft,\fs,\sP)$ given in the construction
of the spliced end $W_\ka$ in Section \ref{subsec:SplicedEnd}.

\begin{thm}
\label{thm:GlobalSplicingDataOverlaps}
Let $\ft$ be a \spinu structure on $X$ and $\fs$ a \spinc structure with $M_\fs\times\Sym^\ell(X)\subset I\sM_\ft$.
For a partition $\sP$  of $N_\ell$,
let $\Si(\ft,\fs,\sP)$ be the subspace defined in \eqref{eq:DefineTrivialStratumXConePoints} of the
gluing data space $\bar\Gl(\ft,\fs,\sP)$ defined in \eqref{eq:GluingData}.
For every partition $\sP$ of $N_\ell$, there is a neighborhood $\sO(\ft,\fs,\sP)$ of
$\Si(\ft,\fs,\sP)$ in   $\bar\Gl(\ft,\fs,\sP)$
such that
the intersection of the images of the crude splicing maps $\bga''_{\ft,\fs,\sP}$ defined in
\eqref{eq:CrudeSplicingMapExplicit},
\begin{equation}
\label{eq:CrudeSplicingMapImageIntersection1}
\bga''_{\ft,\fs,\sP}
\left(
N_{\ft(\ell),\fs}(\delta)\times_{\sG_{\fs}}\sO(\ft,\fs,\sP)
\right)
\cap
\bga''_{\sP'}
\left(
N_{\ft(\ell),\fs}(\delta)\times_{\sG_{\fs}}\sO(\ft,\fs,\sP')
\right),
\end{equation}
is empty if the sets $[\sP<\sP']$ and $[\sP'<\sP]$  defined in \eqref{eq:ConjugateRefinements}
are empty
while if $[\sP<\sP']$ is non-empty, then the intersection is contained in
\begin{align*}
{}& \bga''_{\ft,\fs,\sP}
\left(\rho^{\ft,\fs,d}_{\sP,[\sP']} \left( \tilde
N_{\ft(\ell),\fs}(\delta)\times_{\sG_{\fs}}\sO_1(\ft,\fs,\sP,[\sP']) \right)
\right)
\\
&\quad= \bga''_{\ft,\fs,[\sP<\sP']} \left( \rho^{\ft,\fs,u}_{\sP,[\sP']} \left(
\tilde N_{\ft(\ell),\fs}(\delta)\times_{\sG_{\fs}}\sO_1(\ft,\fs,\sP,[\sP'])
\right) \right),
\end{align*}
where $\sO_1(\ft,\fs,\sP,[\sP'])$ is the
subspace appearing in Proposition \ref{prop:XOverlapControl}
and defined in \eqref{eq:DefineOverlapDomain}.
If $\si\in\fS_\ell$ and $\sP'=\si(\sP)$ as defined by \eqref{eq:ActionOnPartitions},
then images of $\bga''_{\ft,\fs,\sP}$ and $\bga''_{\ft,\fs,\sP'}$  are equal.
\end{thm}

\begin{proof}
The proof follows the same argument as
that of Items \eqref{item:ConstructCollar3} and \eqref{item:ConstructCollar4} in
Lemma \ref{lem:ConstructCollar} for the construction of $W_\ka$.
The final sentence follows from Lemma \ref{lem:CrudeSplicingSymmetricEquiv}.
\end{proof}

We will require the following special case of
Theorem \ref{thm:GlobalSplicingDataOverlaps}.

\begin{lem}
\label{lem:SplicingTrivialOverlap}
Continue the notation of Theorem \ref{thm:GlobalSplicingDataOverlaps}.
Let $\sP<\sP'$ be partitions of $N_\ell$ and let
$T(\ft,\fs,\sP,[\sP'])\subset \bar\Gl(\ft,\fs,\sP,[\sP'])$ be as
defined in \eqref{eq:DefineTrivialInOverlap}.
There is a subspace $T^{\sO}(\ft,\fs,\sP,[\sP'])\subset T(\ft,\fs,\sP,[\sP'])$ defined by
the property
\begin{equation}
\label{eq:DefineSubspaceOfOverlapProductStrata}
\begin{aligned}
{}&\tilde N_{\ft(\ell),\fs}(\delta)\times_{\sG_{\fs}}T^{\sO}(\ft,\fs,\sP,[\sP'])
\\
{}&\quad=
\left(\tilde N_{\ft(\ell),\fs}(\delta)\times_{\sG_{\fs}}T(\ft,\fs,\sP,[\sP'])\right)
\cap
(\rho^{\ft,\fs,d}_{\sP,[\sP']})^{-1}\left(\tilde N_{\ft(\ell),\fs}(\delta)\times_{\sG_{\fs}}\sO(\ft,\fs,\sP)\right),
\end{aligned}
\end{equation}
where $\sO(\ft,\fs,\sP)$ is as defined in Theorem
\ref{thm:GlobalSplicingDataOverlaps}.  Then
$$
(\bga''_{\ft,\fs,\sP}\circ \rho^{\ft,\fs,d}_{\sP,[\sP']})
\left(N_{\ft(\ell),\fs}(\delta)\times T^{\sO}(\ft,\fs,\sP,[\sP']) \right)
\subset
\bga''_{\ft,\fs,\sP'}\left(N_{\ft(\ell),\fs}(\delta)\times \Si(\ft,\fs,\sP') \right).
$$
\end{lem}

\begin{proof}
To see that the property \eqref{eq:DefineSubspaceOfOverlapProductStrata} defines
a subspace of $T(\ft,\fs,\sP,[\sP'])$, we show that
for
$(A_0,\Phi_0)\in \tilde N_{\ft(\ell),\fs}(\delta)$ and
$\bA'\in T(\ft,\fs,\sP,[\sP'])$, the validity of the inclusion
$$
\rho^{\ft,\fs,d}_{\sP,[\sP']}\left([(A_0,\Phi_0),\bA']\right)
\in
\tilde N_{\ft(\ell),\fs}(\delta)\times_{\sG_{\fs}}\sO(\ft,\fs,\sP)
$$
does not depend on $(A_0,\Phi_0)\in \tilde N_{\ft(\ell),\fs}$.
Observe that
the map $\rho^{\ft,\fs,d}_{\sP,[\sP']}$ is defined in
\eqref{eq:XDownwardTransition} by the identity on $\tN_{\ft(\ell),\fs}(\delta)$
and by the map $\rho^{\ft,\fs,d}_{f,\sP,[\sP']}$ (defined in \eqref{eq:DowardsOverlapFiberMap})
on the fibers of $\pi_\Si$.
Because $\rho^{\ft,\fs,d}_{f,\sP,[\sP']}$ does not depend on $(A_0,\Phi_0)\in \tilde N_{\ft(\ell),\fs}(\delta)$,
the validity of the above inclusion does not depend on
$(A_0,\Phi_0)\in \tilde N_{\ft(\ell),\fs}(\delta)$, as required.

The definition of $T(\ft,\fs,\sP,[\sP'])$ and $\rho^{\ft,\fs,u}_{\sP,[\sP']}$ in \eqref{eq:UnparamUpwardsTransition}
imply that
$$
 \rho^{\ft,\fs,u}_{\sP,[\sP']}
\left(N_{\ft(\ell),\fs}(\delta)\times T^{\sO}(\ft,\fs,\sP,[\sP']) \right)
\subset
(N_{\ft(\ell),\fs}(\delta)\times \Si(\ft,\fs,\sP').
$$
The lemma then follows from the commutativity result in Proposition
\ref{prop:XOverlapControl}.
\end{proof}

We now define the \emph{space of global splicing data} to be
\begin{equation}
\label{eq:DefineGlobalGluingDataSpace}
\bar\sM^{\vir}_{\ft,\fs}
:=
\bigsqcup_{\sP} \left. \left( \tilde
N_{\ft(\ell),\fs}(\delta)\times_{\sG_{\fs}}\sO(\ft,\fs,\sP)\right)\right/\sim,
\end{equation}
where $\sO(\ft,\fs,\sP)$ is the space appearing in Theorem
\ref{thm:GlobalSplicingDataOverlaps} and points in
$$
\tilde N_{\ft(\ell),\fs}(\delta)\times_{\sG_{\fs}}\sO(\ft,\fs,\sP)
\quad\text{and}\quad
\tilde N_{\ft(\ell),\fs}(\delta)\times_{\sG_{\fs}}\sO(\ft,\fs,\sP')
$$
are identified by the relation $\sim$ if their images under the
crude splicing maps $\bga''_{\ft,\fs,\sP}$ and $\bga''_{\ft,\fs,\sP'}$
defined in \eqref{eq:CrudeSplicingMapExplicit}
are equal.
Therefore the space
$\bar\sM^{\vir}_{\ft,\fs}$ is, up to the action of $\fS_\ell$ on partitions, the homotopy pushout of the diagram
\eqref{eq:GlobalSplicingCD} in the sense of the following lemma.

\begin{lem}
\label{lem:IdentifyPointsInOverlap}
If $\sP<\sP'$ are partitions of $N_\ell$ and
\[
[A_0,\Phi_0,\bA_{\sP}]\in \tilde N_{\ft(\ell),\fs}(\delta)\times_{\sG_{\fs}}\sO(\ft,\fs,\sP)
\quad\text{and}\quad
[A_0',\Phi_0',\bA_{\sP'}]\in \tilde N_{\ft(\ell),\fs}(\delta)\times_{\sG_{\fs}}\sO(\ft,\fs,\sP'),
\]
then
\begin{equation}
\label{eq:IDentifiedPoints}
[A_0,\Phi_0,\bA_{\sP}] \sim [A_0',\Phi_0',\bA_{\sP'}]
\end{equation}
if and only if $[A_0,\Phi_0]=[A_0',\Phi_0']$ and there is
a point
$[A_0,\Phi_0,\bA']\in \tilde N\times_{\sG_{\fs}}\sO_1(\ft,\fs,\sP,[\sP'])$ satisfying
\begin{equation}
\label{eq:OverlapEquality1}
\begin{aligned}
{}[A_0,\Phi_0,\bA_{\sP}]
{}&= \rho^{\ft,\fs,d}_{\sP,[\sP']}([A_0,\Phi_0,\bA']),
\\
{}[A_0',\Phi_0',\bA_{\sP'}]
{}&= \rho^{\ft,\fs,u}_{\sP,[\sP']}([A_0,\Phi_0,\bA']).
\end{aligned}
\end{equation}
If $\si\in\fS_\ell$ and $\sP'=\si(\sP)$ as defined by \eqref{eq:ActionOnPartitions},
then \eqref{eq:IDentifiedPoints} holds if and only if
$[(A_0',\Phi_0'),\bA_{\sP'}]=[(A_0,\Phi_0),\si(\bA_{\sP})]$,
where $\si(\bA_{\sP})$ is defined by \eqref{eq:SymmetricGroupActionOnSplicingData1}.
\end{lem}

Recall from \cite[p. 43]{MayConciseCourse} that a pair of topological spaces
$(X,A)$ is a \emph{neighborhood deformation retract} (NDR) pair
if there are continuous maps $u:X\to [0,1]$ and $h:X\times [0,1]\to X$
with $A=u^{-1}(0)$, $h(x,0)=x$ for all $x\in X$, $h(a,t)=a$ for all $(a,t)\in A\times [0,1]$,
and $h(x,1)\in A$ if $u(x)<1$.

\begin{lem}
\label{lem:StatifiedSpaceStr}
The space $\bar\sM^{\vir}_{\ft,\fs}$ defined by
\eqref{eq:DefineGlobalGluingDataSpace} is a 
smoothly-stratified space with the property that if $\Si\subset\bar\sM^{\vir}_{\ft,\fs}$
is the complement of the top stratum, then
$(\bar\sM^{\vir}_{\ft,\fs},\Si)$ is an
NDR pair.
\end{lem}

\begin{proof}
By Theorem \ref{thm:GlobalSplicingDataOverlaps} and Lemma \ref{lem:IdentifyPointsInOverlap}, $\bar\sM^{\vir}_{\ft,\fs}$ is the union of open sets
\begin{equation}
\label{eq:OpenSetOfVirModuliGivenByPartition}
\tilde N_{\ft(\ell),\fs}(\delta)\times_{\sG_{\fs}}\sO(\ft,\fs,\sP),
\end{equation}
as $\sP$ varies over partitions, attached to each other by the smoothly-stratified maps,
\begin{equation}
\label{eq:TransitionMapsForVirModuli}
\begin{CD}
\rho^{\ft,\fs,d}_{\sP,[\sP']}
\left(
 \tilde N_{\ft(\ell),\fs}(\delta)\times_{\sG_{\fs}}\sO_1(\ft,\fs,\sP,[\sP'])
\right)
\\
@V \rho^{\ft,\fs,u}_{\sP,[\sP']}\circ (\rho^{\ft,\fs,d}_{\sP,[\sP']})^{-1} VV
\\
\rho^{\ft,\fs,u}_{\sP,[\sP']}
\left(
 \tilde N_{\ft(\ell),\fs}(\delta)\times_{\sG_{\fs}}\sO_1(\ft,\fs,\sP,[\sP'])
\right)
\end{CD}
\end{equation}
or, in the case when $\sP'=\si(\sP)$, by the smoothly-stratified diffeomorphism \eqref{eq:SymmetricGroupActionOnSplicingData1}.
Thus, $\bar\sM^{\vir}_{\ft,\fs}$ is a smoothly-stratified space.

Lemma \ref{lem:SplicedEndNDR} implies that for
$\sU=\tilde N_{\ft(\ell),\fs}(\delta)\times_{\sG_{\fs}}\sO(\ft,\fs,\sP)$,
the pair $(\sU,\Si\cap\sU)$ is  an NDR pair.
A pair $(U,\Si)$ is an NDR pair if and only if  the inclusion $i:\Si\to U$
is a cofibration by \cite[p. 43]{MayConciseCourse} or \cite[Theorem 7.1.10]{SelickIntroToHomotopyTheory}.
Thus, $(\bar\sM^{\vir}_{\ft,\fs},\Si)$ is a local cofibration in the sense that $\bar\sM^{\vir}_{\ft,\fs}$
admits a cover by open sets $\sU$ with the property that the inclusion $\Si\cap \sU\to \sU$ is a cofibration.
Because a local cofibration is a  global one by \cite[Satz 2]{Dold_1968},
the inclusion $\Si\to \bar\sM^{\vir}_{\ft,\fs}$ is a cofibration and hence an NDR pair.
\end{proof}

Because the transition maps
\eqref{eq:TransitionMapsForVirModuli} are $S^1$-equivariant by  Lemmas \ref{lem:UpwardsTransition} and \ref{lem:XUpwardsTransition}, the $S^1$ actions defined in \eqref{eq:S1ActionOnPiece}
on the open sets \eqref{eq:OpenSetOfVirModuliGivenByPartition} determine an $S^1$ action,
\begin{equation}
\label{eq:DefineGlobalS1Action}
S^1\times \bar\sM^{\vir}_{\ft,\fs}
\to \bar\sM^{\vir}_{\ft,\fs}.
\end{equation}
The description of the open cover in the preceding lemma then yields

\begin{cor}
\label{cor:GlobalFibration}
The space $\bar\sM^{\vir}_{\ft,\fs}$ admits an $S^1$-equivariant fibration,
\begin{equation}
\label{eq:GlobalProjectionToN}
\pi_N: \bar\sM^{\vir}_{\ft,\fs} \to
N_{\ft(\ell),\fs}(\delta),
\end{equation}
and an $S^1$-equivariant embedding,
\begin{equation}
\label{eq:GlobalCrudeSplicingMap}
\bga''_{\sM}:\bar\sM^{\vir}_{\ft,\fs}\to \bar\sC_{\ft},
\end{equation}
such that the restriction of $\bga''_{\sM}$ to $\tilde
N_{\ft(\ell),\fs}(\delta)\times_{\sG_{\fs}}\sO(\ft,\fs,\sP)$ is equal to
$\bga''_{\ft,\fs,\sP}$.
\end{cor}

\begin{proof}
Each of the open sets \eqref{eq:OpenSetOfVirModuliGivenByPartition} of $\bar\sM^{\vir}_{\ft,\fs}$
admits an $S^1$-equivariant fibration,
$$
\pi_{N,\sP}:\tilde N_{\ft(\ell),\fs}(\delta)\times_{\sG_{\fs}}\sO(\ft,\fs,\sP) \to N{\ft(\ell),\fs}(\delta).
$$
Because the transition maps \eqref{eq:TransitionMapsForVirModuli} are equal to the identity
on $\tilde N_{\ft(\ell),\fs}(\delta)$ by the definitions \eqref{eq:UnparamUpwardsTransition}
and \eqref{eq:XDownwardTransition}, these local fibrations define the global fibration
\eqref{eq:GlobalProjectionToN}.

Similarly, the embeddings $\bga''_{\ft,\fs,\sP}$ defined on the open sets
\eqref{eq:OpenSetOfVirModuliGivenByPartition} are equal on
the intersections of these open sets by the commutativity of the diagram \eqref{eq:GlobalSplicingCD}.
Hence, the maps $\bga''_{\ft,\fs,\sP}$ define the global map \eqref{eq:GlobalCrudeSplicingMap}.
The global map $\bga''_{\sM}$ is an embedding by the definition of the equivalence relation
following \eqref{eq:DefineGlobalGluingDataSpace}.

Finally, we see that the maps $\pi_N$ and $\bga''_{\sM}$ are $S^1$-equivariant by observing
that their local definitions are $S^1$-equivariant and that the $S^1$ action
\eqref{eq:DefineGlobalS1Action} is defined by the $S^1$ actions on the open sets \eqref{eq:OpenSetOfVirModuliGivenByPartition}.
\end{proof}

\section{Thom--Mather structures on the space of global splicing data}
\label{sec:TMStr}
On each open subspace,
\begin{equation}
\label{eq:DefineUSet}
\bar\sU(\ft,\fs,\sP)
:=
\tN_{\ft(\ell),\fs}(\delta)\times_{\sG_{\fs}}\sO(\ft,\fs,\sP) \subset
\barM^{\vir}_{\ft,\fs},
\end{equation}
there is a projection map,
\begin{equation}
\label{eq:DefineXTubularProj}
\pi(\ft,\fs,\sP): \bar\sU(\ft,\fs,\sP)
\to N_{\ft(\ell),\fs}(\delta)\times\Si(\ft,\fs,\sP),
\end{equation}
defined by $\id_N\times \pi_\Si$, where $\pi_\Si$ is the
restriction of the projection defined in
\eqref{eq:DefineSplicingDataToDiagonalProjection} to
$\sO(\ft,\fs,\sP)$ and $\id_N$ is the identity on $\tilde
N_{\ft(\ell),\fs}(\delta)$.

To prove that the
compatibility conditions \eqref{eq:TMProperties} hold for the
projection maps $\pi(\ft,\fs,\sP)$, we make the following
observations on the overlap maps.
For $\sP<\sP'$, define
\begin{equation}
\label{eq:DefineOverlapFiberAsImage}
U_f(\sP,[\sP'])
:=
\Imag
(\rho^{\ft,\fs,d}_{f,\sP,[\sP']})
\subset
\bar M(\sP),
\end{equation}
where the map on fibers $\rho^{\ft,\fs,d}_{f,\sP,[\sP']}$
is defined in \eqref{eq:DownwardsInclusionFiberMap}.
By the definition of $\rho^{\ft,\fs,d}_{f,\sP,[\sP']}$,
we have
$$
U_f(\sP,[\sP'])
\subseteqq
\bigsqcup_{\sP''\in [\sP<\sP'']}
\prod_{P\in\sP} \sU(\Theta,\sP''_P)/\fS(\sP),
$$
where $\sU(\Theta,\sP''_P)$ is the neighborhood of
$\{[\Theta]\}\times \Si(Z_P(\delta_P),\sP''_P)$
defined in Lemma \ref{lem:TMProjectionOnSplicedEndR4}
and $\fS(\sP)$ acts by permuting the components
of $U_f(\sP,[\sP'])$ in the disjoint union
given by $\sP''\in [\sP<\sP']$.
We use the convention that if $P\in\sP\cap\sP''$, then
$\sU(\Theta,\sP''_P)=\barM^{s,\natural}_{\spl,|P|}(\delta_P)$.

The projection maps,
\begin{align*}
\pi(\Theta,\sP'_P): \sU(\Theta,\sP'_P)\subset
\barM^{s,\natural}_{\spl,|P|}(\delta_P)
&\to \{[\Theta]\}\times
\Si(Z_P(\delta_P),\sP'_P)
\\
&\cong
\Delta^\circ(Z_P(\delta_P),\sP'_P)/\fS(\sP'_P),
\end{align*}
defined in \eqref{eq:SplicedEndR4Projection} are
$\SO(3)\times\SO(4)$-equivariant.
We use the convention that if $P\in\sP\cap \sP'$, then
$\pi(\Theta,\sP'_P)$ is projection to the cone point,
$\pi(\Theta,\sP'_P)([A,F^s,\bx])=[\Theta,F^s,c_P]$.
Together with the
identity maps
on $\tN_{\ft(\ell),\fs}(\delta)$ and on $\Fr(\ft,\fs,\sP)$,
the maps $\pi(\Theta,\sP''_P)$ (for
$\sP''\in [\sP<\sP']$) define a map,
$$
\begin{CD}
\tilde N_{\ft(\ell),\fs}(\delta)\times_{\sG_{\fs}\times S^1}
\Fr(\ft,\fs,\sP) \times_{G(\sP)} \prod\limits_{P\in\sP}
\bigsqcup\limits_{\sP''\in[\sP<\sP']} \sU(\Theta,\sP''_P)
\\
@V p_{\sP,[\sP']} VV
\\
\tilde N_{\ft(\ell),\fs}(\delta)\times_{\sG_{\fs}\times S^1}
\Fr(\ft,\fs,\sP) \times_{G(\sP)} \prod\limits_{P\in\sP}
\bigsqcup\limits_{\sP''\in[\sP<\sP']} \{[\Theta]\}\times
\Si(Z_P(\delta_P),\sP'_P)
\end{CD}
$$
Observe that in the definition of the preceding map,
the subgroup $\fS(\sP'_P)<\fS(\sP)<G(\sP)$ is
contained in
the stabilizer $\Ga(\fS)$ of the partition $\sP$
defined following \eqref{eq:ActionOnPartitions} and thus acts trivially on
$\Fr(\ft,\fs,\sP)$, giving the isomorphisms,
\begin{align*}
{}&
\tilde N_{\ft(\ell),\fs}(\delta)\times_{\sG_{\fs}\times S^1}
\Fr(\ft,\fs,\sP) \times_{G(\sP)} \prod_{P\in\sP}
\bigsqcup_{\sP''\in[\sP<\sP']} \{[\Theta]\}\times
\Delta^\circ(Z_P(\delta_P),\sP'_P)
\\
{}&\quad \cong
\tilde N_{\ft(\ell),\fs}(\delta)\times_{\sG_{\fs}\times S^1}
\Fr(\ft,\fs,\sP) \times_{G(\sP)} \prod_{P\in\sP}
\bigsqcup_{\sP''\in[\sP<\sP']} \{[\Theta]\}\times
\Delta^\circ(Z_P(\delta_P),\sP'_P)/\fS(\sP'_P)
\\
{}&\quad \cong
\tilde N_{\ft(\ell),\fs}(\delta)\times_{\sG_{\fs}\times S^1}
\Fr(\ft,\fs,\sP) \times_{G(\sP)} \prod_{P\in\sP}
\bigsqcup_{\sP''\in[\sP<\sP']} \{[\Theta]\}\times
\Si(Z_P(\delta_P),\sP'_P).
\end{align*}
Thus we can rewrite the range of the map $p_{\sP,[\sP']}$
to give
\begin{equation}
\label{eq:DefineUpperProjectionInLowerNgh}
\begin{CD}
\tilde N_{\ft(\ell),\fs}(\delta)\times_{\sG_{\fs}\times S^1}
\Fr(\ft,\fs,\sP) \times_{G(\sP)} \prod\limits_{P\in\sP}
\bigsqcup\limits_{\sP''\in[\sP<\sP']} \sU(\Theta,\sP''_P)
\\
@V p_{\sP,[\sP']} VV
\\
\tilde N_{\ft(\ell),\fs}(\delta)\times_{\sG_{\fs}\times S^1}
\Fr(\ft,\fs,\sP) \times_{G(\sP)} \prod\limits_{P\in\sP}
\bigsqcup\limits_{\sP''\in[\sP<\sP']} \{[\Theta]\}\times
\Delta^\circ(Z_P(\delta_P),\sP'_P)
\end{CD}
\end{equation}
The inclusions,
\begin{align*}
{}&\bigsqcup_{\sP''\in[\sP<\sP']} \sU(\Theta,\sP''_P) \subset
\barM^{s,\natural}_{\spl,|P|}(\delta_P) \quad \text{and} \quad
\\
{}& \bigsqcup_{\sP''\in[\sP<\sP']} \{[\Theta]\}\times
\Si(Z_P(\delta_P),\sP''_P) \subset
\barM^{s,\natural}_{\spl,|P|}(\delta_P),
\end{align*}
imply that
the map $p_{\sP,[\sP']}$
in \eqref{eq:DefineUpperProjectionInLowerNgh}
is defined on a subspace of
$\tilde N_{\ft(\ell),\fs}(\delta)\times_{\sG_{\fs}}\sO(\ft,\fs,\sP)$.
Observe that the relation,
\begin{equation}
\label{eq:ProjectionRelationOnLowerSet}
(\id_N\times \pi_\Si)\circ p_{\sP,[\sP']} = \id_N\times \pi_\Si,
\end{equation}
holds where $\pi_\Si$ is the projection given in
\eqref{eq:DefineSplicingDataToDiagonalProjection} because
$p_{\sP,[\sP']}$ respects the fibers of $\bar\Gl(\ft,\fs,\sP)$.

For $\sP<\sP'$, the following lemma describes the projection
map $\pi(\ft,\fs,\sP')$ on the image of the overlap map
$\rho^{\ft,\fs,d}_{\sP,[\sP']}$.

\begin{lem}
\label{lem:XProjectionOnOverlaps}
Let $\sP<\sP'$ be partitions of
$N_\ell$.
For the projection $\pi_{\Si,\sP,[\sP']}$ defined in \eqref{eq:ProjectionGluingOverlapToNormalInX}
and resulting map,
$$
\id_N\times \pi_{\Si,\sP,[\sP']}: \tilde
N_{\ft(\ell),\fs}(\delta)\times_{\sG_{\fs}}\bar\Gl(\ft,\fs,\sP,[\sP']) \to
N_{\ft(\ell),\fs}(\delta)\times T(\ft,\fs,\sP,[\sP']),
$$
the following relations hold on $\tN_{\ft(\ell),\fs}(\delta)\times_{\sG_{\fs}}\sO_1(\ft,\fs,\sP,[\sP'])$:
\begin{equation}
\label{eq:XProjectionRelationsOnOverlap}
\begin{aligned}
\pi(\ft,\fs,\sP')\circ \rho^{\ft,\fs,u}_{\sP,[\sP']} &=
\rho^{\ft,\fs,u}_{\sP,[\sP']}\circ (\id_N\times
\pi_{\Si,\sP,[\sP']}),
\\
p_{\sP,[\sP']}\circ
\rho^{\ft,\fs,d}_{\sP,[\sP']} &=
\rho^{\ft,\fs,d}_{\sP,[\sP']}\circ (\id_N\times
\pi_{\Si,\sP,[\sP']}).
\end{aligned}
\end{equation}
\end{lem}

\begin{proof}
The
conclusion
follows from the definition
of the projection maps in
\eqref{eq:DefineXTubularProj},
\eqref{eq:DefineUpperProjectionInLowerNgh},
and \eqref{eq:ProjectionGluingOverlapToNormalInX}
and the construction
of $\rho^{\ft,\fs,d}_{\sP,[\sP']}$ in \eqref{eq:XDownwardTransition}
and of $\rho^{\ft,\fs,u}_{\sP,[\sP']}$ in \eqref{eq:UnparamUpwardsTransition}.
\end{proof}

We now prove that the projection maps $\pi(\ft,\fs,\sP)$
satisfy the first compatibility condition in \eqref{eq:TMProperties}.

\begin{lem}
\label{lem:TM1ForGlobalSplicingData}
If $\sP<\sP'$ are partitions of $N_\ell$, then the following equality holds:
\begin{equation}
\label{eq:TM1ForGlobalSplicing}
\pi(\ft,\fs,\sP)\circ\pi(\ft,\fs,\sP')=\pi(\ft,\fs,\sP)
\quad\text{on }
\bar\sU(\ft,\fs,\sP)\cap \bar\sU(\ft,\fs,\sP').
\end{equation}
\end{lem}

\begin{proof}
By Theorem \ref{thm:GlobalSplicingDataOverlaps}, the intersection
$\bar\sU(\ft,\fs,\sP)\cap\bar\sU(\ft,\fs,\sP')$ is contained in the image
of the overlap maps in the diagram \eqref{eq:GlobalSplicingCD}.
Therefore, using the notation of Lemma \ref{lem:IdentifyPointsInOverlap},
every point $[A_0,\Phi_0,\bA_{\sP}]=[A_0,\Phi_0,\bA_{\sP'}]$ in the intersection is given by
$$
[A_0,\Phi_0,\bA_{\sP}]
= \rho^{\ft,\fs,d}_{\sP,[\sP']}([A_0,\Phi_0,\bA'])
=\rho^{\ft,\fs,u}_{\sP,[\sP']}([A_0,\Phi_0,\bA'])
=[A_0',\Phi_0',\bA_{\sP'}]
$$
for some $[A_0,\Phi_0,\bA']\in \tN_{\ft(\ell),\fs}(\delta)\times_{\sG_{\fs}}\sO_1(\ft,\fs,\sP,[\sP'])$. Therefore,
\begin{align*}
{}&\pi(\ft,\fs,\sP')([A_0',\Phi_0',\bA_{\sP'}])
\\
&\quad= \pi(\ft,\fs,\sP')\circ
\rho^{\ft,\fs,u}_{\sP,[\sP']}([A_0,\Phi_0,\bA'])
\\
&\quad= \rho^{\ft,\fs,u}_{\sP,[\sP']}\circ (\id_N\times
\pi_{\Si,\sP,[\sP']})([A_0,\Phi_0,\bA']) \quad\text{(by \eqref{eq:XProjectionRelationsOnOverlap})}
\\
&\quad= \rho^{\ft,\fs,d}_{\sP,[\sP']}\circ (\id_N\times
\pi_{\Si,\sP,[\sP']})([A_0,\Phi_0,\bA']) \quad\text{(by the definition of $\sim$
following \eqref{eq:DefineGlobalGluingDataSpace})}
\\
&\quad= p_{\sP,[\sP']}\circ \rho^{\ft,\fs,d}_{\sP,[\sP']}([A_0,\Phi_0,\bA']) \quad\text{(by \eqref{eq:XProjectionRelationsOnOverlap})}.
\end{align*}
The preceding equality implies that
\begin{align*}
{}& \pi(\ft,\fs,\sP)\circ \pi(\ft,\fs,\sP')([A_0',\Phi_0',\bA_{\sP'}])
\\
&\quad = \pi(\ft,\fs,\sP)\circ p_{\sP,[\sP']}\circ
\rho^{\ft,\fs,d}_{\sP,[\sP']}([A_0,\Phi_0,\bA'])
\\
&\quad = (\id_N\times\pi_\Si)\circ p_{\sP,[\sP']}\circ
\rho^{\ft,\fs,d}_{\sP,[\sP']}([A_0,\Phi_0,\bA']) \quad\text{(by the
definition of $\pi(\ft,\fs,\sP)$
following \eqref{eq:DefineXTubularProj})}
\\
&\quad= (\id_N\times\pi_\Si) \circ \rho^{\ft,\fs,d}_{\sP,[\sP']}([A_0,\Phi_0,\bA']) \quad\text{(by \eqref{eq:ProjectionRelationOnLowerSet})}
\\
&\quad = \pi(\ft,\fs,\sP)\circ \rho^{\ft,\fs,d}_{\sP,[\sP']}([A_0,\Phi_0,\bA']) \quad\text{(by the
definition of $\pi(\ft,\fs,\sP)$
following \eqref{eq:DefineXTubularProj})}
\\
&\quad = \pi(\ft,\fs,\sP)([A_0,\Phi_0,\bA_{\sP}]),
\end{align*}
and which, together with the equality $[A_0,\Phi_0,\bA_{\sP}]=[A_0',\Phi_0',\bA_{\sP'}]$,
proves \eqref{eq:TM1ForGlobalSplicing}.
\end{proof}

We now construct a tubular distance function on
$\bar\sU(\ft,\fs,\sP)$ in a manner similar to the
definition of $\vec t(X^\ell,g_{\sP})$
in \eqref{eq:DefineDiagonalTubularDistFunction}.
Because the functions $\tilde \la_P$
defined in Lemma
\ref{lem:TMTubularDistanceR4} are $\SO(4)\times\SO(3)$-invariant,
the map
$$
\prod_{P\in\sP}\tilde\la_P: \bar M(\sP) \to [0,1]^{\sP}:=\prod_{P\in\sP} [0,1],
$$
where $\bar M(\sP)$ is the fiber defined in \eqref{eq:DefineGluingDataFiber} of $\bar\Gl(\ft,\fs,\sP)$,
extends from the fiber and defines a smoothly-stratified map
\begin{equation}
\label{eq:DefineTubularDistanceFunction}
\vec t(\ft,\fs,\sP):\bar \sU(\ft,\fs,\sP)
\to
[0,1]^{\sP}/\fS(\sP).
\end{equation}
That is, $\vec t(\ft,\fs,\sP)$ is the unique map whose restriction to each fiber $\bar M(\sP)$ is
given by $\prod_{P\in\sP}\tilde\la_P$.
For $P\in\sP$, we
let $\vec t^{\,P}(\ft,\fs,\sP)$ denote the $P$-th component of the map \eqref{eq:DefineTubularDistanceFunction}.
We now describe the map $\vec t(\ft,\fs,\sP')$ on the intersection
$\bar\sU(\ft,\fs,\sP)\cap\bar\sU(\ft,\fs,\sP')$ where $\sP<\sP'$.

\begin{lem}
\label{lem:TubularDistanceFunctionOnOverlaps2}
Let $\sP<\sP'$ be partitions of $N_\ell$.
Let $\sO$ be the domain of the map $\rho^{\ft,\fs,d}_{f,\sP,[\sP']}$ defined
on the fibers of $\bar\Gl(\ft,\fs,\sP,[\sP'])$ in \eqref{eq:DowardsOverlapFiberMap}.
Then, for the structure group, $G(\sP)$, of the bundle $\Fr(\ft,\fs,\sP)$ defined in \eqref{eq:DefineGluingDataBundleStructureGroup_second_part},
there is a smoothly-stratified, $G(\sP)$-invariant map,
$$
\vec t_f(\sP,[\sP']):
\rho^{\ft,\fs,d}_{f,\sP,[\sP']}\left(\sO\right)
\subset
\prod_{P\in\sP} \barM^{s,\natural}_{\spl,|P|}(\delta_P)\to
\left.\bigsqcup_{\sP''\in [\sP<\sP']} [0,1]^{\sP''}\right/\fS(\sP),
$$
whose extension to $\bar\sU(\ft,\fs,\sP)\cap\bar\sU(\ft,\fs,\sP')$ equals $\vec t(\ft,\fs,\sP')$.
\end{lem}

\begin{proof}
We begin by defining a map,
$$
\vec t_u:
\tilde N_{\ft(\ell),\fs}(\delta)\times_{\sG_{\fs}\times S^1}
\bar\Gl(\ft,\fs,\sP,[\sP'])
\to
\left.\bigsqcup_{\sP''\in [\sP<\sP']} [0,1]^{\sP''}\right/\fS(\sP),
$$
as follows.
For $(A_0,\Phi_0)\in\tilde N_{\ft(\ell),\fs}(\delta)$, and $(F^{T,\fg}_P)_{P\in\sP}\in\Fr(\ft,\fs,\sP)$, and
$v_P\in \Delta^\circ(Z_{|P|}(\delta_P),\sP'_P)$,
and
$[A_Q,F^s_Z,\bx_Q]\in \barM^{s,\natural}_{\spl,|Q|}(\delta_Q)$, where
$Q\in\sP'_P$, we can write
$$
[(A_0,\Phi_0),\bA']=
\left[
(A_0,\Phi_0),\left(F^{T,\fg}_P,v_P,([A_Q,F^s_P,\bx_Q])_{Q\in\sP'_P}\right)_{P\in\sP}
\right]
\in
\tilde N_{\ft(\ell),\fs}(\delta)\times_{\sG_{\fs}\times S^1}\bar\Gl(\ft,\fs,\sP,[\sP']).
$$
Define
$$
\vec t_u([(A_0,\Phi_0),\bA']):= (\tilde \la_Q([A_Q,F^s_P,\bx_Q]))_{Q\in\sP'}.
$$
The definitions of $\vec t(\ft,\fs,\sP')$ and
$\rho^{\ft,\fs,u}_{\sP,[\sP']}$ (see \eqref{eq:DefineTubularDistanceFunction} and \eqref{eq:UnparamUpwardsTransition})
imply that
\begin{equation}
\label{eq:TubularFunctionPulledBackToOverlap}
\vec t(\ft,\fs,\sP')\circ \rho^{\ft,\fs,u}_{\sP,[\sP']}
= \vec t_u.
\end{equation}
Hence, the composition $\vec t(\ft,\fs,\sP')\circ \rho^{\ft,\fs,u}_{\sP,[\sP']}$
is defined by a $G(\sP)$-invariant
map on the fibers of $\bar\Gl(\ft,\fs,\sP,[\sP'])$.
The downwards transition map defined in \eqref{eq:XDownwardTransition}
is given by the $G(\sP)$-equivariant
map $\rho^{\ft,\fs,d}_{f,\sP,[\sP']}$
on the fibers of $\bar\Gl(\ft,\fs,\sP,[\sP'])$
defined in \eqref{eq:DowardsOverlapFiberMap}.
While the map $\rho^{\ft,\fs,d}_{f,\sP,[\sP']}$ in
\eqref{eq:DowardsOverlapFiberMap} is not an embedding
because of symmetric group actions, the map
\eqref{eq:TubularFunctionPulledBackToOverlap} is invariant under
these symmetric group actions and thus defines a
$G(\sP)$-invariant map,
$$
\vec t_{\sO,f}(\sP,[\sP']):
\rho^{\ft,\fs,d}_{f,\sP,[\sP']}(\sO)
\to
\left.\bigsqcup_{\sP''\in [\sP<\sP']} I^{\sP''}\right/\fS(\sP),
$$
where $\sO$ is the domain of $\rho^{\ft,\fs,d}_{f,\sP,[\sP']}$ appearing in the definition
\eqref{eq:DowardsOverlapFiberMap}.
From the definition of $\sO_1(\ft,\fs,\sP,[\sP'])$ in
\eqref{eq:DefineOverlapDomain}, in particular the condition
\eqref{eq:XUpwardsTransitionDeomainRequirement},
the domain of $\rho^{\ft,\fs,d}_{\sP,[\sP']}$ is contained in the union
over the fibers of the map $\pi_\Si$ defined in \eqref{eq:ProjectionGluingOverlapToDiagonalInX}.
Hence, the extension of the map $\vec t_{\sO,f}(\sP,[\sP'])$ from the fibers of $\pi_\Si$
defines a map on the domain of $\rho^{\ft,\fs,d}_{\sP,[\sP']}$.
Lemma \ref{lem:IdentifyPointsInOverlap} implies that
$\bar\sU(\ft,\fs,\sP)\cap\bar\sU(\ft,\fs,\sP')$
equals the image of $\rho^{\ft,\fs,d}_{\sP,[\sP']}$ and so this extension of
$\vec t_{\sO,f}(\sP,[\sP'])$ is the required map.
\end{proof}

The following lemma identifies
the restriction of $\vec t(\ft,\fs,\sP)$ to a subspace of $\bar\sU(\ft,\fs,\sP)$
with $\vec t(X^\ell,g_{\sP})$.

\begin{lem}
\label{lem:TubularDistanceFunctionOnOverlaps3}
Let $T(\ft,\fs,\sP,[\sP'])$ be the subspace defined in
\eqref{eq:DefineTrivialInOverlap} and define
$$
\sO_T(\ft,\fs,\sP,[\sP'])
:=
\sO_1(\ft,\fs,\sP,[\sP']) \cap T(\ft,\fs,\sP,[\sP']).
$$
If $k_{\Si,\sP,[\sP']}$ is the
smoothly-stratified diffeomorphism
\eqref{eq:IdentifyNormal} and $\id_N$ denotes the identity map on
$\tilde N_{\ft(\ell),\fs}(\delta)$, then there is a
smoothly-stratified embedding,
$$
\id_N\times (e(X^\ell,g_{\sP})\circ k_{\Si,\sP,[\sP']}):
\tilde N_{\ft(\ell),\fs}(\delta)\times_{\sG_{\fs}\times S^1}
\sO_T(\ft,\fs,\sP,[\sP'])
\to
N_{\ft(\ell),\fs}(\delta)/S^1\times \Si(X^\ell,\sP').
$$
Let $\pi_2: N_{\ft(\ell),\fs}(\delta)/S^1\times \Si(X^\ell,\sP')\to  \Si(X^\ell,\sP')$
be the projection.  Then  the restriction of the map
$\vec t(\ft,\fs,\sP)\circ \rho^{\ft,\fs,d}_{\sP,[\sP']}$ to
$$
\tilde N_{\ft(\ell),\fs}(\delta)\times_{\sG_{\fs}\times S^1} \sO_T(\ft,\fs,\sP,[\sP']),
$$
is equal to the composition,
$$
\vec t(X^\ell,\sP)\circ \pi_2\circ
\left(\id_N\times e(X^\ell,g_{\sP})\circ k_{\Si,\sP,[\sP']}\right).
$$
\end{lem}

\begin{proof}
The lemma follows immediately from the definition
of $\rho^{\ft,\fs,d}_{\sP,[\sP']}$
in terms of $\bga'_{\Theta,\sP''_P}$
in \eqref{eq:XDownwardTransition},
the
image
of $\bga'_{\Theta,\sP''_P}$ on the domain $\sO_T(\ft,\fs,\sP,[\sP'])$,
the definitions
of $\vec t(\ft,\fs,\sP)$ and $\vec t(X^\ell,\sP)$ in terms of
$\tilde \la_P$ and $t(Z_P)$, in  \eqref{eq:DefineTubularDistanceFunction}
and \eqref{eq:DefineDiagonalTubularDistFunction} respectively, and the equality in
Item \eqref{item:TMTubularDistanceR42} of
Lemma \ref{lem:TMTubularDistanceR4} between $\tilde \la_P$ and $t(Z_P)$.
\end{proof}

The following lemma shows that the Thom--Mather relation
$\vec t(\ft,\fs,\sP)\circ \pi(\ft,\fs,\sP')=\vec t(\ft,\fs,\sP)$
does not hold for the same reasons as the analogous relation
described in the proof of Lemma \ref{lem:SymProjTubDistRelations}.

\begin{lem}
\label{lem:SymProjTubDistRelations0}
Let $\sP<\sP'$ be partitions of $N_\ell$.
For each $\bA\in \bar\sU(\ft,\fs,\sP)\cap\bar\sU(\ft,\fs,\sP')$, the following identities hold:
\begin{enumerate}
\item
\label{item:SymProjTubDistRelations01}
If $P\notin\sP\cap \sP'$, then
$\vec t^{\,P}(\ft,\fs,\sP)(\bA)=\vec t^{\,P}(\ft,\fs,\sP)\circ \pi(\ft,\fs,\sP')(\bA)$.
\item
\label{item:SymProjTubDistRelations02}
If $P\in\sP\cap\sP'$, then $\vec t^{\,P}(\ft,\fs,\sP)\circ \pi(\ft,\fs,\sP')(\bA)=0$.
\item
\label{item:SymProjTubDistRelations03}
If $P\in\sP\cap\sP'$, then $\vec t^{\,P}(\ft,\fs,\sP)(\bA)=\vec t^{\,P}(\ft,\fs,\sP')(\bA)$.
\end{enumerate}
\end{lem}

\begin{proof}
By \eqref{eq:XProjectionRelationsOnOverlap},
we have the equality
$$
\vec t^{\,P}(\ft,\fs,\sP)\circ p_{\sP,[\sP']} \circ \rho^{\ft,\fs,d}_{\sP,[\sP']}
=
\vec t^{\,P}(\ft,\fs,\sP)\circ \pi(\ft,\fs,\sP')\circ \rho^{\ft,\fs,u}_{\sP,[\sP']}
$$
on $\bar\sU(\ft,\fs,\sP)\cap\bar\sU(\ft,\fs,\sP')$.  Thus, computing
$\vec t^{\,P}(\ft,\fs,\sP)\circ p_{\sP,[\sP']}$
will determine $\vec t^{\,P}(\ft,\fs,\sP)\circ \pi(\ft,\fs,\sP')$.
Item \eqref{item:TMTubularDistanceR43} of Lemma \ref{lem:TMTubularDistanceR4},
the definition of $\vec t^{\,}(\ft,\fs,\sP)$ in terms of $\tilde\la_P$ in \eqref{eq:DefineTubularDistanceFunction},
and
the definition of $p_{\sP,[\sP']}$
in \eqref{eq:DefineUpperProjectionInLowerNgh} in terms of
the projection maps $\pi(\Theta,\sP'_P)$ imply that
$\vec t^{\,P}(\ft,\fs,\sP)\circ p_{\sP,[\sP']}$ and $\vec t^{\,P}(\ft,\fs,\sP)$
will be equal for $P\notin \sP\cap\sP'$ while
$\vec t^{\,P}(\ft,\fs,\sP)\circ p_{\sP,[\sP']}$ will be zero for $P\in\sP\cap\sP'$.
This yields Items \eqref{item:SymProjTubDistRelations01} and \eqref{item:SymProjTubDistRelations02}.

Item \eqref{item:SymProjTubDistRelations03}
follows by observing that $\rho^{\ft,\fs,u}_{\sP,[\sP']}$
and $\rho^{\ft,\fs,d}_{\sP,[\sP']}$ are both given by the identity map
on the factor $\barM^{s,\natural}_{\spl,|P|}(\delta)$ of the fibers
\eqref{eq:DefineGluingDataFiber}
of
$\sO(\ft,\fs,\sP)$, and $\sO(\ft,\fs,\sP,[\sP'])$, and
$\sO(\ft,\fs,[\sP<\sP'])$.
\end{proof}

We
obtain a weaker version of the
second
Thom--Mather relation
from \eqref{eq:TMProperties}
in the following lemma and this will be sufficient to construct the
link of $M_\fs\times\Sym^\ell(X)$ in Section \ref{sec:AmbientLink}.

\begin{lem}
\label{lem:TubularDistanceFunctionOnOverlaps1}
Let $\sP<\sP'$ be partitions of $N_\ell$
and let $\eps>\eps'>0$.  Assume
that $\bA\in \bar\sU(\ft,\fs,\sP)\cap \bar\sU(\ft,\fs,\sP')$.
If $\hat D(\sP,\eps)$ denotes
either of the spaces $D(\sP,\eps)$ or $\bar D(\sP,\eps)$
defined in \eqref{eq:DefineSquares}, then the following hold:
\begin{enumerate}
\item
If $\vec t(\ft,\fs,\sP)(\bA)\in \hat D(\sP,\eps)$, then
$\vec t(\ft,\fs,\sP)\circ\pi(\ft,\fs,\sP')(\bA) \in \hat D(\sP,\eps)$.
\item
If $\vec t(\ft,\fs,\sP')(\bA)\in \bar D(\sP',\eps')$
and $\vec t(\ft,\fs,\sP)\circ \pi(\ft,\fs,\sP')(\bA)\in \hat D(\sP,\eps)$,
then $\vec t(\ft,\fs,\sP)(\bA)\in \hat D(\sP,\eps)$.
\end{enumerate}
\end{lem}

\begin{proof}
The
conclusion
follows from Lemma \ref{lem:TubularDistanceFunctionOnOverlaps1}
by the same arguments used to prove Lemma
\ref{lem:SymProjTubDistRelations}
from Lemma \ref{lem:SymProjTubDistRelations1}.
\end{proof}

In our construction of the link of $M_\fs\times\Sym^\ell(X)$ in
$\bar\sM^{\vir}_{\ft,\fs}/S^1$, the subspace
\begin{equation}
\label{eq:Disk1}
\vec t (\ft,\fs,\sP)^{-1}\left(\hat D(\sP,\eps)\right) \subset\bar \sU(\ft,\fs,\sP)
\end{equation}
will play the role of a disk-bundle neighborhood of $N_{\ft(\ell),\fs}(\delta)/S^1\times\Si(X^\ell,\sP)$ in $\bar\sU(\ft,\fs,\sP)$.  Similarly, for partitions $\sP<\sP'$ of $N_\ell$ and the tubular distance function,
$\vec t(X^\ell,g_{\sP}):\sU(X^\ell,\sP)\subset\Sym^\ell(X)\to [0,1]^{\sP}$,
defined in \eqref{eq:DefineDiagonalTubularDistFunction}, the subspace
\begin{equation}
\label{eq:Disk2}
\Si(X^\ell,\sP')
\cap
\left(\vec t(X^\ell,g_{\sP})^{-1}\left(\hat D(\sP,\eps)\right)\right)
\subset
\Si(X^\ell,\sP')
\cap
\sU(X^\ell,\sP)
\end{equation}
can be thought of as a disk-bundle neighborhood of $\Si(X^\ell,\sP)$ in $\Si(X^\ell,\sP')$.
The following corollary, to be used in
Proposition \ref{prop:LinkPieceFiberBundleStructure} to
establish the fiber-bundle structure of a subspace of the link,
describes the intersection of two  subspaces of the form \eqref{eq:Disk1}
in terms of the subspace \eqref{eq:Disk2}
and the projection map $\pi(\ft,\fs,\sP'):\bar\sU(\ft,\fs,\sP')\to N_{\ft(\ell),\fs}/S^1\times\Si(X^\ell,\sP')$.

\begin{cor}
\label{cor:LowerStratumRemoval}
Let $\sP<\sP'$ be partitions of $N_\ell$.
If
$$
D(X^\ell,\sP,[\sP'],\eps)
:=
N_{\ft(\ell),\fs}(\delta)/S^1
\times
\left( \Si(X^\ell,\sP') \cap
      \left(
        \vec t(X^\ell,g_{\sP})^{-1}\left(D(\sP,\eps)\right)
      \right)
\right),
$$
then for any $\eps'\in (0,\eps)$,
we have
\begin{align*}
{}&
\vec t(\ft,\fs,\sP')^{-1}\left(\bar D(\sP',\eps')\right)
\cap
\vec t(\ft,\fs,\sP)^{-1}\left(\bar D(\sP,\eps)\right)
\\
{}&\quad =
\vec t(\ft,\fs,\sP')^{-1}\left(\bar D(\sP',\eps')\right)
\cap
\pi(\ft,\fs,\sP')^{-1}\left(D(X^\ell,\sP,[\sP'],\eps)\right).
\end{align*}
\end{cor}

\begin{proof}
If
$$
T(X^\ell,\sP,[\sP'],\eps)
:=
\left( N_{\ft(\ell),\fs}(\delta)/S^1\times\Si(X^\ell,\sP')\right)
\cap
\vec t(\ft,\fs,\sP)^{-1}\left(\bar D(\sP,\eps))\right),
$$
then Lemma \ref{lem:TubularDistanceFunctionOnOverlaps1} implies that
\begin{align*}
{}&
\left(\vec t(\ft,\fs,\sP')^{-1}\left(\bar D(\sP',\eps')\right)\right)
\cap
\left(\vec t(\ft,\fs,\sP)^{-1}\left(\bar D(\sP,\eps)\right)\right)
\\
{}&\quad =
\left(\vec t(\ft,\fs,\sP')^{-1}\left(\bar D(\sP',\eps')\right)\right)
\cap
\pi(\ft,\fs,\sP)^{-1}\left( T(X^\ell,\sP,[\sP'],\eps)\right).
\end{align*}
The relation between $\vec t(\ft,\fs,\sP)$
and $\vec t(X^\ell,\sP)$ in Lemma \ref{lem:TubularDistanceFunctionOnOverlaps3}
then implies that
$$
T(X^\ell,\sP,[\sP'],\eps)=D(X^\ell,\sP,[\sP'],\eps),
$$
completing the proof of the corollary.
\end{proof}

\section{Global splicing map}
\label{sec:GlobalSplicingMap}
The global crude splicing map $\bga''_{\sM}$ in \eqref{eq:GlobalCrudeSplicingMap} is
unsatisfactory because its image does not include $M_{\fs}\times\Sym^\ell(X)$ --- see Lemma
\ref{lem:CrudeSplicingOnTrivials} and Remark \ref{rmk:LocalSplicingMap}.
To remedy this, we define an isotopy of the embedding $\bga''_{\sM}$ to an
embedding $\bga'_{\sM}$ with the desired image.
The analytical underpinnings required to replace Hypothesis \ref{hyp:Gluing} by a theorem that yields its assertions will be provided in the forthcoming \cite{Feehan_Leness_monopolegluingbook}.

For each partition $\sP$ of $N_\ell$,
the splicing map defined in Section \ref{subsubsec:StdSplicingMap} and the crude splicing map
defined in Section \ref{subsubsec:ConstrCrudeSplice}
give smoothly-stratified embeddings,
$$
\bga'_{\ft,\fs,\sP}:
\bar\sU(\ft,\fs,\sP)
\to \bar\sC_{\ft}
\quad\text{and}\quad
\bga''_{\ft,\fs,\sP}:
\bar\sU(\ft,\fs,\sP)
\to \bar\sC_{\ft},
$$
respectively,
where $\bar\sC_\ft$ is defined in \eqref{eq:IdealPairs}.
By its definition in \eqref{eq:GlobalCrudeSplicingMap}, the restriction of the global crude
splicing map to $\bar\sU(\ft,\fs,\sP)$ equals $\bga''_{\ft,\fs,\sP}$.
From the list of differences between the two splicing maps in the beginning of Section \ref{subsubsec:ConstrCrudeSplice}, one can see that the map $\bga''_{\ft,\fs,\sP}$ equals the map
defined by $\bga'_{\ft,\fs,\sP}$ using a locally flattened Riemannian metric $g_{\sP,\bx}$
(where $\bx\in\Si(X^\ell,\sP)$)
and by replacing the pair $(A_0,\Phi_0)\in\tsC_{\ft(\ell)}$ with the locally flattened pair
$\bTheta(A_0,\Phi_0,\bx)$.
Thus, one can define an isotopy of smoothly-stratified embeddings between
$\bga'_{\ft,\fs,\sP}$ and $\bga''_{\ft,\fs,\sP}$ by applying the construction of
$\bga'_{\ft,\fs,\sP}$ with a convex linear combination of the Riemannian metrics $g$
and $g_{\sP,\bx}$ and of the pairs $(A_0,\Phi_0)$ and $\bTheta(A_0,\Phi_0,\bx)$.

We write this isotopy as
$$
\Ga'_{\ft,\fs,\sP}:\bar\sU(\ft,\fs,\sP)\times [0,1]\to \bar\sC_{\ft},
\quad\text{with}\quad
\Ga'_{\ft,\fs,\sP}(\cdot,s)
=
\begin{cases}
\bga''_{\ft,\fs,\sP} &\text{if } s = 0,
\\
\bga'_{\ft,\fs,\sP} &\text{if } s = 1.
\end{cases}
$$
We now construct the function giving the parameter to give the isotopy.
Let $f_\sP:\Delta^\circ(X^\ell,\sP)\to (0,1]$ be a smooth, $\fS(\sP)$-invariant function
with the property that the rescaled map
$$
\vec t_r(\ft,\fs,\sP)(\cdot)
:=
\vec t (\ft,\fs,\sP)(\cdot) (f_{\sP}\circ \pi(\ft,\fs,\sP))(\cdot),
$$
(where the preceding multiplication means  to multiply each factor
of $\vec t (\ft,\fs,\sP)$ by the scalar
$f_{\sP}\circ \pi(\ft,\fs,\sP)$) takes values in
the subspace $\bar D(\sP,\frac{1}{2})\subset [0,1]^{\sP}$ defined in \eqref{eq:DefineSquares}.
Now, observe that if
$$
\beta: [0,1]^\sP\to [0,1],
$$
is a smooth, $\fS(\sP)$-invariant function, supported in
$D(\sP,\frac{1}{2})$
with $\beta(\bar D(\sP,\frac{1}{4}))=1$, then
$\beta\circ \vec t_r(\ft,\fs,\sP)$
is identically equal to one
on $\bar\sU_1(\ft,\fs,\sP)$ and
is identically equal to zero
on $\bar\sU(\ft,\fs,\sP)\setminus\sU_2(\ft,\fs,\sP)$, where the spaces
$\bar\sU_i(\ft,\fs,\sP)$,
for $i=1,2$,
satisfy
$$
N_{\ft(\ell),\fs}(\delta)/S^1\times \Si(X^\ell,\sP)
\sqsubset
\bar\sU_1(\ft,\fs,\sP)\sqsubset \bar\sU_2(\ft,\fs,\sP)\sqsubset \bar\sU(\ft,\fs,\sP).
$$
Hence, the function,
\begin{equation}
\label{eq:DefineCutOff}
\beta_\sP = \beta\circ \vec t_r(\ft,\fs,\sP),
\end{equation}
defines
a global function on $\bar\sM^{\vir}_{\ft,\fs}$
satisfying
\begin{equation}
\label{eq:SupportOfCutOffs}
\bar\sM^{\vir}_{\ft,\fs}\setminus\bar\sU_2(\ft,\fs,\sP) \subset
\beta_{\sP}^{-1}(0) \quad\text{and}\quad \bar\sU_1(\ft,\fs,\sP)
\subset \beta_{\sP}^{-1}(1).
\end{equation}
We now describe how to
patch
the isotopies $\Ga'_{\ft,\fs,\sP}$
together to define a global isotopy of smoothly-stratified embeddings.
Choose one representative $\sP$ in each conjugacy class $[\sP]$ and
enumerate these representatives in the manner
described in Section \ref{subsec:EnumStrata}.  Assume $\sP_0=\{N_\ell\}$ is the crudest
partition and
$\sP_m$ is the most refined partition.
We define a finite sequence of
smoothly-stratified
embeddings,
$$
\bga'_{\sM,i}:\bar\sM^{\vir}_{\ft,\fs}\to \bar\sC_{\ft},
$$
as follows.  We begin with
\begin{multline}
\label{eq:SplicingIsotopyInduction1}
\bga'_{\sM,0}([(A_0,\Phi_0),\bA])
\\
:=
\begin{cases}
\Ga'_{\ft,\fs,\sP_0}([(A_0,\Phi_0),\bA],\beta_{\sP_0}([(A_0,\Phi_0),\bA]))
{}&
\quad\text{if $[(A_0,\Phi_0),\bA]\in\bar\sU(\ft,\fs,\sP_0)$,}
\\
\bga''_{\sM}([(A_0,\Phi_0),\bA])
{}&
\quad\text{if $[(A_0,\Phi_0),\bA]\notin\bar\sU(\ft,\fs,\sP_0)$.}
\end{cases}
\end{multline}
Observe that the restriction of $\bga'_{\sM,0}$ to $\bar\sU(\ft,\fs,\sP_0)$
is isotopic to $\bga_{\ft,\fs,\sP_0}'$.
We proceed by upwards induction on the enumerated partitions.
Assume inductively that the restriction  of $\bga'_{\sM,i-1}$ to
$\bar\sU(\ft,\fs,\sP_i)$ is isotopic to $\bga''_{\ft,\fs,\sP_i}$.
Then, the construction of the isotopy $\Ga'_{\ft,\fs,\sP_i}$ above
can be modified to give an isotopy of smoothly-stratified embeddings
of $\bar\sU(\ft,\fs,\sP_i)$,
$$
\hat\Ga'_{\ft,\fs,\sP_i}:\bar\sU(\ft,\fs,\sP_i)\times [0,1]\to \bar\sC_{\ft},
\quad\text{with}\quad
\hat\Ga'_{\ft,\fs,\sP_i}(\cdot,s)
=
\begin{cases}
\bga'_{\sM,i-1} &\text{if } s = 0,
\\
\bga'_{\ft,\fs,\sP} &\text{if } s = 1.
\end{cases}
$$
We then continue the induction by defining
\begin{multline}
\label{eq:SplicingIsotopyInduction2}
\bga'_{\sM,i}([(A_0,\Phi_0),\bA])
\\
:=
\begin{cases}
\hat\Ga'_{\ft,\fs,\sP_i}([(A_0,\Phi_0),\bA],\beta_{\sP_i}([(A_0,\Phi_0),\bA]))
{}&
\quad\text{if } [(A_0,\Phi_0),\bA]\in\bar\sU(\ft,\fs,\sP_i),
\\
\bga'_{\sM,i-1}([(A_0,\Phi_0),\bA])
{}&
\quad\text{if } [(A_0,\Phi_0),\bA]\notin\bar\sU(\ft,\fs,\sP_i).
\end{cases}
\end{multline}
If $\sP_m$ is, as defined above, the most refined partition,
we define the \emph{global splicing map} by
\begin{equation}
\label{eq:DefineGlobalSplice}
\bga'_{\sM} := \bga'_{\sM,m}.
\end{equation}
The proof of the following proposition is then technical but straightforward.

\begin{prop}
\label{prop:GlobalSplicingMap}
There is a smoothly-stratified, $S^1$-equivariant embedding,
$$
\bga'_{\sM}:\bar\sM^{\vir}_{\ft,\fs}\to \bar\sC_\ft,
$$
such that the restriction of $\bga'_{\sM}$ to $N_{\ft(\ell),\fs}(\delta)\times\Sym^\ell(X)$ is equal to the
product of the embedding $N_{\ft(\ell),\fs}(\delta)\to\sC_{\ft(\ell)}$ with the
identity on $\Sym^\ell(X)$.
\end{prop}

\section{Projections onto symmetric products}
\label{sec:ProjectionToSym}
We now define a projection map,
\begin{equation}
\label{eq:GlobalProjectionToX}
\pi_X: \bar\sM^{\vir}_{\ft,\fs} \to
\Sym^\ell(X),
\end{equation}
which will be used in describing cohomology classes on
$\bar\sM^{\vir}_{\ft,\fs}$.

We will denote the projection $\pi_\Si$ in
\eqref{eq:DefineSplicingDataToDiagonalProjection} by
$$
\pi_\Si(\ft,\fs,\sP):\bar{\sU}(\ft,\fs,\sP)\to \Si(X^\ell,\sP),
$$
to avoid confusion between the projection maps of the different
strata.  These local projections maps do not agree on the
intersections $\bar\sU(\ft,\fs,\sP)\cap\bar\sU(\ft,\fs,\sP')$, but they can be
related with the aid of the following lemma.

\begin{lem}
\label{lem:LocalXProjectionsRelation}
If $\sP<\sP'$ are partitions
of $N_\ell$, then
\begin{equation}
\label{eq:LocalXProjectionsRelation}
\pi(X^\ell,g_{\sP})\circ \pi_\Si(\ft,\fs,\sP')=\pi_\Si(\ft,\fs,\sP)
\quad\text{on } \bar\sU(\ft,\fs,\sP)\cap \bar\sU(\ft,\fs,\sP').
\end{equation}
\end{lem}

\begin{proof}
The
conclusion
follows immediately from Lemma
\ref{lem:IdentifyPointsInOverlap} and the definitions of the
projection maps.
\end{proof}

We can now define the desired global projection to $\Sym^\ell(X)$.

\begin{lem}
\label{lem:GlobalProjectionToX}
There is an open neighborhood
$\sU_1\subset \bar\sM^{\vir}_{\ft,\fs}$ of
$N_{\ft(\ell),\fs}(\delta)\times\Sym^\ell(X)$ and a projection map,
$$
\pi_X:\sU_1\to \Sym^\ell(X),
$$
such that the restriction of $\pi_X$ to
$\sU_1\cap\bar\sU(\ft,\fs,\sP)$ is homotopic to
$\pi_\Si(\ft,\fs,\sP)$.
\end{lem}

\begin{proof}
Enumerate the strata of $\Sym^\ell(X)$ as described in
Section \ref{subsec:EnumStrata} using partitions $\sP_0,\dots,\sP_n$.
For $0\le i\le n$,
choose open neighborhoods $\sU_1(\ft,\fs,\sP_i)\subset\bar\sU(\ft,\fs,\sP_i)$
of $\Si(\ft,\fs,\sP_i)$ such that, for $j<k$,
$$
\pi(\ft,\fs,\sP_k)(\bar\sU_1(\ft,\fs,\sP_j)\cap \bar\sU_1(\ft,\fs,\sP_k))
\subseteqq
\bigcup_{i<j} T_i,
$$
where $T_i$ is constructed in Lemma \ref{lem:DecomposeSymProduct}.
The
conclusion
then follows from Lemma \ref{lem:GlobalProjToSymCrit}
and the equality \eqref{eq:LocalXProjectionsRelation}.
\end{proof}

\chapter{Obstruction bundle}
\label{chap:obstr}

\section{Introduction}
\label{sec:Introduction_obstruction_bundle}
The gluing map described in Hypothesis \ref{hyp:Gluing}
does not map all points in $\bar\sM^{\vir}_{\ft,\fs}$ to the
moduli space $\bar\sM_{\ft}$.
Instead, the intersection of $\bar\sM_{\ft}$ with the image of the top stratum  of $\bar\sM^{\vir}_{\ft,\fs}$
is diffeomorphic to the zero-locus of a section  of a vector bundle
over the top stratum of $\bar\sM^{\vir}_{\ft,\fs}$.
To describe the intersection of the image of  $\bar\sM^{\vir}_{\ft,\fs}$ under the gluing map
with $\bar\sM_{\ft}$, we require the following notion of a pseudo-bundle.

\begin{defn}
\label{defn:PseudoBundle}
(See Schwartz \cite{Schwartz_Lectures_On_Strat_Objects}.)
Let $Y$ be a stratified space.  Assume that to every stratum
$\Si$ of $Y$, there is an open neighborhood, $U_\Si$, of $\Si$
in $Y$ and a strict deformation retraction, $r_\Si:U_\Si\to \Si$.
A {\em pseudo-bundle\/} over $Y$ is a space $\Upsilon$ with
a surjective, continuous map $\pi_\Upsilon:\Upsilon\to Y$
satisfying the following conditions:
\begin{enumerate}
\item
For every stratum $\Si$ of $Y$, the subspace $\Upsilon_\Si=\pi_\Upsilon^{-1}(\Si)$
is a vector bundle with projection map given by the restriction of $\pi_\Upsilon$.
\item
For every stratum $\Si$ of $Y$, there is an injective  bundle map,
$$
\iota_{\Upsilon,\Si}: r_\Si^*\Upsilon_\Si \to \Upsilon|_{U_\Si},
$$
that is linear on each fiber
and is the identity when restricted to $\Upsilon_\Si$ and
which admits a left inverse,
$$
\Pi_{\Upsilon,\Si}:\Upsilon|_{U_\Si}\to r_\Si^*\Upsilon_\Si,
$$
that is a surjective vector bundle map.
\end{enumerate}
A \emph{section} of a pseudo-bundle $\pi_\Upsilon:\Upsilon\to Y$ is a continuous map
$s:Y\to \Upsilon$ satisfying $\pi_\Upsilon\circ s=\id_Y$.
\end{defn}
In this chapter, we construct a pseudo-bundle
(the \emph{obstruction pseudo-bundle}) over $\bar\sM^{\vir}_{\ft,\fs}$
and then in Hypothesis \ref{hyp:Gluing}
state
the properties of a gluing map (to be proved in the forthcoming \cite{Feehan_Leness_monopolegluingbook}) required to
give a smoothly-stratified
diffeomorphism between the zero-locus
of a section of this obstruction pseudo-bundle
(the \emph{obstruction section})
and an open neighborhood of $M_{\fs}\times\Sym^\ell(X)$
in $\bar\sM_{\ft}$.

For a point $[A',\Phi']$ in the image of the splicing map,
let $D\fS|_{(A',\Phi')}$ denote the linearization of
the $\SO(3)$-monopole equations
at $(A',\Phi')$ with image
in
\begin{equation}
\label{eq:MonopoleEquationsImageSpace}
L^2_{k-1}(\La^+\otimes \fg_{\ft}) \oplus L^2_{k-1}(V^-).
\end{equation}
The fiber of the obstruction pseudo-bundle over $[A',\Phi']$ will be defined by
assigning, $\sG_\ft$-equivariantly, to each pair $(A'',\Phi'')$ in the $\sG_\ft$-orbit
a finite-rank subspace $\Upsilon_{(A'',\Phi'')}$ of the Hilbert space \eqref{eq:MonopoleEquationsImageSpace}
which is close to the cokernel of $D\fS|_{(A'',\Phi'')}$
in the sense that the restriction of
the $L^2$-orthogonal projection onto $\Coker(D\fS|_{(A'',\Phi'')})$ to
$\Upsilon_{(A'',\Phi'')}$ is $L^2$-close to the identity
(see \cite[Proposition 8.19]{FL3}).

Let $[A',\Phi']$ be a point in the image of the splicing map, given by splicing
$\SO(3)$ connections
$A_P$ over $S^4$,
for $P\in\sP$, onto
a background pair $(A_0,\Phi_0)\in\tilde N_{\ft(\ell),\fs}(\delta)$.
In \cite[Section 8]{FL3}, we constructed an embedding
of vector spaces,
$$
\varphi_{\sP}:\Coker(D\fS|_{(A_0,\Phi_0)})
\oplus
\left(
\mathop{\oplus}\limits_{P\in\sP}
\Coker(D\fS|_{A_P,0)})
\right)
\to
L^2_{k-1}(\La^+\otimes \fg_{\ft}) \oplus L^2_{k-1}(V^-),
$$
by multiplying elements of the domain by cut-off functions.
We then showed that $\Imag(\varphi_{\sP})$ had the desired
properties mentioned above and hence
will define the fiber of the obstruction pseudo-bundle over
$[A',\Phi']$.  In this section, we show
how to fit these fibers together as the point
$[A',\Phi']$ and the partition $\sP$ vary to form
a pseudo-bundle.

We refer to the component of $\Imag(\varphi_{\sP})$ given by
$\varphi_{\sP}(\Coker(D\fS|_{(A_0,\Phi_0)}))$ as the {\em background
component\/}.  If $M_{\fs}$ has positive dimension, then
the possibility of spectral flow and the
varying
dimension of the cokernel of $D\fS|_{(A_0,\Phi_0)}$
prevent us from using this cokernel to form a vector bundle.
Instead, we use the obstruction bundle $\Xi_{\ft(\ell),\fs}\to M_{\fs}$
described in Section \ref{subsubsec:ThickenedNeighborhood}.
This obstruction bundle is a stabilization of the cokernel
of $D\fS$ in the sense that it is a trivial bundle which
surjects onto $\Coker(D\fS|_{(A_0,\Phi_0)})$ at each point
$[A_0,\Phi_0]\in M_{\fs}$.
For $[A',\Phi']\in\bar\sU(\ft,\fs,\sP)\cap\bar\sU(\ft,\fs,\sP')$, the
restrictions of the maps $\varphi_{\sP}$ and $\varphi_{\sP'}$
to  $\Xi_{\ft(\ell),\fs}$ differ in the choice of cut-off functions
just as the standard splicing maps $\bga_{\ft,\fs,\sP}'$ differed (see
\eqref{eq:StandardCutOffSection}).
To fit these images together to form a global bundle,
we use the flattening map for pairs, defined in Section \ref{sec:FlatteningPairs},
relying on the consistency property \eqref{eq:ConsistentFlattening}.
The details of this construction appear in Section \ref{sec:SWObstruction}.

We refer to the component of $\Imag(\varphi_{\sP})$ given by
$\varphi_{\sP}(\oplus_{P\in\sP} \Coker(D\fS|_{(A_P,0)}))$ as the {\em instanton component\/}.
The cokernel  of $D\fS|_{(A_P,0)}$
is isomorphic to the cokernel of the twisted
Dirac operator $D_{A_P}:\Om^0(V_{|P|}^+)\to  \Om^0(V_{|P|}^-)$,
where $V_{|P|}$ is the \spinu structure on $S^4$ defined by
the standard spin structure and the rank-two complex bundle
$E_{|P|}\to S^4$.
An argument using the Bochner--Weitzenb\"ock formula for the Dirac operator and the fact that the standard round metric on $S^4$ has positive scalar curvature
(see \cite[Section 4.1]{FL1})
ensures that spectral flow
does not arise
for this
component of the obstruction and that the cokernel of the twisted
Dirac operator forms a vector bundle (the index
bundle) over $M^{s,\natural}_{|P|}(\delta_P)$, as
described in Section \ref{sec:EquivDiracBundle}.  However,
the index of the twisted Dirac operator drops on
the lower-charge moduli spaces
appearing in the lower strata of the
Uhlenbeck compactification
of the moduli space of anti-self-dual connections over $S^4$.
Note that this is a change in the index and not just a
spectral flow problem which could be dealt with by the
stabilization procedure used for the background obstruction.
Hence the image space for the obstruction map
will  only be a pseudo-bundle in the sense of
Definition \ref{defn:PseudoBundle}.  This
change of index
reflects
the failure of $\bar\sM_{\ft}$ to intersect
the lower strata of $\bar\sM^{\vir}_{\ft,\fs}$ in
subspaces of the same codimension.  One can express
this failure in the language of intersection homology \cite{GoreskyMacPhersonIntHomology,GoreskyMacPhersonIntHomologyII,Borel_IntCohom}
by saying that the cycle $\bar\sM_{\ft}$ has non-trivial
perversity.
However, in Section \ref{subsec:RelEulerClasses}, we shall
see that for computations in rational cohomology,
the Euler class of
the obstruction pseudo-bundle behaves in a manner
similar to that of an actual vector bundle and there is
no need to
appeal to the framework
of intersection-homology-valued characteristic classes.

The obstruction to forming a global
obstruction pseudo-bundle from the images
of $\varphi_{\sP}$, restricted to the instanton component
comes from comparing the images of $\varphi_{\sP}$ and $\varphi_{\sP'}$
for different partitions.  We overcome this problem by
adapting the method used in Chapters \ref{chap:SplicedEnd} and \ref{chap:GlobalSplicingData} to compare
the images of the crude splicing maps $\bga_{\ft,\fs,\sP}''$.
In Section \ref{sec:ObstructionSplicing},
we define the {\em spliced-ends obstruction pseudo-bundle\/},
\begin{equation}
\label{eq:SplicedEndsObstr1}
\bar\Upsilon^i_{\spl,\ka}\to\barM^{s,\natural}_{\spl,\ka}(\delta),
\end{equation}
inductively by assuming that the bundle \eqref{eq:SplicedEndsObstr1} exists for $\ka'<\ka$
and
by replacing the index bundle of the twisted Dirac operator
over the subspace $\sU(\Theta,\sP'_P)$ of the
spliced end of $\barM^{s,\natural}_{\spl,\ka}(\delta)$
with the image, under a splicing map similar to $\varphi_{\sP'_P}$,
of
a direct sum of bundles,
$$
\bigoplus_{P\in\sP}\bar\Upsilon^i_{\spl,|P|}\ .
$$
Because this construction is equivariant with respect to
group actions discussed in Section \ref{sec:SpinuAction},
the spliced-ends obstruction pseudo-bundle extends to a pseudo-bundle
over $\bar\sU(\ft,\fs,\sP)$ with an embedding
$\varphi_{\ft,\fs,\sP}$
replacing $\varphi_{\sP}$.
In Section \ref{sec:InstantonObstr}, we
give a proof,
similar to that in Section \ref{sec:OverlapSpacesModuli},
that the intersections of
the images of splicing maps $\varphi_{\ft,\fs,\sP}$ and
$\varphi_{\ft,\fs,\sP'}$
can be
described
by a push-out diagram, thus defining
the instanton obstruction as a global pseudo-bundle.
Finally, in Section \ref{sec:GluingThm}, we state the
properties of the local gluing map and obstruction section required to prove the main results of this 
monograph.

\section{Infinite-rank obstruction pseudo-bundle}
\label{sec:InfiniteDimObstr}
The monopole equations \eqref{eq:PerturbedSO3MonopoleEquations}
define a section, $\fS$,  of
the infinite-rank obstruction pseudo-bundle,
\begin{equation}
\label{eq:InfiniteDimObstr}
\fV_{\ft}
:=
\tsC_{\ft}\times_{\sG_{\ft}}
L^2_{k-1}(\La^+\otimes \fg_{\ft}) \oplus L^2_{k-1}(V^-)
\to
\sC_{\ft}.
\end{equation}
The $S^1$ action on $\sC_{\ft}$ lifts to an $S^1$ action on
the space \eqref{eq:InfiniteDimObstr} given by
the diagonal action of the $S^1$ action on $\tsC_{\ft}$  defined in
\eqref{eq:DefineS1LAction} and scalar multiplication
on the element of $L^2_{k-1}(V^-)$.

We extend the map \eqref{eq:InfiniteDimObstr}
to a map from a union of
spaces
to
$\bar\sC_{\ft}$ by setting
\begin{equation}
\label{eq:ExtendedInfiniteDimObstrBundle}
\bar\fV_{\ft}
:=
\bigsqcup_{\ell=0}^N \
\left(\fV_{\ft(\ell)}\times\Sym^\ell(X)\right)
\to
\bar\sC_{\ft}=
\bigsqcup_{\ell=0}^N \ \left(\sC_{\ft(\ell)}\times\Sym^\ell(X)\right),
\end{equation}
where $N$ is the integer appearing in Theorem \ref{thm:Compactness}
and the space $\bar\sC_{\ft}$
is defined in \eqref{eq:IdealPairs}.
The
space \eqref{eq:ExtendedInfiniteDimObstrBundle}
is endowed with a topology given by the following notion of convergence.

\begin{defn}
\label{defn:ObstrBundleConvergence}
A sequence $\{[A_\alpha,\Phi_\alpha,\bx_\alpha,\Psi_\alpha]\}_{\alpha=1}^\infty\subset\bar\fV_{\ft}$,
where $\{[A_\alpha,\Phi_\alpha,\bx_\alpha]\}_{\alpha=1}^\infty \subset\bar\sC_\ft$
and $\{\Psi_\alpha\}_{\alpha=1}^\infty$ is a subset of the fiber of \eqref{eq:InfiniteDimObstr},
converges to $[A_\infty,\Phi_\infty,\bx_\infty,\Psi_\infty]$ as
$\alpha\to\infty$ if
\begin{enumerate}
\item
$\{[A_\alpha,\Phi_\alpha,\bx_\alpha]\}_{\alpha=1}^\infty$ converges
to $[A_\infty,\Phi_\infty,\bx_\infty]$ in the Uhlenbeck topology on $\bar\sC_\ft$.
\item
If $u_\alpha:V_{\ft}|_{X\setminus\bx} \to V_{\ft(\ell)}|_{X\setminus \bx}$ are the
\spinu bundle isomorphisms appearing in Definition \ref{defn:UhlenbeckConvergence}
such that $\{u_\alpha(A_\alpha,\Phi_\alpha)\}_{\alpha=1}^\infty$ converges to $(A_0,\Phi_0)$ in $L^2_{k,\loc}$
over $X\setminus\bx$ as $\alpha\to\infty$, then $\{u_\alpha(\Psi_\alpha)\}_{\alpha=1}^\infty$ converges to $\Psi_\infty$
in the $L^2_{k-1,\loc}$ topology on $X\setminus \bx_\infty$  as $\alpha\to\infty$.
\end{enumerate}
\end{defn}
The $\SO(3)$-monopole equations then
extend to give a continuous section,
$\fS:\bar\sC_{\ft}\to\bar\fV_{\ft}$.

We define an $S^1$-invariant inner product
on the fibers of $\bar\fV_{\ft}$ by
\begin{equation}
\label{eq:L2FiberNorm}
\langle [A,\Phi,\bx,\Psi_1],[A,\Phi,\bx,\Psi_2] \rangle
:=
(\Psi_1,\Psi_2)_{L^2(X)},
\end{equation}
where $[A,\Phi,\bx]\in \bar\sC_{\ft}$ and
$[A,\Phi,\bx,\Psi_i]\in \bar\fV_{\ft}$.

\section{Background obstruction bundle}
\label{sec:SWObstruction}
Recall from \cite{FL3} or \cite[Section 3.6.1]{FLLevelOne}
that the background component of the obstruction
is constructed by splicing in sections of the obstruction bundle,
$$
\pi_N^*\Xi_{\ft(\ell),\fs}\to N_{\ft(\ell),\fs}(\delta),
$$
where $\pi_N:N_{\ft(\ell),\fs}(\delta)\to M_{\fs}$ is the projection and
$\Xi_{\ft(\ell),\fs}\cong M_{\fs}\times\CC^{r_\Xi}$ is a
product
bundle.
The embedding $N_{\ft(\ell),\fs}(\delta)\to\sC_{\ft(\ell)}$ is covered
by an embedding of vector bundles,
$$
\begin{CD}
\pi_N^*\Xi_{\ft(\ell),\fs} @>>> \fV_{\ft(\ell)}
\\
@VVV @VVV
\\
N_{\ft(\ell),\fs}(\delta) @>>> \sC_{\ft(\ell)}
\end{CD}
$$
The preceding vector-bundle map is $S^1$-equivariant with
respect to the actions described in \cite[Section 3.6.1]{FLLevelOne}
on the domain and the action on $\fV_{\ft(\ell)}$ described
following \eqref{eq:InfiniteDimObstr}.
If $\pi_N^*\tXi_{\ft(\ell),\fs}\to \tilde N_{\ft(\ell),\fs}(\delta)$ is the pullback of $\pi_N^*\Xi_{\ft(\ell),\fs}$
by the projection $\tilde N_{\ft(\ell),\fs}(\delta)\to \tilde N_{\ft(\ell),\fs}(\delta)/\sG_\fs=N_{\ft(\ell),\fs}(\delta)$
(see \eqref{eq:VirtualNormalBundleTotalSpace}), so $\tXi_{\ft(\ell),\fs}/\sG_{\fs}=\Xi_{\ft(\ell),\fs}$
then the preceding vector-bundle map is defined by a $\sG_\fs$ to $\sG_\ft$ equivariant map
$$
\pi_N^*\tXi_{\ft(\ell),\fs} \to \tsC_{\ft(\ell)}
\times
L^2_{k-1}(\La^+\otimes\fg_{\ft(\ell)}\oplus V^-_{\ft(\ell)}),
$$
with respect to the homomorphism $\varrho:\sG_\fs\to\sG_\ft$ defined in \eqref{eq:GaugeGroupInclusion}.
Define
\begin{equation}
\label{eq:DefineBackgroundObstrOverGluing}
\Xi_\fs(\ft,\sP):=
\pi_N^*\tilde\Xi_{\ft(\ell),\fs}
\times_{\sG_{\fs}}
\sO(\ft,\fs,\sP)
\subset
\pi_N^*\tilde\Xi_{\ft(\ell),\fs}
\times_{\sG_{\fs}}
\bar\Gl(\ft,\fs,\sP).
\end{equation}
Because the
$S^1$
action on $\pi_N^*\tilde\Xi_{\ft(\ell),\fs}$
described above commutes with the $\sG_{\fs}$ action,
this $S^1$ action defines an $S^1$ action on $\Xi_\fs(\ft,\sP)$
and
covers the $S^1$ action described in
Lemma \ref{lem:S^1ActionOnSplicingData}.
To embed $\Xi_\fs(\ft,\sP)$ in $\bar\fV_\ft$,
we now apply the crude splicing construction
from Section \ref{sec:CrudeSplicingMap}.
We define the \emph{convex complement} of a function
$\beta$ to be $1-\beta$.
By multiplying the obstruction sections
by the convex complements of the
cut-off functions used to define the flattened section
$\Phi'$ in the definition of the flattening map $\bTheta_{\sP}$,
we define an $S^1$-equivariant vector-bundle embedding,
\begin{equation}
\label{eq:BackObstrEmbedding}
\varphi_{\fs}''(\sP):\Xi_{\fs}(\ft,\sP) \to \bar\fV_{\ft},
\end{equation}
covering the crude splicing map,
\begin{equation}
\label{eq:CrudeObstructionBundleEmbeddingDiagram1}
\begin{CD}
\pi_N^*\tilde\Xi_{\ft(\ell),\fs}
\times_{\sG_{\fs}}
\sO(\ft,\fs,\sP)
@> \varphi_{\fs}''(\sP) >>
\bar\fV_{\ft}
\\
@VVV @VVV
\\
\tilde N_{\ft(\ell),\fs}(\delta)\times_{\sG_{\fs}}\sO(\ft,\fs,\sP)
@>\bga_{\ft,\fs,\sP}'' >> \bar\sC_{\ft}
\end{CD}
\end{equation}
The crude obstruction splicing maps $\varphi_{\fs}''(\sP')$ are
invariant under the action of the symmetric group and thus,
for $\sP<\sP'$, define a crude obstruction splicing map,
$$
\varphi_{\fs}''([\sP<\sP'])
:
\pi_N^*\tXi_{\fs}\times_{\sG_{\fs}}\sO(\ft,\fs,[\sP<\sP'])
\to
\bar\fV_{\ft},
$$
covering the crude splicing map $\bga''_{\ft,\fs,[\sP<\sP']}$
defined in \eqref{eq:UpperStratumCrudeSplicing}.

We can
describe
the overlaps of the images of the embeddings $\varphi_{\fs}''(\sP)$
exactly as was done in the construction of the space of global splicing
data.  The overlap maps, $\rho^{\ft,\fs,u}_{\sP,[\sP']}$
and $\rho^{\ft,\fs,d}_{\sP,[\sP']}$, appearing in the diagram
\eqref{eq:GlobalSplicingCD}
are equal to the identity on the factor  $\tilde N_{\ft(\ell),\fs}(\delta)$ in their domain.
Hence, these overlap maps
are covered by
$S^1$-equivariant bundle maps,
\begin{equation}
\label{eq:SWObstructionOverlap}
\begin{aligned}
{}&
\rho^{\Xi,\ft,\fs,u}_{\sP,[\sP']}:
\pi_N^*\tXi_{\ft(\ell),\fs}\times_{\sG_{\fs}}\sO(\ft,\fs,\sP,[\sP'])
\to
\pi_N^*\tXi_{\ft(\ell),\fs}\times_{\sG_{\fs}}\sO(\ft,\fs,[\sP<\sP']),
\\
{}&
\rho^{\Xi,\ft,\fs,d}_{\sP,[\sP']}:
\pi_N^*\tXi_{\fs}\times_{\sG_{\fs}}\sO(\ft,\fs,\sP,[\sP'])
\to
\pi_N^*\tXi_{\fs}\times_{\sG_{\fs}}\sO(\ft,\fs,\sP),
\end{aligned}
\end{equation}
equal to the identity on $\pi_N^*\tXi_{\ft(\ell),\fs}$.
Exactly as in the proof of Proposition
\ref{prop:XOverlapControl} (specifically, the equality of
the sections $\Phi_u''$ and $\Phi_d''$ appearing in \eqref{eq:OverlapCompositionUpAndDownImages}),
we see that the diagram
\begin{equation}
\label{eq:SWObstructionOverlapDiagram}
\begin{CD}
\pi_N^*\tXi_{\ft(\ell),\fs}\times_{\sG_{\fs}}\sO(\ft,\fs,\sP,[\sP'])
@> \rho^{\Xi,\ft,\fs,u}_{\sP,[\sP']} >>
\pi_N^*\tXi_{\ft(\ell),\fs}\times_{\sG_{\fs}}\sO(\ft,\fs,[\sP<\sP'])
\\
@V \rho^{\Xi,\ft,\fs,d}_{\sP,[\sP']} VV
@V \varphi_{\fs}''([\sP<\sP']) VV
\\
\pi_N^*\tXi_{\fs}\times_{\sG_{\fs}}\sO(\ft,\fs,\sP)
@>\varphi_{\fs}''(\sP)>> \bar\fV_{\ft}
\end{CD}
\end{equation}
commutes and covers the diagram \eqref{eq:GlobalSplicingCD}.

Because the maps $\rho^{\Xi,\ft,\fs,u}_{\sP,[\sP']}$ and $\rho^{\Xi,\ft,\fs,d}_{\sP,[\sP']}$
in the diagram \eqref{eq:SWObstructionOverlapDiagram} are isomorphisms on the fibers of the
bundles and cover open embeddings, we can define a vector bundle
\begin{equation}
\label{eq:DefineGlobalBackgroundObstruction}
\bar\Upsilon^s_{\ft,\fs} \to \bar\sM^{\vir}_{\ft,\fs}
\end{equation}
as the union of the bundles $\Xi_s(\ft,\sP)$ as $\sP$ varies over
the partitions $\sP$ of $N_\ell$, patched together on
the overlaps by the diagram \eqref{eq:SWObstructionOverlapDiagram}.
The embeddings
$\varphi_s''(\sP)$ fit together to define an
$S^1$-equivariant embedding of vector bundles,
\begin{equation}
\label{eq:DefineCrudeGlobalBackgroundObstrEmbedd}
\varphi''_s:\bar\Upsilon^s_{\ft,\fs} \to \bar\fV_{\ft},
\end{equation}
covering the embedding
given by the crude splicing maps, as shown in \eqref{eq:CrudeObstructionBundleEmbeddingDiagram1}.

Like the crude splicing map,
the embedding $\varphi''_s$
in \eqref{eq:DefineCrudeGlobalBackgroundObstrEmbedd}
is
not equal to the identity
over $M_\fs\times\Sym^\ell(X)$.  We therefore define a new
embedding of vector bundles
in a manner similar to that used in defining
the isotopies of Section \ref{sec:GlobalSplicingMap}.
That is, let $\varphi'_s(\sP)$ be the standard slicing map,
$$
\varphi'_s(\sP):\Xi_s(\ft,\fs,\sP) \to \bar\fV_{\ft},
$$
defined by using the cut-off functions vanishing on balls
$B(x_P,4\la_P^{1/3})$.  There is an
isotopy
through vector-bundle embeddings
between
$\varphi'_s(\sP)$ and $\varphi''_s$ given by
taking convex linear
combinations of the two sections parameterized by `time'.
One then follows the inductive argument, using these isotopies through
vector-bundle embeddings, defined like the isotopies in
\eqref{eq:SplicingIsotopyInduction1}
and \eqref{eq:SplicingIsotopyInduction2}, to
define a sequence of embeddings of vector bundles,
$\varphi''_{s,j}:\bar\Upsilon^s_{\ft,\fs}\to\bar\fV_{\ft}$.
In a manner similar to the definition of the global splicing map in \eqref{eq:DefineGlobalSplice}, define
\begin{equation}
\label{eq:DefineGlobalBackgroundObstrEmbedd}
\varphi'_s:\bar\Upsilon^s_{\ft,\fs}\to\bar\fV_{\ft}
\end{equation}
to be the
last embedding of vector bundles in the sequence $\varphi''_{s,j}$ of embeddings of vector bundles.
The embedding of vector bundles $\varphi_s'$
then covers
the global splicing map $\bga'_{\sM}$ in \eqref{eq:DefineGlobalSplice}.

\section{Equivariant Dirac index bundle}
\label{sec:EquivDiracBundle}
To define the instanton component of the obstruction pseudo-bundle \eqref{eq:SplicedEndsObstr1},
we begin by defining the bundle $\Ind(\bD^*_\ka)$.
Let $\fs_{S^4}=(\rho,S^+,S^-)$ be the standard spin structure
on $S^4$ with respect to the round metric.
Let $E_{\ka}\to S^4$ be a complex Hermitian rank-two vector bundle
with $c_2(E_{\ka})=\ka$  and let
$V_\ka^\pm=S^\pm\otimes E_\ka$ be the resulting \spinu structure
on $S^4$.

For $[A]\in M^{s,\natural}_{\ka}(S^4,\delta)$,
we define the Dirac operator
\begin{equation}
\label{eq:DefineS4Dirac}
D_A: L^2_k( V_\ka^+) \to L^2_{k-1}(V_\ka^-)
\end{equation}
to be the Dirac operator defined by tensor product of
the spin connection induced by the Levi-Civita connection
on $\fs_{S^4}$
and of an $\SU(2)$ connection
$A$ on $E_{\ka}$.
Let
$$
\ind(\bD^*_{\ka}) \to M^{s,\natural}_{\ka}(S^4,\delta)
$$
be the vector bundle
defined by the cokernels of the Dirac operators (using $\Coker D_A=\Ker D_A^*$),
\begin{equation}
\label{eq:DefineInstantonDiracBundle}
\ind(\bD^*_{\ka})
:=
\{ (A,F^s,\Psi)\in \tM^{s,\natural}_{\ka}(S^4)\times \Fr(E_{\ka})|_s
\times L^2_{k-1}(S^-\otimes E_{\ka}): D_A^*\Psi=0
\}/\sG_{\ka}.
\end{equation}
Because of the positive scalar curvature
of the standard round metric
on $S^4$,
the estimates in \cite[Section 4.1]{FL1} using a Bochner--Weitzenb\"ock formula for the Dirac operator
imply
that $\Ker D_A=0$ for all $[A,F^s]\in M^{s,\natural}_\ka(S^4,\delta)$.  The fact that the index of the Dirac
operator is constant with respect to $[A,F^s]$ then ensures
that $\ind(\bD^*_\ka)$ defines a vector bundle.
The Chern class of this bundle is
discussed in
\cite{AtiyahJones,ASDiracVector}.
A computation using the Atiyah--Singer Index Theorem shows that $\ind(\bD^*_\ka)$
is a complex, rank-$\ka$ vector subbundle of
the restriction of the Hilbert bundle,
\begin{equation}
\label{eq:S4InfiniteDimObstrBundle}
\fV_\ka
:=
\sA^\natural_\ka(2\delta)\times\Fr(E_\ka)|_s\times_{\sG_\ka} L^2_{k-1}(V^-_\ka)
\to
\sB^s_\ka(S^4,2\delta),
\end{equation}
to $M^{s,\natural}_{\spl,\ka}(\delta)$.
Using the same topology as that described for
the infinite-dimensional obstruction bundle
$\bar\fV_\ft$ in Definition \ref{defn:ObstrBundleConvergence},
we define an extension
of the Hilbert bundle \eqref{eq:S4InfiniteDimObstrBundle}
over the space \eqref{eq:ExtendedFramedQuotientSpace} of ideal framed connections
by
\begin{equation}
\label{eq:ExtendedS4InfiniteDimObstrBundle}
\bar\fV_\ka:=\bigsqcup_{\ell=0}^\ka \fV_{\ka-\ell}\times\Sym^\ell(X).
\end{equation}
If we define a subbundle of the restriction of $\bar\fV_\ka$ to $\barM^{s,\natural}_{\ka}(S^4,\delta)$ by
\begin{equation}
\label{eq:ExtensionOfDiracIndexBundle}
\overline{\ind}(\bD^*_{\ka}):
=
\{[A,F^s,\bx,\Psi]\in\bar\fV_\ka: [A,F^s,\bx]\in \barM^{s,\natural}_{\ka}(S^4,\delta)
\text{ and } D_A^*\Psi=0\},
\end{equation}
then the restriction of $\overline{\ind}(\bD^*_{\ka})$ to $M^{s,\natural}_{\ka}(S^4,\delta)$
is exactly $\ind(\bD^*_{\ka})$.
Because the index of the Dirac operator depends on $\ka$,
the rank of the fiber of the projection map $\overline{\ind}(\bD^*_{\ka})\to \barM^{s,\natural}_{\ka}(S^4,\delta)$
decreases on the lower strata of $\barM^{s,\natural}_{\ka}(S^4,\delta)$ and
$\overline{\ind}(\bD^*_{\ka})$ is only a pseudo-bundle.  The inclusion maps $\iota_\Si$ for this bundle,
as required by Definition \ref{defn:PseudoBundle},
are defined in \cite[Proposition 7.1.32]{DK}.

\section{The action of $\Spinu(4)$}
\label{sec:SpinuAction}
To splice elements of the instanton obstruction pseudo-bundle
over $S^4$ onto $X$, we must use local trivializations
of the \spinu bundles given by frames.  Changing these frames
is equivalent to changing the element of the obstruction pseudo-bundle being spliced in by an action which we now describe.

As in \cite[Section 3.2]{FLLevelOne}, we denote
$$
\Spinu(4)
:=
\left( \U(2)\times \Spinc(4)\right)/S^1.
$$
The standard
homomorphisms, $\Ad_U:\U(2)\to\SO(3)$ and
$\Ad_S:\Spin^c(4)\to\SO(4)$ (see \cite[Equation (D.2)]{LM}), define
a homomorphism,
\begin{equation}
\label{eq:DefineAdu}
\Ad^u:\Spinu(4) \ni [M,S] \mapsto (\Ad_U(M),\Ad_S(S)) \in \SO(3)\times\SO(4),
\end{equation}
where $M\in\U(2)$, $S\in\Spinc(4)$, and $[M,S]\in \Spinu(4)$.
The composition of $\Ad^u$ with projection onto
the factors $\SO(3)$ and $\SO(4)$ define
a pair of surjective homomorphisms (see \cite[Equation (3.13)]{FLLevelOne}),
$$
\Ad^u_{\SO(3)}:\Spinu(4)\to\SO(3),
\quad\text{and}\quad
\Ad^u_{\SO(4)}:\Spinu(4)\to\SO(4).
$$
The action of $\SO(4)\times\SO(3)$ on
$\bar\sB^s_\ka(S^4)$
is defined in \cite[Section 3.3]{FLLevelOne}, with
$\SO(3)$ acting on the frame
and $\SO(4)$ acting by
pulling the rotation action on $\RR^4$ back to $S^4$
by stereographic
projection.
The $\SO(4)$ action on $S^4$
induces an action on connections
on bundles over $S^4$
by pullback.
The action of $\Spinu(4)$ on $\bar\fV_\ka$ is
defined in \cite[Section 3.6.2]{FLLevelOne}  (see also \cite[Section 3.2]{FeehanGeometry}).
The map
$\bar\fV_\ka\to\bar\sB^s_\ka(S^4)$ is equivariant with respect to the
action of $\Spinu(4)$ on $\bar\fV_\ka$ and the action
of $\Spinu(4)$ on $\bar\sB^s_\ka(S^4)$
defined by $\Ad^u$ and the action of $\SO(4)\times\SO(3)$ on $\bar\sB^s_\ka(S^4)$.
The bundle $\Ind(\bD^*_\ka)$ is invariant under this action and
thus is also a $\Spinu(4)$-equivariant bundle.

\section{Pseudo-bundle over the instanton moduli space with spliced ends}
\label{sec:ObstructionSplicing}
The  obstruction pseudo-bundle over the instanton moduli space with spliced ends
\eqref{eq:SplicedEndsObstr1}
will be defined to be
 identical to
the index pseudo-bundle $\bar{\Ind}(\bD^*_\ka)$
defined in \eqref{eq:ExtensionOfDiracIndexBundle}
on the complement
of the spliced ends of $\barM^{s,\natural}_{\spl,\ka}(\delta)$.
On the spliced end, $W_\ka$, the  obstruction pseudo-bundle
will be
constructed as
the image of a splicing map which we now introduce.

Recall from \eqref{eq:EqualsSplicingImage} in Theorem \ref{thm:ExistenceOfSplicedEndsModuli}
that the
spliced
end of
$\barM^{s,\natural}_{\spl,\ka}(\delta)$
near the
stratum containing the product connection,
$$
[\Theta]\times\Si(Z_\ka(\delta),\sP),
$$
where $\sP$ is a partition of $N_\ka$ with $|\sP|>1$,
is given by the image of a splicing map,
$$
\bga'_{\Theta,\sP}:
\sO^{\asd}(\Theta,\sP)\subset
\Delta^\circ(Z_\ka(\delta),\sP)\times_{\fS(\sP)}
\prod_{P\in\sP} \barM^{s,\natural}_{\spl,|P|}(\delta_P)
\to
\barM^{s,\natural}_{\spl,\ka}(\delta).
$$
We will define the
obstruction pseudo-bundle over the instanton moduli space with spliced ends inductively,
just as we defined the instanton moduli space with spliced ends.

An element of the space $\bar\fV_\ka$ in
\eqref{eq:ExtendedS4InfiniteDimObstrBundle}
can be written as
$[A,F^s,\bx,\Psi]$ where $[A,F^s,\bx]\in\bar\sB^{s}_{\ka}(S^4,2\delta)$,
and $\bx\in\Sym^\ell(\RR^4)$,
and $\Psi\in L^2_{k-1}(S^-\otimes E_{\ka-\ell})$.
For each partition
$\sP$ of $N_\ka$ we define an embedding, covering
the splicing map $\bga_{\Theta,\sP}'$, by
\begin{equation}
\label{eq:DefineS4SectionSplicingToTrivial}
\begin{aligned}
{}&
\varphi_{\Theta,\sP}':
\Delta^\circ(Z_\ka(\delta),\sP)
\times_{\fS(\sP)}
\prod_{P\in\sP} \bar\fV_{|P|}
\to
\bar\fV_\ka,
\\
{}&\quad
[x_P,[A_P,F^s_P,v_P,\Psi_P]]_{P\in\sP}
\mapsto
\left(\bga'_{\Theta,\sP}([x_P,[A_P,F^s_P,v_P]_{P\in\sP}),
\sum_{P\in\sP} \chi_{x_P,\frac{1}{8}\la_P^{1/3}} c_{x_P,1}^*\Psi_P
\right).
\end{aligned}
\end{equation}
If $x_P\in\RR^4$, then $\chi_{x_P,\frac{1}{8}\lambda_P^{1/3}}:\RR^4\to [0,1]$ is a smooth cut-off
function that is equal to one on $B(x_P,\frac{1}{8}\lambda_P^{1/3})$ and is supported in $B(x_P,\frac{1}{4}\lambda_P^{1/3})$
(see Section \ref{subsubsec:StdSplicingMap}), where $\la_P=\la([A_P,F^s_P,\bx_P])$.
The map $c_{x_P,1}$ is defined in Lemma \ref{lem:ComposeCenterings}.
The definition \eqref{eq:DefineS4SectionSplicingToTrivial} follows
that in
\cite[Section 8.2 and Equation (8.32)]{FL3}.

Observe that the map \eqref{eq:DefineS4SectionSplicingToTrivial}
is equivariant with respect to the $\Spinu(4)$ action if
$\Spinu(4)$ acts diagonally on the domain, by the projection $\Ad^u_{\SO(4)}$
to $\SO(4)$ on the factor $\Delta^\circ(Z_\ka(\delta),\sP)$.
The properties of the obstruction pseudo-bundle over the instanton moduli space with spliced ends
are given by the

\begin{thm}
\label{thm:ExistenceOfSplicedEndsIndex}
There exists a pseudo-bundle,
\begin{equation}
\label{eq:SplicedEndsIndexBundle}
\bar\Upsilon^i_{\spl,\ka}\to \barM^{s,\natural}_{\spl,\ka}(\delta),
\end{equation}
which is a $\Spinu(4)$-invariant subbundle of the
space
$\bar\fV_\ka\to\bar\sB^{s,\natural}_\ka(S^4,2\delta)$
and
has the following properties:
\begin{enumerate}
\item
\label{item:ExistenceOfSplicedEndsIndex1}
The restriction of $\bar\Upsilon^i_{\spl,\ka}$ to
the top stratum, $M^{s,\natural}_{\spl,\ka}(\delta)$,
is a complex rank-$\ka$ vector bundle.
\item
\label{item:ExistenceOfSplicedEndsIndex2}
For any partition $\sP$ of $\ka$ with $|\sP|>1$,
the restriction of $\bar\Upsilon^i_{\spl,\ka}$ to the
open neighborhood of $[\Theta]\times \Si(Z_\ka(\delta),\sP)$
given in \eqref{eq:EqualsSplicingImage} is equal to the image of
\begin{equation}
\label{eq:SplicedEndIndexBundle}
\Delta^\circ(Z_\ka(\delta),\sP)\times_{\fS(\sP)}
\prod_{P\in\sP}\bar\Upsilon^i_{\spl,|P|},
\end{equation}
under the splicing map $\varphi'_{\Theta,\sP}$.
\item
\label{item:ExistenceOfSplicedEndsIndex3}
If $W_\ka\subset \barM^{s,\natural}_{\spl,\ka}(\delta)$ is as defined in
\eqref{eq:DefineSplicedEnd}, then
the restriction of $\bar\Upsilon^i_{\spl,\ka}$
to $\barM^{s,\natural}_{\spl,\ka}(\delta)\setminus W_\ka$ is given by the
restriction of the pseudo-bundle $\overline{\Ind}(\bD^*_\ka)$
in \eqref{eq:ExtensionOfDiracIndexBundle} to
$\barM^{s,\natural}_{\spl,\ka}(\delta)\setminus W_\ka$.
\item
\label{item:ExistenceOfSplicedEndsIndex4}
The restriction of $\bar\Upsilon^i_{\spl,\ka}$
to the levels
$$
\bar M^{s,\natural}_{\spl,\ka}(\delta)\cap
\left(
M^{s,\natural}_{\spl,\ka-\ell}(\delta)\times\Sym^\ell(\RR^4)
\right),
$$
of $\bar M^{s,\natural}_{\spl,\ka}(\delta)$
is isomorphic to the pullback of $\bar\Upsilon^i_{\spl,\ka-\ell}$
under the projection
$$
M^{s,\natural}_{\spl,\ka-\ell}(\delta)\times\Sym^\ell(\RR^4)
\to
M^{s,\natural}_{\spl,\ka-\ell}(\delta).
$$
\item
\label{item:ExistenceOfSplicedEndsIndex5}
For any $[A',F^s,\bx]\in\barM^{s,\natural}_{\spl,\ka}(\delta)$,
$L^2$-orthogonal projection defines an isomorphism,
$$
\Upsilon^i_{\spl,\ka}|_{[A',F^s,\bx]} \cong \Coker D_{A_t},
$$
for all $t\in [0,1]$, where $A_t=A'+a_t(A')$ and $A=A'+a_1(A')$ is the
point in the image of the gluing map
appearing in \eqref{eq:S4GluingMap} at $A'$.
\end{enumerate}
\end{thm}

\begin{rmk}
\label{rmk:PseudoBundleOnLowerStrata}
Observe that the rank of the fibers of
$\bar\Upsilon^i_{\spl,\ka}\to \bar M^{s,\natural}_{\spl,\ka}(\delta)$
varies with the stratum of $M^{s,\natural}_{\spl,\ka}(\delta)$
by
Item \eqref{item:ExistenceOfSplicedEndsIndex4}.
For this reason,  $\bar\Upsilon^i_{\spl,\ka}$ is a pseudo-bundle.
\end{rmk}

The proof of
Theorem \ref{thm:ExistenceOfSplicedEndsIndex} occupies the rest of this section
and is similar to the construction of the instanton moduli space with spliced ends in 
the proof of Theorem \ref{thm:ExistenceOfSplicedEndsModuli}.

\subsection{Pseudo-bundles and overlap data}
\label{subsec:SplicedEndsObstrOverlap}
If the restriction of the
obstruction pseudo-bundle over the instanton moduli space with spliced ends to the image of the splicing map
$\bga'_{\Theta,\sP}$
is to be given by the image
of the splicing map $\varphi'_{\Theta,\sP}$ restricted to
the subspace \eqref{eq:SplicedEndIndexBundle}
as asserted in Item \eqref{item:ExistenceOfSplicedEndsIndex2} of Theorem \ref{thm:ExistenceOfSplicedEndsIndex},
then we need to describe the overlaps of the images of
the splicing maps $\varphi'_{\Theta,\sP}$ and
$\varphi'_{\Theta,\sP'}$ where $\sP<\sP'$.
This is undertaken by constructing a space of overlap
data analogous to that appearing in
Proposition \ref{prop:CommutingSplicingSymmQuotient}.

We begin by initially considering the pseudo-bundle,
$\tilde\Upsilon(\Theta,\sP,\sP')$, defined by restricting
the pseudo-bundle,
\begin{equation}
\label{eq:DefineOverlapInstantonObstr}
\begin{CD}
\Delta^\circ(Z_\ka(\delta),\sP)
\times
\prod_{P\in\sP}
\left(
\Delta^\circ(Z_{|P|}(\delta_P),\sP'_P)
\times\prod_{Q\in\sP'_P} \bar\Upsilon^i_{\spl,|Q|}
\right)
\\
@VVV
\\
\Delta^\circ(Z_\ka(\delta),\sP)
\times
\prod_{P\in\sP}
\left(
\Delta^\circ(Z_{|P|}(\delta_P),\sP'_P)
\times\prod_{Q\in\sP'_P} \bar M^{s,\natural}_{\spl,|Q|}(\delta_Q)
\right)
\end{CD}
\end{equation}
to the intersection of its base with the
space $\tilde\sO(\Theta,\sP,\sP',\delta)$
defined in Lemma \ref{lem:OverlappingTrivial}.
By induction on $\ka$, the pseudo-bundle $\bar\Upsilon^i_{\spl,|Q|}$
is
defined for all $Q$ in a partition $\sP$
of $N_\ka$ of length greater than one.
For all $P\in\sP$, the splicing maps
$\varphi_{\Theta,\sP'_P}$ define maps,
$$
\varphi_{\Theta,\sP'_P}:
\Delta^\circ(Z_{|P|}(\delta_P),\sP'_P)
\times\prod_{Q\in\sP'_P} \bar\Upsilon^i_{\spl,|Q|}
\to
\bar\Upsilon^i_{\spl,|P|},
$$
covering the splicing map, $\bga'_{\Theta,\sP'_P}$.
Hence, there is a map,
$$
\rho^{\Psi,d}_{\sP,\sP'}:
\tilde\Upsilon(\Theta,\sP,\sP')
\to
\tilde\Upsilon(\Theta,\sP)=
\Delta^\circ(Z_{\ka}(\delta),\sP)\times\prod_{P\in\sP}\bar\Upsilon^i_{\spl,|P|},
$$
defined by
\begin{equation}
\label{eq:DefineDownwardsObstrOverlapOnS4}
\rho^{\Psi,d}_{\sP,\sP'}
:=
\id_\Delta\times\prod_{P\in\sP}\varphi'_{\Theta,\sP'_P},
\end{equation}
where $\id_\Delta$ is the identity map on
$\Delta^\circ(Z_{\ka}(\delta),\sP)$.
By the definition of the overlap
map $\rho^{\Theta,d}_{\sP,\sP'}$
in \eqref{eq:DefineUpwardTransition}, the
map $\rho^{\Psi,d}_{\sP,\sP'}$ fits into the
following commutative diagram of maps of pseudo-bundles,
\begin{equation}
\label{eq:SplicedObstrUpwards}
\begin{CD}
\tilde\Upsilon(\Theta,\sP,\sP')
@> \rho^{\Psi,d}_{\sP,\sP'} >>
\tilde\Upsilon(\Theta,\sP)
\\
@VVV @VVV
\\
\sO^{\asd}(\Theta,\sP,\sP',\delta)
@> \rho^{\Theta,d}_{\sP,\sP'} >>
\sO^{\asd}(\Theta,\sP)
\end{CD}
\end{equation}
The map $\rho^{\Psi,d}_{\sP,\sP'}$ is not injective because
of the action of the symmetric group, but this will not
affect our construction.

We also define a map,
\begin{equation}
\label{eq:UpwardsTransitionObstruction}
\rho^{\Psi,u}_{\sP,\sP'}:
\tilde\Upsilon(\Theta,\sP,\sP')
\to
\tilde\Upsilon(\Theta,\sP'),
\end{equation}
which fits into the
following commutative diagram of maps of pseudo-bundles,
\begin{equation}
\label{eq:SplicedObstrDownwards}
\begin{CD}
\tilde\Upsilon(\Theta,\sP,\sP')
@> \rho^{\Psi,u}_{\sP,\sP'} >>
\tilde\Upsilon(\Theta,\sP')
\\
@VVV @VVV
\\
\sO^{\asd}(\Theta,\sP,\sP',\delta)
@> \rho^{\Theta,u}_{\sP,\sP'} >>
\sO^{\asd}(\Theta,\sP')
\end{CD}
\end{equation}
The following lemma, analogous to
Lemma \ref{lem:CommutingSplicing}, is the key
to the construction of the bundle $\Upsilon^i_{\spl,\ka}$.

\begin{lem}
\label{lem:OverlapObstructionCommuting}
Let $\sP<\sP'$ be partitions of $N_\ka$ with
$|\sP|>1$.
Then the following diagram of pseudo-bundles,
\begin{equation}
\label{eq:OverlapObstructionCommuting}
\begin{CD}
\tilde\Upsilon(\Theta,\sP,\sP')
@> \rho^{\Psi,d}_{\sP,\sP'} >>
\tilde\Upsilon(\Theta,\sP)
\\
@V \rho^{\Psi,u}_{\sP,\sP'} VV
@V \varphi'_{\Theta,\sP} VV
\\
\tilde\Upsilon(\Theta,\sP') @> \varphi'_{\Theta,\sP'} >>
\bar\fV_\ka
\end{CD}
\end{equation}
covers the diagram \eqref{eq:CommutingSplicing},
is $\Spinu(4)$-equivariant, and commutes.
\end{lem}

\begin{proof}
The assertion that the diagram \eqref{eq:OverlapObstructionCommuting}
covers the diagram \eqref{eq:CommutingSplicing} follows
from the definition of the maps
$\rho^{\Psi,d}_{\sP,\sP'}$ and $\rho^{\Psi,u}_{\sP,\sP'}$
(see the diagrams \eqref{eq:SplicedObstrUpwards} and
\eqref{eq:SplicedObstrDownwards}) and the fact that the
splicing maps $\varphi'_{\Theta,\sP}$ cover the
splicing maps $\bga'_{\Theta,\sP}$.

The $\Spinu(4)$ equivariance of the diagram
\eqref{eq:OverlapObstructionCommuting}
follows
from that of the splicing maps $\varphi'_{\Theta,\sP}$.

We now prove the commutativity of the diagram \eqref{eq:OverlapObstructionCommuting}.
Assume that
a point in
$\tilde\Upsilon(\Theta,\sP,\sP')$ is given by data consisting of
points $(x_P)_{P\in\sP}\in\Delta^\circ(Z_{\ka}(\delta),\sP)$,
points $(y_Q)_{Q\in\sP'_P}\in\Delta^\circ(Z_{|P|}(\delta_P),\sP'_P)$,
and $[A_Q,F^s_Q,\bv_Q,\Psi_Q]\in\bar\Upsilon_{\spl,|Q|}$.
We write
\begin{align*}
[A''_u,\Phi''_u,\bx_u,\Psi_u]
{}&=
(\varphi'_{\Theta,\sP'}\circ \rho^{\Psi,u}_{\sP,\sP'})
\left(
    \left( (x_P),\left((y_Q),[A_Q,F^s_Q,\bv_Q,\Psi_Q] \right)_{Q\in\sP'_P}\right)_{P\in\sP}
\right),
\\
[A''_d,\Phi''_d,\bx_d,\Psi_d]
{}&=
(\varphi'_{\Theta,\sP}\circ\rho^{\Psi,d}_{\sP,\sP'})
\left(
    \left( (x_P),\left((y_Q),[A_Q,F^s_Q,\bv_Q,\Psi_Q] \right)_{Q\in\sP'_P}\right)_{P\in\sP}
\right).
\end{align*}
Because the diagram \eqref{eq:OverlapObstructionCommuting}
covers the commuting diagram \eqref{eq:CommutingSplicing}, to show that \eqref{eq:OverlapObstructionCommuting}
commutes, we only need to prove that $\Psi_u=\Psi_d$.

From the definition of $\rho^{\Psi,u}_{\sP,\sP'}$ in \eqref{eq:UpwardsTransitionObstruction} and
the definition of $\rho^{\Theta,u}_{\sP,\sP'}$ in \eqref{eq:R4TrivialSplicingUpwardsOverlapMap},
we have,
$$
 \rho^{\Psi,u}_{\sP,\sP'}
\left(
    \left( (x_P),\left((y_Q),[A_Q,F^s_Q,\bv_Q,\Psi_Q] \right)_{Q\in\sP'_P}\right)_{P\in\sP}
\right)
=
\left( x_{P(Q)}+y_Q,[A_Q,F^s_Q,\bv_Q,\Psi_Q]\right)_{Q\in\sP'},
$$
where, for $Q\in\sP'$, we write $P(Q)\in\sP$ for the unique subset of $N_\ell$ with $Q\subseteq P$.
If $\la_P=\la([A_P,F^s_P,\bx'_P])$, the definition of $\varphi'_{\Theta,\sP}$ then implies that
\begin{equation}
\label{eq:UpwardsSplicingSection1}
\Psi_u=\sum_{Q\in\sP'} \chi_{x_{P(Q)}+y_Q,\frac{1}{8}\lambda_Q^{1/3}}c_{x_{P(Q)}+y_Q,1}^*\Psi_Q.
\end{equation}
From the definition of $\rho^{\Psi,d}_{\sP,\sP'}$ in \eqref{eq:DefineDownwardsObstrOverlapOnS4},
if we write
\begin{align*}
{}&\rho^{\Psi,d}_{\sP,\sP'}
\left(
    \left( (x_P),\left((y_Q),[A_Q,F^s_Q,\bv_Q,\Psi_Q \right)_{Q\in\sP'_P}\right)_{P\in\sP}
\right)
\\{}&\quad=
    \left( (x_P),[A_P,F^s_P,\bx'_P,\Psi_P \right)_{P\in\sP}
\end{align*}
then, for  $\la_Q=\la([A_Q,F^s_Q,\bv_Q])$,
$$\Psi_P=\sum_{Q\in\sP'_P}\chi_{y_Q,\frac{1}{8}\la_Q^{1/3}}c_{y_Q,1}^*\Psi_Q.
$$
If $\la_P=\la([A_P,F^s_P,\bx'_P])$, then the definition of $\varphi_{\Theta,\sP}'$ implies that
$$
\Psi_d
=
\sum_{P\in\sP}\chi_{x_P,\frac{1}{8}\la_P^{1/3}}c_{x_P,1}^*\Psi_P
=
\sum_{P\in\sP}\sum_{Q\in\sP'_P}
\chi_{x_P,\frac{1}{8}\la_P^{1/3}}c_{x_P,1}^*\left(\chi_{y_Q,\frac{1}{8}\la_Q^{1/3}}c_{y_Q,1}^*\Psi_Q\right).
$$
Applying
Lemma \ref{lem:ComposeCenterings} to the preceding
equalities
gives
\begin{equation}
\label{eq:IteratedSplicing1}
\Psi_d=
\sum_{Q\in\sP'}
\chi_{x_{P(Q)},\frac{1}{8}\la_{P(Q)}^{1/3}}\chi_{x_{P(Q)}+y_Q,\frac{1}{8}\la_Q^{1/3}}c^*_{x_Q+y_Q,1}\Psi_Q.
\end{equation}
By the constraint
\eqref{eq:OverlapCondition3},
the points $x_P$ and $y_Q$ in the definition of
the subspace $\tilde\sO(\Theta,\sP,\sP')$
satisfy (see condition \eqref{eq:BallContainment2}),
$$
B(x_{P(Q)}+y_Q,\tquarter\la_Q^{1/3})\Subset B(x_{P(Q)},\teighth\la_{P(Q)}^{1/3}).
$$
Because $\chi_{x,\frac{1}{8}\la^{1/3}}$ is supported on
$B(x,\frac{1}{4} \la^{1/3})$ and equal to one on
$B(x,\teighth \la^{1/3})$,
the preceding inclusion implies that
$\chi_{x_{P(Q)},\frac{1}{8}\la_{P(Q)}^{1/3}}$ is equal to one on the support of
$\chi_{x_{P(Q)}+y_Q,\frac{1}{8}\la_Q^{1/3}}$, so
\eqref{eq:IteratedSplicing1}  and \eqref{eq:UpwardsSplicingSection1} imply that
$$
\Psi_d=
\sum_{Q\in\sP'}
\chi_{x_{P(Q)}+y_Q,\frac{1}{8}\la_Q^{1/3}}c^*_{x_Q+y_Q,1}\Psi_Q=\Psi_u,
$$
as required to prove that the diagram \eqref{eq:OverlapObstructionCommuting} commutes.
\end{proof}

\begin{proof}[Proof of Theorem \ref{thm:ExistenceOfSplicedEndsIndex}]
The proof of this
theorem is largely identical to the
proof of Theorem \ref{thm:ExistenceOfSplicedEndsModuli} so
we omit many of the details.
We construct $\bar\Upsilon^i_{\spl,\ka}$ using induction on $\ka$.

The induction begins with setting
the restriction of $\bar\Upsilon^i_{\spl,1}$
to the top stratum of $\barM^{s,\natural}_{\spl,1}(\delta)$
equal to the bundle $\ind(\bD^*_1)$.  Because there is,
up to dilation and the $\SO(3)$ action on the frame, only
one gauge-equivalence class in $\barM^{s,\natural}_{\spl,1}(\delta)$,
the bundle,
$$
\overline{\ind}(\bD^*_1)
\to
\barM^{s,\natural}_{\spl,1}(\delta)=(0,\delta)\times\SO(3),
$$
is given by
$$
(0,\delta)\times\SU(2)\times_{\pm 1}\CC \to (0,\delta)\times\SO(3).
$$
This bundle extends as a pseudo-bundle over the Uhlenbeck
compactification, $c(\SO(3))$,
of $M^{s,\natural}_{\spl,1}(\delta)$
as
$$
c\left( \SU(2)\times_{\pm 1}\CC\right)
\to
c(\SO(3)),
$$
which satisfies the requirements of
Theorem \ref{thm:ExistenceOfSplicedEndsIndex} for the case $\ka=1$.

Assuming that $\bar\Upsilon^i_{\spl,\ka'}$ has been constructed
for all $\ka'<\ka$,
we define a pseudo-bundle $\bar\Upsilon^i(W_\ka)\to W_\ka$
over the spliced end by the image of the splicing maps
$\varphi'_{\Theta,\sP}$ on the domains in
\eqref{eq:SplicedEndIndexBundle}.
This space is a pseudo-bundle because the overlap
of the images of the maps $\varphi'_{\Theta,\sP}$ is
parameterized
by Lemma \ref{lem:OverlapObstructionCommuting}.

The final property in Theorem \ref{thm:ExistenceOfSplicedEndsIndex}
regarding the surjectivity of $L^2$-orthogonal projection follows
from the results of
\cite[Section 8]{FL3} (see also \cite[Proposition 7.1.32]{DK}).
We can then define an isotopy $R_D$ of
$\bar\Upsilon^i(W_\ka)$ by, for
$t \in [1/2, 3/4]$,
$$
R_D\left(t,[A,F^s,\bx,\Psi]\right)
:=
\left(R(t,[A,F^s,\bx]),
(1-4(t-1/2))\Psi+ 4(t-1/2)\Pi_{A+a(4(t-1/2),A)}\Psi
\right)
$$
where $R$ is the isotopy defined in
\eqref{eq:GluingIsotopy} and for $a(t,A)=a_t(A)$ defined by
the gluing map appearing in \eqref{eq:S4GluingMap},
$$
G(t,[A,F^s,\bx]) = [A+a(t,A),F^s,\bx],
$$
and where $\Pi_{A+a(t,A)}$ denotes $L^2$-orthogonal onto $\Coker D_{A+a(t,A)}$.
For
$t\in(-\infty,1/2]$,
the isotopy $R_D(t,\cdot)$
is equal to the identity and for
$t\in[3/4,1]$,
the isotopy $R_D(\cdot,t)$
is defined by pulling back the section $\Psi$ by the
centering map.

The pseudo-bundle $\bar\Upsilon^i_{\spl,\ka}$ is then defined by
the union of the images $\bar\Upsilon^i(W_\ka)$ under the map
$R_D(\beta(\cdot,)\cdot)$, where $\beta$ is the function defined
in the proof of Theorem \ref{thm:ExistenceOfSplicedEndsModuli},
with $\Ind(\bD^*_\ka)$ restricted to the complement
$$
\barM^{s,\natural}_{\ka}(S^4,\delta)
\setminus
R(1,W_\ka),
$$
as in the definition of $\barM^{s,\natural}_{\spl,\ka}(\delta)$
in the proof of Theorem \ref{thm:ExistenceOfSplicedEndsModuli}.
\end{proof}

\section{Instanton obstruction pseudo-bundle}
\label{sec:InstantonObstr}
We now construct the instanton component of the
obstruction pseudo-bundle and its embedding
into the infinite-rank obstruction pseudo-bundle
$\bar\fV_{\ft}$.

\subsection{The frame bundles}
\label{subsec:Instanton_obstruction_pseudo-bundle_frames}
For a \spinu structure $\ft=(\rho,V)$ over $X$, we defined a frame bundle
in \cite[Section 3.2]{FLLevelOne} by
$$
\Fr_{\CCl(T^*X)}(V)\to \Fr(TX)\to X,
$$
which is a $\U(2)$ bundle over $\Fr(TX)$ and a $\Spinu(4)$ bundle over
$X$.  The fiber over
$F\in\Fr(TX)|_x$ is given by the Clifford module
isomorphisms from the Clifford module $\Delta\otimes\CC^2$
to $V|_x$ with respect to the
isomorphism $\RR^4\cong T^*X|_x$ defined by $F$, for each $x \in X$.
The homomorphism $\Ad^u:\Spinu(4)\to\SO(3)\times\SO(4)$
in \eqref{eq:DefineAdu} defines a bundle map,
$$
\Fr_{\CCl(T^*X)}(V) \to \Fr(\fg_{\ft})\times_X\Fr(TX).
$$
By analogy with the definition of the gluing-data bundle
$\Fr(\ft,\fs,\sP)$ in \eqref{eq:DefineGluingDataBundle},
we define
\begin{equation}
\label{eq:FrameBundleForObstruction}
\Fr(V,\sP,g_{\sP}) \to \Fr(\ft,\fs,\sP),
\end{equation}
by
\begin{multline*}
\Fr(V,\sP,g_{\sP})
:=
\left\{(F^V_1,\dots,F^V_\ell)\in \left.\prod_{i=1}^\ell\Fr_{\CCl(T^*)}(V_{\ft(\ell)})\right|_{\Delta^\circ(X^\ell,\sP)}:\right.
\\
\left. F^V_i=F^V_j \iff \exists\, P\in\sP \text{ with } i,j\in P \right\}.
\end{multline*}
The structure group of the bundle $\Fr(V,\sP,g_{\sP})\to \Delta^\circ(X^\ell,\sP)$ is
\begin{equation}
\tilde G^V(\sP)
=
\left\{(G_1,\dots,G_\ell)\in\prod_{i=1}^\ell\Spinu(4):
G_i=G_j \iff \exists\, P\in\sP \text{ with } i,j\in P \right\}.
\end{equation}
The structure group of the bundle
$\Fr(V,\sP,g_{\sP})\to \Si(X^\ell,\sP)$ is
\begin{equation}
G^V(\sP):=\tilde G^V(\sP)\rtimes\fS(\sP)/\Ga(\sP).
\end{equation}
We then define the instanton obstruction pseudo-bundle for the partition
$\sP$ by
\begin{equation}
\label{eq:DefineInstObstrBundle}
\tilde N_{\ft(\ell),\fs}(\delta)\times_{\sG_{\fs}}\times\Ob(\ft,\fs,\sP)
\to
\tilde N_{\ft(\ell),\fs}(\delta)\times_{\sG_{\fs}}\times\sO(\ft,\fs,\sP),
\end{equation}
where
\begin{equation}
\label{eq:InstantonObstructionBundle}
\Ob(\ft,\fs,\sP)
:=
\Fr(V,\sP,g_{\sP})
\times_{G^V(\sP)}
\prod_{P\in\sP}\bar\Upsilon^{s,\natural}_{\spl,|P|}.
\end{equation}
We define an $S^1$ action on the bundle \eqref{eq:DefineInstObstrBundle}
covering the $S^1$ action on its base defined in
Lemma \ref{lem:S^1ActionOnSplicingData} by the diagonal action on
the factors
$\tilde N_{\ft(\ell),\fs}(\delta)$ and $\bar\Upsilon^{s,\natural}_{\spl,|P|}$
in \eqref{eq:DefineInstObstrBundle}, where $S^1$ acts on
$\tilde N_{\ft(\ell),\fs}(\delta)$ by the action \eqref{eq:S1ZActionOnN}
and on $\bar\Upsilon^{s,\natural}_{\spl,|P|}$ by scalar multiplication
on the fibers of $\bar\Upsilon^{s,\natural}_{\spl,|P|}\to\bar M^{s,\natural}_{\spl,|P|}$.

\subsection{Splicing map}
\label{subsec:Instanton_obstruction_pseudo-bundle_splicing_map}
For each partition $\sP$ of $N_\ell$, we define an
$S^1$-equivariant splicing
map for the instanton obstruction,
\begin{equation}
\label{eq:InstantonObstrSplicingMap}
\varphi_{\ft,\fs,\sP}:
\tN_{\ft(\ell),\fs}(\delta)\times_{\sG_{\fs}}\Ob(\ft,\fs,\sP)
\to
\bar\fV_\ft,
\end{equation}
which will cover the crude splicing map $\bga''_{\ft,\fs,\sP}$.
A choice of frames $(F^V_P)_{P\in\sP}\in\Fr(V,\sP)|_{\by}$,
where $\by=(y_P)_{P\in\sP}$, a \spinu connection $A_0$
on $V_{\ft(\ell)}$,
an $\SO(3)$ connection $A_P$ on $\fg_P$,
and a frame $F^s_P\in\Fr(\fg_P)|_s$
determine  embeddings,
\begin{equation}
\label{eq:CliffordBundleIsomorphism}
\phi(A_0,F^V_P,A_P,F^s_P):V_{|P|}|_{\varphi_n(B(8\la^{1/3}))}
\cong
V_{\ft}|_{B(y_P,8\la^{1/3})},
\end{equation}
as described in \cite[Section 3.4.3]{FLLevelOne}.
If $\chi$ is a cut-off function supported on $\varphi_n(B(8\la^{1/3}))$
and $\Psi$ is a section of $V_{|P|}$, then
we write $\phi(A_0,F^V_P,A_P,F^s_P)(\chi\Psi)$ for the section of
$V_\ft|_{B(y_P,8\la^{1/3})}$ given by the isomorphism
\eqref{eq:CliffordBundleIsomorphism}.

To ensure that the images of the different splicing maps
fit together to form a vector pseudo-bundle, we introduce a flattening map
on \spinu connections just
as
we did for Riemannian metrics on $X$ and connections
on $\fg_{\ft(\ell)}$
in Section \ref{sec:FlatteningPairs}.  That is, given a \spinu
connection $A_0$ on $V_{\ft(\ell)}$ we apply the construction
of Lemma \ref{lem:CompatibleFlattening} and for each
$\by\in \Si(X^\ell,\sP)$, we define a locally flattened \spinu
connection $\bTheta^u(A,\by)$ on $V_{\ft(\ell)}$
satisfying the following properties:
\begin{enumerate}
\item
For $\by=(y_P)_{P\in\sP}$, the connection
$\bTheta^u(A_0,\by)$ is flat on
$$
\bigcup_{P\in\sP} B(y_p,4\tilde s_P(\by)^{1/3}),
$$
where $\tilde s_P$ is the separating function defined in Lemma \ref{lem:CompatibleFlattening}.
\item
If $\sP<\sP'$ and $\by\in \Si(X^\ell,\sP)$ and
$\by'\in\Si(X^\ell,\sP')\cap \sU(X^\ell,\sP)$
satisfy $\pi(X^\ell,\sP)(\by')=\by$, then
$$
\bTheta^u(A_0,\by)=\bTheta^u(A_0,\by').
$$
\end{enumerate}
For $F^V=(F^V_P)_{P\in\sP}\in\Fr(V,\sP)|_{\by}$, define
$$
\phi'(A_0,F^V_P,A_P,F^s_P):=
\phi(\bTheta^u(A_0,\by),F^V_P,A_P,F^s_P),
$$
using the isomorphism \eqref{eq:CliffordBundleIsomorphism}
defined by the locally flattened connection $\bTheta^u(A_0,\by)$ on $V_{\ft(\ell)}$ instead of the connection $A_0$.

For $[A_0,\Phi_0]\in N_{\ft(\ell),\fs}(\delta)\subset\sC_{\ft(\ell)}$
and a point in $\Ob(\ft,\fs,\sP)$
given by the data
$(F^V_P)_{P\in\sP}\in \Fr(V,\sP,g_{\sP})$, and
$[A_P,F^s_P,\bv_P]\in\barM^{s,\natural}_{\spl,|P|}(\delta_P)$, and
$[A_P,F^s_P,\bv_P,\Psi]\in\bar\Upsilon_{\spl,|P|}$,
and
$$
\bga''_{\ft,\fs,\sP}((A_0,\Phi_0),F^V_P,[A_P,F^s_P,\bv_P])_{P\in\sP}
=(A'',\Phi'',\bx),
$$
the instanton obstruction splicing map is then defined by
\begin{equation}
\label{eq:DefineCrudeInstantonObstrSplicing}
\begin{aligned}
\varphi_{\ft,\fs,\sP}
{}&
\left(
\left[(A_0,\Phi_0),(F^V_P,[A_P,F^s_p,\bv_P,\Psi_P])_{P\in\sP}\right]
\right)
\\
{}&:=
\left[
A'',\Phi'',\bx,\sum_{P\in\sP}\phi'(A_0,F^V_P,A_P,F^s_P)(\chi_{n,\frac{1}{8}\la_P^{1/3}}\Psi_P)
\right],
\end{aligned}
\end{equation}
where $\chi_{n,\frac{1}{8}\la_P^{1/3}}:S^4\to [0,1]$ is a smooth cut-off function that is equal to one
on $\varphi_n(B(\frac{1}{8}\lambda_P^{1/3}))$ and is supported in $\varphi_n(B(\frac{1}{4}\lambda_P^{1/3}))$
(see Section \ref{subsubsec:StdSplicingMap}).
This map is $S^1$-equivariant with respect to the
$S^1$ action on its domain described following
\eqref{eq:InstantonObstructionBundle} and the
$S^1$ action on the image $\bar\fV_\ft$ defined
following \eqref{eq:ExtendedInfiniteDimObstrBundle}.

\subsection{Overlap space and overlap maps}
\label{subsec:Instanton_obstruction_pseudo-bundle_overlap_space_maps}
Following the method of previous sections, to
parameterize
the overlap of the images
of the splicing maps $\varphi_{\ft,\fs,\sP}$ and
$\varphi_{\ft,\fs,\sP'}$, we introduce a space of
overlap data and overlap maps.
By analogy with the definition of the overlap space
$\bar\Gl(\ft,\fs,\sP,[\sP'])$ in \eqref{eq:OverlapGluingDataSpace},
we define
$$
\Ob(\ft,\fs,\sP,[\sP'])
\to
\bar\Gl(\ft,\fs,\sP,[\sP'])
$$
by
\begin{equation}
\label{eq:XInstantonObstrOverlapSpace}
\begin{aligned}
{}&
\Ob(\ft,\fs,\sP,[\sP'])
\\
&\quad:=
\Fr(V,\sP,g_{\sP})
\times_{G^V(\sP)}
\bigsqcup_{\sP''\in[\sP<\sP']}
\prod_{P\in\sP}
\left(
\Delta^\circ(Z_{|P|}(\delta_P),\sP''_P)
\times\prod_{Q\in\sP''_P} \bar\Upsilon_{\spl,|Q|}
\right).
\end{aligned}
\end{equation}
Define an $S^1$ action on $\Ob(\ft,\fs,\sP,[\sP'])$
through the diagonal $S^1$ action on the factors of
$\bar\Upsilon_{\spl,|Q|}$ appearing in
\eqref{eq:XInstantonObstrOverlapSpace}.

We define overlap maps,
$$
\rho^{V,d}_{\sP,\sP'}:
\tilde N_{\ft(\ell),\fs}(\delta)\times_{\sG_{\fs}}
\Ob(\ft,\fs,\sP,[\sP'])
\to
\tilde N_{\ft(\ell),\fs}(\delta)\times_{\sG_{\fs}}
\Ob(\ft,\fs,\sP),
$$
covering the overlap maps $\rho^{\ft,\fs,d}_{\sP,[\sP']}$
defined in \eqref{eq:XDownwardTransition} by
\begin{equation}
\label{eq:XDownwardsInstantonOverlapMap}
\rho^{V,d}_{\sP,\sP'}
:=
\id_{\tilde N_{\ft(\ell),\fs}(\delta)}\times
\id_{\Fr(V)}\times
\bigsqcup_{\sP''\in [\sP<\sP']}
\prod_{P\in\sP}\varphi_{\Theta,\sP''_P},
\end{equation}
where $\id_{\Fr(V)}$ is the identity map on
$\Fr(V,\sP,g_{\sP})$.  Observe that $\rho^{V,d}_{\sP,\sP'}$
is $S^1$-equivariant if $S^1$ acts on the domain by
the diagonal action on $\tilde N_{\ft(\ell),\fs}(\delta)$ and
$\Ob(\ft,\fs,\sP,[\sP'])$, where $S^1$ acts on
$\tilde N_{\ft(\ell),\fs}(\delta)$ by the action \eqref{eq:S1ZActionOnN}
and on $\Ob(\ft,\fs,\sP,[\sP'])$ by the action
defined following \eqref{eq:XInstantonObstrOverlapSpace}
and $S^1$ acts on the image
by the action defined
following \eqref{eq:InstantonObstructionBundle}.

For
\begin{equation}
\label{eq:DefineImageOfXUpwardsInstantonOverlapMap}
\begin{aligned}
{}&\tilde N_{\ft(\ell),\fs}(\delta)\times_{\sG_{\fs}}\Ob(\ft,\fs,[\sP<\sP'])
\\
{}&\quad:=
\left.\tilde N_{\ft(\ell),\fs}(\delta)\times_{\sG_{\fs}}
\left(
\bigsqcup_{\sP''\in[\sP<\sP']}
\Ob(\ft,\fs,\sP'')
\right)\right/\fS(\sP),
\end{aligned}
\end{equation}
we define an upwards overlap map,
\begin{equation}
\label{eq:XUpwardsInstantonOverlapMap}
\begin{CD}
\tilde N_{\ft(\ell),\fs}(\delta)\times_{\sG_{\fs}}
\Ob(\ft,\fs,\sP,[\sP'])
\\
@V \rho^{V,u}_{\sP,[\sP']} VV
\\
\left.\tilde N_{\ft(\ell),\fs}(\delta)\times_{\sG_{\fs}}
\left(
\bigsqcup_{\sP''\in[\sP<\sP']}
\Ob(\ft,\fs,\sP'')
\right)\right/\fS(\sP)
\end{CD}
\end{equation}
exactly as was done in \eqref{eq:UnparamUpwardsTransition}.
Recall that the map $\rho^{\ft,\fs,u}_{\sP,[\sP']}$ was defined by
a parallel translation of the frames in
the domain $\sO(\ft,\fs,\sP,[\sP'],g_{\sP})$
and leaving the points in $\tilde N_{\ft(\ell),\fs}(\delta)$ and
the $S^4$ connections unchanged.
Note that this parallel translation is carried out with
respect to the locally flattened connection $\bTheta^u(A_0,\by)$.
We define the map $\rho^{V,u}_{\sP,[\sP']}$
in the same way, using the flattened
\spinu connection $A_0$ to parallel translate the frames
of $\Fr_{\CCl(T^*X)}(V_{\ft(\ell)})$ and leaving the elements
of $\tilde N_{\ft(\ell),\fs}(\delta)$ and
of $\bar\Upsilon_{\spl,|Q|}$ unchanged.
The map $\rho^{V,u}_{\sP,[\sP']}$ is $S^1$-equivariant
when the $S^1$ actions on the range and domain are
as described in the analogous assertion for $\rho^{V,d}_{\sP,[\sP']}$.

By the invariance of the obstruction splicing maps under
the action of the symmetric group, the splicing maps
$\varphi_{\ft,\fs,\sP''}$ for $\sP''\in [\sP<\sP']$
define a splicing map,
\begin{equation}
\label{eq:XCrudeInstantonSplicingMultiple}
\varphi_{\ft,\fs,[\sP<\sP']}:
\tilde N_{\ft(\ell),\fs}(\delta)\times_{\sG_{\fs}}
\Ob(\ft,\fs,[\sP<\sP'])
\to
\bar\fV_{\ft},
\end{equation}
which covers the crude splicing map
$\bga''_{\ft,\fs,[\sP<\sP']}$ in \eqref{eq:UpperStratumCrudeSplicing}.
We then have the following relation between the splicing maps
$\varphi_{\ft,\fs,\sP}$ and $\varphi_{\ft,\fs,[\sP<\sP']}$.

\begin{lem}
\label{lem:InstantonObstrCommDiagr}
Let $\sP<\sP'$ be partitions of $N_\ell$.  Then the diagram
\begin{equation}
\label{eq:InstantonObstrCommDiagr}
\begin{CD}
\tilde N_{\ft(\ell),\fs}(\delta)\times_{\sG_{\fs}}
\Ob(\ft,\fs,\sP,[\sP'])
@>\rho^{V,u}_{\sP,[\sP']} >>
\tilde N_{\ft(\ell),\fs}(\delta)\times_{\sG_{\fs}}
\Ob(\ft,\fs,[\sP<\sP'])
\\
@V\rho^{V,d}_{\sP,[\sP']} VV
@V \varphi_{\ft,\fs,[\sP<\sP']} VV
\\
\tilde N_{\ft(\ell),\fs}(\delta)\times_{\sG_{\fs}}
\Ob(\ft,\fs,\sP)
@> \varphi_{\ft,\fs,\sP} >> \bar\fV_{\ft}
\end{CD}
\end{equation}
commutes,
covers the diagram \eqref{eq:GlobalSplicingCD},
and all the maps in the diagram are $S^1$-equivariant.
\end{lem}

\begin{proof}
That the diagram \eqref{eq:InstantonObstrCommDiagr} covers
the diagram \eqref{eq:GlobalSplicingCD} follows immediately
from the definitions
of the overlap (see \eqref{eq:XDownwardsInstantonOverlapMap} and
the paragraph following \eqref{eq:XUpwardsInstantonOverlapMap})
and obstruction splicing
maps (see \eqref{eq:DefineCrudeInstantonObstrSplicing}).
As in the proof of Proposition
\ref{prop:XOverlapControl}, we assume that the splicing maps
$\varphi_{\ft,\fs,\sP}$ and $\varphi_{\ft,\fs,\sP'}$ are
defined at points $\by\in\Si(X^\ell,\sP)$ and
$\by'\in\Si(X^\ell,\sP')$, respectively, where
$\pi(X^\ell,\sP)(\by')=\by$.  Thus, for
$P'\in\sP'$ and $P'\subseteqq P\in\sP$, the inclusion
$$
B(y_{P'}, 8\tilde s(\by')^{1/3})
\Subset
B(y_P,4 \tilde s(\by)^{1/2})
$$
holds by \eqref{eq:BallInclusions}.
Because the Riemannian metric and \spinu connection used in the splicing
argument are flat on the preceding balls,
the conclusion of the lemma then follows from the argument establishing
Lemma \ref{lem:OverlapObstructionCommuting}.
\end{proof}

The commutativity of the diagram \eqref{eq:InstantonObstrCommDiagr}
and the definition of the space $\bar\sM^{\vir}_{\ft,\fs}$
as the union of the spaces
$$
\tN_{\ft(\ell),\fs}(\delta)\times_{\sG_{\fs}}\sO(\ft,\fs,\sP)
$$
implies that we can define a pseudo-bundle,
\begin{equation}
\bar\Upsilon^i_{\ft,\fs} \to \bar\sM^{\vir}_{\ft,\fs},
\end{equation}
as the union of the pseudo-bundles
$$
\tN_{\ft(\ell),\fs}(\delta)\times_{\sG_{\fs}}\Ob(\ft,\fs,\sP)
\to
\tN_{\ft(\ell),\fs}(\delta)\times_{\sG_{\fs}}\sO(\ft,\fs,\sP),
$$
with the overlaps identified by the diagram
\eqref{eq:InstantonObstrCommDiagr}.  Because the maps in
\eqref{eq:InstantonObstrCommDiagr} are $S^1$-equivariant
and the $S^1$ actions in that diagram cover the $S^1$ actions
on $\bar\sM^{\vir}_{\ft,\fs}$, the $S^1$ actions on
$\tN_{\ft(\ell),\fs}(\delta)\times_{\sG_{\fs}}\Ob(\ft,\fs,\sP)$
define a global $S^1$ action on $\bar\Upsilon^i_{\ft,\fs}$.

The splicing maps $\varphi_{\ft,\fs,\sP}$ define
a vector-bundle embedding of $\bar\Upsilon^i_{\ft,\fs}$
into $\bar\fV_{\ft}$ which covers the embedding of
$\bar\sM^{\vir}_{\ft,\fs}$ into $\bar\sC_{\ft}$ given
by the crude splicing maps.
Finally, we define an embedding,
\begin{equation}
\label{eq:XGlobalInstantonObstrSplicing}
\varphi'_i:\bar\Upsilon^i_{\ft,\fs}\to\bar\fV_{\ft},
\end{equation}
covering the global splicing map $\bga'_{\sM}$ in
\eqref{eq:DefineGlobalSplice} by replacing the
triple $[A'',\Phi'',\bx]$ in the definition
\eqref{eq:DefineCrudeInstantonObstrSplicing} with the triple
$[A',\Phi',\bx]$ given by the global splicing map $\bga'_{\sM}$
of \eqref{eq:DefineGlobalSplice}.

\section{Local gluing hypothesis for $\SO(3)$ monopoles}
\label{sec:GluingThm}
We can now state the relevant gluing hypothesis.

\begin{hyp}[Local gluing hypothesis]
\label{hyp:Gluing}
There is a continuous, $S^1$-equivariant embedding,
\begin{equation}
\label{eq:GluingMap}
\bga_{\sM}:\bar\sM^{\vir}_{\ft,\fs} \to \bar\sC_{\ft},
\end{equation}
which is homotopic
through $S^1$-equivariant, continuous embeddings to the global splicing map $\bga'_{\sM}$,
smooth on each stratum of $\bar\sM^{\vir}_{\ft,\fs}$,
and
equal to the identity on $N_{\ft(\ell),\fs}(\delta)\times\Sym^\ell(X)$.
In addition, there are  $S^1$-equivariant sections $\bchi_s$ and $\bchi_i$
of the pseudo-bundles $\bar\Upsilon^s_{\ft,\fs}$ and
$\bar\Upsilon^i_{\ft,\fs}$ with the following properties:
\begin{enumerate}
\item
\label{item:GluingHyp1}
The restriction of the section $\bar\bchi=\bchi_s\oplus \bchi_i$ of
$\bar\Upsilon^s_{\ft,\fs}\oplus \bar\Upsilon^i_{\ft,\fs}$
to each stratum is smooth.

\item
\label{item:GluingHyp2}
The restriction of the section $\bar\bchi$ to
each stratum of $\bar\sM^{\vir}_{\ft,\fs}$ vanishes
transversely.

\item
\label{item:GluingHyp3}
If we pull back the
fiber metric
\eqref{eq:L2FiberNorm}
to $\bar\Upsilon^s_{\ft,\fs}\oplus \bar\Upsilon^i_{\ft,\fs}$
by the splicing embeddings $\varphi_s'\oplus\varphi_i'$
defined in
\eqref{eq:DefineGlobalBackgroundObstrEmbedd} and \eqref{eq:XGlobalInstantonObstrSplicing},
then the $L^2$ norm of $\bar\bchi$ is lower semi-continuous.

\item
\label{item:GluingHyp4}
The restriction of $\bga_{\sM}$ to the zero-locus
$\bar\bchi^{-1}(0)$ is a homeomorphism between
$\bar\bchi^{-1}(0)$ and
an open neighborhood of $M_{\fs}\times\Sym^\ell(X)$
in $\bar\sM_{\ft}$.
\end{enumerate}
\end{hyp}

\begin{rmk}
Recall from \cite[p. 12]{LiebLoss} that a function $f$ is
lower semi-continuous if $f^{-1}((a,\infty))$ is open for all
real $a$. A lower semi-continuous function $f(x)$ satisfies
$\lim_{x\to x_0} f(x) \ge f(x_0)$ and
we will use this property in
Lemma \ref{lem:ExtendingToThomSection} to obtain lower bounds
on the
fiber norm
of the obstruction section on a relatively
closed subspace of the top stratum of $\bar\sM^{\vir}_{\ft,\fs}$.
\end{rmk}

Although there are some important differences (which we explain further in Section \ref{sec:Notes_on_justification_gluing_hypothesis}), the gluing map and obstruction section
in Hypothesis \ref{hyp:Gluing} are similar
to those
constructed in \cite{FL3}, as we now explain.

The gluing map
in \cite[Section 9]{FL3} is the composition of the standard splicing map
(defined in \cite[Section 3.2]{FL3} and here in Section \ref{subsubsec:StdSplicingMap})
with domain given by the space $\bar\sU(\ft,\fs,\sP)$ defined in Section \ref{eq:DefineUSet}
and a
solution map.
Recall that the solution map,
\begin{equation}
\label{eq:GluingPerturb1}
(A',\Phi')\mapsto (A'+a(A',\Phi'),\Phi+\phi(A',\Phi')),
\end{equation}
is defined on the image of the splicing map and has the property that
$\fS\left(A'+a(A',\Phi'),\Phi+\phi(A',\Phi')\right)$ lies in
a finite-dimensional obstruction space over $(A',\Phi')$, and thus is a solution to the extended $\SO(3)$-monopole equations, while $(A,\Phi)$ is an $\SO(3)$ monopole if and only $\fS(A,\Phi) = 0$, where $\fS$ is as in \eqref{eq:PerturbedSO3MonopoleEquations}.
In \cite[Section 8]{FL3}, this obstruction space over $(A',\Phi')$ is the span of
the small-eigenvalue eigenvectors of the $L^2$ self-adjoint elliptic operator $(D\fS_{(A',\Phi'})(D\fS_{(A',\Phi'})^*$,
where $D\fS_{(A',\Phi'}$ is the linearization of the map $\fS$ at $(A',\Phi')$.
The map
$$
(A',\Phi') \mapsto
\left((A',\Phi'),\fS\left(A'+a(A',\Phi'),\Phi+\phi(A',\Phi')\right)\right)
$$
then defines a section of a finite-rank vector bundle over the image of the slicing map
whose zero locus is mapped to $\bar\sM_\ft$ by the
solution map.

To show that a
solution map
in the sense of \eqref{eq:GluingPerturb1}
exists, we first prove in
\cite[Proposition 5.5]{FL3} that for each pair $(A',\Phi')$ in the image of the splicing map, the error term
$\fS(A',\Phi')$ is suitably small
in the sense of norms defined in \cite{FeehanSlice,TauFrame}.
Using techniques generalizing those introduced by Taubes in
\cite{TauSelfDual, TauIndef, TauFrame} and Donaldson and Kronheimer in \cite[Section 7.2]{DK}
(see also \cite{MorganMrowkaTube, MrowkaThesis}),
we showed that for each point $(A',\Phi')$ in the image of the
splicing map, there is a solution to the system of partial differential equations defining
the
solution map
(see \cite[Section 9.1]{FL3}).

Our proof of existence of the
solution map
in \cite{FL3} should extend to yield the desired gluing map, $\bga_\sM$,
in Hypothesis \ref{hyp:Gluing} and identify the zero locus of the obstruction section, $\bar\bchi$, with an open neighborhood in $\bar\sM_\ft$. The smoothness of $\bga_\sM$ and the obstruction section $\bchi$ asserted
in  Property \eqref{item:GluingHyp1} in Hypothesis \ref{hyp:Gluing} follow by
the construction of the gluing map.
Property \eqref{item:GluingHyp2}, that the
obstruction section vanishes transversely, follows from a formal
argument and the result in \cite{FeehanGenericMetric} that
the $\SO(3)$-monopole map, $\fS$, vanishes transversely
for appropriate choices of generic perturbations of the
$\SO(3)$-monopole equations.

We give proofs of the properties of the gluing map
$\bga_\sM$ and obstruction section $\bchi$ asserted by Hypothesis \ref{hyp:Gluing} in \cite{Feehan_Leness_monopolegluingbook}.
We expect that the continuity of the gluing map with respect to
Uhlenbeck limits will follow along lines similar to the proof of the
same property for the
gluing map for anti-self-dual $\SO(3)$ connections
given in \cite{FLKM1}.  This Uhlenbeck continuity will yield
the Property \eqref{item:GluingHyp3} of the gluing map in Hypothesis \ref{hyp:Gluing}.
One must also show the map is injective and surjective where
by the latter we mean
Property \eqref{item:GluingHyp4} given in
Hypothesis \ref{hyp:Gluing}.
In special cases, proofs of
continuity with respect to Uhlenbeck limits, injectivity, and surjectivity
for gluing maps for anti-self-dual $\SO(3)$ connections
have been given in
\cite[Sections 7.2.5 and 7.2.6]{DK} and \cite{TauSelfDual, TauIndef, TauFrame}.

\section{Notes on the justification of the local gluing hypothesis}
\label{sec:Notes_on_justification_gluing_hypothesis}
The purpose of this section is to summarize the justification of
Hypothesis \ref{hyp:Gluing} as a theorem, whose proof is provided by the authors in \cite{Feehan_Leness_monopolegluingbook}. In our article \cite{FL3}, we discussed analytical issues specific to the problem of gluing $\SO(3)$ monopoles that are not present when gluing anti-self-dual $\SO(3)$ connections and we refer the reader to that discussion. Rather, our goal in this section is to describe modifications to our approach in \cite{FL3} to gluing $\SO(3)$ monopoles. In \cite{FL3}, we had restricted our attention to the case of gluing anti-self-dual connections over $S^4$ and $\SO(3)$ monopoles over $X$ that varied in an open neighborhood with the property that spectral flow for the Laplace operator, $d_{A,\Phi}^1d_{A,\Phi}^{1,*}$, is small. This assumption is valid, for example, when considering small open neighborhoods of a \emph{zero}-dimensional moduli space of Seiberg--Witten monopoles (equivalently, strata of reducible $\SO(3)$ monopoles), since those moduli spaces comprise finite sets of isolated points (because those moduli spaces can be assumed to be compact and smooth). However, in this monograph and our previous articles such as \cite{FL2a} (where no bubbling is allowed) or \cite{FLLevelOne} (where one bubble allowed), we must allow for moduli spaces of Seiberg--Witten moduli spaces that are \emph{positive}-dimensional. In those cases, we cannot assume that spectral flow for $d_{A,\Phi}^1d_{A,\Phi}^{1,*}$ is small and thus our method of gluing \cite{FL3} does not apply without significant modification. Indeed, the rank of the obstruction vector bundle constructed in \cite{FL3}, by an analogue of the small-eigenvalue decomposition that Taubes developed in \cite{TauIndef}, necessarily jumps as eigenvalues of $d_{A,\Phi}^1d_{A,\Phi}^{1,*}$ cross the small-eigenvalue cut-off parameter $\mu \in (0,1]$ used to define the $L^2$-orthogonal projections required to split the $\SO(3)$-monopole equations, following the Kuranishi method \cite{Kuranishi}.

In Section \ref{subsec:Virtual_neighborhood_moduli_space_SO3_monopoles_top_level}, we recall the main ideas from \cite[Section 3]{FL2a} 
(see also Section \ref{subsubsec:ThickenedNeighborhood} in this monograph)
involved in our construction of a virtual neighborhood for the moduli space of $\SO(3)$ monopoles near a possibly singular\footnote{In the sense of Kodaira--Spencer deformation theory, rather than in the sense of bubbling or elliptic regularity theory.} stratum of $\SO(3)$ monopoles\footnote{Reducible $\SO(3)$ monopoles or Seiberg--Witten monopoles in our application.} in the top level. This construction relies on a splitting of the $\SO(3)$-monopole equations defined by a choice of abstract stabilizing bundle, $\Xi$, over an open neighborhood, $\sU$, in the configuration space of pairs. When the stratum is in the top level, there is no bubbling and thus no need to glue in anti-self-dual connections over copies of $S^4$. In Section \ref{subsec:Virtual_neighborhood_moduli_space_anti-self-dual_connections}, we describe two methods of constructing a virtual neighborhood for the moduli space of anti-self-dual connections near a possibly singular stratum of anti-self-dual connections\footnote{Reducible anti-self-dual connections in the context of the proof of the Kotschick--Morgan conjecture \cite{FLKM1, KotschickMorgan}} in the top lower-level; we restrict our attention to the case of anti-self-dual connections so that we may focus on the essential ideas in a familiar and relatively simple case and relate our constructions to those of Donaldson and Kronheimer \cite{DK}, an accessible reference. When the stratum is in a lower level, one has bubbling and thus one must glue in anti-self-dual connections over copies of $S^4$ in order to construct a virtual neighborhood; therefore in Section \ref{subsec:Extrinsic_virtual_neighborhood_moduli_space_anti-self-dual_connections_gluing}, we outline how to extend the construction of Section \ref{subsubsec:Donaldson_Kronheimer_extrinsic_splitting} to the case of gluing two families of anti-self-dual connections over a pair of four-manifolds, $X_1$ and $X_2$, to produce new anti-self-dual connections over a connected sum, $X=X_1\# X_2$. Finally, in Section \ref{subsec:Virtual_neighborhood_moduli_space_SO3_monopoles_lower_level}, we very briefly outline the changes needed to extend our discussion of gluing anti-self-dual $\SO(3)$ connections to gluing $\SO(3)$ monopoles and thus construct a virtual neighborhood for the moduli space of $\SO(3)$ monopoles near a possibly singular stratum of $\SO(3)$ monopoles in a lower-level and complete the proof of Hypothesis \ref{hyp:Gluing} as a theorem.

\subsection{Construction of a virtual neighborhood for the moduli space of $\SO(3)$ monopoles near a top-level singular stratum of $\SO(3)$ monopoles}
\label{subsec:Virtual_neighborhood_moduli_space_SO3_monopoles_top_level}
We begin by recalling the basic outline, in the absence of bubbling, from our article \cite[Section 3]{FL2a} for how to construct a virtual neighborhood for the moduli space of $\SO(3)$ monopoles near a top-level singular stratum of $\SO(3)$ monopoles. The `gluing map' and `obstruction section' are defined in terms of a Kuranishi model that describes how the Seiberg--Witten moduli space, $M_\fs$, is contained in the top level, $\sM_\ft$, of the Uhlenbeck compactification, $\bar\sM_\ft$, of the moduli space of $\SO(3)$ monopoles. 

We can assume that $M_\fs$ contains no zero-section pairs by \cite[Corollary 3.3]{FL2a}. The Banach Lie group, $\sG_\ft$, acts freely on the open subspace, $\tilde{\sC}_\ft^0 \subset \tilde{\sC}_\ft$, of non-zero-section pairs and the quotient, $\sC_\ft^0 = \tilde{\sC}_\ft^0/\sG_\ft$, is a smooth Banach manifold by \cite[Proposition 2.8]{FL1}. We have a smooth Hilbert vector bundle, $\fV = \tilde{\sC}_\ft^0 \times_{\sG_\ft} L^2(X;\Lambda^+\otimes\fg_\ft \oplus W^-\otimes E)$ 
over
$\sC_\ft^0$, as in \cite[Equation (3.39)]{FL2a}, where $\ft=(\rho,W\otimes E)$. We obtain 
a 
smooth splitting of Hilbert vector bundles, $\fV \restriction \iota(\sC_\fs^0) = \fV^t\oplus \fV^n$ for $\iota(\sC_\fs) \subset \tilde{\sC}_\ft$, following \cite[Equation (3.40)]{FL2a}, where $\fs = (\rho,W)$. We let $\sU \subset \sC_\ft^0$ denote an open subset that contains a compact subset of $\iota(\sC_\fs^0)$. By  \cite[Theorem 3.19]{FL2a}, there is a finite-rank, smooth, product vector subbundle, $\Xi \subset \fV$, over $\sU$ with the following properties:
\begin{itemize}
  \item $\Xi$ is $S^1$ equivariant with respect to the circle action specified in \cite[Equation (3.2)]{FL2a} (implied by a splitting $E=\underline{\CC}\oplus L$ and $W\otimes E = W\oplus W\otimes L$, with $\CC^*$ acting trivially on $\underline{\CC}=X\times\CC$ and by scalar multiplication on $L$)
      and given in this monograph by \eqref{eq:DefineS1LAction};
  \item $\Xi \restriction \iota(\sC_\fs^0)$ is a smooth complex vector bundle;
  \item The $\SO(3)$-monopole equations  \eqref{eq:PerturbedSO3MonopoleEquations} (compare \cite[Equation (2.32)]{FL2a}) define a smooth section, $\fS$, of $\fV \to \sC_\ft^0$;
  \item If $\Pi_\Xi:\fV \to \Xi$ is the smooth map of vector bundles defined by $L^2$-orthogonal projection and $\fV \restriction \sU = \Xi\oplus \Xi^\perp$ and $\Pi_{\Xi^\perp} := \id_\fV - \Pi_\Xi$, then
      \[
      \Pi_{\Xi^\perp}: \Ran (D\fS)_{A,\Phi} \to \Xi^\perp_{A,\Phi}
      \]
      is a surjective map of fibers, for each $[A,\Phi] \in \sU$, where
      \[
      \fV_{A,\phi} \subset L^2(X;\Lambda^+\otimes\fg_\ft \oplus W^-\otimes E).
      \]
\end{itemize}
The $\SO(3)$-monopole section, $\fS$, of $\fV \to \sC_\ft^0$ induces smooth sections,
\begin{gather*}
\Pi_\Xi\circ\fS \quad\text{of}\quad \Xi \to \sU,
\\
\Pi_{\Xi^\perp}\circ\fS \quad\text{of}\quad \Xi^\perp \to \sU,
\end{gather*}
with $\Pi_{\Xi^\perp}\circ\fS$ defining the virtual moduli space (denoted by $\sM_\ft(\Xi,\fs)$ in \cite[Equation (3.41)]{FL2a}),
\[
\sM_{\ft,\fs}^\vir = \left(\Pi_{\Xi^\perp}\circ\fS\right)^{-1}(0) \subset \sC_\ft^0.
\]
Furthermore, our \cite[Theorem 3.21]{FL2a} provides that:
\begin{itemize}
\item The Seiberg--Witten moduli space, $M_\fs$, is an $S^1$-invariant, smooth
  submanifold of $\sM_{\ft,\fs}^\vir$ via the smooth embedding $\iota:\sC_\fs^0 \to \sC_\ft^0$;
  \item The restriction of $\fS$ to $\sM_{\ft,\fs}^\vir$ takes values in $\Xi$ and vanishes transversely on $\sM_{\ft,\fs}^\vir - \iota(M_\fs)$;
  \item The smooth vector bundle, $N_{\ft,\fs} \to M_\fs$, given by \cite[Equation (3.42)]{FL2a} (where it is denoted by $N_\ft(\Xi,\fs)$)
      and in \eqref{eq:VirtualNormalandObstBundles} in this monograph, and constructed from the deformation complex for $\fS$ is the normal bundle for the submanifold, $\iota:M_\fs \hookrightarrow \sM^{\vir}_{\ft,\fs}$; the tubular map is equivariant with respect to the circle action on $N_{\ft,\fs}$ given by the trivial action on the base $M_\fs$ and complex multiplication on the fibers, and the circle action on $\sM_{\ft,\fs}^\vir$ induced from the $S^1$ action in \cite[Equation (3.2)]{FL2a}.
\end{itemize}
In the next subsection, we shall further explore certain aspects of the construction of virtual neighborhoods in the simpler setting of the moduli space of anti-self-dual connections.

\subsection{Virtual neighborhoods for the moduli space of anti-self-dual connections}
\label{subsec:Virtual_neighborhood_moduli_space_anti-self-dual_connections}
There are essentially two ways of splitting the non-linear anti-self-dual or Yang--Mills equations in the presence of cokernel obstructions in their linearizations:
\begin{enumerate}
\item
\label{item:Taubes_intrinsic}
Taubes' intrinsic splitting for the Yang--Mills equations \cite[Equation (8.4)]{TauFrame} (compare \cite[Equations (2.7) and (2.8)]{TauIndef} for the anti-self-dual equation), and
\item
\label{item:Donaldson_Kronheimer_extrinsic}
Donaldson--Kronheimer's extrinsic splitting for the anti-self-dual equation \cite[Equations (7.26) or (7.27)]{DK}.
\end{enumerate}
These two approaches are inspired by the method employed by Kuranishi \cite{Kuranishi} in the deformation of complex structures and adapted by Atiyah, Hitchin, and Singer \cite{AHS} to the anti-self-dual equation, but the technical differences between them are significant. In this subsection, we describe these methods in their original context for the anti-self-dual equation but they formally extend, with minor changes, to other non-linear equations arising in gauge theory, such as the Yang--Mills, $\SO(3)$-monopole, and pseudo-holomorphic curve equations. Indeed, in Section \ref{subsec:Virtual_neighborhood_moduli_space_SO3_monopoles_lower_level}, we outline such an extension for the $\SO(3)$-monopole equations, as required by our application in this monograph.

An earlier variant of Method \eqref{item:Taubes_intrinsic} was used by
Taubes in \cite{TauIndef} for his construction of solutions to the anti-self-dual equation and generalizations of his idea are
described, for example, by the authors in \cite{FLKM1, FL3}. Method
\eqref{item:Taubes_intrinsic}, proposed by Taubes \cite{TauFrame} in
the context of the Yang--Mills equations, was adapted to
the anti-self-dual equation by Friedman and Morgan in
\cite[Section 3.4.6]{FrM}.  Method
\eqref{item:Donaldson_Kronheimer_extrinsic} is more abstract and,
unlike Method \eqref{item:Taubes_intrinsic}, does not rely on an eigenvalue
splitting for a Laplace operator depending on a connection, $A$, and
thus is more suitable when spectral flow occurs as the connection $A$
varies in a large neighborhood. In this subsection, we shall review these methods in their original setting of the anti-self-dual equation.

\subsubsection{Taubes' intrinsic splitting}
\label{subsubsec:Taubes_intrinsic_splitting}
Let $G$ be a compact Lie group and $P$ be a smooth principal $G$-bundle over a closed, oriented, smooth Riemannian four-manifold $(X,g)$ and $\tilde\sU \subset \sA(P)$ be an open neighborhood in the pre-configuration space of all $W^{1,p}$ connections\footnote{By analogy with the notation of Adams and Fournier \cite{AdamsFournier}, we define the $W^{k,p}$ Sobolev norm of $a \in C^\infty(X;\Lambda^l\otimes\ad P)$ by requiring that $\nabla{A_0}^ja \in L^p(X;\Lambda^l\otimes\ad P)$ for $j=0,1,\ldots,k$, given a smooth reference connection on $P$, integers $k,l \geq 0$, and constant $p\in[1,\infty]$.} $A$ on $P$ (with $p\in(2,\infty)$). For Method \eqref{item:Taubes_intrinsic}, given a constant $\mu \in (0,1]$, one aims to solve the following pair of equations,
\begin{subequations}
\label{eq:ASD_equation}  
\begin{align}
\label{eq:Extended_ASD_equation}
\Pi_\mu^\perp(A) F^+(A) &= 0 \quad\text{(for $A$ in $\sU$)},
\\
\label{eq:Balancing_ASD_equation}
\Pi_\mu(A)F^+(A) &= 0 \quad\text{(for $A$ in $\sU$ and solving \eqref{eq:Extended_ASD_equation})},
\end{align}
\end{subequations}
where $\Pi_\mu(A) \in \sL(W^{2,p}(X;\Lambda^+\otimes\ad P))$ is the finite-rank $L^2$-orthogonal projection from $W^{2,p}(X;\Lambda^+\otimes\ad P)$ onto the subspace spanned by the eigenvectors of the unbounded operator,
\[
d_A^+d_A^{+,*}: L^2(X;\Lambda^+\otimes\ad P) \to L^2(X;\Lambda^+\otimes\ad P),
\]
with eigenvalues less then $\mu/2$, while $\Pi_\mu^\perp(A) := 1 - \Pi_\mu(A)$; we choose $\mu$ so that no eigenvalue of $d_A^+d_A^{+,*}$ lies in $[\mu/2,\mu]$. In applications of Method \eqref{item:Taubes_intrinsic} and its variant in \cite{TauIndef}, such as in Feehan and Leness \cite{FLKM1}, the point $[A]$ varies in an open neighborhood $\sU = \tilde\sU/\Aut(P)$ in the configuration space, $\sB(P) := \sA(P)/\Aut(P)$, of all $W^{1,p}$ connections modulo the action of the group of $W^{2,p}$ gauge transformations, $\Aut(P)$.  We assume that each $[A] \in \sU$ obeys $F^+(A) < \eps$, for a small $\eps \in (0,1]$. We also assume that $\sU$ is small enough that the eigenvalues of $d_A^+d_A^{+,*}$ do not cross $\mu$ --- in other words, small enough that there is no spectral flow as $[A]$ varies in $\sU$. The infinite-dimensional equation \eqref{eq:Extended_ASD_equation} is often called the \emph{extended anti-self-dual equation} and the finite-dimensional equation \eqref{eq:Balancing_ASD_equation} is often called the \emph{balancing equation}.

\subsubsection{Donaldson--Kronheimer's extrinsic splitting}
\label{subsubsec:Donaldson_Kronheimer_extrinsic_splitting}
This method can be viewed as a more flexible version of Method \eqref{item:Taubes_intrinsic} and is especially appropriate when the open $\Aut(P)$-invariant neighborhood, $\sU \subset \sB(P)$, cannot be assumed to be small and one must contend with spectral flow. If $\Xi$ denotes both a finite-rank, product vector bundle over $\sU$ and (slightly abusing notation) its pullback to $\tilde\sU$ (by the quotient map, $\pi:\sA(P) \to \sB(P)$), we let $\Pi_\Xi(A): L^2(X;\Lambda^+\otimes\ad P) \to \Xi(A)$ and $\Pi_\Xi^\perp(A) := 1 - \Pi_\Xi(A) = \Pi_{\Xi^\perp}(A): L^2(X;\Lambda^+\otimes\ad P) \to \Xi^\perp(A)$ denote the $L^2$-orthogonal projection operators for $A \in \tilde\sU$, where $L^2(X;\Lambda^+\otimes\ad P) = \Xi(A) \oplus \Xi^\perp(A)$. As we discuss in the following paragraph, one can construct the bundle $\Xi$ so that the operators,
\begin{equation}
\label{eq:ASD_stabilizing_bundle_surjectivity}
\Pi_\Xi^\perp(A):\Ran d_A^+ \to \Xi^\perp(A),
\end{equation}
are surjective for all $[A] \in \sU$.

For a small enough neighborhood $\sU$, spectral flow is limited and the bundle $\Xi$ can be constructed by decomposing the spectrum of $d_A^+d_A^{+,*}$ into a set of small non-negative eigenvalues bounded above by $\mu/2$ and a set of eigenvalues bounded below by a uniform positive constant $\mu$, just as in Method \eqref{item:Taubes_intrinsic}. One then defines $\Xi(A) := \Xi_\mu(A) \subset L^2(X;\Lambda^+\otimes\ad P)$, the subspace spanned by the eigenvectors of $d_A^+d_A^{+,*}$ with eigenvalues less than $\mu/2$. For a large open neighborhood, $\sU$, one can instead employ the method described by the authors in \cite[Section 3]{FL2a} to construct $\Xi$ in the context of the $\SO(3)$-monopole equations \eqref{eq:PerturbedSO3MonopoleEquations}.
  
Regardless of how the bundle $\Xi$ is constructed, the smooth map,
\[
\tilde\sU \ni A \mapsto   \Pi_\Xi^\perp F^+(A) := \Pi_\Xi^\perp(A) F^+(A) \in \Pi_\Xi^\perp(A)L^2(X;\Lambda^+\otimes\ad P),
\]
defines a smooth section of the smooth vector bundle, $ \Xi^\perp \to \sU$. Moreover, if $F^+(A) = 0$, then
\[
D(\Pi_\Xi^\perp F^+)(A)a = \Pi_\Xi^\perp(A) d_A^+a, \quad \forall\, a \in W^{1,2}(X;\Lambda^1\otimes\ad P),
\]
is the derivative of the smooth map, $\Pi_\Xi^\perp F^+$, at $A$ in the direction $a$, noting that $d_A^+a$ is the derivative of the smooth map, $A \mapsto F^+(A)$, at any connection $A$ in the direction $a$. In particular, the operators
\[
D(\Pi_\Xi^\perp F^+)(A) \in \sL\left(W^{1,2}(X;\Lambda^1\otimes\ad P), \Pi_\Xi^\perp(A)\right)
\]
are surjective\footnote{And thus admit right inverses, a construction favored by Donaldson and Kronheimer in \cite[Section 7.2.2]{DK} as a way to convert the first-order anti-self-dual partial differential equation into a zeroth-order nonlinear integral equation that may be solved using the Quantitative Inverse Function Theorem.} for any $[A] \in \sU$ with $F^+(A) = 0$. Hence, provided $\sU\cap (F^+)^{-1}(0) \neq \emptyset$, we can arrange that the operators $D(\Pi_\Xi^\perp F^+)(A)$ are surjective for any $[A] \in \sU$, by shrinking $\sU$ if necessary, and the set of solutions $A$ to the extended anti-self-dual equation,
\begin{equation}
\label{eq:Extended_ASD_equation_extrinsic_splitting}
\Pi_\Xi^\perp(A) F^+(A) = 0 \quad\text{(for $A\in\tilde\sU$)},
\end{equation}
is a smooth submanifold of $\tilde\sU$ by the Implicit Function Theorem. The quotient,
\[
M(P,g;\Xi) := (\Pi_\Xi^\perp F^+)^{-1}(0)/\Aut(P),
\]
is a finite-dimensional smooth submanifold of $\sU$ away from points $[A]$ where the stabilizer of $A$ in $\Aut(P)$ is larger than the center of $G$. 

When we solve the balancing equation,
\begin{equation}
\label{eq:Balancing_ASD_equation_extrinsic_splitting}  
\Pi_\Xi(A) F^+(A) = 0 \quad\text{(for $A\in\tilde\sU$ solving \eqref{eq:Extended_ASD_equation_extrinsic_splitting})},
\end{equation}
to produce
\[
\sU \cap M(P,g) = M(P,g;\Xi) \cap (\Pi_\Xi F^+)^{-1}(0),
\]
we recover an open neighborhood in the moduli space $M(P,g)$ of anti-self-dual connections on $P$.

Given an anti-self-dual 
`basepoint' connection, $A^\flat$ on $P$ with $F^+(A^\flat) = 0$, and constant $\eps\in(0,1]$, we define an open ball,
\[
B_{A^\flat}(\eps) := \left\{a \in W_{A^\flat}^{1,p}(X;\Lambda^1\otimes\ad P): \|a\|_{ W_{A^\flat}^{1,p}(X)} < \eps\right\}. 
\]
The Implicit Function Theorem yields a $\Stab(A^\flat)$-equivariant embedding,
\begin{equation}
\label{eq:ASD_Kuranishi_gluing_map_extrinsic_splitting}
\gamma: B_{A^\flat}(\eps) \cap \Ker \left(d_{A^\flat}^* + \Pi_\Xi^\perp(A^\flat) d_{A^\flat}^+\right) \ni a \mapsto \gamma(a) := A^\flat +a+\wp(a) \in \sA(P),
\end{equation}
where $\wp(a)$ is defined by solving \eqref{eq:Extended_ASD_equation_extrinsic_splitting}, that is,
\[
\Pi_\Xi^\perp( A^\flat +a+\wp(a))F^+(A^\flat +a+\wp(a)) = 0.
\] 
Here, $\gamma$ is the analogue of a gluing map in \cite{FrM, TauIndef} and
\begin{equation}
\label{eq:ASD_Kuranishi_obstruction_map_extrinsic_splitting}
\chi = \Pi_\Xi F^+ : B_{A^\flat}(\eps) \cap \Ker\left(d_{A^\flat}^* + \Pi_\Xi^\perp(A^\flat) d_{A^\flat}^+\right) \to
\Pi_\Xi^\perp(A^\flat)L^2(X;\Lambda^+\otimes\ad P)  
\end{equation}
is the analogue of an $\Stab(A^\flat)$-equivariant obstruction map in \cite{FrM, TauIndef}. For small enough $\eps \in (0,1]$, one may replace $\Pi_\Xi^\perp(A^\flat) d_{A^\flat}^+$ by $d_{A^\flat}^+$ in the definitions of $\gamma$ in \eqref{eq:ASD_Kuranishi_gluing_map_extrinsic_splitting} and $\chi$ in \eqref{eq:ASD_Kuranishi_obstruction_map_extrinsic_splitting} and recover the standard Kuranishi model \cite[Proposition 4.2.23]{DK} for an open neighborhood in $M(P,g)$ of a single point $[A^\flat]$. However, the advantage of the model \eqref{eq:ASD_Kuranishi_gluing_map_extrinsic_splitting}, \eqref{eq:ASD_Kuranishi_obstruction_map_extrinsic_splitting}  is that it also adapts to the case where $[A^\flat]$ varies in a smooth stratum $\sS$ of $M(P,g)$ and where there would be spectral flow for the Laplace operator, $d_{A^\flat}^+d_{A^\flat}^{+,*}$, so the bundle $\Ran\Pi_\mu(A^\flat)$ would not have constant rank along $\sS$ for any choice of $\mu \in (0,1]$, whereas --- as we explain in \cite[Section 3]{FL2a} --- we may choose $\Xi$ so that it has constant rank along $\sS$ and \eqref{eq:ASD_stabilizing_bundle_surjectivity} holds. The ball $B_{A^\flat}(\eps)$ would then be replaced by a relatively open tubular neighborhood, $N_\sS(\eps)$, of a closed smooth manifold, $\sS$.

\subsection{Extrinsic virtual neighborhoods for the moduli space of anti-self-dual connections and gluing}
\label{subsec:Extrinsic_virtual_neighborhood_moduli_space_anti-self-dual_connections_gluing}
We shall omit discussions of the more technical analytic aspects of gluing anti-self-dual connections or $\SO(3)$ monopoles described in \cite{FLKM1, FL3}, complex structures, $\CC^*$ or $S^1$ or other group actions, or group equivariances, noting that our gluing constructions are natural and that all group actions and equivariances extend from the setting of \cite{FL2a} (no bubbling) to those considered here in the presence of bubbling. 

We continue the notation and setup of Section \ref{subsubsec:Donaldson_Kronheimer_extrinsic_splitting} and first consider the problem of splicing a pair of families of connections, $A_i \in \tilde\sU_i$, on smooth principal $G$-bundles $P_i$ over closed, oriented, smooth Riemannian four-manifolds, $(X_i,g_i)$ for $i=1,2$, that are connected by a small cylinder (or neck or tube) diffeomorphic to $S^3\times (-1,1)$ and whose neck size is defined by a small neck parameter, $\lambda \in (0,1]$. Our development is inspired by that of Donaldson and Kronheimer \cite[Section 7.2]{DK}.

The open $\Aut(P)$-invariant neighborhoods, $\tilde\sU_i \subset \sA(P_i)$, project under the action of $\Aut(P_i)$ onto open neighborhoods $\sU_i$ in the configuration spaces $\sB(P_i) := \sA(P_i)/\Aut(P_i)$. We assume that each $[A_i] \in \sU_i$ obeys $\|F^+(A_i)\|_{L^2(X_i)} < \eps_i$, for suitably small $\eps_i \in (0,1]$. We let $\Xi_i$ denote finite-rank, product vector bundles over $\sU_i$ and their pullbacks to $\tilde\sU_i$ and let $\Pi_{\Xi_i}(A_i): L^2(X_i;\Lambda^+\otimes\ad P_i) \to \Xi_i(A_i)$ and $\Pi_{\Xi_i}^\perp(A_i) := 1 - \Pi_{\Xi_i}(A_i) = \Pi_{\Xi_i^\perp}(A_i)$ denote $L^2$-orthogonal projections for $i=1,2$, where $L^2(X_i;\Lambda^+\otimes\ad P_i) = \Xi_i(A_i) \oplus \Xi_i^\perp(A_i)$. As we explained in Section \ref{subsubsec:Donaldson_Kronheimer_extrinsic_splitting}, one can construct the $\Xi_i$ so that the operators,
\begin{equation}
\label{eq:ASD_gluing_stabilizing_bundle_surjectivity}
\Pi_{\Xi_i}^\perp(A_i):\Ran d_{A_i}^+ \to \Xi_i^\perp(A_i),
\end{equation}
are surjective for all $[A_i] \in \sU_i$ and $i=1,2$. The smooth maps,
\[
\tilde\sU_i \ni A_i \mapsto \Pi_{\Xi_i}^\perp(A_i) F^+(A_i) \in \Pi_{\Xi_i}^\perp(A_i)L^2(X_i;\Lambda^+\otimes\ad P_i),
\]
define smooth sections of the smooth vector bundles, $ \Xi_i^\perp \to \sU_i$. We can arrange that the operators $D(\Pi_{\Xi_i}^\perp F^+)(A_i)$ are surjective for any $[A_i] \in \sU_i$, by shrinking $\sU_i$ if necessary, and the zero locus,
\[
(\Pi_{\Xi_i}^\perp F^+)^{-1}(0) = \left\{A_i \in \sU_i: \Pi_{\Xi_i}^\perp F^+(A_i) = 0 \right\},
\]
is a finite-dimensional, smooth submanifold of $\sU_i$.

If $X := X_1\# X_2$ denotes the connected sum of $X_1$ and $X_2$ along a small cylinder, then we may use splicing to construct approximate solutions, $A' := A_1\#A_2$, of the extended anti-self-dual equation \eqref{eq:Extended_ASD_equation_extrinsic_splitting}
on the smooth principal $G$-bundle $P = P_1\# P_2$ over the closed, oriented, smooth Riemannian four-manifold $(X,g)$, where $g$ coincides with the $g_i$ outside the neck region, $\Xi := \Xi_1\oplus\Xi_2$ is the induced finite-rank product vector bundle over $X$, the operators
\[
\Pi_\Xi^\perp(A):\Ran d_A^+ \to \Xi^\perp(A)
\]
are surjective (by construction) for all $[A] \in \sU$, where $L^2(X;\Lambda^+\otimes\ad P) = \Xi(A)\oplus \Xi^\perp(A)$ and $\Pi_\Xi^\perp(A) := 1 - \Pi_\Xi(A) = \Pi_{\Xi^\perp}(A)$, and $\sU \subset \sB(P)$ is an open neighborhood containing $\sU_1\#\sU_2$, the set of all gauge-equivalence classes of the spliced connections, $A_1\# A_2$, and $\tilde\sU \subset \sA(P)$ is its $\Aut(P)$-invariant pullback. We now apply the (Quantitative) Implicit Function Theorem to construct solutions $a := \wp(A') \in W_{A'}^{1,p}(X;\Lambda^1\otimes\ad P)$ (for small enough parameters $\eps_1, \eps_2, \lambda$) to the extended anti-self-dual equation \eqref{eq:Extended_ASD_equation_extrinsic_splitting}, namely
\[
\Pi_\Xi^\perp F^+(A' + a) = 0 \quad\text{(for $A' \in \tilde\sU$)}.
\] 
Just as in the case of the component sections, $\Pi_{\Xi_i}^\perp F^+$, the zero locus,
\[
M(P,g;\Xi) := \left\{[A] \in \sU: \Pi_\Xi^\perp F^+(A) = 0 \right\},
\]
is a finite-dimensional, smooth submanifold of $\sU$. By solving the balancing equation \eqref{eq:Balancing_ASD_equation_extrinsic_splitting}, namely
\[
\Pi_\Xi F^+(A) = 0 \quad\text{(for $A\in\tilde\sU$ solving \eqref{eq:Extended_ASD_equation_extrinsic_splitting})},
\]
to produce
\[
\sU \cap M(P,g) = M(P,g;\Xi) \cap (\Pi_\Xi F^+)^{-1}(0),
\]
we recover an open neighborhood in the moduli space $M(P,g)$ of anti-self-dual connections on $P$.

Suppose now that we are given a pair of anti-self-dual basepoint connections, $A_i^\flat$ on $P_i$ with $F^+(A_i^\flat) = 0$, and constants $\eps_i\in(0,1]$ for $i=1,2$. For small enough $\eps_1,\eps_2,\lambda$, the Quantitative Implicit Function Theorem yields a $\Stab(A_1^\flat)\times \Stab(A_2^\flat)$-equivariant embedding,
\begin{multline}
\label{eq:ASD_tube_gluing_map_extrinsic_splitting}
\gamma: B_{A_1^\flat}(\eps_1) \cap \Ker \left(d_{A_1^\flat}^* + \Pi_{\Xi_1}^\perp(A_1^\flat) d_{A_1^\flat}^+\right)
\times B_{A_2^\flat}(\eps_2) \cap \Ker \left(d_{A_2^\flat}^* + \Pi_{\Xi_2}^\perp(A_2^\flat) d_{A_2^\flat}^+\right)
\\
\times \Isom_G(P_{x_1},P_{x_2}) \ni (a_1,a_2,h) \mapsto A'+a+\wp(a) \in \sA(P),
\end{multline}
where $\wp(A')$ is defined by solving \eqref{eq:Extended_ASD_equation_extrinsic_splitting}, that is,
\[
\Pi_\Xi^\perp F^+(A'+\wp(A')) = 0,
\]
and the approximately anti-self-dual connection, $A'$ on $P$, is defined by splicing $A_i^\flat+a_i$ for $i=1,2$. Here, $\gamma$ is the gluing map and
\begin{multline}
\label{eq:ASD_tube_obstruction_map_extrinsic_splitting}
\chi: B_{A_1^\flat}(\eps_1) \cap \Ker \left(d_{A_1^\flat}^* + \Pi_{\Xi_1}^\perp(A_1^\flat) d_{A_1^\flat}^+\right)
\times B_{A_2^\flat}(\eps_2) \cap \Ker \left(d_{A_2^\flat}^* + \Pi_{\Xi_2}^\perp(A_2^\flat) d_{A_2^\flat}^+\right)
\\
\times \Isom_G(P_{x_1},P_{x_2}) \ni (a_1,a_1,h)
\\
\mapsto \Pi_\Xi F^+(A'+\wp(A')) \in \Pi_\Xi^\perp(A'+\wp(A'))L^2(X;\Lambda^+\otimes\ad P)
\end{multline}
is the $\Stab(A_1^\flat)\times \Stab(A_2^\flat)$-equivariant obstruction map. The existence and properties of $\gamma$ in \eqref{eq:ASD_tube_gluing_map_extrinsic_splitting} and $\chi$ in \eqref{eq:ASD_tube_obstruction_map_extrinsic_splitting} are established by Donaldson and Kronheimer in \cite[Theorems 7.2.62 and 7.2.63]{DK}. The gluing map $\gamma$ in \eqref{eq:ASD_tube_gluing_map_extrinsic_splitting} descends to the quotient by $\Stab(A_1^\flat)\times \Stab(A_2^\flat)$ to give a continuous embedding onto a relatively open subset of $\sB(P)$ (respectively, smooth embedding away from singularities) and restricts to a continuous embedding from $\chi^{-1}(0)$ onto an open subset of $M(P,g)$ (respectively, smooth embedding away from singularities). 

In the case where $[A_i^\flat]$ varies in smooth strata $\sS_i$ of $M(P_i,g)$ and where there would be spectral flow for the Laplace operators, $d_{A_i^\flat}^+d_{A_i^\flat}^{+,*}$, so the bundles $\Ran\Pi_{\mu_i}(A_i^\flat)$ would not have constant rank along $\sS_i$ for any choice of $\mu_i \in (0,1]$, we may choose $\Xi \cong \Xi_1\oplus \Xi_2$ so that it has constant rank along $\sS_1\times\sS_2$ and \eqref{eq:ASD_stabilizing_bundle_surjectivity} holds. The balls $B_{A_i^\flat}(\eps)$ would then be replaced by a relatively open tubular neighborhoods, $N_{\sS_i}(\eps)$, of a closed smooth manifolds $\sS_i$.

In applications to the proof of the Kotschick--Morgan Conjecture in this monograph, the above-mentioned strata, $\sS$ or $\sS_1, \sS_2$, in this context comprise isolated points corresponding to reducible anti-self-dual connections. However, in our application of these ideas to the proof of the Witten Conjecture, we must allow for positive-dimensional strata of reducible $\SO(3)$ monopoles, that is, positive-dimensional moduli spaces of Seiberg--Witten monopoles.

When the neck parameter, $\lambda$, tends to zero, the pair $(\gamma,\chi)$ converges continuously to $(\gamma_1\times\gamma_2, \chi_1\times\chi_2)$. In typical applications (for example, see Donaldson and Kronheimer \cite[Theorem 8.2.3 and Proposition 8.2.4]{DK} or our application to the proof of the Kotschick--Morgan Conjecture), we choose $(X_2,g_2)$ to be $(S^4,g_\round)$, the four-sphere equipped with its standard round metric, $g_\round$, of radius one. The gluing theory summarized in this subsection yields an open neighborhood of a boundary point $([A_1,q_1],[A_2,q_2])$ in the \emph{bubble-tree compactification} of $M(P_1,g_1)$, where $q_1 \in P_1|_{x_1}$ and $q_2 \in P_2|_{s}$ and the point $x_1 \in X_1$ is identified with the south pole $s \in S^4$.
By `forgetting' the additional data on $S^4$ (aside from $c_2(P_2)$ when $G=\SU(2)$, which defines the multiplicity of the point $x_1$)
and the frame $q_1 \in P_1|_{x_1}$ when the neck parameter $\lambda$ becomes zero, the gluing theory yields an open neighborhood of a boundary point $([A_1],x_1)$ in the \emph{Uhlenbeck compactification} of $M(P,g)$.

Lastly, the gluing theory for solutions to the anti-self-dual equation on the principal $G$-bundle $P=P_1\# P_2$ over the connected sum $X=X_1\#X_2$ summarized in this subsection extends to the case of \emph{multiple connected sums} (see Donaldson and Kronheimer  \cite[Chapter 8]{DK} and Feehan and Leness \cite{FL3}) and \emph{tree connected sums} (see Feehan \cite{FeehanGeometry}, Peng \cite{Peng_1995, Peng_1996}, and Taubes \cite{TauFrame}) required when building parameterizations for open neighborhoods of points in the boundary of the bubble-tree (and hence Uhlenbeck) compactification of $M(P_1,g_1)$.

\subsection{Construction of a virtual neighborhood for the moduli space of $\SO(3)$ monopoles near a lower-level singular stratum of $\SO(3)$ monopoles}
\label{subsec:Virtual_neighborhood_moduli_space_SO3_monopoles_lower_level}
The gluing theory described in Section \ref{subsec:Extrinsic_virtual_neighborhood_moduli_space_anti-self-dual_connections_gluing} for solutions to the anti-self-dual equation generalizes \emph{formally} to the case of the $\SO(3)$-monopole equations \eqref{eq:PerturbedSO3MonopoleEquations}. However, the analysis required to prove that the gluing map $\bgamma_\sM$ and obstruction map $\bar{\bchi}$ have all the properties asserted by Hypothesis \ref{hyp:Gluing} is considerable and the details of that analysis do \emph{not} extend in a straightforward manner from those previously encountered in the case of the anti-self-dual equation (for example, Donaldson and Kronheimer \cite{DK}, Feehan \cite{FeehanGeometry}, Feehan and Leness \cite{FLKM1}, Morgan and Mrowka \cite{MorganMrowkaTube}, Mrowka \cite{MrowkaThesis}, and Taubes \cite{TauSelfDual, TauIndef, TauFrame}). We summarized the principal new analytical difficulties in our article \cite{FL3} and address them fully in \cite{Feehan_Leness_monopolegluingbook}.

In this subsection, we outline a method for extending the construction of the stabilizing bundle $\Xi$ for $\Ran d_A^+$ to give a stabilzing bundle $\Xi$ for $\Ran d_{A,\Phi}^1$ and hence define a splitting of the $\SO(3)$-monopole equations \eqref{eq:PerturbedSO3MonopoleEquations}. If $(X,g) = (S^4,g_\round)$ and $P$ is a principal $\SU(2)$-bundle over $S^4$ and $(\rho,W)$ is the standard \spinc structure over $S^4$, we recall from \cite[Remark 4.6]{FL1} that if $(A,\Phi)$ is an $\SO(3)$ monopole on $(P,W^+)$, then $\Phi \equiv 0$ and $F_A^+ = 0$. The linearization, $d_{A,0}^1 = (D\fS)(A,0)$, of the $\SO(3)$-monopole equations, $\fS(A,\Phi) = 0$, in  \eqref{eq:PerturbedSO3MonopoleEquations} at an $\SO(3)$ monopole $(A,0)$ over $S^4$ is given by \cite[Equation (2.49)]{FL2a}, namely,
\[
(D\fS)(A,0) = (d_A^+, D_A),
\]
where $D_A:C^\infty(S^4,W^+\otimes E) \to C^\infty(S^4,W^-\otimes E)$ is the coupled Dirac operator (with $L^2$-adjoint $D_A^*:C^\infty(S^4,W^-\otimes E) \to C^\infty(S^4,W^+\otimes E)$) and, as usual, $d_A^+:C^\infty(S^4;\Lambda^1\otimes\ad P) \to C^\infty(S^4;\Lambda^+\otimes\ad P)$ is the linearization of the anti-self-dual equation, $F_A^+ = 0$, at an anti-self-dual connection, $A$, and $E = P\times_{\std}\CC^2$. According to \cite[Lemma 8.12]{FL3}, when $F_A^+ = 0$, the kernel of $D_A^*D_A$ (and thus $D_A$) is zero, the real dimension of the kernel of $D_AD_A^*$ (and thus $D_A^*$) is equal to twice the complex index of $D_A^*$, namely $\Ind_\CC D_A^* = \dim_\CC\Ker D_A^* - \dim_\CC\Ker D_A = c_2(E)$ by \cite[Proof of Proposition 2.28]{FL1}, and the least positive eigenvalue of the Dirac Laplacian, $D_AD_A^*$, is equal to $3$. Hence, as $A$ varies over $\tilde M(P,g_\round)$, the kernels, $\Ker D_A^*$ (equivalently, the cokernels, $\Coker D_A$), define a complex vector bundle, $\Ker \bD^* \cong \Coker \bD$, of complex rank $c_2(E)$ over $M(P,g_\round)$. As usual, $\Coker d_A^+ = \{0\}$ when $A$ is an anti-self-dual connection over $S^4$ by \cite[Proof of Theorem 6.1]{AHS} and thus $\Coker d_{A,0}^1 \cong \Coker D_A$.

We may extend the basic setting of Section \ref{subsec:Extrinsic_virtual_neighborhood_moduli_space_anti-self-dual_connections_gluing} to the case of $\SO(3)$ monopoles and consider splicing a family of $\SO(3)$ monopoles, $(A_1,\Phi_1)$, on $(P_1,W_1^+)$ over $(X_1,g_1) = (X,g)$, and a family of zero-section $\SO(3)$ monopoles, $(A_2,0)$, on $(P_2,W_2^+)$ over $(X_2,g_2) = (S^4,g_\round)$ to form a family of approximate $\SO(3)$ monopoles, $(A',\Phi') = (A_1\#A_2,\Phi_1\# 0)$, on $(P,W^+) = (P_1\#P_2,W_1^+\# W_2^+)$ over the Riemannian connected sum, $(X\# S^4,g\# g_\round)$. We may choose the stabilizing bundle $\Xi_1$ for $\Ran d_{A_1,\Phi_1}^1$ described in Section \ref{subsec:Virtual_neighborhood_moduli_space_SO3_monopoles_top_level} and choose $\Xi_2 = \Coker \bD$ to form a stabilizing bundle, $\Xi = \Xi_1\oplus\Xi_2$, for $\Ran d_{A',\Phi'}^1$ over the connected sum, $X\# S^4$. 

As in Section \ref{subsec:Extrinsic_virtual_neighborhood_moduli_space_anti-self-dual_connections_gluing}, the preceding construction extends to the case of multiple connected and tree connected sums. One difference in the case that all summands in a tree except $(X_1,g_1) = (X,g)$ comprise copies of $(S^4,g_\round)$ is that we may use the Taubes small-eigenvalue decomposition to construct obstruction bundles corresponding to connected sums of copies of $(S^4,g_\round)$ since there is no spectral flow for the Dirac Laplacians, $D_AD_A^*$, over $(S^4,g_\round)$.

\chapter{Link of an ideal Seiberg--Witten moduli space}
\label{chap:Link}
We now use the Thom--Mather structure defined in
Section \ref{sec:TMStr} to construct a
\emph{virtual link},
$\bar\bL^{\vir}_{\ft,\fs}$, given by the boundary of a
neighborhood of $M_{\fs}\times\Sym^\ell(X)$ in
$\bar\sM^{\vir}_{\ft,\fs}/S^1$.
The \emph{link}, $\bar\bL_{\ft,\fs}$, appearing in \eqref{eq:RawCobordismSum},
of $M_{\fs}\times\Sym^\ell(X)$ in
$\bar\sM_\ft/S^1$
is then the intersection of the
virtual link $\bar\bL^{\vir}_{\ft,\fs}$ with
the zero locus of the obstruction section $\bar\bchi$ appearing in
Hypothesis \ref{hyp:Gluing}.

We construct the link in Section \ref{sec:AmbientLink} using the
tubular distance functions $\vec t(\ft,\fs,\sP_i)$
defined in \eqref{eq:DefineTubularDistanceFunction}. We use our understanding of
the overlap maps to show that the
virtual link
$\bar\bL^{\vir}_{\ft,\fs}$ can be decomposed into
closed subspaces,
enumerated by the strata of $\Sym^\ell(X)$.
Each of these closed subspaces is a smoothly-stratified
space.
From this definition and a discussion of the orientation of the link,
we prove Theorem \ref{thm:CobordismThm}
and thus the cobordism sum \eqref{eq:RawCobordismSum}.

Our understanding of the overlaps of the spaces $\bar\sU(\ft,\fs,\sP)$
allows us to prove in Section \ref{sec:FiberBundle}
that each of the aforementioned closed subspaces of the
virtual link admits a fiber bundle structure.
In Section \ref{sec:Boundaries}, we describe the intersections
of these subspaces of the link and the interaction of the intersections
with the fiber bundle structures.

\section{Definition of the link of an ideal Seiberg--Witten moduli space}
\label{sec:AmbientLink}
We define the link $\bL_{\ft,\fs}$
of an ideal Seiberg--Witten moduli space
by first
constructing a
virtual link in Section \ref{subsec:AmbientLink}.
This virtual link is the union of closed subspaces, as described in
\eqref{eq:AmbientLink}, \eqref{eq:InstantonLinkComponent}, and \eqref{eq:DefineLinkStratum}.
We discuss the orientations of the link in Section
\ref{subsubsec:Orient}
as needed to define the intersection numbers appearing in the terms of the cobordism sum in
\eqref{eq:RawCobordismSum} in Theorem
\ref{thm:CobordismThm}.
Then in Section \ref{sec:FiberBundle}, we
give a more detailed description of the closed subspaces of the link defined in
\eqref{eq:DefineLinkStratum},
showing how they admit a fiber bundle structure which will be used in
the pushforward-pullback computations of Chapter \ref{chap:Comp}.

\subsection{The virtual link of an ideal Seiberg--Witten moduli space}
\label{subsec:AmbientLink}
The
virtual link is defined as the
boundary of a neighborhood of $M_{\fs}\times\Sym^\ell(X)$ in
$\bar\sM^{\vir}_{\ft,\fs}/S^1$. The subspace
$M_{\fs}\times\Sym^\ell(X)$ of $\bar\sM^{\vir}_{\ft,\fs}/S^1$ is
the intersection of the zero-locus of the function defined by
a norm on the fibers of the vector bundle $N_{\ft(\ell),\fs}(\delta)\to
M_{\fs}$,
\begin{equation}
\label{eq:NTubularDist} t_N:\bar\sM^{\vir}_{\ft,\fs}/S^1\to
[0,\delta],
\end{equation}
with the union of the zero-loci of the tubular distance function
$\vec t(\ft,\fs,\sP_j)$ defined on $\bar\sU(\ft,\fs,\sP_j)$ in
\eqref{eq:DefineTubularDistanceFunction}. We will define the
virtual link as the union of two
subspaces,
\begin{equation}
\label{eq:AmbientLink}
\bar\bL^{\vir}_{\ft,\fs} =
\bar\bL^{\vir,s}_{\ft,\fs}\cup \bar\bL^{\vir,i}_{\ft,\fs}.
\end{equation}
The first subspace, $\bar\bL^{\vir,s}_{\ft,\fs}$, is defined to be
the
codimension-one subspace,
\begin{equation}
\label{eq:DefineSWComponentOfLink} \bar\bL^{\vir,s}_{\ft,\fs} =
t_N^{-1}(\delta) \subset \bar\sM^{\vir}_{\ft,\fs}/S^1,
\end{equation}
in the sense that the intersection of $\bar\bL^{\vir,s}_{\ft,\fs}$ with each smooth stratum of
$\bar\sM^{\vir}_{\ft,\fs}/S^1$ is a codimension-one submanifold.
To define the instanton component of the link,
$\bar\bL^{\vir,i}_{\ft,\fs}$,
first enumerate the strata of $\Sym^\ell(X)$
with partitions $\sP_0,\dots,\sP_r$ of $N_\ell$ as done in \S \ref{subsec:EnumStrata}.
Then, define
\begin{equation}
\label{eq:InstantonLinkComponent}
\bar\bL^{\vir,i}_{\ft,\fs} :=
\bigcup_{j=0}^r \ \bar\bL^{\vir,i}_{\ft,\fs}(\sP_j)
\end{equation}
with each subspace $\bar\bL^{\vir,i}_{\ft,\fs}(\sP_j)$ given as the pre-image
of a codimension-one submanifold with corners
(in the sense of \cite{Joyce_2012})
in $[0,1]^{|\sP_j|}$
under $\vec t(\ft,\fs,\sP_j)$.
We now define the
subspaces, $\bar\bL^{\vir,i}_{\ft,\fs}(\sP_j)$, of the instanton link.

Assign a small,
generic constant $\eps_j$  to each such partition $\sP_j$ with
$\eps_j>\eps_k$ for $j<k$ in a manner to be specified in Lemma \ref{lem:BoundaryOfNeigh}.
Then define a
subspace of the instanton component of the link by
\begin{equation}
\label{eq:DefineLinkStratum}
\bar\bL^{\vir,i}_{\ft,\fs}(\sP_j)
:=
\bar\sU(\ft,\fs,\sP_j)/S^1\cap \vec t(\ft,\fs,\sP_j)^{-1}(\rd \bar D(\sP,\eps_j))
\setminus \bigcup_{k\neq j} \vec t(\ft,\fs,\sP_k)^{-1}(D(\sP,\eps_k)).
\end{equation}
Recall from \eqref{eq:DefineTubularDistanceFunction} that the
map $\vec t(\ft,\fs,\sP_j)$
is defined by the
functions $\tilde \la_P$ on the connections over $S^4$ making up the fiber of $\bar\Gl(\ft,\fs,\sP)\to \Si(X^\ell,\sP)$.
Because of the separating condition \eqref{eq:SeparatingCondition1} on the scales of these connections over $S^4$
in the definition of $\sU(\ft,\fs,\sP_j)$, the inclusion
$$
\pi(\ft,\fs,\sP_j)^{-1}(\bx)\cap \vec t(\ft,\fs,\sP_j)^{-1}(\bar D(\sP,\eps_j))
\subseteq
\pi(\ft,\fs,\sP_j)^{-1}(\bx)\cap\sU(\ft,\fs,\sP_j)
$$
need not hold as $\bx\in\Si(X^\ell,\sP_j)$ approaches lower strata in $\Sym^\ell(X)$ and hence
the subspace $\vec t(\ft,\fs,\sP_j)^{-1}(\bar D(\sP,\eps_j))$ of $\bar\sU(\ft,\fs,\sP_j)$
need not define a fiber bundle.
The following lemma
addresses this issue.

\begin{lem}
\label{lem:SmallConst}
Let $\sP_j$ be one of the partitions of $N_\ell$ chosen before \eqref{eq:InstantonLinkComponent}.
For any compact
subset
$K\Subset \Si(X^\ell,\sP_j)$, there is a
positive
constant
$\eps_j$ such that for all
positive
$\eps\le \eps_j$,
the set
$$
\vec t(\ft,\fs,\sP_j)^{-1}(\bar D(\sP_j,\eps) \cap
\pi(\ft,\fs,\sP_j)^{-1}\left( N_{\ft(\ell),\fs}(\delta)/S^1\times K\right)
$$
is a proper subspace of
$$
\bar\sU(\ft,\fs,\sP_j)/S^1 \cap \pi(\ft,\fs,\sP_j)^{-1}\left(
N_{\ft(\ell),\fs}(\delta)/S^1\times K\right).
$$
\end{lem}

\begin{proof}
The lemma follows immediately from the compactness of $K$.
\end{proof}

We now describe how to chose the decreasing sequence of
constants appearing in
\eqref{eq:DefineLinkStratum}.
We say that a codimension-one submanifold $N$ of a manifold $M$
is \emph{smoothly collared} if there is a smooth embedding, $N\times (-\eps,\eps)\to M$, whose restriction to $N\times \{0\}$ is
the inclusion $N\to M$.

\begin{lem}
\label{lem:BoundaryOfNeigh}
There is a decreasing sequence of positive
constants $\eps_0>\eps_1>\cdots>\eps_r$ such that
the
virtual link
$\bar\bL^{\vir}_{\ft,\fs}$ defined in \eqref{eq:AmbientLink},
\eqref{eq:DefineSWComponentOfLink}, and
\eqref{eq:InstantonLinkComponent}
by these constants has the following properties:
\begin{enumerate}
\item
\label{item:BoundaryOfNeigh1}
$\bar\bL^{\vir}_{\ft,\fs}$ is the boundary of a
neighborhood of $M_{\fs}\times \Sym^\ell(X)$ in
$\bar\sM^{\vir}_{\ft,\fs}/S^1$.
\item
\label{item:BoundaryOfNeigh2}
The intersection of $\bar\bL^{\vir}_{\ft,\fs}$
with
each stratum $\sS$ of $\bar\sM^{\vir}_{\ft,\fs}/S^1$ is  a
Whitney-stratified space
whose top stratum is a codimension-one, collared submanifold of $\sS$.
\item
\label{item:BoundaryOfNeigh3}
The intersection of $\bar\bL^{\vir}_{\ft,\fs}$ with the top stratum of $\bar\sM^{\vir}_{\ft,\fs}/S^1$
is a topological manifold.
\end{enumerate}
\end{lem}

\begin{proof}
We work by induction on the partitions $\sP_j$
of $N_\ell$ used to enumerate the strata of $\Sym^\ell(X)$
before \eqref{eq:DefineLinkStratum}.
The lowest stratum, $\Si(X^\ell,\sP_0)$, of
$\Sym^\ell(X)$ is necessarily compact and so there is a small positive constant $\eps_0$ such that Lemma \ref{lem:SmallConst} applies with $K = \Si(X^\ell,\sP_0)$.

Assume that we have defined $\eps_0,\dots,\eps_{k-1}$
such  that for $0\le j\le k-1$,
$\eps_j$ is a constant sufficiently small that Lemma \ref{lem:SmallConst}
applies to  $K_j\Subset \Si(X^\ell,\sP_j)$ where
$K_j$ is constructed as in
Lemmas \ref{lem:DecomposeSymProduct} and \ref{lem:CompactSubsetsOfSi}
using the constants
$\eps_i$ with $i<j$.
If we define $K_k\subset \Si(X^\ell,\sP_k)$ as in
Lemmas \ref{lem:DecomposeSymProduct} and \ref{lem:CompactSubsetsOfSi} applied to the constants
$\eps_0,\dots,\eps_{k-1}$, then Lemma \ref{lem:CompactSubsetsOfSi}
implies that $K_k$ is compact. We now
select $\eps_k<\eps_{k-1}$ so Lemma \ref{lem:SmallConst}
applies to $K_k$.

We note that for any $\eps>0$, the set
$\vec t(\ft,\fs,\sP)^{-1}(\bar D(\sP,\eps))$ is a
closed neighborhood
of $N_{\ft(\ell),\fs}(\delta)/S^1\times \Si(X^\ell,\sP)$ in $\bar\sU(\ft,\fs,\sP)$.
Hence,
\begin{multline}
\label{eq:LocalNghOfReducibles}
\bigcup_j\ \vec t(\ft,\fs,\sP_j)^{-1}(\bar D(\sP_j,\eps_j))
\\
=
\bigcup_j
\left(
\vec t(\ft,\fs,\sP_j)^{-1}(\bar D(\sP_j,\eps_j))
\left.\setminus
\bigcup_{i\neq j}
\vec t(\ft,\fs,\sP_i)^{-1}(D(\sP_i,\eps_i))
\right.
\right)
\end{multline}
defines a
closed neighborhood of $N_{\ft(\ell),\fs}(\delta)/S^1\times \Sym^\ell(X)$ in
$\cup_j\bar\sU(\ft,\fs,\sP_i)/S^1$.  The link \eqref{eq:AmbientLink}
is the boundary of the preceding neighborhood, proving Item \eqref{item:BoundaryOfNeigh1}.

We prove that the intersection of the space \eqref{eq:LocalNghOfReducibles} with each stratum of
$\bar\sM^{\vir}_{\ft,\fs}/S^1$ is a smooth manifold with corners using the following criterion.
If $Y$ is a smooth manifold and $f_j:Y\to \RR$ is a smooth function for $1\le j\le p$, then
the subspace,
$$
Z:=\{y\in Y: f_j(y)\le \eps_j\ \text{for $j=1,\dots,p$}\},
$$
for given constants $\eps_1,\dots,\eps_p$ will be a manifold with corners if for every collection of
indices $1\le j_1,\dots,j_a\le p$ and every $y\in Y$ satisfying $f_{j_k}(y)=\eps_{j_k}$, the derivatives
$(df_{j_1})_y,\dots, (df_{j_a})_y$ are linearly independent.

The level sets playing the role of $f_j^{-1}(\eps_j)$ for  the space \eqref{eq:LocalNghOfReducibles} are
the subspaces of $\bar\sU(\ft,\fs,\sP)/S^1$ defined by
\begin{equation}
\label{eq:BoundaryHypersurface}
H_P(\sP,\eps):=
\vec t(\ft,\fs,\sP)^{-1}\left( \{(\la_Q)_{Q\in\sP}: \la_P=\eps\}\subset\RR^{|\sP|}\right).
\end{equation}
From the definition  \eqref{eq:DefineTubularDistanceFunction}
of $\vec t(\ft,\fs,\sP)$ in terms of the functions $\tilde\la_P$, the intersection of
$H_P(\sP,\eps)$ with open sets of $\sU(\ft,\fs,\sP)$ of the form
\begin{equation}
\label{eq:OnePartitionOpenSet}
N\times\sO_P\times \prod_{P\in\sP}\barM^{s,\natural}_{\spl,|P|}(\delta),
\end{equation}
where $N\subset N_{\ft(\ell),\fs}(\delta)$ and $\sO_\sP\subset\Delta^\circ(X^\ell,\sP)$ are open subsets,
is given by the pre-image of $\tilde\la_P^{-1}(\eps)\subset \barM^{s,\natural}_{\spl,|P|}(\delta)$ under
the obvious projection.   Hence, for generic values of $\eps$, the intersection of $H_P(\sP,\eps)$
with each stratum of $\bar\sM^{\vir}_{\ft,\fs}/S^1$ will be a smooth, collared, codimension-one submanifold
of that stratum.
For any collection  $P_1,\dots,P_s\in\sP$,  at any point in the intersection of
\begin{equation}
\label{eq:OnePartitionIntersection}
H_{P_1}(\sP,\eps)\cap\dots\cap H_{P_s}(\sP,\eps)
\end{equation}
with an open set of the form \eqref{eq:OnePartitionOpenSet},
the pullbacks of the derivatives $d\tilde\la_{P_1},\dots, d\tilde\la_{P_s}$ will be linearly independent at points in
\eqref{eq:OnePartitionIntersection} by the product structure of \eqref{eq:OnePartitionOpenSet}.

To prove that the above criterion applies to more general intersections of the form
\begin{equation}
\label{eq:HypersurfaceIntesections}
H_{P_{k_1}}(\sP_{k_1},\eps_{k_1})
\cap
\dots
\cap
H_{P_{k_a}}(\sP_{k_a},\eps_{k_a}),
\end{equation}
we argue by induction on $a$. For $a=1$, the result follows by the argument above
given for intersections of the form \eqref{eq:OnePartitionIntersection}.
Now assume $a>1$.
After re-ordering the indices, we can assume that $\sP_{k_1}\le\sP_{k_2}\le\dots\le\sP_{k_a}$.
The case where $\sP_{k_1}=\sP_{k_2}=\dots=\sP_{k_a}$ has been done in \eqref{eq:OnePartitionIntersection}.
Then assume that $\sP_{k_1}=\dots=\sP_{k_b}<\sP_{k_{b+1}}$.
By Theorem \ref{thm:GlobalSplicingDataOverlaps}, points in $\sU(\ft,\fs,\sP_{k_b})\cap\sU(\ft,\fs,\sP_{k_{b+1}})$ have
neighborhoods which are smoothly-stratified diffeomorphic to an open subspace of (see \eqref{eq:OverlapGluingDataSpace})
\begin{equation}
\label{eq:LocalPicOfOverlap}
N\times\sO \times \prod_{P\in\sP_{k_b}}
\left(
 \Delta^\circ(Z_{|P|}(\delta_P),(\sP_{k_{b+1}})_P) \times \prod_{P'\in(\sP_{k_{b+1}})_P}  \barM^{s,\natural}_{\spl,|P'|}(\delta)
\right),
\end{equation}
where $N\subset N_{\ft(\ell),\fs}(\delta)$ and $\sO_\sP\subset\Delta^\circ(X^\ell,\sP_{k_b})$ are open subsets.
The space \eqref{eq:LocalPicOfOverlap} is smoothly-stratified diffeomorphic to a product,
$\sU_\Si\times\sU_M$, where
\begin{align*}
\sU_\Si & :=N\times\sO \times \prod_{P\in\sP_{k_b}} \Delta^\circ(Z_{|P|}(\delta_P),(\sP_{k_{b+1}})_P),
\\
\sU_M & := \prod_{P\in\sP_{k_b}}\prod_{P'\in (\sP_{k_{b+1}})_P}  \barM^{s,\natural}_{\spl,|P'|}(\delta).
\end{align*}
We write $p_\Si$ and $p_M$ for the projections maps from the space \eqref{eq:LocalPicOfOverlap} to $\sU_\Si$ and $\sU_M$ respectively.

From Lemma \ref{lem:TubularDistanceFunctionOnOverlaps2},  for $u>b$ the restriction of $\vec t(\ft,\fs,\sP_{k_u})$ to an open subspace of the form \eqref{eq:LocalPicOfOverlap} equals the pullback
by $p_M$ of a
map on $\sU_M$.
By induction, the derivatives of the components of $\vec t(\ft,\fs,\sP_{k_u})$ defining $H_{P_{k_u}}(\sP_{k_u},\eps_{k_u})$
are linearly independent at points in the intersection
\begin{equation}
\label{eq:HypersurfaceIntesections1}
H_{P_{k_{b+1}}}(\sP_{k_{b+1}},\eps_{k_{b+1}})
\cap
\dots
\cap
H_{P_{k_a}}(\sP_{k_a},\eps_{k_a}).
\end{equation}
If $P_{k_v}\notin\sP_{k_{b+1}}$ for  all $1\le v\le b$, then from
Item \eqref{item:TMTubularDistanceR43} of Lemma \ref{lem:TMTubularDistanceR4}
and the second equality of \eqref{eq:XProjectionRelationsOnOverlap},
the
maps $\vec t(\ft,\fs,\sP_{k_v})$ are equal to the pullbacks of
functions on $\sU_\Si$ by $p_\Si$.
Again by induction, the derivatives of the components of $\vec t(\ft,\fs,\sP_{k_v})$ defining $H_{P_{k_v}}(\sP_{k_v},\eps_{k_v})$
are linearly independent at points in the intersection
\begin{equation}
\label{eq:HypersurfaceIntesections2}
H_{P_{k_1}}(\sP_{k_1},\eps_{k_1})
\cap
\dots
\cap
H_{P_{k_b}}(\sP_{k_a},\eps_{k_b}).
\end{equation}
Consequently,  the derivatives of the functions defining the codimension-one submanifolds
in \eqref{eq:HypersurfaceIntesections}
are linearly independent at points in the intersection
\eqref{eq:HypersurfaceIntesections} by the product structure $\sU_\Si\times\sU_M$ of
\eqref{eq:LocalPicOfOverlap}.

If $P_{k_v}\in \sP_{k_{u}}$ and $P_{k_v}=P_{k_u}$ for  $1\le v\le b<u$ then,
because $\eps_{k_v}\neq \eps_{k_u}$,
the intersection \eqref{eq:HypersurfaceIntesections} is empty.

If $P_{k_v}\in \sP_{k_{u}}$ and $P_{k_v}\neq P_{k_u}$ for  $1\le v\le b<u$ then,
on the open set \eqref{eq:LocalPicOfOverlap},
the $P_{k_v}$-th components of $\vec t(\ft,\fs,\sP_{k_v})$ and $\vec t(\ft,\fs,\sP_{k_u})$ are equal.  Thus, we can
replace $H_{P_{k_v}}(\sP_{k_v},\eps_{k_v})$ with $H_{P_{k_v}}(\sP_{k_u},\eps_{k_v})$.
Continue doing these replacements until the intersection is reduced to the previous case considered.
This completes the induction and shows that the intersections \eqref{eq:HypersurfaceIntesections} satisfy our criterion,
proving that each stratum of the space
\eqref{eq:LocalNghOfReducibles} is a smooth manifold with corners.

The
virtual link
$\bar\bL^{\vir}_{\ft,\fs}$  is the boundary of the space \eqref{eq:LocalNghOfReducibles}.
Hence, the intersection of $\bar\bL^{\vir}_{\ft,\fs}$
with the top stratum of $\bar\sM^{\vir}_{\ft,\fs}/S^1$ is the boundary of a smooth manifold with corners
and hence a topological manifold, proving Item \eqref{item:BoundaryOfNeigh3}.

From the preceding argument, we see that for generic, decreasing
values of $\eps_1,\dots,\eps_r$,
the intersection of
\begin{equation}
\label{eq:CompactDiskBundle}
\bar\sU(\ft,\fs,\sP_j)/S^1\cap \vec t(\ft,\fs,\sP_j)^{-1}(\bar D(\sP,\eps_j))
\left.\setminus \bigcup_{k\neq j} \vec t(\ft,\fs,\sP_k)^{-1}(D(\sP,\eps_k))\right.
\end{equation}
with each stratum of $\bar\sM^{\vir}_{\ft,\fs}/S^1$ is a smooth manifold with corners.
Then, considering $\rd\bar D(\sP_j,\eps_j)$ as a Whitney-stratified subspace of $\RR^{|\sP_j|}$
with the lower strata being given by the intersections $\cap_\ell H_{P_\ell}(\sP_j,\eps_j)$, the function
$\vec t(\ft,\fs,\sP_j)$ is transverse to $\rd\bar D(\sP_j,\eps_j)$ in the sense of
\cite[Definition 1.3.1]{GorMacPh} and so $\vec t(\ft,\fs,\sP_j)^{-1}(\rd\bar D(\sP_j,\eps_j))$ is
a Whitney-stratified subspace of the space \eqref{eq:CompactDiskBundle}, whose top stratum
has a smooth collaring by \cite[Section 1.5]{GorMacPh}, proving Item \eqref{item:BoundaryOfNeigh2}.
\end{proof}

\subsection{The link of an ideal Seiberg--Witten moduli space}
\label{subsec:Link}
We can now give the

\begin{defn}
\label{defn:DefineLink}
The \emph{link of an ideal Seiberg--Witten moduli space}, $\bar\bL_{\ft,\fs}$, of
$M_{\fs}\times\Sym^\ell(X)$ in $\bar\sM_{\ft}/S^1$ is
$$
\bar\bL_{\ft,\fs}=\bar\bchi^{-1}(0)\cap \bar\bL^{\vir}_{\ft,\fs},
$$
where $\bar\bchi$ is the obstruction section in
Hypothesis \ref{hyp:Gluing} and $\bar\bL^{\vir}_{\ft,\fs}$ is the
virtual link defined in \eqref{eq:AmbientLink}.
\end{defn}

From Item \eqref{item:GluingHyp2} of Hypothesis \ref{hyp:Gluing},
the intersection of $\bar\bchi^{-1}(0)$ with each stratum of
$\bar\sM^{\vir}_{\ft,\fs}/S^1$ is a smooth submanifold of that stratum.
For generic values of $\eps$, the subspaces $H_P(\sP,\eps)$ from \eqref{eq:BoundaryHypersurface} defining
the
virtual link $\bar\bL^{\vir}_{\ft,\fs}$ will intersect $\bar\bchi^{-1}(0)$ transversally.
Hence, Lemma \ref{lem:BoundaryOfNeigh} yields:

\begin{lem}
\label{lem:DefiningLink}
For generic values of the
constants
$\delta,\eps_i$ used to define the virtual link in Lemma \ref{lem:BoundaryOfNeigh}, the following hold:
\begin{enumerate}
  \item
  \label{item:DefiningLink1}

  $\bar\bL_{\ft,\fs}$ is the boundary of a
 closed neighborhood of
 $M_{\fs}\times\Sym^\ell(X)$ in $\bar\sM_{\ft}/S^1$.

\item
\label{item:DefiningLink2}

    The intersection of $\bar\bL_{\ft,\fs}$ with each stratum
    $\sS$ of $\bar\sM_{\ft}/S^1$ is  a Whitney-stratified space
    whose top stratum is a codimension-one, collared submanifold of $\sS$.

\item
\label{item:DefiningLink3}

  The intersection of $\bar\bL^{\vir}_{\ft,\fs}$ with the top stratum of $\bar\sM_{\ft}/S^1$ is a topological manifold.
\end{enumerate}
\end{lem}

The proof of the following lemma then translates immediately from
that of \cite[Lemma 3.9]{FLLevelOne}.

\begin{lem}
\label{lem:GRIntersect0}
Assume $w\in H^2(X;\ZZ)$ is such that
$w\pmod{2}$ is \emph{good} in the sense of Definition \ref{defn:Good}.
Given a \spinu structure $\ft$ on $X$ and  a \spinc structure
$\fs$ on $X$ satisfying $\ell(\ft,\fs)\ge 0$ and $w_2(\ft)\equiv
w\pmod{2}$, there are positive constants $\eps_0$ and $\delta_0$
such that the following hold for all generic choices of positive constants
$\eps_i\leq \eps_0$ and $\delta\leq\delta_0$ defining
$\bar\bL^{\vir}_{\ft,\fs}$:
\begin{enumerate}
\item $\bar\bL_{\ft,\fs}$ is disjoint from $\bar M^w_\ka$
and $\bar\sM_{\ft}^{\red}$ in the stratification
\eqref{eq:StratificationCptPU(2)Space} of $\bar\sM_{\ft}/S^1$.

\item For all $z\in\AAA(X)$ and positive integers $\eta$
satisfying \eqref{eq:DimensionCondition}, and for the geometric representatives $\bar\sV(z)$ and
$\bar\sW$ discussed in Section
\ref{sec:Cohomology}, the intersection,
\begin{equation}
\label{eq:Int}
\bar\sV(z)\cap\bar\sW^{\eta}\cap\bar\bL_{\ft,\fs},
\end{equation}
is a finite collection of points contained in the top stratum
$\bL_{\ft,\fs}$ of $\bar\bL_{\ft,\fs}\subset\bar\sM_{\ft}/S^1$.
\end{enumerate}
\end{lem}

\subsection{A subspace of the virtual link of an ideal Seiberg--Witten moduli space}
The deformation retraction $N_{\ft(\ell),\fs}(\delta)\to M_{\fs}$ defines a
deformation retraction of $\bar\bL^{\vir,i}_{\ft,\fs}$ to the
subspace,
\begin{equation}
\label{eq:InstantonLinkComponentBase}
\bar{\bB\bL}^{\vir}_{\ft,\fs}
:=
\bigcup_{j=0}^r \ \bar{\bB\bL}^{\vir}_{\ft,\fs}(\sP_j),
\quad\text{where}\quad
\bar{\bB\bL}^{\vir}_{\ft,\fs}(\sP_j):=t_N^{-1}(0)\cap \bL^{\vir,i}_{\ft,\fs}(\sP_j),
\end{equation}
where the function $t_N$ is defined in \eqref{eq:NTubularDist}.
Essentially, $\bar{\bB\bL}^{\vir}_{\ft,\fs}(\sP_i)$ is defined by
replacing $\tilde N_{\ft(\ell),\fs}(\delta)$ with $\tilde M_{\fs}$ in the
definition \eqref{eq:DefineLinkStratum}.
In Section \ref{sec:ReduceToBase}, we will show that the intersections
with $\bar\bL_{\ft,\fs}$
in \eqref{eq:RawCobordismSum} can be replaced by intersections
with $\bar{\bB\bL}^{\vir}_{\ft,\fs}$.

\subsection{Orientations of the link of an ideal Seiberg--Witten moduli space}
\label{subsubsec:Orient}
To define an intersection number using the
intersection \eqref{eq:Int}, it is necessary to discuss
orientations of the link. An orientation for $\sM_{\ft}$
determines one for $\bL_{\ft,\fs}$ through the convention
introduced in \cite[Equations (2.16), (2.16) and (2.25)]{FL2b} by
considering $\bL_{\ft,\fs}$ as a boundary of
$(\sM_{\ft}\setminus\bga_{\sM}(\bar\sM^{\vir}_{\ft,\fs}))/S^1$. Specifically, at a
point $[A,\Phi]\in \bL_{\ft,\fs}$, if
\begin{itemize}
\item $\vec r\in T\sM_{\ft}^{*,0}$ is an outward-pointing radial
vector with respect to the open neighborhood
$\sM_{\ft}\cap\bga(\bar\sM^{\vir}_{\ft,\fs})$
 and spans a subspace
 complementary to the tangent space of
$\bL_{\ft,\fs}$,

\item $v_{S^1}\in T\sM_{\ft}^{*,0}$ is tangent to
the orbit of $[A,\Phi]$ under the (free) circle action (where
$S^1\subset\CC$ has its usual orientation), and

\item
$\la_{\sM}\in\det(T\sM_{\ft}^{*,0})$ is an orientation for
$T\sM_{\ft}$ at $[A,\Phi]$,
\end{itemize}
then we define an orientation $\la_{L}$ for
the tangent space  $T_{[A,\Phi]}\bL_{\ft,\fs}$ when $[A,\Phi]$ is a point in the top stratum of $\bar \sM_{\ft}$ by
\begin{equation}
\label{eq:QuotientBoundaryOrientation} \la_{\sM} =
-v_{S^1}\wedge\vec r \wedge \tilde{\la}_{L},
\end{equation}
where the lift $\tilde\la_{L}\in
\Lambda^{\max-2}(T\sM_{\ft}^{*,0})$ at
the point $[A,\Phi] \in \sM_{\ft}^{*,0}$
of
\[
\la_{L}\in\det(T\bL_{\ft,\fs})\subset
\Lambda^{\max-1}(T(\sM_{\ft}^{*,0}/S^1)),
\]
obeys
$\pi_*\tilde\la_{L}=\la_{L}$, where $\pi:\sM_{\ft}\to\sM_{\ft}/S^1$
denotes the quotient map.

\begin{defn}
\label{defn:BdryOrnLink} If $O$ is an orientation for $\sM_{\ft}$,
we call the orientation for $\bL_{\ft,\fs}$ related to $O$ by
equation \eqref{eq:QuotientBoundaryOrientation} the {\em boundary
orientation defined by $O$\/}.
\end{defn}

Given a homology orientation $\Om$, that is, an orientation for
$H^+(X;\RR)\oplus H^1(X;\RR)$ as defined prior to \cite[Definition 2.1]{KMStructure}
and an integral lift $w\in H^2(X;\ZZ)$ of $w_2(\ft)$, an orientation
$O^{\asd}(\Om,w)$ for $\sM_\ft$ is defined in \cite[Definition 2.3]{FL2b}.

The {\em standard orientation\/}
 for $\bar\bL_{\ft,\fs}$ is defined
in \cite[Definition 3.12]{FLLevelOne}.  The standard orientation
arises from orientations of spaces used in the definition of $\bar\bL^{\vir}_{\ft,\fs}$ and
thus is more natural for computations on this space. The standard
and boundary orientations are related in the following

\begin{lem}
\label{lem:OrientationFactor}
(See Feehan and Leness \cite[Lemmas 3.13 and 3.14]{FLLevelOne}.)
Let $\ft$ and $\ft(\ell)$ be \spinu structures on
$X$ satisfying
$$
p_1(\ft)=p_1(\ft(\ell))-4\ell,\quad c_1(\ft)=c_1(\ft(\ell)),
\quad\text{and}\quad w_2(\ft)=w_2(\ft(\ell)).
$$
Let $\Om$ be a homology orientation and let
$w \in H^2(X;\ZZ)$ be an integral
lift of $w_2(\ft) \in H^2(X;\ZZ/2\ZZ)$.  If $\ft(\ell)$ admits a splitting
$\ft(\ell)=\fs\oplus \fs\otimes L$, then the standard orientation
for $\bar\bL_{\ft,\fs}$ and the boundary orientation for
$\bar\bL_{\ft,\fs}$ defined through the orientation
$O^{\asd}(\Om,w)$ for $\sM_{\ft}$ differ by a factor of
\begin{equation}
\label{eq:OrientChangeFactor}
(-1)^{o_{\ft}(w,\fs)},
\quad\text{where}\quad o_{\ft}(w,\fs) :=
\frac{1}{4}
(w-c_1(L))^2.
\end{equation}
\end{lem}

We note some alternative expressions for the change
of orientation formula in \eqref{eq:OrientChangeFactor},
in particular one which  matches
that appearing in Conjecture \ref{conj:WC}.
Recall from \cite[Definition 1.2.8(c)]{GompfStipsicz} that $w\in H^2(X;\ZZ)$ is \emph{characteristic}
if $w\cdot \alpha\equiv \alpha^2\pmod 2$ for all $\alpha\in H^2(X;\ZZ)$.

\begin{lem}
\label{eq:SWSeriesOrientationFactor}
Continue the assumptions and notation of Lemma
\ref{lem:OrientationFactor}.  Then,
$$
o_\ft(w,\fs)
\equiv
\frac{1}{2}\left(w^2 + c_1(\fs)\cdot(w-c_1(\ft))\right)
+\frac{1}{2}\left(\sigma - w^2\right) \pmod 2.
$$
If $w$ is
characteristic and $c_1(\ft)\cdot c_1(\fs)=0$, then
$$
o_\ft(w,\fs)
\equiv
\frac{1}{2}\left(w^2 + c_1(\fs)\cdot w)\right)  \pmod 2.
$$
\end{lem}

\begin{proof}
The first equality
is given by
\cite[Equation 4.62]{FL2b}.  The
second equality follows immediately from the first and the additional
assumptions.
\end{proof}

\subsection{An equality of intersection numbers provided by the $\SO(3)$-monopole cobordism}
For all $z\in\AAA(X)$ and positive integers $\eta$ satisfying
\eqref{eq:DimensionCondition}, and for
the geometric representatives  $\bar\sV(z)$ and $\bar\sW$ discussed in Section
\ref{sec:Cohomology},
Lemma \ref{lem:GRIntersect0} implies that
we can define the intersection
number,
\begin{equation}
\label{eq:IntNumber}
\#\left(\bar\sV(z)\cap\bar\sW^{\eta}\cap\bar\bL_{\ft,\fs}\right),
\end{equation}
to be the oriented count of points in the intersection
\eqref{eq:Int}, using the standard orientation of
$\bar\bL_{\ft,\fs}$. This yields the following
equality of intersection numbers provided by the $\SO(3)$-monopole cobordism between links of the moduli space of anti-self-dual connections and links of the ideal moduli spaces of Seiberg--Witten moduli spaces.

\begin{thm}
\label{thm:CobordismThm}
Let $\ft$ be a \spinu structure on a
closed, oriented, smooth four-manifold $X$.  Let $z\in\AAA(X)$ and $\eta$
be a non-negative integer satisfying
\begin{equation}
\label{eq:DimensionCondition}
\deg(z)+2\eta=\dim\sM^{*,0}_{\ft}-2.
\end{equation}
Assume that there is a class $w\in H^2(X;\ZZ)$ satisfying $w_2(\ft)\equiv
w\pmod 2$ and which is \emph{good} in the sense of Definition
\ref{defn:Good}. Then the intersection numbers
\eqref{eq:IntNumber} obey
\begin{equation}
\label{eq:RawCobordismSum1}
\#\left(\bar\sV(z)\cap \bar\sW^{\eta-1}\cap \bar\bL^w_{\ft,\ka}\right)
=
-\sum_{\fs\in\Spinc(X)} (-1)^{o_{\ft}(w,\fs)}\ \# \left(
\bar\sV(z)\cap \bar\sW^{\eta-1}\cap \bar\bL_{\ft,\fs}\right),
\end{equation}
where $\bar\bL^w_{\ft,\ka}$ is the link of the
moduli space of
anti-self-dual connections in $\bar\sM_{\ft}/S^1$
specified in
\cite[Definition 3.7]{FL2a}.
\end{thm}

\section[Fiber bundle structure of the instanton component of the link]{Fiber bundle structure of the instanton component of the link of an ideal Seiberg--Witten moduli space}
\label{sec:FiberBundle}
We now describe the fiber bundle structure of the
subspaces $\bar\bL^{\vir,i}_{\ft,\fs}(\sP_j)$ (
defined in \eqref{eq:DefineLinkStratum}) of the
link of an ideal Seiberg--Witten moduli space.

Recall from Lemma \ref{lem:TubularDistanceFunctionOnOverlaps2} that
for non-negative integers
$j\le k\le r$,
the restriction of $\vec t(\ft,\fs,\sP_k)$ to
$\bar\sU(\ft,\fs,\sP_j)\cap\bar\sU(\ft,\fs,\sP_k)$ is given by the extension from the fibers of $\bar\sU(\ft,\fs,\sP_j)\to N_{\ft(\ell),\fs}(\delta)\times\Si(X^\ell,\sP_j)$ of the function
$$
\vec t_f(\sP_j,[\sP_k]):
\bar M(\sP_j)
\to
\bigsqcup_{\sP''\in [\sP_j<\sP_k]} I^{\sP''}.
$$
(For $j=k$, recall that $[\sP_j<\sP_k]=\{\sP_j\}$.)
We note that the subset
$$
\bigsqcup_{\sP''\in [\sP_j<\sP_k]} D(\sP'',\eps_k)
\subset
\bigsqcup_{\sP''\in [\sP_j<\sP_k]} I^{\sP''}
$$
is
closed under the action of $\fS(\sP)$.
For the constants
$\beps=(\eps_0,\dots,\eps_r)$  used in the definition
of the link in Lemma \ref{lem:DefiningLink}, define
\begin{equation}
\label{eq:DefineLinkFiber}
\begin{aligned}
\bar M(\sP_j,\beps)
&:=
\bar M(\sP_j)\cap \vec t_f(\sP_j,[\sP_j])^{-1}(\rd \bar D(\sP_j,\eps_j))
\\
&\qquad
\left.\setminus
\bigcup_{k>j}\  \vec t_f(\sP_j,[\sP_k])^{-1}
\left(\bigsqcup_{\sP''\in [\sP_j<\sP_k]} D(\sP'',\eps_k)\right)\right. .
\end{aligned}
\end{equation}
We then have the following fiber bundle structure on
the subspace $\bL^{\vir,i}_{\ft,\fs}(\sP_j)$ of the link.

\begin{prop}
\label{prop:LinkPieceFiberBundleStructure}
Let $\beps=(\eps_0,\dots,\eps_r)$ be the generic constants
defining the link
in Lemma \ref{lem:DefiningLink}.
Let $K_j\Subset\Si(X^\ell,\sP_j)$ be the compact
subset defined by those
constants as in Lemma \ref{lem:CompactSubsetsOfSi}.
For $0\le j\le r$,
the subspace $\bL^{\vir,i}_{\ft,\fs}(\sP_j)$ of
$\bar\sU(\ft,\fs,\sP_j)/S^1$ defined in
\eqref{eq:DefineLinkStratum}
admits a fiber bundle structure,
\begin{equation}
\label{LocalInstantonLink0}
\begin{CD}
\barM(\sP_j,\beps) @>>> \tilde N_{\ft(\ell),\fs}(\delta)\times_{\sG_{\fs}\times
S^1}
\Fr(\ft,\fs,\sP_j)|_{K_j}\times_{G(\sP_j)}\barM(\sP_j,\beps)
\\
@. @V \pi(\ft,\fs,\sP_j) VV
\\
@. N_{\ft(\ell),\fs}(\delta)/S^1\times K_j
\end{CD}
\end{equation}
where $\barM(\sP_j,\beps)$  is defined in \eqref{eq:DefineLinkFiber}.
\end{prop}

\begin{proof}
The proposition follows from our previously established
results on the intersection,
$$
\left(\vec t(\ft,\fs,\sP_k)^{-1}(D(\sP_k,\eps_k))\right)
\cap
\left( \vec t(\ft,\fs,\sP_j)^{-1}(\rd \bar D(\sP_j,\eps_j))\right),
$$
as we now describe.  For $j<k$, Lemma \ref{lem:TubularDistanceFunctionOnOverlaps2}
equates $\vec t(\ft,\fs,\sP_k)$ and the extension of $\vec t_f(\sP_j,[\sP_k])$
and thus shows that the preceding intersection is given by the definition in \eqref{eq:DefineLinkFiber}.
For $k<j$, Corollary \ref{cor:LowerStratumRemoval}
and the construction of $K_j$ in
Lemma \ref{lem:CompactSubsetsOfSi}
shows that removing the subset
$$
\bigcup_{j<k}\ \vec t(\ft,\fs,\sP_k)^{-1}(D(\sP_k,\eps_k))
$$
from $\vec t(\ft,\fs,\sP_j)^{-1}(\rd \bar D(\sP_j,\eps_j))$ is equivalent
to restricting the fiber bundle \eqref{LocalInstantonLink0} to $K_j$.
\end{proof}

The following lemma will be used to define quotients of the link $\bL^{\vir}_{\ft,\fs}$ in
Chapter \ref{chap:Comp}.

\begin{lem}
\label{lem:S1ActionFreeOnFibers}
The restriction of the $\SO(3)$ action on $\bar M(\sP)$
defined by the diagonal action on the frames,
\begin{equation}
\label{eq:DiagonalFrameAction}
(([A_P,F^s_P,\bx_P])_{P\in \sP},M)\mapsto
(([A_P,F^s_PM,\bx_P])_{P\in \sP}),
\end{equation}
to $\bar M(\sP,\beps)$
(defined in \eqref{eq:DefineLinkFiber})
is free, where $M\in\SO(3)$ and $[A_P,F^s_P,\bx_P]\in \barM^{s,\natural}_{\spl,|P|}(\delta)$.
\end{lem}

\begin{proof}
The $\SO(3)$ action on $\barM^{s,\natural}_{\spl,|P|}(\delta)$ given by
the action on the frame is free at all points except for points
of the form $[\Theta,F^s,v_P]$, where $\Theta$ is the product connection.
Thus a fixed point of the action \eqref{eq:DiagonalFrameAction} will
be of the form, $A_{\sP}=(([\Theta,F^s_P,\bx_P])_{P\in \sP})$.
If $A_{\sP}\in \bar M(\sP,\beps)$, then
the requirement that $\vec t_f(\sP,[\sP])(A_\sP)\in \rd \bar D(\sP,\eps)$ in
\eqref{eq:DefineLinkFiber}
implies that not all the points $\bx_P$ can be cone points, $c_P$.
Hence, for such a point $A_{\sP}\in \bar M(\sP,\beps)$ and for any $F_{\sP}\in\Fr(\ft,\fs,\sP)$,
the point $[F_{\sP},A_{\sP}]$ will lie in the subspace
$T(\ft,\fs,\sP,[\sP'])$ defined in \eqref{eq:DefineTrivialInOverlap}
for some partition $\sP'$ with $\sP<\sP'$.  By a suitable choice
of $F_{\sP}$ (that is, one which lies over a point $\bx\in\Si(X^\ell,\sP)$ sufficiently
far from the lower strata in $\Sym^\ell(X)$), we can ensure that
$[F_{\sP},A_{\sP}]\in T^{\sO}(\ft,\fs,\sP,[\sP'])$, where
$T^{\sO}(\ft,\fs,\sP,[\sP'])$ is defined in Lemma \ref{lem:SplicingTrivialOverlap}.
By Lemma \ref{lem:SplicingTrivialOverlap},
Theorem \ref{thm:GlobalSplicingDataOverlaps},
and the definition of the equivalence relation $\sim$ in
\eqref{eq:DefineGlobalGluingDataSpace},
there is a point $[F_{\sP'},A_{\sP'}]\in\Si(\ft,\fs,\sP')$,
where $\Si(\ft,\fs,\sP')\subset\bar\Gl(\ft,\fs,\sP')$ is defined in \eqref{eq:DefineTrivialStratumXConePoints},
 with
$$
[(A_0,\Phi_0),F_{\sP'},A_{\sP'}]
\sim
[(A_0,\Phi_0),F_{\sP},A_{\sP}],
$$
for any
$[A_0,\Phi_0]\in N_{\ft(\ell),\fs}(\delta)$.
Because $[F_{\sP'},A_{\sP'}]\in\Si(\ft,\fs,\sP')$, we have
\[
\vec t(\ft,\fs,\sP')([(A_0,\Phi_0),F_{\sP},A_{\sP}])=\vec 0,
\]
where $\vec 0$ is the
vector with all coordinates equal to zero.
However, Lemma
\ref{lem:TubularDistanceFunctionOnOverlaps2} then implies that
$\vec t_f(\sP,[\sP'])(A_{\sP})=\vec 0$, contradicting
the assumption that $A_{\sP}\in \bar M(\sP,\beps)$.
\end{proof}

As we describe in the following lemma,
the subspaces $\bar{\bB\bL}^{\vir}_{\ft,\fs}(\sP_j)$
admit a fiber bundle structure similar to that provided by
Proposition \ref{prop:LinkPieceFiberBundleStructure}.

\begin{lem}
\label{lem:LocalBaseLinkFiberBundle}
Continue the hypotheses of Proposition \ref{prop:LinkPieceFiberBundleStructure}.
Then the space
\begin{equation}
\label{eq:DefineLocalBaseLink}
\bar{\bB\bL}^{\vir}_{\ft,\fs}(\sP_j)
= t_N^{-1}(0)\cap \bar\bL^{\vir,i}_{\ft,\fs}(\sP_j)
\end{equation}
is the restriction of the fiber bundle \eqref{LocalInstantonLink0}
to the subspace
$$
M_{\fs}\times K_j \subset N_{\ft(\ell),\fs}(\delta)/S^1\times K_j
$$
of the base and is thus given by
\begin{equation}
\label{eq:LocalBaseLinkFibration}
\begin{CD}
\barM(\sP_j,\beps) @>>> \tilde M_{\fs}\times_{\sG_{\fs}\times S^1} \Fr(\ft,\fs,\sP_j,g_{\sP_j})\times_{G(\sP_j)}\barM(\sP,\beps)
\\
@. @V \pi(\ft,\fs,\sP_j) VV
\\
@. M_{\fs}\times K_j
\end{CD}
\end{equation}
\end{lem}

\begin{proof}
The conclusion of the lemma follows immediately from Proposition
\ref{prop:LinkPieceFiberBundleStructure}.
\end{proof}

\section{Boundaries of components of links of ideal Seiberg--Witten moduli spaces}
\label{sec:Boundaries}
In Chapter \ref{chap:Cohom} we will show that the intersection number \eqref{eq:IntNumber}
can be written as a cohomological pairing
and in Lemma \ref{lem:ReduceToBase} we will show that this pairing can be written as a pairing
with a fundamental class defined by the subspace $\bar{\bB\bL}^{\vir}_{\ft,\fs}$.
To apply the technique for computing such pairings  described in Section \ref{subsubsec:IntroComp},
we will define quotients of the intersections of the subspaces $\bar{\bB\bL}^{\vir}_{\ft,\fs}(\sP_j)$
in Chapter \ref{chap:Comp}.  To ensure these quotients satisfy the conditions described in Section \ref{subsubsec:IntroComp}, we now prove that these intersections
are subbundles of the fiber bundles \eqref{eq:LocalBaseLinkFibration}.

Define
\begin{equation}
\label{eq:DefineBoundaryOne}
\rd_{k_1}\cdots\rd_{k_p}\bar{\bB\bL}^{\vir}_{\ft,\fs}(\sP_j)
:=
 \bar{\bB\bL}^{\vir}_{\ft,\fs}(\sP_j)
 \cap
 \bigcap_{u=1}^p \vec t(\ft,\fs,\sP_{k_u})^{-1}(\rd\bar D(\sP_{k_u},\eps_{k_u})).
\end{equation}
From the definitions \eqref{eq:DefineLocalBaseLink} and \eqref{eq:DefineLinkStratum},
we have the equality,
\begin{equation}
\label{eq:MatchingBoundaries}
\rd_k\bar{\bB\bL}^{\vir}_{\ft,\fs}(\sP_j)
=
\bar{\bB\bL}^{\vir}_{\ft,\fs}(\sP_j) \cap
\bar{\bB\bL}^{\vir}_{\ft,\fs}(\sP_k)
=
\rd_j\bar{\bB\bL}^{\vir}_{\ft,\fs}(\sP_k).
\end{equation}
More generally,
$$
\bigcap_{j=0}^p \bar{\bB\bL}^{\vir}_{\ft,\fs}(\sP_{i_j})
=
\rd_{i_1}\cdots\rd_{i_p}\bar{\bB\bL}^{\vir}_{\ft,\fs}(\sP_{i_0}).
$$
If
$$
\cl\left( \Si(X^\ell,\sP_j)\right)\cap \Si(X^\ell,\sP_k) =\emptyset
\quad\text{and}\quad \Si(X^\ell,\sP_j)\cap\cl\left(\Si(X^\ell,\sP_k)\right)
=\emptyset,
$$
then the boundaries in \eqref{eq:MatchingBoundaries} are empty.

To describe the intersection \eqref{eq:DefineBoundaryOne}, we
introduce the following notation.
For $j<i_1<i_2<\dots <i_v$, define
\begin{equation}
\label{eq:DefineLinkFiberBoundary}
\rd_{i_1}\rd_{i_2}\cdots\rd_{i_v}
\bar M(\sP_j,\beps)
:=
\bar M(\sP_j,\beps)
\cap
\bigcap_{u=1}^v \vec t_f(\sP_j,[\sP_{i_u}])^{-1} (\rd\bar D(\sP_{i_u},\eps_{i_u})).
\end{equation}
For
non-negative integers
$k_1< k_2<\dots< k_p<j \leq r$, we define
\begin{equation}
\label{eq:DefineBoundaryOfCompactum}
\rd_{k_1}\cdots\rd_{k_p}K_j
:=
K_j\cap
\bigcap_{u=1}^p \vec t(X^\ell,g_{\sP_{k_u}})^{-1}(\rd\bar D(\sP_{k_u},\eps_{k_u}))
\subset \Si(X^\ell,\sP_j),
\end{equation}
where the function $\vec t(X^\ell,g_{\sP_{k_u}})$ is defined in
\eqref{eq:DefineDiagonalTubularDistFunction}.

\begin{prop}
\label{prop:LinkCornerdescription}
For
non-negative integers
$k_1< k_2<\dots< k_p<j$ and $j<i_1<i_2<\dots <i_v\le r$,
the
intersection,
$$
\rd_{k_1}\cdots\rd_{k_p}\rd_{i_1}\cdots\rd_{i_v}\bar{\bB\bL}^{\vir}_{\ft,\fs}(\sP_j),
$$
defined in \eqref{eq:DefineBoundaryOne}
admits a description as a fiber bundle,
\begin{equation}
\label{eq:LinkCornerFiberBundle}
\begin{CD}
\rd_{i_1}\rd_{i_2}\cdots\rd_{i_v}\bar M(\sP_j,\beps)
@>>>
\rd_{k_1}\cdots\rd_{k_p}\rd_{i_1}\cdots\rd_{i_v}\bar{\bB\bL}^{\vir}_{\ft,\fs}(\sP_j)
\\
@. @V \pi(\ft,\fs,\sP_j) VV
\\
@. M_{\fs}\times \rd_{k_1}\cdots\rd_{k_p}K_j
\end{CD}
\end{equation}
arising from the
equality,
\begin{align*}
{}&
\rd_{k_1}\cdots\rd_{k_p}\rd_{i_1}\cdots\rd_{i_v}\bar{\bB\bL}^{\vir}_{\ft,\fs}(\sP_j)
\\
{}&\quad
=
\tilde M_{\fs}\times_{\sG_{\fs}\times S^1}
\Fr(\ft,\fs,\sP_j)|_{\rd_{k_1}\cdots\rd_{k_p}K_j}\times_{G(\sP_j)}
\rd_{i_1}\rd_{i_2}\cdots\rd_{i_v}\bar M(\sP_j,\beps).
\end{align*}
\end{prop}

\begin{proof}
For
non-negative integers
$j<i_u \leq r$,
Lemma \ref{lem:TubularDistanceFunctionOnOverlaps2} implies that
the intersection of $\bar{\bB\bL}^{\vir}_{\ft,\fs}(\sP_j)$
with $\vec t(\ft,\fs,\sP_{i_u})^{-1} (\rd\bar D(\sP_{i_u},\eps_{i_u}))$
is given by the intersection of the fiber $\bar M(\sP_j,\beps)$
with
\[
\vec t_f(\sP_j,[\sP_{i_u}])^{-1} (\rd\bar D(\sP_{i_u},\eps_{i_u})).
\]
For
non-negative integers
$k<j \leq r$, we claim that the intersection of
$\bar{\bB\bL}^{\vir}_{\ft,\fs}(\sP_j)$
with
\[
\vec t(\ft,\fs,\sP_k)^{-1} (\rd\bar D(\sP_k,\eps_k))
\]
is
given by the inverse image, under the projection
$\pi(\ft,\fs,\sP_j)$, of the intersection of $M_{\fs}\times K_j$
with
\[
M_{\fs}\times \vec t(X^\ell,g_{\sP_k})^{-1}(\rd\bar D(\sP_k,\eps_k)).
\]
To prove this claim, we will show that for $\bA\in \bar{\bB\bL}^{\vir}_{\ft,\fs}(\sP_j)$
and $\bA' := \pi(\ft,\fs,\sP_j)(\bA)$ we have
\begin{equation}
\label{eq:RelationOnTub}
\vec t(\ft,\fs,\sP_k)(\bA)\in \rd\bar D(\sP_k,\eps_k)
\quad\text{if and only if}\quad
\vec t(\ft,\fs,\sP_k)(\bA')\in \rd\bar D(\sP_k,\eps_k).
\end{equation}
The latter condition is true if and only if $\bA'\in M_{\fs}\times \rd_k K_j$
by Lemma \ref{lem:TubularDistanceFunctionOnOverlaps3}
and the definition of $\rd_k K_j$.
Thus it suffices to prove \eqref{eq:RelationOnTub} in order to complete the proof of Proposition \ref{prop:LinkCornerdescription}.
To this end, we first collect the following observations.

Items \eqref{item:SymProjTubDistRelations01} and
\eqref{item:SymProjTubDistRelations02} of Lemma \ref{lem:SymProjTubDistRelations0} imply that
for all $Q\in\sP_k$ we have
\begin{equation}
\label{eq:TubVal1a}
\vec t^Q(\ft,\fs,\sP_k)(\bA')\le \vec t^Q(\ft,\fs,\sP_k)(\bA)
\end{equation}
and
\begin{equation}
\label{eq:TubVal1b}
\vec t^Q(\ft,\fs,\sP_k)(\bA)\neq \vec t^Q(\ft,\fs,\sP_k)(\bA')
\quad
\text{if and only if $Q\in\sP_j\cap \sP_k$ }
\end{equation}
and in this case ($Q\in\sP_j\cap \sP_k$), we see that $\vec t^Q(\ft,\fs,\sP_k)(\bA')=0$.
We further note that because
$\bA \in \bar{\bB\bL}^{\vir}_{\ft,\fs}(\sP_j)$, we have
\begin{equation}
\label{eq:TubVal1}
\vec t(\ft,\fs,\sP_j)(\bA) \in \bar D(\sP_j,\eps_j)
\quad\text{and thus}\quad
\vec t^Q(\ft,\fs,\sP_j)(\bA)\le \eps_j, \quad\text{for all $Q\in \sP_j$}.
\end{equation}
Item \eqref{item:SymProjTubDistRelations03} in Lemma \ref{lem:SymProjTubDistRelations0}
and \eqref{eq:TubVal1}
imply that
\begin{equation}
\label{eq:TubVal1c}
\vec t^Q(\ft,\fs,\sP_k)(\bA)=\vec t^Q(\ft,\fs,\sP_j)(\bA)\le \eps_j<\eps_k
\quad
\text{for $Q\in\sP_j\cap\sP_k$}.
\end{equation}
If $\vec t(\ft,\fs,\sP_k)(\bA)\in \rd\bar D(\sP_k,\eps_k)$,
then $\vec t^Q(\ft,\fs,\sP_k)(\bA)\le\eps_k$ for all $Q\in\sP_k$
and there is at least one $Q_0\in\sP_k$ with $\vec t^{Q_0}(\ft,\fs,\sP_k)(\bA)=\eps_k$.
By \eqref{eq:TubVal1c}, we have $Q_0\notin\sP_k\cap\sP_j$ so
\eqref{eq:TubVal1b} implies that $\vec t^{Q_0}(\ft,\fs,\sP_k)(\bA')=\eps_k$.
This equality and \eqref{eq:TubVal1a} then imply that
$\vec t(\ft,\fs,\sP_k)(\bA')\in \rd\bar D(\sP_k,\eps_k)$.

If $\vec t(\ft,\fs,\sP_k)(\bA')\in \rd\bar D(\sP_k,\eps_k)$,
then $\vec t^Q(\ft,\fs,\sP_k)(\bA')\le\eps_k$ for all $Q\in\sP_k$
and there is at least one $Q_0\in\sP_k$ with $\vec t^{Q_0}(\ft,\fs,\sP_k)(\bA')=\eps_k$.
Because $\eps_k\neq 0$, then \eqref{eq:TubVal1} implies that  $Q_0\notin\sP_k\cap\sP_j$ and
$\vec t^{Q_0}(\ft,\fs,\sP_k)(\bA)=\eps_k$.
For all other $Q\in\sP_k$, if $Q\notin \sP_k\cap\sP_j$, then
\eqref{eq:TubVal1b} implies that
$$
\vec t^Q(\ft,\fs,\sP_k)(\bA)
=
\vec t^Q(\ft,\fs,\sP_k)(\bA')\le\eps_k,
$$
while for $Q\in\sP_k\cap\sP_j$, we see that \eqref{eq:TubVal1c} implies
$\vec t^Q(\ft,\fs,\sP_k)(\bA)<\eps_k$.  Hence,
$\vec t(\ft,\fs,\sP_k)(\bA)\in \rd\bar D(\sP_k,\eps_k)$,
completing the proof of \eqref{eq:RelationOnTub}. This completes the proof of Proposition \ref{prop:LinkCornerdescription}.
\end{proof}

Because $U_f(\sP_i,[\sP_j])$ is defined
in \eqref{eq:DefineOverlapFiberAsImage}
as the image of the injective (up to symmetric group action)
splicing map $\rho^{\ft,\fs,d}_{f,\sP,[\sP']}$, there is a
homeomorphism between
$U_f(\sP_i,[\sP_j])$ and an open subspace of
$$
\left.
\bigsqcup_{\sP''\in [\sP_i<\sP_j]} \prod_{P\in\sP_i} \left(
\Delta^\circ(Z_{P}(\delta_P),\sP''_P) \times\prod_{Q\in\sP''_P}
\barM^{s,\natural}_{\spl,|Q|}(\delta_Q) \right) \right/\fS(\sP_i)  .
$$
Composing this homeomorphism with the obvious projection map defines a surjective map,
\begin{equation}
\label{eq:FiberBoundaryProjectionMap} U_f(\sP_i,[\sP_j]) \to
\bigsqcup_{\sP''\in [\sP_i<\sP_j]} \prod_{P\in\sP_i}
\Delta^\circ(Z_{P}(\delta_P),\sP''_P)/\fS(\sP_i).
\end{equation}
In addition, because
\begin{align*}
\left.
\left(
\bigsqcup_{\sP''\in [\sP_i<\sP_j]} \prod_{P\in\sP_i}
\prod_{Q\in\sP''_P} \barM^{s,\natural}_{\spl,|Q|}(\delta_Q)
\right) \right/\fS(\sP_i)
&=
\left.
\left(
\bigsqcup_{\sP''\in [\sP_i<\sP_j]} \prod_{Q\in\sP''}
\barM^{s,\natural}_{\spl,|Q|}(\delta_Q)
\right)\right/\fS(\sP_i)
\\
&=
\left.
\left(
\bigsqcup_{\sP''\in [\sP_i<\sP_j]}\barM(\sP'')
\right)\right/\fS(\sP_i),
\end{align*}
where $\barM(\sP'')$ is defined in \eqref{eq:DefineGluingDataFiber}, there is
a $G(\sP_i)$-equivariant  map,
\begin{equation}
\label{eq:FiberBoundaryProjectionMapToConn}
c_{j,i}:U_f(\sP_i,[\sP_j]) \to
\left.\left(
\bigsqcup_{\sP''\in [\sP_i<\sP_j]} \barM(\sP'')\right)\right/\fS(\sP_i),
\end{equation}
which will be used in Chapter \ref{chap:Comp}.  In that
chapter,
we will also require the characterization of
a group action on the image of $c_{j,i}$ as described
in the following two lemmas.

\begin{lem}
\label{eq:PseudoFreeActions}
Let $S^1$ act  on a
topological
space
$\sX$.
Assume that there is
a non-zero
$N\in\NN$ such that the order of every isotropy group of this action divides $N$.
If there is an action of a finite group $G$ on
$\sX$ which commutes with the action of $S^1$, then
the $S^1$ action on $\sX$ descends to an $S^1$ action
on $\sX/G$ such that if $\la[x]=[x]$
for $\la\in S^1$, $x\in \sX$, and $[x]\in \sX/G$, then $\la$
is a $|G|N$-th root of unity.
\end{lem}

\begin{proof}
Denote
the $S^1$ action on $\sX/G$ by $(\la,[x])\mapsto [\la x]$.
Because the $S^1$ and $G$ actions on $\sX$ commute, this
action is well-defined.
If $\la [x]=[x]$, then there is a group element, $g\in G$, such that
$\la x=gx$.  Then $\la^{|G|}x=g^{|G|}x=ex=x$
so $\la^{|G|}$ is in the isotropy group of $x$.  By the assumption on the order of
the isotropy groups, $\la^{|G|N}=e$, as asserted.
\end{proof}

We say that an action of a group $H$ on a
topological
space
$\sX$ has {\em finite isotropy\/} if there
is a finite subgroup $K<H$ such that for all
$x\in \sX$, we have $H_x<K$,
where $H_x \subset H$ denotes the isotropy subgroup of $x \in \sX$.

\begin{lem}
\label{lem:FreeActionOnCijImage}
The restriction of the diagonal $\SO(3)$ action on the frames
in \eqref{eq:DiagonalFrameAction} to
$$
c_{j,i}\left(U_f(\sP_i,[\sP_j])\cap \bar M(\sP_i,\eps) \right)
\subset
\left.\left(\bigsqcup_{\sP''\in [\sP_i<\sP_j]} \barM(\sP'')\right)\right/\fS(\sP_i)
$$
has finite isotropy.
\end{lem}

\begin{proof}
For the product connection, $\Theta$, on $S^4\times\so(3)$, we
define
$$
T(\sP):=\{(([A_P,F^s_P,\bx_P])_{P\in\sP})\in\bar M(\sP): A_P=\Theta\
\text{for all $P\in\sP$}\}.
$$
We first claim that the image of the restriction of
$c_{j,i}$ to $U_f(\sP_i,[\sP_j])\cap \bar M(\sP_i,\eps)$
is disjoint from $T(\sP'')$ for all $\sP''\in [\sP_i<\sP_j]$.
The proof of Lemma \ref{lem:S1ActionFreeOnFibers}
implies that $\bar M(\sP_i,\eps)\cap T(\sP_i)$ is empty.
If $A(\sP_i)\in U_f(\sP_i,[\sP_j])\cap \bar M(\sP_i,\eps)$
and $c_{j,i}(A(\sP_i))\in T(\sP'')$, then the
definition
of $c_{j,i}$ in \eqref{eq:FiberBoundaryProjectionMapToConn} and the
definition of the map
$\rho^{\ft,\fs,d}_{f,\sP_i,[\sP_j]}$ in \eqref{eq:DownwardsInclusionFiberMap}
would imply that
$A(\sP_i)\in T(\sP_i)$, a contradiction.

The proof of Lemma \ref{lem:S1ActionFreeOnFibers} implies that
the action \eqref{eq:DiagonalFrameAction}
is free on $\bar M(\sP'')\setminus T(\sP'')$.  The action on the
$\fS(\sP_i)$ quotient then  has finite isotropy by
Lemma \ref{eq:PseudoFreeActions}.
This completes the proof of Lemma \ref{lem:FreeActionOnCijImage}.
\end{proof}

\chapter{Cohomology and duality}
\label{chap:Cohom}
In this chapter, we will prove that the intersection number,
$$
\#\left(\bar\sV(z)\cap\bar\sW^{\eta} \cap
\bar\bL_{\ft,\fs}\right),
$$
appearing in \eqref{eq:RawCobordismSum},
can be expressed as the pairing of a cohomology class with a
homology class.  Our aim in doing this is to rewrite the above
intersection number in a form to which we can apply the pushforward-pullback argument
described in Section \ref{subsubsec:IntroComp}.

\section{Introduction}
\label{sec:Cohom_intro}
Explicitly, in this chapter we will prove the following analogue
of \cite[Proposition 5.2]{FLLevelOne}.

\begin{prop}
\label{prop:Duality}
For $\beta\in H_\bullet(X;\RR)$, let
$\barmu_p(\beta), \barmu_c \in H^\bullet(\bar\sM^{\vir,*}_{\ft,\fs}/S^1;\RR)$
be the cohomology classes in Definition \ref{defn:ExtendedCohomologyClasses}.
Let $\bar e_s=e(\bar\Upsilon^s_{\ft,\fs}/S^1)$ be the Euler class of the
background obstruction bundle in \eqref{eq:DefineGlobalBackgroundObstruction}
and let $\bar e_I$ be the extension of the Euler class of the instanton obstruction
bundle in Definition \ref{defn:ExtendedCohomologyClasses}. If
$d(\ft)=\dim\sM^{*,0}_{\ft,\fs}$
and $\deg(z)+2\eta=d(\ft)-2$,
then let
$[\bar{\bL}^{\vir}_{\ft,\fs}]\in
H_{d(\ft)-2}(\bar\sM^{\vir,*}_{\ft,\fs}/S^1;\RR)$ be the homology
class defined in \eqref{eq:DefineAmbientLinkFund}. Then
\begin{equation}
\label{eq:Duality}
\#\left(\bar\sV(z)\cap\bar\sW^{\eta} \cap
\bar{\bL}_{\ft,\fs}\right) = \langle \barmu_p(z)\smile
\barmu_c^{\eta}\smile \bar e_I \smile\bar e_s,
[\bar\bL^{\vir}_{\ft,\fs}] \rangle.
\end{equation}
\end{prop}

Let $\hat\bL$ be any codimension-zero, compact submanifold of
the top stratum of $\bL^{\vir}_{\ft,\fs}$ with the property
that its boundary lies in a neighborhood $U$ of the
lower strata, where $U$ retracts to the lower strata.
The existence of such a neighborhood, $U$, follows
immediately from Lemma \ref{lem:StatifiedSpaceStr}.
The fundamental class of $\hat\bL$ is a relative
class and we can define the fundamental class of
$\bL^{\vir}_{\ft,\fs}$ through the exact sequence
of the pair $(\bL^{\vir}_{\ft,\fs},\bL^{\vir}_{\ft,\fs}\cap U)$.

A geometric representative $\sV$ on a stratified
space $\sM$, as described in Definition \ref{defn:GeomRepresentative},
naturally defines a relative cohomology class in
$H^\bullet(\sM,\sM\setminus\sV;\RR)$.
If $\sV$ were a smooth submanifold of a manifold $\sM$, then the
relative cohomology class would be given by the Thom class of the
normal bundle of $\sV$.
The geometric representatives $\bar\sV(z)$ and $\bar\sW$ define such
cohomology classes.  In addition, the restriction of the zero locus of the
obstruction section to the top stratum also defines such a relative
cohomology class, namely the relative Euler class.  We note that
these classes are, initially, defined only on the top stratum of
$\bar\sM^{\vir,*}_{\ft,\fs}/S^1$.  The intersection number
in \eqref{eq:Duality} equals the pairing of the product of these
relative cohomology classes with the relative fundamental class of
the manifold with boundary $\hat\bL$.

To rewrite this pairing of relative  classes in a form
to which we can apply the pushforward-pullback argument
described in Section \ref{subsubsec:IntroComp}, we
relate the relative cohomology classes mentioned above
with absolute cohomology classes defined on $\bar\sM^{\vir,*}_{\ft,\fs}/S^1$
and compute these cohomology classes in terms of cohomology classes on
$M_{\fs}\times\Sym^\ell(X)$ and the Chern class defined by the
$S^1$ action.

The relative cohomology classes defined by $\bar\sV(z)$, $\bar\sW$,
and the obstruction section extend to relative cohomology
classes on the complement of a small subspace of $\bar\sM^{\vir,*}_{\ft,\fs}/S^1$.
If an embedded surface $\Si\subset X$ satisfies
$q[\Si]=h\in H_2(X;\RR)$ for $q\in\RR$,
then
the
relative cohomology class for the geometric representative $\bar\sV(h)$
defines such a relative cohomology class in
$$
H^\bullet\left(\bar\sM^{\vir,*}_{\ft,\fs}/S^1\setminus\sI(\nu(\Si)),
\bar\sM^{\vir,*}_{\ft,\fs}/S^1\setminus \left(\sI(\nu(\Si))\cup\bar\sV(h)\right);\RR\right),
$$
where $\sI(\nu(\Si))$ is the subspace of points $[A,\Phi,\bx]$ where
the support of $\bx$ intersects $\nu(\Si)$, a tubular neighborhood
of $\Si$. The relative cohomology
class for $\bar\sW$ extends to a relative cohomology class on the complement
of $\sI(\nu(x))$.  The relative Euler class of the obstruction section
extends to the complement of the intersection of the zero locus of the
obstruction section with the lower strata.
We then compute that the image of these extended relative cohomology
classes, in the exact sequence of the relevant pair, equals the
restriction of the cohomology classes of the desired form from
$\bar\sM^{\vir,*}_{\ft,\fs}/S^1$ to the subspaces described above.
With these equalities established, the proof of Proposition
\ref{prop:Duality} is then largely a formal manipulation.

After defining subspaces of and cohomology classes on
$\bar\sM^{\vir,*}_{\ft,\fs}/S^1$ in Section \ref{sec:CohomDefn},
we define the relevant fundamental class in Section \ref{sec:FundClassOfAmbLink}.
In Section \ref{sec:ComputingMuClasses}, we define the relative cohomology
classes corresponding to the geometric representatives $\bar\sV(z)$
and $\bar\sW$, define their extensions, and compute the corresponding
absolute cohomology classes in terms of the cohomology classes
defined in Section \ref{sec:CohomDefn}.
We discuss relative Euler classes and carry out a similar program
for a geometric dual of $\bar\bchi^{-1}(0)$ in Section \ref{sec:RelEulerClassObstr}.
The proof of Proposition \ref{prop:Duality} appears
in Section \ref{sec:Duality}.  Finally, in Section \ref{sec:ReduceToBase},
we show how to replace the pairing
in \eqref{eq:Duality} with a pairing with the fundamental class
of $\bar{\bB\bL}^{\vir}_{\ft,\fs}$, the space defined in
\eqref{eq:InstantonLinkComponentBase}.

\section{Definitions}
\label{sec:CohomDefn}

\subsection{Subspaces and maps}
We begin by defining subspaces of $\bar\sM^{\vir}_{\ft,\fs}$.
First, let
\begin{equation}
\label{eq:TopStratumOfVirtualInclusion}
\iota:\sM^{\vir}_{\ft,\fs}\to\bar\sM^{\vir,*}_{\ft,\fs}
\end{equation}
be the inclusion of the top stratum,
which is mapped to $\sC_{\ft}$ by
the splicing and gluing maps $\bga'_{\sM}$ and $\bga_{\sM}$,
and let
\begin{equation}
\label{eq:RedComplementInVirtual}
\bar\sM^{\vir,*}_{\ft,\fs}
=
\bar\sM^{\vir}_{\ft,\fs} \setminus \left(M_{\fs}\times\Sym^\ell(X)\right).
\end{equation}
The $S^1$ action \eqref{eq:DefineGlobalS1Action} is free on the subspace
$\bar\sM^{\vir,*}_{\ft,\fs}$.
The map
$$
\pi_N:\bar\sM^{\vir}_{\ft,\fs}\to N_{\ft(\ell),\fs}(\delta)
$$
was defined in \eqref{eq:GlobalProjectionToN}
and we define
\begin{equation}
\label{eq:ProjectionToSW}
\pi_{\fs}:\bar\sM^{\vir}_{\ft,\fs}\to M_{\fs}
\end{equation}
as the composition of $\pi_N$ with the projection
$N_{\ft(\ell),\fs}(\delta)\to M_{\fs}$.
Recall that the projection
$$
\pi_X:\bar\sM^{\vir}_{\ft,\fs}\to \Sym^\ell(X)
$$
was defined in \eqref{eq:GlobalProjectionToX}.
We will also write
\begin{align}
\label{eq:ProductOfProj1}
\pi_{\fs,X}:\bar\sM^{\vir}_{\ft,\fs}
\to M_{\fs}\times\Sym^\ell(X),
\\
\label{eq:ProductOfProj2}
\pi_{N,X}:\bar\sM^{\vir}_{\ft,\fs}
\to N_{\ft(\ell),\fs}(\delta)\times\Sym^\ell(X),
\end{align}
for the projection maps.
We use the same notation to denote the projections
from $\bar\sM^{\vir}_{\ft,\fs}/S^1$ (although we note that
$\pi_N$ maps this quotient to $N_{\ft(\ell),\fs}(\delta)/S^1$).

For any subspace $Y\subset X$, define
\begin{equation}
\label{eq:DefineSingularSupportLocus}
\sI(Y)\subset \bar\sM^{\vir,*}_{\ft,\fs}/S^1
\end{equation}
to be the subset of triples $[A,\Phi,\bx]$ where the support of $\bx$ contains
a point of $Y$.

\subsection{The incidence locus}
We now define subspaces of $\Sym^\ell(X)\times X$ and their Poincar\'e duals
to be used in computing expressions for cohomology classes dual to the geometric
representatives $\bar\sV(z)$ appearing in \eqref{eq:Duality}.

The incidence locus,
\begin{equation}
\label{eq:Incidence}
\sV_\ell(\Delta)\subset \Sym^\ell(X)\times X,
\end{equation}
is defined to be the set of points $(\bx,y)\in \Sym^\ell(X)\times X$ such
that $y$ is in the support of $\bx$.  Alternately, one can
describe $\sV_\ell(\Delta)$ by first defining
\begin{equation}
\label{eq:IncidenceInProduct}
\tilde\sV_\ell(\Delta) = \left( \bigcup_i
(\pi_i\times\id_X)^{-1}(\Delta_2)\right)
\subset X^\ell \times X,
\end{equation}
where $\pi_i:X^\ell\to X$ is projection onto the $i$-th factor and
$\Delta_2\subset X\times X$ is the diagonal and observing that
$\sV_\ell(\Delta)=\tilde\sV_\ell(\Delta)/\fS_\ell$. We will write
\begin{equation}
\label{eq:DefineIncidenceStratum}
\begin{aligned}
\tilde\sV_{\sP}(\Delta)&=(\Delta^\circ(X^\ell,\sP)\times X)\cap
\tilde\sV_\ell(\Delta),
\\
\sV_{\sP}(\Delta)&= (\Si(X^\ell,\sP) \times X)\cap
\sV_\ell(\Delta).
\end{aligned}
\end{equation}
Observe that $\sV_{\sP}(\Delta)=\sV_{\sP'}(\Delta)$ if the
partitions $\sP$ and $\sP'$ are in the same orbit of the
$\fS_\ell$ action on partitions of $N_\ell$.

For any smooth submanifold $Y\subseteqq X$, we define analogous
subspaces of $\Sym^\ell(X)\times Y$,
\begin{equation}
\label{eq:DefineIncidenceSubmanifold}
\begin{aligned}
\sV_\ell(\Delta,Y) &=\sV_\ell(\Delta)\cap\left(\Sym^\ell(X)\times Y\right),
\\
\tilde\sV_\ell(\Delta,Y)
&=\tilde\sV_\ell(\Delta)\cap\left(X^\ell\times Y\right),
\\
\sV_{\sP}(\Delta,Y) &=\sV_{\sP}(\Delta)\cap
\left(\Si(X^\ell,\sP)\times Y\right),
\\
\tilde\sV_{\sP}(\Delta,Y) &=\sV_{\sP}(\Delta)\cap
\left(\Delta^\circ(X^\ell,\sP)\times Y\right).
\end{aligned}
\end{equation}
We have the following description of the preceding subspaces.

\begin{lem}
\label{lem:ComponentsOfIncidence}
Let $Y\subseteqq X$ be an
oriented, smooth
submanifold with an oriented normal bundle and let $\sP$ be a partition of $N_\ell$. The subspace
$\tilde\sV_{\sP}(\Delta,Y)$ defined in
\eqref{eq:DefineIncidenceSubmanifold} is then the disjoint union
of components,
$$
\tilde\sV_\sP(\Delta,Y) = \bigsqcup_{P\in\sP} \tilde\sV_\sP(\Delta,Y,P),
$$
where
\begin{equation}
\label{eq:DefineIncidenceComponent} 
\tilde\sV_\sP(\Delta,Y,P): =
\{\left(x_P)_{P\in\sP},y\right)\in \tilde\sV_\sP(\Delta,Y):
x_P=y\}.
\end{equation}
Each component \eqref{eq:DefineIncidenceComponent} is a
codimension-four, oriented, smooth
submanifold of $\Delta^\circ(X^\ell,\sP)\times Y$ with an oriented normal bundle.
Under the action of $\fS(\sP)$ on $\Delta^\circ(X^\ell,\sP)\times
T$, the components $\tilde\sV_\sP(\Delta,Y,P)$ and
$\tilde\sV_\sP(\Delta,Y,P')$ are identified if $|P|=|P'|$.
\end{lem}

\begin{proof}
The fact that $\tilde\sV_\sP(\Delta,Y)$ is given by the union in the statement of the lemma
is clear.  The fact that the given union is disjoint follows by observing
that for any $\bx=(x_P)_{P\in\sP}\in\Delta^\circ(X^\ell,\sP)$, the points $x_P$ and $x_{P'}$
are always unequal for $P\neq P'$, so for $(\bx,y)\in
\tilde\sV_{\sP}(\Delta,Y)$, the equality $x_P=y=x_{P'}$ is
impossible.

To see that each component \eqref{eq:DefineIncidenceComponent} is a
codimension-four, oriented, smooth submanifold we
observe that
$\Delta_2(Y):=\{(y,y)\in X\times Y:y\in Y\}$ is a smooth submanifold of $X\times Y$, diffeomorphic to $Y$, with orientable normal bundle given by
the direct sum of orientable bundles,
$$
\{(w^\perp_Y,0)\in T(X\times Y)|_{\Delta_2(Y)}: w^\perp_Y\in\nu(Y)\}
\oplus
\{(w_Y,-w_Y)|_{\Delta_2(Y)}: w_Y\in TY\},
$$
where $\nu(Y)$ is the normal bundle of $Y$ in $X$.
Because the restriction of the projection
maps $\pi_i,\pi_j:X^\ell\to X$ to $\Delta^\circ(X^\ell,\sP)$ are equal for any $i,j\in P$,
the component
\eqref{eq:DefineIncidenceComponent} is the pre-image of
$\Delta_2(Y)$ under the map
$$
\pi_i\times\id_Y:\Delta^\circ(X^\ell,\sP)\times Y \to X\times Y
$$
for any $i\in P$.  Because the map $\pi_i\times\id_Y$ is a
submersion, this pre-image is then also a codimension-four, smooth
submanifold with normal bundle given by pullback of the normal bundle of $\Delta_2(Y)$.
Because the  normal bundle of $\Delta_2(Y)$ is orientable, its pullback by $\pi_i\times\id_Y$
is also orientable, as required.

Lastly, the proof of the observation on the action of $\fS(\sP)$ is straightforward.
\end{proof}

We will use the following computation of the cohomology of the
complement of $\sV_{\sP}(\Delta,Y)$.

\begin{lem}
\label{lem:CohomologyOfIncidenceSubmanifoldComplement}
Let $\sP$ be a partition of $N_\ell$, let $Y\subset X$ be a smooth
submanifold, and let $\sV_{\sP}(\Delta,Y)\subset
\Si(X^\ell,\sP)\times Y$ be the subspace defined in
\eqref{eq:DefineIncidenceSubmanifold}.  Then,
\begin{equation}
H_k\left( \Si(X^\ell,\sP)\times Y,\left(\Si(X^\ell,\sP)\times Y\right)\setminus \sV_{\sP}(\Delta,Y);\RR\right)
\cong
\begin{cases}
0 & \text{if $k<4$,}
\\
(\oplus_{P\in\sP}\RR)/\fS(\sP) & \text{if $k=4$.}
\end{cases}
\end{equation}
The group $H_4(\Si(X^\ell,\sP)\times Y,(\Si(X^\ell,\sP)\times Y)\setminus\sV_{\sP}(\Delta,Y);\RR)$ is generated by
$\phi_{\sP,P}:=(\tilde\pi_\ell\times \id_Y)\circ\tilde\phi_{\sP,P}$,
where
\begin{equation}
\label{eq:DiskTransvToIncid}
\tilde\phi_{\sP,P}:D^4\to
\Delta^\circ(X^\ell,\sP)\times Y
\end{equation}
is any map
intersecting $\tilde\sV_\sP(\Delta,Y,P)$ transversely at the origin
and $\tilde\pi_\ell:X^\ell\to\Sym^\ell(X)$ is the projection.
\end{lem}

\begin{proof}
The conclusion follows immediately from
the Thom isomorphism,
\cite[Theorem 5.7.10]{Spanier},
\cite[Theorem 10.4]{MilnorStasheff},
or \cite[Equation VIII.11.3]{Dold}, and
the assertion in Lemma
\ref{lem:ComponentsOfIncidence} that $\sV_{\sP}(\Delta,Y)$ is a
disjoint union of smooth, codimension-four submanifolds
with oriented normal bundles,
enumerated by $\sP/\fS(\sP)$, of $\Si(X^\ell,\sP)\times Y$.
\end{proof}

\subsection{Cohomology classes}
\label{subsec:DefiningCohomClasses}
We now define some cohomology classes on
$\bar\sM^{\vir,*}_{\ft,\fs}/S^1$
to express our computations of the duals of the geometric representatives of
$\bar\sV(z)$ and $\bar\sW$.

\begin{defn}
\label{defn:DefineNu}
Let $\nu\in H^2(\bar\sM^{\vir,*}_{\ft,\fs}/S^1;\ZZ)$ be the
first Chern class of the $S^1$ bundle,
\begin{equation}
\bar\sM^{\vir,*}_{\ft,\fs}\to \bar\sM^{\vir,*}_{\ft,\fs}/S^1,
\end{equation}
where the $S^1$ action is defined in \eqref{eq:DefineGlobalS1Action}.
If
\begin{equation}
\label{eq:LineBundleForS1ZAction}
\LL_\nu=\bar\sM^{\vir,*}_{\ft,\fs}\times_{S^1}\CC,
\end{equation}
is  the complex line bundle associated to this $S^1$ action,
then $c_1(\LL_\nu)=\nu$.
\end{defn}

\begin{defn}
\label{defn:LineBundleForMus}
For $\beta\in H_\bullet(X;\RR)$,
let $\mu_{\fs}(\beta)\in H^{2-\bullet}(M_{\fs};\ZZ)$ be
the $\mu$-class defined in
\eqref{eq:SWMuMap}.
We will use the same notation for the pullback by the
projection $\pi_{\fs}$ of these classes to
$\bar\sM^{\vir}_{\ft,\fs}/S^1$.  We will also use the notation
$$
\LL_{\fs}
=
\tM_{\fs}\times_{\sG_{\fs}}(X\times\CC)
\to
M_{\fs}\times X
$$
for the restriction of the bundle defined in \eqref{eq:DefineSWUniversal}.
\end{defn}

By \cite[Theorem II.19.2]{BredonSheafTheory}, there is an isomorphism,
\begin{equation}
\label{eq:SymmetricProductCohomology}
H^\bullet(\Sym^\ell(X);\RR)
\cong
H^\bullet(X^\ell;\RR)^{\fS_\ell}.
\end{equation}
We can thus make the following definition.
\begin{defn}
\label{defn:DefnSymmBeta}
If $\beta\in H_\bullet(X;\RR)$ and $\pi_i:X^\ell\to X$
is projection onto the $i$-th factor,
define
$$
\tilde S^\ell(\beta):=\sum_i\pi_i^*\PD[\beta]\in H^{4-\bullet}(X^\ell;\RR)
$$
and let $S^\ell(\beta)\in H^{4-\bullet}(\Sym^\ell(X);\RR)$ be the
cohomology class satisfying $\tilde\pi_\ell^*S^\ell(\beta)=\tilde S^\ell(\beta)$,
where
$\varpi:X^\ell\to\Sym^\ell(X)$ is the projection,
as specified by Equation \eqref{eq:SymmetricProductCohomology}, and $\PD$ denotes the
Poincar\'e duality isomorphism.

For $\alpha\in H^\bullet(X;\RR)$, define
$$
\tilde S^\ell(\alpha):=\sum_i\pi_i^*\alpha\in H^{\bullet}(X^\ell;\RR)
$$
and let $S^\ell(\alpha)\in H^{\bullet}(\Sym^\ell(X);\RR)$ be
the cohomology class satisfying
$\varpi^*S^\ell(\alpha)=\tilde S^\ell(\alpha)$.
\end{defn}

Although the incidence locus $\sV_\ell(\Delta)$ defined in \eqref{eq:Incidence} is not a submanifold
of $\Sym^\ell(X)\times X$, we can define a cohomology class which will serve the same role as
a Poincar\'e dual for the purposes of counting intersections as follows.
Recall that if $Y$ is a smooth, closed oriented $m$-dimensional submanifold of an oriented, closed $n$-dimensional manifold $Z$ then
by \cite[Proposition VIII.11.18]{Dold},
the Poincar\'e dual of the fundamental class of $Y$, namely $[Y]\in H_m(Z;\ZZ)$, and the Thom class of
the normal bundle of $Y$, namely $\Th(Y)\in H^{n-m}(Z,Z\setminus Y;\ZZ)$, are related by
\begin{equation}
\label{eq:PDualAndThomClass}
\PD_M([Y])=\jmath^*\Th(Y),
\end{equation}
where $\jmath:(Z,\emptyset)\to (Z,Z\setminus Y$) is the inclusion map.

Because the diagonal $\Delta_2\subset X\times X$ is a smooth submanifold,
it has a Thom class and a Poincar\'e dual,
$$
\Th(\Delta_2)\in H^4(X\times X,X\times X\setminus\Delta_2;\ZZ)
\quad\text{and}\quad
\PD[\Delta_2]\in H^4(X\times X;\ZZ),
$$
satisfying \eqref{eq:PDualAndThomClass}.
The relation \eqref{eq:IncidenceInProduct} between
the incidence locus
$\sV_\ell(\Delta)$ and $\Delta_2$ implies that
it makes sense to define
\begin{subequations}
\begin{align}
\label{eq:DefineThomClassOfIncidenceInProduct}
T_{\tilde\sV}
{}&:=
\sum_i (\pi_i\times\id_X)^*\Th(\Delta_2)
\in H^4(X^\ell\times X,X^\ell\times X \setminus \tilde\sV_\ell(\Delta);\RR),
\\
\label{eq:DefinePDOfIncidenceInProduct}
\PD[\tilde\sV]
{}&:=\sum_i (\pi_i\times\id_X)^*\PD[\Delta_2]\in H^4(X^\ell\times X;\RR),
\end{align}
\end{subequations}
where $\pi_i:X^\ell\to X$ is projection onto the $i$-th factor
and we take images of integer cohomology classes in real cohomology.

The isomorphism \eqref{eq:SymmetricProductCohomology} and
the $\fS_\ell$ symmetry of $T_{\tilde\sV}$ and
$\PD[\tilde\sV]$ allow us to define
\begin{subequations}
\label{eq:DefinePDIncid_and_ThomSubmanifoldIncidence}
\begin{align}
\label{eq:DefinePDIncid}
\PD[\sV]&\in H^4(\Sym^\ell(X)\times X;\RR),
\\
\label{eq:DefineThomSubmanifoldIncidence}
T_{\sV}&\in H^4(\Sym^\ell(X)\times X,\Sym^\ell(X)\times X \setminus \sV_\ell(\Delta);\RR)
\end{align}
\end{subequations}
to be the unique cohomology classes satisfying
\begin{equation}
\label{eq:PullbackToProductDefinitionOfPDandThomOfIncidence}
(\varpi_\ell\times\id_X)^*\PD[\sV]=\PD[\tilde\sV]
\quad\text{and}\quad
(\varpi_\ell\times\id_X)^*T_{\sV}=T_{\tilde\sV},
\end{equation}
where $\varpi_\ell:X^\ell\to \Sym^\ell(X)$ denotes the
projection map.
There is a commutative diagram
$$
\begin{CD}
\left(X^\ell\times X,\emptyset \right)
@> \tilde \jmath >>
\left( X^\ell\times X,X^\ell\times X \setminus \tilde\sV_\ell(\Delta) \right)
\\
@V \varpi_\ell\times\id_X VV @V \varpi_\ell\times\id_X   VV
\\
\left(\Sym^\ell(X)\times X,\emptyset \right)
@>  \jmath >>
\left( \Sym^\ell(X)\times X,\Sym^\ell(X)\times X \setminus \sV_\ell(\Delta) \right)
\end{CD}
$$
where the horizontal maps are inclusions.  By the commutativity of this diagram,
the definitions \eqref{eq:DefinePDIncid_and_ThomSubmanifoldIncidence},
and the relation between the Thom class and Poincar\'e dual \eqref{eq:PDualAndThomClass} satisfied by
$\PD[\Delta_2]$ and $\Th(\Delta_2)$, we have
\begin{equation}
\tilde \jmath^*T_{\tilde\sV}=\PD[\tilde \sV]
\quad\text{and}\quad
\jmath^*T_{\sV}=\PD[\sV].
\end{equation}
For a closed, oriented, smooth submanifold, $Y\subseteqq X$,
let $\iota_\sP:\Si(X^\ell,\sP)\to \Sym^\ell(X)$ and $\iota_Y:Y\to X$
denote the inclusion maps.  There is an inclusion map of pairs,
$$
\begin{CD}
(\Si(X^\ell,\sP)\times Y,\Si(X^\ell,\sP)\times Y\setminus\sV_\Si(\sP,Y))
\\
@V \iota_{\sP}\times\iota_Y VV
\\
(
\Sym^\ell(X)\times X,
\Sym^\ell(X)\times X\setminus\sV_\ell(\Delta)
)
\end{CD}
$$
The following computation, together with
Lemma \ref{lem:CohomologyOfIncidenceSubmanifoldComplement}, will be useful
in characterizing $\PD[\sV]$.

\begin{lem}
\label{eq:IncidenceLocusThomClassEval}
Let $Y\subseteqq X$ be a closed, oriented, smooth submanifold, $T_{\sV}$ be the cohomology class defined in
\eqref{eq:DefineThomSubmanifoldIncidence}, and $[\phi_{\sP,P}]$ be as defined prior to
equation \eqref{eq:DiskTransvToIncid}. Then,
$$
\langle (\iota_{\sP}\times\iota_Y)^*T_{\sV},[\phi_{\sP,P}]
\rangle
=
|P|,
$$
where $\iota_\sP$ and $\iota_Y$ are the inclusion maps defined above.
\end{lem}

\begin{proof}
By construction, the map $\phi_{\sP,P}$ is covered by the map
$\tilde\phi_{\sP,P}:D^4\to\Delta^\circ(X^\ell,\sP)\times Y$
defined in \eqref{eq:DiskTransvToIncid}
in the sense that $(\tilde\pi_\ell\times \id_Y)\circ \tilde\phi_{\sP,P}=\phi_{\sP,P}$.
This gives the equality
$$
(\tilde\pi_\ell\times\id_Y)_*[\tilde\phi_{\sP,P}]
=
[\phi_{\sP,P}].
$$
The inclusion $\iota_{\sP}\times\iota_Y$ is thus covered by
an inclusion map
$$
\tilde\iota_{\sP}\times\iota_Y:
\Delta^\circ(X^\ell,\sP)\times Y \to X^\ell\times X,
$$
in the sense that $(\tilde\pi_\ell\times \id_X)\circ (\tilde\iota_{\sP}\times\iota_Y)=\iota_{\sP}\times\iota_Y$.
We can
therefore compute that
\begin{align*}
\left\langle (\iota_{\sP}\times\iota_Y)^*T_{\sV},[\phi_{\sP,P}]
\right\rangle
{}&=
\left\langle (\iota_{\sP}\times\iota_Y)^*T_{\sV},(\tilde\pi_\ell\times\id_Y)_*[\tilde\phi_{\sP,P}]
\right\rangle
\\
{}&=
\left\langle
(\varpi_\ell\times\id_Y)^*(\iota_{\sP}\times\iota_Y)^*T_{\sV},[\tilde\phi_{\sP,P}]
\right\rangle
\\
{}&=
\left\langle
(\tilde\iota_{\sP}\times\iota_Y)^*(\tilde\pi_\ell\times\id_X)^*T_{\sV},[\tilde\phi_{\sP,P}]
\right\rangle
\\
{}&=
\left\langle
(\tilde\iota_{\sP}\times\iota_Y)^*T_{\tilde\sV},[\tilde\phi_{\sP,P}]
\right\rangle
\quad\text{(by \eqref{eq:PullbackToProductDefinitionOfPDandThomOfIncidence})}
\\
{}&=
\left\langle
\sum_i(\pi_i\times\id_X)^*\Th(\Delta_2),(\tilde\iota_{\sP}\times\iota_Y)_*[\tilde\phi_{\sP,P}]
\right\rangle
\quad\text{(by \eqref{eq:DefineThomClassOfIncidenceInProduct})}
\\
{}&=
\sum_i
\left\langle \Th(\Delta_2),(\pi_i\circ\iota_\sP)\times\iota_Y)_*[\tilde\phi_{\sP,P}]\right\rangle.
\end{align*}
By the argument in Lemma \ref{lem:ComponentsOfIncidence}, the image of
$\Delta^\circ(X^\ell,\sP)\times Y$ under the map $\pi_i\times\iota_Y$
intersects the diagonal $\Delta_2$ transversely.  Hence, the image of the map
$$
(\pi_i\circ\iota_\sP)\times\iota_Y)\circ\tilde\phi_{\sP,P}: D^4\to X\times X
$$
intersects the diagonal transversely at one point if
$i\in P$ and is disjoint from the diagonal if $i\notin P$.
This yields the equality
$$
\langle \Th(\Delta_2),(\pi_i\circ\iota_\sP)\times\iota_Y)_*[\tilde\phi_{\sP,P}]\rangle
=
\begin{cases}
1 & \text{if $i\in P$,} \\
0 & \text{if $i\notin P$,}
\end{cases}
$$
and this in turn gives the equality in the statement of the lemma.
\end{proof}

The equality on Poincare duality given in
\cite[Theorem 30.6]{GreenbergHarper}, the definition of
$\PD[\sV]$ following \eqref{eq:DefinePDIncid},
the definition of $S^\ell(\beta)$ in
Definition \ref{defn:DefnSymmBeta},
and the isomorphism \eqref{eq:SymmetricProductCohomology},
yield the following

\begin{lem}
\label{lem:PDual}
For any $\beta\in H_\bullet(X;\RR)$,
$$
\PD[\sV]/\beta= S^\ell(\beta),
$$
where $S^\ell(\beta)$ is  defined in
Definition \ref{defn:DefnSymmBeta}
and $\PD[\sV]$ in \eqref{eq:DefinePDIncid}.
\end{lem}

\section{Fundamental class of the virtual link of the ideal moduli space of Seiberg--Witten monopoles}
\label{sec:FundClassOfAmbLink}
We now define the fundamental class of the virtual link appearing in
\eqref{eq:Duality},
$$
[\bar\bL^{\vir}_{\ft,\fs}]\in
H_{d(\ft)-2}\left(\bar\sM^{\vir,*}_{\ft,\fs}/S^1;\RR\right),
$$
where $d(\ft)$ is the dimension of $\sM^{\vir,*}_{\ft,\fs}$.
The union of the lower strata in $\bar\bL^{\vir}_{\ft,\fs}$ will be
denoted $\bL^{\sing}_{\ft,\fs}\subset\bar\bL^{\vir}_{\ft,\fs}$.

\begin{lem}
\label{lem:NghOfEnd}
Let $\bar\sV(z)$ and $\bar\sW$ be the geometric representatives
appearing in \eqref{eq:Duality} and
let $\bar\bchi$ be the obstruction section appearing in
Hypothesis \ref{hyp:Gluing}.
For any neighborhood $V$ of
$\bL^{\sing}_{\ft,\fs}$ in the lower strata of
$\bar\sM^{\vir,*}_{\ft,\fs}/S^1$ such that
$\cl(V)\cap\bar\sV(z)\cap\bar\sW^{\eta}\cap\bar\bchi^{-1}(0)$
is empty, there is a neighborhood $U$ of the singular strata in
$\bar\sM^{\vir,*}_{\ft,\fs}/S^1$ such that the following hold:
\begin{enumerate}
\item
\label{item:NghOfEnd1}
There is a deformation retraction $r$ of $U$ onto $V$ which
respects the level sets of the function $t_N$ given in
\eqref{eq:NTubularDist}.
\item
\label{item:NghOfEnd2}
There is a compact topological manifold with
boundary, $\hat\bL$,
given as the boundary of a smooth manifold with corners
with
$\bL^{\vir}_{\ft,\fs}\setminus U\subset\hat\bL\subset
\bar\bL^{\vir}_{\ft,\fs}$.
\item
\label{item:NghOfEnd3}
The intersection
$\bar\sV(z)\cap\bar\sW^{\eta}\cap\bar\bchi^{-1}(0)\cap U$ is
empty.
\end{enumerate}
\end{lem}

\begin{proof}
Let
$\sM^{\sing,*}_{\ft,\fs}/S^1\subset\bar\sM^{\vir,*}_{\ft,\fs}/S^1$
be the union of the singular strata. The local Whitney stratified
structure of the spliced-ends moduli space, specifically the
property described in Lemma \ref{lem:SplicedEndNDR}, implies that
the union of the lower strata is locally a neighborhood
deformation retraction  and thus a neighborhood deformation
retraction (NDR)
(combine the results \cite[Proposition 5]{Goresky_Triang_Strat_Objects} and \cite[Corollary 3.3.11]{Spanier}).
The fact that this deformation retraction respects the level
sets of $t_N$ follows by observing that the deformation
retractions defined in Lemma \ref{lem:SplicedEndNDR} are invariant
with respect to the frame action.  The NDR structure gives a
cofinal sequence of neighborhoods of $\sM^{\sing,*}_{\ft,\fs}/S^1$
which deformation retract onto $\sM^{\sing,*}_{\ft,\fs}/S^1$.
Because
$\cl(V)\cap\bar\sV(z)\cap\bar\sW^{\eta}\cap\bar\bchi^{-1}(0)$
is empty and $\bar\sM^{\vir,*}_{\ft,\fs}/S^1$ is normal, we
can find a neighborhood deformation retraction, $r:N\to
\sM^{\sing,*}_{\ft,\fs}/S^1$ onto $\sM^{\sing,*}_{\ft,\fs}/S^1$, such
that $U=r^{-1}(V)$ satisfies the
conclusions \eqref{item:NghOfEnd1} and \eqref{item:NghOfEnd3} of
the lemma. Define $\hat\bL$ to be the complement in
$\bL^{\vir}_{\ft,\fs}$ of tubular neighborhoods (in the Thom--Mather
sense) of the local strata of $\sM^{\sing,*}_{\ft,\fs}/S^1$.  By
choosing these tubular neighborhoods to be contained in $N$, the
resulting topological manifold with boundary will satisfy the
conclusion \eqref{item:NghOfEnd2}.
\end{proof}

For the subspaces  $U$ and $V$ appearing in Lemma \ref{lem:NghOfEnd},
the  manifold with
boundary $\hat\bL$
defined in Lemma \ref{lem:NghOfEnd} has a fundamental class,
\begin{equation}
\label{eq:FundClass1}
[\hat\bL,\rd\hat\bL]\in
H_{d(\ft)-2}(\bL^{\vir}_{\ft,\fs}\cup U^{\circ},U^{\circ};\ZZ),
\end{equation}
where $U^{\circ}:=U\setminus\sM^{\sing,*}_{\ft,\fs}/S^1$.  Note that there is an
excision isomorphism,
$$
\iota_* :
H_\bullet(\bL^{\vir}_{\ft,\fs}\cup U^{\circ},U^{\circ};\ZZ)
\cong
H_\bullet(\bar\bL^{\vir}_{\ft,\fs}\cup U,U;\ZZ).
$$
Because $V$ is contained
in the codimension-four subspace $\sM^{\sing,*}_{\ft,\fs}/S^1$, the set
$V$ has dimension three less than $\bL^{\vir}_{\ft,\fs}$.
Because $U$ deformation retracts onto $V$, we have $H_k(U;\ZZ)\cong H_k(V;\ZZ)=0$ for
$k> \dim\hat\bL -3$. Hence,
the inclusion map of pairs,
\begin{equation}
\label{eq:CompactJMap} \bar\jmath: \left(
\bar\bL^{\vir}_{\ft,\fs}\cup U,\emptyset \right) \to \left(
\bar\bL^{\vir}_{\ft,\fs}\cup U,U \right),
\end{equation}
induces an isomorphism,
$$
\bar\jmath_*: H_{d(\ft)-2}(\bar\bL^{\vir}_{\ft,\fs}\cup U;\ZZ) \cong
H_{d(\ft)-2}(\bar\bL^{\vir}_{\ft,\fs}\cup U,U;\ZZ).
$$
We now define
\begin{equation}
\label{eq:DefineAmbientLinkFund}
[\bar\bL^{\vir}_{\ft,\fs}]
:=
\jmath_*^{-1}\iota_*[\hat\bL,\rd\hat\bL] \in
H_{d(\ft)-2}(\bar\bL^{\vir}_{\ft,\fs}\cup U;\ZZ)
\end{equation}
to be the fundamental class of the virtual link.

\section{Computation of the $\mu$-classes}
\label{sec:ComputingMuClasses}
The geometric representatives $\bar{\sV}(z)$ and $\bar\sW$
are dual to the cohomology classes $\mu_p(z)$ and $\mu_c$
defined in Section \ref{sec:Cohomology}.  We now compute
the pullbacks $\bga^*_{\sM}\mu_p(z)$ and $\bga^*_{\sM}\mu_c$
in terms of the cohomology classes defined in
Section \ref{subsec:DefiningCohomClasses}.

\subsection{Geometric representatives and cocycles}
\label{subsec:GRepsAndCocycles}
We first recall the

\begin{defn}
\label{defn:GeomRepresentative}
(See Kronheimer and Mrowka \cite[p. 588]{KMStructure}.)
Let $Z$
be a smooth manifold.  A {\em geometric
representative\/} for $\mu\in H^d(Z;\RR)$
is a closed, smoothly stratified subspace $\sV$ of $Z$
together with a real coefficient, $q$, the {\em multiplicity\/},
with the following properties:
\begin{enumerate}
\item
\label{item:GeomRepresentative1}
The top stratum $\sV_0$ of $\sV$ is a codimension-$d$, smooth submanifold of
$Z$ with an oriented normal bundle.
\item
\label{item:GeomRepresentative2}
All other strata of $\sV$ have codimension $d+2$ or greater in $Z$.
\item
\label{item:GeomRepresentative3}
The pairing of $\mu$ with a homology class $h$ of
dimension $d$ is obtained by choosing a smooth singular cycle
$\sigma$ representing $h$, transverse to all strata of $\sV$,
then taking $q$ times the count (with signs) of the intersection
points between the cycle and the top stratum of $\sV$,
$$
\langle \mu,h\rangle = q\cdot \#(\sV_0\cap\sigma).
$$
\end{enumerate}
\end{defn}
A choice of a geometric representative $\sV$ for a cohomology class $\mu$
specifies a singular cocycle $c$ representing $\mu$, supported
on $\sV$, as described in the following lemma.

\begin{lem}
\label{lem:DefineSingularCocycle}
Let $Z$ be a smooth manifold
and let $(\sV,q)$ be a geometric representative for a real
cohomology class $\mu$ of dimension $d$ on $Z$. Then there is a class
$c\in H^d(Z,Z\setminus\sV;\RR)$ such that $\jmath_{\sV}^*c=\mu$,
where
$$
\jmath_{\sV}:\left(Z,\emptyset\right)\to \left(Z,Z\setminus V\right)
$$
is the inclusion map of pairs.
\end{lem}

\begin{proof}
Let $\Delta^{\8}_\bullet(Z)$ denote
the chain complex of smooth singular chains \cite[p. 291]{BredonTopGeom} with
$Z^{\8}_d(Z)\subset \Delta^{\8}_d(Z)$ and $B^{\8}_d(Z)\subset \Delta^{\8}_d(Z)$
denoting the submodules of
cycles and boundaries respectively.
As we are considering real cohomology, we shall view $\Delta^{\8}_\bullet(Z)$ as
a real vector space rather than as a $\ZZ$-module.
Property \eqref{item:GeomRepresentative3} in Definition
\ref{defn:GeomRepresentative} implies that a geometric representative
$(\sV,q)$ for $\mu\in H^d(Z;\RR)$ defines an element of
$\Hom(Z^{\infty}_d(Z);\RR)$ which is a representative of the cohomology class $\mu$.
However,
the intersection number $\#(\sV_0\cap\sigma)$ in Property \eqref{item:GeomRepresentative3}
is well-defined if
$\si$ is a smooth, singular chain whose boundary lies in $Z\setminus\sV$,
that is,
$\si\in\rd^{-1}(\Delta^{\infty}_{d-1}(Z\setminus\sV))$.
The intersection number $\#(\sV_0\cap\sigma)$
vanishes if $\si$ is a boundary or an element of
$\Delta^{\infty}_d(Z\setminus\sV)$.
Hence,
a geometric representative $(\sV,q)$ of $\mu\in H^d(Z;\RR)$
defines
an element of
\begin{equation}
\label{eq:SmoothRelCochains}
\Hom\left(
\rd^{-1}(\Delta^{\infty}_{d-1}(Z\setminus\sV))/\left(B^{\8}_d(Z)+\Delta^{\8}_d(Z\setminus\sV)\right);\RR\right).
\end{equation}
We now argue that \eqref{eq:SmoothRelCochains} is isomorphic to
$H^d(Z,Z\setminus\sV;\RR)$.
By the de Rham Theorem
(see the discussion in \cite[p. 291]{BredonTopGeom}),
there is a functorial isomorphism between $H^{\bullet}(Z;\RR)$
and the homology of the complex $\Hom(\Delta^{\8}_{\bullet}(Z);\RR)$.
The Universal Coefficient Theorem then identifies $H_{\bullet}(Z;\RR)$
with the homology of the complex $\Delta^{\8}_{\bullet}(Z)$.
Thus, we have isomorphisms,
\begin{equation}
H_{\bullet}(Z;\RR)
\cong
H_{\bullet}(\Delta^{\8}_{\bullet}(Z);\RR),
\quad
H_{\bullet}(Z\setminus\sV;\RR)\cong H_{\bullet}(\Delta^{\8}_{\bullet}(Z\setminus\sV);\RR).
\end{equation}
The Five Lemma then identifies $H_{\bullet}(Z,Z\setminus\sV;\RR)$ with the homology of
the complex
$$
\Delta^{\8}_{\bullet}(Z)/\Delta^{\8}_{\bullet}(Z\setminus\sV),
$$
and thus,
$$
H_d(Z,Z\setminus\sV;\RR)
\cong
\rd^{-1}(\Delta^{\infty}_{d-1}(Z\setminus\sV))/\left(B^{\8}_d(Z)+\Delta^{\8}_d(Z\setminus\sV)\right).
$$
The Universal Coefficient Theorem then yields the required
isomorphism between \eqref{eq:SmoothRelCochains} and
$H^d(Z,Z\setminus\sV;\RR)$.
\end{proof}

\subsection{Cocycles as pullbacks}
\label{subsubsec:CocyclesPullbacks}
For $z=\beta_1\cdots\beta_r\in\AAA(X)$, where
$\beta_i\in H_{\dim\beta_i}(X;\RR)$,
we now review the definition of the cocycles defined by
$\bar\sV(\beta_i)$ and $\bar\sW$.

Let $T\subset X$ be a smooth submanifold with
$q[T]=\beta$ for $q\in\RR$ and let $\nu(\beta)\subset X$ be a tubular
neighborhood of $T$.
Let $\sB^*(\nu(\beta))$ be the space of irreducible connections
on the bundle $\fg_\ft|_{\nu(\beta)}$
(see \cite[Section 2 (ii)]{KMStructure} or \cite[Section 3.2.2]{FL2b})
and (abusing notation slightly) let
$$
r_\beta:\sM^{\vir}_{\ft,\fs}/S^1 \ni [A,\Phi] \mapsto [\hat A|_{\nu(\beta)}] \in \sB^*(\nu(\beta)),\quad
$$
be the composition of the gluing map $\bga_{\sM}$ and the
restriction map.
For $\sI(\nu(\beta))$ as defined in \eqref{eq:DefineSingularSupportLocus},
let
\begin{equation}
\label{eq:ExtendedRestrictionMap}
\bar r_{\beta}:
\bar\sM^{\vir,*}_{\ft,\fs}/S^1 \setminus\sI(\nu(\beta))
\ni [A,\Phi,\bx]
\mapsto
[\hat A|_{\nu(\beta)}] \in
\sB^*(\nu(\beta))
\end{equation}
be the extension of $r_\beta$.
We will denote the image of $r_\beta$ by $\sM_\beta$ and
the image of $\bar r_\beta$ by $\bar\sM_{\beta}$.

Let $\sV_T\subset \sB^*(\nu(\beta))$ be the geometric
representative for the cohomology class $\mu_p(\beta)$ defined in
\cite[Section 2 (ii)]{KMStructure}.  From the definition in \cite[Section 3.2.3]{FL2b},
it follows that
$$
\bar\sV(\beta)\cap
\left(
\bar\sM^{\vir,*}_{\ft,\fs}/S^1 \setminus\sI(\nu(\beta))
\right)=
\bar r_\beta^{-1}(\sV_T).
$$
Note that in \cite[Section 3.2.3]{FL2b}, a tubular neighborhood of
the union of $T$ and certain loops (a ``suitable'' neighborhood
in the sense of \cite{KMStructure})
were used instead of
the tubular neighborhood $\nu(\beta)$ in defining $\bar\sV(\beta)$.
However, by \cite[Lemma 5.5]{FLLevelOne} the geometric representative $\bar\sV(\beta)$
can be replaced with
geometric representatives pulled back from tubular rather than
suitable neighborhoods of $T$,
when computing the intersection
numbers in \eqref{eq:Duality}.
Let
\begin{equation}
\label{eq:CocycleOnRestrictedSpace}
[c_{T,\beta}]\in H^{4-\dim\beta}\left(\bar\sM_{\beta},\bar\sM_{\beta}\setminus\sV_T;\RR\right)
\end{equation}
be the cohomology class defined by the geometric representative
$\sV_T$ as described in
Lemma \ref{lem:DefineSingularCocycle}
and define
\begin{equation}
\label{eq:DefineRelMuClass}
[c_\beta]:=\bar r_\beta^*[c_{T,\beta}]
\in
H^{4-\bullet}
\left(
\bar\sM^{\vir,*}_{\ft,\fs}/S^1\setminus\sI(\nu(\beta)),
\bar\sM^{\vir,*}_{\ft,\fs}/S^1\setminus \left(\sI(\nu(\beta))\cup\bar\sV(\beta)\right);\RR
\right).
\end{equation}
A similar construction defines a cocycle representing $\bar\sW$ of the
cohomology class $\mu_c$.  For
the configuration space $\sC^{*,0}_{\ft}(\nu(x))$ of pairs on the restriction of the \spinu structure
$\ft$ to $\nu(x)$ (see \cite[Section 3.2.2]{FL2b}), we define
\begin{equation}
\label{eq:RestrictionForConfigSpace}
\bar r_x: \bar\sM^{\vir,*}_{\ft,\fs}/S^1 \setminus\sI(\nu(x))
\to
\sC^{*,0}_{\ft}(\nu(x))/S^1
\end{equation}
to be the composition of the gluing map $\bga_{\sM}$ with the restriction map.
In \cite[Equation (3.14)]{FL2b}, the geometric representative
$\bar\sW$
is defined by
$$
\bar\sW\cap
\left( \bar\sM^{\vir,*}_{\ft,\fs}/S^1 \setminus\sI(\nu(x))\right)
=
\bar r_x^{-1}(\sW_T),
$$
where $\sW_T$ is the zero-locus of a section of an appropriate line
bundle.
Define
\begin{equation}
\label{eq:DefinePreCocycleForWClass0}
[c_{T,W}]
\in H^2(\sC^{*,0}_{\ft}(\nu(x))/S^1,\sC^{*,0}_{\ft}(\nu(x))/S^1\setminus\sW_T;\RR)
\end{equation}
to be the cocycle obtained from the geometric representative  $\sW_T$ by
Lemma \ref{lem:DefineSingularCocycle} and define
\begin{equation}
\label{eq:DefinePreCocycleForWClass}
[c_\sW]:=\bar r_x^*[c_{T,W}]\in
H^2(\bar\sM^{\vir,*}_{\ft,\fs}/S^1\setminus\sI(\nu(x)),\bar\sM^{\vir,*}_{\ft,\fs}/S^1\setminus\bar\sW;\RR)
\end{equation}
to be its pullback.

\subsection{Computations of cocycles}
\label{subsubsec:ComputCocycles}
In Section \ref{subsec:InitialDuality}, we will write the intersection number
appearing in \eqref{eq:Duality}
in terms of a cup product of $[c_\sW]$ and
$[c_{\beta_i}]$. We will then compute this product
in terms of the cohomology classes
defined in Section \ref{subsec:DefiningCohomClasses}.
For this purpose, we will need to relate the relative cohomology classes
$[c_\sW]$ and $[c_{\beta_i}]$ to absolute cohomology
classes. Specifically,
we will need to compute the cohomology class $\bar\jmath_\beta^*\bar r_\beta^*[c_{T,\beta}]$
in terms of the cohomology classes defined in Section \ref{subsec:DefiningCohomClasses},
where
\begin{equation}
\label{eq:DefinePairsInclusionForBeta}
\bar\jmath_\beta:
\left(\bar\sM^{\vir,*}_{\ft,\fs}/S^1 \setminus\sI(\nu(\beta)),
\emptyset\right)
\to
\left(\bar\sM^{\vir,*}_{\ft,\fs}/S^1 \setminus\sI(\nu(\beta)),
\bar\sM^{\vir,*}_{\ft,\fs}/S^1 \setminus \left( \sI(\nu(\beta))\cup\bar\sV_\beta\right)\right)
\end{equation}
is the inclusion of pairs.
To that end, we consider the bundle
\begin{equation}
\label{eq:TubNghBundle}
\FF_T
:=
\sA^*_{\nu(T)}\times_{\sG_{\nu(T)}} \fg_{\ft}|_T
\to
\sB^*_{\nu(T)}\times T.
\end{equation}
If $\bar\sM_{\beta}$ is as defined in the paragraph following \eqref{eq:ExtendedRestrictionMap}, then the construction of the geometric representative $\sV_T$ implies that
$$
\jmath_T^*[c_{T,\beta}]=p_1(\FF_{T})/[T],
$$
where
$$
\jmath_T:
\left(\bar\sM_{\beta}\times T, \emptyset\right)
\to
\left(\bar\sM_{\beta}\times T,\bar\sM_{\beta}\times T\setminus\sV_T\right)
$$
is the inclusion of pairs.
To compute $\bar r_{\beta}^*\jmath_T^*[c_{T,\beta}]$,
we shall compare the Pontrjagin classes of the bundles
\begin{equation}
\label{eq:ComparingPontrjClasses}
(\bar r_{\beta}\times\id_T)^*\FF_T\quad\text{and}\quad (\id_{\sM}\times\iota_T)^*\FF^{\vir}_{\ft,\fs},
\end{equation}
where
$$
\id_{\bar\sM}\times\iota_T:
\bar\sM^{\vir,*}_{\ft,\fs}/S^1\times T
\to
\bar\sM^{\vir,*}_{\ft,\fs}/S^1\times X
$$
is the inclusion map
and $\FF^{\vir}_{\ft,\fs}\to \bar\sM^{\vir,*}_{\ft,\fs}/S^1\times X$
is the bundle defined by restricting
\begin{equation}
\label{eq:VirtualUniversal}
\FF^{\vir}_{\ft,\fs}=
\bigcup_{\sP}\left(
\tN_{\ft(\ell),\fs}(\delta)\times\bar\Gl(\ft,\fs,\sP)\times_{\sG_{\fs}\times S^1} \fg_{\ft(\ell)}\right)
\end{equation}
to $\bar\sM^{\vir,*}_{\ft,\fs}/S^1\times X$.
We begin with the computation of $p_1(\FF^{\vir}_{\ft,\fs})$.

\begin{lem}
\label{lem:PontrjaginOfVirtualUniv}
Assume the \spinu structure $\ft(\ell)$ admits a splitting
$\ft(\ell)=\fs\oplus \fs\otimes L$ as described in Section \ref{subsubsec:SpincuStr},
where $L\to X$ is a complex line bundle.
Let $\FF^{\vir}_{\ft,\fs}\to \bar\sM^{\vir,*}_{\ft,\fs}/S^1\times X$
be the bundle defined in \eqref{eq:VirtualUniversal} and $\LL_\fs$ the
line bundle given in Definition \ref{defn:LineBundleForMus}.
Let $\pi_{X,2}:\sM^{\vir,*}_{\ft,\fs}/S^1\times X\to X$
be the projection onto the second factor.  Then
\begin{equation}
\label{eq:PontrjaginOfVirtualUniv}
p_1(\FF^{\vir}_{\ft,\fs})
=
\left(
2(\pi_{\fs}\times\id_X)^*c_1(\LL_\fs)-\nu+\pi_{X,2}^*c_1(L)
\right)^2.
\end{equation}
\end{lem}

\begin{proof}
Equation \eqref{eq:PontrjaginOfVirtualUniv} follows from \cite[Lemma 4.7]{FLLevelOne}.
\end{proof}

To compare the Pontrjagin classes in \eqref{eq:ComparingPontrjClasses},
we first describe a subspace on which the bundles
are isomorphic.
Let
$\sV_\sP(\Delta,T)\subset\Si(X^\ell,\sP)\times T$
be the incidence locus
defined in \eqref{eq:DefineIncidenceSubmanifold}
(with $Y=T$).
Let
\begin{equation}
\label{eq:DefineLocalProjToSymm}
\pi(T,\sP):
\bar\sU(\ft,\fs,\sP)\subseteqq
\tN_{\ft(\ell),\fs}(\delta)\times_{\sG_{\fs}\times S^1}\bar\Gl(\ft,\fs,\sP)
\to \Si(X^\ell,\sP)
\end{equation}
be the composition of the projection $\pi(\ft,\fs,\sP)$
defined in \eqref{eq:DefineXTubularProj} and the projection
$$
N_{\ft(\ell),\fs}(\delta)\times\Si(\ft,\fs,\sP,g_{\sP})
\cong
N_{\ft(\ell),\fs}(\delta)\times
\Si(X^\ell,\sP)
\to
\Si(X^\ell,\sP),
$$
where the identification
$\Si(\ft,\fs,\sP,g_{\sP})\cong \Si(X^\ell,\sP)$
is given in \eqref{eq:IdentifyConePointStratum}.
Define
\begin{equation}
\label{eq:LocalIncidenceComplement}
\sO(T,\sP) = \left(\pi(T,\sP)\times\id_T\right)^{-1}
\left(\Si(X^\ell,\sP)\times T\setminus\sV_\sP(\Delta,T)\right).
\end{equation}
We now have the following

\begin{lem}
\label{lem:SplicedBundleIsomOnIncidenceComplement}
The restrictions of the bundles
$(\bar r_{\beta}\times\id_T)^*\FF_T$
and
$(\id_{\sM}\times\iota_T)^*\FF^{\vir}_{\ft,\fs}$ to
$$
\sO(T):=\bigcup_{\sP} \sO(T,\sP)
$$
are isomorphic.
\end{lem}

\begin{proof}
The construction of the bundle $\fg_{\ft}$ by splicing
together the background bundle $\fg_{\ft(\ell)}$ and
the bundles $\fg_P\to S^4$ in \eqref{eq:ConstructionOfSplicedBundle}
and the definition of
$\sO(\ft,\fs,\sP)$
in \eqref{eq:LocalIncidenceComplement}
as the complement of the points
where any of the splicing points lie on $T$ give
an isomorphism of the restriction of the bundles
to each subspace $\sO(\ft,\fs,\sP)$.  These
isomorphisms agree on the overlaps
$\sO(\ft,\fs,\sP)\cap\sO(\ft,\fs,\sP')$
and hence define the desired global isomorphism.
\end{proof}
Lemma \ref{lem:SplicedBundleIsomOnIncidenceComplement} implies that
the difference,
\begin{equation}
\label{eq:DifferenceOfPClasses}
(\bar r_{\beta}\times\id_T)^*p_1(\FF_T)
-
(\id_{\sM}\times\iota_T)^*\FF^{\vir}_{\ft,\fs},
\end{equation}
lies in the image of the map $\bar\jmath_{\sO,\beta}^*$,
where $\bar\jmath_{\sO,\beta}$ is the inclusion map of pairs,
$$
\begin{CD}
\left(
\left(\bar\sM^{\vir,*}_{\ft,\fs}/S^1 \setminus \sI(\beta)\right)\times T,
\emptyset
\right)
\\
@V \bar\jmath_{\sO,\beta} VV
\\
\left(
\left(\bar\sM^{\vir,*}_{\ft,\fs}/S^1 \setminus \sI(\beta)\right)\times T,
\sO(T)
\right)
\end{CD}
$$
We can thus calculate generators of
the image of $\bar\jmath_{\sO,\beta}^*$ in the next
lemma.

\begin{lem}
\label{lem:RelHomOfIncidenceHom}
Let $T$ be a
closed, connected,
oriented,
smooth submanifold of $X$
with $[T]=\beta\in H_\bullet(X;\RR)$
and let $\sO(T)$ be as defined in Lemma \ref{lem:SplicedBundleIsomOnIncidenceComplement}.
Then the relative homology,
$$
H_k\left(
\left(\bar\sM^{\vir,*}_{\ft,\fs}/S^1 \setminus \sI(\beta)\right)\times T,
\sO(T);\ZZ
\right),
$$
is trivial for $k<4$, while for $k=4$ it is generated by the images of
disks,
$$
\psi_{\sP,P}:D^4\to \left(\bar\sU(\ft,\fs,\sP) \setminus \sI(\beta)\right)\times T
\subset \bar\sM^{\vir,*}_{\ft,\fs}/S^1\times T,
$$
for a partition $\sP$ of $N_\ell$ and $P\in\sP$, satisfying
\begin{enumerate}
\item
The image of $\psi_{\sP,P}$ is in the top stratum of $\bar\sM^{\vir,*}_{\ft,\fs}/S^1\times T$,
\item
For the map $\phi_{\sP,P}$ defined
in Lemma \ref{lem:CohomologyOfIncidenceSubmanifoldComplement}
and the map $\pi(T,\sP)$ defined in \eqref{eq:DefineLocalProjToSymm},
\begin{equation}
\label{eq:CoveringDisk}
(\pi(T,\sP)\times \id_T)\circ \psi_{\sP,P}
=
\phi_{\sP,P}.
\end{equation}
\end{enumerate}
\end{lem}

\begin{proof}
By the Whitney stratification of $\bar\sU(\ft,\fs,\sP)\times T$,
we can require any relative homology class to be represented
by a chain which intersects $(\pi(T,\sP)\times \id_T))^{-1}(\sV_\sP(\Delta,T))$
transversely (compare the proof of \cite[Lemma 6.2]{FrM}) and thus has the same
image in the relative homology group as $\psi_{\sP,P}$.
The conclusion follows by observing that the codimension of
the intersection of
$(\pi(T,\sP)\times \id_T)^{-1}\sV_\sP(\Delta,T)$ with all lower strata
is greater than four.
\end{proof}

We then compute.

\begin{lem}
\label{lem:PClassDifference}
Continue the notation of Lemmas
\ref{lem:RelHomOfIncidenceHom} and
\ref{lem:SplicedBundleIsomOnIncidenceComplement}.
For any partition $\sP$ of $N_\ell$, the following equality holds:
\begin{equation}
\label{eq:PClassDifference}
\langle (\bar r_{\beta}\times\id_T)^* p_1(\FF_T)
- (\id_{\sM}\times\iota_T)^*p_1(\FF^{\vir}_{\ft,\fs}), (\bar\jmath_{\sO,\beta})_*[\psi_{\sP,P}]
\rangle =-4|P|.
\end{equation}
\end{lem}

\begin{proof}
We begin by giving an explicit construction of
a map $\psi_{\sP,P}:D^4\to\sU(\ft,\fs,\sP)\times T$
generating the relative homology as in Lemma \ref{lem:RelHomOfIncidenceHom}.
Let $f:D^4\to X$ be a smooth embedding with $f(0)=x_T\in T$.
For $Q\in \sP$ with $Q\neq P$, fix points $y_{Q}\in X$ with $y_{Q}\notin f(D^4)$.
Define the map $\psi_{\sP,P}$,
for $v\in D^4$, by
\begin{equation}
\label{eq:ConstructTransvDisk}
\psi_{\sP,P}(v)
:=
\left(\left[ (A_0,\Phi_0),\left[ (F^{\fg}_Q(v),F^T_Q(v)),[A_Q,F^s_Q]\right]_{Q\in\sP}\right],x_T\right),
\end{equation}
where points in $\sU(\ft,\fs,\sP)$ are described in Notation \ref{notn:SplicingDataFactors} and
\begin{enumerate}
\item
$(A_0,\Phi_0)\in\tN_{\ft(\ell)}(\delta)$, and $[A_Q,F^s_Q]\in \barM^{s,\natural}_{\spl,|Q|}$, and $x_T\in T$
do not depend on $v$,
\item
For $Q \neq P$, then $F^{\fg}_{Q}(v)\in\Fr(\fg_{\ft(\ell)})|_{y_{Q}}$ and $F^T_{Q}(v)\in\Fr(TX)|_{y_{Q}}$
are fixed frames, independent of $v\in D^4$,
\item
The map $D^4\to \Fr(\fg_{\ft(\ell)})\times_X\Fr(TX)$ defined by
$v\mapsto (F^{\fg}_P(v),F^T_P(v))$
is a smooth lift of $f:D^4\to X$.
\end{enumerate}
Then the map $\phi_{\sP,P}=(\pi(T,\sP)\times \id_T)\circ\psi_{\sP,P}$ as in \eqref{eq:CoveringDisk}
is transverse to the
incidence locus $\sV_\sP(\Delta,T)$
defined in \eqref{eq:DefineIncidenceSubmanifold}
and thus is one of the generators of homology appearing
in Lemma \ref{lem:CohomologyOfIncidenceSubmanifoldComplement}.  Hence, $\psi_{\sP,P}$ satisfies the conditions
of Lemma \ref{lem:RelHomOfIncidenceHom}.

Let $\pi_\sU:\sU(\ft,\fs,\sP)\times T\to\sU(\ft,\fs,\sP)$ and $\pi_T:\sU(\ft,\fs,\sP)\times T\to T$
be the projection maps.
Define
\begin{align*}
\FF_T(x_T)&:=\FF_T|_{\sB^*_{\nu(T)}\times \{x_T\}},
\\
\FF^{\vir}_{\ft,\fs}(x_T)&:=
\bigcup_{\sP}\left( \tN_{\ft(\ell),\fs}(\delta)\times\bar\Gl(\ft,\fs,\sP)\times_{\sG_{\fs}\times S^1} \fg_{\ft(\ell)}|_{x_T}\right).
\end{align*}
Because $\pi_T\circ\psi_{\sP,P}=x_T$ from \eqref{eq:ConstructTransvDisk}, then the definitions \eqref{eq:TubNghBundle} and \eqref{eq:VirtualUniversal} yield
\begin{subequations}
\begin{align}
\label{eq:PullbackToDiskT}
\FF_T(D^4)&:=\psi_{\sP,P}^*(\bar r_\beta\times\id_T)^*\FF_T
=
(\pi_{\sU}\circ\psi_{\sP,P})^* \bar r_\beta^* \FF_T(x_T),
\\
\label{eq:PullbackToDiskVir}
\FF^{\vir}(D^4)&:=\psi_{\sP,P}^*(\id_{\sM}\times\iota_T)^*\FF^{\vir}_{\ft,\fs}
=
(\pi_{\sU}\circ\psi_{\sP,P})^*\FF^{\vir}_{\ft,\fs}(x_T).
\end{align}
\end{subequations}
Parallel translation
with respect to the connection $\hat A_0$
of the frame $F^\fg_P(v)$ appearing in \eqref{eq:ConstructTransvDisk} from
$f(v)$ to $x_T$ then defines a trivialization of $\FF^{\vir}(D^4)$.
If $\sO(T,\sP)$ is as defined in \eqref{eq:LocalIncidenceComplement}, then
$\psi_{\sP,P}(D^4\setminus\{0\})\subset\sO(T,\sP)$ and so
the isomorphism of Lemma \ref{lem:SplicedBundleIsomOnIncidenceComplement} gives an
isomorphism $\FF_T(D^4)|_{\rd D^4}\cong \FF^{\vir}(D^4)|_{\rd D^4}$.
Therefore, to compute the difference \eqref{eq:PClassDifference}, we only
need to compute the Pontrjagin class of $\FF_T(D^4)$ relative to the trivialization
over $\rd D^4$ given by the trivialization of $\FF^{\vir}(D^4)$ and
the isomorphism between the two bundles over $\rd D^4$.

Let $\fg_\ft(v)\to X$ be the
$\SO(3)$ bundle constructed in \eqref{eq:ConstructionOfSplicedBundle}
using the data $(\pi_\sU\circ\psi_{\sP,P})(v)$. The definition \eqref{eq:ConstructTransvDisk} implies
that $\fg_\ft(v)$ is constructed by attaching bundles $\fg_Q\to S^4$ to $\fg_{\ft(\ell)}$ at points
$y_{P'}$ for $Q=P'\neq P$ and at $f(v)$ for $Q=P$.

As in Section \ref{subsubsec:StdSplicingMap},
let
$\varphi_n(B(\la)) \subset S^4\setminus\{s\}$ be the image of a ball of radius $\la$ centered at the origin
under the conformal diffeomorphism $\varphi_n:\RR^4\to S^4\setminus\{s\}$ defined by inverse stereographic projection.
Consider the map $F:D^4\times B(\la)\to X$ defined by the property that $F(v,\cdot):B(\la)\to X$ is
the smooth embedding defined by the gluing data $\pi_{\sU}\circ \psi_{\sP,P}(v)$, identifying the
ball $B(\la)$ with a ball around $f(v)$ using the exponential map
around $f(v)$
used in the construction of the crude splicing map in Section \ref{subsubsec:ConstrCrudeSplice}
(given by the flattened metric and tangent frames specified by $\psi_{\sP,P}(v)$)
and
$\varphi_n^{-1}$.
Hence, if $n\in S^4$ denotes the North Pole, then $F(v,n)=f(v)$ for all $v\in D^4$.
By the Implicit Function Theorem, there is a smooth map
$h:D^4\to B(\la)$ with the property that $F(v,h(v))=x_T$.
The equality $F(v,n)=f(v)$ implies that $dF_{(0,n)}(\delta v,0)=df_0(\delta v)$.
The definition of $F(0,\cdot):B(\la)\to X$ as the embedding given by composing $\varphi_n^{-1}$
with the exponential map specified above 
implies that $dF_{(0,n)}(0,\cdot): T_nS^4\to T_{x_T}X$ is an orientation-preserving isomorphism.
The definition $F(v,h(v))=x_T$ implies that
$$
0
=dF_{(0,n)}(\delta v,dh_0(\delta v))
=df_0(\delta v)+dF_{(0,n)}(0,dh_0(\delta v)),
$$
and thus, because $T_0 D^4$ is even-dimensional, $dh_0$ is an orientation-preserving isomorphism.
Hence, $h:(D^4,D^4\setminus \{0\})\to (B(\la),B(\la)\setminus\{n\})$ is a degree-one map.

By the construction of $\fg_\ft(v)$ and the isomorphism \eqref{eq:PullbackToDiskT},
the fiber $\fg_\ft(v)|_{x_T}$ is identified with $\fg_P|_{h(v)}$, giving
\begin{equation}
\label{eq:RestrictedDiskBundle2}
\FF_T(D^4)
\cong
h^* \fg_P.
\end{equation}
Recall that we trivialized $\FF^{\vir}(D^4)$ by parallel translating $F^\fg_P(v)$ from $f(v)$
to $x_T$.  In the construction of $\fg_{\ft}(v)$, this parallel translation of $F^\fg_P(v)$
is identified with the parallel translation of the frame
$F^s\in\Fr(\fg_{P})|_s$, where $s\in S^4$ is the South Pole, from $s$ to $h(v)$.  Therefore, the trivialization
of $\FF_T(D^4)|_{\rd D^4}$ given by the trivialization of $\FF^{\vir}(D^4)$ and the isomorphism
$\FF_T(D^4)|_{\rd D^4}\cong \FF^{\vir}(D^4)|_{\rd D^4}$ equals the
pullback by $h$ of the trivialization of $\fg_P$ given by
parallel translations of the frame $F^s$.
Consequently, the relative Pontrjagin class of $h^*\fg_P$ equals
the pullback by $h$ of the  Pontrjagin class
of $\fg_P$ relative to the trivialization given by parallel translations of the frame $F^s$
and this in turn equals the absolute Pontrjagin class of $\fg_P$.
Because $h$ is a degree-one map,
the relative Pontrjagin class of the bundle \eqref{eq:PullbackToDiskT}
equals $p_1(\fg_P)=-4|P|$, as required.
\end{proof}

\begin{cor}
\label{cor:DiffOfVirtualBundles}
For any submanifold $T$ of $X$, the equality,
\begin{equation}
\label{eq:DiffOfVirtualBundles}
(\bar r_{\beta}\times\id_T)^*p_1(\FF_T)
-
(\id_{\sM}\times\iota_T)^*p_1(\FF^{\vir}_{\ft,\fs})
=
-4(\pi_X\times\id_T)^*\PD[\sV],
\end{equation}
holds in real cohomology.
\end{cor}

\begin{proof}
By Lemmas \ref{lem:PClassDifference} and
\ref{eq:IncidenceLocusThomClassEval}, the cohomology classes appearing
on the two sides of equation \eqref{eq:DiffOfVirtualBundles} agree on
the generators of the relative homology described in
Lemma \ref{lem:RelHomOfIncidenceHom}.  Hence, by the Universal
Coefficient Theorem the equality holds in real cohomology.
\end{proof}

Let $\ga_1,\dots,\ga_b\in H_1(X;\ZZ)$ be a basis for $H_1(X;\ZZ)/\Tor$
with dual basis $\bga_1^*,\dots,\bga_b^*\in H^1(X;\RR)$
satisfying $\langle \bga_i^*,\bga_j\rangle
=\delta^{ij}$ (the Kronecker delta).
Let $\bga^{J,*}_1,\dots,\bga^{J,*}_b\in H^1(M_\fs;\RR)$
be the related basis as defined in
\cite[Definition 2.21 and Corollary 2.22]{FL2a} (denoted there by $\br^*\bga^{J,*}_i$).
We now define expressions in the cohomology classes from
Section \ref{subsec:DefiningCohomClasses} which will be related to
the cocycles previously defined.

\begin{defn}
\label{defn:ExtendedCohomologyClasses}
For $x\in H_0(X;\RR)$, define
$\barmu_p(x)\in H^{4}(\bar\sM^{\vir,*}_{\ft,\fs}/S^1;\RR)$ by
\begin{equation}
\label{eq:DefineExtendedMuP4}
\barmu_p(x)
:=
-\frac{1}{4} (2\mu_\fs-\nu)^2+S^\ell(x),
\end{equation}
where $\nu$ is as in Definition \ref{defn:DefineNu}.
For $\ga\in H_1(X;\RR)$ written as $\ga=\sum_i q_i\ga_i$, where $q_i\in\RR$ and $\ga_i\in H_1(X;\ZZ)$ are the generators above, define
\begin{equation}
\label{eq:DefineExtendedMuP1}
\barmu_p(\ga)
:=-\frac{1}{2}\sum_i q_i(2\mu_\fs-\nu)\bga^{J,*}_i+ S^\ell(\ga)\in H^{3}(\bar\sM^{\vir,*}_{\ft,\fs}/S^1;\RR).
\end{equation}
For $\beta\in H_2(X;\RR)$, define
$\barmu_p(\beta)\in H^{2}(\bar\sM^{\vir,*}_{\ft,\fs}/S^1;\RR)$
by
\begin{equation}
\label{eq:DefineExtendedMuP2}
\barmu_p(\beta)
:=-\frac{1}{2}(2\mu_\fs-\nu)\langle c_1(\ft)-c_1(\fs),\beta\rangle
-2\sum_{i<j}\mu_{\fs}(\bga_i\bga_j)\langle \bga_i^*\smile\bga_j^*,\beta\rangle
+S^\ell(\beta),
\end{equation}
where $S^\ell(\beta)$ is as in Definition \ref{defn:DefnSymmBeta}.

For $Y\in H_3(X;\RR)$, define
$\barmu_p(Y) \in H^{1}(\bar\sM^{\vir,*}_{\ft,\fs}/S^1;\RR)$ by
\begin{equation}
\label{eq:DefineExtendedMuP3}
\barmu_p(Y)
:=-\sum_i \left\langle(c_1(\ft)-c_1(\fs)\smile \bga^*_i,[Y]\right\rangle \mu_{\fs}(\ga_i)+S^\ell([Y]).
\end{equation}
Define
\begin{equation}
\label{eq:DefineExtendedMuC}
\barmu_c:=-\nu\in H^2(\bar\sM^{\vir,*}_{\ft,\fs}/S^1;\RR).
\end{equation}
\end{defn}
We define
\begin{equation}
\label{eq:PartialInclusions}
\begin{aligned}
{}&
\iota_{\beta,1}:
\sM^{\vir,*}_{\ft,\fs}/S^1
\to
\bar\sM^{\vir,*}_{\ft,\fs}/S^1 \setminus \sI(\nu(\beta)),
\\
{}&
\iota_{\beta,2}:
\bar\sM^{\vir,*}_{\ft,\fs}/S^1 \setminus \sI(\nu(\beta)) \to
\bar\sM^{\vir,*}_{\ft,\fs}/S^1,
\end{aligned}
\end{equation}
to be the inclusion maps.
Standard computations (compare \cite[Lemma 4.10]{FLLevelOne}
and \cite[Corollary 4.7]{FL2b}),
Lemma \ref{lem:PontrjaginOfVirtualUniv},
Corollary \ref{cor:DiffOfVirtualBundles},
Lemma \ref{lem:PDual}, and the equality,
$$
c_1(\LL_{\fs})=\mu_\fs\times 1 +\sum_{i=1}^{b_1} \bga^{J,*}_i\times \bga^*_i,
$$
from \cite[Lemma 2.24]{FL2a}, where $\LL_{\fs}$ is as in Definition \ref{defn:LineBundleForMus},
then yield the

\begin{cor}
\label{cor:CohomologyClasses}
Assume that the \spinu structure $\ft(\ell)$ admits a splitting
$\ft(\ell)=\fs\oplus \fs\otimes L$,
so $c_1(L)=c_1(\ft)-c_1(\fs)$.
For $\om$ equal to
$x\in H_0(X;\RR)$,
$\ga\in H_1(X;\RR)$,
$\beta\in H_{2}(X;\RR)$, or  $[Y]\in H_3(X;\RR)$,
let $[c_{\om}]$ be the relative cohomology
class defined in \eqref{eq:DefineRelMuClass}
and $\iota_{\om,2}$ the inclusion defined in
\eqref{eq:PartialInclusions}.
Then
\begin{equation}
\label{eq:PullbacksOfMuAndNu}
\begin{aligned}
\bar\jmath_x^*[c_x]
&=
\iota^*_{x,2}\barmu_p(x),
\qquad
\bar\jmath_{\ga}^*[c_{\ga}]
=
\iota^*_{\ga,2}\barmu_p(\ga),
\\
\bar\jmath_\beta^*[c_{\beta}]
&=
\iota^*_{\beta,2}\barmu_p(\beta),
\qquad
\bar\jmath_Y^*[c_{Y}]
=
\iota^*_{Y,2}\barmu_p(Y),
\end{aligned}
\end{equation}
where $\bar\jmath_x$, $\bar\jmath_\beta$, $\bar\jmath_{\ga_i}$, and $\bar\jmath_Y$
are the inclusion maps of pairs defined in
equation \eqref{eq:DefinePairsInclusionForBeta}.
\end{cor}

The analogue below of Corollary \ref{cor:CohomologyClasses}
for $[c_\sW]$ follows immediately from the observation that $\bar\sW$ is
the zero locus of the complex line bundle $\LL_\nu$ defined
in \eqref{eq:LineBundleForS1ZAction} and the computation
in \cite[Lemma 4.2]{FLLevelOne}.

\begin{cor}
\label{cor:CohomClassC1}
If
$$
\bar\jmath_\nu:
\left(
    \bar\sM^{\vir,*}_{\ft,\fs}/S^1\setminus\sI(\nu(x)),\emptyset
\right)
\to
\left(
    \bar\sM^{\vir,*}_{\ft,\fs}/S^1\setminus\sI(\nu(x)),\bar\sM^{\vir,*}_{\ft,\fs}/S^1\setminus\bar\sW
\right)
$$
is the inclusion map of the pair, then
\begin{equation}
\label{eq:PullbacksOfNu}
\bar\jmath_\nu^*[c_\sW]=-\iota_{x,2}^*\nu,
\end{equation}
where $[c_\sW]$ is defined in \eqref{eq:DefinePreCocycleForWClass}
and $\iota_{x,2}$ is defined in \eqref{eq:PartialInclusions}.
\end{cor}

\section{Relative Euler class of the obstruction pseudo-bundle}
\label{sec:RelEulerClassObstr}
We now use the framework of relative Euler classes to
construct a cocycle to represent the geometric representative
defined by the zero locus of the obstruction section,
$\bchi^{-1}(0)\subset \bar\sM^{\vir,*}_{\ft,\fs}/S^1$,
appearing in Hypothesis \ref{hyp:Gluing}.
The restriction of the obstruction pseudo-bundle to the top stratum of $\bar\sM^{\vir}_{\ft,\fs}/S^1$
is a bundle which we will refer to as the \emph{obstruction bundle}.
We begin by showing that the Euler class of this obstruction bundle
is the restriction to the top stratum of $\bar\sM^{\vir,*}_{\ft,\fs}/S^1$ of a cohomology class
\begin{equation}
\label{eq:ExtendedEulerClass}
\bar e_I\smile e(\bar\Upsilon^s_{\ft,\fs}/S^1)\in H^\bullet(\bar\sM^{\vir,*}_{\ft,\fs}/S^1;\RR),
\end{equation}
given by an expression in the cohomology classes appearing in Section \ref{subsec:DefiningCohomClasses}.
We use the relative Euler class of the obstruction bundle to construct a cocycle corresponding to
the geometric representative $\bchi^{-1}(0)$ in the sense of Lemma \ref{lem:DefineSingularCocycle}
and prove that that cocycle extends over the subspace of $\bar\sM^{\vir,*}_{\ft,\fs}/S^1$ given
by the complement of a neighborhood of the intersection of $\bar\bchi^{-1}(0)$ with the lower strata
of $\bar\sM^{\vir,*}_{\ft,\fs}/S^1$.  Finally, we prove that this extension of the cocycle given
by the relative Euler class is mapped to the cohomology class \eqref{eq:ExtendedEulerClass}
in the exact sequence of the pair.

\subsection{Euler class of the Seiberg--Witten component of the obstruction pseudo-bundle}
\label{subsec:SWObstrEuler}
The computation of the Euler class of the Seiberg--Witten, or background,
component of the obstruction is identical to that given in
\cite[Lemma 4.11]{FLLevelOne}.

\begin{lem}
\label{lem:SWObstruction}
Let $r_\Xi$ be the complex rank of the background obstruction bundle,
$\bar\Upsilon^s_{\ft,\fs}\to \bar\sM^{\vir,*}_{\ft,\fs}$.  Then
\begin{equation}
e(\bar\Upsilon^s_{\ft,\fs}/S^1)=\iota^*(-\nu)^{r_\Xi}.
\end{equation}
\end{lem}

Because the background obstruction bundle is defined over
all of $\bar\sM^{\vir,*}_{\ft,\fs}$, there is no need to
define an extension of $e(\bar\Upsilon^s_{\ft,\fs}/S^1)$.

\subsection{Local Euler class  of the instanton component of the obstruction bundle}
\label{subsec:InstObstrEulerLocal}
We now give a description of the Euler class
of the bundle defined by the
restriction of the instanton component of the
obstruction pseudo-bundle, $\Upsilon^i_{\ft,\fs}/S^1\to\sM^{\vir}_{\ft,\fs}/S^1$,
to the top stratum, $\sU(\ft,\fs,\sP)$, of $\bar\sU(\ft,\fs,\sP)$.

This description requires the introduction of an additional cohomology class.
Let $\pi_{\fs,i}:M_{\fs}\times X^\ell\to M_{\fs}\times X$ be
defined by the identity map on $M_{\fs}$ and projection map from $X^\ell$
onto the $i$-th factor of $X$.
Let $\LL_\fs\to M_\fs\times X$ be the line bundle in Definition \ref{defn:LineBundleForMus}.
The Chern classes of the bundle
$$
\bigoplus_{i=1}^\ell \pi_{\fs,i}^*\LL_{\fs}
\to
M_{\fs}\times X^\ell
$$
are invariant under the obvious action of the symmetric group
$\fS_\ell$ and thus define cohomology classes,
\begin{equation}
\label{eq:ChernClassesOfSymSWUniv}
c_{\fs,\ell,j}\in H^{2j}(M_{\fs}\times \Sym^\ell(X);\RR).
\end{equation}
In addition, we let
$$
\iota_{\sP}: \Si(X^\ell,\sP)\to\Sym^\ell(X)
\quad\text{and}\quad
\tilde\iota_{\sP}:\Delta^\circ(X^\ell,\sP)\to X^\ell
$$
denote the inclusion maps.
Finally, we let
\begin{equation}
\label{eq:SpinuCharacter}
c(\ft)\in H^\bullet(X;\RR)
\end{equation}
denote any real characteristic class of the bundle
$\Fr_{\CCl(T^*X)}(V)\to X$; those characteristic classes include $c_1(\ft)$, $p_1(\ft)$, $p_1(X)$,
and $e(X)$.

If $\om\in H^\bullet(\cdot;\ZZ)$, we refer to the image
of $\om$ in $H^\bullet(\cdot;\RR)$ under the functor
$H^\bullet(\cdot;\ZZ)\to H^\bullet(\cdot;\RR)$ defined
by the homomorphism $\ZZ\to \RR$
as the image of $\om$ in real cohomology.

\begin{prop}
\label{prop:LocalInstantonEuler}
Let $\sU(\ft,\fs,\sP)\subset\bar\sU(\ft,\fs,\sP)$ be
the top stratum.  The image in real cohomology of the
Euler class of the restriction of the obstruction bundle
$\Upsilon^i_{\ft,\fs}/S^1$ to $\sU(\ft,\fs,\sP)$ is given
by a polynomial in the cohomology classes
\begin{equation}
\nu, \quad
\pi_X^*\iota_{\sP}^*S^\ell(c(\ft)), \quad\text{and}\quad
\pi_{X,\fs}^*\iota_{\sP}^*c_{\fs,\ell,j},
\end{equation}
with coefficients depending only on the partition
$\ell=|P_1|+\dots + |P_r|$, where $\sP=\{P_1,\dots,P_r\}$.
\end{prop}

To prove Proposition
\ref{prop:LocalInstantonEuler}, we first show that the restriction
of the instanton obstruction bundle to  $\sU(\ft,\fs,\sP)$
admits a direct sum decomposition in Lemma \ref{lem:LocalInstantonPullbackToCover}
and then in Lemmas \ref{lem:RelateInstantonObstructionFactors}
and \ref{lem:XFactorOfLocalInstantonEuler}
prove that the components of this direct sum satisfy the
conclusions of Proposition \ref{prop:LocalInstantonEuler}.

Let $\Gl(\ft,\fs,\sP)\subset\bar\Gl(\ft,\fs,\sP)$ denote the top stratum,
$$
\Gl(\ft,\fs,\sP)
=
\Fr(\ft,\fs,\sP)\times_{G(\sP)}\prod_{P\in\sP} M^{s,\natural}_{\spl,|P|}(\delta_P).
$$
Because $\tilde N_{\ft(\ell),\fs}(\delta)$ retracts $\sG_\fs$-equivariantly onto
$\tM_{\fs}$,
the space $\sU(\ft,\fs,\sP)$ retracts onto
the intersection of $\bar\sU(\ft,\fs,\sP)$ with the subspace
\begin{equation}
\label{eq:SWTopStratumSubspace}
\tM_{\fs}\times_{\sG_{\fs}\times S^1}\Gl(\ft,\fs,\sP).
\end{equation}
We compute the pullback of the restriction of the
Euler class of the instanton obstruction to
the subspace \eqref{eq:SWTopStratumSubspace}
by the map,
\begin{equation}
\label{eq:LocalCoveringMap}
\tM_{\fs}\times_{\sG_{\fs}\times S^1}\widetilde{\Gl}(\ft,\fs,\sP)
\to
\tM_{\fs}\times_{\sG_{\fs}\times S^1}\Gl(\ft,\fs,\sP),
\end{equation}
defined by the projection
$$
\widetilde{\Gl}(\ft,\fs,\sP)
:=
\Fr(\ft,\fs,\sP)\times_{\tG(\sP)}\prod_{P\in\sP} M^{s,\natural}_{\spl,|P|}(\delta_P)
\to
\Gl(\ft,\fs,\sP),
$$
where the group $\tG(\sP)\le G(\sP)$ is defined in
\eqref{eq:DefineGluingDataBundleStructureGroup}
by omitting the symmetric group factor in $G(\sP)$.
Observe that the preceding map appears in the diagram,
$$
\begin{CD}
\tilde M_{\fs}\times_{\sG_{\fs}\times S^1}
\widetilde{\Gl}(\ft,\fs,\sP)
@>>>
\tilde M_{\fs}\times_{\sG_{\fs}\times S^1}
\Gl(\ft,\fs,\sP)
\\
@V \tilde\pi_{X,\fs} VV @V \pi_{X,\fs} VV
\\
M_{\fs}\times \Delta^\circ(X^\ell,\sP) @>>> M_{\fs}\times \Si(X^\ell,\sP)
\end{CD}
$$
For each $P\in\sP$, we define
\begin{equation}
\label{eq:DefineSingleInstantonSplicingDataBundle}
\Gl(\ft,\fs,P):=\Fr(\fg_{\ft(\ell)})\times_X\Fr(TX)\times_{\SO(3)\times\SO(4)}
M^{s,\natural}_{\spl,|P|}(\delta_P),
\end{equation}
and observe that
the projection $\widetilde{\Gl}(\ft,\fs,\sP)\to\Delta^\circ(X^\ell,\sP)$
fits into a diagram
\begin{equation}
\label{eq:DefineProjectionOfGluingParamBundleOntoFactor}
\begin{CD}
\tM_{\fs}\times_{\sG_{\fs}\times S^1}
\widetilde{\Gl}(\ft,\fs,\sP) @> \pi_{\widetilde{\Gl},P} >>
\tM_{\fs}\times_{\sG_{\fs}\times S^1}
\Gl(\ft,\fs,P)
\\
@VVV @VVV
\\
M_{\fs}\times \Delta^\circ(X^\ell,\sP) @> \id_{M_{\fs}}\times \pi_P >> M_{\fs}\times X
\end{CD}
\end{equation}
We then have the following

\begin{lem}
\label{lem:LocalInstantonPullbackToCover}
Let $\ft(\ell)=(\rho,V_\ell)$ be the \spinu structure defined in
\eqref{eq:DefineLowerChargeSpinuStr}.
We abbreviate $\fg_{\ell}=\fg_{\ft(\ell)}$ and $V_\ell=V_{\ft(\ell)}$.
The pullback of the restriction of
$\Upsilon^i_{\ft,\fs}/S^1$ to $\sU(\ft,\fs,\sP)$
by the covering map \eqref{eq:LocalCoveringMap}
splits into a direct sum,
$$
\bigoplus_{P\in\sP}\pi_{\widetilde{\Gl},P}^*\Upsilon^i_{\ft,\fs}(P)/S^1,
$$
where for $\Gl(\ft,\fs,P)$ defined in \eqref{eq:DefineSingleInstantonSplicingDataBundle},   the bundle
\begin{equation}
\label{eq:FactorOfInstantonEuler}
\Upsilon^i_{\ft,\fs}(P)/S^1
\to
\tM_{\fs}\times_{\sG_{\fs}\times S^1}
\Gl(\ft,\fs,P)
\end{equation}
is defined by
$$
\begin{CD}
\tM_{\fs}\times_{\sG_{\fs}\times S^1}
\Fr_{\CCl(T^*X)}(V_\ell)\times_{\Spinu(4)}\Upsilon^i_{\spl,|P|}
\\
@VVV
\\
\tM_{\fs}\times_{\sG_{\fs}\times S^1}
\Gl(\ft,\fs,P)
\end{CD}
$$
and the map $\pi_{\widetilde{\Gl},P}$ is defined in the diagram
\eqref{eq:DefineProjectionOfGluingParamBundleOntoFactor}.
\end{lem}

The $S^1$ action in \eqref{eq:FactorOfInstantonEuler} on the space
\begin{equation}
\label{eq:FactorOfInstantonEuler1}
\tM_{\fs}\times_{\sG_{\fs}}
\Fr_{\CCl(T^*X)}(V_\ell)\times_{\Spinu(4)}
\Upsilon^i_{\spl,|P|}
\end{equation}
is defined by the action of $S^1$ on
the infinite-dimensional obstruction space described following
\eqref{eq:InfiniteDimObstr}.
That is, for $(A_0,\Phi_0)\in\tM_{\fs}$, and $\tilde F_u\in\Fr_{\CCl(T^*X)}(V_\ell)|_x$, and
$[A,F_s,\Psi]\in\Upsilon^i_{\spl,|P|}$, the action is given by
\begin{equation}
\label{eq:S1ObstructionAction,SingleInstanton1}
\left(
    e^{i\theta},
    \left[ (A_0,\Phi_0), \left( \tilde F_u,[A,F_s,\Psi]\right) \right]
\right)
\mapsto
\left[ (A_0,\Phi_0 e^{i\theta}), \left( \tilde F_u,[A,F_s,\Psi e^{i\theta}]\right) \right].
\end{equation}
Recall that $\varrho_{\fs}(e^{i\theta})$ is defined in \eqref{eq:RelateS1Actions}.
Because the action of scalar multiplication on $V_\ell$ differs from the action of
$\varrho_\fs(e^{i\theta})$
by
an element of the stabilizer of $(A_0,\Phi_0)\in\tM_{\fs}$, the action \eqref{eq:S1ObstructionAction,SingleInstanton1}
is the same as the action
\begin{equation}
\label{eq:S1ObstructionAction,SingleInstanton2}
\left(
    e^{i\theta},
    \left[ (A_0,\Phi_0), \left( \tilde F_u,[A,F_s,\Psi]\right) \right]
\right)
\mapsto
\left[ (A_0,\Phi_0 ), \left( \varrho_{\fs}(e^{-i\theta})\tilde F_u,[A,F_s,\Psi e^{i\theta}]\right) \right].
\end{equation}
Lemma \ref{lem:LocalInstantonPullbackToCover}
implies that to compute the restriction
of the instanton obstruction bundle to the subspace
\eqref{eq:SWTopStratumSubspace},
it suffices to characterize the Euler class of
the bundle \eqref{eq:FactorOfInstantonEuler}.
The base of the bundle $\Upsilon^i_{\ft,\fs}(P)/S^1$
admits a product decomposition,
$$
\tM_{\fs}\times_{\sG_{\fs}\times S^1}\Gl(\ft,\fs,P)
\cong
M_{\fs}\times S^1\backslash\Gl(\ft,\fs,P).
$$
Define
\begin{equation}
\label{eq:XFactorObstruction}
\Upsilon^i_{X,P}/S^1
:=
M_{\fs}\times
S^1\backslash\left(\Fr_{\CCl(T^*X)}(V_\ell)\times_{\Spinu(4)}\Upsilon^i_{\spl,|P|}\right)
\to
M_{\fs}\times S^1\backslash\Gl(\ft,\fs,P).
\end{equation}
As in the discussion of the $S^1$ action on
\eqref{eq:FactorOfInstantonEuler1}, the $S^1$ action in \eqref{eq:XFactorObstruction}
is given by the diagonal action of
$\varrho_{\fs}(e^{-i\theta})$ on $\Fr_{\CCl(T^*X)}(V_\ell)$
and scalar multiplication on the fibers of
$\Upsilon^i_{\spl,|P|}$.

\begin{lem}
\label{lem:RelateInstantonObstructionFactors}
Let $\LL_{\fs}\to M_{\fs}\times X$ be the restriction of the
universal Seiberg--Witten
line bundle defined in \eqref{eq:DefineSWUniversal}
and $\pi_{X,\fs}:M_{\fs}\times S^1\backslash\Gl(\ft,\fs,P)\to M_{\fs}\times X$
be the projection map.  Then
\begin{equation}
\label{eq:RelateInstantonObstructionFactors}
\Upsilon^i_{\ft,\fs}(P)/S^1
\cong
\pi_{X,\fs}^*\LL_{\fs}^*\otimes\Upsilon^i_{X,P}/S^1.
\end{equation}
\end{lem}

\begin{proof}
By  \cite[Lemma 3.27]{FL2a},
the tensor product in \eqref{eq:RelateInstantonObstructionFactors}
can be described by the fibered product,
\begin{equation}
\label{eq:FiberProductToTensor}
\begin{aligned}
{}&\left(
\left(
\tM_{\fs}\times_{\sG_{\fs}} X\times S^1
\right)\times_{M_{\fs}\times X}\Upsilon^i_{X,P}/S^1
\right)/S^1
\\
{}&\quad
=
\left(
\left(
\tM_{\fs}\times_{\sG_{\fs}} X\times S^1
\right)
\times_{M_{\fs}\times X}
\left(
    M_{\fs}\times S^1\backslash\Fr_{\CCl(T^*X)}(V_\ell)\times_{\Spinu(4)}\Upsilon^i_{\spl,|P|}
\right)
\right)/S^1,
\end{aligned}
\end{equation}
where the final $S^1$ acts diagonally on the factor of $S^1$ in $\tM_{\fs}\times_{\sG_{\fs}} X\times S^1$
 and
on the fiber of $\Upsilon^i_{X,P}/S^1$.  The fibered product \eqref{eq:FiberProductToTensor} admits
a simplification as
\begin{equation}
\label{eq:FiberProductToTensor2}
\begin{aligned}
{}&
\left(
\left(
\tM_{\fs}\times_{\sG_{\fs}} X\times S^1
\right)
\times_{M_{\fs}\times X}
\left(
    M_{\fs}\times S^1\backslash\Fr_{\CCl(T^*X)}(V_\ell)\times_{\Spinu(4)}\Upsilon^i_{\spl,|P|}
\right)
\right)/S^1
\\
{}&\quad\cong
\tM_{\fs}\times_{\sG_{\fs}} S^1\backslash\left(\Fr_{\CCl(T^*X)}(V_\ell)\times_{\Spinu(4)}\Upsilon^i_{\spl,|P|}\right),
\end{aligned}
\end{equation}
given by the map (using the same notation as in \eqref{eq:S1ObstructionAction,SingleInstanton2}),
$$
\left(
\left( (A_0,\Phi_0),x,e^{i\theta}\right),\left[[A_0,\Phi_0], \left[\tilde F_u,[A,F_s,\Psi]\right]\right]
\right)
\mapsto
\left(
(A_0,\Phi_0),\left[\varrho_\fs(e^{-i\theta})\tilde F_u,[A,F_s,\Psi]\right]
\right).
$$
This map is $S^1$-equivariant with respect to the $S^1$ action on the domain given by the diagonal action
on the factor $S^1$ and scalar multiplication on the fibers of $\Upsilon^i_{\spl,|P|}$, if $S^1$
acts on the image by $\varrho_{\fs}(e^{-i\theta})$ on $\Fr_{\CCl(T^*X)}(V_\ell)$ and by
scalar multiplication on the fibers of $\Upsilon^i_{\spl,|P|}$.  This is the same $S^1$ action as
that appearing in \eqref{eq:S1ObstructionAction,SingleInstanton2}.
The conclusion of the lemma then follows by observing that the presentation for the $S^1$ group action in
\eqref{eq:S1ObstructionAction,SingleInstanton2} identifies the bundle on the right-hand-side of
\eqref{eq:FiberProductToTensor2} with $\Upsilon^i_{\ft,\fs}(P)/S^1$ as required.
\end{proof}

Next, we characterize the (rational)
Euler class of $\Upsilon^i_{X,P}/S^1$.  To this end, we
define the cohomology class,
\begin{equation}
\label{eq:XFactorS1Chern}
\nu_{\Gl}\in H^2(M_{\fs}\times S^1\backslash\Gl(\ft,\fs,P);\RR),
\end{equation}
(denoted by $\nu_\ft$ in \cite[Definition 5.22]{FLLevelOne})
to be the first Chern class of the $S^1$ bundle,
$$
M_{\fs}\times \Gl(\ft,\fs,P)
\to
M_{\fs}\times S^1\backslash\Gl(\ft,\fs,P).
$$
The equality,
\begin{equation}
\label{eq:RelateNus}
\nu=\nu_{\Gl}+2\pi_{X,\fs}^*c_1(\LL_{\fs}),
\end{equation}
follows from \cite[Lemma 5.23]{FLLevelOne}.
The advantage of the bundle \eqref{eq:XFactorS1Chern}
over the $S^1$ bundle defining $\nu$
in \eqref{eq:LineBundleForS1ZAction}
is that the former pulls back
by the projection $M_{\fs}\times S^1\backslash\Gl(\ft,\fs,P)
\to S^1\backslash\Gl(\ft,\fs,P)$.

\begin{lem}
\label{lem:XFactorOfLocalInstantonEuler}
The images of the Chern classes of $\Upsilon^i_{X,P}/S^1$
in real cohomology are given by polynomials in
 $\pi_X^*c(\ft)$ and
$\nu_{\Gl}$ with coefficients depending only on $|P|$.
\end{lem}

\begin{proof}
Observe that the bundle $\Upsilon^i_{X,P}/S^1$ is the pullback by the projection
$$
M_{\fs}\times S^1\backslash\Gl(\ft,\fs,P)
\to
S^1\backslash\Gl(\ft,\fs,P)
$$
of the bundle
\begin{equation}
\label{eq:SimplifyBundle}
S^1\backslash\Fr_{\CCl(T^*X)}(V_\ell)\times_{\Spinu(4)}\Upsilon^i_{\spl,|P|}
\to
S^1\backslash\Gl(\ft,\fs,\sP).
\end{equation}
Thus, it suffices to prove that the Chern class of the
bundle \eqref{eq:SimplifyBundle} satisfies the conclusion of the lemma.

The kernel of the projection
\begin{equation}
\label{eq:GroupProjection}
(\Ad^u_{\SO(3)},\Ad^u_{\SO(4)}):
\Spinu(4)\to\SO(3)\times\SO(4)
\end{equation}
is the central $S^1$ in $\Spinu(4)$.  By the identity
(see \cite[Equation (3.14)]{FLLevelOne})
$$
\Fr_{\CCl(T^*X)}(V_\ell)/S^1=\Fr(\fg_\ell)\times_X\Fr(TX),
$$
(where the $S^1$ is the kernel of the homomorphism \eqref{eq:GroupProjection}),
we can rewrite the space $\Gl(\ft,\fs,P)$ as
$$
\Gl(\ft,\fs,P)
=
\Fr_{\CCl(T^*X)}(V_\ell)\times_{\Spinu(4)}M^{s,\natural}_{\spl,|P|}(\delta_P),
$$
where $\Spinu(4)$ acts on  $M^{s,\natural}_{\spl,|P|}(\delta_P)$
via the projection \eqref{eq:GroupProjection}
and the action of $\SO(3)\times\SO(4)$ on
$M^{s,\natural}_{\spl,|P|}(\delta_P)$.
Let $\tilde\iota_u:\Fr_{\CCl(T^*X)}(V_\ell)\to\ESpinu(4)$ be the classifying
map appearing in the diagram,
$$
\begin{CD}
\Fr_{\CCl(T^*X)}(V_\ell) @> \tilde\iota_u >> \ESpinu(4)
\\
@VVV @VVV
\\
X @> \iota_u >> \BSpinu(4)
\end{CD}
$$
Recall from \cite[Section 4]{AtiyahJones},
that the map $[A]\to [D_A]$,
where $D_A$ is the Dirac operator defined
by the connection $A$, defines a continuous map,
$$
f_\ka: \sB^s_\ka\to\sF_\ka,
$$
where $\sF_\ka$ is the space of Fredholm operators of index $\ka$.
The space $\sF_\ka$ is homotopic to $\BU$, the classifying space
for the stable unitary group $\U:=\lim_{n\to\infty}\U(n)$
(see \cite[Equation (4.7]{AtiyahJones}).
Although there is no index bundle defined for the Dirac operators
parameterized by $\sB^s_\ka$, because $\sB^s_\ka$ is not compact, the
restrictions of the pullback of the universal Chern
classes in $\BU$ by $f_\ka$
to the subspace $\sK^s_\ka\subset\sB^s_\ka$,
where the
kernel of the Dirac operator vanishes,
are equal to the Chern class of the
vector bundle defined by
$\ind(\bD^*_\ka)$
in \eqref{eq:DefineInstantonDiracBundle}
on this subspace.
By Property
\eqref{item:ExistenceOfSplicedEndsIndex5}
in the conclusion of Theorem \ref{thm:ExistenceOfSplicedEndsIndex},
the Chern classes of the bundle,
$$
\Upsilon^i_{\spl,|P|}
\to
M^{s,\natural}_{\spl,|P|},
$$
are given by the
pullbacks of the universal Chern classes
in $H^\bullet(\BU)$ by the map $f_{|P|}$.
By the $\Spinu(4)$-equivariance of the Dirac operator, the
map $f_\ka$ extends to a map
$$
F_\ka: \ESpinu(4)\times_{\Spinu(4)}\sB^s_\ka \to\sF_\ka,
$$
such that the restriction of $F_\ka$ to
$$
\ESpinu(4)\times_{\Spinu(4)}\sK^s_\ka
$$
is the composition of a classifying map for the index bundle
with the inclusion $\BU(\ka)\to\BU$.
We conclude that the classifying map for the bundle
$$
\Upsilon^i_{X,P}
=
\Fr_{\CCl(T^*X)}(V_\ell)\times_{\Spinu(4)}\Upsilon^i_{\spl,|P|}
\to
\Fr_{\CCl(T^*X)}(V_\ell)\times_{\Spinu(4)}M^{s,\natural}_{\spl,|P|}(\delta_P)
$$
factors through the composition
\begin{equation}
\label{eq:ClassifyingMapFactoring}
\begin{CD}
\Fr_{\CCl(T^*X)}(V_\ell)\times_{\Spinu(4)}M^{s,\natural}_{\spl,|P|}(\delta_P)
@>\tilde\iota_u\times\iota_M >>
\ESpinu(4)\times_{\Spinu(4)}\sB^s_{|P|}
@> F_{|P|}>>
\sF_{|P|},
\end{CD}
\end{equation}
where $\iota_M:M^{s,\natural}_{\spl,|P|}(\delta_P)\to\sB^s_{|P|}$
is the inclusion map.  The bundle map
$\tilde\iota_u\times\iota_M$
is $S^1$-equivariant and thus
descends to the $S^1$ quotients,
$$
\tilde\iota_u\times\iota_M:
S^1\left\backslash
\left(
\Fr_{\CCl(T^*X)}(V_\ell)\times_{\Spinu(4)}M^{s,\natural}_{\spl,|P|}(\delta_P)
\right)\right.
\to
S^1\left\backslash
\left(\ESpinu(4)\times_{\Spinu(4)}\sB^s_{|P|}\right)\right. .
$$
Because the rational cohomology of $\sB^s_{|P|}$
is trivial (see \cite[Lemma 5.1.14]{DK}), the rational
cohomology of the space
$$
S^1\left\backslash
\left(\ESpinu(4)\times_{\Spinu(4)}\sB^s_{|P|}\right)\right.
$$
is generated by the first Chern class of the $S^1$ action
and cohomology classes pulled back by the projection
$$
S^1\left\backslash
\left(\ESpinu(4)\times_{\Spinu(4)}\sB^s_{|P|}\right)\right.
\to
\BSpinu(4).
$$
Under the map $\tilde\iota_u\times\iota_M$,
the first Chern class of the $S^1$ action pulls back to
$\nu_{\Gl}$, while the cohomology classes pulled back
from $\BSpinu(4)$ pull back to characteristic classes of
$\Fr_{\CCl(T^*X)}(V_\ell)$, which we are denoting by $c(\ft)$.
The conclusion of the lemma then follows from the
factorization
of the classifying map for the bundle
\eqref{eq:SimplifyBundle} given in
\eqref{eq:ClassifyingMapFactoring}.
\end{proof}

\begin{rmk}
Unlike the case $\ell=1$, the action of $\sG_{\fs}$ on
$S^1\backslash\Gl(\ft,\fs,\sP)$ is not trivial because there
can be more than one frame in the definition of $\Gl(\ft,\fs,\sP)$.
That is, if $F_1$ and $F_2$ are frames over separate points
$x_1$ and $x_2$
in $X$, the action of $S^1$ identifies $(F_1,F_2)$ with
$(e^{i\theta}F_1,e^{i\theta}F_2)$.
But there are gauge transformations, $u\in\sG_\fs$,
with $u(x_1)\neq u(x_2)$ and for such $u\in\sG_\fs$, the pair
$(uF_1,uF_2)$ would not be identified with $(F_1,F_2)$ by the $S^1$ action.
\end{rmk}

\begin{proof}[Proof of Proposition \ref{prop:LocalInstantonEuler}]
By the Splitting Principle, we can assume that there are line bundles
$N_{P,i}$ such that $\Upsilon^i_{X,P}/S^1\cong\oplus_{i\in P}N_{P,i}$.
By Lemmas \ref{lem:LocalInstantonPullbackToCover},
\ref{lem:RelateInstantonObstructionFactors}, and
\ref{lem:XFactorOfLocalInstantonEuler}, we see that the pullback by the map
\eqref{eq:LocalCoveringMap} of the restriction of $\Upsilon^i_{\ft,\fs}/S^1$
to the space \eqref{eq:SWTopStratumSubspace} is given by
\begin{align*}
\bigoplus_{P\in\sP}
\pi_{\widetilde{\Gl},P}^* \Upsilon^i_{\ft,\fs}(P)/S^1
{}&=
\bigoplus_{P\in\sP}
\pi_{\widetilde{\Gl},P}^*  \left(\pi_{X,\fs}^*\LL_{\fs}^*\otimes\Upsilon^i_{X,P}/S^1\right)
\\
{}&=
\bigoplus_{P\in\sP}\bigoplus_{i\in P}
\pi_{\widetilde{\Gl},P}^*  \left(\pi_{X,\fs}^*\LL_{\fs}^*\otimes N_{P,i}\right).
\end{align*}
The Euler class of the bundle
\[
\bigoplus_{P\in\sP}
\pi_{\widetilde{\Gl},P}^* \Upsilon^i_{\ft,\fs}(P)/S^1
\]
on the left-hand side of the preceding sequence of equalities is thus
$$
\prod_{P\in\sP}\prod_{i\in P}
\pi_{\widetilde{\Gl},P}^*\left( -\pi_{X,\fs}^*c_1(\LL_{\fs}) + c_1(N_{P,i})\right).
$$
Observe that $\pi_{X,\fs}\circ\pi_{\widetilde{\Gl},P}=\pi_{\fs,i}\tilde\pi_{X,\fs}$
(where $\pi_{\fs,i}$ is defined prior to \eqref{eq:ChernClassesOfSymSWUniv}).
Hence,
the preceding cohomology class can be expressed in terms of  symmetric
polynomials in
$\tilde\pi_{X,\fs}^*\pi_{\fs,i}^*c_1(\LL_{\fs})$
and
$\pi_{\widetilde{\Gl},P}^*c_1(N_{P,i})$.
Symmetric polynomials in $\tilde\pi_{X,\fs}^*\pi_{\fs,i}^*c_1(\LL_{\fs})$
are given by the Chern classes $c_{\fs,\ell,j}$ defined in
\eqref{eq:ChernClassesOfSymSWUniv}, while symmetric polynomials
in $\pi_{\widetilde{\Gl},P}^*c_1(N_{P,i})$ are given by the Chern
classes of the bundle
$$
\bigoplus_{P\in\sP}\pi_{\widetilde{\Gl},P}^*\Upsilon^i_{X,P}/S^1.
$$
By Lemma \ref{lem:XFactorOfLocalInstantonEuler}, the Chern classes of
the preceding bundle are given by polynomials in $\nu$ and $\pi_X^*c(\ft)$,
with coefficients depending only on the partition
$\ell=|P_1|+\dots + |P_r|$.
\end{proof}

\subsection{Global Euler class of the instanton component of the obstruction bundle}
\label{subsec:InstObstrEulerGlobal}
We now piece together the local computations
of Proposition \ref{prop:LocalInstantonEuler}
to give a global characterization of the instanton
obstruction bundle, $\Upsilon^i_{\ft,\fs}/S^1\to\sM^{\vir}_{\ft,\fs}/S^1$.

\begin{prop}
\label{prop:GlobalInstantonEulerClass}
The Euler class of the instanton obstruction
bundle,
$\Upsilon^i_{\ft,\fs}/S^1\to\sM^{\vir}_{\ft,\fs}/S^1$,
is given by an element in $H^{2\ell}(\sM^{\vir}_{\ft,\fs}/S^1;\RR)$ that is a polynomial
in $\iota^*\nu$, and
$\iota^*\pi_X^*S^\ell(c(\ft))$, and
$\iota^*\pi_{X,\fs}^*c_{\fs,\ell,i}$,
with coefficients which are independent of $X$.
\end{prop}

\begin{proof}
The construction of the bundle
$\Upsilon^i_{\spl,\ka}\to M^{s,\natural}_{\spl,\ka}(\delta)$
and the $\Spinu(4)$ equivariance of the diagram
\eqref{eq:OverlapObstructionCommuting} imply that
the particular classifying maps of the local instanton obstruction
bundles which admit the factorization
\eqref{eq:ClassifyingMapFactoring}
can be chosen to be equal on the overlap
given in Lemma \ref{lem:InstantonObstrCommDiagr}.
Hence, the local equality of the cohomology classes given
in Proposition \ref{prop:LocalInstantonEuler} is an equality
of cocycles and thus an equality of global cohomology classes.
\end{proof}

\begin{defn}
Let
\begin{equation}
\label{eq:ExtendedInstEuler}
\bar e_I\in H^{2\ell}(\bar\sM^{\vir,*}_{\ft,\fs}/S^1;\RR)
\end{equation}
be the extension of $e(\Upsilon^i_{\ft,\fs}/S^1)\in H^{2\ell}(\sM^{\vir,*}_{\ft,\fs}/S^1;\RR)$
obtained by omitting the pullback $\iota^*$ from the expression in
Proposition \ref{prop:GlobalInstantonEulerClass}.
\end{defn}

\begin{rmk}
Note that the vector bundle $\Upsilon^i_{\ft,\fs}/S^1$ does
not extend from $\sM^{\vir,*}_{\ft,\fs}/S^1$ to
$\bar\sM^{\vir,*}_{\ft,\fs}/S^1$.  Indeed, the
cohomology class $u_i$ from \cite[Definition 8.3.16]{DK}
can be seen as an obstruction to such an extension.
\end{rmk}

The following lemma describing the extension $\bar e_I$
will be used in Section \ref{sec:Duality} to relate
intersection numbers with a pairing of cohomology and homology
classes.

\begin{lem}
\label{lem:ExtendedObstructionClassOnLowerStrata}
Let $i_\Si:\Si\subset \bar\sM^{\vir,*}_{\ft,\fs}/S^1$ be
the inclusion map of any
smooth stratum.
Let $e_I(\Si)$ be the Euler class of the vector bundle
defined by the restriction of the pseudo-vector bundle
$\bar\Upsilon^i_{\ft,\fs}$ to $\Si$.
For $\bar e_I$ as defined in \eqref{eq:ExtendedInstEuler}, then
$$
i_\Si^*(\bar e_I\smile e(\bar\Upsilon^s_{\ft,\fs}/S^1))
=e_I(\Si)_\RR\smile \om_i(\Si)\smile i_\Si^*(e(\bar\Upsilon^s_{\ft,\fs}/S^1)),
$$
where $e_I(\Si)_\RR$ is the image of the integral cohomology
class $e_I(\Si)$ in real cohomology and $\om_i(\Si)$
is a real cohomology class.
\end{lem}

\begin{proof}
If $N$ is a sufficiently small neighborhood of $\Si$
with a deformation retraction $\pi:N\to\Si$,
then the description of $\bar\Upsilon^i_{\ft,\fs}$ in
Theorem \ref{thm:ExistenceOfSplicedEndsIndex}
and \cite[Proposition 7.2.32]{DK} imply that there is an inclusion of
vector bundles,
$$
\pi^*(\bar\Upsilon^i_{\ft,\fs}|_\Si)|_{N\cap \sM^{\vir,*}_{\ft,\fs}/S^1}
\subset
\bar\Upsilon^i_{\ft,\fs}|_{N\cap \sM^{\vir,*}_{\ft,\fs}/S^1}.
$$
Therefore, the restriction of $\bar e_I$ to $N\cap \sM^{\vir,*}_{\ft,\fs}/S^1$
admits the
factorization
asserted in the lemma.
\end{proof}

\subsection{Relative Euler classes}
\label{subsec:RelEulerClasses}
An extensive discussion of relative characteristic classes
appears in \cite{Kervaire} but we shall only use a
small and self-contained portion of that theory here.
For any oriented vector bundle, $V\to Y$, and section, $s:Y\to V$, a {\em
relative Euler class\/} can be defined as follows. The section $s$
defines a map of pairs,
$$
s:(Y,Y \setminus s^{-1}(0)) \to (V,V \setminus Y).
$$
If $\Th(V)\in H^r(V,V\setminus Y)$ is the Thom class of $V$, where
$r=\rank_\RR(V)$, then the relative Euler class of the section $s:Y\to
V$ is defined to be
\begin{equation}
\label{eq:DefineRelEulerClass}
e(V,s):=s^*\Th(V)\in
H^r(Y,Y \setminus s^{-1}(0);\RR).
\end{equation}
We note the following properties of the relative Euler class.

\begin{lem}
\label{lem:RelEulerClass}
Let $V\to Y$ be an
oriented, real vector bundle of rank $r$ and
$s:Y\to V$ be a section.  Then the following hold.
\begin{enumerate}
\item
\label{item:RelEulerClass1}
If $\jmath_s:\left(Y,\emptyset\right)\to \left(Y,Y\setminus s^{-1}(0)\right)$
is the inclusion map, then
\begin{equation}
\label{eq:RelEulerSatisfies} \jmath_s^*e(V,s)=e(V),
\end{equation}
where $e(V)\in H^r(Y;\ZZ)$ is the Euler class of $V$.
\item
\label{item:RelEulerClass2}
If $s=s_1\oplus s_2$ with respect to
a decomposition $V=W_1\oplus W_2$
as a direct sum of oriented real vector bundles,
then
$$
e(V,s)=e(W_1,s_1)\smile e(W_2,s_2).
$$
\end{enumerate}
\end{lem}

\begin{proof}
Item \eqref{item:RelEulerClass1} follows from \cite[Equation (11.2)]{Kervaire} or
\cite[p. 98]{MilnorStasheff}.  Item \eqref{item:RelEulerClass2} follows from
the product formula for Thom classes (see \cite[Proposition VIII.11.26]{Dold}).
\end{proof}

\begin{rmk}
If $s:Y\to V$ is the zero section, then
Item \eqref{item:RelEulerClass1} of Lemma \ref{lem:RelEulerClass},
reduces to the well-known result that $s^*\Th(V)\in H^r(Y,\emptyset;\ZZ)$
is the absolute Euler class. If one of the sections $s_i$ in
Item \eqref{item:RelEulerClass2} of Lemma \ref{lem:RelEulerClass} is the zero section,
then the resulting product formula decomposes $e(V,s)$ as a cup product
of a relative and an absolute class.
\end{rmk}

Recalling that $\sM^{\vir}_{\ft,\fs}\subset\bar\sM^{\vir}_{\ft,\fs}$ is the top stratum,
we denote
$$
\bar\Upsilon_{\ft,\fs}/S^1
:=
\bar\Upsilon^i_{\ft,\fs}/S^1
\oplus\bar\Upsilon^s_{\ft,\fs}/S^1
\quad\text{and}\quad
\Upsilon_{\ft,\fs}/S^1
:=
\bar\Upsilon_{\ft,\fs}/S^1|_{\sM^{\vir}_{\ft,\fs}/S^1}.
$$
For any $\eps>0$, let $D_\eps:=D_\eps(\Upsilon_{\ft,\fs}/S^1)
\subset \Upsilon_{\ft,\fs}/S^1$
be the open $\eps$-disk subbundle defined by the $L^2$-norm on
the fibers of $\Upsilon_{\ft,\fs}/S^1$.
Let $r$ be the real rank of $\Upsilon_{\ft,\fs}/S^1$.
For the Eilenberg-MacLane space $K(\QQ,r)$,
let $V(\QQ,r)$ be a contractible neighborhood of a basepoint $*\in K(\QQ,r)$.
Let the cocycle
\begin{equation}
\label{eq:DefineUnivCocycle}
\iota(\QQ,r)\in Z^r(K(\QQ,r),V(\QQ,r);\QQ)
\end{equation}
represent the universal class.  There is a map of triples
\begin{equation}
\label{eq:ThomClassifyingMap}
k_\eps:
\left(
\Upsilon_{\ft,\fs}/S^1,\Upsilon^\circ_{\ft,\fs}/S^1,\Upsilon_{\ft,\fs}/S^1\setminus D_\eps
\right)
\to
\left(
K(\QQ,r),V(\QQ,r),*
\right),
\end{equation}
where $\Upsilon^\circ_{\ft,\fs}/S^1$ is the complement of the zero-section,
such that $[k_\eps^*\iota(\QQ,r)]$ is the Thom
class of $\Upsilon_{\ft,\fs}/S^1$.

The following lemma constructs an extension of the
relative Euler class $e(\Upsilon_{\ft,\fs}/S^1,\bchi)$
to a larger subspace of $\bar\sM^{\vir,*}_{\ft,\fs}/S^1$.

\begin{lem}
\label{lem:ExtendingToThomSection}
Let $\bar\bchi$ be the obstruction section of $\bar\Upsilon_{\ft,\fs}/S^1$
appearing in Hypothesis \ref{hyp:Gluing}.
For any open neighborhood $\sU$  of $\bar\bchi^{-1}(0)$
in $\bar\sM^{\vir}_{\ft,\fs}/S^1$, there
is a constant $c>0$ such that $\bar\bchi':=c \bar\bchi$ obeys the following:
\begin{enumerate}
\item
\label{item:ExtendingToThomSection1}
$\bar\bchi'(\sM^{\vir,*}_{\ft,\fs}/S^1\setminus\sU)
\subset \Upsilon_{\ft,\fs}/S^1\setminus D_\eps$.
\item
\label{item:ExtendingToThomSection2}
For the map $k_\eps$ defined in \eqref{eq:ThomClassifyingMap},
$k_\eps\circ \bar\bchi'(\sM^{\sing}_{\ft,\fs}/S^1\setminus\sU)=*$.
\item
\label{item:ExtendingToThomSection3}
The composition $k_\eps\circ \bar\bchi'$ extends as a map of pairs,
\begin{equation}
\label{eq:DefineBarKObstr}
\bar k_{\bchi}:
\left( \bar\sM^{\vir,*}_{\ft,\fs}/S^1\setminus \sM^{\sing}\cap \sU,
\bar\sM^{\vir,*}_{\ft,\fs}/S^1\setminus \sU,
\right)
\to
\left( K(\QQ,r),* \right)
\end{equation}
which is constant on a neighborhood of $\sM^{\sing}_{\ft,\fs}/S^1\setminus\sU$.
\end{enumerate}
\end{lem}

\begin{proof}
By Hypothesis \ref{hyp:Gluing}, we have
$M_{\fs}\times\Sym^\ell(X)\subset \bar\bchi^{-1}(0)$.
Hence, $\bar\sM^{\vir}_{\ft,\fs}/S^1\setminus\sU=\bar\sM^{\vir,*}_{\ft,\fs}/S^1\setminus\sU$
and we can assume that $\bar\sM^{\vir,*}_{\ft,\fs}/S^1\setminus\sU$ is
compact.  By the lower semi-continuity of the $L^2$ norm of $\bar\bchi$
given by \eqref{eq:L2FiberNorm}
on the fibers of $\bar\Upsilon_{\ft,\fs}/S^1$
(see Hypothesis \ref{hyp:Gluing}), there is a minimum value, $\eps_0$, of
the $L^2$ norm of $\bar\bchi$ on $\bar\sM^{\vir,*}_{\ft,\fs}/S^1\setminus\sU$.
Therefore, if $\bar\bchi'=(\eps/\eps_0)\bar\bchi$, then
$\bar\bchi'(\sM^{\vir,*}_{\ft,\fs}/S^1\setminus\sU)
\subset \Upsilon_{\ft,\fs}/S^1\setminus D_\eps$, proving Property \eqref{item:ExtendingToThomSection1}.
Properties \eqref{item:ExtendingToThomSection2} and \eqref{item:ExtendingToThomSection3} follow from
Property \eqref{item:ExtendingToThomSection1} and the definition
\eqref{eq:ThomClassifyingMap}
of $k_\eps$.
\end{proof}

Lemma \ref{lem:ExtendingToThomSection} gives an extension of
the relative Euler class $e(\Upsilon_{\ft,\fs}/S^1,\bchi)$
to $ \bar\sM^{\vir,*}_{\ft,\fs}/S^1\setminus \sM^{\sing}\cap \sU$.
We now compare this extension with
the extension of the Euler class defined in
\eqref{eq:ExtendedInstEuler}.

\begin{lem}
\label{lem:ExtendedRelativeEulerClass} Continue the notation of
Lemma \ref{lem:ExtendingToThomSection}. Define
\begin{equation}
\label{eq:DefineExtendedRelativeEuler}
\bar e(\bar\Upsilon_{\ft,\fs}/S^1,\bar\bchi')
:=
[\bar k_{\bchi}^*\iota(\QQ,r)],
\end{equation}
where the map $\bar k_{\bchi}$ is defined in \eqref{eq:DefineBarKObstr}
and the cocycle $\iota(\QQ,r)$ is defined in \eqref{eq:DefineUnivCocycle}.
If we abbreviate
$\bar\sM(\sU):=\bar\sM^{\vir,*}_{\ft,\fs}/S^1\setminus\left(\sM^{\sing}\cap \sU\right)$ then
$$
\bar e(\bar\Upsilon_{\ft,\fs}/S^1,\bar\bchi')
\in
H^r\left( \bar\sM(\sU),\bar\sM(\sU)\setminus\bar\bchi^{-1}(0);\RR\right).
$$
As an element of
real cohomology, $\bar e(\bar\Upsilon_{\ft,\fs}/S^1,\bar\bchi')$
satisfies
\begin{equation}
\label{eq:ExtendedRelEulerEqualities1}
\iota_{\bchi,1}^*\bar e(\bar\Upsilon_{\ft,\fs}/S^1,\bar\bchi')
=e(\Upsilon_{\ft,\fs}/S^1,\bar\bchi)
\quad\text{and}\quad \jmath_{\bchi,2}^*
\bar e(\bar\Upsilon_{\ft,\fs}/S^1,\bar\bchi')) = \iota_{\bchi,2}^*(\bar e_I\smile
\bar e_s),
\end{equation}
where $\iota_{\bchi,1}$, $\iota_{\bchi,2}$, and $\jmath_{\bchi,2}$
are the inclusion maps,
$$
\begin{CD}
\sM^{\vir}_{\ft,\fs}/S^1 @> \iota_{\bchi,1} >> \bar\sM(\sU) @>
\iota_{\bchi,2} >> \bar\sM^{\vir,*}_{\ft,\fs}/S^1
\end{CD}
$$
and
$$
\jmath_{\bchi,2}:\left( \bar\sM(\sU),\emptyset \right) \to \left(
\bar\sM(\sU), \bar\sM(\sU)\setminus\bar\bchi^{-1}(0) \right),
$$
and $\bar e_I$ is the extension of the
real Euler class of
$\Upsilon_{\ft,\fs}/S^1$ from \eqref{eq:ExtendedInstEuler} and
$\bar e_s:=e(\Upsilon^s_{\ft,\fs}/S^1)$.
\end{lem}

\begin{proof}
The defining property of $k_\eps$ in \eqref{eq:ThomClassifyingMap}, namely, $[k_\eps^*\iota(\QQ,r)]=\Th(\Upsilon_{\ft,\fs}/S^1)$, 
the definition of the relative Euler class in \eqref{eq:DefineRelEulerClass},
and the homotopy between $\bar\bchi$ and $\bar\bchi'$ give
the first equality in  \eqref{eq:ExtendedRelEulerEqualities1}.

By \cite[Theorem 8.1.15]{Spanier}, the cohomology classes
$\jmath_{\bchi,2}^* \bar e(\bar\Upsilon_{\ft,\fs}/S^1,\bar\bchi')$ and
$\iota_{\bchi,2}^*(\bar e_I\smile \bar e_s)$ are determined by
pointed homotopy classes $[\bar k_\bchi]$ and $[k_2]$, respectively,
in $[\bar\sM(\sU),K(\QQ,r)]$.
Let $\sM^{\sing}(U)\subset \bar\sM(\sU)$ be the complement of the top stratum
and let $N^{\sing}\subset \bar\sM(\sU)$ be a neighborhood of $\sM^{\sing}(\sU)$.
Because the inclusion $\sM^{\sing}(\sU)\subset \bar\sM(\sU)$
is a cofibration, we can assume that $N^{\sing}$ deformation retracts
onto $\sM^{\sing}(\sU)$.

If $N^{\sing}$ is sufficiently small, the construction of $\bar k_{\bchi}$
in Lemma \ref{lem:ExtendingToThomSection}
implies that $k_{\bchi}(N^{\sing})=*$, where $*$ is the basepoint
in $K(\QQ,r)$.
The description of the restriction of $\bar e_I$ to a lower
stratum $\Si\subset\bar\sM(\sU)$ in Lemma
\ref{lem:ExtendedObstructionClassOnLowerStrata},
$$
\bar e_I|_\Si=
e_I(\Si)_\QQ\smile \om_i(\Si),
$$
and the argument of Lemma \ref{lem:ExtendingToThomSection}
imply that the homotopy class
$[k_2]$ also admits a representative
$k_2$ satisfying $k_2(\sM^{\sing}(U))=*$.
Because $N^{\sing}$ retracts to $\sM^{\sing}(U)$, we can assume that
$k_2(N^{\sing})=*$.
Thus,
\begin{equation}
\label{eq:ExtendedInImage}
\bar e_I\smile \bar e_s
\in
\Imag\left(
H^r(\bar\sM(\sU),N^{\sing};\QQ)\to H^r(\bar\sM(\sU);\QQ)
\right).
\end{equation}
By the first equality in  \eqref{eq:ExtendedRelEulerEqualities1}
and the relation between relative and absolute Euler classes in
\eqref{eq:RelEulerSatisfies},
$$
\iota_{\bchi,1}^*
\left(
\jmath_{\bchi,2}^*\bar e(\bar\Upsilon_{\ft,\fs}/S^1,\bar\bchi')
-
\iota_{\bchi,2}^*(\bar e_I\smile\bar e_s)
\right)
=0.
$$
Thus,
\begin{equation}
\label{eq:EulerClassDifference}
\jmath_{\bchi,2}^*\bar e(\bar\Upsilon_{\ft,\fs}/S^1,\bar\bchi')
-
\iota_{\bchi,2}^*(\bar e_I\smile\bar e_s)
\in
\Imag
\left(
H^r(\bar\sM(\sU),\sM^{\vir}_{\ft,\fs}/S^1;\QQ)\to H^r(\bar\sM(\sU);\QQ)
\right).
\end{equation}
By excision,
$$
H^r(\bar\sM(\sU),\sM^{\vir}_{\ft,\fs}/S^1;\QQ)
\cong
H^r(N^{\sing},N^{\sing}-\sM^{\sing}(\sU);\QQ),
$$
and by the construction of $k_{\bchi}$ and by
\eqref{eq:ExtendedInImage}, both $\bar e(\bar\Upsilon_{\ft,\fs}/S^1,\bar\bchi')$
and $\bar e_I\smile \bar e_s$ vanish on cycles in $N^{\sing}$.
Hence, the difference \eqref{eq:EulerClassDifference}
vanishes, completing the proof of the
second equality in \eqref{eq:ExtendedRelEulerEqualities1}.
\end{proof}

\section{Duality and the link of an ideal Seiberg--Witten moduli space}
\label{sec:Duality}
We now perform the computation proving the equality
\eqref{eq:Duality}.

\subsection{The initial duality}
\label{subsec:InitialDuality}
We begin by showing how the intersection number in
\eqref{eq:Duality} can be written as a pairing of relative
cohomology classes with the fundamental class
$[\hat\bL,\rd\hat\bL]$.

The relative cohomology classes describing the geometric
representatives are given by the cocycles defined in Section
\ref{subsubsec:CocyclesPullbacks}. For
$z=\beta_1\cdots\beta_s\in\AAA(X)$ and $\beta_i\in
H_\bullet(X;\RR)$, denote
\begin{equation}
\sK(z,\eta):=\sV(z)\cap \sW^\eta \quad \text{and}\quad \bar\sK(z,\eta) :=
\bar\sV(z)\cap\bar\sW^\eta.
\end{equation}
The inclusion map $\iota_{\beta,1}$  defined in
\eqref{eq:PartialInclusions}  also defines an inclusion of pairs,
$$
\begin{CD}
\left( \sM^{\vir,*}_{\ft,\fs}/S^1 , \sM^{\vir,*}_{\ft,\fs}/S^1
\setminus\sV(\beta) \right)
\\
@V \iota_{\beta,1} VV
\\
\left( \bar\sM^{\vir,*}_{\ft,\fs}/S^1 \setminus \sI(\nu(\beta)),
\bar\sM^{\vir,*}_{\ft,\fs}/S^1 \setminus \sI(\nu(\beta))\setminus\bar\sV(\beta)
\right)
\end{CD}
$$
A similar inclusion of pairs also holds for $\bar\sW$ in place of
$\bar\sV(\beta)$. We then define
\begin{equation}
\label{eq:DefineProductOfCocycles}
[c(z,\eta)]:=\iota_{\beta_1,1}^*[c_{\beta_1}]\smile\dots\smile
\iota_{\beta_r,1}^*[c_{\beta_r}]\smile \iota_{x,1}^*[c_{\sW}]^\eta,
\end{equation}
where $[c_{\beta_i}]$ is as defined in \eqref{eq:DefineRelMuClass}
and $[c_{\sW}]$ is as defined in
\eqref{eq:DefinePreCocycleForWClass}.

With these relative cohomology classes defined, we observe that
dimension-counting arguments yield
\begin{equation}
\bar\sK(z,\eta)\cap\bar\bchi^{-1}(0)\cap\sM^{\sing,*}_{\ft,\fs}/S^1=\emptyset.
\end{equation}
Hence, for a sufficiently small neighborhood $U$ of
$\sM^{\sing}_{\ft,\fs}/S^1$ as constructed in Lemma
\ref{lem:NghOfEnd}, we have
\begin{equation}
\label{eq:RequirementOnNghOfSingStrata}
\bar\sK(z,\eta)\cap\bar\bchi^{-1}(0)\cap U=\emptyset.
\end{equation}
For $\beta\in H_\bullet(X;\RR)$, let $\nu(\beta)$ be the tubular
neighborhood described in Section \ref{subsubsec:CocyclesPullbacks} and let
$\iota_{\beta,1}$ and $\iota_{\beta,2}$ be the inclusion maps defined
in \eqref{eq:PartialInclusions}. For
the complement $U^{\circ}\subset U$ of the lower strata in $U$,
as defined following \eqref{eq:FundClass1},
the cup product
\begin{equation}
[c(z,\eta)] \smile e(\Upsilon_{\ft,\fs}/S^1,\bchi) \in
H^d(\bL^{\vir}_{\ft,\fs}\cup U^{\circ}, \bL^{\vir}_{\ft,\fs}\cup
U^{\circ}\setminus\sK(z,\eta)\cap\bchi^{-1}(0);\RR)
\end{equation}
can be paired with the homology class $[\hat\bL,\rd\hat\bL]$ from
\eqref{eq:FundClass1}.
By the definition of the cocycles $[c_\beta]$ and $[c_\sW]$
as dual to the geometric representatives $\sV(\beta)$ and $\sW$,
the definition of $e(\Upsilon_{\ft,\fs}/S^1,\bchi)$ in
terms of a Thom class, and because $\bchi$ vanishes
transversely (see Hypothesis \ref{hyp:Gluing}), we can write the intersection number as
\begin{equation}
\label{eq:Duality1}
\begin{aligned}
{}&\#\left(\bar\sV(z)\cap\bar\sW^{\eta} \cap
\bar{\bL}_{\ft,\fs}\right)
\\
{}&\quad = \left\langle
\iota_{\beta_1,1}^*[c_{\beta_1}]\smile\dots\smile
\iota_{\beta_r,1}^*[c_{\beta_r}]\smile
\iota_{x,1}^*[c_{\sW}]^\eta\smile e(\Upsilon_{\ft,\fs}/S^1,\bchi),
[\hat\bL,\rd\hat\bL] \right\rangle.
\end{aligned}
\end{equation}
The first equality in Equation
\eqref{eq:ExtendedRelEulerEqualities1} allows us to rewrite
Equation \eqref{eq:Duality1} as
\begin{equation}
\label{eq:Duality2}
\begin{aligned}
{}&\#\left(\bar\sV(z)\cap\bar\sW^{\eta} \cap
\bar{\bL}_{\ft,\fs}\right)
\\
{}&\quad = \left\langle
\iota_{\beta_1,1}^*[c_{\beta_1}]
\smile\dots\smile
\iota_{\beta_r,1}^*[c_{\beta_r}]
\smile
\iota_{x,1}^*[c_{\sW}]^\eta
\smile
\iota_{\bchi,1}^*
\bar e(\bar\Upsilon_{\ft,\fs}/S^1,\bar\bchi'), [\hat\bL,\rd\hat\bL] \right\rangle.
\end{aligned}
\end{equation}
To replace the pairing of the relative classes in
\eqref{eq:Duality1} by pairings of absolute classes with the
fundamental class of $\bar\bL^{\vir}_{\ft,\fs}$, we need to write
the cohomology class $[c(z,\eta)]\smile
e(\Upsilon_{\ft,\fs},\bchi)$ as
$$
[c(z,\eta)]\smile e(\Upsilon_{\ft,\fs},\bchi)= \iota^*[\bar
c(z,\eta)]\smile \bar e(\Upsilon_{\ft,\fs},\bchi),
$$
where
$$
\jmath^*[\bar c(z,\eta)]\smile \bar e(\Upsilon_{\ft,\fs},\bchi) =
\barmu_p(z)\smile \barmu_c^{\eta}\smile \bar e_I \smile\bar
e_s.
$$
We will show how to accomplish this in the following sections.

\subsection{Extension of the cocycles}
\label{subsec:Extending_cocycles}
The relative cohomology classes,
\begin{align*}
[c_\beta]
{}&\in H^\bullet(\bar\sM^{\vir,*}_{\ft,\fs}/S^1\setminus\sI(\nu(\beta)),
\bar\sM^{\vir,*}_{\ft,\fs}/S^1 \setminus\left(\sI(\nu(\beta)) \cup\bar\sV(\beta)\right);\RR),
\\
[c_{\sW}]
{}& \in H^\bullet(\bar\sM^{\vir,*}_{\ft,\fs}/S^1\setminus\sI(\nu(x)),
\bar\sM^{\vir,*}_{\ft,\fs}/S^1 \setminus\left(\sI(\nu(x)) \cup\bar\sV(\beta)\right);\RR),
\\
\bar e(\bar\Upsilon_{\ft,\fs}/S^1,\bar\bchi')
{}& \in H^\bullet\left(\bar\sM^{\vir,*}_{\ft,\fs}/S^1\setminus\left(\sM^{\sing}\cap U'\right),
\bar\sM^{\vir,*}_{\ft,\fs}/S^1 \setminus \left(\sM^{\sing}\cap U'\cup\bchi^{-1}(0)\right);\RR\right),
\end{align*}
defined in \eqref{eq:DefineRelMuClass},
\eqref{eq:DefinePreCocycleForWClass}, and Lemma
\ref{lem:ExtendedRelativeEulerClass}, respectively, are extensions
only over subspaces of $\bar\sM^{\vir,*}_{\ft,\fs}/S^1$. We now
discuss how these cohomology classes can be extended to relative
cohomology classes on $\bar\sM^{\vir,*}_{\ft,\fs}/S^1$.

For any pair of spaces $(A,B)$, let $S^p(A,B;\RR)=\Hom(S_p(A)/S_p(B);\RR)$ denote the complex of
real-valued relative simplicial cochains  and let
$Z^p(A,B;\RR)$ be the cycles of this complex. We can consider
$S^p(A,B;\RR)$ to be the set of real-valued cochains on $A$ which vanish
when restricted to $B$. If $\fU$ is an open cover of $A$, then a
chain is {\em small of order $\fU$\/} (see \cite[Remark 15.8]{GreenbergHarper})
if it is a sum of simplices
each of which has image contained in some element of $\fU$. By \cite[Theorem
15.9]{GreenbergHarper}, any element of $H_p(A,B;\RR)$ can be
represented by a cycle which is small of order $\fU$. By the
Universal Coefficient Theorem, to define an element of
$H^p(A,B;\RR)$ it then suffices to define it on cycles which are
small of order $\fU$.

Lemma \ref{lem:ExtensionPerturbationV} below shows that we can change the cocycle $c_\beta$
to one that extends over $\bar\sM^{\vir,*}_{\ft,\fs}/S^1$ by adding
an exact cocycle supported near $\sI(\nu(\beta))$.  In the proof of
Proposition \ref{prop:Duality}, we will show that such a change does not change the
pairing \eqref{eq:Duality2}.

\begin{lem}
\label{lem:ExtensionPerturbationV}
For $\beta\in H_\bullet(X;\RR)$, let
$\sU(\beta)\subset\bar\sM^{\vir,*}_{\ft,\fs}/S^1$ be any open
neighborhood of $\sI(\nu(\beta))$ satisfying
$\sI(\nu(\beta))\Subset \sU(\beta)$.
Let $c_\beta$ be
the relative cocycle defined in \eqref{eq:DefineRelMuClass}.
Then there are a cochain,
$$
\theta_\beta\in S^{\deg(\beta)-1} \left(
\bar\sM^{\vir,*}_{\ft,\fs}/S^1\setminus\sI(\nu(\beta)),
\bar\sM^{\vir,*}_{\ft,\fs}/S^1\setminus\sU(\beta);\RR \right)
$$
and a cocycle
$$
\bar c_\beta\in Z^{\deg(\beta)}\left(\bar\sM^{\vir,*}_{\ft,\fs}/S^1,
\bar\sM^{\vir,*}_{\ft,\fs}/S^1\setminus\bar\sV(\beta)\setminus\sU(\beta);\RR\right)
$$
satisfying
$$
\jmath_{\sU(\beta)}^*[\bar c_\beta]=\barmu_p(\beta)\quad
\text{and}\quad
\iota_{\beta,2}^*\bar c_\beta=c_\beta+\delta^*\theta_\beta,
$$
where the map $\iota_{\beta,2}$ is
defined in \eqref{eq:PartialInclusions}  and
$$
\jmath_{\sU(\beta)}: \left(
\bar\sM^{\vir,*}_{\ft,\fs}/S^1,\emptyset \right) \to \left(
\bar\sM^{\vir,*}_{\ft,\fs}/S^1,\bar\sM^{\vir,*}_{\ft,\fs}/S^1\setminus\bar\sV(\beta)\setminus\sU(\beta)
\right)
$$
is the inclusion map.
\end{lem}

\begin{proof}
The conclusion of the lemma follows from an argument similar to the proof of \cite[Lemma
5.14]{FLLevelOne}.
If
\begin{equation}
\label{eq:AdditionalInclusionsForV}
\jmath_{\beta,2}:
\left(
    \bar\sM^{\vir,*}_{\ft,\fs}/S^1 \setminus\sI(\nu(\beta)),\emptyset
\right)
\to
\left(
    \bar\sM^{\vir,*}_{\ft,\fs}/S^1 \setminus\sI(\nu(\beta)),
    \bar\sM^{\vir,*}_{\ft,\fs}/S^1\setminus\sI(\nu(\beta))\setminus\bar\sV(\beta)) \right)
\end{equation}
is the inclusion map, then the equality,
$$
\iota_{\beta,2}^*\barmu_p(\beta)=\jmath_{\beta,2}^*[c_\beta],
$$
provided by Corollary \ref{cor:CohomologyClasses}
implies that there are a cocycle $\bar c_\beta' \in
Z^{\deg(\beta)}\left(\bar\sM^{\vir,*}_{\ft,\fs}/S^1;\RR\right)$ and a
cochain
$$
\theta_0\in S^{\deg(\beta)-1}
\left(\bar\sM^{\vir,*}_{\ft,\fs}/S^1\setminus\sI(\nu(\beta));\RR\right),
$$
satisfying
$$
\barmu_p(\beta)=[\bar c_\beta'] \quad\text{and}\quad
\iota_{\beta,2}^*\bar c_\beta'= j_{\beta,2}^*c_\beta
+\delta^*\theta_0.
$$
Because $\sI(\nu(\beta))\Subset\sU(\beta)$, there is a closed
neighborhood, $\sU'(\beta)$, of $\sI(\nu(\beta))$ satisfying
$\sI(\nu(\beta))\Subset\sU'(\beta)\subset\sU(\beta)$. Then
$$
\fU=\left\{\bar\sM^{\vir,*}_{\ft,\fs}/S^1\setminus \sU'(\beta),\sU(\beta)\right\}
$$
is an open cover of $\bar\sM^{\vir,*}_{\ft,\fs}/S^1$ and using the
argument of the proof of the
Excision Lemma (see \cite[Proposition 2.21]{Hatcher}), we can assume that
all cochains in this discussion are defined only on chains
which are small with respect to this open cover.

Because the pairs appearing in the diagram are
excisive couples (see \cite[Theorem 4.6.3]{Spanier} and \cite[p. 218]{Spanier}), the map
$$
\begin{CD}
S^p\left( \bar\sM^{\vir,*}_{\ft,\fs}/S^1\setminus\sI(\nu(\beta)),
\bar\sM^{\vir,*}_{\ft,\fs}/S^1\setminus \sU(\beta);\RR\right) \oplus
S^p\left(
\bar\sM^{\vir,*}_{\ft,\fs}/S^1\setminus\sI(\nu(\beta)),\sU'(\beta);\RR
\right)
\\
@VVV
\\
S^p\left( \bar\sM^{\vir,*}_{\ft,\fs}/S^1\setminus\sI(\nu(\beta)),
\left(\bar\sM^{\vir,*}_{\ft,\fs}/S^1\setminus \sU(\beta)\right)\cap\sU'(\beta);\RR
\right)
\end{CD}
$$
is surjective.  In addition, because   the intersection
$$
\left(\bar\sM^{\vir,*}_{\ft,\fs}/S^1\setminus \sU(\beta)\right)\cap
\sU'(\beta)
$$
is empty, the map is actually a map to the space of absolute
cochains,
$$
S^p(\bar\sM^{\vir,*}_{\ft,\fs}/S^1\setminus\sI(\nu(\beta)),\emptyset;\RR).
$$
Thus we can write $\theta_0=\theta_\beta+\theta_p$, where
\begin{align*}
\theta_\beta\in & S^{\deg(\beta)-1} \left(
\bar\sM^{\vir,*}_{\ft,\fs}/S^1\setminus\sI(\nu(\beta)),\bar\sM^{\vir,*}_{\ft,\fs}/S^1\setminus \sU(\beta);
\RR \right),
\\
\theta_p \in & S^{\deg(\beta)-1} \left(
\bar\sM^{\vir,*}_{\ft,\fs}/S^1\setminus\sI(\nu(\beta)),\sU'(\beta);\RR
\right).
\end{align*}
By definition, $\theta_p$ vanishes on all simplices in
$\sU'(\beta)$.
Because we have assumed that we only need to
define cochains on chains which are small with respect to the open
cover $\fU$ and because $\sI(\nu(\beta))\Subset \sU'(\beta)$, then
$\theta_p$ defines a cochain,
$$
\bar\theta_p \in S^{\deg(\beta)-1} \left(
\bar\sM^{\vir,*}_{\ft,\fs}/S^1,\sU'(\beta);\RR \right),
$$
with $\iota_{\beta,2}^*\bar\theta_p=\theta_p$. Define $\bar
c_\beta:=\bar c'_\beta-\delta^*\bar \theta_p$ and observe that
$$
\iota_{\beta,2}^*\bar c_\beta = \iota_{\beta,2}^*\bar
c_\beta'-\iota_{\beta,2}^*\delta^*\bar\theta_p =
\iota_{\beta,2}^*\bar c_\beta'-\delta^*\theta_p =
\jmath_{\beta,2}^* c_\beta+\delta^*\theta_\beta,
$$
as required.  The preceding equation implies that the cocycle
$\bar c_\beta$ vanishes when restricted to
$\bar\sM^{\vir,*}_{\ft,\fs}/S^1\setminus\left(\bar\sV(\beta)\cup\sU\right)$ and thus
defines the relative cohomology class asserted by the lemma. Finally, we compute
that
$$
\jmath_{\sU(\beta)}^*[\bar c_\beta] = [\bar
c_\beta'-\delta^*\bar\theta_p] = [\bar c_\beta'] =
\barmu_p(\beta),
$$
completing the proof of Lemma \ref{lem:ExtensionPerturbationV}.
\end{proof}

The
argument of Lemma \ref{lem:ExtensionPerturbationV}
gives corresponding extension results for
$[c_{\sW}]$ and $e(\Upsilon,\hat\bchi)$.
The extension result for $[c_{\sW}]$ is given by

\begin{lem}
\label{lem:ExtensionPerturbationW}
Let $c_{\sW}$ be the cochain
defined in \eqref{eq:DefinePreCocycleForWClass}. For any
neighborhood $\sU(x)\subset\bar\sM^{\vir,*}_{\ft,\fs}/S^1$ of
$\sI(\nu(x))$ satisfying $\sI(\nu(x))\Subset \sU(x)$, there are a
cochain,
$$
\theta_{\sW}\in
S^1\left(\bar\sM^{\vir,*}_{\ft,\fs}/S^1\setminus\sI(\nu(x)),
\bar\sM^{\vir,*}_{\ft,\fs}/S^1\setminus \sU(x);\RR\right),
$$
and a cocycle,
$$
\bar c_{\sW}\in
Z^2(\bar\sM^{\vir,*}_{\ft,\fs}/S^1,\bar\sM^{\vir,*}_{\ft,\fs}/S^1\setminus\bar\sW\setminus \sU(x);\RR),
$$
satisfying
$$
\jmath_{\sU(x)}^*[\bar c_{\sW}]=\barmu_c \quad\text{and}\quad
\iota_{x,2}^*\bar c_{\sW}=c_{\sW}+\delta^*\theta_{\sW},
$$
where
$$
\jmath_{\sU(x)}: \left( \bar\sM^{\vir,*}_{\ft,\fs}/S^1,\emptyset
\right) \to \left(
\bar\sM^{\vir,*}_{\ft,\fs}/S^1,\bar\sM^{\vir,*}_{\ft,\fs}/S^1\setminus\bar\sW\setminus \sU(x)
\right)
$$
is the inclusion map and $\iota_{x,2}$ is defined
in \eqref{eq:PartialInclusions}.
\end{lem}

Finally, we record
the extension result for
$\bar e(\bar\Upsilon_{\ft,\fs}/S^1,\bar\bchi')$.

\begin{lem}
\label{lem:ExtensionPerturbationChi} Continue the notation of
Lemma \ref{lem:ExtendedRelativeEulerClass}. For any neighborhood
$U_{\bchi}\subset\bar\sM^{\vir,*}_{\ft,\fs}/S^1$ of
$\sM^{\sing}\cap \sU$ with $\sM^{\sing}\cap \sU\Subset U_{\bchi}$,
there are a cochain,
$$
\theta_\bchi\in S^{r-1}\left( \bar\sM(\sU),
\bar\sM^{\vir,*}_{\ft,\fs}/S^1 \setminus U_{\bchi}\right),
$$
and a cocycle,
$$
\bar e_\bchi \in Z^r\left(\bar\sM^{\vir,*}_{\ft,\fs}/S^1,
\bar\sM^{\vir,*}_{\ft,\fs}/S^1\setminus\bchi^{-1}(0)\setminus U_{\bchi};\RR\right),
$$
satisfying
$$
\jmath_{\bchi}^*[\bar e_\bchi] =\bar e_I\smile \bar e_{\fs}
\quad\text{and}\quad \iota_{\bchi,2}^* \bar e_{\bchi}=
\bar e(\bar\Upsilon_{\ft,\fs}/S^1,\bar\bchi')+\delta^*\theta_\bchi,
$$
where
$$
\jmath_{\bchi}: \left( \bar\sM^{\vir,*}_{\ft,\fs}/S^1,\emptyset
\right) \to \left(
\bar\sM^{\vir,*}_{\ft,\fs}/S^1,\bar\sM^{\vir,*}_{\ft,\fs}/S^1\setminus\bchi^{-1}(0)\setminus U_{\bchi}
\right)
$$
is the inclusion map.
\end{lem}

We can now give the

\begin{proof}[Proof of Proposition \ref{prop:Duality}]
The proof is analogous to that of \cite[Proposition 5.2]{FLLevelOne},
using relative Euler classes of the obstruction
sections and representatives of the cohomology with compact
support along the geometric representatives.

Let $\sU$ be the neighborhood of $\bar\bchi^{-1}(0)$ appearing
in Lemma \ref{lem:ExtendingToThomSection}.
Because the intersection,
\begin{equation}
\label{eq:IntersectionOfLowerStrat}
\sI(\nu(\beta_1))\cap\dots\cap\sI(\nu(\beta_k)) \cap
\sI(\nu(x_1))\cap\dots\cap\sI(\nu(x_\eta)) \cap \sU \cap\sM^{\sing},
\end{equation}
is empty, we can find open neighborhoods $U(\beta_i)$ of
$\sI(\nu(\beta_i))$, and $U_{\sW,j}$ of $\sI(\nu(x_j))$, and $U_\bchi$
of $\sU$, and $U$ of $\sM^{\sing}$, such that the intersection
obtained by replacing any of $\sI(\nu(\beta_i))$ with
$U(\beta_i)$, or $\sI(\nu(x_j))$ with $U_{\sW,j}$, or $\sU$ with
$U_\bchi$, or $\sM^{\sing}$ with $U$ in
\eqref{eq:IntersectionOfLowerStrat} is still empty.
For pairs of spaces and subspaces, $(Y,A)$ and $(Y,B)$, cup product defines a map of relative cohomology,
$$
H^p(Y,A;\RR)\otimes H^q(Y,B;\RR)
\to
H^{p+q}(Y,A\cup B;\RR).
$$
Therefore, if we let $\theta_{\beta_i}$, $\theta_{\nu}$, and $\theta_{\bchi}$
be the cochains
defined in Lemmas \ref{lem:ExtensionPerturbationV},
\ref{lem:ExtensionPerturbationW}, and
\ref{lem:ExtensionPerturbationChi} respectively, for the
neighborhoods $U(\beta_i)$, $U_{\sW,j}$ and $U_{\bchi}$,
and we replace one or more of
$[c_{\beta_1}]$ with $\delta^*\theta_{\beta_i}$, or
$[c_{\sW}]$ with $\delta^*\theta_{\nu}$,
or
$\bar e(\bar\Upsilon_{\ft,\fs}/S^1,\bar\bchi')$ with $\delta^*\theta_{\bchi}$,
in
$$
\iota_{\beta_1,1}^*[c_{\beta_1}]\smile\dots\smile
\iota_{\beta_r,1}^*[c_{\beta_r}]\smile
\iota_{x,1}^*[c_{\sW}]^\eta\smile \iota_{\bchi,1}^*
\bar e(\bar\Upsilon_{\ft,\fs}/S^1,\bar\bchi'),
$$
then the resulting cup product will vanish.
Hence, we
obtain the following equality for the pairing \eqref{eq:Duality2},
\begin{equation}
\label{eq:Duality3}
\begin{aligned}
{}& \left\langle \iota_{\beta_1,1}^*[c_{\beta_1}]\smile\dots\smile
\iota_{\beta_r,1}^*[c_{\beta_r}]\smile
\iota_{x,1}^*[c_{\sW}]^\eta\smile
\iota_{\bchi,1}^*\bar e(\bar\Upsilon_{\ft,\fs}/S^1,\bar\bchi'), [\hat\bL,\rd\hat\bL] \right\rangle
\\
{}&\quad = \left\langle
\iota_{\beta_1,1}^*[c_{\beta_1}+\delta^*\theta_{\beta_1}] \smile
\dots \smile
\iota_{\beta_r,1}^*[c_{\beta_s}+\delta^*\theta_{\beta_s}] \right.
\\
{}&\qquad \left.\smile \iota_{x,1}^*[c_{\sW}+\delta^*\theta_\nu]^\eta
\smile
\iota_{\bchi,1}^*(\bar e(\bar\Upsilon_{\ft,\fs}/S^1,\bar\bchi')+\delta^*\theta_\bchi),
[\hat\bL,\rd\hat\bL] \right\rangle.
\end{aligned}
\end{equation}
Recall that for $j=1,2$, the inclusion maps $\iota_{\beta,j}$ and
$\iota_{\bchi,j}$ are defined in
\eqref{eq:PartialInclusions} and
Lemma \ref{lem:ExtendedRelativeEulerClass}, respectively.
Applying the equalities
$$
\iota=\iota_{\beta,2}\circ\iota_{\beta,1}=\iota_{\bchi,2}\circ\iota_{\bchi,1}
$$
and Lemmas \ref{lem:ExtensionPerturbationV},
\ref{lem:ExtensionPerturbationW}, and
\ref{lem:ExtensionPerturbationChi} to equations
\eqref{eq:Duality2} and \eqref{eq:Duality3} then implies that
for the cocycle $\bar c_{\beta_1}$ defined in Lemma \ref{lem:ExtensionPerturbationV},
the cocycle $\bar c_{\sW}$ defined in Lemma \ref{lem:ExtensionPerturbationW},
and the cocycle $\bar e_{\bchi}$ defined in Lemma \ref{lem:ExtensionPerturbationChi},
we obtain the equalities,
\begin{equation}
\label{eq:Duality4}
\begin{aligned}
\#\left(\bar\sV(z)\cap\bar\sW^{\eta} \cap
\bar{\bL}_{\ft,\fs}\right)
{}&= \left\langle \iota^*\left( [\bar c_{\beta_1}]\smile\dots [\bar
c_{\beta_s}] \smile [\bar c_{\sW}]^\eta \smile [\bar e_{\bchi}]
\right), [\hat\bL,\rd\hat\bL] \right\rangle
\\
{}&= \left\langle [\bar c_{\beta_1}]\smile\dots [\bar c_{\beta_s}]
\smile [\bar c_{\sW}]^\eta \smile [\bar e_{\bchi}],
\bar\jmath_*[\bar\bL^{\vir}_{\ft,\fs}] \right\rangle,
\end{aligned}
\end{equation}
where the map $\bar\jmath$ is defined in \eqref{eq:CompactJMap}
and we have used the definition of the homology class
$[\bar\bL^{\vir}_{\ft,\fs}]$ given in
\eqref{eq:DefineAmbientLinkFund}. If $\jmath_{\sU(\beta_i)}$,
$\jmath_{\sU(x)}$, and $\jmath_{\bchi}$
are the inclusion maps defined in Lemmas \ref{lem:ExtensionPerturbationV},
\ref{lem:ExtensionPerturbationW}, and
\ref{lem:ExtensionPerturbationChi}, then
\begin{align*}
{}&\#\left(\bar\sV(z)\cap\bar\sW^{\eta} \cap
\bar{\bL}_{\ft,\fs}\right)
\\
{}&= 
\left\langle \bar\jmath^*\left([\bar c_{\beta_1}]\smile\dots \smile[\bar c_{\beta_s}] \smile [\bar c_{\sW}]^\eta 
\smile [\bar e_{\bchi}]\right), [\bar\bL^{\vir}_{\ft,\fs}] \right\rangle
\\
{}&= \left\langle \jmath_{\sU(\beta_1)}^*[\bar c_{\beta_1}]
\smile\dots\smile \jmath_{\sU(\beta_s)}^*[\bar c_{\beta_s}] \smile
\jmath_{\sU(x)}^*[\bar c_{\sW}]^\eta \smile \jmath_{\bchi}^*[\bar
e_{\bchi}], [\bar\bL^{\vir}_{\ft,\fs}] \right\rangle
\\
&\qquad\text{(by \cite[Theorem 5.6.8]{Spanier})}
\\
{}&= \left\langle \barmu_p(\beta_1)\smile\dots\smile\barmu_p(\beta_s) \smile
\barmu_c^\eta \smile \bar e_I\smile\bar e_s,
[\bar\bL^{\vir}_{\ft,\fs}] \right\rangle,
\end{align*}
where the final equality follows from Lemmas
\ref{lem:ExtensionPerturbationV},
\ref{lem:ExtensionPerturbationW}, and
\ref{lem:ExtensionPerturbationChi}. This completes the proof of Proposition \ref{prop:Duality}.
\end{proof}

\section{Reduction to the subspace $\bar{\bB\bL}^{\vir}_{\ft,\fs}$}
\label{sec:ReduceToBase}
We now show how to rewrite pairings
with the homology class $[\bar\bL^{\vir}_{\ft,\fs}]$ defined in
\eqref{eq:DefineAmbientLinkFund} as pairings with the homology
class of a subspace
and thus eliminate the topology of the normal
bundle $N_{\ft(\ell),\fs}(\delta)\to M_{\fs}$.

We first give a global description of the relation between the space
$\bar{\bB\bL}^{\vir}_{\ft,\fs}$ defined in
\eqref{eq:InstantonLinkComponentBase} and the space
$\bar\bL^{\vir,i}_{\ft,\fs}\subset\bar\bL^{\vir}_{\ft,\fs} $
defined in \eqref{eq:AmbientLink}, \eqref{eq:InstantonLinkComponent}, and \eqref{eq:DefineLinkStratum}.
Let $\tilde\iota_\nu$ and $\iota_\nu$ be the maps in the classifying
diagram for the
cohomology class $\nu$ of Definition \ref{defn:DefineNu},
\begin{equation}
\label{eq:S1ClassifyingMap}
\begin{CD}
\bar{\sM}^{\vir,*}_{\ft,\fs} @> \tilde \iota_\nu >> \ES^1
\\
@VVV @VVV
\\
\bar{\sM}^{\vir,*}_{\ft,\fs}/S^1 @> \iota_\nu >> \BS^1
\end{CD}
\end{equation}
The proof of the following lemma is identical to that of \cite[Lemma
5.19]{FLLevelOne}.

\begin{lem}
\label{lem:BaseSpaceNormalBundle}
Let $S^1$ act on the fibers of $N_{\ft(\ell),\fs}(\delta)\to M_{\fs}$
by scalar multiplication.  Then
the deformation retraction,
$\bar\bL^{\vir,i}_{\ft,\fs}\to \bar{\bB\bL}^{\vir}_{\ft,\fs}$, is the
pullback of the disk bundle,
$$
N_{\ft(\ell),\fs}(\delta)\times_{S^1}\ES^1 \to M_{\fs}\times\BS^1,
$$
by the map
$$
\pi_N\times\tilde\iota_\nu: \bar\sM^{\vir,*}_{\ft,\fs}/S^1 \to
N_{\ft(\ell),\fs}(\delta)\times_{S^1}\ES^1,
$$
where $\bar\sM^{\vir,*}_{\ft,\fs}$ is defined in \eqref{eq:RedComplementInVirtual},
the map $\pi_N$ is defined in \eqref{eq:GlobalProjectionToN}, and the map
$\tilde\iota_\nu$ in \eqref{eq:S1ClassifyingMap}.
\end{lem}

Let $r_N$ denote the complex rank of $N_{\ft(\ell),\fs}\to M_{\fs}$ and
define
\begin{equation}
\label{eq:DefineThomClass}
\Th(t_N) \in
H^{2r_N}\left(\bar\sM^{\vir,*}_{\ft,\fs}/S^1,\bar\sM^{\vir,*}_{\ft,\fs}/S^1 \setminus t_N^{-1}(0);\RR\right)
\end{equation}
to be the pullback of the Thom class of the bundle
$$
N_{\ft(\ell),\fs}(\delta)\times_{S^1}\ES^1 \to M_{\fs}\times\BS^1
$$
by the map $\pi_N\times\tilde\iota_\nu$ appearing in
Lemma \ref{lem:BaseSpaceNormalBundle}. A fundamental class,
\begin{equation}
\label{eq:DefineBaseBundClass}
[\bar{\bB\bL}^{\vir}_{\ft,\fs}] \in
H_{d(\ft)-2r_N-2}\left(\bar\bL^{\vir}_{\ft,\fs}\cap t_N^{-1}(0);\RR\right),
\end{equation}
is defined exactly as in \eqref{eq:DefineAmbientLinkFund}.
Namely, let
\begin{equation}
\label{eq:BaseFundClass1}
[\widehat{\bB\bL},\rd\bB\bL]
\end{equation}
be the homology class defined by the fundamental class of the
manifold with corners given by the intersection of $t_N^{-1}(0)$
with the manifold defining the homology class
$[\hat\bL,\rd\hat\bL]$ from \eqref{eq:FundClass1}. Now define
$[\bar{\bB\bL}^{\vir}_{\ft,\fs}]$ to be the homology class
satisfying
\begin{equation}
\label{eq:DefineBFundClass}
\iota_* [\widehat{\bB\bL},\rd\bB\bL] =
\bar\jmath_*[\bar{\bB\bL}^{\vir}_{\ft,\fs}],
\end{equation}
where $\iota$ and $\bar\jmath$ are the maps defined in Equations
\eqref{eq:TopStratumOfVirtualInclusion} and \eqref{eq:CompactJMap},
respectively. Note that the restriction of $\bar\jmath$ to
$t_N^{-1}(0)$ also defines an isomorphism on homology of the appropriate
dimension because the deformation retraction $r$
in Lemma \ref{lem:NghOfEnd}
preserves the level sets of $t_N$.

Let $U$ be the neighborhood of $\sM^{\sing}_{\ft,\fs}/S^1$
used in \eqref{eq:FundClass1} to
define the homology class $[\hat\bL,\rd\hat\bL]$.
Because we can choose the boundaries of the space $U$ to be
generic, we can assume that $t_N$ vanishes transversely on these
boundaries.  Hence, we have the equality,
\begin{equation}
[\widehat{\bB\bL},\rd\bB\bL]= \iota^*\Th(t_N)\cap [\hat\bL,\rd\hat\bL].
\end{equation}
We then have

\begin{lem}
\label{lem:ReduceToBase}
If $b:=\dim \bB\bL^{\vir}_{\ft,\fs}$, then for any $\om\in H^b(\bar\sM^{\vir,*}_{\ft,\fs}/S^1)$ we have
$$
\left\langle
\iota_{\bB}^*\om,[\bar{\bB\bL}^{\vir}_{\ft,\fs}]
\right\rangle
=
\left\langle
\om\smile \jmath_\bB^*\Th(t_N),[\bL^{\vir}_{\ft,\fs}]
\right\rangle,
$$
where $\Th(t_N)$ is defined in \eqref{eq:DefineThomClass}
and
$$
\jmath_{\bB}:
\left(
\bar\sM^{\vir,*}_{\ft,\fs}/S^1,\emptyset
\right)
\to
\left(
\bar\sM^{\vir,*}_{\ft,\fs}/S^1,\bar\sM^{\vir,*}_{\ft,\fs}/S^1 \setminus t_N^{-1}(0)
\right)
$$
and
$$
\iota_{\bB}:t_N^{-1}(0)\to \bar\sM^{\vir,*}_{\ft,\fs}/S^1
$$
are the inclusion maps.
\end{lem}

\begin{proof}
The proof is similar to that of \cite[Equation (5.59)]{FLLevelOne}.

Let $U$ be the neighborhood of $\sM^{\sing}_{\ft,\fs}/S^1$ as
described prior to the statement of the lemma.
Abbreviate $\bar\bL^{\vir}_{\ft,\fs}\cup U$ by $\bL(U)$
and write $t_N^{-1}(0)=\bB$.
Then
for $d(\ft)=\dim \sM^{\vir}_{\ft,\fs}$ and
$d(\ft)-2=b+2r_N$, where $r_N$
is the complex rank of the vector bundle $N_{\ft(\ell),\fs}\to M_{\fs}$,
we have
$$
\om\smile\Th(t_N)\in
H^{b+2r_N}(\bL(U),\bL(U)\setminus \bB;\RR).
$$
Consider the following commutative diagram of inclusions:
\begin{equation}
\label{eq:RelativeIncl1}
\begin{CD}
@.
\left(\bL(U),\emptyset \right)
@>\jmath_{\bB}>>
\left(\bL(U),\bL(U)\setminus \bB\right)
\\
@. @V\bar\jmath VV @V\jmath_{\bB,U}VV
\\
@. \left(\bL(U),U\right)
@>\jmath_{U,\bB}>>
\left(\bL(U),U\cup(\bL(U)\setminus \bB)\right)
\end{CD}
\end{equation}
Define $U^{\sing}:=U\cap\sM^{\sing}_{\ft,\fs}/S^1$.
Because $U$ retracts onto $U^{\sing}$ and $\dim U^{\sing}\le d(\ft)-1-4$,
the inclusion map
$\jmath_{\bB,U}$ defines an isomorphism,
$$
\jmath_{\bB,U}^*:{}
H^{b+2r_N}\left(\bL(U),U\cup(\bL(U)\setminus \bB);\RR\right)
\cong
H^{b+2r_N}\left(\bL(U),\bL(U)\setminus \bB;\RR\right).
$$
Consequently, there is a unique class
$\Om\in H^{b+2r_N}\left(\bL(U),U\cup(\bL(U)\setminus \bB);\RR\right)$
such that
\begin{equation}
\label{eq:DefineOm}
\jmath_{\bB,U}^*\Om=\om\smile\Th(t_N).
\end{equation}
We now calculate that
\begin{align*}
\left\langle\jmath_{\bB}^*(\om\smile\Th(t_N)),[\bar\bL^{\vir}_{\ft,\fs}] \right\rangle
{}&=
\left\langle\jmath_{\bB}^*\jmath_{\bB,U}^*\Om,[\bar\bL^{\vir}_{\ft,\fs}] \right\rangle
\quad\text{(by \eqref{eq:DefineOm})}
\\
{}&=
\left\langle \bar\jmath^*\jmath_{U,\bB}^*\Om,[\bar\bL^{\vir}_{\ft,\fs}] \right\rangle
\quad\text{(by commutativity of \eqref{eq:RelativeIncl1})}
\\
{}&=
\left\langle \jmath_{U,\bB}^*\Om,\bar\jmath_*[\bar\bL^{\vir}_{\ft,\fs}] \right\rangle
\\
{}&=
\left\langle \jmath_{U,\bB}^*\Om,\iota_*[\bL,\rd\bL] \right\rangle
\quad\text{(by \eqref{eq:DefineAmbientLinkFund}),}
\end{align*}
that is,
\begin{equation}
\label{eq:Reduction1}
\left\langle\jmath_{\bB}^*(\om\smile\Th(t_N)),[\bar\bL^{\vir}_{\ft,\fs}] \right\rangle
=
\left\langle \jmath_{U,\bB}^*\Om,\iota_*[\bL,\rd\bL] \right\rangle.
\end{equation}
To get a better description of $\Om$,
we now reduce the preceding pairing to one on the subspace,
$$
\bL^i(U):=\bL(U)\setminus\Int(\bar\bL^{\vir,s}_{\ft,\fs}),
$$
where $\bar\bL^{\vir,s}_{\ft,\fs}$ is defined in
\eqref{eq:DefineSWComponentOfLink}.
Observe that
the top stratum, $\bL^i(U)\cap \bL^{\vir}_{\ft,\fs}$,
is an oriented topological manifold with boundary by the argument giving Item \eqref{item:BoundaryOfNeigh3}
of Lemma \ref{lem:BoundaryOfNeigh}
and thus has a relative fundamental
class
(see, for example, \cite[p. 253]{Hatcher})
which we denote by
$$
[\bL^i,\rd\bL^i]
\in
H_{d(\ft)-2}(\bL^i(U),U\cup (\bL^i\setminus \bB);\RR).
$$
Because $\bar\bL^{\vir,s}_{\ft,\fs}\subset \bL(U)\setminus \bB$,
the inclusion of pairs,
$$
\iota_{U,i}:
\left(
\bL^i(U),U\cup (\bL^i\setminus \bB)
\right)
\to
\left(\bL(U),U\cup(\bL(U)\setminus \bB)\right),
$$
is an excision map and thus induces an isomorphism on cohomology.
Consequently,
\begin{equation}
\label{eq:RelFundClassEq}
(\iota_{U,i})_*[\bL^i,\rd\bL^i]=(\jmath_{U,\bB})_*\iota_*[\bL,\rd\bL],
\end{equation}
because both $[\bL^i,\rd\bL^i]$ and $[\bL,\rd\bL]$
are fundamental classes.
We compute that the right-hand side of
Equation \eqref{eq:Reduction1} is given by
\begin{align*}
\left\langle \jmath_{U,\bB}^*\Om,\iota_*[\bL,\rd\bL] \right\rangle
{}&=
\left\langle \Om,(\jmath_{U,\bB})_*\iota_*[\bL,\rd\bL]\right\rangle
\\
{}&=
\left\langle \Om,(\iota_{U,i})_*[\bL^i,\rd\bL^i]\right\rangle
\quad\text{(by \eqref{eq:RelFundClassEq})}
\\
{}&=
\left\langle \iota_{U,i}^*\Om,[\bL^i,\rd\bL^i]\right\rangle,
\end{align*}
that is,
\begin{equation}
\label{eq:Reduction2}
\left\langle \jmath_{U,\bB}^*\Om,\iota_*[\bL,\rd\bL] \right\rangle
=
\left\langle \iota_{U,i}^*\Om,[\bL^i,\rd\bL^i]\right\rangle.
\end{equation}
We now compute $\iota_{U,i}^*\Om$ in terms of
similar restrictions of $\om$ and $\Th(t_N)$.
The $(\sG_{\fs}\times S^1)$-equivariant
deformation retraction, $\tilde N_{\ft(\ell),\fs}\to \tilde M_{\fs}$,
defines a deformation retraction
$$
r_B: \bL^i(U) \to \bB.
$$
Observe that $r_B$ is a homotopy equivalence and also
defines a homotopy equivalence of pairs,
$$
r_{\bB}:(\bL^i(U),U)\to (\bB,\bB\cap U).
$$
Because the retraction of $U$ in Lemma \ref{lem:NghOfEnd} onto the lower strata respects
the level sets of $t_N$,
the subspace $\bB\cap U$ retracts onto the lower strata of
$\bB$ which have codimension greater than or equal to four
in $\bB$.  Therefore, the map
$$
\bar\jmath^*:H^b(\bB,\bB\cap U;\RR)\to H^b(\bB;\RR)
$$
is an isomorphism.
Because the map $r_{\bB}$ induces
isomorphisms $r_{\bB}^*:H^b(\bB,\bB\cap U;\RR)\to H^b(\bL^i(U),U;\RR)$
and $r_{\bB}^*:H^b(\bB;\RR)\to H^b(\bL^i(U);\RR)$, there is a commutative
diagram,
$$
\begin{CD}
H^b(\bB,\bB\cap U;\RR)
@>\bar\jmath^* >>
H^b(\bB;\RR)
\\
@V r_{\bB}^* VV @Vr_{\bB}^* VV
\\
H^b(\bL^i(U),\bL^i(U)\cap U;\RR)
@>\bar\jmath^* >>
H^b(\bL^i(U);\RR)
\end{CD}
$$
Hence, there is a unique
$y_{\bB}\in H^b(\bB,U\cap\bB;\RR)$ such that
\begin{equation}
\label{eq:yBDefining}
\iota_i^*\om= \bar\jmath^*r_{\bB}^*y_{\bB}
\quad\text{and}\quad
\iota_{\bB}^*\om=\bar\jmath^* y_{\bB},
\end{equation}
where $\iota_i:\bL^i(U)\to\bL(U)$
is the inclusion map.
We now claim that
\begin{equation}
\label{eq:RestrOm}
\iota_{U,i}^*\Om=(r_{\bB}^*y_{\bB})\smile\iota_i^*\Th(t_N).
\end{equation}
To see this, we consider the following commutative diagram:
$$
\begin{CD}
H^b(\bL^i(U),U;\RR)
\otimes
H^{2r_N}(\bL^i,\bL^i\setminus \bB;\RR)
@> \bar\jmath^*\otimes\id >>
H^b(\bL^i(U);\RR)
\otimes
H^{2r_N}(\bL^i(U),\bL^i\setminus \bB;\RR)
\\
@V \smile VV @V \smile VV
\\
H^{b+2r_N}(\bL^i(U),U\cup\bL^i\setminus \bB;\RR)
@>\jmath_{\bB,U}^* >>
H^{b+2r_N}(\bL^i(U),\bL^i\setminus \bB;\RR)
\\
@A \iota_i^* AA  @A \iota_i^* AA
\\
H^{b+2r_N}(\bL(U),U\cup (\bL(U)\setminus B);\RR)
@> \jmath_{\bB,U}^* >>
H^{b+2r_N}(\bL(U),\bL(U)\setminus B;\RR)
\end{CD}
$$
The class $r_{\bB}^*(y_{\bB})\otimes\iota_i^*\Th(t_N)$ is in the upper-left
entry of the preceding diagram, while $\Om$ is in the lower-left entry
of the preceding diagram.
To prove that \eqref{eq:RestrOm} holds, we must show that these two classes have the same
image in the center-left entry.
By \eqref{eq:DefineOm}, the class $\iota_i^*\Om$ is mapped
to $\iota_i^*\om\smile\iota_i^*\Th(t_N)$ in
the center-right entry.  By \eqref{eq:yBDefining}, the class
$r_{\bB}^*(y_{\bB})\smile\iota_i^*\Th(t_N)$ is also mapped
to $\iota_i^*\om\smile\iota_i^*\Th(t_N)$ in
the center-right entry, so
$$
\jmath_{\bB,U}^*\iota_i^*\Om=\jmath_{\bB,U}^*(r_{\bB}^*(y_{\bB})\smile\iota_i^*\Th(t_N)).
$$
Because $U$ retracts onto a set of codimension
greater than or equal to two and $b+2r_N=d(\ft)-2$ equals
the dimension of the top stratum of $\bL(U)$, we observe that both
arrows labeled $\jmath_{\bB,U}^*$ in the preceding diagram
are isomorphisms.
Hence, the preceding equality yields Equation \eqref{eq:RestrOm}, as claimed.

By applying Equation \eqref{eq:RestrOm} to
the right-hand side of
Equation
\eqref{eq:Reduction2}, we obtain
\begin{align*}
\left\langle \iota_{U,i}^*\Om,[\bL^i,\rd\bL^i]\right\rangle
{}&=
\left\langle r_{\bB}^*(y_{\bB})\smile\iota_i^*\Th(t_N), [\bL^i,\rd\bL^i]\right\rangle
\\
{}&=
\left\langle r_{\bB}^*y_{\bB},\iota_i^*(\Th(t_N))\cap[\bL^i,\rd\bL^i]\right\rangle,
\end{align*}
that is,
\begin{equation}
\label{eq:Reduction3}
\left\langle \iota_{U,i}^*\Om,[\bL^i,\rd\bL^i]\right\rangle
=
\left\langle r_{\bB}^*y_{\bB},\iota_i^*(\Th(t_N))\cap[\bL^i,\rd\bL^i]\right\rangle.
\end{equation}
From \cite[p. 371, Equation (1)]{BredonTopGeom}, we have the equality
\begin{equation}
\label{eq:RelThomIsom}
\iota_i^*\Th(t_N)\cap [\bL^i,\rd\bL^i]
=
(\iota_{\bB,i})_*\iota_*[\bB,\rd\bB],
\end{equation}
where  $\iota_{\bB,i}:(\bB,\bB\cap U)\to (\bL^i(U),\bL^i\setminus B)$
is the inclusion map.
The right-hand side of
Equation \eqref{eq:Reduction3} is then equal to
\begin{align*}
\left\langle r_{\bB}^*y_{\bB},\iota_i^*(\Th(t_N))\cap[\bL^i,\rd\bL^i]\right\rangle
{}&=
\left\langle r_{\bB}^*y_{\bB},(\iota_{\bB,i})_*\iota_*[\bB,\rd \bB]\right\rangle
\\
{}&=
\left\langle (\iota_{\bB,i})^*r_{\bB}^*y_{\bB},\iota_*[\bB,\rd \bB]\right\rangle
\\
{}&=
\left\langle y_{\bB},\iota_*[\bB,\rd\bB]\right\rangle
\quad\text{(because $r_{\bB}$ is a retraction, so $r_{\bB}\circ \iota_{\bB,i}=\id$)}
\\
{}&=
\left\langle y_{\bB},\jmath_*[\bar{\bB\bL}^{\vir}_{\ft,\fs}]\right\rangle
\quad\text{(by \eqref{eq:DefineBFundClass})}
\\
{}&=
\left\langle \jmath^*y_{\bB},[\bar{\bB\bL}^{\vir}_{\ft,\fs}]\right\rangle
\\
{}&=
\left\langle \iota_{\bB}^*\om,[\bar{\bB\bL}^{\vir}_{\ft,\fs}]\right\rangle
\quad\text{(by \eqref{eq:yBDefining})}.
\end{align*}
This completes the proof of Lemma \ref{lem:ReduceToBase}.
\end{proof}

Lemma \ref{lem:ReduceToBase} can be used to reduce the pairing
in \eqref{eq:Duality} to one with $[\bar{\bB\bL}^{\vir}_{\ft,\fs}]$
when we have an expression for division by the Thom class
$\Th(t_N)$.  The following proposition gives that result.

\begin{prop}
\label{prop:ReductionToBase}
Continue to denote $d_s=\dim M_{\fs}$.
For $0\le j\le [d_s/2]$, let $s_j(N)\in H^{2j}(M_{\fs};\RR)$ be the Segre classes of the
complex-rank-$r_N$ vector bundle $N_{\ft(\ell),\fs}\to M_{\fs}$.
Let $k$ and $m$ be non-negative integers satisfying
$k+2m=\dim \bL^{\vir}_{\ft,\fs}$.
For any $\alpha\in H^k(M_{\fs}\times \Sym^\ell(X);\RR)$
and the Chern class $\nu\in H^2(\bar\bL^{\vir}_{\ft,\fs};\RR)$ defined in Definition \ref{defn:DefineNu},
we have the equality
\begin{equation}
\label{eq:ReduceToBase}
\begin{aligned}
{}& \left\langle \nu^m\smile
\pi_{X,\fs}^*\alpha,[\bar\bL^{\vir}_{\ft,\fs}] \right\rangle
\\
{} &\quad  = \sum_{j=0}^{d_s/2} (-1)^{r_N+j} \left\langle \nu^{m-r_N-j}
\smile \pi_{\fs}^*s_j(N) \smile \pi_{X,\fs}^*\alpha,
[\bar{\bB\bL}^{\vir}_{\ft,\fs}] \right\rangle.
\end{aligned}
\end{equation}
\end{prop}

\begin{proof}
See \cite[Section 5.2]{FLLevelOne}.
\end{proof}

\begin{rmk}
\label{rmk:GeneralSegreClasses} The Segre classes $s_i(N)$ have
been computed under some assumptions on $H^1(X;\RR)$ in
\cite[Lemma 4.11]{FL2a}. From \cite[Theorem 3.29]{FL2a}, one can
see that in general, these Segre classes will be given by a
universal polynomial in $\mu_{\fs}(x)$ and $\mu_{\fs}(\bga_i)$
with coefficients depending only on the indices $n_s'$ and $n_s''$
appearing in \eqref{eq:NormalComponentDims}
and which in turn only depend on $p_1(\ft(\ell))$ and $c_1(\ft(\ell))$.
\end{rmk}

\chapter{Computation of the intersection numbers}
\label{chap:Comp}

\section{Introduction}
\label{sec:Comp_intro}
In this chapter, we perform the computation leading to a proof of
the following theorems which form the technical heart of this monograph.

\begin{thm}
\label{thm:LinkPairing}
Let $X$ be a closed, connected, oriented, smooth, Riemannian four-manifold with $b_1(X)=0$
and $\ft$ be a \spinu structure on $X$. Suppose that $\ell\geq 0$ is an integer and $\fs$ is a \spinc structure on $X$ such that $M_{\fs}\times\Sym^\ell(X)$ is a subset of the space of
gauge-equivalence classes of ideal
monopoles $I\sM_{\ft}$ defined in \eqref{eq:idealmonopoles}.
Let $\bar\bL_{\ft,\fs}$ be the link of the
stratum of gauge-equivalence classes of
reducible $\SO(3)$ monopoles determined by $\fs\in\Spinc(X)$
as in Definition \ref{defn:DefineLink}.
For non-negative integers $m, \delta, \eta$ satisfying
$$
\delta-2m+2\eta=\dim\sM_{\ft}-2,
$$
a class $h \in H_2(X;\RR)$, and a generator
$x \in H_0(X;\ZZ)$,
let $z=h^{\delta-2m}x^m\in\AAA(X)$ and $\bar\sV(z)$ and $\bar\sW^{\eta}$ be the geometric
representatives defined in Section \ref{sec:Cohomology}.
Then
\begin{equation}
\label{eq:LinkIntersectionFormula}
\begin{aligned}
{}&
\#\left( \bar\sV(z)\cap\bar\sW^{\eta}\cap\bar\bL_{\ft,\fs}\right)
\\
&\quad
=
\SW_X(\fs)
\sum_{i=0}^{ \min(\ell,\left[\frac{\delta-2m}{2}\right])}
\left(
q_{\delta,\ell,m,i}(c_1(\fs)-c_1(\ft),c_1(\ft))Q_X^i
\right) (h),
\end{aligned}
\end{equation}
where $q_{\delta,\ell,m,i}$ are degree $\delta-2m-2i$
homogeneous polynomials which are universal functions
of the constants given in Theorem \ref{thm:MainThm}.
\end{thm}

We note that a similar result can also be achieved
by the methods of this monograph without
the assumptions that
$b_1(X)=0$ and $z=h^{\delta-2m}x^m$,
but the resulting expression \eqref{eq:LinkIntersectionFormula} becomes considerably more
complicated.
While
an explicit formula
for the intersection number in
\eqref{eq:LinkIntersectionFormula} is still unknown
in general, many special cases have been computed in \cite{FL2b,FLLevelOne,FL6, FL7, FL8} under additional hypotheses.
However,
the following important result,
referred to as the \emph{Multiplicity Conjecture} in
\cite{FL2b,FKLM, FLLevelOne}, holds without any
additional hypotheses. When $b_1(X)>0$, we write
$$
\SW_{X,\fs}(\om)
:=
\langle \mu_{\fs}(\om),[M_{\fs}]\rangle,
$$
for all $\om\in\AAA_2(X)$.

\begin{thm}
\label{thm:Multiplicity}
Assume the hypotheses of Theorem \ref{thm:LinkPairing}, but allow $b_1(X)>0$.
If $\SW_{X,\fs}(\om)$ vanishes for all $\om\in\AAA_2(X)$, then
$$
\#\left( \bar\sV(z)\cap\bar\sW^{\eta}\cap\bar\bL_{\ft,\fs}\right)=0.
$$
\end{thm}

To prove Theorems \ref{thm:LinkPairing} and \ref{thm:Multiplicity},
we begin by observing that
Proposition \ref{prop:Duality} gives the equality,
\begin{equation}
\label{eq:IntersectionToCohom1}
\#\left(\bar\sV(z)\cap\bar\sW^{\eta} \cap \bar\bL_{\ft,\fs}\right)
= \left\langle \barmu_p(z)\smile \barmu_c^{\eta}\smile \bar e_I \smile\bar e_s,
 [\bar\bL^{\vir}_{\ft,\fs}] \right\rangle.
\end{equation}
The expression for the $\mu$-classes in Corollary \ref{cor:CohomologyClasses},
the expression for $\bar e_s$ in Lemma \ref{lem:SWObstruction},
the expression for $\bar e_I$ in Proposition
\ref{prop:GlobalInstantonEulerClass}, and the equality between
the homology classes $[\bar\bL^{\vir}_{\ft,\fs}]$ and
$[\bar{\bB\bL}^{\vir}_{\ft,\fs}]$ in
Proposition \ref{prop:ReductionToBase} implies that the
intersection number in \eqref{eq:IntersectionToCohom1} can be
reduced to a linear combination of pairings of the form
\begin{equation}
\label{eq:ReducedPairing1}
\langle
\nu^a\smile\pi_X^*S^\ell(\beta)\smile \pi_{\fs}^*\mu_s(z_2),
[\bar{\bB\bL}^{\vir}_{\ft,\fs}]
\rangle
\end{equation}
where $\beta\in H_\bullet(X;\RR)$ and $z_2\in \AAA_2(X)$.

We wish to compute the pairings \eqref{eq:ReducedPairing1} by
writing them as sums over pairings with the subspaces
$\bar{\bB\bL}^{\vir}_{\ft,\fs}(\sP_i)$
defined in \eqref{eq:InstantonLinkComponentBase}
and then applying a pushforward-pullback
argument to the fiber-bundle
structure in Lemma \ref{lem:LocalBaseLinkFiberBundle} of each of these subspaces.
To do this, we replace the cohomology classes in
\eqref{eq:ReducedPairing1} by cohomology classes with compact support
away from the boundaries $\rd_j\bar{\bB\bL}^{\vir}_{\ft,\fs}(\sP_i)$ defined in
\eqref{eq:DefineBoundaryOne}.
Such a replacement requires a choice of cocycle in the cohomology class
and this choice must be done consistently on each subspace.
We encode these choices in the following geometric data.  We will define a quotient,
${\bQ\bL}^{\vir}_{\ft,\fs}$, of
$\bar{\bB\bL}^{\vir}_{\ft,\fs}$ by replacing the boundaries
$\rd_j\bar{\bB\bL}^{\vir}_{\ft,\fs}(\sP_i)$  with spaces
of codimension greater than or equal to two.
The cohomology classes in \eqref{eq:ReducedPairing1} pull back
from cohomology classes on this quotient.  Thus, we may chose cocycles
representing the cohomology classes in \eqref{eq:ReducedPairing1}
which pull back from the quotient.
Because the image of the boundaries in the quotient has
codimension greater than or equal to two, the cocycles pulled
back from the quotient will have compact support away
from the boundaries.  Using these cocycles is
equivalent to computing \eqref{eq:ReducedPairing1} using the
fundamental class of ${\bQ\bL}^{\vir}_{\ft,\fs}$
which, because the image of the boundaries has codimension two,
can be written as a sum of fundamental classes of subspaces.

One can compare this method of computation to the method
of ``adding caps''
employed by Ozsv{\'a}th in \cite{OzsvathBlowUp} as follows.
In the case where there are only two subspaces, if one
adds the cap to each subspace defined by the mapping
cone of the restriction of the quotient map to the
boundary of that subspace, the resulting compactification
of  each subspace $\bar{\bB\bL}^{\vir}_{\ft,\fs}(\sP_i)$ would be
homotopic to the image of $\bar{\bB\bL}^{\vir}_{\ft,\fs}(\sP_i)$
in the quotient.  Thus our method differs from that of
Ozsv{\'a}th in that there are more
than two open subspaces in this computation and in that there
are no correction terms arising from using different compactifications
of the subspaces.  Ozsv{\'a}th has informed us that he has an
extension of his method to the case of arbitrarily many open sets,
so our method offers no advantage in that regard.  Correction terms arise in
\cite{OzsvathBlowUp} because the natural compactifications to each
subspace define different quotients of the common boundary.
No correction terms appear in this computation because on
each boundary, $\rd_j\bar{\bB\bL}^{\vir}_{\ft,\fs}(\sP_i)$,
the natural quotient for the top stratum (defined by
extending the fiber bundle from $\Si_j$ to $\cl(\Si_j)$)
respects the fiber-bundle structure of the lower stratum.
Hence, the absence of correction terms in this computation
is not an indication of a better method but merely an
exploitation of a simpler situation.

In Section \ref{sec:Quotient}, we construct the quotient space
${\bQ\bL}^{\vir}_{\ft,\fs}$.  This construction shows
that the image of each subspace $\bar{\bB\bL}^{\vir}_{\ft,\fs}(\sP_i)$
in the quotient has a fiber-bundle structure identical to that
of $\bar{\bB\bL}^{\vir}_{\ft,\fs}(\sP_i)$ but with a compact
fiber and base.
In Section \ref{sec:HomOfQuot}, we show how to use the quotient
to obtain homology classes with the desired properties.
In Section \ref{sec:BundlesPushforwards}, we describe the
fiber-bundle properties of the images of the subspaces in the quotient.
In Section \ref{sec:Computations}, we perform the computation
giving the needed characterization of the intersection
pairing \eqref{eq:ReducedPairing1}.
Finally in Section \ref{sec:Proofs}, we give the proofs of
Theorems \ref{thm:LinkPairing} and \ref{thm:Multiplicity}.

\section{Quotient space of $\bar{\bB\bL}^{\vir}_{\ft,\fs}$}
\label{sec:Quotient}
We wish to use the fiber bundle structure,
\begin{equation}
\label{eq:FiberBundle}
\bar{\bB\bL}^{\vir}_{\ft,\fs}(\sP_j) \to M_{\fs}\times K_j
\subset
M_\fs\times \Si(X^\ell,\sP_j),
\end{equation}
given in Lemma \ref{lem:LocalBaseLinkFiberBundle}
to compute cohomological pairings with
$\bar{\bB\bL}^{\vir}_{\ft,\fs}$. To do this, we must write the
pairing in \eqref{eq:ReducedPairing1} with the homology class
$[\bar{\bB\bL}^{\vir}_{\ft,\fs}]$ as a sum over pairings with
homology classes representing the subspaces
$\bar{\bB\bL}^{\vir}_{\ft,\fs}(\sP_j)$ defined in \eqref{eq:InstantonLinkComponentBase}. However,  the subspaces
$\bar{\bB\bL}^{\vir}_{\ft,\fs}(\sP_j)$ have boundaries as
described in Section \ref{sec:Boundaries} and so would only define
relative homology classes. To overcome this difficulty, we construct a
quotient of $\bar{\bB\bL}^{\vir}_{\ft,\fs}$ obtained by leaving
the interiors of the subspaces $\bar{\bB\bL}^{\vir}_{\ft,\fs}(\sP_j)$
unchanged and replacing the boundaries
$\rd_k\bar{\bB\bL}^{\vir}_{\ft,\fs}(\sP_j)$ with quotients
obtained by deleting part of the gluing data as
described in the following

\begin{prop}
\label{prop:GlobalQuotient}
There is a surjective, continuous map,
$$
Q:\bar{\bB\bL}^{\vir}_{\ft,\fs} \to {\bQ\bL}^{\vir}_{\ft,\fs}
$$
onto a smoothly-stratified space ${\bQ\bL}^{\vir}_{\ft,\fs}$ with the
following properties:
\begin{enumerate}
\item
\label{item:GlobalQuotient_1}
The
map $Q$ is injective on the complement
of $\cup_{j\neq i}\rd_j\bar{\bB\bL}^{\vir}_{\ft,\fs}(\sP_i)$.
\item
\label{item:GlobalQuotient_2}
The image of each intersection,
$Q(\rd_j\bar{\bB\bL}^{\vir}_{\ft,\fs}(\sP_i))$, has codimension
greater than
two in ${\bQ\bL}^{\vir}_{\ft,\fs}$.
\item
\label{item:GlobalQuotient_3}
For the space $\barM(\sP_i,\beps_i)$ defined in \eqref{eq:DefineLinkFiber},
there is a continuous, surjective, $G(\sP_i)$-equivariant map, $q_i: \barM(\sP_i,\beps_i)\to \check M(\sP_i,\beps_i)$,
onto the smoothly-stratified space $\check M(\sP_i,\beps_i)$ given in Definition \ref{eq:DefineFiberQuotient}
such that
$$
Q(\bar{\bB\bL}^{\vir}_{\ft,\fs}(\sP_i))
\cong
M_{\fs}\times_{\sG_{\fs}\times S^1}
\bar\Fr(\ft,\fs,\sP_i)\times_{G(\sP_i)}\check M(\sP_i,\beps_i),
$$
where $\bar\Fr(\ft,\fs,\sP_i)$ is defined in
\eqref{eq:DefineExtendedBundle}
and the restriction of the map $Q$ to $\bar{\bB\bL}^{\vir}_{\ft,\fs}(\sP_i)$
is given, for $(A_0,\Phi_0)\in \tilde M_{\fs}$, and $F\in \Fr(\ft,\fs,\sP_i)$,
and $\vec A\in \bar M(\sP_i,\beps_i)$ by
$$
Q\left( \left[(A_0,\Phi_0),F,[\vec A]\right]\right)
:=
\left[(A_0,\Phi_0),R_i(F),q_i([\vec A])\right],
$$
where $R_i$ is the map defined in Lemma \ref{lem:FiberBundleQuotient}.
\item
\label{item:GlobalQuotient_4}
There is a
positive integer $b$
such that
the restriction of the bundle $\LL_\nu^{\otimes b}$,
where $\LL_\nu$ was defined in
\eqref{eq:LineBundleForS1ZAction}, to
$\bar{\bB\bL}^{\vir}_{\ft,\fs}$ is pulled back
from an $S^1$ bundle over ${\bQ\bL}^{\vir}_{\ft,\fs}$.
\item
\label{item:GlobalQuotient_5}
There is a map, $\pi_{Q,X}: {\bQ\bL}^{\vir}_{\ft,\fs}\to\Sym^\ell(X)$,
satisfying $\pi_X=\pi_{Q,X}\circ Q$, where $\pi_X$ is the
map defined in Lemma \ref{lem:GlobalProjectionToX}.
\item
\label{item:GlobalQuotient_6}
There is a map,
$$
\check\pi_{\fs}: {\bQ\bL}^{\vir}_{\ft,\fs}\to M_{\fs},
$$
such that $\check\pi_{\fs}\circ Q=\pi_{\fs}$, where
$\pi_{\fs}$ is the
map
defined in \eqref{eq:ProjectionToSW}.
\end{enumerate}
\end{prop}

Because the images of the
boundaries, $Q(\rd_j\bar{\bB\bL}^{\vir}_{\ft,\fs}(\sP_i))$,
have codimension greater than two,
the fundamental class of ${\bQ\bL}^{\vir}_{\ft,\fs}$
can be written as a sum of the fundamental classes
of the images of the subspaces.  This allows us to apply a
pushforward-pullback argument to the pairings of the
cohomology classes with these subspaces.

We will construct the quotient ${\bQ\bL}^{\vir}_{\ft,\fs}$
as follows.   Let
\begin{equation}
\label{eq:DefineExtendedBundle}
\bar\Fr(\ft,\fs,\sP_j)\to\Delta(X^\ell,\sP_j)
\end{equation}
be the extension of the fiber bundle
$\Fr(\ft,\fs,\sP_j)\to\Delta^\circ(X^\ell,\sP_j)$ defined in \eqref{eq:DefineGluingDataBundle}
from the open
diagonal $\Delta^\circ(X^\ell,\sP_j)$ to its closure,
$\Delta(X^\ell,\sP_j)$.
In Lemma \ref{lem:FiberBundleQuotient}, we define a surjective,
$G(\sP_j)$-equivariant map,
$$
R_j: \Fr(\ft,\fs,\sP_j)|_{K_j}\to\bar\Fr(\ft,\fs,\sP_j),
$$
which is injective over the complement of the boundaries $\cup_{i<j}\rd_i K_j$
defined in \eqref{eq:DefineBoundaryOfCompactum}.
We will define $R_j$  as a composition of maps
$R_{j,i}$, defined in Lemma \ref{lem:PartialTranslMap}.

Each of the maps $R_{j,i}$ defines a quotient $Q_{j,i}$ of the boundary
$\rd_i\bar{\bB\bL}^{\vir}_{\ft,\fs}(\sP_j)$.
In Lemma
\ref{lem:DescribeUpwardsBoundaryQuotient}, we show that for $j<k$
the quotient $Q_{k,j}$ is defined by a quotient of the fiber of
$\bar{\bB\bL}^{\vir}_{\ft,\fs}(\sP_j)\to M_{\fs}\times K_j$.
For $j<k<r$, the fibers of the quotient $Q_{r,j}$ contain
those of $Q_{k,j}$ in the corner
$\rd_r\rd_k \bar{\bB\bL}^{\vir}_{\ft,\fs}(\sP_j)$ and hence
one can `compose' these quotients.  The result of applying
the quotients $Q_{j+1,j},Q_{j+2,j},\dots Q_{n,j}$ successively
to $\bar{\bB\bL}^{\vir}_{\ft,\fs}(\sP_j)$ will be to replace
the fiber $\barM(\sP_j,\beps)$ with the quotient
$\check M(\sP_j,\beps_j)$ appearing in Proposition
\ref{prop:GlobalQuotient}.
For $i<j$, the quotient map $Q_{j,i}$ respects the fibers
of the bundle $\bar{\bB\bL}^{\vir}_{\ft,\fs}(\sP_j)\to M_{\fs}\times K_j$
and thus can be applied to the quotient obtained by replacing
the fiber $\barM(\sP_j,\beps)$ with the quotient
$\check M(\sP_j,\beps_j)$.
In Lemma \ref{lem:DownwardsQuotients1}, we show that one can
take the quotients $Q_{j,j-1}, Q_{j,j-2},\dots,Q_{j,0}$ successively
and obtain the quotient of $\bar{\bB\bL}^{\vir}_{\ft,\fs}(\sP_j)$
described in Proposition
\ref{prop:GlobalQuotient}.

Finally, we show that the quotient maps $Q_{j,i}$  on
$\rd_i\rd_j\bar{\bB\bL}^{\vir}_{\ft,\fs}(\sP_k)$
only identify points already identified
by quotient maps $Q_{k,r}$ or $Q_{r,k}$
with $r=i$ or $r=j$.
Thus, applying all the quotient maps $Q_{j,i}$ to the space
$\bar{\bB\bL}^{\vir}_{\ft,\fs}$ defines the quotient space
${\bQ\bL}^{\vir}_{\ft,\fs}$ described in the proposition.

\subsection{Quotient maps}
\label{subsec:Quotient_maps}
We construct the quotient maps $Q_{j,i}$
by defining a quotient of the frame bundle
$\Fr(\ft,\fs,\sP_k)|_{K_k}$ in \eqref{eq:DefineGluingDataBundle}.
As described in the Introduction to Section \ref{sec:Quotient},
this quotient will be given by the extension
\eqref{eq:DefineExtendedBundle}
of the fiber bundle
$\Fr(\ft,\fs,\sP)\to\Delta^\circ(X^\ell,\sP)$
to the closed diagonal, $\Delta(X^\ell,\sP)$.
Note that the structure group $G(\sP)/\Ga(\sP)$ does not act freely on the
extension $\bar\Fr(\ft,\fs,\sP)$
due to the presence of diagonals corresponding to cruder
partitions in $\Delta(X^\ell,\sP)$.  For this reason, we will
discuss $G(\sP)$-equivariant maps rather than bundle maps.

For partitions $\sP_i<\sP_k$ of $N_\ell$, let
$\nu(X^\ell,\sP_i\to [\sP_k])$ be the normal bundle of
$\Si(X^\ell,\sP_i)$ in $\Si(X^\ell,\sP_k)$ defined in
\eqref{eq:EndOfUpperStratum}.  Let $\sO(X^\ell,\sP_i\to [\sP_k],g_{\sP_i})\subset
\nu(X^\ell,\sP_i\to [\sP_k])$ be the neighborhood of $\Si(X^\ell,\sP_i)$
on which the exponential map $e(X^\ell,g_{\sP_i})$ of \eqref{eq:VaryingMetricExponentialMap} is defined.
Let $\sU(X^\ell,\sP_i\to [\sP_k],g_{\sP_i})\subset\Si(X^\ell,\sP_k)$ be the image of
$\sO(X^\ell,\sP_i\to [\sP_k],g_{\sP_i})$ under $e(X^\ell,g_{\sP_i})$.
For each $(A_0,\Phi_0)\in\tilde N_{\ft(\ell),\fs}(\delta)$, there is
a $G(\sP_k)$-equivariant homeomorphism,
\begin{equation}
\label{eq:EndBundleTriv}
\begin{CD}
\bar\Fr(\ft,\fs,\sP_k)\times_{\Si(X^\ell,\sP_i)}
\tilde\sO(X^\ell,\sP_i\to [\sP_k],g_{\sP_i})
\\
@V T(\sP_i,\sP_k)(A_0) VV
\\
\bar\Fr(\ft,\fs,\sP_k)|_{\sU(X^\ell,\sP_i\to [\sP_k],g_{\sP_i})}
\end{CD}
\end{equation}
defined by parallel translation of frames with respect to
the locally flattened connection $A'_0$, where $(A'_0,\Phi_0')=\Theta_{\sP_i}(A_0,\Phi_0)$,
and the locally flattened metric $g_{\sP_i,\bx}$, where the frames
in $\bar\Fr(\ft,\fs,\sP_k)$ lie over $\bx\in\Si(X^\ell,\sP_i)$.
This is the same
parallel translation used in the definition
of the upwards transition map in \eqref{eq:XUpwardsTransitionMap}.

\begin{lem}
\label{lem:PartialTranslMap}
Let $\sP_i<\sP_k$ be partitions of $N_\ell$.
Let $\vec t(X^\ell,g_{\sP_i})$ be the tubular distance function defined in
\eqref{eq:DefineDiagonalTubularDistFunction}
and let $\bar D(\sP_i,\eps_i)$ be the cube defined in \eqref{eq:DefineSquares}.
Then there is a $G(\sP_k)$-equivariant map,
$R_{k,i}:\bar\Fr(\ft,\fs,\sP_k)\to \bar\Fr(\ft,\fs,\sP_k)$,
with the following properties:
\begin{enumerate}
\item
\label{item:PartialTranslMap1}
$R_{k,i}$ is injective on the restriction of $\Fr(\ft,\fs,\sP_k)$
to the complement of the neighborhood
$\vec t(X^\ell,g_{\sP_i})^{-1}(\bar D(\sP_i,\eps_i))$.
\item
\label{item:PartialTranslMap2}
For $j<i$, $R_{k,i}$ maps the restriction of
$\bar\Fr(\ft,\fs,\sP_i)$ to $\vec t(X^\ell,g_{\sP_j})^{-1}(\bar D(\sP_j,\eps_j))$
to itself.
\item
\label{item:PartialTranslMap3}
If $K_{k,i}=K_k\cup (\cup_{j=k-1}^i T_{k,j})$,
where $T_{k,j}$ is defined in Lemma \ref{lem:CompactSubsetsOfSi}, then
$R_{k,i}$ maps $\bar \Fr(\ft,\fs,\sP_k)|_{K_{k,i+1}}$ onto
$\bar\Fr(\ft,\fs,\sP_k)|_{K_{k,i}}$.
\item
\label{item:PartialTranslMap4}
If $F_1,F_2\in \Fr(\ft,\fs,\sP_k)$, then
$R_{k,i}(F_1)=R_{k,i}(F_2)$ if and only if there are pairs
$(F_3,v_a)$, for $a=1,2$, in the domain of the map $T(\sP_i,\sP_k)$
defined in \eqref{eq:EndBundleTriv} with
$v_a\in \sO_1(X^\ell,\sP_i\to [\sP_k])$ and
$T(\sP_i,\sP_k)(F_3,v_a)=F_a$.
\item
\label{item:PartialTranslMap5}
$R_{k,i}$ covers the restriction of the map $r_i$
defined in Lemma \ref{lem:SymmetricProductRetractions} to $\Si(X^\ell,\sP_k)$.
\end{enumerate}
\end{lem}

\begin{proof}
For a partition $\sP_i$ of $N_\ell$ with $\sP_i<\sP_k$,
Item \eqref{item:SymmetricProductRetractions_3} of Lemma \ref{lem:SymmetricProductRetractions}
implies that
the map $r_i:\Sym^\ell(X)\to\Sym^\ell(X)$ pulls back,
by $e(X^\ell,g_{\sP_i})$, to a bundle map,
$$
r_{k,i}:
\sO(X^\ell,\sP_i\to [\sP_k],g_{\sP_i})
\to
\sO(X^\ell,\sP_i\to [\sP_k],g_{\sP_i}).
$$
The map $r_{k,i}$ and the $G(\sP_k)$-equivariant homeomorphism \eqref{eq:EndBundleTriv}
define a $G(\sP_k)$-equivariant map of
$\bar\Fr(\ft,\fs,\sP_k)|_{\sU(X^\ell,\sP_i\to [\sP_k],g_{\sP_i})}$ to itself.
Because $r_i$ is the identity on the complement of
the subspace $\sU(X^\ell,g_{\sP_i})$,
we can extend this map as the identity on the complement of
the restriction
$\bar \Fr(\ft,\fs,\sP_k)|_{\sU(X^\ell,g_{\sP_i})}$.
Let $R_{k,i}$ be this extension.  Item \eqref{item:PartialTranslMap5} follows immediately from
this definition.

Item \eqref{item:PartialTranslMap1} follows from this definition of $R_{k,i}$
and Item
\eqref{item:SymmetricProductRetractions_9}
of Lemma \ref{lem:SymmetricProductRetractions}.
Item \eqref{item:PartialTranslMap2} follows from Item \eqref{item:SymmetricProductRetractions_6} of Lemma \ref{lem:SymmetricProductRetractions}.
Item \eqref{item:PartialTranslMap3}  follows from Item \eqref{item:SymmetricProductRetractions_8} of Lemma \ref{lem:SymmetricProductRetractions}.

The $G(\sP_k)$-equivariance of $R_{k,i}$ implies that $R_{k,i}(F_1)=R_{k,i}(F_2)$
only if $F_1$ and $F_2$ lie over distinct points identified by $r_{k,i}$.
Item \eqref{item:PartialTranslMap4}  then follows from the construction of $R_{k,i}$.
\end{proof}

The following lemma  will be used to
compare the fibers of the quotient maps $Q_{k,i}$.

\begin{lem}
\label{lem:LowerEqualImpliesHigher}
Let $\sP_i<\sP_j<\sP_k$ be partitions of $N_\ell$.
Let $\rd_i\rd_j K_k$ be the corner of $K_k$ defined in \eqref{eq:DefineBoundaryOfCompactum}.
Let $R_{k,i}$ and $R_{k,j}$ be the $G(\sP_k)$-equivariant maps
defined in Lemma \ref{lem:PartialTranslMap}.
If $F_1,F_2\in \Fr(\ft,\fs,\sP_k)|_{\rd_i\rd_j K_k}$
and $R_{k,j}(F_1)=R_{k,j}(F_2)$, then
$R_{k,i}(F_1)=R_{k,i}(F_2)$.
\end{lem}

\begin{proof}
By Item \eqref{item:PartialTranslMap4}  of Lemma \ref{lem:PartialTranslMap},
and the assumption that $R_{k,j}(F_1)=R_{k,j}(F_2)$,
there are pairs $(F_3,v_a)$, for
$a=1,2$, in the domain of $T(\sP_j,\sP_k)(A_0)$
with $T(\sP_j,\sP_k)(A_0)(F_3,v_a)=F_a$.  The frame $F_3$
is itself a parallel translation of a frame
$F_4\in \bar\Fr(\ft,\fs,\sP_k)|_{\Si(X^\ell,\sP_i)}$.
Because the parallel translations defining these maps
are performed with respect
to a locally flattened connection and metric
which have no holonomy along these paths,
there are points
$w_a\in \tilde\nu(X^\ell,\sP_i \to [\sP_k])$ with
$T(\sP_i,\sP_k)(A_0)(F_4,w_a)=F_a$.  Hence, Item \eqref{item:PartialTranslMap3}
of Lemma \ref{lem:PartialTranslMap} implies that
$R_{k,i}(F_1)=R_{k,i}(F_2)$.
\end{proof}

\begin{lem}
\label{lem:FiberBundleQuotient}
Let $\sP_k$ be a partition of
$N_\ell$ and let $K^\circ_k$ be the interior of $K_k$. Then there is a
surjective, $G(\sP)$-equivariant map,
$$
R_k:
\Fr(\ft,\fs,\sP_k)|_{K_k}
\to
\bar\Fr(\ft,\fs,\sP_k),
$$
which is injective on $\Fr(\ft,\fs,\sP_k)|_{K_k^\circ}$.
For $F_1,F_2\in\Fr(\ft,\fs,\sP_k)|_{\rd_jK_k-\cup_{i<j}\rd_i K_k}$, we have
$R_k(F_1)=R_k(F_2)$ if and only if there are pairs
$(F_3,v_a)$, for $a=1,2$, in the domain of the map $T(\sP_j,\sP_k)$ defined in \eqref{eq:EndBundleTriv}
satisfying $T(\sP_j,\sP_k)(F_3,v_a)=F_a$.
In addition, $R_k$ covers the map $m_{k-1}$ defined in the proof of
Lemma \ref{lem:GlobalProjToSymCrit}.
\end{lem}

\begin{proof}
For the map $R_{k,i}$ defined in Lemma \ref{lem:PartialTranslMap},
define $R_0$ to be the identity map and for $k>0$, set
\begin{equation}
\label{eq:RComposition}
R_k:=R_{k,0}\circ R_{k,1}\circ \dots \circ R_{k,k-1}.
\end{equation}
Because $R_{k,i}$ covers $r_i$ by Item \eqref{item:PartialTranslMap5} of Lemma \ref{lem:PartialTranslMap}, the definition
\eqref{eq:RComposition} implies that $R_k$ covers $m_{k-1}=r_0\circ r_1\circ\dots\circ r_{k-1}$ as required.
Because it is the composition of $G(\sP_k)$-equivariant maps,
$R_k$ is $G(\sP_k)$-equivariant.
Let $T_{k,i}$ be as defined in Lemma \ref{lem:CompactSubsetsOfSi}.
Item \eqref{item:PartialTranslMap3} from Lemma \ref{lem:PartialTranslMap}
and the equality,
$$
\cl\left( \Si(X^\ell,\sP_k)\right)
=
K_k \cup \bigcup_{i<k} T_{k,i},
$$
from Item \eqref{item:CompactSubsetsOfSi_1} of Lemma \ref{lem:CompactSubsetsOfSi}
imply that $R_k$ maps
$\Fr(\ft,\fs,\sP_k)|_{K_k}$ onto $\bar\Fr(\ft,\fs,\sP_k)$.
It remains only to check the injectivity properties.

Let $F_1,F_2\in \Fr(\ft,\fs,\sP_k)|_{K_k}$ satisfy
$R_k(F_1)=R_k(F_2)$.  There is an index $u$  such that if
$\tilde F_a:=(R_{k,u+1}\circ \dots \circ R_{k,k-1})(F_a)$
for $a=1,2$, then $\tilde F_1\neq \tilde F_2$ while
$R_{k,u}(\tilde F_1)=R_{k,u}(\tilde F_2)$.
By Item \eqref{item:PartialTranslMap1} of  Lemma \ref{lem:PartialTranslMap}, $\tilde F_1$ and $\tilde F_2$
lie over $\vec t(X^\ell,g_{\sP_u})^{-1}(\bar D(\sP_u,\eps_u))$.
By Item \eqref{item:PartialTranslMap2} of Lemma \ref{lem:PartialTranslMap}, $F_1$ and $F_2$ also
lie over $\vec t(X^\ell,g_{\sP_u})^{-1}(\bar D(\sP_u,\eps_u))$.
Thus, $F_1$ and $F_2$ lie over
$K_k\cap \vec t(X^\ell,g_{\sP_u})^{-1}(\bar D(\sP_u,\eps_u))=\rd_u K_k$
and the restriction of $R_k$ to $\Fr(\ft,\fs,\sP_k)|_{K_k^\circ}$ is injective.

If $F_1$ and $F_2$ lie over $\rd_jK_k-\cup_{i<j}\rd_i K_k$, and $u$
and $\tilde F_a$ are as in the preceding paragraph, then because
the frames $F_a$ lie over $\rd_uK_k$, we must have $u\ge j$.
By Lemma \ref{lem:LowerEqualImpliesHigher},
the equality
$R_{k,u}(\tilde F_1)=R_{k,u}(\tilde F_2)$ would imply that
$R_{k,j}(\tilde F_1)=R_{k,j}(\tilde F_2)$.
By Item \eqref{item:PartialTranslMap4} of Lemma \ref{lem:PartialTranslMap},
there are pairs $(F_3,v_a)$ in the domain of $T(\sP_j,\sP_k)$
with $T(\sP_j,\sP_k)(F_3,v_a)=\tilde F_a$ for $a=1,2$. By construction of the maps
$R_{k,i}$, the frames $F_a$ are parallel translations along piecewise continuous
paths of the frames $\tilde F_a$.  Because the maps
$T(\sP_i,\sP_k)$ are defined by parallel translation with respect
to a locally-flattened connection and metric
which have no holonomy along these paths,
we therefore obtain pairs $(F_3,w_a)$ in the domain of
$T(\sP_j,\sP_k)$ with $T(\sP_j,\sP_k)(F_3,v_a)=F_a$ for $a=1,2$, completing the proof of the lemma.
\end{proof}

\begin{defn}
\label{defn:LocalQuotient}
For $\sP_i<\sP_k$ and the subspace $\rd_i \bar{\bB\bL}^{\vir}_{\ft,\fs}(\sP_k)$ defined in
\eqref{eq:InstantonLinkComponentBase},
define a quotient map $Q_{k,i}$ on $\bar{\bB\bL}^{\vir}_{\ft,\fs}(\sP_k)$
by identifying points,
$$
[(A_a,\Phi_a),F_a,[\vec A_a]]\in \rd_i \bar{\bB\bL}^{\vir}_{\ft,\fs}(\sP_k)
\quad\text{for $a=1,2$},
$$
where $(A_a,\Phi_a)\in\tilde M_{\fs}$, and $F_a\in\Fr(\ft,\fs,\sP_k)$,
and $[\vec A_a]\in\bar M(\sP_k,\beps)$ (the fiber of $\bL^{\vir}_{\ft,\fs}(\sP_k)$ defined in
\eqref{eq:DefineLinkFiber}),
if $R_{k,i}(F_1)=R_{k,i}(F_2)$,
where $R_{k,i}$ is the map defined
in Lemma \ref{lem:PartialTranslMap}.
\end{defn}

Therefore, by Item \eqref{item:PartialTranslMap1} of Lemma
\ref{lem:PartialTranslMap}, the quotient map $Q_{k,i}$
is injective on the complement of the boundary
$\rd_i\bar{\bB\bL}^{\vir}_{\ft,\fs}(\sP_k)$ and on that boundary
is defined by the projection map,
$$
\Fr(\ft,\fs,\sP_k)|_{\rd_iK_k} \to \bar \Fr(\ft,\fs,\sP_k)|_{K_i},
$$
defined by the parallel translation map in \eqref{eq:EndBundleTriv}
and the projection map,
$$
\bar\Fr(\ft,\fs,\sP_k)\times_{\Si(X^\ell,\sP_i)}
\tilde\sO(X^\ell,\sP_i\to [\sP_k],g_{\sP_i})
\to
\bar\Fr(\ft,\fs,\sP_k)|_{\Si(X^\ell,\sP_i)}.
$$
We note the following relation between quotient maps
$Q_{k,i}$ and $Q_{k,j}$.

\begin{lem}
\label{lem:DownwardsQuotients1}
Let $\sP_i<\sP_j<\sP_k$ be partitions of $N_\ell$.
If $\bA_1,\bA_2\in \rd_i\rd_j\bar{\bB\bL}^{\vir}_{\ft,\fs}(\sP_k)$,
the subspace defined in
\eqref{eq:InstantonLinkComponentBase},
and $Q_{k,j}(\bA_1)=Q_{k,j}(\bA_2)$ then  $Q_{k,i}(\bA_1)=Q_{k,j}(\bA_i)$.
\end{lem}

\begin{proof}
This follows immediately from the corresponding properties of
$R_{k,i}$ and $R_{k,j}$ described in Lemma \ref{lem:LowerEqualImpliesHigher}.
\end{proof}

Lemma \ref{lem:DownwardsQuotients1} implies that for each $k$, one can apply the
quotient maps $Q_{k,i}$ in order of descending $i$ as was done for the maps
$R_{k,i}$ in Equation \eqref{eq:RComposition}.  The resulting quotient map is $Q_k$.

\subsection{Construction of the local quotient}
\label{subsec:Constructing_local_quotient}
We now construct a quotient of the subspace
$\bar{\bB\bL}^{\vir}_{\ft,\fs}(\sP_i)$
using the maps $Q_{j,i}$.
We have described
the quotient map $Q_{j,i}$
on $\rd_i\bar{\bB\bL}^{\vir}_{\ft,\fs}(\sP_j)$ when $\sP_i<\sP_j$.
By \eqref{eq:MatchingBoundaries}, we have $\rd_j\bar{\bB\bL}^{\vir}_{\ft,\fs}(\sP_i)=\rd_i\bar{\bB\bL}^{\vir}_{\ft,\fs}(\sP_j)$
so $Q_{j,i}$ also defines a quotient of the boundary
$\rd_j\bar{\bB\bL}^{\vir}_{\ft,\fs}(\sP_i)$.
We begin by showing that for $\sP_i<\sP_j$,
the quotient map $Q_{j,i}$ on $\rd_j\bar{\bB\bL}^{\vir}_{\ft,\fs}(\sP_i)$ is defined
by a quotient of the fiber of
$\bar{\bB\bL}^{\vir}_{\ft,\fs}(\sP_i)\to M_{\fs}\times K_i$.

\begin{lem}
\label{lem:DescribeUpwardsBoundaryQuotient}
If $\sP_i<\sP_j$ are partitions of
$N_\ell$
and $\rd_j\barM(\sP_i,\beps_i)$ is the boundary defined in \eqref{eq:DefineLinkFiberBoundary},
then the restriction of the map
$Q_{j,i}$ to $\rd_j\bar{\bB\bL}^{\vir}_{\ft,\fs}(\sP_i)$
is given by the map,
$$
\begin{CD}
M_{\fs}\times_{\sG_{\fs}\times S^1}
\Fr(\ft,\fs,\sP_i)\times_{G(\sP_i)}\rd_j\barM(\sP_i,\beps_i)
\\
@VVV
\\
M_{\fs}\times_{\sG_{\fs}\times S^1}
\Fr(\ft,\fs,\sP_i)\times_{G(\sP_i)}
\bigsqcup_{\sP''\in [\sP_i<\sP_j]} \barM(\sP'')/\fS(\sP_i)
\end{CD}
$$
determined by the map $c_{j,i}$ in \eqref{eq:FiberBoundaryProjectionMapToConn}
on the fibers.
\end{lem}

\begin{proof}
The upwards overlap map,
$$
\begin{CD}
\tilde M_{\fs}\times_{\sG_{\fs}}
\Fr(\ft,\fs,\sP_i)\times_{G(\sP_i)}
\bigsqcup_{\sP''\in[\sP_i<\sP_j]}
\prod_{P\in\sP_j} \left( \Delta^\circ(Z_{|P|}(\delta_P),\sP''_P) \times \barM(\sP''_P)\right)
\\
@V \rho^{\ft,\fs,u}_{\sP_i,[\sP_i]} VV
\\
\left.
\left(\bigsqcup_{\sP''\in[\sP_i<\sP_j]}
\left(
\tilde M_{\fs}\times_{\sG_{\fs}}\Fr(\ft,\fs,\sP'')\times_{G(\sP'')} \barM(\sP'')
\right)\right)\right/\fS(\sP_i)
\\
@V \cong VV
\\
\tilde M_{\fs}\times_{\sG_{\fs}}\Fr(\ft,\fs,\sP_j)\times_{G(\sP_j)} \barM(\sP_j)
\end{CD}
$$
is defined in \eqref{eq:UpwardTransitionBasePoints} by the exponential map and parallel translation with
respect to the same connections used to define $R_{j,i}$.
Hence, pullback of the quotient $Q_{j,i}$ by the map
$\rho^{\ft,\fs,u}_{\sP_i,[\sP_i]}$ is defined by the map on the fibers of the
domain of $\rho^{\ft,\fs,u}_{\sP_i,[\sP_i]}$ given by the projection,
$$
\begin{CD}
\bigsqcup_{\sP''\in[\sP_i<\sP_j]}
\prod_{P\in\sP_j} \left( \Delta^\circ(Z_{|P|}(\delta_P),\sP''_P) \times \barM(\sP''_P)\right)
\\
@VVV
\\
\bigsqcup_{\sP''\in[\sP_i<\sP_j]}
\prod_{P\in\sP_j}   \barM(\sP''_P)
\\
@V = VV
\\
\bigsqcup_{\sP''\in[\sP_i<\sP_j]}
\barM(\sP'')
\end{CD}
$$
The overlap map $\rho^{\ft,\fs,d}_{\sP_i,[\sP_i]}$ is defined by the $G(\sP_i)$-equivariant
map  $\rho^{\ft,\fs,d}_{f,\sP_i,[\sP_j]}$
of fibers from \eqref{eq:DownwardsInclusionFiberMap}
and thus pushes the preceding projection forward to the projection
described in the statement of the lemma.
\end{proof}

To compare the quotient maps $Q_{i,k}$ and $Q_{j,k}$,
for $k<i<j$, Lemma \ref{lem:DescribeUpwardsBoundaryQuotient}
shows that
it suffices to compare the maps $c_{i,k}$ and $c_{j,k}$ and
which we now do.

\begin{lem}
\label{lem:MultiOverlapQuotientOfFiber}
Let $\sP_k<\sP_i<\sP_j$ be partitions of $N_\ell$.
If $[\vec A_1],[\vec A_2]\in \rd_i\rd_j \barM(\sP_k,\beps)$, the
intersection defined in
\eqref{eq:DefineLinkFiberBoundary},
satisfy $c_{i,k}([\vec A_1])=c_{i,k}([\vec A_2])$, then
$c_{j,k}([\vec A_1])=c_{j,k}([\vec A_2])$.
\end{lem}

\begin{proof}
The assumption that $[\vec A_a]\in \rd_i\rd_j \barM(\sP_k,\beps)$ for $a=1,2$
implies that $[\vec A_a]$ lies in the images of the maps of the fibers
$\rho^{\ft,\fs,d}_{f,\sP_k,[\sP_i]}$ and
$\rho^{\ft,\fs,d}_{f,\sP_k,[\sP_j]}$
defined in \eqref{eq:DowardsOverlapFiberMap}.
Thus, by the definition of
the maps $\rho^{\ft,\fs,d}_{f,\sP_k,[\sP_i]}$ and
$\rho^{\ft,\fs,d}_{f,\sP_k,[\sP_j]}$,
if we write $[\vec A_a]=([A_{P,a}])_{P\in\sP_k}$, where
$[A_{P,a}]\in \barM^{s,\natural}_{\spl,|P|}(\delta_P)$, then
$$
[A_{P,a}] \in\Imag(\bga'_{\Theta,\sP'_P}) \cap \Imag(\bga'_{\Theta,\sP''_P})
\quad\text{for all $P\in\sP_k$,}
$$
for some $\sP'\in [\sP_k<\sP_i]$ and $\sP''\in [\sP_k<\sP_j]$,
where $\bga'_{\Theta,\sP'_P}$ is the splicing map defined in \eqref{eq:SplicedGenConn}.
Because the preceding intersection is non-empty, we have $\sP_k<\sP'<\sP''$.
Lemma \ref{lem:ConstructCollar}
and the commutativity of the diagram \eqref{eq:CommutingSplicingQuotient}
imply that
$$
A_{P,a}
=
\left(\bga'_{\Theta,\sP'_P}\circ \rho^{\Theta,d}_{\sP'_P,[\sP''_P]}\right)(B_P)
=
\left(\bga'_{\Theta,\sP''_P}\circ \rho^{\Theta,u}_{\sP'_P,[\sP''_P]}\right)(B_P),
$$
for some $B_P$ in the following overlap space defined in \eqref{eq:SplicedModuliOverlapR4} (omitting the symmetric group
as it is not relevant to this discussion),
\begin{align*}
{}&\tilde\sO^{\asd} (\Theta,\sP'_P,[\sP''_P],\delta)
\\
{}&\quad
\subseteqq
\Delta^\circ(Z_P(\delta_P),\sP_P')
\times
\prod_{Q'\in\sP'_P}
\left(
\Delta^\circ( Z_{Q'}(\delta_{Q'}),\sP''_{Q'})
\times
\prod_{Q''\in\sP''_{Q'}}\barM^{s,\natural}_{\spl,|Q''|}(\delta_{Q''})
\right).
\end{align*}
The map,
$$
c_{j,k}
\circ
\left(\prod_{P\in\sP_k}\bga'_{\Theta,\sP''_P}\circ \rho^{\Theta,u}_{\sP'_P,[\sP''_P]}
\right),
$$
is given by the projection,
$$
\begin{CD}
\prod_{P\in\sP_k}
\left(\Delta^\circ(Z_P(\delta_P),\sP_P')
\times
\prod_{Q'\in\sP'_P}
\left(
\Delta^\circ( Z_{Q'}(\delta_{Q'}),\sP''_{Q'})
\times
\prod_{Q''\in\sP''_{Q'}}\barM^{s,\natural}_{\spl,|Q''|}(\delta_{Q''})
\right)
\right)
\\
@VVV
\\
\prod_{P\in\sP_k}
\prod_{Q'\in\sP'_P}
\prod_{Q''\in\sP''_{Q'}}\barM^{s,\natural}_{\spl,|Q''|}(\delta_{Q''})
\end{CD}
$$
Moreover, the map,
$$
c_{i,k}\circ
\left(\prod_{P\in\sP_k}
\bga'_{\Theta,\sP'_P}\circ \rho^{\Theta,d}_{\sP'_P,[\sP''_P]}
\right),
$$
is given by the projection,
$$
\begin{CD}
\prod_{P\in\sP_k}
\left(
\Delta^\circ(Z_P(\delta_P),\sP_P')
\times
\prod_{Q'\in\sP'_P}
\left(
\Delta^\circ( Z_{Q'}(\delta_{Q'}),\sP''_{Q'})
\times
\prod_{Q''\in\sP''_{Q'}}\barM^{s,\natural}_{\spl,|Q''|}(\delta_{Q''})
\right)
\right)
\\
@VVV
\\
\prod_{P\in\sP_k}
\prod_{Q'\in\sP'_P}
\left(
\Delta^\circ( Z_{Q'}(\delta_{Q'}),\sP''_{Q'})
\times
\prod_{Q''\in\sP''_{Q'}}\barM^{s,\natural}_{\spl,|Q''|}(\delta_{Q''})
\right)
\end{CD}
$$
The desired equality, $c_{j,k}(A_1)=c_{j,k}(A_2)$,
follows by comparing the preceding projections
and using the commutativity of the diagram \eqref{eq:CommutingSplicingQuotient}.
\end{proof}

Lemma \ref{lem:MultiOverlapQuotientOfFiber}
allows us to make the following

\begin{defn}
\label{eq:DefineFiberQuotient}
Define $\check M(\sP_k,\beps)$ to be the quotient of
the fiber $\barM(\sP_k,\beps)$ of $\bar\bL^{\vir}_{\ft,\fs}(\sP_k)$ defined in \eqref{eq:DefineLinkFiber}
by applying the quotient maps
$c_{k+1,k},c_{k+2,k},\dots,c_{n,k}$, defined in \eqref{eq:FiberBoundaryProjectionMapToConn}, successively.
Let $q_k: \barM(\sP_k,\beps)\to \check M(\sP_k,\beps)$ denote the resulting quotient map.
\end{defn}

\begin{lem}
\label{lem:CharacterizeFiberQuotient}
Let $q_k: \barM(\sP_k,\beps)\to \check M(\sP_k,\beps)$
be the quotient map in Definition \ref{eq:DefineFiberQuotient}.
Then $\check M(\sP_k,\beps)$ is a smoothly-stratified
space, with a smoothly-stratified $G(\sP_k)$ action satisfying the following properties:
\begin{enumerate}
\item
\label{item:CharacterizeFiberQuotient_1}
For all $i>k$, the image $q_k(\rd_i \barM(\sP_k,\beps))$ has
codimension greater than or equal to two in
$\check M(\sP_k,\beps)$.
\item
\label{item:CharacterizeFiberQuotient_2}
The action of $S^1<\SO(3)$ on
$$
\barM(\sP_k,\beps) \subseteqq
\prod_{P\in\sP_k} \barM^{s,\natural}_{\spl,|P|}(\delta_P)
$$
given by \eqref{eq:DiagonalFrameAction}
pulls back from an action of $\SO(3)$ on
$\check M(\sP_k,\beps)$ with finite isotropy
as defined prior to Lemma \ref{lem:FreeActionOnCijImage}.
\item
\label{item:CharacterizeFiberQuotient_3}
The map $q_k$ is $G(\sP_k)$-equivariant.
\end{enumerate}
\end{lem}

\begin{proof}
Recall from \eqref{eq:DefineLinkFiber} that $\bar M(\sP_k,\beps)$ is the subspace of
$\bar M(\sP_k)$ where the maps $\vec t_f(\sP_k,[\sP_j])$ with $k\le j$ take values
in smoothly-stratified subspaces of $[0,1]^{|\sP_j|}$.
Hence,  the smooth strata of $\bar M(\sP_k,\beps)$ are given by taking the intersection
of smooth strata of $\bar M(\sP_k)$ with the pre-image
under $\vec t_f(\sP_k,[\sP_j])$
of these smooth strata of $[0,1]^{|\sP_j|}$.
The quotients given by the maps $c_{j,k}$ are defined by replacing the union of
strata $\rd_j \bar M(\sP_k,\beps)$ with the smoothly-stratified space given
by the image of $c_{j,k}$.  Hence, $\check M(\sP_k,\beps)$ is a smoothly-stratified space.
The maps $\vec t_f(\sP_k,[\sP_j])$ are invariant under the smooth, stratum-preserving
$G(\sP_k)$ action on $\bar M(\sP_k)$, so
the $G(\sP_k)$ action also preserves the pre-images  of the strata of
$[0,1]^{|\sP_j|}$ under the map $\vec t_f(\sP_k,[\sP_j])$.
The quotient maps $c_{j,k}$ are $G(\sP_k)$-equivariant and so the action of $G(\sP_k)$
on $\bar M(\sP_k,\beps)$ descends to define an action on $\check M(\sP_k,\beps)$.
The equivariance of $q_k$ in Item \eqref{item:CharacterizeFiberQuotient_3} follows
immediately from this construction.

The assertion regarding the
codimensions of the images of the boundaries
follows immediately from the construction.
The assertion regarding the $\SO(3)$ action follows
from the equivariance of each map $c_{i,k}$
with respect to this action and from the
fact that $\SO(3)$ acts with finite isotropy on the image
of each map $c_{i,k}$ by Lemma \ref{lem:FreeActionOnCijImage}.
That is, if
$e^{i\theta}\in S^1<\SO(3)$ and $[[\vec A]] \in \check M(\sP_k,\beps)$ is the image of $[\vec A]\in\barM(\sP_k,\beps)$
with $e^{i\theta}[[\vec A]]=[e^{i\theta}[\vec A]]=[[\vec A]]$,
then there is a smallest integer $j$ such that $c_{j,k}(e^{i\theta}[\vec A])=c_{j,k}([\vec A])$.
Lemma \ref{eq:PseudoFreeActions} implies that $e^{i\theta}$ is therefore
a root of unity, of order $o_j$ determined by $j$.  Hence, if
$[\vec A]\in\barM(\sP_k,\beps)$ and $e^{i\theta}\in S^1<\SO(3)$ satisfy
$e^{i\theta}[[\vec A]]=[[\vec A]]$, then
$e^{i\theta}$ must be a root of unity
whose order divides the product $o_{k+1}o_{k+2}\dots o_n$ of those orders.
\end{proof}

By Lemma \ref{lem:DescribeUpwardsBoundaryQuotient},
the result of applying the quotient maps
$$
Q_{k+1,k},Q_{k+2,k},\dots,Q_{n,k}
$$ to
$\bar{\bB\bL}^{\vir}_{\ft,\fs}(\sP_k)$ is the same
as that of applying the quotient maps
$$
c_{k+1,k},c_{k+2,k},\dots,c_{n,k}
$$ 
to
the fiber $\barM(\sP_k,\beps)$.
From Lemma \ref{lem:CharacterizeFiberQuotient}, we see
that the quotient of $\bar{\bB\bL}^{\vir}_{\ft,\fs}(\sP_k)$ by the quotient maps
$Q_{k+1,k},Q_{k+2,k},\dots,Q_{n,k}$ is equal to
\begin{equation}
\label{eq:PartialQuotient}
\widehat{\bB\bL}^{\vir}_{\ft,\fs}(\sP_k)
:=
\tilde M_{\fs}\times_{\sG_{\fs}\times S^1}
\Fr(\ft,\fs,\sP_k)|_{K_k}\times_{G(\sP_k)}\check M(\sP_k,\beps).
\end{equation}
For integers $j<k$, the quotient maps $Q_{k,j}$ preserve the fibers
of $\bar{\bB\bL}^{\vir}_{\ft,\fs}(\sP_k)\to M_{\fs}\times K_k$
and hence define quotients of $\widehat{\bB\bL}^{\vir}_{\ft,\fs}(\sP_k)$.
We characterize the result of taking these further quotients as
follows.

\begin{lem}
\label{lem:QuotientOfPieces}
Let $\sP_k$ be a partition of $N_\ell$.  Then
after applying the quotient maps
$$
Q_{k+1,k},Q_{k+2,k},\dots,Q_{n,k},
Q_{k,k-1},Q_{k,k-2},\dots,Q_{k,1}
$$
successively, we obtain the quotient map,
$$
\begin{CD}
\tilde M_{\fs}\times_{\sG_{\fs}\times S^1}
\Fr(\ft,\fs,\sP_k)|_{K_k}\times_{G(\sP_k)}\barM(\sP_k,\beps)
\\
@V Q_k VV
\\
\tilde M_{\fs}\times_{\sG_{\fs}\times S^1}
\bar\Fr(\ft,\fs,\sP_k)\times_{G(\sP_k)}\check M(\sP_k,\beps)
\end{CD}
$$
defined for $(A_0,\Phi_0)\in\tilde M_{\fs}$,
and $F\in\Fr(\ft,\fs,\sP_k)|_{K_k}$, and
$[\vec A ]\in \barM(\sP_k,\beps)$
by
$$
Q_k\left([(A_0,\Phi_0),F,[\vec A ]]\right)
:= \left[(A_0,\Phi_0), R_k(F),q_k([\vec A ])\right],
$$
where $R_k$ is defined in Lemma \ref{lem:FiberBundleQuotient}
and $q_k$ is defined in Lemma \ref{lem:CharacterizeFiberQuotient}.
\end{lem}

\begin{proof}
By applying the quotient maps
$Q_{k+1,k},Q_{k+2,k},\dots,Q_{n,k}$,
we obtain the space
$\widehat{\bB\bL}^{\vir}_{\ft,\fs}(\sP_k)$ as described
above in \eqref{eq:PartialQuotient}.
From the definition of the quotient maps $Q_{k,j}$ in Definition \ref{defn:LocalQuotient},
the result of applying the quotient maps $Q_{k,k-1},Q_{k,k-2},\dots,Q_{k,1}$
to $\widehat{\bB\bL}^{\vir}_{\ft,\fs}(\sP_k)$ is the same
as that of applying the quotient maps
$R_{k,k-1},R_{k,k-2},\dots,R_{k,1}$ to
$\Fr(\ft,\fs,\sP_k)|_{K_k}$.
By Lemma \ref{lem:FiberBundleQuotient},
applying the quotient maps
$R_{k,k-1},R_{k,k-2},\dots,R_{k,1}$ to
$\Fr(\ft,\fs,\sP_k)|_{K_k}$
gives the quotient
of $\Fr(\ft,\fs,\sP_k)|_{K_k}$
defined by the
surjective map, $R_k:\Fr(\ft,\fs,\sP_k)|_{K_k}\to\bar\Fr(\ft,\fs,\sP_k)$.
\end{proof}

We note the following

\begin{lem}
\label{lem:ProjectionToSigFromQuotientPiece}
Continue the notation of Lemma \ref{lem:QuotientOfPieces}.
The pullback by $Q_k$ of the projection,
\begin{equation}
\label{eq:LocalQuotientProjToX}
\check\pi_{X,k}:
\tilde M_{\fs}\times_{\sG_{\fs}\times S^1}
\bar\Fr(\ft,\fs,\sP_k)\times_{G(\sP_k)}\check M(\sP_k,\beps),
\to
\Si(X^\ell,\sP_k)
\end{equation}
is equal to the restriction of the map $\pi_X$ to $\bar{\bB\bL}^{\vir}_{\ft,\fs}(\sP_k)$.
\end{lem}

\begin{proof}
The conclusion follows immediately from the assertion that
$R_k$ covers the map $m_{k-1}$ in Lemma \ref{lem:FiberBundleQuotient}
and the construction of the map $\pi_X$ in
Lemmas  \ref{lem:GlobalProjToSymCrit} and \ref{lem:GlobalProjectionToX}.
\end{proof}

\subsection{Global quotient of $\bar{\bB\bL}^{\vir}_{\ft,\fs}$}
\label{subsec:Global_quotient}
Lemma \ref{lem:QuotientOfPieces} describes a quotient
of each subspace $\bar{\bB\bL}^{\vir}_{\ft,\fs}(\sP_j)$ of
$\bar{\bB\bL}^{\vir}_{\ft,\fs}$.  We have used Lemma
\ref{lem:DescribeUpwardsBoundaryQuotient} to ensure
that we are applying the same quotient to both
boundaries in the equality
\eqref{eq:MatchingBoundaries},
$$
\rd_i\bar{\bB\bL}^{\vir}_{\ft,\fs}(\sP_j)
=
\rd_j\bar{\bB\bL}^{\vir}_{\ft,\fs}(\sP_i).
$$
We will define the global quotient appearing in the statement of
Proposition \ref{prop:GlobalQuotient} by applying
all the quotient maps $Q_{j,i}$ in an appropriate order.
To ensure that the resulting global quotient satisfies Proposition \ref{prop:GlobalQuotient}, we
need
the following
lemma
to show
that no
identifications, other than those already
in the statement of
Lemma \ref{lem:QuotientOfPieces}, appear in the global
quotient.

\begin{lem}
\label{lem:GlobalQuotientNotWorse}
For $i<j$, if
$\bA_1,\bA_2 \in \rd_i\rd_j\bar{\bB\bL}^{\vir}_{\ft,\fs}(\sP_r)$
satisfy $Q_{j,i}(\bA_1)=Q_{j,i}(\bA_2)$, then the following hold:
\begin{enumerate}
\item
\label{item:GlobalQuotientNotWorse_1}
If $i<j<r$, then $Q_{r,i}(\bA_1)=Q_{r,i}(\bA_2)$.
\item
\label{item:GlobalQuotientNotWorse_2}
If $r<i<j$, then $Q_{j,r}(\bA_1)=Q_{j,r}(\bA_2)$.
\item
\label{item:GlobalQuotientNotWorse_3}
If $i<r<j$, then there is an element
$\bA'\in \rd_i\rd_j \bar{\bB\bL}^{\vir}_{\ft,\fs}(\sP_r)$
such that
$Q_{r,i}(\bA_1)=Q_{r,i}(\bA')$,
and
$Q_{j,r}(\bA_1)=Q_{j,r}(\bA')$.
\end{enumerate}
\end{lem}

\begin{proof}
Items \eqref{item:GlobalQuotientNotWorse_1} and \eqref{item:GlobalQuotientNotWorse_2} follow immediately from
Lemmas \ref{lem:MultiOverlapQuotientOfFiber}
and \ref{lem:DownwardsQuotients1}, respectively.
We now prove Item \eqref{item:GlobalQuotientNotWorse_3}.
Because $\rd_i\rd_j\bar{\bB\bL}^{\vir}_{\ft,\fs}(\sP_r)=\rd_i\rd_r\bar{\bB\bL}^{\vir}_{\ft,\fs}(\sP_j)$
by \eqref{eq:MatchingBoundaries}, we can consider the points $\bA_1$ and $\bA_2$ to
be elements of $\bar{\bB\bL}^{\vir}_{\ft,\fs}(\sP_j)$.
By Definition \ref{defn:LocalQuotient},
points
$\bA_1,\bA_2\in \rd_i \bar{\bB\bL}^{\vir}_{\ft,\fs}(\sP_j)$
are identified by $Q_{j,i}$ if they can be expressed as
$$
\bA_a=\left[(A_0,\Phi_0),T(\sP_i,\sP_j)(A_0)(F_3,v_a),([A_Q])_{Q\in\sP_j}\right],
\quad a=1,2,
$$
where $(A_0,\Phi_0)\in\tilde M_{\fs}$,
and $(F_3,v_a)$ is in the domain of $T(\sP_i,\sP_j)$ given in
\eqref{eq:EndBundleTriv}, and
$([A_Q])_{Q\in\sP_j}\in \barM(\sP_j,\beps)$.
Note that the frames $T(\sP_i,\sP_j)(A_0)(F_3,v_a)$ lie over points in $\rd_i\rd_r K_j$.
Since
\begin{equation}
\label{eq:FrameBundleInclusion}
\bar\Fr(\ft,\fs,\sP_j)|_{\Si(X^\ell,\sP_i)}
\cong
\Fr(\ft,\fs,\sP_i)\times_{\tilde G(\sP_i)}\tilde G(\sP_j),
\end{equation}
we can assume that $F_3\in \bar\Fr(\ft,\fs,\sP_i)$.
Because the parallel translations used to define the maps $T(\sP_i,\sP_j)(A_0)$ are performed
with the locally flattened connection, $A'_0$, we can write
the parallel translation of $F_3$ given by $T(\sP_i,\sP_j)(A_0)(F_3,v_a)$
as a composition of parallel translations, first the
parallel translation $T(\sP_i,\sP_r)(A_0)(\cdot,u_a)$ to a point
in $\rd_i K_r$, and then
the parallel translation
$T(\sP_r,\sP_j)(A_0)(\cdot,w_a)$ from
that point to the point in $\rd_i\rd_r K_j$:
\begin{equation}
\label{eq:CompositionOfParallelTranslation1}
T(\sP_i,\sP_j)(A_0)(F_3,v_a)
=
T(\sP_j,\sP_r)(A_0)\left(T(\sP_i,\sP_r)(A_0)\left( F_3,v_a'\right),w_a\right).
\end{equation}
By the identification in \eqref{eq:IdentifyingNormalBundleOfDiagonalInSymmProduct1},
we can write the elements $(F_3,v_a')$ of the domain of $T(\sP_i,\sP_r)(A_0)$ as
$(F_3, (v'_{P,a})_{P\in\sP''})$, where $\sP''\in [\sP_r<\sP_j]$
and $v'_{P,a}\in \Delta^\circ(Z_{|P|}(\delta_P),\sP''_P)$.
Then the definition of the upwards transition map in \eqref{eq:UnparamUpwardsTransition} implies that
$\bA_a=\rho^{\ft,\fs,u}_{\sP_r,[\sP_j]}(\bA^r_a)$, where
$$
\bA^r_a=
\left[(A_0,\Phi_0),T(\sP_i,\sP_r)(A_0)(F_3,w_a),(v'_{P,a},(A_Q)_{Q\in\sP''_P})_{P\in\sP_r}\right],
$$
and $T(\sP_i,\sP_r)(A_0)(F_3,w_a)\in\Fr(\ft,\fs,\sP_r)$, and $\sP''\in [\sP_r<\sP_j]$,
$v'_{P,a}$ and $A_Q$ are as above.

If we define
$$
\bA':=
\rho^{\ft,\fs,u}_{\sP_r,[\sP_j]}
\left(
\left[(A_0,\Phi_0),
       T(\sP_i,\sP_r)(A_0)(F_3,w_2),
        (v'_{P,1},(A_Q)_{Q\in\sP''_P})_{P\in\sP_r}
\right]
\right),
$$
then $Q_{j,r}(\bA_2)=Q_{j,r}(\bA')$.
The latter equality holds because
if we write $F_3'=T(\sP_i,\sP_r)(A_0)(F_3,w_a)$, then by \eqref{eq:CompositionOfParallelTranslation1}
we can express $\bA_2$ and $\bA'$ as
\begin{align*}
\bA_2&=\left[ (A_0,\Phi_0), T(\sP_j,\sP_r)(A_0)(F_3',v'_{P,2}),(A_Q)_{Q\in\sP_j}\right],
\\
\bA'&=\left[ (A_0,\Phi_0), T(\sP_j,\sP_r)(A_0)(F_3',v'_{P,1}),(A_Q)_{Q\in\sP_j}\right].
\end{align*}
To see that $Q_{r,i}(\bA_1)=Q_{r,i}(\bA')$, we apply the downwards overlap map $\rho^{\ft,\fs,d}_{\sP_r,[\sP_j]}$
defined in \eqref{eq:XDownwardTransition},
so that
\begin{align*}
\bA_1
{}&=
\rho^{\ft,\fs,d}_{\sP_r,[\sP_j]}
\left(
\left[(A_0,\Phi_0),T(\sP_i,\sP_r)(A_0)(F_3,w_1),(v'_{P,1},(A_Q)_{Q\in\sP''_P})_{P\in\sP_r}\right]
\right)
\\
{}&=
\left[(A_0,\Phi_0),T(\sP_i,\sP_r)(A_0)(F_3,w_1),(A'_P)_{P\in\sP_r}\right],
\\
\bA'
{}&=
\rho^{\ft,\fs,d}_{\sP_r,[\sP_j]}
\left(
\left[(A_0,\Phi_0),
       T(\sP_i,\sP_r)(A_0)(F_3,w_2),
        (v'_{P,1},(A_Q)_{Q\in\sP''_P})_{P\in\sP_r}
\right]
\right)
\\
{}&=
\left[(A_0,\Phi_0),T(\sP_i,\sP_r)(A_0)(F_3,w_2),(A'_P)_{P\in\sP_r}\right],
\end{align*}
where, for $P\in\sP_r$ and $\sP''\in [\sP_r<\sP_j]$, we denote
$$
A'_P := \bga'_{\Theta,\sP''_P}\left((v'_{P,1},(A_Q)_{Q\in\sP''_P}) \right).
$$
Because the preceding expressions for $\bA_1$ and $\bA'$ differ only in the ``$w_a$'' coordinates,
the desired equality $Q_{r,i}(\bA_1)=Q_{r,i}(\bA')$ follows.
\end{proof}

We can now give the

\begin{proof}[Proof of Proposition \ref{prop:GlobalQuotient}]
Let $\sP_0,\dots,\sP_n$ be the partitions of $N_\ell$ chosen, as in \S \ref{subsec:EnumStrata} ,
to enumerate the strata of $\Sym^\ell(X)$.
We perform the quotients $Q_{j,i}$ in the following order:
\begin{align*}
{}& Q_{n,n-1},
\\
{}& Q_{n-1,n-2},\ Q_{n,n-2},\
\\
{}&\ \vdots
\\
{}& Q_{j+1,j},\ Q_{j+2,j},\dots, Q_{n,j},
\\
{}&\ \vdots
\\
{}& Q_{1,0},\ Q_{2,0},\dots, Q_{n,0}.
\end{align*}
Lemma \ref{lem:GlobalQuotientNotWorse} implies that
the application of these
quotients to the subspace $\bar\bB\bL^{\vir}_{\ft,\fs}(\sP_r)$
yields the same quotient as that defined by Lemma \ref{lem:QuotientOfPieces}.
This establishes Items \eqref{item:GlobalQuotient_1}--\eqref{item:GlobalQuotient_3}.

The restriction of the bundle $\LL_\nu$ to the
subspace
$\bar{\bB\bL}^{\vir}_{\ft,\fs}(\sP_r)$ is the complex line
bundle associated to the $S^1$ bundle,
\begin{equation}
\label{eq:RestrictionOfS1BundleToPiece}
\begin{CD}
\tilde M_{\fs}\times_{\sG_{\fs}}\Fr(\ft,\fs,\sP_r)|_{K_r}
\times_{G(\sP_r)}\bar M(\sP_r,\beps)
\\
@VVV
\\
\tilde M_{\fs}\times_{\sG_{\fs}\times S^1}\Fr(\ft,\fs,\sP_r)|_{K_r}
\times_{G(\sP_r)\times S^1}\bar M(\sP_r,\beps)
\end{CD}
\end{equation}
where the $S^1$ action is defined in \eqref{eq:S1ActionOnPiece}.
This $S^1$ action is trivial on $\tilde M_{\fs}$ and acts diagonally
on the frames
in $\Fr(\ft,\fs,\sP_r)$ as described in
Section \ref{subsubsec:GrpActionOnFrameBundle}.
Because of the action of $G(\sP_r)$ on the base of \eqref{eq:RestrictionOfS1BundleToPiece},
this $S^1$ action equals
the diagonal action of $S^1$ on the frames in $\bar M(\sP_r,\beps)$ as described in
\eqref{eq:DiagonalFrameAction}.
By Item \eqref{item:CharacterizeFiberQuotient_2} of Lemma \ref{lem:CharacterizeFiberQuotient}, this
action on $\bar M(\sP_r,\beps)$ pulls back from an action
with finite isotropy group,
$\Ga_r$, on $\check M(\sP_r,\beps)$.
By Lemma \ref{eq:PseudoFreeActions},
the orders of the isotropy groups of the $S^1$ action on
\begin{equation}
\label{eq:S1ActionOnPieceOfQuotient}
\tilde M_{\fs}\times \bar\Fr(\ft,\fs,\sP_r)
\times_{G(\sP_r)}\check M(\sP_r,\beps),
\end{equation}
given
by this $S^1$ action on $\check M(\sP_r,\beps)$ will divide
$b_r:=|\Ga_r| |\fS(\sP_r)|$.
If we define $b:=b_0\cdots b_N$, then there
is a complex line bundle over ${\bQ\bL}^{\vir}_{\ft,\fs}$,
with restriction to $Q(\bar{\bB\bL}^{\vir}_{\ft,\fs}(\sP_r))$
given by
$$
\tilde M_{\fs}\times \bar\Fr(\ft,\fs,\sP_r)
\times_{G(\sP_r)}\check M(\sP_r,\beps)
\times_{(S^1,\times b)}\CC,
$$
whose pullback to $\bar{\bB\bL}^{\vir}_{\ft,\fs}$
is $\LL_\nu^{\otimes b}$.
By $(S^1,\times b)$ in the preceding expression, we mean that $S^1$ acts with multiplicity $b$
on $\CC$ to ensure that the stabilizers act trivially on $\CC$.  This completes the proof of Item
\eqref{item:GlobalQuotient_4}.

Define the map $\pi_{X,Q}$ by the projection
\eqref{eq:LocalQuotientProjToX} on $Q(\bar{\bB\bL}^{\vir}_{\ft,\fs}(\sP_j))$.
Then
Item \eqref{item:GlobalQuotient_5} follows
from Lemma \ref{lem:ProjectionToSigFromQuotientPiece}.
On each subspace $Q(\bar{\bB\bL}^{\vir}_{\ft,\fs}(\sP_j))$, the projection map
$$
\check\pi_{\fs,j}:
\tilde M_{\fs}\times_{\sG_{\fs}\times S^1}
\bar\Fr(\ft,\fs,\sP_j) \times_{G(\sP_i)}\check M(\sP_j,\beps)
\to
M_{\fs}
$$
satisfies $\pi_{\fs}=\check\pi_{\fs,j}\circ Q$.
By the definition of the quotient map $Q$ and by the proof of Corollary \ref{cor:GlobalFibration},
we see that
$\check\pi_{\fs,i}=\check\pi_{\fs,j}$
on the appropriate intersection of images of the subspaces $Q(\bar{\bB\bL}^{\vir}_{\ft,\fs}(\sP_i))$
in the quotient $\bQ\bL^{\vir}_{\ft,\fs}$.  Hence, the maps $\check\pi_{\fs,i}$ fit together to
define a global map $\check\pi_\fs$ on $\bQ\bL^{\vir}_{\ft,\fs}$ satisfying Item \eqref{item:GlobalQuotient_6}.
\end{proof}

\section{Homology and cohomology classes of the quotient}
\label{sec:HomOfQuot}
To use the quotient
space ${\bQ\bL}^{\vir}_{\ft,\fs}$ in Proposition
\ref{prop:GlobalQuotient} to compute the pairing
\eqref{eq:ReducedPairing1}, we must verify that
cohomology classes appearing in \eqref{eq:ReducedPairing1}
pull back from the quotient.

\begin{lem}
\label{lem:QuotientS1Class}
Continue the notation of Proposition \ref{prop:GlobalQuotient}.
Then there are
a positive integer $n$ and a class,
\begin{equation}
\label{eq:QuotientNu}
\check\nu\in H^2({\bQ\bL}^{\vir}_{\ft,\fs};\RR),
\end{equation}
such that $Q^*\check\nu=n\nu$, where $\nu$ is the cohomology
class appearing in \eqref{eq:ReducedPairing1}.
\end{lem}

\begin{proof}
Choose $\check\nu$ to be the first Chern class of the line
bundle over ${\bQ\bL}^{\vir}_{\ft,\fs}$ defined in
Item \eqref{item:GlobalQuotient_6} of Proposition \ref{prop:GlobalQuotient}.
\end{proof}

\begin{lem}
Continue the notation of Proposition \ref{prop:GlobalQuotient}.
Define
$$
{\bQ\bL}^{\vir}_{\ft,\fs}(\sP_i)
:=
Q(\bar{\bB\bL}^{\vir}_{\ft,\fs}(\sP_i)).
$$
Then there are homology classes,
\begin{equation}
\label{eq:QuotientPiecesFund}
[{\bQ\bL}^{\vir}_{\ft,\fs}(\sP_i)]
\in
H_\bullet ({\bQ\bL}^{\vir}_{\ft,\fs}(\sP_i);\ZZ),
\end{equation}
satisfying
\begin{equation}
Q_*[\bar{\bB\bL}^{\vir}_{\ft,\fs}]
=
\sum_i (\check\iota_i)_*[{\bQ\bL}^{\vir}_{\ft,\fs}(\sP_i)],
\end{equation}
where $\check\iota_i: {\bQ\bL}^{\vir}_{\ft,\fs}(\sP_i)
\to {\bQ\bL}^{\vir}_{\ft,\fs}$ is the inclusion map.
\end{lem}

\begin{proof}
From the definition of $\bar{\bB\bL}^{\vir}_{\ft,\fs}(\sP_i)$ in
\eqref{eq:InstantonLinkComponentBase}  and \eqref{eq:DefineLinkStratum}
and the definition of $\rd_j\bar{\bB\bL}^{\vir}_{\ft,\fs}(\sP_i)$ in
\eqref{eq:DefineBoundaryOne}, we see that the boundary of
$\bar{\bB\bL}^{\vir}_{\ft,\fs}(\sP_i)$ is given by
$$
\rd \bar{\bB\bL}^{\vir}_{\ft,\fs}(\sP_i)
=
\bigcup_{j\neq i}\ \rd_j\bar{\bB\bL}^{\vir}_{\ft,\fs}(\sP_i).
$$
We then define
$$
b{\bQ\bL}^{\vir}_{\ft,\fs}(\sP_i)
:=
Q(\rd \bar{\bB\bL}^{\vir}_{\ft,\fs}(\sP_i)).
$$
Then there is an equality of cycles,
\begin{equation}
\label{eq:DecomposeFund}
[\bar{\bB\bL}^{\vir}_{\ft,\fs}]
=
\sum_i
(\iota_i)_*[\bar{\bB\bL}^{\vir}_{\ft,\fs}(\sP_i), \rd \bar{\bB\bL}^{\vir}_{\ft,\fs}(\sP_i)],
\end{equation}
where $\iota_i:\bar{\bB\bL}^{\vir}_{\ft,\fs}(\sP_i)\to
\bar{\bB\bL}^{\vir}_{\ft,\fs}$ is the inclusion map and
$$
[\bar{\bB\bL}^{\vir}_{\ft,\fs}(\sP_i), \rd \bar{\bB\bL}^{\vir}_{\ft,\fs}(\sP_i)]
\in
H_\bullet(
\bar{\bB\bL}^{\vir}_{\ft,\fs}(\sP_i), \rd \bar{\bB\bL}^{\vir}_{\ft,\fs}(\sP_i);\RR)
$$
is a relative fundamental class.
Item \eqref{item:GlobalQuotient_2} of Proposition \ref{prop:GlobalQuotient}
implies that the map
in the exact sequence of the pair
$\left({\bQ\bL}^{\vir}_{\ft,\fs}(\sP_i), b{\bQ\bL}^{\vir}_{\ft,\fs}(\sP_i)\right)$,
$$
(\check\jmath_i)*:
H_{\max}({\bQ\bL}^{\vir}_{\ft,\fs}(\sP_i);\RR)
\to
H_{\max}({\bQ\bL}^{\vir}_{\ft,\fs}(\sP_i),
b{\bQ\bL}^{\vir}_{\ft,\fs}(\sP_i);\RR),
$$
is an isomorphism.  Hence,
there is a homology class,
$$
[{\bQ\bL}^{\vir}_{\ft,\fs}(\sP_i)]
\in
H_\bullet ({\bQ\bL}^{\vir}_{\ft,\fs}(\sP_i);\RR),
$$
as in
\eqref{eq:QuotientPiecesFund} which maps to
$Q_*([\bar{\bB\bL}^{\vir}_{\ft,\fs}(\sP_i),
\rd\bar{\bB\bL}^{\vir}_{\ft,\fs}(\sP_i)])$
under $(\check\jmath_i)*$.
Applying $Q_*$ to the equality
\eqref{eq:DecomposeFund} yields the conclusion.
\end{proof}

\section{Fiber bundles and pushforwards}
\label{sec:BundlesPushforwards}
We now show how to use the fiber bundle structure
of ${\bQ\bL}^{\vir}_{\ft,\fs}(\sP_i)$ in cohomological
pairings. The
subspaces of the quotient ${\bQ\bL}^{\vir}_{\ft,\fs}(\sP_i)$
are not fiber bundles over $M_{\fs}\times\cl(\Si(X^\ell,\sP_i))$
but do pull back to a fiber bundle over $\Delta(X^\ell,\sP_i)$
as we record in the following lemma.

\begin{lem}
\label{lem:CoverOfQuotientPiece}
For the subgroup $\tG(\sP_i)$ of
the structure group $G(\sP_i)$ defined in \eqref{eq:DefineGluingDataBundleStructureGroup_first_part},
define
$$
\widetilde{\bQ\bL}^{\vir}_{\ft,\fs}(\sP_i)
:=
M_{\fs}\times_{\sG_{\fs}\times S^1}
\bar\Fr(\ft,\fs,\sP_i)\times_{\tG(\sP_i)}
\check M(\sP_i,\beps).
$$
Then there are a branched covering,
$c_i: \widetilde{\bQ\bL}^{\vir}_{\ft,\fs}(\sP_i) \to
{\bQ\bL}^{\vir}_{\ft,\fs}(\sP_i)$,
and a homology class,
$$
[\widetilde{\bQ\bL}^{\vir}_{\ft,\fs}(\sP_i)]
\in
H_\bullet(\widetilde{\bQ\bL}^{\vir}_{\ft,\fs}(\sP_i);\RR),
$$
satisfying $(c_i)_*[\widetilde{\bQ\bL}^{\vir}_{\ft,\fs}(\sP_i)]
=e_i [{\bQ\bL}^{\vir}_{\ft,\fs}(\sP_i)]$ for some $e_i\in\ZZ$
and a commutative diagram,
$$
\begin{CD}
\check M(\sP_i,\beps)
@>>>
\widetilde{\bQ\bL}^{\vir}_{\ft,\fs}(\sP_i)
@> c_i >>
{\bQ\bL}^{\vir}_{\ft,\fs}(\sP_i)
\\
@. @V\tilde\pi_{\bQ,i} VV @V \check \pi_i VV
\\
@. M_{\fs}\times\Delta(X^\ell,\sP_i)
@>>>
M_{\fs}\times \cl (\Si(X^\ell,\sP_i))
\end{CD}
$$
where the map
$\widetilde{\bQ\bL}^{\vir}_{\ft,\fs}(\sP_i)\to M_{\fs}\times\Delta(X^\ell,\sP_i)$
is a fiber bundle with structure group $\sG_{\fs}\times\tG(\sP_i)$.
\end{lem}

We record the behavior of the cohomology class $\check\nu$
under the branched cover $c_i$.

\begin{lem}
\label{lem:BranchedCoverS1Bundle}
Continue the notation of Lemma \ref{lem:CoverOfQuotientPiece}.
Let $\check\nu$ be the cohomology class defined in
\eqref{eq:QuotientNu}.  Then $c_i^*\iota_i^*\check\nu=\tilde \nu_i$,
where $\tilde\nu_i$ is the first Chern class of the $S^1$-bundle
given by
\begin{equation}
\label{eq:DefineBranchedCoverS1Bundle}
M_{\fs}\times_{\sG_{\fs}}
\bar\Fr(\ft,\fs,\sP_i)\times_{\tG(\sP_i)}
\check M(\sP_i,\beps)
\to
M_{\fs}\times_{\sG_{\fs}\times S^1}
\bar\Fr(\ft,\fs,\sP_i)\times_{\tG(\sP_i)}
\check M(\sP_i,\beps).
\end{equation}
\end{lem}

\begin{proof}
The conclusion follows immediately from the definitions
$\check\nu$ and of the map $c_i$.
\end{proof}

We write
\begin{equation}
\label{eq:BranchedCoverProjections}
\tilde\pi_{X,i}:\widetilde{\bQ\bL}^{\vir}_{\ft,\fs}(\sP_i)
\to
\Delta(X^\ell,\sP)
\quad\text{and}\quad
\tilde\pi_{\fs,i}:\widetilde{\bQ\bL}^{\vir}_{\ft,\fs}(\sP_i)
\to M_{\fs}
\end{equation}
for the composition of $\tilde\pi_{\bQ,i}$ with the obvious projection maps.

We now analyze the structure group of the bundle defined
in Lemma \ref{lem:CoverOfQuotientPiece}.
The equality
$\Delta(X^\ell,\sP)=\prod_{P\in\sP}X$,
implies that we can write the bundle $\bar \Fr(\ft,\fs,\sP)\to \Delta(X^\ell,\sP)$
as
$$
\bar \Fr(\ft,\fs,\sP)
=
\prod_{P\in\sP} \Fr(\fg_{\ft(\ell)})\times_X\Fr(TX).
$$
The preceding identity and
the diagonal inclusion map, $\tilde M_{\fs}\to\prod_{P\in\sP}\tilde M_{\fs}$,
define an inclusion of $S^1$ bundles,
\begin{equation}
\label{eq:DiagonalBundleInclusion1}
\begin{CD}
\tilde M_{\fs}\times_{\sG_{\fs}}\bar\Fr(\ft,\fs,\sP)
\times_{\tG(\sP_i)}
\check M(\sP_i,\beps)
\\
@VVV
\\
\prod_{P\in\sP}
\left(
\tilde M_{\fs}\times_{\sG_{\fs}} \Fr(\fg_{\ft(\ell)})\times_X\Fr(TX)
\right)
\times_{\tG(\sP_i)}
\check M(\sP_i,\beps)
\end{CD}
\end{equation}
covering the inclusion of base spaces,
$$
\begin{CD}
\tilde M_{\fs}\times_{\sG_{\fs}\times S^1}\bar\Fr(\ft,\fs,\sP)
\times_{\tG(\sP_i)}
\check M(\sP_i,\beps)
\\
@VVV
\\
\prod_{P\in\sP}
\left(
\tilde M_{\fs}\times_{\sG_{\fs}} \Fr(\fg_{\ft(\ell)})\times_X\Fr(TX)
\right)/S^1
\times_{\tG(\sP_i)}
\check M(\sP_i,\beps)
\end{CD}
$$
We now show how to reduce the structure group of the preceding bundle.
The $\SO(3)$ bundle,  $\Fr(\fg_{\ft(\ell)})\to X$, admits a reduction
to an $S^1$ bundle,
\begin{equation}
\label{eq:S1Reduction}
\Fr(\fg_{\ft(\ell)})
\cong
Q_{\ft,\fs}\times_{S^1}\SO(3),
\end{equation}
where $c_1(Q_{\ft,\fs})=c_1(\ft)-c_1(\fs)$, as described in the
beginning of Section \ref{sec:Reducibles} (see also
\cite[Lemma 2.9]{FL2a}).

\begin{lem}
\label{lem:CharClassOfGsTwistedBundle}
Define an $S^1$-bundle $P_{\ft,\fs}$ by
$$
P_{\ft,\fs}:=\tilde M_{\fs}\times_{\sG_{\fs}} Q_{\ft,\fs} \to M_{\fs}\times X,
$$
where $\sG_{\fs}$ acts on
$Q_{\ft,\fs}$ by $(q,u)\mapsto q u(p_Q(q))^2$, for $q\in Q_{\ft,\fs}$ and $u\in\sG_{\fs}$. Then
$$
\tilde M_{\fs}\times_{\sG_{\fs}}Q_{\ft,\fs}
\to
M_{\fs}\times X
$$
is a principal $S^1$-bundle with first Chern class equal to
$2c_1(\LL_{\fs})+\pi_2^*c_1(Q_{\ft,\fs})$, where
$\pi_2:M_{\fs}\times X\to X$ is the projection.
\end{lem}

\begin{proof}
The conclusion is given by \cite[Equation (3.68)]{FL2a}.
\end{proof}

We define
$$
\tH(\sP_i)
:=
\prod_{P\in\sP} S^1\times\SO(4).
$$
We can thus
rewrite the bundle
\begin{align*}
{}&
\prod_{P\in\sP_i}
\left(
\tilde M_{\fs}\times_{\sG_{\fs}} \Fr(\fg_{\ft(\ell)})\times_X\Fr(TX)
\right)
\times_{\tG(\sP_i)}
\check M(\sP_i,\beps)
\\
{}&\quad
\cong
\prod_{P\in\sP_i}
\left(
\tilde M_{\fs}\times_{\sG_{\fs}} Q_{\ft,\fs}\times_X\Fr(TX)
\right)
\times_{\tH(\sP_i)}
\check M(\sP_i,\beps)
\\
{}&\quad
\cong
\prod_{P\in\sP_i}
\left(
P_{\ft,\fs}\times_X\Fr(TX)
\right)
\times_{\tH(\sP_i)}
\check M(\sP_i,\beps).
\end{align*}
If we denote $\Fr_{\tH}(\sP_i):=\prod_{P\in\sP}P_{\ft,\fs}\times_X\Fr(TX)$, then
composing the preceding
reduction of the structure group with the inclusion
\eqref{eq:DiagonalBundleInclusion1}
yields an inclusion of $S^1$ bundles,
\begin{equation}
\label{eq:DiagonalInclusionOfQuotientPiece}
\begin{CD}
\tilde M_{\fs}\times_{\sG_{\fs}}\bar\Fr(\ft,\fs,\sP)
\times_{\tG(\sP_i)}
\check M(\sP_i,\beps)
@> \tilde\Delta_{\ft,\fs}^i >>
\Fr_{\tH}(\sP_i)
\times_{\tH(\sP_i)}
\check M(\sP_i,\beps)
\\
@VVV @VVV
\\
\tilde M_{\fs}\times_{\sG_{\fs}\times S^1}\bar\Fr(\ft,\fs,\sP)
\times_{\tG(\sP_i)}
\check M(\sP_i,\beps)
@> \Delta_{\ft,\fs}^i>>
\Fr_{\tH}(\sP_i)
\times_{\tH(\sP_i)}
\check M(\sP_i,\beps)/S^1
\end{CD}
\end{equation}
where the $S^1$ action on $\check M(\sP_i,\beps)$ is given
by the diagonal $S^1$ action on frames defined in
\eqref{eq:DiagonalFrameAction}.

Let $f_{\ft,\fs}:M_{\fs}\times X\to \BS^1\times\BSO(4)$
be the classifying map for the bundle
$$
F_{\ft,\fs}:=
P_{\ft,\fs}\times_X\Fr(TX)
\to
M_{\fs}\times X.
$$
By the identities,
$$
\EtH(\sP_i)
=
\prod_{P\in\sP_i} \ES^1\times\ESO(4)
\quad\text{and}\quad
\BtH(\sP_i)
=
\prod_{P\in\sP_i} \BS^1\times\BSO(4),
$$
if we define $f_{\ft,\fs}(\sP):=\prod_{P\in\sP} f_{\ft,\fs}$ and $\tilde f_{\ft,\fs}(\sP)$ to be the bundle maps
in the diagram,
\begin{equation}
\label{eq:BundleMapForFts}
\begin{CD}
\prod_{P\in\sP} F_{\ft,\fs} @> \tilde f_{\ft,\fs}(\sP) >> \EtH(\sP)
\\
@VVV @VVV
\\
\prod_{P\in\sP} M_{\fs}\times X @> f_{\ft,\fs}(\sP) >> \BtH(\sP)
\end{CD}
\end{equation}
and if we write $\id_M$ and $\id_{M/S^1}$ for the identify maps on
$\check M(\sP_i,\beps)$ and $\check M(\sP_i,\beps)/S^1$ respectively,
then
we have a commutative diagram of $S^1$ bundles
and $\tH(\sP_i)$ bundles:
\begin{equation}
\label{eq:BranchedCoverToClassifyingSpace}
\begin{CD}
\prod_{P\in\sP_i}
F_{\ft,\fs}
\times_{\tH(\sP_i)}
\check M(\sP_i,\beps)
@> \tilde f_{\ft,\fs}(\sP_i)\times \id_M >>
\EtH(\sP_i)
\times_{\tH(\sP_i)}
\check M(\sP_i,\beps)
\\
@VVV @VVV
\\
\prod_{P\in\sP_i}
F_{\ft,\fs}
\times_{\tH(\sP_i)}
\check M(\sP_i,\beps)/S^1
@> \tilde f_{\ft,\fs}(\sP_i)\times \id_{M/S^1} >>
\EtH(\sP_i)\times_{\tH(\sP_i)}
\check M(\sP_i,\beps)/S^1
\\
@VVV @V p_i VV
\\
\prod_{P\in\sP_i} (M_{\fs}\times X)
@> f_{\ft,\fs}(\sP_i) >>
\BtH(\sP_i)
\end{CD}
\end{equation}
We then have the following lemma.

\begin{lem}
\label{lem:DefineEquivariantS1Bundle}
Define
$$
\nu_{H,i}
\in
H^2(\EtH(\sP_i)\times_{\tH(\sP_i)}
\check M(\sP_i,\beps)/S^1;\RR)
$$
to be the first Chern class of the bundle
$$
\EtH(\sP_i)\times_{\tH(\sP_i)}
\check M(\sP_i,\beps)
\to
\EtH(\sP_i)\times_{\tH(\sP_i)}
\check M(\sP_i,\beps)/S^1.
$$
If $\tilde\nu_i$ is the cohomology class defined in
\eqref{eq:DefineBranchedCoverS1Bundle},
$\Delta_{\ft,\fs}^i$ the map  defined in
\eqref{eq:DiagonalInclusionOfQuotientPiece},
and we define
\begin{equation}
\label{eq:ComposingDiagonalWithFtsBundleMap}
\tilde g_i:= (\tilde f_{\ft,\fs}(\sP_i)\times \id_{M/S^1})\circ \Delta_{\ft,\fs}^i,
\end{equation}
then
$\tilde g_i^*\nu_{H,i}=\tilde\nu_i$.
\end{lem}

\begin{proof}
The lemma follows immediately from the definitions
of $\nu_{H,i}$ and $\tilde\nu_i$ and the
commutativity of the diagrams
\eqref{eq:DiagonalInclusionOfQuotientPiece}
and \eqref{eq:BranchedCoverToClassifyingSpace}.
\end{proof}

We will require the following computations of pushforwards.

\begin{lem}
\label{lem:PushforwardComputation}
Continue the notation of Lemma \ref{lem:DefineEquivariantS1Bundle}.
For $\sP=\{P_1,\dots,P_r\}$, define
$k(\sP_i)=\sum_{i=1}^r \left(8|P_i|-4\right)-2$.
Let $p_i: \EtH(\sP_i)\times_{\tH(\sP_i)}\check M(\sP_i,\eps)/S^1
\to \BtH(\sP_i)$ be the projection map appearing in
the diagram \eqref{eq:BranchedCoverToClassifyingSpace}.
Then there is a class,
$$
m_{|\sP_i|,\delta}\in H^{2\delta-k(\sP_i)}(\BtH(\sP_i);\RR),
$$
such that
$$
(p_i)_* (\nu_{H,i}^\delta)=m_{|\sP_i|,\delta}.
$$
\end{lem}

\begin{proof}
The conclusion follows by observing that $k(\sP_i)$ is the dimension of
the fiber $\check M(\sP_i,\beps)/S^1$.
\end{proof}

\section{Computations of intersection numbers on $\bar{\bL}_{\ft,\fs}$}
\label{sec:Computations}
We can now begin to compute the pairings
appearing in \eqref{eq:ReducedPairing1}.
First, we show how to reduce such pairings
to
pairings with the
subspaces of the
quotient space.

\begin{lem}
\label{lem:ReduceToLocalPairings}
For any  $\mu_{\fs}(z_2)\in H^\bullet(M_{\fs};\RR)$, and
$S^\ell(\beta)\in H^\bullet(\Sym^\ell(X);\RR)$,
and $z_2\in\AAA_2(X)$,
the equality,
\begin{equation}
\label{eq:ReduceToLocalPairings}
\begin{aligned}
{}&
\left\langle \nu^\delta\smile \pi_{\fs}^*\mu_{\fs}(z_2)\smile\pi_X^*S^\ell(\beta),
[\bar{\bB\bL}^{\vir}_{\ft,\fs}]\right\rangle
\\
&\quad=
\sum_i
n^{-\delta}
\left\langle  \iota_i^*\left(\check\nu^\delta \smile \check\pi_{\fs}^*\mu_{\fs}(z_2)
\smile \check\pi_{X,i}^*S^\ell(\beta)\right),
[{\bQ\bL}^{\vir}_{\ft,\fs}(\sP_i)]\right\rangle,
\end{aligned}
\end{equation}
holds, where we defined $n$ and $\check\nu$ in Lemma
\ref{lem:QuotientS1Class}, and $\check \pi_{X,i}$ in
\eqref{eq:LocalQuotientProjToX}, and
$\check\pi_{\fs}$ in Item \eqref{item:GlobalQuotient_6} of Proposition
\ref{prop:GlobalQuotient}.
\end{lem}

\begin{proof}
From Items \eqref{item:GlobalQuotient_4}, \eqref{item:GlobalQuotient_5}, and \eqref{item:GlobalQuotient_6} of Proposition
\ref{prop:GlobalQuotient} and from
Lemma \ref{lem:QuotientS1Class} we have
\begin{align*}
{}&\left\langle \nu^\delta\smile \pi_{\fs}^*\mu_{\fs}(z_2)\smile\pi_X^*S^\ell(\beta),
[\bar{\bB\bL}^{\vir}_{\ft,\fs}]\right\rangle
\\
{}&\quad
=
n^{-\delta}
\left\langle (Q^*\check\nu)^\delta \smile Q^*\check\pi_{\fs}^*\mu_{\fs}
\smile Q^*\pi_{Q,X}^*S^\ell(\beta),
[\bar{\bB\bL}^{\vir}_{\ft,\fs}]\right\rangle
\\
{}&
\quad
=
n^{-\delta}
\left\langle Q^*\left( \check\nu^\delta \smile \check\pi_{\fs}^*\mu_{\fs}
\smile \check\pi_{Q,X}^*S^\ell(\beta)\right),
[\bar{\bB\bL}^{\vir}_{\ft,\fs}]\right\rangle
\\
{}&
\quad
=
n^{-\delta}
\left\langle \check\nu^\delta \smile \check\pi_{\fs}^*\mu_{\fs}
\smile \check\pi_{Q,X}^*S^\ell(\beta),
Q_*[\bar{\bB\bL}^{\vir}_{\ft,\fs}]\right\rangle
\\
{}&\quad
=
n^{-\delta}
\sum_i
\left\langle  \check\nu^\delta \smile \check\pi_{\fs}^*\mu_{\fs}
\smile \check\pi_{X,i}^*S^\ell(\beta),
(\iota_i)_*[{\bQ\bL}^{\vir}_{\ft,\fs}(\sP_i)]\right\rangle,
\end{align*}
where the final equality follows from Lemma
\ref{eq:QuotientPiecesFund} and the definition of
$\check\pi_{Q,X}$ in terms of $\check\pi_{X,i}$.
\end{proof}

The description of the fiber-bundle structures in the
preceding section leads to the following expression
for the pairing on the right-hand side of
equation \eqref{eq:ReduceToLocalPairings}.

\begin{lem}
\label{lem:UseOfFiberBundle}
Continue the notation of Lemma \ref{lem:ReduceToLocalPairings}.
For $\sP_i=\{P_1,\dots,P_r\}$, let $k(\sP_i)$ be the constant
defined in Lemma \ref{lem:PushforwardComputation}.
Write $\mu_{\fs}(z_2)$ and $\tilde S^\ell(\beta)$ for
$\mu_{\fs}(z_2)\times 1$ and $1\times\tilde S^\ell(\beta)$ respectively.
Then there are universal polynomials,
$$
m_{|\sP_i|,\delta}(\ft,\fs)\in H^{2\delta-k(\sP_i)}(M_{\fs}\times\Delta(X^{\ell},\sP_i);\RR),
$$
in the cohomology classes $c_1(\LL_{\fs})$ and $\tilde S^\ell(\alpha)$,
where $\alpha$ may denote $c_1(\ft)-c_1(\fs)$, $e(X)$, or $p_1(X)$,
with coefficients depending only on $\delta$ and the partition
$\ell=|P_1|+\dots + |P_r|$,
such that
\begin{align*}
{}&
\left\langle  \iota_i^*\left(\check\nu^\delta \smile \check\pi_{\fs}^*\mu_{\fs}(z_2)
\smile \check\pi_{X,i}^*S^\ell(\beta)\right),
[{\bQ\bL}^{\vir}_{\ft,\fs}(\sP_i)]\right\rangle
\\
{}&\quad
=
\left\langle m_{|\sP_i|,\delta}(\ft,\fs)\smile
\mu_{\fs}(z_2)\smile \tilde S^\ell(\beta),[M_{\fs}\times \Delta(X^{\ell},\sP_i)]\right\rangle.
\end{align*}
\end{lem}

\begin{proof}
Lemmas \ref{lem:CoverOfQuotientPiece}
and \ref{lem:BranchedCoverS1Bundle} allow us to write
\begin{align*}
{}&
\left\langle
\iota_i^*\left(\check\nu^\delta \smile \check\pi_{\fs}^*\mu_{\fs}(z_2)
\smile \check\pi_{X,i}^*S^\ell(\beta)\right),
[{\bQ\bL}^{\vir}_{\ft,\fs}(\sP_i)]\right\rangle
\\
{}&\quad =
e_i^{-1}
\left\langle
\iota_i^*\left(\check\nu^\delta \smile \check\pi_{\fs}^*\mu_{\fs}(z_2)
\smile \check\pi_{X,i}^*S^\ell(\beta)\right),
(c_i)_*\left[\widetilde{\bQ\bL}^{\vir}_{\ft,\fs}(\sP_i)\right]\right\rangle
\\
{}&\quad =
e_i^{-1}
\left\langle
c_i^*\iota_i^*\left(\check\nu^\delta \smile \check\pi_{\fs}^*\mu_{\fs}(z_2)
\smile \check\pi_{X,i}^*S^\ell(\beta)\right),
\left[\widetilde{\bQ\bL}^{\vir}_{\ft,\fs}(\sP_i)\right]\right\rangle
\\
{}&\quad =
e_i^{-1}
\left\langle
\tilde\nu_i^\delta \smile
\tilde\pi_{\bQ,i}^*\left(\mu_{\fs}(z_2)\smile \tilde S^\ell(\beta)\right),
\left[\widetilde{\bQ\bL}^{\vir}_{\ft,\fs}(\sP_i)\right]\right\rangle,
\end{align*}
where in the final step we have adopted the convention of
writing $\mu_{\fs}(z_2)$ for $\mu_{\fs}(z_2)\times 1$
and $\tilde S^\ell(\beta)$ for $1\times \tilde S^\ell(\beta)$
mentioned in the statement of the lemma; we record the conclusion as
\begin{multline}
\label{eq:PullPairingsToCover}
\left\langle
\iota_i^*\left(\check\nu^\delta \smile \check\pi_{\fs}^*\mu_{\fs}(z_2)
\smile \check\pi_{X,i}^*S^\ell(\beta)\right),
[{\bQ\bL}^{\vir}_{\ft,\fs}(\sP_i)]\right\rangle
\\
=
e_i^{-1}
\left\langle
\tilde\nu_i^\delta \smile
\tilde\pi_{\bQ,i}^*\left(\mu_{\fs}(z_2)\smile \tilde S^\ell(\beta)\right),
\left[\widetilde{\bQ\bL}^{\vir}_{\ft,\fs}(\sP_i)\right]\right\rangle.
\end{multline}
We then use a pushforward-pullback argument
to compute the pairing in
the right-hand side of \eqref{eq:PullPairingsToCover}.
Define
$$
\Delta_M: M_{\fs}\times\prod_{P\in\sP}X
\to
\prod_{P\in\sP}M_{\fs}\times \Delta(X^\ell,\sP),
$$
by the diagonal inclusion, $M_{\fs}\to\prod_{P\in\sP} M_{\fs}$,
and the identification
$\Delta(X^\ell,\sP)=\prod_{P\in\sP}X$
(where $\Delta(X^\ell,\sP)\subset X^\ell$ is defined in \eqref{eq:DefineDelta}).
Composing the diagrams
\eqref{eq:DiagonalInclusionOfQuotientPiece}
and
\eqref{eq:BranchedCoverToClassifyingSpace}
gives the diagram,
\begin{equation}
\label{eq:PushPullDiagram}
\begin{CD}
\widetilde{\bQ\bL}^{\vir}_{\ft,\fs}(\sP_i)
@> \tilde g_i  >>
\EtH(\sP_i)\times_{\tG(\sP_i)} \check M(\sP_i,\beps)/S^1
\\
@V \tilde\pi_{\bQ,i} VV @V p_i VV
\\
M_{\fs}\times \Delta(X^\ell,\sP_i)
@> g_i >>
\BtH(\sP_i)
\end{CD}
\end{equation}
where
$\tilde g_i$ is defined in \eqref{eq:ComposingDiagonalWithFtsBundleMap} 
and $g_i:=(\prod_{P\in\sP} f_{\ft,\fs})\circ \Delta_M$ (and $f_{\ft,\fs}$ is defined prior to \eqref{eq:BundleMapForFts}).
Lemma \ref{lem:DefineEquivariantS1Bundle} then implies that
we can rewrite
the right-hand side of
\eqref{eq:PullPairingsToCover} as
\begin{equation}
\label{eq:PushPull}
\begin{aligned}
{}&
e_i^{-1}
\left\langle
\tilde\nu_i^\delta \smile
\tilde\pi_{\bQ,i}^*\left(\mu_{\fs}(z_2)\smile \tilde S^\ell(\beta)\right),
\left[\widetilde{\bQ\bL}^{\vir}_{\ft,\fs}(\sP_i)\right]\right\rangle
\\
{}&\quad
=
e_i^{-1}
\left\langle
\tilde g_i^*\nu_{H,i}^\delta \smile
\tilde\pi_{\bQ,i}^*\left(\mu_{\fs}(z_2)\smile \tilde S^\ell(\beta)\right),
\left[\widetilde{\bQ\bL}^{\vir}_{\ft,\fs}(\sP_i)\right]
\right\rangle.
\end{aligned}
\end{equation}
We then apply a pushforward-pullback argument to
the diagram \eqref{eq:PushPullDiagram}
to rewrite the right-hand side of \eqref{eq:PushPull} as:
\begin{align*}
{}&
e_i^{-1}
\left\langle
\tilde g_i^*\nu_{H,i}^\delta \smile
\tilde\pi_{\bQ,i}^*\left(\mu_{\fs}(z_2)\smile \tilde S^\ell(\beta)\right),
\left[\widetilde{\bQ\bL}^{\vir}_{\ft,\fs}(\sP_i)\right]
\right\rangle
\\
{}&\quad
=
e_i^{-1}
\left\langle
(\tilde\pi_{\bQ,i})_*
\left(\tilde g_i^*\nu_{H,i}^\delta \smile
\tilde\pi_{\bQ,i}^*\left(\mu_{\fs}(z_2)\smile \tilde S^\ell(\beta)\right)\right),
[M_{\fs}\times \Delta(X^\ell,\sP_i)]
\right\rangle
\\
{}&\quad
=
e_i^{-1}
\left\langle
(\tilde\pi_{\bQ,i})_*\left(\tilde g_i^*\nu_{H,i}^\delta \right)
\smile \mu_{\fs}(z_2)\smile \tilde S^\ell(\beta),
[M_{\fs}\times \Delta(X^\ell,\sP_i)]
\right\rangle
\\
{}&\quad
=
e_i^{-1}
\left\langle
\left(\tilde g_i^*(p_i)_*\nu_{H,i}^\delta\right)
\smile \mu_{\fs}(z_2)\smile \tilde S^\ell(\beta),
[M_{\fs}\times \Delta(X^\ell,\sP_i)]
\right\rangle
\\
{}&\quad
=
e_i^{-1}
\left\langle
\tilde g_i^*m_{|\sP_i|,\delta}
\smile \mu_{\fs}(z_2)\smile \tilde S^\ell(\beta),
[M_{\fs}\times \Delta(X^\ell,\sP_i)]
\right\rangle,
\end{align*}
where the final equality follows from
Lemma \ref{lem:PushforwardComputation}.
By the definition of $g_i$ in terms of the maps $f_{\ft,\fs}$
and by the definition of the maps $f_{\ft,\fs}$ as classifying maps
for the bundles
$P_{\ft,\fs}$
(defined in Lemma \ref{lem:CharClassOfGsTwistedBundle})
and $\Fr(TX)$, the pullback
$g_i^*m_{|\sP_i|,\delta}$ is given by a universal polynomial
in the characteristic classes of $P_{\ft,\fs}$ and $\Fr(TX)$.
By Lemma \ref{lem:CharClassOfGsTwistedBundle}, these characteristic
classes are $2c_1(\LL_{\fs})+c_1(\ft)-c_1(\fs)$, and $e(X)$, and
$p_1(X)$.  This completes the proof of Lemma \ref{lem:UseOfFiberBundle}.
\end{proof}

\section{Proofs of the main theorems}
\label{sec:Proofs}
We now prove Theorems
\ref{thm:LinkPairing} and \ref{thm:Multiplicity} and Theorem
\ref{thm:MainThm}.

\begin{proof}[Proof of Theorem \ref{thm:LinkPairing}]
Recall Equation \eqref{eq:IntersectionToCohom1}:
\begin{equation}
\label{eq:IntersectionToCohom1a}
\#\left(\bar\sV(h^{\delta-2m}x^m)\cap\bar\sW^{\eta} \cap \bar\bL_{\ft,\fs}\right)
=\left\langle \barmu_p(h^{\delta-2m}x^m)\smile \barmu_c^{\eta}\smile \bar e_I \smile\bar e_s,
 [\bar\bL^{\vir}_{\ft,\fs}] \right\rangle.
\end{equation}
For convenience, we abbreviate
$L_h:=-\frac{1}{2}\langle c_1(\ft)-c_1(\fs),h\rangle$
and $\alpha:=2\mu_\fs-\nu$. By Definition \ref{defn:ExtendedCohomologyClasses}
we can then rewrite \eqref{eq:IntersectionToCohom1a} as
\begin{align*}
{}&\#\left(\bar\sV(h^{\delta-2m}x^m)\cap\bar\sW^{\eta} \cap \bar\bL_{\ft,\fs}\right)
\\
{}&\quad
=
\left\langle
\left(
\alpha L_h
+\pi_X^*S^\ell(h)\right)^{\delta-2m}
\smile
\left( -\frac{1}{4} \alpha^2+S^\ell(x)\right)^m
\smile
\nu^\eta
\smile(-\nu)^{r_\Xi}
\smile
\bar e_I,
[\bar\bL^{\vir}_{\ft,\fs}]
\right\rangle
\\
{}&\quad=
\sum_{r=0}^{\delta-m}
\sum_{i+2j=r}
(-1)^{r_\Xi} C_{\delta,m}(i,j)L_h^i
\left\langle
\alpha^r\smile S^\ell(h)^{\delta-2m-i}\smile S^\ell(x)^{m-j}
\smile\nu^{\eta+r_\Xi}
\smile\bar e_I,
[\bar\bL^{\vir}_{\ft,\fs}]
\right\rangle,
\end{align*}
that is,
\begin{multline}
\label{eq:LinkPairing1}
\#\left(\bar\sV(h^{\delta-2m}x^m)\cap\bar\sW^{\eta} \cap \bar\bL_{\ft,\fs}\right)
\\
=
\sum_{r=0}^{\delta-m}
\sum_{i+2j=r}
(-1)^{r_\Xi} C_{\delta,m}(i,j)L_h^i
\left\langle
\alpha^r\smile S^\ell(h)^{\delta-2m-i}\smile S^\ell(x)^{m-j}
\smile\nu^{\eta+r_\Xi}
\smile\bar e_I,
[\bar\bL^{\vir}_{\ft,\fs}]
\right\rangle,
\end{multline}
where
$$
C_{\delta,m}(i,j):=\binom{\delta-2m}{i}\binom{m}{j}
(-1)^j 2^{-r}.
$$
The expression for $\bar e_I$ in Proposition
\ref{prop:GlobalInstantonEulerClass}
implies that
we can write the terms in the right-hand side of
\eqref{eq:LinkPairing1}
 as
\begin{equation}
\label{eq:LinkPairing2}
\begin{aligned}
{}&
L_h^i
\left\langle
\alpha^r\smile S^\ell(h)^{\delta-2m-i}\smile S^\ell(x)^{m-j}
\smile\nu^{\eta+r_\Xi}
\smile\bar e_I,
[\bar\bL^{\vir}_{\ft,\fs}]
\right\rangle
\\
{}&\quad
=
\sum_{j=0}^\ell
L_h^i
\left\langle\alpha^r\smile
S^\ell(h)^{\delta-2m-i}\smile S^\ell(x)^{m-j}
\smile \nu^{\eta+r_{\Xi}+j}\smile n_{\ell,j},
[\bar\bL^{\vir}_{\ft,\fs}]
\right\rangle,
\end{aligned}
\end{equation}
where $n_{\ell,j}\in H^{2\ell-2j}(M_{\fs}\times\Sym^\ell(X);\RR)$
is a polynomial in $\mu_\fs(\cdot)$ and $S^\ell(c(\ft))$ with
real coefficients
depending only on 
$\ell$ and $j$.
Proposition \ref{prop:ReductionToBase} and the expression for the
Segre classes $s_i(N)$ in \cite[Lemma 4.11]{FL2b}
imply that the terms
in the sum
on the left-hand side of
\eqref{eq:LinkPairing2}
can be rewritten as
\begin{equation}
\label{eq:LinkPairing3}
\begin{aligned}
{}&
L_h^i
\left\langle\alpha^r\smile
S^\ell(h)^{\delta-2m-i}\smile S^\ell(x)^{m-j}
\smile \nu^{\eta+r_{\Xi}+j}\smile n_{\ell,j},
[\bar\bL^{\vir}_{\ft,\fs}]\right\rangle
\\
{}&\quad
=
L_h^i
\left\langle
S^\ell(h)^{\delta-2m-i}\smile S^\ell(x)^{m-j}
\smile s_{\ft,\fs,r}(\nu,\mu_s,\ft,X)
,
[\bar{\bB\bL}^{\vir}_{\ft,\fs}]
\right\rangle,
\end{aligned}
\end{equation}
where $s_{\ft,\fs,r}(\nu,\mu_s,\ft,X)$ is a polynomial in $\nu$,
$\pi_{\fs}^*\mu_{\fs}(\cdot)$, and $\pi_X^*S^\ell(\alpha)$
where $\alpha$ is $c_1(\ft)$, $c_1(\fs)$, $e(X)$, or $p_1(X)$
and the coefficients of $s_{\ft,\fs,r}(\nu,\mu_s,\ft,X)$ depend
only on the homotopy type of $\ft$, $\fs$, and $X$.
Applying  Lemmas
\ref{lem:ReduceToLocalPairings} and
\ref{lem:UseOfFiberBundle} to
\eqref{eq:LinkPairing3} yields the equality,
\begin{equation}
\label{eq:LinkPairing4}
\begin{aligned}
{}&
\left\langle
S^\ell(h)^{\delta-2m-i}\smile S^\ell(x)^{m-j}
\smile s_{\ft,\fs,r}(\nu,\mu_s,\ft,X)
,
[\bar{\bB\bL}^{\vir}_{\ft,\fs}]
\right\rangle
\\
{}&\quad
=
\sum_{k=0}^r
L_h^i
\left\langle
q(\sP_k)\smile
S^\ell(h)^{\delta-2m-i}\smile S^\ell(x)^{m-j},
[M_{\fs}\times \Delta(X^\ell,\sP_k]
\right\rangle,
\end{aligned}
\end{equation}
where $q(\sP_k)\in H^\bullet(M_{\fs}\times\Delta(X^\ell,\sP_k)$
is a polynomial in $\mu_{\fs}(\cdot)$ and $S^\ell(\alpha)$, and
where $\alpha$ can equal $c_1(\ft)$, $c_1(\fs)$,
$(c_1(\ft)-c_1(\fs)^2$, $e(X)$, or $p_1(X)$.
Because $b_1(X)=0$, the equality
$$
[M_{\fs}\times\Delta(X^\ell,\sP_k)]
=
[M_{\fs}]\times [\Delta(X^\ell,\sP_k)]
$$
implies that the expressions
in \eqref{eq:LinkPairing4} will vanish unless they contain
a factor $\mu_{\fs}(x)^{d_s(\fs)/2}$.  Hence, the terms in the sum
on the right-hand side of
\eqref{eq:LinkPairing4} can be rewritten as
\begin{equation}
\label{eq:LinkPairing5}
\begin{aligned}
{}&
L_h^i
\left\langle
q(\sP_k)\smile
S^\ell(h)^{\delta-2m-i}\smile S^\ell(x)^{m-j},
[M_{\fs}\times \Delta(X^\ell,\sP_k]
\right\rangle
\\
{}&\quad
=
\SW_X(\fs)\la_h^i
\left\langle
S^\ell(h)^{\delta-2m-i}\smile S^\ell(x)^{m-j}\smile \alpha_{\ft,\fs,}(\delta,m,j,k),
[\Delta(X^\ell,\sP_k]
\right\rangle,
\end{aligned}
\end{equation}
where $\alpha_{\ft,\fs,}(\delta,m,j,k)\in H^\bullet(\Delta(X^\ell,\sP_k)$
is a polynomial in $S^\ell(\alpha)$ for $\alpha$ equal to
$c_1(\ft)$, $c_1(\fs)$, $(c_1(\ft)-c_1(\fs)^2$, $p_1(X)$, or $e(X)$.
Because $\Delta(X^\ell,\sP_k)$ is a product of $|\sP_k|$ copies of
the four-dimensional
space $X$ and because $S^\ell(h)$ is a two-dimensional cohomology class,
terms in \eqref{eq:LinkPairing5}
will vanish unless each factor of $S^\ell(h)$ is paired
either with another factor of $S^\ell(h)$ or a factor $S^\ell(\alpha_2)$,
where $\alpha_2=c_1(\ft)$ or $c_1(\ft)-c_1(\fs)$.
Pairing $S^\ell(h)$ with itself will give
powers of $Q_X(h)$, pairing $S^\ell(h)$ with $S^\ell(\alpha_2)$
will give powers of $\langle c_1(\ft),h\rangle$
or $\langle c_1(\ft)-c_1(\fs),h\rangle$.
Thus, recalling the definition of $L_h$
from before \eqref{eq:LinkPairing1},
we see that the
expression on the right-hand side of
\eqref{eq:LinkPairing5} must be a sum over
terms of the form
$$
\left\langle c_1(\ft)-c_1(\fs),h\right\rangle^i
\left\langle c_1(\ft),h\right\rangle^j
Q_X(h)^a,
$$
where $i+j+2a=\delta-2m$.  We note that the power $a$
appearing in $Q_X(h)^a$ must be less than  or equal to the length of the partition
$\sP_k$, so $a\le \ell$.  This completes the proof of Theorem \ref{thm:LinkPairing}.
\end{proof}

We record a particular case of Theorem \ref{thm:LinkPairing}
in a form
which is used in the proof of Witten's Conjecture \cite{FL7,FL8}.

\begin{prop}
\label{prop:theorem_LinkPairing_dim_Ms_zero}
Continue the notation and hypotheses of
Theorem \ref{thm:LinkPairing}.  In addition,
assume that $\dim M_{\fs}=0$ and 
$$
2(\delta+\eta)=\dim \sM_{\ft}-2.
$$
Writing $\fs_h:=\langle c_1(\fs),h\rangle$ and
$\la_h:=\langle c_1(\ft),h\rangle$, then
\begin{equation}
\label{eq:ReducedForm}
\begin{aligned}
{}&
\left\langle
\nu^{\delta+\eta-2m-i}\smile \pi_X^*S^\ell(h)^i\smile \barmu_p(x)^m
\smile \bar e_I\smile e_s,
[\bar\bL^{\vir}_{\ft,\fs}]
\right\rangle
\\
{}&\quad=
\SW_X(\fs)
\sum_{j=0}^{\min\left( \ell,[i/2]\right)}
f_{\delta,\ell,m,\eta,i,j}(\fs_h,\la_h) Q_X(h)^j,
\end{aligned}
\end{equation}
where
$$
f_{\delta,\ell,m,\eta,i,j}(\fs_h,\la_h)
=
\sum_{k=0}^{i-2j} a_{\delta,\ell,m,\eta,i,j,k}\fs_h^{i-2j-k}\la_h^k,
$$
and the coefficients $a_{\delta,\ell,m,\eta,i,j,k}$ depend on
$\delta$, $\ell$, $m$, $\eta$, $i$, $j$, $k$, $\chi$, $\si$,
$c_1(\ft)^2$, $c_1(\fs)\cdot c_1(\ft)$, and $c_1(\fs)^2$.
\end{prop}

\begin{proof}
Under the assumption
$\dim M_{\fs}=0$,
the left-hand side of
\eqref{eq:LinkPairing2} reduces
to the expression on
the left-hand side of \eqref{eq:ReducedForm}. Hence,
the conclusion of Proposition \ref{prop:theorem_LinkPairing_dim_Ms_zero} follows from the
proof of
Theorem \ref{thm:LinkPairing}  from
\eqref{eq:LinkPairing2} to the end.
\end{proof}

We now give the

\begin{proof}[Proof of Theorem \ref{thm:Multiplicity}]
The proof of Theorem  \ref{thm:LinkPairing}
shows that the intersection number
in the statement of the theorem
equals a sum of
expressions of the form on the left-hand-side of
\eqref{eq:LinkPairing5}.
Without the assumption that $b_1(X)=0$, the left-hand side of \eqref{eq:LinkPairing5}
equals a sum over terms that each contain a factor of the form $\langle \mu_\fs(\om),[M_{\fs}]\rangle$,
for some $\om\in\AAA_2(X)$.  These vanish by hypothesis and hence the intersection number vanishes.
\end{proof}

The definition of the Donaldson invariant for $X$ in
\eqref{eq:DefineDonaldson} is given in terms of
the blow-up, $\tilde X=X\#\bar{\CC\PP}^2$, of $X$.
Hence, before commencing the proof of Theorem \ref{thm:MainThm}
we will need to restate Theorem \ref{thm:LinkPairing}
in a form useful for computations on $\tilde X$.

Let $e\in H_2(\tilde X;\ZZ)$ be the fundamental class
of the exceptional curve and define $e^*:=\PD[e]$.
We will identify elements of the homology and cohomology
of $X$ with elements of the homology and cohomology
of $\tilde X$ by means of the obvious inclusions.

By the discussions in \cite[Section 12.4]{SalamonSWBook} or
\cite[Section 4.6]{NicolaescuSWNotes},
given $\fs\in\Spinc(X)$, there is a unique \spinc structure $\tilde\fs_k\in\Spinc(\tilde X)$
such that $c_1(\tilde\fs_k)=c_1(\fs)+(2k-1)e^*$ and all
such \spinc structures on $\tilde X$ are so obtained.  The dimensions of the
corresponding Seiberg--Witten moduli space satisfy
$$
d_s(\tilde\fs_k)=d_s(\fs)-k(k-1).
$$
We then have the

\begin{lem}
\label{lem:BlowUpLinkPairing}
Continue the notation and hypotheses of Theorem \ref{thm:LinkPairing}.
Let $\tilde X=X\#\bar{\CC\PP}^2$ be the blow-up of $X$ and
let $e\in H_2(\tilde X;\ZZ)$ be the homology class of the exceptional curve.
Let $\tilde\ft$ be a \spinu structure over $\tilde X$
satisfying
$$
p_1(\tilde\ft)= p_1(\ft)-1,
\quad
c_1(\tilde\ft)=c_1(\ft),
\quad\text{and}\quad
w_2(\tilde\ft)\equiv w_2(\ft)+\PD[e]\pmod 2.
$$
Then
\begin{equation}
\label{eq:BlownUpLinkIntersectionFormula}
\begin{aligned}
{}&
\#\left( \bar\sV(ze)\cap\bar\sW^{\eta}\cap\bar\bL_{\tilde\ft,\tilde\fs_k}\right)
\\
&\quad
=
\SW_X(\fs)
\sum_{i=0}^{ \min(\ell,\left[\frac{\delta-2m}{2}\right])}
\left(
\tilde q_{\delta,\ell,m,i,k}(c_1(\fs)-c_1(\ft),c_1(\ft))Q_X^i
\right) (h),
\end{aligned}
\end{equation}
where $\tilde q_{\delta,\ell,m,i,k}$ are degree $\delta-2m-2i$
homogeneous polynomials
whose coefficients are universal functions
of the constants given in Theorem \ref{thm:MainThm}.
\end{lem}

\begin{proof}
The proof follows that of Theorem \ref{thm:LinkPairing}
once one has made the following observations.
By \eqref{eq:ReducibleLevel},
the levels of $\fs$ and $\tilde\fs_k$
are related by
\begin{align*}
4\ell(\tilde\ft,\tilde\fs_k)
{}&=
\left( c_1(\tilde \fs_k) -c_1(\tilde\ft)\right)^2 - p_1(\tilde\ft)
\\
{}&=
\left( c_1(\fs) -c_1(\tilde\ft)\right)^2-(2k-1)^2 -p_1(\ft)+1
\\
{}&=
4\ell(\ft,\fs)-4k(k-1).
\end{align*}
Hence, $\ell(\tilde\ft,\tilde\fs_k)\le \ell(\ft,\fs)$.
Next,
by \cite[Theorem 4.6.8]{NicolaescuSWNotes},
$$
\SW_{\tilde X}(\tilde\fs_k)
=
\begin{cases}
\pm \SW_X(\fs) & \text{if $d(\fs)\ge k(k-1)$,}
\\
0 &\text{if $d(\fs)< k(k-1)$.}
\end{cases}
$$
The conclusion of the lemma now follows from the proof of Theorem \ref{thm:LinkPairing}.
\end{proof}

\begin{proof}[Proof of Theorem \ref{thm:MainThm}]
To compute the Donaldson polynomial appearing on
the left-hand side of \eqref{eq:MainEquation},
we select a \spinu structure $\tilde\ft$ over $\tilde X$
with $\La=c_1(\tilde \ft)$ and $p_1(\tilde\ft)$ determined
by
$$
\delta+1=-p_1(\tilde\ft)-\frac{3}{4}\left(\chi(X)+\si(X)\right),
$$
where $\delta$ is the constant appearing in \eqref{eq:MainEquation}.
To apply \eqref{eq:ASDPairing}, we need to verify that
$n_a(\tilde \ft)>0$.  From the definition of $n_a(\tilde\ft)$
in \eqref{eq:Transv}
and for $i(\La)=\La^2-\frac{1}{4}(3\chi+7\si)$ as  in the statement
of Theorem \ref{thm:MainThm}, we have
\begin{align*}
4n_a(\tilde\ft)
{}&=
p_1(\tilde\ft) +\La^2 -i(\tilde X)
\\
{}&=
-\delta-1-\frac{3}{4}\left(\chi(X)+\si(X)\right)+\La^2- i(\tilde X)+1
\\
{}&=
-\delta+\La^2-\quarter\left(3\chi(X)+7\si(X)\right)
\\
{}&=
i(\La)-\delta.
\end{align*}
Hence, the hypothesis of the theorem
that $\delta<i(\La)$
implies that $n_a(\tilde\ft)>0$.  Equation \eqref{eq:ASDPairing}
and the definition of the Donaldson invariant in \eqref{eq:DefineDonaldson}
then imply that for $z=h^{\delta-2m}x^m$, and $n_a=n_a(\tilde \ft)$,
and $\tilde w=w+e^*$, and $\tilde\ka=-(1/4)p_1(\tilde\ft)$ we have
\begin{equation}
\label{eq:MTProof1}
\begin{aligned}
D^w_X(z)
{}&=
\#\left( \bar\sV(ze)\cap \bar M^{\tilde w}_{\tilde\ka}\right)
\\
{}&
=
2^{1-n_a}
\#\left( \bar\sV(ze)\cap \bar\sW^{n_a-1}\cap \bar\bL^{\tilde w}_{\tilde\ft,\tilde\ka}\right).
\end{aligned}
\end{equation}
The cobordism formula \eqref{eq:RawCobordismSum1} implies that
\begin{equation}
\label{eq:MTProof2}
\#\left( \bar\sV(ze)\cap \bar\sW^{n_a-1}\cap \bar\bL^{\tilde w}_{\tilde\ft,\tilde\ka}\right)
=
-\sum_{\tilde\fs\in\Spinc(\tilde X)}
(-1)^{o_{\ft}(\tilde w,\tilde\fs)}
\#\left( \bar\sV(ze)\cap \bar\sW^{n_a-1}\cap\bar\bL_{\tilde \ft,\tilde\fs}\right).
\end{equation}
For any four-manifold $Y$ with $b_1(Y)=0$, let $B(Y)$ denote the set of Seiberg--Witten basic classes,
$$
B(Y)
:=
\{c_1(\fs):\fs\in\Spinc(Y)\ \text{and}\ \SW_Y(\fs)\neq 0\}.
$$
By \cite[Theorem 4.6.8]{NicolaescuSWNotes}, the Seiberg--Witten basic classes of $\tilde X$ and
$X$ are related by
$$
B(\tilde X)
=
\{c_1(\tilde \fs_k): c_1(\fs)\in B(X)\ \text{and}\ d_s(\fs)-k(k-1)\ge 0\}.
$$
Theorem \ref{thm:Multiplicity} then implies that
the
right-hand side of \eqref{eq:MTProof2}
equals
\begin{equation}
\label{eq:MTProof3}
\begin{aligned}
{}&
-\sum_{\tilde\fs\in B(\tilde X)}
(-1)^{o_{\tilde\ft}(\tilde w,\tilde\fs)}
\#\left( \bar\sV(ze)\cap \bar\sW^{n_a-1}\cap\bar\bL_{\tilde \ft,\tilde\fs}\right)
\\
{}&\quad
=
-\sum_{\fs\in B(X)}
\sum_{\begin{subarray}{c} k^2-k\le d_s(\fs), \\ k \geq 0\end{subarray}}
(-1)^{o_{\ft}(\tilde w,\tilde\fs_k)}
\#\left( \bar\sV(ze)\cap \bar\sW^{n_a-1}\cap\bar\bL_{\tilde \ft,\tilde\fs_k}\right).
\end{aligned}
\end{equation}
The definition of $\fo_{\tilde\ft}(\tilde w,\tilde\fs_k)$ in
\eqref{eq:OrientChangeFactor}
implies that
$$
o_{\tilde\ft}(\tilde w,\tilde\fs_k)
=
\frac{1}{4}
\left( w+e^*-\La+c_1(\fs)+(2k-1)e^* \right)^2
=
o_{\ft}(w,\fs)-k.
$$
Thus, we can rewrite the sum
on the right-hand side of
\eqref{eq:MTProof3} as
\begin{equation}
\label{eq:MTProof4}
-\sum_{\fs\in B(X)}
(-1)^{o_{\ft}(w,\fs)}
\sum_{\begin{subarray}{c} k^2-k\le d_s(\fs), \\ k \geq 0\end{subarray}}
(-1)^{k}
\#\left( \bar\sV(ze)\cap \bar\sW^{n_a-1}\cap\bar\bL_{\tilde \ft,\tilde\fs_k}\right).
\end{equation}
Theorem \ref{thm:MainThm} then follows from
Equations \eqref{eq:MTProof1}--\eqref{eq:MTProof4}
and Lemma \ref{lem:BlowUpLinkPairing}.
\end{proof}

\chapter{Kotschick--Morgan Conjecture}
\label{chap:KMconj}
We now show how the
proof of the
Kotschick--Morgan Conjecture,
stated here as Conjecture \ref{conj:KoM} and described in
Section \ref{subsec:IntroKMConj},
can be reduced to
the proof of
a gluing theorem analogous to Hypothesis \ref{hyp:Gluing}.
The methods for this reduction are exactly those that lead to the proof of
Theorem \ref{thm:LinkPairing}.  Hence, there will be little new work in this section;
instead, we will reference the analogous discussion from earlier sections. Throughout this chapter, we will
assume that $X$ is a smooth, closed, and oriented four-manifold with
$b^+(X)=1$ and $b_1(X)=0$.

\section{Cobordisms and reducible connections}
\label{sec:Cobordisms_reducible_connections}
As described in Section \ref{subsec:IntroKMConj},
the Donaldson invariant of a manifold with $b^+=1$ depends on the Riemannian metric.
The Kotschick--Morgan Conjecture \ref{conj:KoM},
describes how the Donaldson invariant
changes with the Riemannian metric.  These changes are given by intersection
numbers on the link of
gauge-equivalence classes of
reducible connections in the cobordism
\eqref{eq:DefineModuliCobordism}.  We therefore now review
how such gauge-equivalence classes appear in the cobordism.

For a complex, rank-two Hermitian vector bundle $E\to X$ with $c_1(E)=w$, denote
$\fg^w_\ka:=\su(E)$, where $-4\ka=p_1(\fg^w_\ka)$.
Let $g_I$ denote a smooth path of Riemannian metrics on $X$ given by $I=[-1,1] \ni t\mapsto g_t$.
Let $M^w_\ka(g_t)$ be
the moduli space of
$\SO(3)$ connections on $\fg^w_\ka$ which are anti-self-dual with
respect to the metric $g_t$, as defined in \eqref{eq:ASDModuliSpace}, and
let $\bar M^w_\ka(g_t)$ be its Uhlenbeck compactification.
The smooth path of Riemannian metrics, $g_I$, on $X$ defines a
cobordism of moduli spaces,
\begin{equation}
\label{eq:DefineModuliCobordism}
\bar M^w_\ka(g_I)
=
\{[A,\bx,t]: [A,\bx]\in \bar M^w_\ka(g_t)\ \text{and $t\in [-1,1]$}\}.
\end{equation}
To simplify the discussion, we shall assume that $w$ is good
in the sense of Definition \ref{defn:Good}.  By the Morgan--Mrowka blow-up trick
\cite{MorganMrowkaPoly}, this assumption can always be satisfied
and so there is no loss of generality.

A reducible connection $A$ on $\fg^w_\ka\to X$ can be written
as $A=\Theta_\RR\oplus A_L$
with respect to a splitting $\fg^w_\ka\cong i\ubarRR\oplus L$,
where $\Theta_\RR$ is the product connection on the product bundle
$\ubarRR = X\times\RR$ and $A_L$ is a unitary connection on a complex line bundle $L$
with $c_1(L)\equiv w\pmod 2$, and $c_1(L)^2=p_1(\fg^w_\ka)=-4\ka$
(see \cite[Proposition 4.2.15]{DK}).
The connection
$\Theta_\RR\oplus A_L$ will be anti-self-dual if and only if
$A_L$ is anti-self-dual.  Let $\om^+(g)$ denote the unique (once
an orientation for $H^+(X)$ is specified) unit-length self-dual 
harmonic two-form for the metric $g$.  Then, by
\cite[p. 147]{DK}, the connection $A_L$ will be anti-self-dual if
and only if
\begin{equation}
\label{eq:HarmonicPerp}
\om^+(g)\smile F_{A_L}=0.
\end{equation}
Hence, $M^w_\ka(g)$ will contain a gauge-equivalence class of a reducible connection of the form $\Theta_\RR\oplus A_L$
if and only if $\om^+(g)$ lies on the codimension-one `wall'
in the positive cone,
$$
\Om_X:=\{\beta\in H^2(X;\RR): \beta^2>0\}/\RR^+,
$$
of $H^2(X;\RR)/\RR^+$ defined
by the zero-locus of the map $\Om \ni [h]\mapsto  h\smile c_1(L)$.
If \eqref{eq:HarmonicPerp} holds and $b^1(X)=0$, then by
 \cite[Proposition 2.2.6]{DK}  there is a unique
gauge-equivalence class in $M^w_\ka(g)$ of a connection giving the splitting $\fg^w_\ka\cong i\ubarRR\oplus L$.
We summarize the above discussion in the following

\begin{lem}
\label{lem:ReductionConditions}
The compactification $\barM^w_\ka(g_0)$ contains a gauge-equivalence class of
an ideal reducible connection $A$ in the sense that
$$
[A]\times \Sym^\ell(X) \subset \barM^w_\ka(g_0),
$$
for $A=\Theta_\RR\oplus A_L$
with respect to a reduction $\fg^w_{\ka-\ell}\cong \ubarRR\oplus L$
for a complex line bundle $L$ with $\ell\ge 0$ \emph{if and only if}
\begin{equation}
\label{eq:ReductionConditions}
c_1(L)^2=-4(\ka-\ell)\ge 0 \quad\text{and}\quad
c_1(L)\equiv w\pmod 2,
\end{equation}
and \eqref{eq:HarmonicPerp} holds with $g=g_0$.
\end{lem}

The wall defined by $\alpha\in H^2(X;\ZZ)$ is a \emph{$(w,\ka)$-wall}
if there is a complex line bundle $L\to X$
satisfying \eqref{eq:ReductionConditions} with $\alpha=c_1(L)$.

\section{Cohomology classes on the cobordism}
\label{sec:ASDCohomOnCobord}
Let $\sA^w_\ka$ be the space of $L^2_k$ ($k\ge 2$)
$\SO(3)$ connections on $\fg^w_\ka$ and for
a rank-two, complex Hermitian bundle $E\to X$ with $\su(E)\cong \fg^w_\ka$,
let $\sG_E$ be the
Hilbert Lie group of $L^2_{k+1}$, determinant-one, unitary gauge
transformations of $E$.  The group $\sG_E$ acts through the
adjoint representation on $\su(E)\cong \fg^w_\ka$. Let
$\sB^w_\ka:=\sA^w_\ka/\sG_E$ denote the quotient space of connections
and define
$$
\bar\sB^w_\ka:=\bigsqcup_{\ell=0}^{[\ka]}\ \sB^w_{\ka-\ell}\times\Sym^\ell(X),
$$
with topology defined by Uhlenbeck convergence as in
Definition \ref{defn:UhlenbeckConvergence}.
Then there is an embedding, $\iota_M$, defined by the composition,
$$
\iota_M: \bar M^w_\ka(g_I) \to \bar\sB^w_\ka\times [-1,1]
\to \bar\sB^w_\ka,
$$
of the inclusion with the projection.
The set of gauge-equivalence classes of reducible connections in $\bar\sB^w_\ka$ are given by
$$
\bar R^w_\ka=\{[A,\bx]\in \bar\sB^w_\ka: \text{$A$ is reducible}\},
$$
and we define
$$
\sB^{w,*}_\ka:=\sB^w_\ka \setminus\bar R^w_\ka,\quad
\bar\sB^{w,*}_\ka:=\bar\sB^w_\ka \setminus\bar R^w_\ka,
\quad\text{and}\quad
\bar M^{w,*}_\ka(g_I):=\bar M^w_\ka(g_I)\setminus\iota_M^{-1}(\bar R^w_\ka).
$$
For $z\in \AAA(X)$, let $\mu(z)\in H^{\deg(z)}(\sB^{w,*}_\ka;\RR)$ be
the cohomology class described in \cite[Definition 5.1.12]{DK}
or \cite[Section 2.2]{KMStructure}.  We note that the cohomology class
$\mu_p(z)$ defined in \eqref{eq:MuPMap} is the pullback of $\mu(z)$ by
the projection map $\sC^*_\ft\to\sB^{w,*}_\ka$.
In \cite[Lemma 4.1.2]{KotschickMorgan}, Kotschick and Morgan define a cohomology class,
$$
\barmu(z)\in H^{\deg(z)}(\bar\sB^{w,*}_\ka;\RR),
$$
whose restriction to $\sB^{w,*}_\ka$ equals $\mu(z)$.
For Riemannian metrics $g$ with $\bar M^w_\ka(g)\cap\bar R^w_\ka$ empty,
the Donaldson invariant for $g$ is given by
$$
D^w_X(g)(z)
:=
\langle \barmu(z),[\bar M^w_\ka(g)]\rangle.
$$
From \cite[Theorem 3.0.1]{KotschickMorgan}, this invariant only depends on the
`chamber' in which the self-dual harmonic two-form $\om^+(g)$
appearing in \eqref{eq:HarmonicPerp} lies.
By \emph{chamber}, we mean a connected component of the complement in
the positive cone $\Om_X$
of all the walls of type $(w,\ka)$.
As described in Section \ref{subsec:IntroKMConj},
if $g_I$ is a smooth path of Riemannian metrics,
the failure of the cobordism
given by $\bar M^w_\ka(g_I)$ to give the equality
$D^w_X(g_{-1})(z)=D^w_X(g_1)(z)$ as in  \eqref{eq:DInvariantWallCross}
arises from the presence of gauge-equivalence classes of ideal reducible connections, $[A,\bx]\in \iota_M(\bar M^w_\ka(g_I))$.

For each
complex line bundle $L\to X$ with $c_1(L)$ defining a wall of type $(w,\ka)$,
we will construct in Section \ref{subsec:DefineReducibleNghInParamModuli}
(see \eqref{eq:DefineNghOfReduciblesInParamModuli})
a closed neighborhood $\bar U^w_\ka(L)$ of
$[A]\times\Sym^\ell(X)$ in $\iota_M(\bar M^w_\ka(g_I))$.
The {\em difference term\/} corresponding to this gauge-equivalence class of a reducible connection is
\begin{equation}
\label{eq:DifferenceTerm}
\delta^w_{L,\ka}(z):=\left\langle \iota_M^*\barmu(z),[\rd\bar U^w_\ka(L)]\right\rangle.
\end{equation}
Computations of the difference term for $\ell=0$ appeared in
\cite{DonHCobord,KotschickBPlus1}, for $\ell=1$ in \cite{Yang,LenessWC},
and for $\ell=2$ in \cite{LenessWC}, where $\ell$ is the level of
the reducible connection as it appears in Lemma \ref{lem:ReductionConditions}.
Details of the construction of $\bar U^w_\ka(L)$ appearing
in the following sections show that the boundary of this neighborhood
admits a fundamental class.
The difference between the Donaldson invariants for the Riemannian metrics $g_1$ and $g_{-1}$
is given by
$$
D^w_X(g_1)(z) - D^w_X(g_{-1})(z)=\sum_L \delta^w_{L,\ka}(z),
$$
where the sum is over all complex line bundles $L$
for which $c_1(L)$ defines a wall of type $(w,\ka)$
which $\om^+(g_t)$ crosses as $t$ moves from $-1$ to $1$.
We can now state the \emph{Kotschick--Morgan Conjecture}.

\begin{conj}
\label{conj:KoM}
(See Kotschick and Morgan \cite[Conjectures 6.2.1 and 6.2.2]{KotschickMorgan}.)
Let $X$ be a smooth, connected, oriented four-manifold with $b^+(x)=1$ and $b_1(X)=0$.
Let $\ka\in\tquarter\ZZ$ with $\ka<0$ and let $w\in H^2(X;\ZZ)$ satisfy $w^2=-4\ka\pmod 4$.
Let $\delta$ and $m$ be non-negative integers satisfying
$2\delta=8\ka-(3/2)(\chi(X)+\si(X))$ and $2m\le \delta$.
For a complex line bundle $L\to X$ with $c_1(L)$ defining a wall of type $(w,\ka)$,
let $\ell=(c_1(L)^2+4\ka)/4$.
Then,
for $h\in H_2(X;\RR)$ and a generator $x\in H_0(X;\ZZ)$,
$$
\delta^w_{L,\ka}(h^{\delta-2m}x^m)
=
\sum_{i=0}^{\min( \ell,[(\delta-2m)/2] )} a_i \langle c_1(L),h\rangle^{\delta-2m-2i} Q_X(h)^i,
$$
where the coefficients $a_i$ depend only on $e(X)$, $\si(X)$, $m$,
$\delta$, and $\ell$.
\end{conj}

\section[Neighborhoods of ideal reducible connections]{Neighborhoods of gauge-equivalence classes of ideal reducible connections}
To compute the difference term \eqref{eq:DifferenceTerm}, we show
how to define and parameterize the neighborhood $\bar U^w_\ka(L)$ of
$\{[A]\}\times\Sym^\ell(X)$ in
$\bar M^w_\ka(g_I)$,
where $[A]$ is the gauge-equivalence class of a reducible anti-self-dual connection.

\subsection{Kuranishi model for a neighborhood of a reducible connection}
We first record a description of a neighborhood of a gauge-equivalence class of a reducible
connection in the parameterized moduli space $M^w_{\ka-\ell}(g_I)$.

\begin{lem}
\label{lem:ParamModuliKuran}
Let $[A]\in M^w_{\ka-\ell}(g)$ be a gauge-equivalence class of a reducible connection,
$A=\Theta_\RR\oplus A_L$, with respect to a splitting
$\fg^w_{\ka-\ell}\cong \ubarRR\oplus L$.  Then for generic smooth paths of metrics $g_I$
with $g=g_0$
and for integers $n, r\in\ZZ$ satisfying $2n=\dim M^w_{\ka-\ell}+1$ and $n+r\ge 0$,
there are
\begin{enumerate}
\item
\label{item:ParamModluliKuran1}
An open ball, $B(\delta)\subset\CC^{n+r}$, with center at the origin and radius $\delta$,
\item
\label{item:ParamModluliKuran2}
A smooth, $S^1$-equivariant map $\chi_A:B(\delta)\to \CC^r$ (where $S^1$ acts by scalar multiplication on the domain and range),
and
\item
\label{item:ParamModluliKuran3}
An $S^1$-equivariant embedding $\tilde\varphi_A:B(\delta)\to \sA^w_{\ka-\ell}$ (where $S^1$ acts on $\sA^w_{\ka-\ell}$ as
the stabilizer of $A$),
\end{enumerate}
such that $\tilde\varphi_A$ covers a smoothly-stratified embedding,
$$
\varphi_A: B(\delta)/S^1\to \sB^w_{\ka-\ell},
$$
with the property that
$\varphi_A(0)=[A]$ and
$\varphi_A(\chi_A^{-1}(0)/S^1)$ is a neighborhood of $[A]$ in $M^w_{\ka-\ell}(g_I)$.
\end{lem}

\begin{proof}
The conclusion (compare \cite[Proposition 4.2.23]{DK}) follows from a standard argument using
the Kuranishi Lemma \cite[Proposition 4.2.19]{DK} and
the transversality of the period map $g\mapsto \om^+(g)$
from \cite[Proposition 4.3.14]{DK}.  An extended discussion
of this model appears in \cite[Section 2.6]{FLKM1}.
\end{proof}

\begin{rmk}
If the index $n$ in Lemma \ref{lem:ParamModuliKuran}
is strictly positive, then by
adapting
\cite[Theorem 4.19]{FU} to the case $b^+(X)=1$,
one can assume that  for a generic path of metrics,
there is no obstruction in
the Kuranishi model, that is, $r=0$.
The index $n$ will be negative when $c_1(L)^2=-1$ and as this case
can be handled by the techniques developed in this monograph, we include it for completeness.
\end{rmk}

\subsection{Crude splicing maps}
We now combine the description of a neighborhood of a gauge-equivalence class of a reducible connection,
$[A]$ in $\iota_M(M^w_{\ka-\ell}(g_I))$, from Lemma \ref{lem:ParamModuliKuran}
with the construction of the space of global splicing data from Chapter \ref{chap:GlobalSplicingData} to parameterize the neighborhood
$\bar U^w_\ka(g_I)$ appearing in \eqref{eq:DifferenceTerm}
of the level $[A]\times \Sym^\ell(X)$ in $\iota_M(\bar M^w_\ka(g_I))$.

Let $\sP$ be a partition of $N_\ell$ with $\Si=\Si(X^\ell,\sP)$.
Although \spinu structures play no direct role in this construction, we introduce
them here to allow us to use the results from earlier parts of this monograph.
Let $\ft=(\rho,V_\ft)$ be a \spinu structure on $X$ with $\fg_\ft\cong \fg^w_\ka$
and let $\fs=(\rho_\fs,W)$ be a \spinc structure on $X$ so that the
\spinu structure $\ft(\ell)=(\rho_\ell,V_{\ft(\ell)})$ admits a splitting
$V_{\ft(\ell)}=W\oplus W\otimes L$ as in Section \ref{subsec:RedPU2Monopole}
and $\fg_{\ft(\ell)}\cong i\ubarRR\oplus L$ as in \eqref{eq:EndReduction}.
We note the equality,
$$
c_1(L)=c_1(\ft)-c_1(\fs),
$$
and define
$$
\Fr(L,\sP):=\Fr(\ft,\fs,\sP),\quad
\bar\Gl(L,\sP):=\bar\Gl(\ft,\fs,\sP),
\quad
\sO(L,\sP):=\sO(\ft,\fs,\sP),
$$
where the frame bundle $\Fr(\ft,\fs,\sP)$ is defined in \eqref{eq:DefineGluingDataBundle},
the bundle of gluing data $\bar\Gl(\ft,\fs,\sP)$ is defined in \eqref{eq:GluingData},
and the subspace $\sO(\ft,\fs,\sP)\subset \bar\Gl(\ft,\fs,\sP)$ is defined in \eqref{eq:DefineCrudeSplicingDomain}.
The reduction
$\fg^w_{\ka-\ell}\cong \ubarRR\oplus L$
and the $S^1$ action on $L$ defined by scalar multiplication
define an $S^1$ action on $\Fr(\fg^w_{\ka-\ell})$. This $S^1$ action in turn
defines an $S^1$ action on $\bar\Gl(\ft,\fs,\sP)$,
\begin{equation}
\label{eq:S1ActionOnASDFrameBundle}
S^1\times\bar\Gl(L,\sP) \to \bar\Gl(L,\sP),
\end{equation}
given by the diagonal $S^1$ action on the factors of $\Fr(\fg^w_{\ka-\ell})$
making up $\bar\Gl(\ft,\fs,\sP)$.
We will define the maps parameterizing
a neighborhood of $[A]\times \Si(X^\ell,\sP)$ in $\bar M^w_\ka(g_I)$ on
a subspace of
$$
B(\delta)\times_{S^1}\bar\Gl(L,\sP),
$$
where $B(\delta)$ is as defined in Lemma \ref{lem:ParamModuliKuran}
and
where $S^1$ acts diagonally by the action in Lemma \ref{lem:ParamModuliKuran} on
$B(\delta)$ and by the action \eqref{eq:S1ActionOnASDFrameBundle} on $\bar\Gl(L,\sP)$.

For
$\zeta\in B(\delta)$,
let $g(\zeta)$ be the Riemannian metric with respect to which
$\tilde\varphi_A(\zeta)$
is anti-self-dual, where $\tilde\varphi_A$ is the map given in Lemma \ref{lem:ParamModuliKuran}.
By applying the flattening construction of Lemma \ref{lem:FlatteningMetrics2} to the metrics
$g(\zeta)$ with
$\zeta\in B(\delta)$, we obtain a family of
locally-flattened metrics parameterized by $B(\delta)/S^1\times \Si$.
Define $\sO(L,\sP):=\sO(\ft,\fs,\sP)$ as in \eqref{eq:DefineCrudeSplicingDomain}.
The crude splicing map from \eqref{eq:CrudeSplicingMapExplicit}, applied with this family of
locally-flattened Riemannian metrics and with the background section $\Phi_0$ set equal to zero,
defines a smoothly-stratified embedding,
\begin{equation}
\label{eq:KMLocalCrudeSplicingMap}
\bga_{L,\sP}'':
B(\delta)\times_{S^1}\sO(L,\sP)
\to \bar\sB^w_\ka.
\end{equation}
We will construct the space $\bar U^w_\ka(L)$ as the union, over partitions $\sP$ of $N_\ell$,
of the images of the zero locus of an obstruction section under these maps.

\subsection{Overlap spaces and maps}
The next step in constructing the link $\rd\bar U^w_\ka(L)$ is
to define a space of overlap data and upwards and downwards transition maps
as in Section \ref{sec:OverlapSpacesModuli}.
Let $\sP<\sP'$ be partitions of $N_\ell$.
Following \eqref{eq:OverlapGluingDataSpace}, we define
$$
\Gl(L,\sP,[\sP'])
=
\Fr(L,\sP)\times_{G(\sP)}
\bigsqcup_{\sP''\in[\sP<\sP']} \prod_{P\in\sP} \left(
\Delta^\circ(Z_{|P|}(\delta_P),\sP''_P) \times
\barM(\sP''_P)\right).
$$
To define the upwards overlap map, we follow
\eqref{eq:SimultUpwardsOverlapSpaces} and
introduce the space
$$
\Gl(L,[\sP<\sP']) =\left. \left(
\bigsqcup_{\sP''\in[\sP<\sP']}
\Fr(L,\sP'')\times_{\tilde G(\sP'')}
\prod_{Q\in\sP''}\barM^{s,\natural}_{\spl,|Q|} \right)\right/\fS(\sP).
$$
If we define an open subspace $\sO(L,\sP,[\sP'])\subset \Gl(L,\sP,[\sP'])$
by
condition \eqref{eq:XUpwardsTransitionDomain1},
then we can define an $S^1$-equivariant upwards transition map,
$$
\rho^{u,L}_{\sP,[\sP']}:
B(\delta)\times\sO(L,\sP,[\sP'])
\to
B(\delta)\times\bar\Gl(L,[\sP<\sP']),
$$
exactly as was done for the upwards transition map
$\rho^{\ft,\fs,u}_{\sP,[\sP']}$ in \eqref{eq:UnparamUpwardsTransition}.
By $S^1$-equivariant above, we mean equivariant with respect to the
$S^1$ action on the domain given by the
diagonal action on $B(\delta)$ and on $\Fr(L,\sP)$
and the same action on the image.

If we further assume that the open subspace $\sO(L,\sP,[\sP'])$
satisfies the condition in \eqref{eq:XUpwardsTransitionDeomainRequirement},
we can define an $S^1$-equivariant downward transition map,
$$
\rho^{d,L}_{\sP,[\sP']}:
B(\delta)\times\sO(L,\sP,[\sP'])
\to
B(\delta)\times\Gl(L,\sP),
$$
by the construction in \eqref{eq:XDownwardTransition} and
\eqref{eq:DownwardsInclusionFiberMap}.

As described in the proof of Proposition \ref{prop:XOverlapControl},
if the neighborhoods $\sO(L,\sP)$ and $\sO(L,\sP'')$ are sufficiently
small, then we can find an open set $\sO(L,\sP,[\sP'])$ such
that the overlap of the images of the crude splicing maps is
``controlled'' by the overlap space in the sense that
\begin{equation}
\label{eq:ASDOverlapNgh}
\begin{aligned}
{}&
\bga_{L,\sP}''\left(B(\delta)\times_{S^1}\sO(L,\sP) \right)
\cap
\bga_{L,\sP'}''\left(B(\delta)\times_{S^1}\sO(L,\sP') \right)
\\
{}&\quad
=
\left(\bga_{L,\sP}''\circ \rho^{d,L}_{\sP,[\sP']}\right)
\left(B(\delta)\times_{S^1}\sO(L,\sP,[\sP']) \right).
\end{aligned}
\end{equation}
More precisely, we have the following

\begin{prop}
\label{prop:ASDOverlapSplicing}
Let $\sP<\sP'$ be partitions of $N_\ell$.
Assume that the families of metrics $g_{\sP''}$ and $g_{\sP}$
used to define the crude splicing maps $\bga_{L,\sP}''$
and $\bga_{L,\sP''}''$ satisfy the
conditions in Section \ref{subsubsec:FiberBundleMetrics}.
Then there are open neighborhoods $\sO(L,\sP)$ of
$\Si(X^\ell,\sP)$ in $\bar\Gl(L,\sP)$ and an $\fS(\sP)$-invariant family
of neighborhoods $\sO(L,\sP'')$  in $\bar\Gl(L,\sP'')$ of $\Si(X^\ell,\sP'')$
such that if the open set $\sO(L,\sP,[\sP'])$ satisfies
the analogues of the conditions
\eqref{eq:XUpwardsTransitionDomain1},
\eqref{eq:XUpwardsTransitionDeomainRequirement},
and those following \eqref{eq:DefineOverlapDomain},
then the following diagram commutes and \eqref{eq:ASDOverlapNgh} holds:
\begin{equation}
\label{eq:ASDGlobalSplicingCD}
\begin{CD}
B(\delta)\times_{S^1}\sO(L,\sP,[\sP'])
@> \rho^{L,u}_{\sP,[\sP']} >>
B(\delta)\times_{S^1}\sO(L,[\sP<\sP'])
\\
@V \rho^{L,d}_{\sP,[\sP']} VV @V
\bga''_{L,[\sP<\sP']} VV
\\
\tilde B(\delta)\times_{S^1}\sO(L,\sP) @>
\bga''_{L,\sP} >> \bar\sB^w_\ka
\end{CD}
\end{equation}
\end{prop}

\begin{proof}
The proof is identical to that of Proposition
\ref{prop:XOverlapControl}.
\end{proof}

\subsection{Definition of the neighborhood of a gauge-equivalence class of an ideal reducible connection}
\label{subsec:DefineReducibleNghInParamModuli}
We define the space of global splicing data by
\begin{equation}
\label{eq:KoMGlobalSplicing}
\widetilde\sU^w_\ka(L)
:=
\left.\left( \bigsqcup_{\sP} B(\delta) \times \sO(L,\sP) \right)\right/ \sim,
\end{equation}
where the relation $\sim$ is defined as follows.
If  $\sP_1<\sP_2$ are partitions of $N_\ell$,
$B(\delta)\subset\CC^n$ is the ball appearing in Lemma \ref{lem:ParamModuliKuran},
$\sO(L,\sP_i)$ is the space of gluing data defined prior to  \eqref{eq:KMLocalCrudeSplicingMap}, and
$(z_i,\bA_i)\in B(\delta)\times\sO(L,\sP_i)$ for $i=1,2$,
then  $(z_1,\bA_1)\sim (z_2,\bA_2)$ if $z_1=z_2$ and
there is a point $\bA_{12}\in\sO(L,\sP_i,[\sP_j])$ with
$\rho^{L,u}_{\sP_1,[\sP_2]}(z_2,\bA_{12})=(z_2,\bA_2)$
and
$\rho^{L,d}_{\sP_1,[\sP_2]}(z_1,\bA_{12})=(z_1,\bA_1)$.
While
Proposition \ref{prop:ASDOverlapSplicing} implies that this relation satisfies the transitivity
property up to $S^1$ quotients, the proof without the $S^1$ quotient requires the following.

\begin{lem}
The relation $\sim$ defined above is transitive.
\end{lem}

\begin{proof}
Assume that $(z_1,\bA_1)\sim (z_2,\bA_2)$ and $(z_2,\bA_2)\sim (z_3,\bA_3)$.
Thus, $z_1=z_2=z_3$.
By Proposition \ref{prop:ASDOverlapSplicing},
$\bga''_{L,\sP_1}([z_1,\bA_1])=\bga''_{L,\sP_2}([z_2,\bA_2])$
and  $\bga''_{L,\sP_2}([z_2,\bA_2])=\bga''_{L,\sP_3}([z_3,\bA_3])$, so
$\bga''_{L,\sP_1}([z_1,\bA_1])=\bga''_{L,\sP_3}([z_3,\bA_3])$,
where $[z_i,\bA_i]\in B(\delta)\times_{S^1}\sO(L,\sP_i)$ denotes the image
of $(z_i,\bA_i)\in B(\delta)\times\sO(L,\sP_i)$ in the $S^1$ quotient.
Assume that $\sP_1<\sP_2<\sP_3$ (the proofs for the other cases are identical).
Proposition \ref{prop:ASDOverlapSplicing},
the equality $\bga''_{L,\sP_1}([z_1,\bA_1])=\bga''_{L,\sP_3}([z_3,\bA_3])$,
and the equalities $z_1=z_2=z_3$
imply that there is a point $\bA_{13}\in\sO(L,\sP_1,[\sP_3])$
with
$[\rho^{L,u}_{\sP_1,[\sP_3]}(z_1,\bA_{13})]=[z_1,\bA_3]=[z_3,\bA_3]$
and
$[\rho^{L,d}_{\sP_1,[\sP_2]}(z_1,\bA_{13})]=[z_1,\bA_1]$.
Consequently, there is a $\la\in S^1$ such that $(z_1,\bA_1)\sim \la (z_1,\bA_3)$.
Hence, $z_1=\la z_1$ and so if $z_1\neq 0$, then $\la=1$ and $(z_1,\bA_1)\sim  (z_1,\bA_3)$.

If $z_1=0$, then we still have $\bA_{13}\in\sO(L,\sP_1,[\sP_3])$
with $[\rho^{L,u}_{\sP_1,[\sP_3]}(0,\bA_{13})]=[0,\bA_3]$
and
$[\rho^{L,d}_{\sP_1,[\sP_2]}(0,\bA_{13})]=[0,\bA_1]$.
Because $\rho^{L,d}_{\sP_1,[\sP_2]}$ and $\rho^{L,d}_{\sP_1,[\sP_2]}$ are open embeddings,
we can find a sequence $\{w_\alpha\}\subset B(\delta)\setminus\{0\}$
converging to $0$ and a sequence
$\{\bA_{13}(\alpha)\}\subset\sO(L,\sP_1,[\sP_3])$ with
$[\rho^{L,u}_{\sP_1,[\sP_3]}(w_\alpha,\bA_{13}(\alpha))]=[w_\alpha,\bA_3]$
and
$[\rho^{L,d}_{\sP_1,[\sP_2]}(w_\alpha,\bA_{13}(\alpha))]=[w_\alpha,\bA_1]$.
Since $w_\alpha\neq 0$, the preceding argument implies that
$\rho^{L,u}_{\sP_1,[\sP_3]}(w_\alpha,\bA_{13}(\alpha))=(w_\alpha,\bA_3)$
and
$\rho^{L,d}_{\sP_1,[\sP_3]}(w_\alpha,\bA_{13}(\alpha))=(w_\alpha,\bA_1)$.
Because $\rho^{L,d}_{\sP_1,[\sP_2]}$ and $\rho^{L,d}_{\sP_1,[\sP_2]}$ are open embeddings,
$\bA_{13}(\alpha)$ converges to $\bA_{13}$ and so the continuity of these maps implies
that
$\rho^{L,u}_{\sP_1,[\sP_3]}(0,\bA_{13})=(0,\bA_3)$
and
$\rho^{L,d}_{\sP_1,[\sP_3]}(0,\bA_{13})=(0,\bA_1)$.
Hence $(0,\bA_1)\sim (0,\bA_3)$ as required.
\end{proof}

For partitions $\sP_1<\sP_2$ of $N_\ell$, the overlap maps
$\rho^{L,u}_{\sP_1,[\sP_2]}$ and $\rho^{L,d}_{\sP_1,[\sP_2]}$ are
$S^1$ equivariant, so
the $S^1$ actions on $B(\delta)\times\sO(L,\sP_1)$
and $B(\delta)\times\sO(L,\sP_2)$ respect the equivalence relation $\sim$.
Hence, these $S^1$ actions define an $S^1$ action on $\widetilde\sU^w_\ka(L)$.
By Proposition \ref{prop:ASDOverlapSplicing}, we then obtain
$$
\bga_{\sU,L}'': \widetilde\sU^w_\ka(L)/S^1 \to \bar\sB^w_\ka,
$$
a global crude splicing map.

\subsection{Thom--Mather structures on $\widetilde\sU^w_\ka(L)/S^1$}
We consider $\sU(L,\sP)=B(\delta)\times\sO(L,\sP)$
as a subspace of $\widetilde\sU^w_\ka(L)$ and define maps
$$
\pi(L,\sP):\sU(L,\sP)\to B(\delta)\times\Si(X^\ell,\sP),
\quad
\vec t(L,\sP): \sU(L,\sP)\to [0,1]^{\sP},
$$
as was done in \eqref{eq:DefineXTubularProj} and
\eqref{eq:DefineTubularDistanceFunction}.

The proof of Lemma \ref{lem:TM1ForGlobalSplicingData} then
implies that for $\sP<\sP'$ on $\sU(L,\sP)\cap\sU(L,\sP')$,
we have the equality
$$
\pi(L,\sP)\circ \pi(L,\sP')=\pi(L,\sP).
$$
The analogues of Lemmas \ref{lem:SymProjTubDistRelations0}
and \ref{lem:TubularDistanceFunctionOnOverlaps1} hold
for $\vec t(L,\sP)$.

\subsection{Global projection map for $\widetilde\sU^w_\ka(L)/S^1$}
By the arguments of Lemma \ref{lem:GlobalProjectionToX},
we can define a map,
$$
\pi_{L,X}: \widetilde\sU^w_\ka(L)/S^1 \to \Sym^\ell(X),
$$
whose restriction to $\sU(L,\sP)$ is homotopic to $\pi(L,\sP)$.

\subsection{Global splicing map on $\widetilde\sU^w_\ka(L)/S^1$}
The construction of the global splicing map $\bga_{\sM}'$
of \eqref{eq:DefineGlobalSplice} translates to give
a global splicing map,
$$
\bga_L':\widetilde\sU^w_\ka(L)/S^1\to\bar\sB^w_\ka,
$$
whose restriction to $\sU(L,\sP)$ is homotopic to the
crude splicing map $\bga_{L,\sP}''$.

\subsection{Obstruction bundle on $\widetilde\sU^w_\ka(L)/S^1$}
\label{subsec:ObstructionBundle}
Let $n$ and $r$ be the non-negative integers appearing in Lemma \ref{lem:ParamModuliKuran}.
Define a complex-rank-$r$ vector bundle,
$\widetilde\Xi^w_\ka(L)\to\widetilde\sU^w_\ka(L)$, by
$$
\widetilde\Xi^w_\ka(L):=\left.\left( \bigsqcup_{\sP} (B(\delta)\times\CC^r\times \sO(L,\sP))\right)\right/\sim_\Xi
$$
where the equivalence relation $\sim_\Xi$ is defined, for $(z_i,\zeta_i,\bA_i)\in B(\delta)\times\CC^r\times \sO(L,\sP_i)$, $i=1,2$,
by
$$
(z_1,\zeta_1,\bA_1)\sim_\Xi (z_2,\zeta_2,\bA_2)
$$
if $\zeta_1=\zeta_2$ and $(z_1,\bA_1)\sim (z_2,\bA_2)$, where $\sim$ is the equivalence relation
used in the definition of $\widetilde\sU^w_\ka(L)$ in
\eqref{eq:KoMGlobalSplicing}.
The diagonal action of $S^1$ on $B(\delta)\times\CC^r\times \sO(L,\sP)$ given by $S^1$ acting on $B(\delta)$
and $\CC^r$ by scalar multiplication and on $\sO(L,\sP)$ as in \eqref{eq:S1ActionOnASDFrameBundle}
is compatible with the relation $\sim_\Xi$ and so $\widetilde \Xi^w_\ka(L)$ is an $S^1$-equivariant bundle
over $\widetilde\sU^w_\ka(L)$.
If we define
$$
\widetilde\sU^{w,*}_\ka(L):=\widetilde\sU^w_\ka(L)\setminus\left([A]\times\Sym^\ell(X)\right),
$$
then the $S^1$ action on $\widetilde\sU^w_\ka(L)$ is free on $\widetilde\sU^{w,*}_\ka(L)$ and the $S^1$-quotients
define a vector bundle,
\begin{equation}
\label{eq:KoMGlobalObstructionBundle}
\widetilde \Xi^w_\ka(L)/S^1 \to \widetilde\sU^{w,*}_\ka(L)/S^1.
\end{equation}
This is the desired obstruction bundle.

\subsection{Gluing hypothesis}
We can now state the analogue for  parameterized anti-self-dual connections of the gluing hypothesis for $\SO(3)$ monopoles.

\begin{hyp}
There is a smoothly-stratified, continuous embedding,
$$
\bga_L: \widetilde\sU^w_\ka(L)/S^1\to \bar \sB^w_\ka,
$$
and a smoothly-stratified section $\bchi_L$ of the vector bundle $\widetilde\Xi^w_\ka(L)/S^1\to\widetilde\sU^w_\ka(L)/S^1$
defined in \eqref{eq:KoMGlobalObstructionBundle}
with the following properties:
\begin{enumerate}
\item
The restriction of $\bga_L$ to
$[A]\times\Sym^\ell(X)$ is the
identity map.
\item
There is a homotopy, through smoothly-stratified embeddings,
between $\bga_L'$ and $\bga_L$.
\item
The image $\bga_L(\bchi_L^{-1}(0))$ is an open neighborhood of
$[A]\times\Sym^\ell(X)$ in $\bar M^w_\ka(g_I)$.
\item
The restriction of the section $\bchi_L$ to each stratum of $\sU^{w,*}_\ka(L)/S^1$ vanishes transversely.
\end{enumerate}
\end{hyp}

We then set
\begin{equation}
\label{eq:DefineNghOfReduciblesInParamModuli}
\bar U^w_\ka(L):=
\bga_L\left(\bchi_L^{-1}(0)\cap\widetilde\sU^w_\ka(L)/S^1 \right),
\end{equation}
to define the neighborhood $\bar U^w_\ka(L)$ appearing in \eqref{eq:DifferenceTerm}.

\section{Cohomology classes on the space of global splicing data}
Define
$$
\nu\in H^2(\widetilde\sU^{w,*}_\ka(L)/S^1;\ZZ)
$$
to be the first Chern class of the $S^1$ bundle
$\widetilde\sU^{w,*}_\ka(L)\to \widetilde\sU^{w,*}_\ka(L)/S^1$.
From \cite[Lemma 4.7.4]{KotschickMorgan}, we have
the following computation of the pullback of the $\mu$-classes
by the gluing map.

\begin{lem}
\label{lem:KoMmuClasses}
For $h\in H_2(X;\RR)$ and a generator $x\in H_0(X;\ZZ)$, then the following hold:
\begin{enumerate}
\item
$\bga_L^*\barmu(h)=\frac{1}{2}\langle c_1(L),h\rangle \nu + \pi_{L,X}^*S^\ell(h)$,
\item
$\bga_L^*\barmu(x)=-\frac{1}{4}\nu^2 + \pi_{L,X}^*S^\ell(x)$,
\end{enumerate}
where the cohomology classes $S^\ell(h)$ and $S^\ell(x)$ are defined in
Definition \ref{defn:DefnSymmBeta}.
\end{lem}

\begin{lem}
\label{lem:KoMObstructionEuler}
The Euler class of the obstruction bundle \eqref{eq:KoMGlobalObstructionBundle} is given by
$$
e\left(\widetilde \Xi^w_\ka(L)/S^1\right)=(-\nu)^r.
$$
\end{lem}

\begin{proof}
The result follows immediately from \cite[Lemma 3.27]{FL2a} and the observation that the $S^1$ action
on $\widetilde \Xi^w_\ka(L)$
is diagonal, so the obstruction
bundle is the direct sum of $r$ copies of the line bundle given by the negative of the $S^1$ action defining the cohomology
class $\nu$.
\end{proof}

\section{Definition of the link of a gauge-equivalence class of an ideal reducible connection}
We define a
virtual
link $\rd \widetilde\sU^w_\ka(L)/S^1$
following the construction of Section \ref{subsec:AmbientLink}.
The standard norm on $\CC^n$ defines an
$S^1$-invariant norm on $B(\delta)$ and thus a map,
$$
t_B: \widetilde\sU^{w,*}_\ka(L)/S^1 \to [0,\delta].
$$
Enumerate the strata of $\Sym^\ell(X)$ by partitions
$\sP_0$,\dots,$\sP_n$ in the manner described in
Section \ref{subsec:EnumStrata}.
Choose a small, generic constant $\eps_i$ for each
stratum as was done in Section \ref{subsec:AmbientLink}
(with $\eps_i>\eps_j$ for $i<j$).
Following \eqref{eq:DefineLinkStratum}, define
$$
\rd\sU(L,\sP_i)
:=
\vec t(L,\sP_i)^{-1}(\rd\bar D(\sP_i,\eps_i))
\setminus
\bigcup_{j\neq i}
\vec t(L,\sP_j)^{-1}(D(\sP_j,\eps_j)).
$$
Then, following
\eqref{eq:AmbientLink},
\eqref{eq:DefineSWComponentOfLink}, and
\eqref{eq:InstantonLinkComponent}, we define the virtual link
of $[A_L]\times\Sym^\ell(X)$
by
$$
\rd \widetilde\sU^w_\ka(L)
:=
t_B^{-1}(\delta)
\cup
\left(
\bigcup_i \rd\sU(L,\sP_i)
\right).
$$
The link of $[A_L]\times\Sym^\ell(X)$ is then defined by
$$
\rd\bar U^w_\ka(L):=
\bga_L\left(\bchi_L^{-1}(0)\cap\rd\widetilde\sU^w_\ka(L)/S^1 \right).
$$
The existence of a fundamental class for the virtual link $\rd \widetilde\sU^w_\ka(L)$
follows from the discussion in Section \ref{sec:FundClassOfAmbLink}.
Up to a sign depending on the orientation discussed in
\cite[p. 450--451]{KotschickMorgan},
\begin{equation}
\label{eq:KMDifference}
\begin{aligned}
\delta^w_{L,\ka}(h^{\delta-2m}x^m)
{}&=
\left\langle
\left( \frac{1}{2}\langle c_1(L),h\rangle \nu + \pi_{L,X}^*S^\ell(h)
\right)^{\delta-2m} \right.
\\
{}&\qquad
\smile\left.
\left(
-\frac{1}{4}\nu^2 + \pi_{L,X}^*S^\ell(x)
\right)^m,
e\left(\widetilde \Xi^w_\ka(L)/S^1\right)\frown
\left[\rd \widetilde\sU^w_\ka(L)/S^1\right]
\right\rangle
\end{aligned}
\end{equation}
is the difference term.

\section{Computations of the difference term}
\label{sec:KoMComputations}
We can now prove

\begin{thm}
Conjecture \ref{conj:KoM} is true.
\end{thm}

\begin{proof}
Using the argument in Lemma \ref{lem:ReduceToBase} and
Proposition \ref{prop:ReductionToBase}, the pairing in \eqref{eq:KMDifference}
can be reduced to a pairing with $t_B^{-1}(0)\cap \rd \widetilde\sU^w_\ka(L)$.
Observe that
$$
t_B^{-1}(0)\cap \rd\sU(L,\sP_i)
\subset
S^1\backslash\sO(L,\sP_i),
$$
that is, $t_B^{-1}(0)$ is given by $\{0\}\in B(\delta)$.
Following Proposition \ref{prop:GlobalQuotient}, we
construct a quotient
map
$$
Q:
t_B^{-1}(0)\cap \rd \widetilde\sU^w_\ka(L)/S^1
\to Q^w_\ka(L)
$$
satisfying the following properties:
\begin{enumerate}
\item
There is a map $\check\pi_X: Q^w_\ka(L) \to \Sym^\ell(X)$
satisfying $\check\pi_X\circ Q=\pi_{L,X}$.
\item
There is a class $\check\nu\in H^2(Q^w_\ka(L))$ and a
non-zero integer $r$
with $Q^*\check\nu=r\nu$.
\item
If $Q^w_\ka(L,\sP_i)=Q(t_B^{-1}(0)\cap \rd\sU(L,\sP_i))$, then
$$
Q^w_\ka(L,\sP_i)
\cong
S^1\backslash\bar\Fr(L,\sP_i)\times_{G(\sP_i)} \check M(\sP_i,\eps_i).
$$
\end{enumerate}
We can write the pairing
\eqref{eq:KMDifference} in terms of a sum over pairings of the form
\begin{equation}
\label{eq:LocalKoMPairing}
\left\langle \check\nu^i\smile\check\pi_X^*(S^\ell(h)^j\smile S^\ell(x)^k),
[Q^w_\ka(L,\sP_u)]\right\rangle.
\end{equation}
The conclusion now follows by an argument analogous to that of the proofs of
Lemma \ref{lem:UseOfFiberBundle} and Theorem \ref{thm:LinkPairing}.
\end{proof}

\begin{rmk}
If $b^1(X)>0$, then the gauge-equivalence classes of reducible, anti-self-dual, ideal
connections in $\barM^w_\ka(g_I)$ would no longer appear in families of the form
$[A]\times\Sym^\ell(X)$ but rather in families of the form
$$
\left( H^1(X;i\RR)/H^1(X;2i\ZZ)
\right)
\times \Sym^\ell(X).
$$
The techniques of this chapter extend to construct a link $\rd\bar U^w_\ka(L)$ of such
families but one must replace the spaces $B(\delta)\times_{S^1}\sO(L,\sP)$ appearing in
\eqref{eq:KMLocalCrudeSplicingMap} with spaces of the form
$$
B(\delta)\times H^1(X;i\RR)
\times_{S^1\times H^1(X;2i\ZZ)}
\sO(L,\sP),
$$
where $S^1\times H^1(X;2i\ZZ)$ acts on $\Gl(L,\sP)$
through the inclusion $\bh_p:S^1\times H^1(X;i\ZZ)\to\Map(X,S^1)$
given by the harmonic gauge transformations described in
\cite[Section 2.4.1]{FL2a}.
The construction of the space $\widetilde\sU^w_\ka(L)$ for manifolds with $b^1=0$
works when $b^1>0$ as does the construction of the quotient space $Q^w_\ka(L)$
given in Section \ref{sec:KoMComputations}.  However,
the expression for the difference term $\delta^w_{L,\ka}$ will be more
complicated when $b^1(X)>0$ because of the presence of
the Chern class of the line bundle,
\begin{equation}
\label{eq:HarmonicGTLineBundle}
H^1(X;i\RR)\times X\times_{H^1(X;i\ZZ)}\CC
\to
\left( H^1(X;i\RR)/H^1(X;i\ZZ)\right)\times X,
\end{equation}
computed in \cite[Lemma 2.23]{FL2a}
among the characteristic classes appearing in the analogue of Lemma \ref{lem:UseOfFiberBundle}
for pairings of the form \eqref{eq:LocalKoMPairing}.
In addition, the Chern class of the line bundle \eqref{eq:HarmonicGTLineBundle}
will appear in the computations of the $\mu$-classes in
Lemma \ref{lem:KoMmuClasses}
and change those expressions
to a form similar to that of the cohomology classes appearing in Definition \ref{defn:ExtendedCohomologyClasses}.

For this reason, the difference term $\delta^w_{L,\ka}$, although still only depending on
the homotopy type of $X$, $\fg^w_\ka$, $\ell$, and $L$, will not be a polynomial only in $c_1(L)$
and $Q_X$ but will also contain terms given by elements of $H^1(X;\ZZ)$.
\end{rmk}

\twocolumn
\chapter*{Glossary of Notation}
\bigskip

$\bA,\bA'$
\hfill
Eq. \eqref{eq:NotationForXOverlapSplicingData}

$[\vec A]$
\hfill Lemma
\ref{defn:LocalQuotient}

$\sA_{\ft}$
\hfill
\S \ref{subsubsec:PU2Monopoles}

$\sA_{\fs}$
\hfill
\S \ref{subsec:SWMonopoles}

$\AAA_2(X)$
\hfill
Eq.
\eqref{eq:DefineAAA2}

$\AAA(X)$
\hfill
Eq.
\eqref{eq:DefineAAA}

$\Ad^u_{\hfill\SO(3)}$, $\Ad^u_{\hfill\SO(4)}$, $\Ad^u$
\hfill
Eq. \eqref{eq:DefineAdu}

$\sB_\ka$, $\sB^s_\ka$
\hfill
before Eq. \eqref{eq:IdealConnS4}

$\bar\sB^{s}_{|P|}(\eps)$,
$\bar\sB^{s,\natural}_{|P|}(\eps)$
\hfill
Eq.
\eqref{eq:ExtendedFramedQuotientSpace}

$\bar\sB_\ka$
\hfill
Eq.
\eqref{eq:IdealConnS4}

$\bar{\bB\bL}^{\vir}_{\ft,\fs}$
\hfill
Eq. \eqref{eq:InstantonLinkComponentBase}

$\bar{\bB\bL}^{\vir}_{\ft,\fs}(\sP_j)$
\hfill
Eq. \eqref{eq:DefineLocalBaseLink}

$[\bar{\bB\bL}^{\vir}_{\ft,\fs}]$
\hfill
Eq. \eqref{eq:DefineBaseBundClass}

$\widehat{\bB\bL}^{\vir}_{\ft,\fs}(\sP_k)$
\hfill
Eq. \eqref{eq:PartialQuotient}

$\rd_{k_1}\dots\rd_{k_r}\bar{\bB\bL}^{\vir}_{\ft,\fs}(\sP_j)$
\hfill Eq. \eqref{eq:DefineBoundaryOne}

$c_P$
\hfill Eq. \eqref{eq:DefineR4DiagonalConePoint}

$\sC_{\fs}$, $\tsC_{\fs}$
\hfill
Eq.
\eqref{eq:SpincPreConfig}

$\tsC_{\ft}$
\hfill
Eq.
\eqref{eq:SpinUPreConfiguration}

$\sC_{\ft}$
\hfill
Eq.
\eqref{eq:SpinUPreConfiguration}

$\sC_{\ft}^*$, $\sC^0_{\ft}$, $\sC_{\ft}^{*,0}$
\hfill
following Eq. \eqref{eq:SpinUPreConfiguration}

$c_1(\fs)$
\hfill
Eq. \eqref{eq:DefineChernClassOfSpinc}

$c_1(\ft)$
\hfill
Eq.
\eqref{eq:SpinUCharacteristics}

$[c_{T,\beta}]$
\hfill Eq. \eqref{eq:CocycleOnRestrictedSpace}

$[c_\beta]$
\hfill Eq. \eqref{eq:DefineRelMuClass}

$[c_{T,W}]$
\hfill Eq. \eqref{eq:DefinePreCocycleForWClass0}

$[c_\sW]$
\hfill Eq. \eqref{eq:DefinePreCocycleForWClass}

$\bar c_\beta$
\hfill Lemma \ref{lem:ExtensionPerturbationV}

$[c(z,\eta)]$
\hfill Eq. \eqref{eq:DefineProductOfCocycles}

$\bar c_{\sW}$
\hfill Lemma \ref{lem:ExtensionPerturbationW}

$c_{\fs,\ell,j}$
\hfill Eq. \eqref{eq:ChernClassesOfSymSWUniv}

$c(\ft)$
\hfill Eq. \eqref{eq:SpinuCharacter}

$c_{j,i}$
\hfill Eq. \eqref{eq:FiberBoundaryProjectionMapToConn}

$\bar\bchi$, $\bchi_s$, $\bchi_i$
\hfill
Hypothesis \ref{hyp:Gluing}

$D(\sP,\eps)$, $\bar D(\sP,\eps)$, $\rd\bar D(\sP,\eps)$
\hfill Eq. \eqref{eq:DefineSquares}

$d(\ft),d_a(\ft)$
\hfill
Eq.
\eqref{eq:Transv}

$\Delta^\circ(X^\ell,\sP)$,
$\Delta(X^\ell,\sP)$
\hfill Eq. \eqref{eq:DefineDelta}

$\Delta^\circ(Z_P,\sP'_P)$
\hfill Eq. \eqref{eq:DefineCenteredR4Diagonal}

$\Delta^\circ(Z_\ka(\delta),\sP)$
\hfill Eq. \eqref{eq:DiagonalsOfZP}

$d_s(\fs)$
\hfill
Eq.
\eqref{eq:DimSW}

$D^w_X$
\hfill \S \ref{sec:Donaldsonseries}

$\bar e_I$
\hfill Eq. \eqref{eq:ExtendedInstEuler}

$\bar e(\bar\Upsilon_{\ft,\fs}/S^1,\bar\bchi')$
\hfill Eq. \eqref{eq:DefineExtendedRelativeEuler}

$\bar e_\bchi$
\hfill Lemma \ref{lem:ExtensionPerturbationChi}

$e(Z_P,\sP'_P)$
\hfill Eq. \eqref{eq:NormalOfR4Diagonal}

$e(X^\ell,g_{\sP})$
\hfill Eq. \eqref{eq:VaryingMetricExponentialMap}

$\FF_{\ft}$
\hfill
Eq.
\eqref{eq:UniversalU2Bundle}

$\FF_T$
\hfill Eq. \eqref{eq:TubNghBundle}

$\FF^{\vir}_{\ft,\fs}$
\hfill Eq. \eqref{eq:VirtualUniversal}

$\Fr(\ft,\fs,\sP)$
\hfill Eq. \eqref{eq:DefineGluingDataBundle}

$\bar\Fr(\ft,\fs,\sP_j)$
\hfill Eq. \eqref{eq:DefineExtendedBundle}

$\Fr(V,\sP,g_{\sP})$
\hfill Eq. \eqref{eq:FrameBundleForObstruction}

$\fg_{\ft}$
\hfill Eq.
\eqref{eq:SpinAssociatedBundles}

$\sG_{\fs}$
\hfill
\S \ref{subsec:SWMonopoles}

$\sG_{\ft}$
\hfill
\S \ref{subsubsec:PU2Monopoles}

$\bar\Gl(\ft,\fs,\sP)$
\hfill Eq. \eqref{eq:GluingData}

$\Gl(\ft,\fs,P)$
\hfill
Eq. \eqref{eq:DefineSingleInstantonSplicingDataBundle}

$\bar\Gl(\ft,\fs,\sP,[\sP'])$
\hfill Eq. \eqref{eq:OverlapGluingDataSpace}

$\bar\Gl(\ft,\fs,[\sP<\sP'])$
\hfill Eq. \eqref{eq:SimultUpwardsOverlapSpaces}

$\tilde G(T,\sP)$
\hfill Eq. \eqref{eq:TgBundleDiagonalStructureGrp}

$\tG(\sP)$, $G(\sP)$
\hfill Eq. \eqref{eq:DefineGluingDataBundleStructureGroup}

$g_{\sP}$, $g_{\sP,\bx}$
\hfill
Lemma \ref{lem:VaryingMetricExponentialMap}

$\bga'_{\Theta,\sP}$
\hfill Eqs. \eqref{eq:R4TrivialSplicing} \&
\eqref{eq:SplicedGenConn}

$\bga'_{\Theta,\sP,[\sP']}$
\hfill Eq. \eqref{eq:SplicingOverSeveralPartitions}

$\bga''_{\ft,\fs,\sP}$
\hfill Eqs. \eqref{eq:DefineCrudeSplicing} \& \eqref{eq:CrudeSplicingMapExplicit}

$\bga_{\ft,\fs,[\sP<\sP']}''$
\hfill Eq. \eqref{eq:UpperStratumCrudeSplicing}

$\bga''_{\sM}$
\hfill Eq. \eqref{eq:GlobalCrudeSplicingMap}

$\bga'_{\sM}$
\hfill Eq. \eqref{eq:DefineGlobalSplice}

$\bga_{\sM}$
\hfill Eq. \eqref{eq:GluingMap}

$\iota$
\hfill Eq. \eqref{eq:TopStratumOfVirtualInclusion}

$\iota_{\beta,1}$, $\iota_{\beta,2}$
\hfill Eq. \eqref{eq:PartialInclusions}

$\ind(\bD^*_{\ka})$
\hfill Eq. \eqref{eq:DefineInstantonDiracBundle}

$\sI(Y)$
\hfill Eq. \eqref{eq:DefineSingularSupportLocus}

$\bar\jmath_\beta$
\hfill Eq. \eqref{eq:DefinePairsInclusionForBeta}

$K_k$
\hfill
Lemma \ref{lem:CompactSubsetsOfSi}

$\rd_{k_1}\dots\rd_{k_r}K_j$
\hfill Eq. \eqref{eq:DefineBoundaryOfCompactum}

$\LL_{\fs}$
\hfill
Eq.
\eqref{eq:DefineSWUniversal}

$\LL_\nu$
\hfill Eq. \eqref{eq:LineBundleForS1ZAction}

$\bar\bL^w_{\ft,\ka}$
\hfill Eq. \eqref{eq:DefineASDLink}

$\bar\bL_{\ft,\fs}$
\hfill
Definition \ref{defn:DefineLink}

$\bar\bL^{\vir}_{\ft,\fs}$
\hfill Eq. \eqref{eq:AmbientLink}

$\bar\bL^{\vir,s}_{\ft,\fs}$
\hfill Eq. \eqref{eq:DefineSWComponentOfLink}

$\bar\bL^{\vir,i}_{\ft,\fs}(\sP_j)$
\hfill Eq. \eqref{eq:DefineLinkStratum}

$[\bar\bL^{\vir}_{\ft,\fs}]$
\hfill Eq. \eqref{eq:DefineAmbientLinkFund}

$\ell(\ft,\fs)$
\hfill Eq. \eqref{eq:ReducibleLevel}

$\la[A,\bx]$
\hfill Eq. \eqref{eq:ConnectionScale}

$\tilde\la_\ka$
\hfill Eq. \eqref{eq:TMScaleR4}

$\sM_{\ft}$
\hfill
before Eq.\eqref{eq:PU2MonopoleSubspaces}

$\bar\sM_{\ft}$
\hfill
before Eq.\eqref{eq:idealmonopoles}

$\bar\sM^{\red}_{\ft}$
\hfill
Eq.
\eqref{eq:StratificationCptPU(2)Space}

$\sM^{*,0}_{\ft}$
\hfill
Eq.
\eqref{eq:PU2MonopoleSubspaces}

$\bar\sM^{*,0}_{\ft}$
\hfill
following Eq.\eqref{eq:StratificationCptPU(2)Space}

$M^w_\ka$
\hfill
Eq.
\eqref{eq:ASDModuliSpace}

$M_{\fs}$
\hfill
\S \ref{subsec:SWMonopoles}

$\barM^{s,\natural}_{\spl,\ka}(\delta)$
\hfill
Theorem \ref{thm:ExistenceOfSplicedEndsModuli}

$\barM(\sP)$
\hfill Eq. \eqref{eq:DefineGluingDataFiber}

$\check M(\sP_k,\beps)$
\hfill
Definition \ref{eq:DefineFiberQuotient}

$\bar\sM^{\vir}_{\ft,\fs}$
\hfill Eq. \eqref{eq:DefineGlobalGluingDataSpace}

$\sM^{\vir}_{\ft,\fs}$
\hfill Eq. \eqref{eq:TopStratumOfVirtualInclusion}

$\bar\sM^{\vir,*}_{\ft,\fs}$
\hfill Eq. \eqref{eq:RedComplementInVirtual}

$\rd_{i_1}\rd_{i_2}\dots\rd_{i_v}\bar M(\sP_j,\beps)$
\hfill Eq. \eqref{eq:DefineLinkFiberBoundary}

$\mu_p$
\hfill
Eq.
\eqref{eq:MuPMap}

$\mu_{\fs}$
\hfill
Eq.
\eqref{eq:SWClass}

$\mu_c$
\hfill
Eq.
\eqref{eq:DefineMuC1}

$\barmu_c$
\hfill Eq. \eqref{eq:DefineExtendedMuC}

$\barmu_p(\beta)$
\hfill Definition \ref{defn:ExtendedCohomologyClasses}

$N_\ell$
\hfill
\S \ref{subsubsec:SymGrps}

$n_a(\ft)$
\hfill
Eq.
\eqref{eq:Transv}

$n_s(\ft,\fs)$
\hfill
Eq.
\eqref{eq:NormalComponentDims}

$N_{\ft,\fs}$
\hfill Eq. \eqref{eq:VirtualNormalandObstBundles}

$\nu$
\hfill Eq. \eqref{defn:DefineNu}

$\nu_{\Gl}$
\hfill Eq. \eqref{eq:XFactorS1Chern}

$\check\nu$
\hfill Eq. \eqref{eq:QuotientNu}

$\tilde\nu_i$
\hfill
Lemma \ref{lem:BranchedCoverS1Bundle}

$\nu_{\Gl}$
\hfill Eq. \eqref{eq:XFactorS1Chern}

$\tilde\nu(X^\ell,\sP)$
\hfill Eq. \eqref{eq:DiagonalNormal}

$\nu(X^\ell,\sP)$
\hfill
Lemma \ref{lem:Normal}

$\tilde\nu(X^\ell,\sP\to\sP')$
\hfill Eq. \eqref{eq:DoubleDiagonalNormal}

$\nu(X^\ell,\sP\to\sP')$
\hfill Eq. \eqref{eq:EndOfUpperStratum}

$\tilde\nu(X^\ell,\sP \to \sP')$
\hfill Eq. \eqref{eq:OverlapNormalBundle}

$\tilde\nu(X^\ell,\sP,[\sP'],g_{\sP})$
\hfill Eq. \eqref{eq:OverlapNormalBundleSymm}

$\nu(X^\ell,\sP\to [\sP'])$
\hfill Eq. \eqref{eq:EndOfUpperStratum}

$\tilde\nu(Z_{\ka}(\delta),\sP)$
\hfill Eq. \eqref{eq:R4DiagonalNormal}

$\tilde\nu(Z_{\ka}(\delta),\sP,\sP')$
\hfill Eq. \eqref{eq:DefineR4NormalPP'}

$\tilde\nu(Z_{\ka},[\sP<\sP'])$
\hfill Eq. \eqref{eq:R4DiagonalOverlap}

$\tilde\nu(\Theta,\sP,\sP',\delta)$
\hfill Eq. \eqref{eq:R4OverlapSpace}

$\nu(\Theta,\sP,[\sP'],\delta)$
\hfill Eq. \eqref{eq:SplicingOverlapSpaceR4SymmQuotient}

$\tilde\sO(X^\ell,\sP\to\sP',g_{\sP})$
\hfill
Lemma \ref{lem:DiagonalInDiagonalNormal}

$\tilde\sO(X^\ell,\sP,[\sP'],g_{\sP})$
\hfill Eq. \eqref{eq:TubNghBundleOverlaps}

$\sO(X^\ell,\sP,[\sP'],g_{\sP})$
\hfill Eq. \eqref{eq:TubNghBundleOverlaps}

$\sO(X^\ell,[\sP<\sP'],g_{\sP''})$
\hfill Eq. \eqref{eq:XlTubNghOverlapUpwardsSym}

$\tilde\sO(Z_{\ka}(\delta),\sP \to \sP')$,
\hfill Eq. \eqref{eq:DefineTubNghR4Normal}

$\sO(Z_P,\sP'_P)$
\hfill
Lemma \ref{lem:RelatingTrivialStrata}

$\tilde\sO(X^\ell,g_{\sP})$
\hfill
Lemma \ref{lem:VaryingMetricExponentialMap}

$\tilde\sO(Z_{\ka}(\delta),\sP)$
\hfill
Lemma \ref{lem:R4DiagonalsNormal}

$\sO(Z_\ka(\delta),[\sP<\sP'])$
\hfill
Lemma \ref{lem:SymmProdNeigh}

$\tilde\sO(\Theta,\sP,\delta)$
\hfill Eq. \eqref{eq:DefineSplicingDataR4}

$\sO(\Theta,\sP,\delta)$
\hfill Eq. \eqref{eq:R4TrivialSplicing}

$\tilde\sO(\Theta,\sP,\sP',\delta)$
\hfill Eqs. \eqref{eq:OverlapCondition1}
\& \eqref{eq:OverlapCondition4}

$\tilde\sO^{\asd}(\Theta,\sP,\delta)$,
$\sO^{\asd}(\Theta,\sP,\delta)$
\hfill Eq. \eqref{eq:DefineASDSplicingDomain}

$\tilde\sO(\Theta,\sP,[\sP'],\delta)$
\hfill Eq. \eqref{eq:SplicedModuliOverlapR4}

$\sO(\Theta,\sP,[\sP'],\delta)$
\hfill Eq. \eqref{eq:SplicedModuliOverlapR4}

$\sO^{\asd}(\Theta,[\sP<\sP'],\delta)$
\hfill Eq. \eqref{eq:UpperPartitionR4ASDNormal}

$\sO_d(\ft,\fs,\sP,[\sP'])$
\hfill Eq. \eqref{eq:XUpwardsTransitionDeomainRequirement}

$\Ob(\ft,\fs,\sP)$
\hfill Eq. \eqref{eq:InstantonObstructionBundle}

$\Ob(\ft,\fs,\sP,[\sP'])$
\hfill Eq. \eqref{eq:XInstantonObstrOverlapSpace}

$p_1(\ft)$
\hfill
Eq.
\eqref{eq:SpinUCharacteristics}

$\sP$ \hfill
\S \ref{subsubsec:SymGrps}

$[\sP]$ \hfill
Lemma \ref{lem:PreImageOfSi}

$\sP'_P$
\hfill
Eq.
\eqref{eq:DefineRefinementPartition}

$P_{\ft,\fs}$
\hfill
Lemma \ref{lem:CharClassOfGsTwistedBundle}

$[\sP<\sP']$
\hfill Eq.
\eqref{eq:ConjugateRefinements}

$\PD[\sV]$
\hfill Eq. \eqref{eq:DefinePDIncid}

$\pi(Z_P,\sP'_P)$
\hfill Eq. \eqref{eq:R4DiagonalProjection}

$\pi(X^\ell,g_{\sP})$
\hfill Eq. \eqref{eq:DefineXlTubularProjection}

$\pi(\Theta,\sP)$
\hfill Eq. \eqref{eq:SplicedEndR4Projection}

$\pi_N$
\hfill Eq. \eqref{eq:GlobalProjectionToN}

$\pi(\ft,\fs,\sP)$
\hfill Eq. \eqref{eq:DefineXTubularProj}

$\pi_X$
\hfill Eq. \eqref{eq:GlobalProjectionToX}

$\pi_{\fs}$
\hfill Eq. \eqref{eq:ProjectionToSW}

$\pi_{\fs,X}$
\hfill Eq. \eqref{eq:ProductOfProj1}

$\pi_{N,X}$
\hfill Eq. \eqref{eq:ProductOfProj2}

$\varphi_{\Theta,\sP}'$
\hfill Eq. \eqref{eq:DefineS4SectionSplicingToTrivial}

$\varphi_{\fs}''(\sP)$
\hfill Eq. \eqref{eq:BackObstrEmbedding}

$\varphi''_s$
\hfill Eq. \eqref{eq:DefineCrudeGlobalBackgroundObstrEmbedd}

$\varphi'_s$
\hfill Eq. \eqref{eq:DefineGlobalBackgroundObstrEmbedd}

$\varphi_{\ft,\fs,\sP}$
\hfill Eq. \eqref{eq:InstantonObstrSplicingMap}

$\varphi_{\ft,\fs,[\sP<\sP']}$
\hfill Eq. \eqref{eq:XCrudeInstantonSplicingMultiple}

$\varphi'_i$
\hfill Eq. \eqref{eq:XGlobalInstantonObstrSplicing}

$Q$, ${\bQ\bL}^{\vir}_{\ft,\fs}$
\hfill
Proposition \ref{prop:GlobalQuotient}

$Q_{k,i}$
\hfill
Definition \ref{defn:LocalQuotient}

$q_k$
\hfill Lemma \ref{lem:CharacterizeFiberQuotient}

$[{\bQ\bL}^{\vir}_{\ft,\fs}(\sP_i)]$
\hfill Eq. \eqref{eq:QuotientPiecesFund}

$\widetilde{\bQ\bL}^{\vir}_{\ft,\fs}(\sP_i)$
\hfill Lemma \ref{lem:CoverOfQuotientPiece}

$r_P$
\hfill Eq. \eqref{eq:ConeRetractionMap}

$r(X^\ell,\sP)$
\hfill
Lemma \ref{lem:XDefRetraction}

$r_i$
\hfill
Lemma \ref{lem:SymmetricProductRetractions}

$\bar r_x$
\hfill Eq. \eqref{eq:RestrictionForConfigSpace}

$\rho^{X,d}_{\sP,\sP'}$
\hfill Eq. \eqref{eq:DefineXlDownwardsOverlap}

$\rho^{X,u}_{\sP,\sP'}$
\hfill Eq. \eqref{eq:DefineXlUpwardsOverlap}

$\rho^d_{\sP,[\sP']}$
\hfill Eq. \eqref{eq:DefineXellDownwardsOverlapMap}

$\rho^u_{\sP,[\sP']}$
\hfill Eq. \eqref{eq:XlUpwardsOverlapMapSym}

$\rho^{\Theta,u}_{\sP,\sP'}$
\hfill Eq. \eqref{eq:R4TrivialSplicingUpwardsOverlapMap}

$\rho^{\Theta,d}_{\sP,\sP'}$
\hfill Eq. \eqref{eq:DefineUpwardTransition}

$\rho^{\Theta,u}_{\sP,[\sP']}$
\hfill Eq. \eqref{eq:SymmetricProductNormalBundlesOfHigherStratum}

$\rho^{\Psi,d}_{\sP,\sP'}$
\hfill Eq. \eqref{eq:DefineDownwardsObstrOverlapOnS4}

$\rho^{\Psi,u}_{\sP,\sP'}$
\hfill Eq. \eqref{eq:UpwardsTransitionObstruction}

$\rho^{\ft,\fs,u}_{\sP,[\sP']}$
\hfill Eq. \eqref{eq:XUpwardsTransitionMap}

$\rho^{\ft,\fs,d}_{\sP,[\sP']}$
\hfill Eq. \eqref{eq:XDownwardTransition}

$\rho^{\ft,\fs,d}_{f,\sP,[\sP']}$
\hfill Eq. \eqref{eq:DownwardsInclusionFiberMap}

$\rho^{\Xi,\ft,\fs,u}_{\sP,[\sP']}$,
$\rho^{\Xi,\ft,\fs,d}_{\sP,[\sP']}$
\hfill Eq. \eqref{eq:SWObstructionOverlap}

$\rho^{V,d}_{\sP,\sP'}$
\hfill Eq. \eqref{eq:XDownwardsInstantonOverlapMap}

$\rho^{V,u}_{\sP,[\sP']}$
\hfill Eq. \eqref{eq:XUpwardsInstantonOverlapMap}

$\bar r_{\beta}$
\hfill Eq. \eqref{eq:ExtendedRestrictionMap}

$R_{k,i}$
\hfill Lemma \ref{lem:PartialTranslMap}

$R_k$
\hfill Eq. \eqref{eq:RComposition}

$\fs, \fs\otimes L$
\hfill
\S \ref{subsubsec:SpincuStr}

$\Si(X^\ell,\sP)$
\hfill
Eq. \eqref{eq:DefineStratumOfPartition}

$\Si(Z_\ka(\delta),\sP)$
\hfill Eq. \eqref{eq:DiagonalsOfZP}

$\Si(\ft,\fs,\sP)$
\hfill Eq. \eqref{eq:DefineTrivialStratumXConePoints}

$\fS(\sP),\Ga(\sP)$, $W(\sP)$
\hfill
\S \ref{subsubsec:SymGrps}

$\SW_X(\fs)$
\hfill Eq. \eqref{eq:DefSW}

$S^\ell(\beta)$, $\tilde S^\ell(\beta)$
\hfill
Definition \ref{defn:DefnSymmBeta}

$\ft$
\hfill
\S \ref{subsubsec:SpincuStr}

$\ft(\ell)$
\hfill
Eq. \eqref{eq:DefineLowerChargeSpinuStr}

$T\Delta^\circ(X^\ell,\sP)$
\hfill Eq. \eqref{eq:TangentToDiagonal}

$T(\Theta,\sP,\delta)$
\hfill Eq. \eqref{eq:DefineTrivialStrataInR4GluingSpace}

$T(\Theta,\sP,\sP',\delta)$
\hfill Eq. \eqref{eq:DefineTrivialOverlap}

$T(\Theta,\sP,[\sP'],\delta)$
\hfill Eq. \eqref{eq:TrivialStrataOverlapSymmQuotient}

$T(\ft,\fs,\sP)$
\hfill Eq. \eqref{eq:DefineTrivialStratumX}

$T(\ft,\fs,\sP,[\sP'])$
\hfill Eq. \eqref{eq:DefineTrivialInOverlap}

$\bTheta_{\sP}$
\hfill Eq. \eqref{eq:CompatibleFlattening}

$t_{P,0}$
\hfill Eq. \eqref{eq:DefineInitialConeParameter}

$t(Z_P)$
\hfill
Lemma \ref{lem:R4TMScale}

$\vec t(X^\ell,g_{\sP})$
\hfill Eq. \eqref{eq:DefineDiagonalTubularDistFunction}

$\vec t(\ft,\fs,\sP)$, $\vec t^{\,P}(\ft,\fs,\sP)$
\hfill Eq. \eqref{eq:DefineTubularDistanceFunction}

$\vec t_f(\sP,[\sP'])$
\hfill Lemma \ref{lem:TubularDistanceFunctionOnOverlaps2}

$t_N$
\hfill Eq. \eqref{eq:NTubularDist}

$T_i$
\hfill
Lemma \ref{lem:DecomposeSymProduct}

$T_{k,j}$
\hfill
Lemma \ref{lem:CompactSubsetsOfSi}

$T(\sP_i,\sP_k)(A_0)$
\hfill Eq. \eqref{eq:EndBundleTriv}

$\sU(Z_P,\sP'_P)$
\hfill
Lemma \ref{lem:RelatingTrivialStrata}

$\tilde\sU(X^\ell,g_{\sP})$
\hfill
Lemma \ref{lem:VaryingMetricExponentialMap}

$\tilde\sU(Z_{\ka}(\delta),\sP)$
\hfill following Lemma \ref{lem:R4DiagonalsNormal}

$\sU(\Theta,\sP)$
\hfill
Lemma \ref{lem:TMProjectionOnSplicedEndR4}

$\bar\sU(\ft,\fs,\sP)$
\hfill Eq. \eqref{eq:DefineUSet}

$U_f(\sP,[\sP'])$
\hfill Eq. \eqref{eq:DefineOverlapFiberAsImage}

$\bar\Upsilon^i_{\spl,\ka}$
\hfill Eq. \eqref{eq:SplicedEndsIndexBundle}

$\tilde\Upsilon(\Theta,\sP,\sP')$
\hfill Eq. \eqref{eq:DefineOverlapInstantonObstr}

$\bar\Upsilon^s_{\ft,\fs}$
\hfill Eq. \eqref{eq:DefineGlobalBackgroundObstruction}

$\Upsilon^i_{X,P}/S^1$
\hfill Eq. \eqref{eq:XFactorObstruction}

$\Upsilon^i_{\ft,\fs}(P)/S^1$
\hfill Eq. \eqref{eq:FactorOfInstantonEuler}

$V$, $V^+$, $V^-$
\hfill
\S \ref{subsubsec:SpincuStr}

$\sV(z),\bar\sV(z)$
\hfill \S \ref{sec:Cohomology}

$\fV_\ka$
\hfill Eq. \eqref{eq:S4InfiniteDimObstrBundle}

$\bar\fV_\ka$
\hfill Eq. \eqref{eq:ExtendedS4InfiniteDimObstrBundle}

$\fV_{\ft}$
\hfill Eq. \eqref{eq:InfiniteDimObstr}

$\bar\fV_{\ft}$
\hfill Eq. \eqref{eq:ExtendedInfiniteDimObstrBundle}

$\sV_\ell(\Delta)$
\hfill Eq. \eqref{eq:Incidence}

$\tilde\sV_{\sP}(\Delta)$,
$\sV_{\sP}(\Delta)$
\hfill Eq. \eqref{eq:DefineIncidenceStratum}

$\sV_\ell(\Delta,T)$,
$\sV_{\sP}(\Delta,T)$
\hfill Eq. \eqref{eq:DefineIncidenceSubmanifold}

$\sW$, $\bar\sW$
\hfill
following Eq.\eqref{eq:DefineMuC1}

$W_\ka$
\hfill Eq. \eqref{eq:DefineSplicedEnd}

$w$ \hfill
before Eq. \eqref{eq:ASDModuliSpace}

$\Xi_{\ft,\fs}$
\hfill Eq. \eqref{eq:VirtualNormalandObstBundles}

$\Xi_\fs(\ft,\sP)$
\hfill Eq. \eqref{eq:DefineBackgroundObstrOverGluing}

$Z_P$
\hfill Eq. \eqref{eq:DefineZ}

$Z(\sP)$
\hfill Eq. \eqref{eq:DiagonalNormalFiber}

$Z(\sP,\sP')$
\hfill Eq. \eqref{eq:DoubleDiagNormalFiber}

$z[A,\bx]$
\hfill Eq. \eqref{eq:CenterOfMass}

\onecolumn

\begin{theindex}
\item Anti-self-dual connection, \pageref{eq:ASDModuliSpace}

\item Associativity of splicing maps, \pageref{subsec:IntroAssocSplic}

\item Circle action on space of $\SO(3)$ pairs, \pageref{Circle_action_on_SO(3)_pairs}, \pageref{subsec:S1Actions}

\item Clifford module, \pageref{subsubsec:SpincuStr}

\item Cohomology classes on configuration space of pairs, \pageref{sec:Cohomology}

\subitem Cocycles and their pullbacks, \pageref{subsubsec:CocyclesPullbacks}
\subsubitem Computation of cocycles, \pageref{subsubsec:ComputCocycles}
\subsubitem Extension of cocycles, \pageref{subsec:Extending_cocycles}
\subitem Cohomology classes on the virtual moduli space of $\SO(3)$ monopoles,
\pageref{subsec:DefiningCohomClasses}

\subitem Computation, \pageref{sec:ComputingMuClasses}
\subitem Dual geometric representatives, \pageref{sec:ComputingMuClasses}, \pageref{sec:Duality}

\item Computation of intersection numbers, \pageref{chap:Comp}

\item Configuration space of $\SO(3)$ pairs, \pageref{Configuration_space_SO(3)_pairs}

\item Connections over the four-dimensional sphere, \pageref{subsec:SpaceOfConn}
\subitem Center of mass, \pageref{eq:CenterOfMass}
\subitem Composition of splicing maps, \pageref{subsec:OverlapSpliceMapsR4}
\subsubitem Definition of the overlap data, \pageref{subsec:Definition_overlap data}
\subsubitem Equality of splicing maps, \pageref{subsec:Equality of splicing maps}
\subitem Scale, \pageref{eq:ConnectionScale}
\subitem Splicing onto the product connection, \pageref{subsec:SplicingMapsR4}
\subitem Strata containing the product connection, \pageref{subsec:R4Diagonals}

\item Crude splicing map, \pageref{intro_crude_splicing_maps}, \pageref{global_splicing_data_crude_splicing_map}, \pageref{sec:CrudeSplicingMap}, \pageref{subsubsec:ConstrCrudeSplice}
\subitem Properties, \pageref{subsec:Properties of the crude splicing map}

\item Dirac operator, \pageref{Dirac_operator}

\item Donaldson invariants, \pageref{sec:Donaldsonseries}

\item Donaldson series, \pageref{eq:DefineDonaldsonSeries}

\item Gauge transformations on space of $\SO(3)$ pairs, \pageref{Gauge_transformations_on_SO(3)_pairs}

\item Geometric representatives for cohomology classes, \pageref{geometric_representatives}
\subitem Cocycles, \pageref{subsec:GRepsAndCocycles}

\item Global splicing data space, \pageref{subsec:SpaceOfGlobalSplicingData}, \pageref{chap:GlobalSplicingData}, \pageref{sec:SplicingData}, \pageref{subsubsection:GrpActionOnFrames}, \pageref{sec:GlobalSplicingData}
\subitem Background pairs, \pageref{subsec:Background pairs}
\subitem Flattening map on pairs, \pageref{sec:FlatteningPairs}
\subitem Gluing data bundle, \pageref{eq:DefineGluingDataBundle}
\subitem Overlap space, \pageref{subsubsec:OverlapSpaceModuli}
\subsubitem Downwards overlap map, \pageref{subsec:Global_splicing_data_downwards overlap map}
\subsubitem Equality of splicing maps, \pageref{subsec:Global_splicing_data_quality_splicing_maps}
\subsubitem Upwards overlap map, \pageref{subsec:Global_splicing_data_upwards overlap map}
\subitem Projections onto symmetric products, \pageref{sec:ProjectionToSym}
\subitem Thom--Mather structures, \pageref{sec:TMStr}

\item Global splicing map, \pageref{sec:GlobalSplicingMap}

\item Instanton moduli space with spliced ends, \pageref{subsec:Introduction_spliced-ends_moduli_space}, \pageref{chap:SplicedEnd}, \pageref{subsec:SplicedEnd}
\subitem As a Whitney stratified space, \pageref{Instanton_moduli_space_spliced_end_Whitney_stratified_space}
\subitem Isotopy, \pageref{subsec:Collar}
\subitem Properties, \pageref{subsec:PropsOfSplicedEnd}
\subitem Tubular neighborhoods, \pageref{subsec:TMOfSplicedEnd}

\item Kotschick--Morgan conjecture, \pageref{chap:KMconj}, \pageref{Introduction_Kotschick-Morgan_conjecture},
\pageref{subsec:IntroKMConj}

\item Kronheimer--Mrowka basic classes, \pageref{Kronheimer-Mrowka_basic_classes}

\item Kronheimer--Mrowka simple type, \pageref{Kronheimer-Mrowka_simple_type}

\item Kronheimer--Mrowka structure theorem for Donaldson invariants, \pageref{Kronheimer-Mrowka_structure_theorem}

\item Link of ideal Seiberg--Witten moduli space, \pageref{chap:Link}, \pageref{subsec:Link}
\subitem Boundaries of components, \pageref{sec:Boundaries}
\subitem Instanton component, \pageref{eq:InstantonLinkComponent}
\subsubitem Fiber bundle structure, \pageref{sec:FiberBundle}
\subitem Orientation of link of ideal Seiberg--Witten moduli space, \pageref{subsubsec:Orient}
\subitem Virtual link of ideal Seiberg--Witten moduli space, \pageref{subsec:AmbientLink}

\item Locally flattened Riemannian metric, \pageref{defn:Locally_flattened_metric}, \pageref{sec:TMMapsOnSymX}

\item Moduli space
\subitem of ideal $\SO(3)$ monopoles, \pageref{ideal_SO(3)_monopoles}
\subitem of $\SO(3)$ monopoles, \pageref{Introduction_moduli_space_SO(3)_monopoles}, \pageref{moduli_space_SO(3)_monopoles}
\subitem of Seiberg--Witten monopoles, \pageref{Introduction_moduli_space_Seiberg-Witten_monopoles}

\item Obstruction bundle, \pageref{chap:obstr}
\subitem Background obstruction bundle, \pageref{Obstruction_bundle_Seiberg-Witten_moduli_space}, \pageref{sec:SWObstruction}
\subitem Equivariant Dirac index bundle, \pageref{sec:EquivDiracBundle}
\subitem Infinite-dimensional obstruction pseudo-bundle, \pageref{sec:InfiniteDimObstr}
\subitem Instanton obstruction pseudo-bundle, \pageref{sec:InstantonObstr}
\subitem Instanton component of the local obstruction bundle
\subsubitem Euler class, \pageref{subsec:InstObstrEulerLocal}
\subsubitem Frame bundles, \pageref{subsec:Instanton_obstruction_pseudo-bundle_frames}
\subitem Instanton component of the global obstruction bundle
\subsubitem Euler class, \pageref{subsec:InstObstrEulerGlobal}
\subitem Local gluing hypothesis for $\SO(3)$ monopoles, \pageref{sec:GluingThm}
\subitem Obstruction pseudo-bundle over instanton moduli space with spliced ends, \pageref{sec:ObstructionSplicing}
\subsubitem Overlap obstruction pseudo-bundle, \pageref{subsec:SplicedEndsObstrOverlap}
\subsubitem Overlap space and overlap maps, \pageref{subsec:IntroOverlapSpaceAndMaps}, \pageref{subsec:Instanton_obstruction_pseudo-bundle_overlap_space_maps}
\subsubitem Splicing map, \pageref{subsec:Instanton_obstruction_pseudo-bundle_splicing_map}
\subitem Relative Euler class of the obstruction bundle, \pageref{sec:RelEulerClassObstr}
\subitem Seiberg--Witten obstruction bundle, \pageref{Obstruction_bundle_Seiberg-Witten_moduli_space}
\subitem Seiberg--Witten component of the obstruction bundle
\subsubitem Euler class, \pageref{subsec:SWObstrEuler}

\item Problem of overlaps, \pageref{subsec:Overlaps}

\item Pidstrigach--Tyurin $\SO(3)$-monopole cobordism program, \pageref{Pidstrigach-Tyurin_SO(3)-monopole_cobordism}

\item Pseudo-bundle, \pageref{defn:PseudoBundle}

\item Partial tubular neighborhood structure, \pageref{partial_tubular_neighborhood_structure}

\subitem Vector-valued tubular distance function, \pageref{vector-valued_tubular_distance_function}

\item Perturbation parameters for $\SO(3)$ monopole equations, \pageref{perturbation_parameters}

\item Property P for knots, \pageref{Property_P_knots}

\item Reducible $\SO(3)$ monopoles, \pageref{subsec:RedPU2Monopole}

\item Relative Euler classes, \pageref{subsec:RelEulerClasses}

\item Seiberg--Witten invariant, \pageref{subsec:SWInvariants}

\item Seiberg--Witten monopole, \pageref{subsec:SWMonopoles}
\subitem Seiberg--Witten monopole equations, \pageref{eq:SeibergWitten}

\item Smoothly stratified space, \pageref{smoothly_stratified_space}
\subitem Compatibility condition, \pageref{compatibility_condition}
\subitem Condition of the frontier, \pageref{condition_frontier}
\subitem Control condition, \pageref{control_condition}
\subitem Thom--Mather stratification, \pageref{Thom-Mather_stratification}
\subitem Whitney stratification, \pageref{Whitney_stratification}

\item $\SO(3)$ monopole, \pageref{SO(3)_monopole}
\subitem $\SO(3)$ monopole equations, \pageref{eq:PerturbedSO3MonopoleEquations}

\item $\SO(3)$-monopole cobordism formula, \pageref{thm:MainThm}
\subitem As equality of intersection pairings with links, \pageref{eq:RawCobordismSum}, \pageref{thm:CobordismThm}

\item \spinc structure, \pageref{Introduction_spinc_structure}, \pageref{spinc_structure}

\item \spinu structure, \pageref{Introduction_spinu_structure}, \pageref{spinu_structure}

\item Standard splicing map, \pageref{subsubsec:StdSplicingMap}

\item Stratum of anti-self-dual (zero-section) $\SO(3)$-monopoles, \pageref{sec:ASDsingularities}

\item Stratum of Seiberg--Witten (reducible) $\SO(3)$-monopoles, \pageref{sec:Reducibles}

\item Superconformal simple type conjecture, \pageref{Superconformal_simple_type_conjecture}

\item Symmetric group, \pageref{subsubsec:SymGrps}
\subitem Partition, \pageref{subsubsec:SymGrps}
\subitem Partition length, \pageref{subsubsec:SymGrps}
\subitem Partition refinement, \pageref{defn:DefineRefinement}

\item Symmetric products of manifolds, \pageref{chap:Diagonals}
\subitem Big diagonal, \pageref{subsec:DefiningTheDiagonals}
\subitem Center of mass map, \pageref{eq:R4CenterOfMassProjection}
\subitem Cone parameter function, \pageref{eq:DefineInitialConeParameter}
\subitem Cone point, \pageref{eq:DefineR4DiagonalConePoint}
\subitem Diagonal of symmetric product of manifolds, \pageref{subsec:DefiningTheDiagonals}
\subitem Incidence relations among diagonals and strata, \pageref{subsec:IncidenceRelations}
\subitem Normal bundle of diagonal, \pageref{sec:DiagonalNormals}
\subitem Normal bundle of stratum, \pageref{lem:NormalInSymmProd}, \pageref{sec:NormalBundlesOfSymlXStrata}
\subitem Overlap maps, \pageref{sec:OverlapMaps}
\subsubitem Commuting overlap maps, \pageref{subsec:Commuting_overlap_map}
\subsubitem Downwards overlap map, \pageref{eq:DefineXlDownwardsOverlap}
\subsubitem Upwards overlap map, \pageref{eq:DefineXlUpwardsOverlap}
\subitem Partial ordering of strata, \pageref{subsec:EnumStrata}
\subitem Small diagonal, \pageref{subsec:DefiningTheDiagonals}
\subitem Stratum of symmetric product defined by partition, \pageref{eq:DefineStratumOfPartition}
\subitem Tubular distance function, \pageref{tubular_distance_function}, \pageref{sec:SymlXTubDistance}
\subitem Tubular neighborhood projection, \pageref{tubular_neighborhood_projection}
\subitem Tubular neighborhood structure, \pageref{tubular_neighborhood_structure}

\item Transversality for the moduli space of  $\SO(3)$ monopoles, \pageref{thm:Transv}

\item Uhlenbeck compactification for the moduli space of $\SO(3)$ monopoles, \pageref{Introduction_Uhlenbeck_compactification}, \pageref{thm:Compactness}

\item Uhlenbeck topology for the moduli space of $\SO(3)$ monopoles, \pageref{Uhlenbeck_topology}

\item Universal Seiberg--Witten line bundle, \pageref{eq:DefineSWUniversal}

\item Virtual fundamental class of moduli space of $\SO(3)$ monopoles, \pageref{Virtual_fundamental_class}

\item Virtual link of the ideal moduli space of Seiberg--Witten monopoles,
\subitem Fundamental class, \pageref{sec:FundClassOfAmbLink}

\item Virtual moduli space of Seiberg--Witten monopoles, \pageref{Virtual_Seiberg-Witten_moduli_space}

\item Virtual normal bundle of Seiberg--Witten moduli space, \pageref{subsubsec:ThickenedNeighborhood}

\item Witten's conjecture, \pageref{conj:WC}
\end{theindex}

\bibliography{/Users/pfeehan/Dropbox/LATEX/Bibinputs/master,/Users/pfeehan/Dropbox/LATEX/Bibinputs/mfpde}
\bibliographystyle{amsplain-nodash}

\end{document}